\documentclass[11pt]{amsart}

\usepackage{parskip}
\allowdisplaybreaks

\usepackage{geometry} 
\usepackage{framed} 

\usepackage{tikz}
\tikzset{every picture/.style={line width=1pt}} 

\usepackage{setspace}
\setdisplayskipstretch{2.5}

\usepackage{amsthm} 
\usepackage{amsmath}
\usepackage{amssymb}
\usepackage{mathrsfs}

\usepackage{hyperref}

\usepackage{amsfonts} 
\newcommand\A{\mathcal{A}}
\newcommand\B{\mathcal{B}}

\newcommand\Bfrak{\mathfrak{B}}
\newcommand\C{\mathbb{C}}
\newcommand\CC{\mathcal{C}}

\newcommand\BC{\mathbf{C}}
\newcommand\Cbf{\mathbf{C}}
\newcommand\BD{\mathbf{D}}
\newcommand\Dbf{\mathbf{D}}

\newcommand\G{\mathcal{G}}
\renewcommand\H{\mathcal{H}}

\newcommand\N{\mathbb{N}}

\newcommand\R{\mathbb{R}}
\renewcommand\S{\mathcal{S}}
\newcommand\sfrak{\mathfrak{s}}
\newcommand\Z{\mathbb{Z}}

\newcommand\one{\mathbf{1}}

\newcommand\w{\omega}
\newcommand\vphi{\varphi}
\renewcommand\phi{\vphi}
\newcommand\eps{\varepsilon}

\usepackage{stmaryrd} 

\newcommand\id{\textnormal{id}}
\newcommand\sing{\textnormal{sing}}
\newcommand\reg{\textnormal{reg}}
\newcommand\spt{\textnormal{spt}}
\newcommand\Lip{\textnormal{Lip}}
\newcommand\dist{\textnormal{dist}}
\newcommand\graph{\textnormal{graph}}
\newcommand\ext{\mathrm{d}}
\newcommand\del{\partial}
\newcommand\vartan{\textnormal{VarTan}}

\newcommand{\res}{\mathbin{\hspace{0.1em}\vrule height 1.3ex depth 0pt width 0.13ex\vrule height 0.13ex depth 0pt width 1.0ex}} 

\newcommand{\weakly}{\rightharpoonup}
\renewcommand{\div}{\textnormal{div}}

\makeatletter
\def\@tocline#1#2#3#4#5#6#7{\relax
  \ifnum #1>\c@tocdepth 
  \else
    \par \addpenalty\@secpenalty\addvspace{#2}%
    \begingroup \hyphenpenalty\@M
    \@ifempty{#4}{%
      \@tempdima\csname r@tocindent\number#1\endcsname\relax
    }{%
      \@tempdima#4\relax
    }%
    \parindent\z@ \leftskip#3\relax \advance\leftskip\@tempdima\relax
    \rightskip\@pnumwidth plus4em \parfillskip-\@pnumwidth
    #5\leavevmode\hskip-\@tempdima
      \ifcase #1
       \or\or \hskip 1em \or \hskip 2em \else \hskip 3em \fi%
      #6\nobreak\relax
    \dotfill\hbox to\@pnumwidth{\@tocpagenum{#7}}\par
    \nobreak
    \endgroup
  \fi}
\makeatother

\newtheoremstyle{newtheoremstyle}
{3pt}
{3pt}
{\itshape}
{\parindent}
{\bfseries}
{.}
{0.5em}
{} 
\newtheoremstyle{newtheoremstyledefn}
{3pt}
{3pt}
{}
{\parindent}
{\bfseries}
{.}
{0.5em}
{} 

\theoremstyle{newtheoremstyle}
\newtheorem{theorem}{Theorem}
\newtheorem*{theorem*}{Theorem}
\newtheorem{lemma}[theorem]{Lemma}
\newtheorem{prop}[theorem]{Proposition}
\newtheorem{corollary}[theorem]{Corollary}

\newtheorem{thmx}{Theorem}

\theoremstyle{newtheoremstyledefn}
\newtheorem{defn}[theorem]{Definition}
\newtheorem{example}[theorem]{Example}
\newtheorem{remark}[theorem]{Remark}

\newtheorem{remarkx}{Remark}

\numberwithin{equation}{section} 
\numberwithin{theorem}{section}

\usepackage{fancyhdr}
\begin{document}

\title{On the Nature of Stationary Integral Varifolds near Multiplicity 2 Planes}

\author{
	Spencer Becker-Kahn
	\and
	Paul Minter
	\and
	Neshan Wickramasekera
}

\address{\textnormal{St Paul’s Girls’ School, Brook Green, Hammersmith, London
W6 7BS}}
\email{spencer.becker-kahn@spgs.org}
\address{\textnormal{Department of Pure Mathematics and Mathematical Statistics, University of Cambridge}\\
and\\
\textnormal{Department of Mathematics, Stanford University, Building 380, Stanford, CA 94305, USA}}
\email{pdtwm2@cam.ac.uk; pminter@stanford.edu}
\address{\textnormal{Department of Pure Mathematics and Mathematical Statistics, University of Cambridge}}
\email{N.Wickramasekera@dpmms.cam.ac.uk}

\begin{abstract}
	We study stationary integral $n$-varifolds $V$ in the unit ball $B_1(0)\subset\R^{n+k}$. Allard's regularity theorem establishes the existence of $\eps = \eps(n,k)\in (0,1)$ for which if $V$ is $\eps$-close in the varifold sense to the plane $P_0 = \{0\}^k\times\R^n$ with multiplicity 1 then, in $B_{1/2}(0)$, $V$ is represented by a single $C^{1,\alpha}$ minimal graph. However, when instead $P_0$ occurs with multiplicity $Q\in \{2,3,\dotsc\}$, simple examples show that this conclusion, now as a multi-valued graph, may fail, even if $V$ corresponds to an area-minimising rectifiable current.
	
	In the present work we investigate the structure of such $V$ which are close to planes with multiplicity $Q>1$, focusing primarily on the case $Q=2$. We show that an $\eps$-regularity theorem holds when $V$ is close, as a varifold, to $P_0$ with multiplicity $2$, \emph{provided} $V$ satisfies a certain topological structural condition on the part of its support where the density of $V$ is $<2$. The conclusion then is that, in $B_{1/2}(0)$, $V$ is represented by the graph of a Lipschitz $2$-valued function over $P_0$ with small Lipschitz constant; in fact, the function is $C^{1,\alpha}$ in a precise generalised sense, and satisfies estimates, implying that all tangent cones at singular points in $B_{1/2}(0)$ are unique and comprised of stationary unions of $4$ half-planes (which may form a union of two distinct planes or a single multiplicity $2$ plane). The theorem does not require any additional assumption on the part of $V$ with density $\geq 2$ (which a priori may be a relatively large set in $\H^n$-measure, possibly with high topological complexity).
	
	The topological structural condition is a requirement on $V$ in so-called $\beta$\emph{-coarse gaps}, which are vertical round cylinders over $P_0$ in which the density of $V$ is $<2$ and, at the scale of the cylinder, $V$ is $\beta$-close as a varifold to some multiplicity $2$ plane parallel to $P_0$. We say that $V$ satisfies the topological structural condition for some $\beta\in(0,1)$ if in any $\beta$-coarse gap on the corresponding sub-cylinder with half the radius $V$ is represented by two (disjoint) minimal graphs. Our $\eps$-regularity theorem then applies to $V$ satisfying the topological structural condition for some choice of $\beta$ (with $\eps$ depending only on $n$, $k$ and $\beta$).
	
	Consequently, if $\mathcal{V}$ is a class of stationary integral $n$-varifolds closed under ambient homotheties, translations, and rotations, one may view our $\eps$-regularity theorem as saying that in order to establish an $\eps$-regularity theorem for varifolds in $\mathcal{V}$ close to a multiplicity $2$ plane, it suffices to establish the $\eps$-regularity theorem assuming additionally that the density of the varifold is $<2$ everywhere. 
    
    As a corollary, we show that our $\eps$-regularity theorem applies unconditionally to stationary $2$-valued Lipschitz graphs with arbitrary Lipschitz constant, yielding improved regularity and uniform a priori estimates. In particular, no other types of singularity, such as those modelled on Lawson--Osserman cones, can accumulate near branch points of stationary $2$-valued Lipschitz graphs. We also recover known regularity results near density $2$ branch points $Z$ of codimension $1$ stationary integral varifolds $W$ with no triple junction singularities and having finite Morse index (of the regular part of $W$) near $Z$. Moreover, our results give a topological property equivalent to $W$ having infinite Morse index in every neighbourhood of $Z$, analogous to locally infinite genus in the 2-dimensional case. In general, our result applies locally about any density 2 branch point of an \emph{arbitrary} stationary integral varifold, provided the branch point is not a limit point of triple junction singularities or multiplicity one regular points where the tangent plane is vertical relative to a given tangent plane at the point.
	
    Our $\eps$-regularity theorem follows from a more general excess decay \emph{tri}chotomy for arbitrary stationary integral varifolds close to a multiplicity $2$ plane, which is reminiscent of the main excess decay dichotomy in the work of L.~Simon on stationary varifolds in multiplicity 1 classes satisfying an integrability condition.
    We also give analogous results when the varifold is instead close to a union of two (distinct) multiplicity 1 planes or a stationary union of $4$ (distinct) half-planes. All the results of the paper extend to integral varifolds with generalised mean curvature in $L^p$ for some $p>n$, as well as to Riemannian ambient spaces.
\end{abstract} 

\maketitle

\tableofcontents

\part{\centering Introduction}\label{part:intro}

Let $V$ be a stationary integral $n$-varifold in the open unit ball $B_1(0)\subset\R^{n+k}$ (or, more generally, a smooth Riemannian manifold). It has long been known from the seminal work of Allard \cite{All72} that close to a point in $\spt\|V\|$ where one tangent cone is a plane with multiplicity $1$, $V$ corresponds to the graph of a single $C^{1,\alpha}$ function with controlled $C^{1,\alpha}$ norm. Elliptic PDE theory can then be invoked to improve on this regularity and give further estimates. In particular, such points are regular points of $V$, namely they have a neighbourhood in which $V$ is represented by a smoothly embedded submanifold.

In the general setting, namely concerning the nature of a stationary integral varifold near a non-immersed (singular) point where at least one tangent cone is a plane of (necessarily integer) multiplicity $Q\geq 2$, prior to the present work, as far as we are aware, little has been known beyond Almgren's Lipschitz approximation theorem \cite[Corollary 3.10]{Alm00}. Understanding such singular points, which we shall refer to as \emph{branch points}, is arguably the most basic open question in the regularity theory of stationary varifolds. This question is paramount to understanding the size and structure of the singular set; for instance, it is not known whether it is possible for the $\H^n$-measure of the singular set to be positive (and throughout this work we allow for this possibility). Almgren's Lipschitz approximation theorem implies that, near such a point, $V$ agrees with a $Q$-valued Lipschitz graph over the plane outside a set of small measure. One of the questions we investigate here is when one can guarantee that, near a branch point, $V$ is \emph{exactly} equal to a $Q$-valued Lipschitz graph over its tangent plane and, if this is the case, whether the graph has any higher regularity and quantitative estimates (for instance, control of a $C^{1,\alpha}$-nature, which would in particular give uniqueness of the tangent plane). Our main result  will give an answer to this question when $Q=2$.

One should first observe that the analogous $\eps$-regularity theorem to that in Allard's work \cite{All72} \emph{fails} when $Q=2$: simple examples, such as the standard $1$-parameter family of catenoids in $\R^3$, illustrate this (cf.~ Example \ref{ex:catenoids}). Moreover, even the weaker claim described above, i.e.~a (possibly multi-valued) Lipschitz graphical structure around a singular point where one tangent cone is a plane of multiplicity $2$, fails (cf.~Example~\ref{ex:graph}).\footnote{In the case of bounded mean curvature this failure is more dramatic, as can be seen from a well-known example of Brakke \cite[Section 6.1]{Bra78} (which has been quantified and further strengthened in the work of Kolasiński--Menne \cite[Example 10.3]{KM17}). We give a brief description of Brakke's example in Example \ref{ex:brakke}; it has a singular set with positive $\H^n$-measure and illustrates a key obstruction to regularity in that every singular point is a limit point of neck-like regions modelled on the catenoid. We emphasise that whilst one can guarantee, by scaling, that the mean curvature in Brakke's example is arbitrarily small (in $L^\infty$), it is not known whether there is a version of Brakke's example where the mean curvature vanishes everywhere; that is, where the surface is minimal.}

If however one additionally assumes that the density $\Theta_V$ of $V$ satisfies $\Theta_V\geq 2$ almost everywhere, then one can apply Allard's regularity theorem to $V$ when it is close to a multiplicity $2$ plane, by considering the (stationary integral) varifold with the same support as $V$ but with half the density. This and the example of the catenoid suggest that to establish an $\eps$-regularity theorem for stationary integral varifolds near a multiplicity $2$ plane, an appropriate additional condition just on the (relatively open) part $\{\Theta_{V} < 2\}$ may suffice. The main result of this work (Theorem~\ref{thm:main-intro} and Theorem~\ref{thm:main}) will provide such an assumption; indeed, no additional assumption is needed on the set $\{\Theta_V\geq 2\}$, which a priori may be a relatively large set in $\H^n$-measure, possibly with high topological complexity. More generally, our main technical lemma (Theorem \ref{thm:main-gap-intro}) provides a trichotomy which must be satisfied by \emph{any} stationary integral varifold sufficiently close to a plane of multiplicity $2$.

Let us start by describing our main $\eps$-regularity theorem and the assumption needed on the region $\{\Theta_V<2\}$. For this, fix a parameter $\beta\in (0,1)$. We call a cylinder $C_\rho(x):= \R^k\times B^n_\rho(x)$ a \emph{$\beta$-coarse gap} (for $V$) if it satisfies the following conditions:
\begin{itemize}
	\item $C_\rho(x)\subset \{\Theta_V<2\}$;
	\item $\frac{3}{2}\leq (\w_n\rho^n)^{-1}\|V\|(C_\rho(x)) \leq \frac{5}{2}$;
	\item For some $z\in \R^k$ we have
	$$\frac{1}{\rho^{n+2}}\int_{C_\rho(x)}\dist^2(y, \{z\}\times\R^n)\, \ext\|V\|(y) < \beta^2.$$
\end{itemize}
Thus, a $\beta$-coarse gap is one where, at the scale of the cylinder $C_\rho(x)$, $V$ is $\beta$-close (in $L^2$) to a plane parallel to $\{0\}^k\times\R^n$ with multiplicity $2$ yet the density of $V$ is $<2$ everywhere in $C_\rho(x)$. Knowing that the density of $V$ is $<2$ heavily restricts the possible singularities which can occur in $V$; this will be discussed further in Remark \ref{remark:b}.

The assumption we will need to make on the set $\{\Theta_V<2\}$ is that the varifold satisfies a separation property in $\beta$-coarse gaps, described as follows:
\begin{defn}\label{defn:top-str-con}
	Suppose $C_\rho(x)$ is a $\beta$-coarse gap for $V$. We say $V$ satisfies the \emph{topological structural condition in $C_\rho(x)$} if $\spt\|V\|\cap C_{\rho/2}(x)$ has at least two connected components, and furthermore both of these connected components intersect $C_{\rho/4}(x)$.
\end{defn}

\begin{remark}\label{remark:sheeting-intro}
We will see in Proposition \ref{prop:sheeting-equivalence} that $V$ satisfying the topological structural condition in $\beta$-coarse gaps $C_\rho(x)$ is equivalent (up to changing $\beta$ by a controlled amount) to requiring that $V\res C_{\rho/2}(x)$ is the sum of two minimal graphs; this in particular guarantees that $V\res C_{\rho/2}(x)$ is topologically simple as, topologically, it is the union of two (disjoint) disks. Notice that this is not the case for the catenoid, where this intersection would be an annulus.
\end{remark}

\begin{remark}
    In Remark \ref{remark:g} we will see that (again, up to changing $\beta$ by a controlled amount) that $V$ satisfying the topological structural condition in a $\beta$-coarse gap is also equivalent to $V$ not containing any: (i) triple junction singularities; and (ii) multiplicity one regular points where the tangent plane is ``vertical'' relative to the fixed reference plane ($\{0\}^k\times\R^n$ in the above).
\end{remark}

Our main $\eps$-regularity result can then be informally described as follows (see Theorem \ref{thm:main} for the precise statement).

\begin{thmx}\label{thm:main-intro}
	Fix $\beta>0$. Then, there exists $\eps = \eps(n,k,\beta)\in (0,1)$ such that the following is true. Suppose that $V$ is a stationary integral $n$-varifold in $B_2(0)\subset\R^{n+k}$ such that in the cylinder $\R^k\times B^n_1(0),$ $V$ is $\eps$-close, in the sense of varifolds, to the disk $\{0\}^k\times B^n_1(0)$ with multiplicity $2$, i.e.
	\begin{itemize}
		\item $\frac{3}{2}\leq \w_n^{-1}\|V\|(B_1(0))\leq \frac{5}{2}$;
		\item $\hat{E}_V<\eps$, where
		$$\hat{E}_V^2 := \int_{\R^k\times B^n_1(0)}\dist^2(x,\{0\}^k\times\R^n)\, \ext\|V\|(x).$$
	\end{itemize}
	Then, provided $\widetilde{V} := V\res B_{3/2}(0)$ satisfies the topological structural condition in any $\beta$-coarse gap, we have:
	\begin{enumerate}
		\item [(A)] $\widetilde{V}\res (\R^k\times B^n_{1/2}(0))$ is equal to (the varifold associated to) the graph of a $2$-valued Lipschitz function $f$ over the disk $\{0\}^k\times B^n_{1/2}(0)$. Moreover, $f$ is $C^{1,\alpha}$ in a certain generalised sense (see Definition \ref{defn:gc-1}), with quantitative estimates, including
		$$\sup_{B^n_{1/2}(0)}|f| + \Lip(f) \leq C\hat{E}_V;$$
		\item [(B)] At every singular point of $\widetilde{V}$ in $\R^k\times B^n_{1/2}(0)$, the tangent cone is unique and comprised of $4$ half-planes (which includes the possibility of the tangent cone being the sum of two multiplicity $1$ planes or a single multiplicity $2$ plane);
		\item [(C)] The tangent cones along the singular set of $\widetilde{V}$ vary $\alpha$-Hölder continuously.
	\end{enumerate}
	Here, $C = C(n,k)\in (0,\infty)$ and $\alpha = \alpha(n,k)\in (0,1)$.
\end{thmx}
In fact, one readily checks that the (unique) tangent cones in conclusion (B) must be one of three types:
\begin{enumerate}
	\item [(a)] a plane with multiplicity $2$;
	\item [(b)] the union of two multiplicity $1$ planes which intersect along a subspace of dimension $\in \{0,1,2,\dotsc,n-1\}$;
	\item [(c)] a \emph{twisted cone}, i.e.~a stationary union of four (distinct) multiplicity $1$ half-planes all with a common boundary (which is a subspace of dimension $n-1$), but which do not form a pair of planes.
\end{enumerate}
If there are no twisted cones arising as tangent cones, then in fact the above `generalised' notion of $C^{1,\alpha}$ regularity for the function $f$ reduces to being $C^{1,\alpha}$ in the usual sense of $2$-valued functions. In particular, it is simple to check that in the special case $k=1$ (i.e.~the codimension is $1$) there are no stationary twisted cones and so in this case the conclusion in (A) becomes that $V\res (\R^k\times B^n_{1/2}(0))$ is equal to the varifold associated to a $2$-valued $C^{1,\alpha}$ graph. The same is true when $V$ corresponds to an area minimising current. Improved regularity properties for $C^{1,\alpha}$ $2$-valued graphs have been previously studied in \cite{SW16, KW21}, therefore providing further information.

In light of Remark \ref{remark:sheeting-intro}, one illustrative viewpoint which the reader may take for Theorem \ref{thm:main-intro} is the following: the validity of the topological structural condition in $\beta$-coarse gaps is morally saying that one already knows an $\eps$-regularity theorem when $V$ is both close to a multiplicity 2 plane \emph{and additionally has $\Theta_V<2$ everywhere}. The density being $<2$ restricts the types of singularities in $V$, meaning this special case should be strictly easier to establish, and indeed in certain situations we know this is the case (cf.~Theorem \ref{thm:main-2-intro} and Remark \ref{remark:b}). Thus, Theorem \ref{thm:main-intro} says that if $\mathcal{V}$ is a class of stationary integral varifolds closed under homotheties and rotations for which one has already established the special case of Theorem \ref{thm:main-intro} when $\Theta_V<2$ everywhere, \emph{then} one can subsequently establish the general case of the theorem where there is \emph{no restriction} on the density $\Theta_V$.

When $V$ corresponds to the graph of a (stationary) $2$-valued Lipschitz function over some plane (with arbitrary Lipschitz constant), one can directly check that the topological structural condition is satisfied in $\beta$-coarse gaps (for \emph{any} choice of $\beta$). Therefore, Theorem \ref{thm:main-intro} directly applies in this situation, providing improved regularity and quantitative estimates for stationary $2$-valued Lipschitz graphs with small $L^2$-norm. This corollary is summarised in the following (cf.~Theorem \ref{thm:main-2} for the precise statement, in particular Section \ref{sec:q-valued-notation} for the definition of $\mathbf{v}(f)$).

\begin{thmx}\label{thm:main-2-intro}
	Fix $L\in (0,\infty)$. Suppose that $f:B^n_1(0)\to \A_2(\R^k)$ is a Lipschitz $2$-valued function with $\Lip(f)\leq L$, such that the (integral) varifold $V := \mathbf{v}(f)$ associated to the graph of $f$ is stationary in $\R^{k} \times B_{1}^{n}(0)$. Then, for any $\beta\in (0,1)$, $V$ satisfies the topological structural condition in any $\beta$-coarse gap.
	
	In particular, Theorem \ref{thm:main-intro} applies to $V$ (with e.g.~$\beta=1/2$). Hence, there exists $\eps = \eps(n,k)\in (0,1)$ such that if $\|f\|_{L^2(B^n_1(0))}<\eps$, then $f|_{B^n_{1/2}(0)}$ is $C^{1,\alpha}$ in a generalised sense (as in Definition \ref{defn:gc-1})), with quantitative estimates. In particular, for some $C = C(n,k)$ we have
	$$\sup_{B^n_{1/2}(0)}|f| + \Lip(\left.f\right|_{B^n_{1/2}(0)}) \leq C\|f\|_{L^2(B_1^n(0)).}$$
\end{thmx}

A direct consequence of Theorem \ref{thm:main-2-intro} is that the only types of singularity which can limit onto a (density $2$) branch point in a stationary $2$-valued Lipschitz graph are those modelled on unions of $4$ half-planes. In particular, more complicated singularities such as those modelled on the Lawson--Osserman cones \cite{LO77} cannot converge to such branch points.

\begin{remark}\label{remark:special-case-stable}
	Another situation where one can verify the topological structural condition in $\beta$-coarse gaps (this time for suitably small $\beta = \beta(n)\in (0,1)$) is when $k=1$ and $V$ corresponds to a stationary integral varifold with stable regular part and no triple junction singularities (the latter condition holds, for example, if $V$ is associated to a current with vanishing boundary). In fact, this particular case of Theorem \ref{thm:main-intro} has already been established by the second and third authors \cite[Theorem D]{MW24}. The verification of the topological structural condition follows directly from the Sheeting Theorem of \cite{Wic14}, or from combining singularity analysis in the region $\{\Theta_V<2\}$ with the regularity theorem of Schoen--Simon \cite{SS81} (see also \cite{Bel25}).
\end{remark}

\begin{remark}\label{remark:area-min}
    In \cite{Cha88}, S.~Chang showed that the local structure about (necessarily isolated) singular points of a 2-dimensional area minimising current is given by a finite union of $C^{3,\beta}$ branched disks. Arguably one of the key aspects of Chang's work is showing that one can verify that a suitable version of our topological structural condition holds on a neighbourhood of a branch point; this is achieved (more or less) by arguing that one can produce suitable tangent maps which are non-trivial, average-free, homogeneous multi-valued Dirichlet energy minimisers. Such functions obey a separation property away from the singular point, which can be passed back down to the level of the current and used to deduce the local structure (this is essentially what happens about density $2$ branch points, and is achieved about higher density branch points using a more complicated blow-up procedure relative to a \emph{branched center manifold}). We will explain this further in Remark \ref{remark:f}, where we will also explain how the same principle applies to area minimisers mod $4$ to determine their local structure about branch points (which necessarily have density $2$ in this case). It therefore appears that all known regularity results which determine the local structure about density $2$ branch points can be explained by the verification of the topological structural condition in $\beta$-coarse gaps (which in particular have small $L^2$ excess to a plane).
\end{remark}

We stress that conclusion (B) of Theorem \ref{thm:main-intro} implies, in this setting, uniqueness of multiplicity $2$ tangent planes arising as tangent cones, and so in particular for Lipschitz $2$-valued stationary graphs. Uniqueness of multiplicity $\geq 2$ tangent planes remains largely an open question for \emph{any} class of integral varifolds with generalised mean curvature in $L^p$ (when $p>n$), and in particular for stationary integral varifolds. The work here appears to be the first instance of showing such a uniqueness result everywhere, in arbitrary dimension and codimension, with the \emph{only} variational assumption being control on the generalised mean curvature (or, equivalently, the first variation of area). We remark that for area minimising and semi-calibrated integral currents some significant results concerning uniqueness of tangent cones are now known (see \cite{Whi83, KW23a, KW23b, KW23c, DLS23a, DLS23b, DLMS23, DLMS24, MPSS24} for more details). As a possible long term goal, one might hope that the results of the present work suggest that integral varifolds satisfying certain `natural' geometric conditions have improved regularity, and moreover that the conditions which imply this are checkable when the varifold arises in a `natural' fashion, such as the limit of more regular surfaces satisfying certain structural restrictions or having `bounded' geometry or topology. We refer the reader to Remark \ref{remark:b} in Section \ref{sec:remarks} for further discussion of this.

In general, when a stationary integral varifold $V$ is sufficiently close, in the sense of varifolds, to a multiplicity $2$ plane, Theorem \ref{thm:main-intro} trivially gives the following dichotomy: either $V$ is represented by a $2$-valued Lipschitz graph, or there is a $\beta$-coarse gap within which $V$ \emph{fails} the topological structural condition. This in itself provides a new structural condition which must be satisfied about a branch point of locally infinite Morse index in a codimension $1$ stationary varifold (see Corollary \ref{cor:main}). This dichotomy (with the full regularity assumption replaced by an excess decay statement) is reminiscent of the main lemma in Simon's influential work \cite[Lemma 1]{Sim93}. However, the reader will note that, unlike in Simon's work, this does not provide a \emph{fixed} lower bound on the radius of the $\beta$-coarse gap; having a fixed lower bound on such a region would allow one to deduce further structural information. It turns out that we \emph{are} able to prove there is a fixed lower bound on the size of the $\beta$-coarse gap in which $V$ fails the topological structural condition, \emph{provided} certain more restrictive $\beta$-coarse gaps, which we call $(\beta,\gamma)$\emph{-fine gaps}, always satisfy the topological structural condition (with \emph{no} size restrictions). Thus, we can establish a general \emph{tri}chotomy for stationary integral varifolds near multiplicity $2$ planes.

To describe this trichotomy, we first define a $(\beta,\gamma)$\emph{-fine gap}. For this, fix $\beta,\gamma>0$. We then call a cylinder $C_\rho(x) = \R^k\times B^n_\rho(x)$ a $(\beta,\gamma)$\emph{-fine gap} if it is a $\beta$-coarse gap with the additional property that there is a cone $\BC$ which is the sum of $4$ distinct half-planes with a common boundary or two disjoint planes
such that
\begin{align*}
	\frac{1}{\rho^{n+2}}\int_{C_\rho(x)}&\dist^2(y,\BC)\, \ext\|V\|(y)\\
	& + \frac{1}{\rho^{n+2}}\int_{\BC\cap C_\rho(x)\setminus \{x:\dist(x,S(\BC))<\rho/8\}}\dist^2(y,\spt\|V\|)\, \ext\H^n(y)\\
	& \hspace{12em} < \gamma^2\inf_{z\in \R^k}\frac{1}{\rho^{n+2}}\int_{C_\rho(x)}\dist^2(y,\{z\}\times\R^n)\, \ext\|V\|(y),
\end{align*}
i.e.~$V$ is closer (by an amount determined by $\gamma$) to $\BC$ than it is to any plane parallel to $\{0\}^k\times\R^n$. We can then ask whether $V$ satisfies the topological structural condition in a $(\beta,\gamma)$-fine cylinder. Notice that if $\BC = P_1\cup P_2$ for two planes $P_1,P_2$ with $\inf_{x_i\in P_i}|x_1-x_2| =: \eta>0$, then for $\gamma = \gamma(\beta,\eta)$ sufficiently small, the topological structural condition is automatically satisfied in $(\beta,\gamma)$-fine gaps. For area minimisers, this is also true under the weaker condition that (cf.~\cite[Theorem 1.2]{KW23b} or \cite[Theorem 3.2]{DLMS23})
$$\frac{1}{\rho^{n+2}}\int_{C_\rho(x)}\dist^2(y,P_1\cup P_2)\, \ext\|V\|(y) < \eta\inf_{x_i\in P_i}|x_1-x_2|.$$
The general trichotomy that we establish is informally described as follows (see Theorem \ref{thm:main-3} for more details).

\begin{thmx}\label{thm:main-gap-intro}
	Fix $\beta>0,\gamma>0$. Then, there exists $\eps = \eps(n,k,\beta,\gamma)\in (0,1)$ such that the following is true. Suppose that $V$ is a stationary integral $n$-varifold in the cylinder $\R^k\times B^n_1(0)$ and that $V$ is $\eps$-close, in the sense of varifolds, to the disk $\{0\}^k\times B^n_1(0)$ with multiplicity $2$, i.e.
	\begin{itemize}
		\item $\frac{3}{2}\leq \w_n^{-1}\|V\|(B_1(0))\leq \frac{5}{2}$;
		\item $\hat{E}_V<\eps$.
	\end{itemize}
	Suppose also that $\Theta_V(0)\geq 2$. Then, one of the following alternatives must occur:
	\begin{enumerate}
		\item [(i)] there exists a cone $\BC$ comprised of $4$ half-planes (which includes the case where $\BC$ is the sum of $2$ full planes or a single multiplicity $2$ plane) and a scale $\theta = \theta(n,k,\beta,\gamma)\in (0,1/4)$ such that
		$$\frac{1}{\theta^{n+2}}\int_{C_{\theta}(0)}\dist^2(x,\spt\|\BC\|)\, \ext\|V\|(x) \leq \frac{1}{4}\int_{C_1(0)}\dist^2(x,\{0\}\times \R^n)\, \ext\|V\|(x)$$
		and
		$$\dist_\H(\spt\|\BC\|\cap B_1(0),\{0\}^k\times B^n_1(0)) \leq C\left(\int_{C_1(0)}\dist^2(x,\{0\}\times\R^n)\, \ext\|V\|(x)\right)^{1/2}$$
		for some $C = C(n,k,\beta,\gamma)$;
		\item [(ii)] there is a $\beta$-coarse gap for $V$ of radius $\geq \eta = \eta(n,k,\beta,\gamma)\in (0,1)$ which fails the topological structural condition;
		\item [(iii)]  there is a $(\beta,\gamma)$-fine gap for $V$ which fails the topological structural condition.
	\end{enumerate}
\end{thmx}
Under the assumption that there are no $\beta$-coarse gaps of \emph{any} size (and hence no $(\beta,\gamma)$-fine gaps for any $\gamma\in (0,1)$), one can more or less iterate the excess improvement conclusion (i) at all scales and centred at any density $\geq 2$ point in a manner which is by now well-established in the field. The result of doing so gives Theorem \ref{thm:main-intro}. One might ask whether it is possible to also have a uniform lower bound on the size of the $(\beta,\gamma)$-fine gap in Theorem \ref{thm:main-gap-intro}(ii); at present we do not know if this is possible or not.

Before continuing, let us give some examples to demonstrate occurrences of situations which exhibit $\beta$-coarse gaps (for suitable $\beta$) failing the topological structural condition.

\begin{example}\label{ex:catenoids}
	Consider, for $k=1,2,\dotsc$, the sequence of catenoids in $\R^3$ described in cylindrical coordinates $(r,\theta,z)\in [0,\infty)\times [0,2\pi)\times\R$ by:
	$$C_k := \{(r,\theta,z):z=k^{-1}\cosh^{-1}(kr)\}.$$
	These converge (in the sense of varifolds) to the multiplicity two plane $2|\{z=0\}|$. This shows that a stationary integral varifold (in fact, a smooth minimal surface) can become arbitrary close to a multiplicity two plane without ever becoming a $2$-valued graph over the plane. Furthermore, in this example $C_1(0)$ is not only a coarse gap which fails the topological structural condition, but it is in fact a fine gap (for suitable $\beta,\gamma>0$) which fails the topological structural condition.
\end{example}

\begin{example}\label{ex:truncated-catenoid}
	Using the same notation as in Example \ref{ex:catenoids}, let $C \equiv C_1 := \{(r,\theta,z):z=\cosh^{-1}(r)\}$ be the standard catenoid in $\R^3$. Let $h>0$ be the (unique) positive number for which the upward-pointing unit normal vector to $C$ at points on the circle $C\cap \{z=h\}$ makes an angle of $30^\circ$ with the horizontal plane $\{z=0\}$. Let
	$$\Sigma_\pm := \{x\in \R^3: x\pm (0,0,h)\in C\cap \{z>h\}\}\ \ \ \ \text{and}\ \ \ \ \Sigma_0 := \{z=0\}\cap \{r\leq \cosh(h)\}.$$
	Then,
	$$V:= |\Sigma_+| + |\Sigma_-| + |\Sigma_0|$$
	is a stationary integral varifold in $\R^3$ with $\Theta_V<2$ everywhere (this $V$ is sometimes called the \emph{truncated catenoid} or the $Y$\emph{-catenoid}). At each singular point $y\in \spt\|V\|\cap \{r=\cosh(h)\}$, $V$ has density $\Theta_V(y) = 3/2$ and the unique tangent cone to $V$ at $y$ is equal to three half-planes meeting at equal angles. Now, consider ``blowing down'' $V$, i.e.~look at the sequence $V_k := k^{-1}V$ for $k=1,2,\dotsc$. It is easy to show that $V_k\to 2|\{z=0\}|$ in the sense of varifolds. However, $V$ is never a $2$-valued graph over $\{z=0\}$; notice also that the support of $V$ remains connected in any cylinder centred at the origin. Again, these also provide examples of fine gaps (for suitable $\beta,\gamma>0$) which fail the topological structural condition.
\end{example}

\begin{remark}
	Notice that the varifolds in Example \ref{ex:truncated-catenoid} cannot be made into a current with zero boundary, as the triple junction singularities would introduce interior boundary points for the current.
\end{remark}

In the application of the main $\eps$-regularity theorem to Lipschitz $2$-valued stationary graphs, the Lipschitz assumption is crucial for the conclusion. This is seen from the following example.

\begin{example}\label{ex:graph}
	Let $V$ denote the stationary integral varifold associated to the holomorphic variety $\{(z_1,z_2,z_3)\in \C^3: z_1^2 = z_2^3 z_3\}$. Notice that $V$ has a singular point at the origin with $2|\{z_1=0\}|$ being the (unique) multiplicity $2$ tangent plane there. In a neighbourhood of the origin, $V$ \emph{is} a $2$-valued graph over its tangent plane. However, the $2$-valued graph is only Hölder continuous; it is the graph of $(z_2,z_3)\mapsto \pm z_2^{3/2}z_3^{1/2}$, which is $C^{0,1/2}$ as a $2$-valued function. We stress that our topological structural condition fails in cylinders centred at points on $\{z_2\neq 0,z_3=0\}$ (which are \emph{regular} points of the varifold). As this $V$ is also naturally an area minimising current, this example also shows that one cannot hope to prove that area minimising currents are locally given by \emph{Lipschitz} multi-valued graphs over their tangent planes about branch points. Note also that $0$ has ``planar frequency'' $2$, and is the limit of points with planar frequency $3/2$ (see \cite{KW23a} for the notion of planar frequency).
\end{example}

\textbf{A Guide to the Reader.} We felt it was necessary to keep the presentation of the paper reasonably self-contained, and with enough details to ensure readability. The length of the paper is a reflection of this, as well as the complexity of the arguments involved. The main technical parts of the proof appear in Part \ref{part:fine-reg}. On a first reading, the reader may wish to take the main result of Part \ref{part:fine-reg}, namely Theorem \ref{thm:fine-reg}, as a black-box and read Part \ref{part:intro}, Part \ref{part:coarse-blow-ups}, and Part \ref{part:blow-up-reg} in that order. The expert may also find it sufficient to only skim over Part \ref{part:coarse-blow-ups}, and go directly from Part \ref{part:intro} to Part \ref{part:blow-up-reg}.

As a very quick word on the proof for the reader acquainted with the proof of Allard's regularity theorem: the key difficulty in our proof is understanding the regularity of coarse blow-ups of sequences of (certain) stationary integral varifolds converging to a multiplicity $2$ plane. In the multiplicity $1$ case, one can readily use the weak $W^{1,2}$ convergence to show that the (single-valued) coarse blow-up must be a (weakly) harmonic function, in which case the Weyl lemma and elliptic PDE theory give both the regularity and estimates for coarse blow-ups. In the multiplicity $>1$ case, one has a multi-valued blow-up, for which the regularity theory (which we develop here subject to our topological structural condition on the varifolds) is significantly more involved. 

\subsection{Remarks and Future Directions}\label{sec:remarks}

We now discuss, through a series of remarks, how the results and techniques of the present work will be used in future work, as well as thoughts on other possible ways in which they might be used.

\begin{remarkx}[Structure of $2$-Valued Lipschitz Stationary Graphs]\label{remark:a}
	We begin by commenting on how our main theorem can be used to control the size of the singular set of a stationary $2$-valued Lipschitz graph. First, note that in the case of non-zero mean curvature, in general the singular set could have positive $\H^n$-measure; the union of the graphs of two $C^2$ functions which touch on a large, nowhere dense, closed set provides a simple example illustrating this. However, if we look at the case of zero mean curvature, i.e.~the stationary case, one can give a precise description of the singular set. Indeed, work of the second and third authors with B.~Krummel \cite{KMW25} shows that the branch set has Hausdorff dimension at most $n-2$ in this case. In fact, the combined results of the present work with \cite{KMW25, BK17} give the following:
	
	\begin{theorem}\label{thm:dimension}
		Suppose that $u:B^n_1(0)\to \A_2(\R^k)$ is a Lipschitz $2$-valued function with the property that $V:= \mathbf{v}(u)$ is a stationary integral $n$-varifold in $\R^k\times B^n_1(0)$. Then, there exists $\alpha = \alpha(n,k)\in (0,1)$ such that the following is true. Write $U\subset B^n_1(0)$ for the (open) subset of $B^n_1(0)$ consisting of the points locally about which $u$ is generalised-$C^{1,\alpha}$. Then we have
		$$\dim_\H(B^n_1(0)\setminus U)\leq n-3.$$
		In fact, we have the disjoint union
		$$\spt\|V\| = \mathscr{R}\cup \mathscr{B}\cup\mathscr{C}\cup \mathscr{K}$$
		and that $U = \pi(\mathscr{R}\cup\mathscr{B}\cup\mathscr{C})\equiv B^n_1(0)\setminus\pi(\mathscr{K})$, where $\pi:\R^{k}\times\R^{n}\to \R^n$ is the orthogonal projection map and:
		\begin{enumerate}
			\item [(i)]
			\begin{itemize}
				\item $\mathscr{R}$ is the set of regular points of $V$; 
				\item $\mathscr{B}$ is the set of singular points of $V$ where at least one tangent cone is a plane of multiplicity $2$ (by Theorem \ref{thm:main-intro}, this is then the unique tangent cone and $V$ locally has a generalised-$C^{1,\alpha}$ structure);
				\item $\mathscr{C}$ is the set of singular points of $V$ where at least one tangent cone consists of either a sum of $2$ distinct planes or $4$ distinct half-planes with a common boundary (by \cite{BK17} this is then the unique tangent cone and $V$ locally has a generalised-$C^{1,\alpha}$ structure);
				\item $\mathscr{K}$ is a closed subset, which consists of the singular points of $V$ where no tangent cone can be written as the sum of $4$ half-planes (including the possibility of two planes which may or may not coincide).
			\end{itemize}
			\item [(ii)]
			\begin{itemize}
				\item $\dim_\H(\mathscr{B})\leq n-2$ (\cite{KMW25}, relying on Theorem~\ref{thm:main-intro});
				\item $\dim_\H(\mathscr{C})\leq n-1$ (in fact, $\mathscr{C}$ is locally contained within an $(n-1)$-dimensional $C^{1,\alpha}$ submanifold by \cite{BK17});
				\item $\dim_\H(\mathscr{K})\leq n-3$;
                \item either $\Theta_{V}(x) = 1$ for every $x \in \mathscr{R}$ or $u(x) = 2\llbracket u_{1}(x) \rrbracket$ for every $x \in B_{1}^{n}(0)$ and some single-valued Lipschitz function $u_{1} \, : \, B_{1}^{n}(0) \to \R^{k}$ which is smooth away from a relatively closed set of Hausdorff dimension $\leq n-4$ (see \cite{Bar79} for this latter fact).
			\end{itemize}
			\end{enumerate}
            In particular, as $\sing(V) = \mathscr{B}\cup\mathscr{C}\cup\mathscr{K}$, we have $\dim_\H(\sing(V))\leq n-1$.
	\end{theorem}
	The proof of this theorem utilises the improved regularity and estimates provided in the present work (through Theorem \ref{thm:main-intro}) to show approximate monotonicity of the \emph{planar frequency function} (introduced by the third author and B.~Krummel \cite{KW23a}) for a stationary $2$-valued Lipschitz graph. We stress that establishing this approximate monotonicity requires stationarity, i.e.~the mean curvature to be $0$. Theorem \ref{thm:dimension} also uses the work of the first author \cite{BK17} to understand the local structure about non-planar tangent cones formed of the union of $4$ distinct half-planes. Similarly to \cite{Bar79} in the single-valued Lipschitz setting, it might be possible to improve the bound on the Hausdorff dimension of $\mathscr{K}$ in Theorem \ref{thm:dimension} (and hence of $B^n_1(0)\setminus U$) to $\leq n-4$; already in codimension one this seems to be an interesting question.
	
	We note that recent work by Hirsch--Spolaor \cite{HS24} gives that the singular set of a stationary $2$-valued Lipschitz graph has Hausdorff dimension $\leq n-1$ (namely, the last conclusion in Theorem \ref{thm:dimension} on the Hausdorff dimension of $\sing(V)$). It should be noted that \cite{HS24} establishes this using \emph{only} the Lipschitz regularity of the stationary graph, \emph{without} proving or using any higher regularity or uniqueness of tangent cones. They introduce \emph{multi-valued} generalised gradient Young measures to analyse the blow-up procedure, and show that these objects are sufficient for one to utilise Almgren's program with its center manifold construction in its general form, including the use of intervals of flattening and changing center manifolds (cf.~\cite{Alm00, DLS16a}). 
    
    To improve upon this result of \cite{HS24} and to prove the optimal dimension bound on the branch set (as asserted in the above Theorem~\ref{thm:dimension}), one then needs to analyse the possible occurrence of branch points of ``planar frequency $1$'': this question is studied in the present work via Theorem \ref{thm:fine-reg} (in Section~\ref{part:fine-reg}), and is ultimately settled in 
    Theorem~\ref{thm:main-intro} through the $C^{1,\alpha}$-style decay estimates at every branch point, implying that planar frequency $1$ branch points cannot occur. The work \cite{KMW25} uses this improved regularity information, namely through the locally uniform estimates provided, to then prove the optimal dimension bound on the branch set. As mentioned above, the decay estimate at branch points is used to show approximate monotonicity of the planar frequency function for a stationary $2$-valued Lipschitz graph. 
    An added benefit of this approach (in addition to providing detailed information on the local structure of the stationary Lipschitz graph including uniqueness of tangent cones near branch points) is that it allows one to significantly simplify (much like in the recent work \cite{KW23a, KW23b, KW23c} on area-minimising integral currents) the argument pertaining to Almgren's center manifold, which along with the higher regularity is used to prove the desired dimension bound ($\leq n-2$) on the branch set (and hence also the dimension bound of $\leq n-1$ on the full singular set). Indeed, for branch points which do not have \emph{exactly} quadratic decay (i.e.~those with planar frequency $\neq 2$) one can prove the branch set dimension bound using the planar frequency directly, \emph{without} using a center manifold. For branch points of quadratic (or higher) decay, one can then construct a \emph{single} center manifold which touches the varifold at \emph{all} nearby such points, which removes the need for both changing center manifolds and consideration of intervals of flattening. Ultimately, this allows one to conclude that in fact the branch set has codimension $\geq 2$, which is optimal; the remaining conclusions of Theorem \ref{thm:dimension} follow from the regularity theory developed by the first author in \cite{BK17}.
	
	Additionally, we already know from the work of L.~Simon and the third author \cite{SW16} that the Hausdorff dimension of the branch set is $\leq n-2$ for a stationary $2$-valued $C^{1,\alpha}$ graph (i.e.~when the graph is $C^{1,\alpha}$ rather than generalised-$C^{1,\alpha}$). In fact, in this case the branch set is either empty or has positive $\H^{n-2}$-measure; moreover, it is countably $(n-2)$-rectifiable and has unique tangent maps at $\H^{n-2}$-a.e.~branch point by work of the third author and B.~Krummel \cite{KW21}. Notice that if one assumes (or knows a priori, which happens in codimension one or when $V$ is associated to an area minimising current, as already noted) that no twisted cones occur as tangent cones, then generalised-$C^{1,\alpha}$ regularity for a $2$-valued stationary graph reduces to $C^{1,\alpha}$ regularity in the usual sense of $2$-valued functions. In such a situation, \cite{SW16, KW21} therefore already provide improved regularity results in our main theorem. It would be interesting to know if one can establish these improved regularity results in general, i.e.~allowing for twisted cones as tangent cones.
\end{remarkx}

\begin{remarkx}[Regions of Density $<2$ and Limits of Smoothly Embedded Minimal Surfaces]\label{remark:b}
	Here we discuss properties of the region $\{\Theta_V<2\}$ in $V$, which includes $\beta$-coarse gaps in which failure of the topological structural condition may in general take place. This will be particularly informative when $V$ arises as a varifold limit of smoothly embedded minimal surfaces.
	
	By general stratification theorems, one can always write the singular set of a stationary integral varifold as a disjoint union
	$$\sing(V) = \S_{n-1}\cup \mathcal{B}$$
	where $\dim_\H(\S_{n-1})\leq n-1$ and where at each point of $\mathcal{B}$ there is at least one tangent cone equal to a plane of multiplicity $\{2,3,\dotsc\}$. Thus, as $\beta$-coarse gaps are contained within $\{\Theta_V<2\}$, we have that $\sing(V)$ agrees with $\S_{n-1}$ in any such region. We therefore always have $\dim_\H(\sing(V)\cap \{\Theta_V<2\})\leq n-1$. In general, if we write
	$$\spt\|V\| \cap B_1(0) = \S^+\cup \S^-$$
	where $\S^+ := \{x\in \spt\|V\|\cap B_1(0):\Theta_V(x)\geq 2\}$ and $\S^-:= \{x\in \spt\|V\|\cap B_1(0):\Theta_V(x)<2\}$, then by upper semi-continuity of the density we know that $\S^-$ is a relatively open subset of $\spt\|V\|\cap B_1(0)$. As far as anyone at present knows, $\S^-$ may be a non-empty set whose $n$-dimensional Hausdorff measure is proportionately very small compared to $\spt\|V\|\cap B_1(0)$ (we stress that, in our hypotheses, conclusions, and analysis, we allow for the possibility that $\H^n(\sing(V))>0$). Nonetheless, our work shows that a certain topological control on parts of $\S^-$ (namely through topologically simple $\beta$-coarse gaps) is all that is additionally necessary for $V$ to correspond, near a multiplicity $2$ plane, to a $2$-valued Lipschitz graph over the plane with small Lipschitz constant and unique tangent cones. Consequently (in view of Theorem \ref{thm:dimension} in Remark \ref{remark:a}), assuming the topological structural condition and closeness to a multiplicity $2$ plane, we can conclude that either $\S^- = \emptyset$ or $\S^+ \equiv \spt\|V\|\setminus \S^-$ has Hausdorff dimension $\leq n-1$.
	
	By Allard's regularity theorem \cite{All72} together with the work of Simon \cite{Sim93}, straightforward tangent cone analysis (when $\Theta_V<2$) allows us to deduce that $\sing(V)\cap \S^-$ is equal to the disjoint union $\S_Y\cup \S_{n-2}$, where $\S_Y$ consists of points $z\in \spt\|V\|$ about which $V$ is represented by the sum of three disjoint embedded $C^{1,\alpha}$ submanifolds-with-boundary meeting at $120^\circ$ angles along a common $(n-1)$-dimensional free boundary containing $z$ (so-called \emph{triple junction} singularities), and $\dim_\H(\S_{n-2})\leq n-2$.
	
	Furthermore, we can write $\S_{n-2} = \S_{n-3}\cup \widetilde{\S}_{n-2}$, where $\dim_\H(\S_{n-3})\leq n-3$ and $\widetilde{\S}_{n-2}$ has the property that each point in $\widetilde{\S}_{n-2}$ is a limit point of $\S_Y$ (indeed, this follows from analysing the singular points in $2$-dimensional stationary cones with vertex density $<2$). As such, if $\S_Y = \emptyset$ (which would be the case if, for example, $V$ is the varifold limit of a sequence of embedded orientable minimal surfaces, or was otherwise known to be a current with vanishing boundary) then we would in fact have that $\sing(V)\cap \S^- = \S_{n-3}$ and so
	$$\dim_\H(\sing(V)\cap \S^-) \leq n-3.$$
	In particular, when $n=2$ and $\S_Y = \emptyset$, we have $\S^-\subset\reg(V)$. This provides a natural setting in which the key hypotheses under which Theorem \ref{thm:main-intro} holds is a hypothesis only on the regular part of the varifold. Furthermore, it would be interesting to know whether it is possible to show that the validity of the topological structural condition in $\beta$-coarse gaps where there are \emph{no} singularities \emph{implies} the topological structural condition in $\beta$-coarse gaps (perhaps with a different value of $\beta$) where the singular set is known to have Hausdorff dimension $\leq n-3$. If this were true then when $V$ does not contain triple junction singularities one can reduce the hypothesis in Theorem \ref{thm:main-intro} to being only on the regular part of $V$. Notice that Example \ref{ex:truncated-catenoid} above shows that one cannot do this when $V$ contains triple junction singular points.
	
	The assumption $\S_Y = \emptyset$, i.e.~that the varifold contains no triple junction singularities, occurs naturally when studying the closure in the varifold topology of smoothly embedded (or immersed) minimal surfaces in a Riemannian manifold under a uniform mass bound. Indeed, if the limiting (stationary integral) varifold $V$ has a triple junction singularity, the $\eps$-regularity theorem of Simon \cite{Sim93} implies that the approximating sequence $(V_j)_j$ must also exhibit triple junction singularities for sufficiently large $j$, contradicting their embeddedness.
	
	Taking this a step further, still assuming that $\S_Y = \emptyset$ so that $\sing(V)\cap \S^- = \S_{n-3}$, we may further write this as a disjoint union $\sing(V)\cap \S^- = \S_{n-4}\cup \widetilde{\S}_{n-3}$, where $\dim_\H(\S_{n-4})\leq n-4$ and each point $x\in \widetilde{\S}_{n-3}$ has the following property: at least one tangent cone to $V$ at $x$ takes the form $\BC_0\times\R^{n-3}$, where $\BC_0$ is a $3$-dimensional multiplicity $1$ cone in $\R^{k+3}$ with an isolated singularity at $0$. In particular, the \emph{link} of $\BC_0$, namely $\Sigma:=\BC_0\cap S^{k+2}$, is a smooth $2$-dimensional minimal surface in $S^{k+2}$. When $k=1$ (i.e.~the codimension is $1$), the work of Almgren \cite{Alm66} implies that the link $\Sigma$ \emph{cannot} have the topology of the sphere, and so must have genus $\geq 1$. Furthermore, when also $n=3$ the work of Simon \cite{Sim83b} gives that such a tangent cone must be unique, with the local structure of $V$ given by a $C^1$ perturbation of the cone. Hence, when $n=3$ and $k=1$, whenever $y\in \sing(V)$ has a tangent cone which is a plane of multiplicity $2$, the $\beta$-coarse gaps in which we must check the topological structural condition either have that $V$ is smoothly embedded, or there is a singular point belonging to $\widetilde{\S}_{n-3}$ ($\equiv \widetilde{\S}_0$ as $n=3$), where in the latter case the local structure of $V$ about the singular point is modelled by a cone whose link is a topologically non-trivial surface. It remains to be seen whether the latter situation can occur infinitely often if $V$ arises as the varifold limit of a sequence of smoothly embedded minimal surfaces with some type of topological control.
\end{remarkx}

\begin{remarkx}[Generalisations to Arbitrary Multiplicity]\label{remark:c}
    A natural question is why the main theorem of the present work, Theorem \ref{thm:main-intro}, is restricted only to the case $Q=2$, i.e.~when the plane has multiplicity $2$, rather than any arbitrary (integer) multiplicity $Q\geq 2$. Indeed, many of the techniques here apply for general $Q$ (and, for instance, the analysis in Part \ref{part:coarse-blow-ups} is performed for general $Q$). There are several ``topological structural conditions'' which one could assume for general $Q$, each assuming that a different amount of `splitting' occurs. Indeed, in cylindrical regions where the density of the varifold $V$ is everywhere $<Q$ and $V$ is close to a multiplicity $Q$ plane, one could impose one of the following graphical splitting conditions on a smaller cylinder:
	\begin{itemize}
		\item \emph{complete splitting}: $V$ is equal to the sum of $Q$ (possibly intersecting) smooth minimal graphs;
		\item \emph{partial splitting}: $V$ is equal to the sum of a $Q_1$-valued stationary Lipschitz graph and a $Q_2$-valued stationary Lipschitz graph (again, possibly intersecting) for some $Q_1,Q_2\in \Z_{\geq 1}$ with $Q_1+Q_2=Q$ (and so $V$ could have branch points of density $\leq \max\{Q_1,Q_2\}<Q$);
		\item \emph{no splitting}: $V$ is equal to a $Q$-valued (stationary) Lipschitz graph (which, as the density is $<Q$, has no branch points of density $Q$ but could have branch points of any density $<Q$).
	\end{itemize}
	Notice that these assumptions are progressively weaker and that when $Q=2$ all of these notions coincide. Indeed, a $2$-valued stationary Lipschitz graph (over a ball) with density $<2$ everywhere has no coincidence points and decomposes as the sum of two single-valued stationary Lipschitz graphs (cf.~Lemma \ref{lemma:Lip-2}). Subsequently, since we have estimates for single-valued Lipschitz minimal graphs which are close to a plane (by, for example, Allard regularity), we get improved estimates in these regions.
        
	In the above list `no splitting' is the most general assumption. Note that there is no immediate a priori reason why a graph as described in the `no splitting' situation should decompose as the sum of two multi-valued graphs each with $<Q$ values, i.e.~as in the `partial splitting' assumption. The `no splitting' scenario is the condition one naturally has, for instance, when studying stationary Lipschitz $Q$-valued functions for arbitrary $Q$. 
    
    In future work we will address the question of how one may generalise the results of the present work to general $Q$ under these various assumptions. For now, we note that generalising Theorem \ref{thm:main-intro} to arbitrary $Q$ under the `complete splitting' assumption should be a fairly routine (but somewhat involved) modification of the present work and known earlier work in settings of higher multiplicity (see \cite[Theorem 9.1]{BW19} and \cite{Wic14, MW24}). Under the `partial splitting' situation significant new difficulties occur, such as the need to perform iteratively finer and finer blow-up procedures in the form of a ``blow-up cascade'' as seen in \cite{Min21a, Min22}. The future work will focus on the most general `no splitting' case described above.
\end{remarkx}

\begin{remarkx}[Regularity Theorems near Unions of (Half-)Planes]\label{remark:d}
	Our main theorem is concerned with stationary integral varifolds close to a multiplicity $2$ plane. A key part of the proof (cf.~Theorem \ref{thm:fine-reg}) is establishing a special case of the theorem where the varifold is not only close to a multiplicity $2$ plane, but is much closer to a cone comprised of either a pair of distinct (multiplicity one) planes, or a twisted cone. A highly relevant and natural question for further analysis of singularities in stationary integral varifolds is understanding when one can establish an $\eps$-regularity theorem for stationary integral varifolds close to a cone which is either a sum of two distinct multiplicity one planes, or a twisted cone. Using the ideas of the present work (in fact, in a simpler fashion more akin to those by the first author in \cite{BK17}), analogues of Theorem \ref{thm:main-intro} and Theorem \ref{thm:main-gap-intro} can be established in this case also, provided one isolates the corresponding gap regions and topological structural conditions. Indeed, the `degenerate' case of this theorem, established as Theorem \ref{thm:fine-reg} of the present work, is the harder case of this result to establish. We discuss these results in Section \ref{sec:non-degenerate}.
\end{remarkx}

\begin{remarkx}[Non-Zero Mean Curvature and General Ambient Spaces]\label{remark:e}
	All the results of the present work have adaptations to the case where the generalised mean curvature is in $L^p$ for some $p>n$, and when the ambient space is a general (sufficiently regular) Riemannian manifold. These adaptations are discussed in Section \ref{sec:generalisation}. However, all the main difficulties and technicalities of the proof are already present in the stationary case. Indeed, one can deduce the more general $L^p$ mean curvature case by directly modifying the proof of the stationary case more or less in a manner analogous to that in the proof of Allard's regularity theorem  for $L^{p}$ mean curvature \cite{All72}. As such, most of the present paper will focus entirely on the case of \emph{stationary} integral varifolds in \emph{Euclidean space}.
\end{remarkx}

\begin{remarkx}[Area Minimisers]\label{remark:f}
    Suppose we are in the special case where $V$ is the varifold associated to an area minimising current with zero boundary. In general, it is still not known whether the tangent cone at a branch point is always unique, even for density $2$ branch points. The best general result known to date is uniqueness except on a set $\S_{\mathcal{B}}$ satisfying:
    \begin{enumerate}
        \item $\H^{n-2}(\S_{\B})=0$;
        \item each point in $\S_{\B}$ is a limit point of `large density gaps';
    \end{enumerate}
    (see \cite{KW23b, DLMS23}) and furthermore:
    \begin{enumerate}
        \item [(3)] along $\S_{\B}$ the current satisfies a `weak, locally uniform, approximation property' saying that for every $p\in \S_{\B}$, the current at any nearby point with density 
        $\geq \Theta_V(p)$ and at all sufficiently small scales (closeness to $p$ and smallness of scales both depending on $p$) is significantly closer to a non-planar cone than to any plane \cite[Theorem 1.1]{KW23a}.
    \end{enumerate}
    The arguments in the present work allow one to say more about the structure of the current within these large `density gaps' at points in $\S_{\B}$ except for a set of Hausdorff dimension $\leq n-3$, namely that they must be $(\beta_j,\gamma_j)$-fine 
    gaps (with  cones therein being planes intersecting along an axis of dimension $n-2$) which fail the topological structural condition, where $\beta_j,\gamma_j\to 0$. 
    
    We stress that many arguments in the present work are drastically simplified if one works in the special case of area minimisers. For instance, a significant part of the present work is establishing that coarse blow-ups have a monotone frequency function, which is already known to hold for general coarse blow-ups of area minimisers. Moreover, since area minimisers cannot be close to cones which are twisted or are two planes intersecting along an $(n-1)$-dimensional axis, the analysis related to fine blow-ups in Part \ref{part:fine-reg} can also be simplified. We also stress that, when considering blow-ups of area minimisers, the key analytic information that our topological structural condition gives is that the frequency of the blow-up must be $\geq 1$ everywhere. Points of frequency $<1$ being present in the blow-up necessarily give rise to $\beta$-coarse gaps failing the topological structural condition.

    Let us now elaborate further on Remark \ref{remark:area-min}. Suppose $T$ is a $2$-dimensional area minimising current in $B^{2+k}_2(0)$ and that $0\in \sing(T)$ is a branch point with $\Theta_T(0) = 2$. By \cite{Whi83}, we know that $T$ has unique tangent cones and moreover there is a $C^{1,\alpha}$ decay-rate of $T$ to any given tangent cone.
    
    We claim that for $r>0$ sufficiently small, one can verify the topological structural condition holds (relative to the \emph{fixed} unique tangent plane) on $B_{\rho}(0)\setminus B_{\rho/2}(0)$ for all $\rho\in (0,r)$. The basic idea behind this is to argue by contradiction, utilising that the local structure of $T$ near $0$ should be determined by a tangent map of $T$ at $0$ (i.e.~a coarse blow-up of suitable sequences $(\eta_{0,\rho_j})_\#T$ with $\rho_j\to 0$). Indeed, from such a sequence we may generate a $2$-valued Dir-minimising tangent map $u$ which is non-trivial, has zero-average, and homogeneous\footnote{Notice that by the decay rate to the (unique) tangent cone provided by \cite{Whi83} the planar frequency function from \cite{KW23a} is always approximately monotone. Thus, if the planar frequency is $<2$, this map can be generated by simply blowing-up off the (unique) tangent plane. Otherwise, if the planar frequency is $\geq 2$ one can generate such a map by blowing-up off a (single) center manifold.} (note that we do not claim uniqueness of this tangent map). In particular, the coincidence set of $u$ must just be $\{0\}$, since the singular set of such a function is isolated. This means that away from $0$ the two values of $u$ are distinct from one another, and that locally away from $0$ we know $u$ decomposes as the sum of two harmonic functions. But this separation of $u$ on $B_1(0)\setminus B_{1/2}(0)$ forces that the sequence $(\eta_{0,\rho_j})_\#T$ generating it must also separate (see \cite[Theorem 1.2]{KW23b} or \cite[Theorem 3.2]{DLMS23}) on this region as well, which proves the claimed separation property. However, the separation property on annuli at all scales combined with Allard regularity implies that the local structure of $T$ about $0$ is given by a $2$-valued $C^{1,\alpha}$ graph with an isolated singularity at $0$.

    Thus, it is in this way that one can use a frequency function and a separation property to determine the structure of $T$ by the behaviour of its tangent maps; notice that this argument also allows us to verify our topological structural condition needed to apply Theorem \ref{thm:main-intro} to determine the local structure, although the full strength of Theorem \ref{thm:main-intro} is in the end not needed in this instance. Nonetheless, it illustrates that a key aspect in determining the local structure about a singular point is through a suitable verification of a separation property like our topological structural condition. In this setting, we are not aware of a way in which one can directly verify the topological structural condition without deducing it from the structure of tangent maps and thus using an argument based on frequency functions (the reader should note that the homogeneity of the tangent map in the above argument was key to deducing separation of the tangent map away from $0$). Of course the above argument is specific to density $2$ branch points, as for density $\geq 3$ branch points it can be the case that the tangent maps have multiplicity as well. For this one needs a more complicated blow-up procedure to ``split'' the multiplicity and again reduce to the multiplicity one setting, and this is done by Chang \cite{Cha88} via building a \emph{branched center manifold}. 

    We end by noting that the above principle can also be applied about arbitrary branch points in $2$-dimensional area minimising mod $4$ currents $T$. Indeed, we know (from \cite{DLMS24}) that the tangent cone is a unique multiplicity $2$ plane and there is a $C^{1,\alpha}$ decay estimate of the current towards the (unique) tangent plane. Arguing as above, one can produce a tangent map $u$ which is non-trivial, average-free, and homogeneous, with the only branch point singularity being at $0$ (due to the homogeneity and the planar frequency at branch points being $>1$). However, for mod $4$ minimisers the blow-up could have (density $2$) classical singularities along rays. We then know from \cite[Remark 1.4]{DLMS24} that the local structure of the current near such points is determined by a $C^{1,\alpha}$ perturbation of these classical singularities. Away from the classical singularities, we again have the separation property for $T$, and so we know again the structure of $T$ on all sufficiently small annuli about the singular point, which allows us to determine the local structure as before. Alternatively, one can verify that our topological structural condition holds in a neighbourhood of the branch point and directly apply Theorem \ref{thm:main-intro}. It should be noted that, unlike in the case of area minimising integral currents discussed above, in both arguments here we need the fine $\eps$-regularity theorem established in Theorem \ref{thm:fine-reg} (in the much simpler form as in the similar results originally established in \cite{Wic14, MW24}, as ``density gaps'' do not occur in this setting) in order to say that the singular points of the tangent map away from $0$ really induce classical singularities in the current. (For the area minimising mod $4$ setting, this simpler version of the fine $\eps$-regularity result in Theorem \ref{thm:fine-reg} also follows from the results in \cite{DLMS24}, which in turn utilise an adaptation of the techniques used in \cite{Wic14, MW24} together with the key height bound established in \cite[Theorem 1.2]{KW23b} or \cite[Theorem 3.2]{DLMS23}.)
\end{remarkx}

\begin{remarkx}[Alternative to Structural Condition]\label{remark:g}
    There is an alternative structural condition one can make in a $\beta$-coarse gap which is equivalent to the topological structural condition in Definition \ref{defn:top-str-con} (up to allowing the parameter $\beta$ to change, similar to Proposition \ref{prop:sheeting-equivalence}). Indeed, suppose $V$ is a stationary integral $n$-varifold in $B^{n+k}_2(0)$ and $C_\rho(x)$ is a $\beta$-coarse gap for $V$. This alternative condition is then:
    \begin{equation}\tag{$\|$}\label{E:vertical-tangents}
        \text{In $C_{\rho/2}(x)$, $V$ has no triple junction singularities or regular points $x$ with $\pi_{P_0}(T_xV) \neq P_0$.}
    \end{equation}
    Note that because $\Theta_V<2$ in a $\beta$-coarse gap, necessarily all regular points there must have multiplicity one. The condition that $\pi_{P_0}(T_xV) = P_0$ is equivalent to requiring that $T_xV$ can be written as the graph of a linear function over $P_0$, and thus the second requirement above is saying that the tangent plane at regular points is not ``vertical'' relative to $P_0$. In the codimension one setting, the condition $\pi_{P_0}(T_xV) \neq P_0$ is equivalent to saying that the unit normals to $T_xV$ and $P_0$ are orthogonal.
    
    With this viewpoint, Theorem \ref{thm:main-intro} gives as a corollary that the only density $2$ branch points in an \emph{arbitrary} stationary integral varifold $V$ about which Theorem \ref{thm:main-intro} does \emph{not} directly apply are those which are limit points of either (i) triple junction singularities, or (ii) multiplicity one regular points with tangent planes ``vertical'' to $P_0$ (as defined above); here $P_0$ is a choice of the (possibly non-unique) tangent plane at the branch point. This corollary should be compared with \cite[Corollary 3.6]{Wic14}. Notice that as a special case, at any point $x\in \spt\|V\|$ where the (unique) tangent cone is a plane and the map $(\reg(V)\cap \{\Theta_V<2\})\cup\{x\}$ sending $y\mapsto T_yV$ is continuous at $x$, then one can check the above conditions are satisfied and thus the conclusions of Theorem \ref{thm:main-intro} hold locally about $x$.
    
    In particular, as triple junction singularities do not occur in area-minimising currents, the only manner in which the assumptions of Theorem \ref{thm:main-intro} can fail to hold locally about a branch point in that setting are when there are multiplicity one regular points converging to the branch point which have tangent planes vertical to the given one. Notice this is exactly what happens in Example \ref{ex:graph}.

    Let us now justify the equivalence between \eqref{E:vertical-tangents} and the topological structural condition in Definition \ref{defn:top-str-con}. The fact that the topological structural condition implies \eqref{E:vertical-tangents} can be seen as a direct consequence of Proposition \ref{prop:sheeting-equivalence} (again, up to changing the parameter $\beta$ to ensure the correct radius of the cylinder for the conclusions) as the graph structure rules out both triple junction singularities and regular points with vertical tangent planes.

    Now let us assume \eqref{E:vertical-tangents} for some choice of $\beta$, and suppose $C_\rho(x)$ is a $\tilde{\beta}$-coarse gap for $V$, where $\tilde{\beta}\in (0,\beta)$ is suitably small to be specified. By Remark \ref{remark:b}, since $V$ has no triple junction singularities in $C_{\rho/2}(x)$ we know that $\dim_\H(\sing(V)\cap C_{\rho/2}(x))\leq n-3$. Now consider $A:= \pi_{P_0}(\sing(V)\cap C_{\rho/2}(x))$: we know that this has $\H^{n-2}(A)=0$, and thus $B^n_{\rho/2}(x)\setminus A$ is simply-connected (cf.~\cite[Appendix]{SW16}). Now consider the function $f:B^n_{\rho/2}(x)\setminus A\to \{0,1,2,\dotsc\}$ given by
    $$f(x) := \#\{y\in \R^k: (y,x)\in \spt\|V\|\} \equiv \#\pi^{-1}_{P_0}(x)\cap \spt\|V\|.$$
    Notice that for $x\in B^n_{\rho/2}(x)\setminus A$, by definition of $A$ we know that if $(y,x)\in \spt\|V\|$ then $(y,x)$ is a regular point of $V$, which as $C_\rho(x)$ is a $\tilde{\beta}$-coarse gap we have $\Theta_V<2$ here, and thus $(y,x)$ necessarily is a multiplicity one regular point of $V$. Our assumption that no regular points have tangent planes vertical to $P_0$ therefore gives that $f$ is a locally constant function, and thus by connectedness of $B^n_{\rho/2}(x)\setminus A$ we know that $f$ is constant. Our mass assumption on $V$ in the $\tilde{\beta}$-coarse gap then implies that we must have either $f\equiv 1$ or $f\equiv 2$; the former case cannot happen provided we choose (as we may) $\tilde{\beta} = \tilde{\beta}(n,k)$ sufficiently small (indeed, this follows from Almgren's Lipschitz approximation Theorem \ref{thm:Lipschitz-approx}, as this implies for sufficiently small $L^2$ norm there must be regions of $P_0$ where $V$ has two points again a given one in $P_0$, and thus points where $f = 2$). Hence, $f\equiv 2$. Thus, as $B^n_{\rho/2}(x)\setminus A$ is simply-connected, there are smooth functions $u_1,u_2:B^n_{\rho/2}(x)\setminus A\to \R^k$ such that
    $$V\res (\R^k\times (B^n_{\rho/2}(x)\setminus A)) = \mathbf{v}(u_1) + \mathbf{v}(u_2).$$
    Writing $V_1 = \mathbf{v}(u_1)$ and $V_2 = \mathbf{v}(u_2)$, we then know that these are stationary integral varifolds in $C_{\rho/2}(x)\setminus (\R^k\times A)$. In fact, they extend as stationary integral varifolds in $C_{\rho/2}(x)\setminus \sing(V)$ (this follows from looking at the connected components of $\reg(V)$ which contain $V_1$ and $V_2$, which must be disjoint from one another by the density upper bound and the assumption on the tangent planes not being vertical). Since $\H^{n-1}(C_{\rho/2}(x)\cap \sing(V)) = 0$, a standard argument (using the volume bounds from $V$) gives that they extend as stationary integral varifolds in $C_{\rho/2}(x)$. But then we must have $\spt\|V_1\|\cap \spt\|V_2\|=\emptyset$, since if $x\in \spt\|V_1\|\cap \spt\|V_2\|$ then $\Theta_V(x) = \Theta_{V_1}(x) + \Theta_{V_2}(x) \geq 2$, a contradiction to the fact that $\Theta_V<2$ in the $\tilde{\beta}$-coarse gap. Thus, we have that $V\res C_{\rho/2}(x) = V_1+V_2$, where $V_1$ and $V_2$ are stationary integral varifolds with disjoint supports. Since $V_1$ and $V_2$ are graphical on $\R^k\times (B^n_{\rho/2}(x)\setminus A)$ and $\H^{n-2}(A)=0$, it follows that each of $V_1$, $V_2$ have non-empty intersection with $C_{\rho/4}(x)$, thus proving the topological structural condition holds under \eqref{E:vertical-tangents}, thereby proving the equivalence.
    \qed
\end{remarkx}

\subsection{Overview of the Proof}

As stated above, from now until Section \ref{sec:generalisation} we will focus entirely on the case of \emph{stationary} integral $n$-varifolds in $B^{n+k}_2(0) \subset \R^{n+k}$.

The starting point for our theory is Almgren's (weak) $Q$-valued Lipschitz approximation. This gives that if $(V_j)_{j=1}^\infty$ is a sequence of stationary integral $n$-varifolds in $B_1(0)$ which converge to the multiplicity $Q$ disk $Q|\{0\}^k\times B^n_1(0)|$ then, for all sufficiently large $j$, we can approximate $V_j$ by the graph of a $Q$-valued Lipschitz function $u_j$. It is then possible to show (see Part \ref{part:coarse-blow-ups}) that one can take a (subsequential) limit
$$v_j:=\hat{E}_{V_j}^{-1}u_j\to v$$
in a suitable topology, where
$$\hat{E}_{V_j}:=\left(\int_{\R^k\times B^n_1(0)}\dist^2(x,P_0)\, \ext\|V_j\|(x)\right)^{1/2}$$
is the $L^2$ \emph{height excess} of $V_j$ relative to $P_0 = \{0\}^k\times\R^n$. Here, $v$ is a $Q$-valued $W^{1,2}_{\text{loc}}(B^n_1(0))$ function which we call a \emph{coarse blow-up}. The function $v$ however does \emph{not} naturally minimise a functional, nor is it even clear whether or not $v$ is stationary for the Dirichlet energy. Therefore, in proving properties about coarse blow-ups, our principle means is using the stationarity of the varifolds themselves. One of the key difficulties in gaining regularity information about coarse blow-ups is the possibility of what we call ``concentration of tilt-excess'', namely it is unknown whether or not the tilt-excess of the varifold, that is $\int_{\R^k\times B^n_1(0)}\|\pi_{P_0}-\pi_{T_xV}\|^2\, \ext\|V\|(x)$, can concentrate in the non-graphical part of $V$ (here, $\pi_L$ denotes the natural orthogonal projection onto a subspace $L$ and $\|\cdot\|$ the Hilbert--Schmidt norm). This issue is morally saying we do not know whether or not we have energy convergence of the sequence $v_j$ to the coarse blow-up $v$. Slightly more precisely, the concern is that the quantity
$$\hat{E}_{V_j}^{-2}\int_{\text{non-graphical part of $V$}}\|\pi_{P_0}-\pi_{T_xV}\|^2\, \ext\|V\|(x)$$
might fail to be small. In spite of this, we prove a new non-concentration estimate for the gradient \emph{of the coarse blow-up}. To state this cleanly, write $v_f:=\sum^Q_{\alpha=1}\llbracket v^\alpha(x)-v_{a}(x)\rrbracket$ for the average-free-part of $v$, where $v_a(x):= Q^{-1}\sum^Q_{\alpha=1}v^\alpha(x)$ is the average-part of $v$. Furthermore, write $v^\kappa_f:= v_f\cdot e_\kappa$ where $\{e_1,\dotsc,e_k\}$ is the standard orthonormal basis for $\R^k$. Then, our estimate says that for all $\delta\in (0,1)$ and $\kappa\in \{1,\dotsc,k\}$, we have
$$\sum_{\alpha=1}^{Q}\int_{B^n_{1/2}(0)\cap \{|v_f^{\kappa,\alpha}|<\delta\}}|Dv_f^{\kappa,\alpha}|^2 \leq C\delta$$
where $C = C(n,k,Q)>0$. We stress that this estimate \emph{only requires stationarity of the varifolds} $(V_j)_{j=1}^\infty$, and it holds in \emph{any dimension and codimension} as well as for \emph{any multiplicity}; it therefore \emph{does not require any additional structural hypothesis on the varifolds} and so in particular it does not require the aforementioned `topological structural condition' to prove it. For our purposes, this estimate keeps control of the energy coming from a neighbourhood of potential branch points in $v$ when checking the so-called `squeeze' identity for a generalised-$C^1$ coarse blow-up $v$ (which is a key step in establishing monotonicity of the frequency function for general $v$). This theorem and other basic properties of coarse blow-ups are proved in Part \ref{part:coarse-blow-ups}. 

Next, if we consider a coarse blow-up $v$ which is generalised-$C^1$, it may have some points in its graph at which the tangent cone is a twisted cone or the union of a pair of distinct planes. Near these points, one might hope to gain a better understanding of $v$ from the sequence of varifolds blowing-up to it; in particular, we would like to prove the aforementioned `squeeze' identity near such points. For this, we will use an $\eps$-regularity theorem for $v$. This will come from a so-called \emph{fine} $\eps$-regularity theorem for the varifolds $(V_j)_j$. To establish this, we will again need to utilise the topological structural condition in suitable $\beta$-coarse gaps (in fact, it suffices to know the topological structural condition holds in suitable $(\beta,\gamma)$-fine gaps). The scenario one finds themselves in is that not only are the varifolds $V_j$ converging to a multiplicity $2$ plane, but simultaneously each $V_j$ is significantly closer to a non-planar cone $\BC_j$ than it is to $2|P_0|$, where $\BC_j$ is either a twisted cone or the union of two distinct planes. In such a situation, the idea is that one wishes to instead blow-up $V_j$ off $\BC_j$ in order to establish an $\eps$-regularity theorem. This general procedure is known as a \emph{fine blow-up} and was first introduced by the third author in \cite{Wic14}. To carry out this procedure requires $L^2$ estimates which are inspired by -- but far more involved than -- the main estimates in \cite[Section 3]{Sim93}. Establishing them requires ideas similar to those seen in \cite[Section 10]{Wic14}, \cite[Theorem 3.1]{MW24}, and \cite[Section 3.2]{BK17} (see also \cite{Kru14b} for highly relevant similar ideas in the Dir-minimising setting), and indeed we will need to unify and develop all of these ideas further in order to achieve this. Establishing these estimates, and then subsequently the fine $\eps$-regularity theorem, takes up a significant portion of this work; this is done in Part \ref{part:fine-reg}. We stress that a major difficulty when compared with the works \cite{Sim93, Wic14, MW24} is that we may now \emph{not} have points of ``good density'' converging to every point on the axis of the cone; this not only significantly complicates the analysis used to establish the main $L^2$ estimates, but also allows (as seen in the simpler setting of \cite{BK17}) the axis dimension of the new cone to decrease in the excess decay lemma. 

To elaborate on this further, let $\hat{F}_{V,\BC}$ denote the \emph{two-sided fine excess} of $V$ relative to such a cone $\BC$. This is formally defined in \eqref{E:fine-excess}, but roughly speaking it has the form
$$\hat{F}_{V,\BC}:=\left(\int\dist^2(x,\spt\|\BC\|)\, \ext\|V\|(x) + \int\dist^2(x,\spt\|V\|)\, \ext\|\BC\|(x)\right)^{1/2}.$$
The analysis in Part \ref{part:fine-reg} then takes place in the regime where the coarse excess is small and the fine excess is much smaller than the coarse excess, i.e.~$\hat{E}_V\ll 1$ and
\begin{equation}\label{E:intro-1}
	\inf_{\BC}\hat{F}_{V,\BC}\ll \hat{E}_V.
\end{equation}
The main fine $\eps$-regularity theorem in Part \ref{part:fine-reg} (see Theorem \ref{thm:fine-reg}) then roughly says that these hypotheses imply the full regularity conclusions (A), (B), (C) stated in Theorem \ref{thm:main-intro}. One may therefore think of this as saying that, in our strategy, in order to prove our main theorem in its full generality it is necessary to fully handle this special case of it first. As mentioned, key to this are $L^2$ estimates for fine blow-ups which are based on the remainder term in the monotonicity formula. Indeed they show, for instance, that in the above regime, if $\Theta_V(0)\geq 2$ then
$$\int_{B_{1/2}(0)}\frac{|x^{\perp_{T_xV}}|^2}{|x|^{n+2}}\, \ext\|V\|(x) \leq C\int_{B_1(0)}\dist^2(x,\spt\|\BC\|)\, \ext\|V\|(x).$$
Notice that this estimate is \emph{much} easier to prove when \eqref{E:intro-1} fails, as then it suffices to bound the left-hand side by $\hat{E}_V^2$ which is a much more standard argument.

We stress however that, even when \eqref{E:intro-1} holds (and not just $\hat{E}_V\ll 1$), our topological structural condition is \emph{still} necessary to show conclusions (A), (B), and (C) of Theorem \ref{thm:main-intro}. This is illustrated by the following example.

\begin{example}[Scherk surfaces with small angle]\label{ex:scherk}
	For a fixed choice of parameter $\phi\in [0,\pi)$, consider the Weierstrass representation given by $g(z)=iz$ and
$$f(z) = \frac{4}{(z^2-e^{-2i\phi})(z^2-e^{2i\phi})}.$$
This defines a Scherk surface with the parameter $\phi$ being half the angle between the two planes whose union form the tangent cone at infinity to the surface. Thus, if we let $\phi\to 0$ then these (smoothly embedded) minimal surfaces converge to a multiplicity two plane. Hence, given $\eps_j,\gamma_j\downarrow 0$, by appropriately rescaling and rotating we can produce a sequence $(V_j)_{j=1}^\infty$ of stationary integral $2$-varifolds, each corresponding to a smoothly embedded surface in $\R^3$, such that
$$\hat{E}_{V_j}<\eps_j \qquad \inf_{\BC}\hat{F}_{V_j,\BC}<\gamma_j\hat{E}_{V_j}$$
but such that they never form a 2-valued graph over any ball $B_r^n(0)$, $r>0$.
\end{example}

In Part \ref{part:blow-up-reg}, we then establish regularity of coarse blow-ups under the assumptions of the main theorem, namely in the specific case when $Q=2$ and assuming the varifolds obey our topological structural condition in any $\beta$-coarse gap. This requires using the previous fine $\eps$-regularity theorem as well as the energy non-concentration estimate for coarse blow-ups described previously. We will use them to show that coarse blow-ups $v$ which are (qualitatively) generalised-$C^1$ enjoy monotonicity of Almgren's frequency function, i.e.~of the quantity
$$\rho\ \longmapsto\ \frac{\rho^{2-n}\int_{B_\rho(0)}|Dv|^2}{\rho^{1-n}\int_{\del B_\rho(0)}|v|^2}.$$
To then understand general coarse blow-ups (without any a priori regularity assumption) we classify the coarse blow-ups which are homogeneous of degree one. It is within this classification that one is able to leverage the previous case in which one \emph{assumes} qualitative generalised-$C^1$ regularity, and which ultimately gives that coarse blow-ups are, a posteriori, always generalised-$C^{1,\alpha}$ and that they satisfy uniform decay estimates. 

In Part \ref{part:extension} we give extensions of the results in the present work. We first prove that, a posteriori, we \emph{do} have energy convergence, and hence strong $W^{1,2}_{\text{loc}}(B_1^n(0))$ convergence, in this setting to the coarse blow-up. Then, we extend the main theorem to the case where the varifold has mean curvature which instead lies in $L^p$ for some $p>n$, as well as to more general ambient Riemannian manifolds. As already noted, the case of non-vanishing mean curvature is particularly useful to bear in mind for the results of the paper due to Brakke's example, which we quickly summarise.

\begin{example}[{Brakke \cite[Section 6.1]{Bra78}}]\label{ex:brakke}
	One can construct a $2$-dimensional integral varifold $V$ in $\R^3$ with bounded mean curvature such that there is a set $A\subset\R^3$ of positive $\H^2$-measure for which no element $x\in A$ has a neighbourhood in which $V$ can be represented by the graph of a multi-valued function. To construct $V$, we use a model based on the catenoid in $\R^3$ as follows. The catenoid has zero mean curvature. Taking any radius $r>0$ and $B>0$, one can take a catenoid with a very small central hole, and gradually deform the two sheets together away from the hole. The result is an integral varifold with mean curvature uniformly bounded by $B$, has a hole, and agrees with a multiplicity $2$ plane outside $B_r(0)$. To construct $V$ from this model, start with a multiplicity $2$ plane in $\R^3$, and remove a disjoint collection of disks whose union is dense in the plane, yet leaves behind a set $A$ of positive $\H^2$-measure. Replace each disk with a section of the aforementioned bent catenoid with a hole so that the edges match smoothly. The resulting integral varifold $V$ has integer density, is smoothly immersed, and has mean curvature bounded uniformly by $B$. Yet, if $x\in A$, then $V$ has holes in every neighbourhood of $A$ and hence cannot be represented by the graph of a multi-valued function. Notice also that the (density 2) branch set of $V$ is exactly $A$, and thus has positive $\H^2$-measure. This example clearly fails the topological structural condition, as the holes are regions where the density is $<2$, yet the support is connected. We also remark that $V$ arises as the varifold limit of a sequence of smoothly immersed surfaces which have uniformly bounded area and mean curvature uniformly bounded by $B$.
\end{example}

\begin{remark}\label{remark:cascade}
As a remark for the experts, we wish to describe a key moral difference in understanding the regularity of different types of blow-ups (for us in the present work, coarse blow-ups and fine blow-ups). The difference lies in when the cone one wishes to blow-up relative to occurs with multiplicity and a multi-valued approximation must be made instead of a tuple of single-valued functions. Again, in the present work coarse blow-ups need to use a multi-valued (Lipschitz) approximation, whilst for the fine blow-up procedure we are able to use single-valued functions. In a situation where one can use single-valued functions, one can rotate the components of the cone (2 full planes, or 4 half-planes with a common boundary, in the present work) independently of one another (keeping a common axis in the case of 4 half planes not forming 2 planes) whilst performing a blow-up procedure, which ultimately allows one the freedom to subtract different linear functions \emph{independently} from the components of the blow-up, whilst still maintaining key estimates for the blow-up (this is also key for the analysis in \cite{BK17, Min21a}, the latter doing this with multiplicity also). This allows one more freedom with the homogeneous degree one elements that can be subtracted off from the blow-up in question, which makes some aspects of the analysis of their regularity easier. However, if one must use a multi-valued approximation (as we must for our coarse blow-ups) we are only able to subtract the \emph{same} linear function from the values of the multi-valued blow-up whilst remaining in the class of such blow-ups. The difference here is a natural artefact arising from the geometric nature of the blow-up procedure in question, which either generates certain freedoms or restrictions, namely whether the components of the cone (counting multiplicity) can be rotated independently of each other or not during the blow-up process.

The consequence is that when the blow-up has a (genuinely) multi-valued part, one needs a \emph{fine} $\eps$-regularity theorem for the varifolds in order to understand the regularity of the original blow-up; this is ultimately the reason why in the present work we need a fine $\eps$-regularity theorem for coarse blow-ups (which are multi-valued) and not for fine blow-ups. The finer blow-up procedure `splits' the multiplicity, which inductively one expects to be easier to understand. This should be compared with the setting in \cite{Min21a, Min22}, where since the fine blow-up could still have a multi-valued part, a further, even finer, blow-up procedure was needed (called an \emph{ultra fine blow-up} in \cite{Min21a}) in order to understand the regularity of fine blow-ups. As the ultra fine blow-ups in \cite{Min21a} were single-valued, no subsequent finer blow-up procedures were needed. The upshot of this is that one only needs $\eps$-regularity theorems for blow-ups when they are genuinely multi-valued, and in principle these can be established via the means of the `blow-up cascade' idea mentioned previously in Remark \ref{remark:c} of Section \ref{sec:remarks} (in this language, the fine blow-up is a blow-up cascade of length $1$, with the fine blow-up/ultra fine blow-up procedure in \cite{Min21a} a blow-up cascade of length $2$). We expect these ideas to be useful in the future. Another thing to note is that to establish a full $\eps$-regularity theorem one needs to rule out bad gap regions of \emph{all} (arbitrarily small) scales; this explains why we cannot give a fixed lower bound on the $(\beta,\gamma)$-fine gaps in Theorem \ref{thm:main-gap-intro} (but only on the $\beta$-coarse gaps failing the topological structural condition) as we need the full fine $\eps$-regularity theorem to understand the regularity of coarse blow-ups, which requires us to rule out \emph{all} $(\beta,\gamma)$-fine gaps which fail the topological structural condition.

It should be noted however that even when the blow-up is comprised only of single-valued functions (and so the regularity can be understood using the freedom to subtract of separate linear functions, as described above), one could presumably also understand the regularity using the blow-up cascade idea. For instance, in the present setting, to understand regularity of fine blow-ups relative to a union of planes with spine dimension $n-1$, one needs to allow blow-ups off unions of planes with spine dimension $\leq n-2$ also. However, one could first try to establish (through a fine blow-up procedure) a fine $\eps$-regularity theorem when $V$ is \emph{much} closer to a union of planes with spine dimension $\leq n-2$ than it is to a union of planes with spine dimension $n-1$; if this were established first, then one would only need to consider blow-ups off unions of planes with spine dimension $n-1$, as whenever $V$ is much closer to a cone of lower spine dimension we would immediately have the regularity conclusion, and otherwise we can replace the cone of lower spine dimension by one with spine dimension $n-1$.
\end{remark}

\textbf{Acknowledgements.} Part of this research was conducted during the period PM was a Clay Research Fellow. PM and NW would like to thank Homerton College, Cambridge, for their hospitality and for providing a pleasant environment where many discussions relating to this work took place.

\section{Notation and Preliminaries}\label{sec:prelim}

We start by setting out the basic notation and terminology that is common to all of the paper. We refer the reader to \cite{MW24} for more explanation on the language where necessary. Throughout the paper we fix integers $n\geq 2$, $k\geq 1$. In the following, $x_0\in \R^{n+k}$ and $\rho>0$.
\begin{itemize}
	\item $B_\rho(x_0):= \{x\in \R^{n+k}:|x-x_0|<\rho\}$ is the open ball of radius $\rho$ centred at $x_0$ in $\R^{n+k}$. If we want to emphasise the dimension of a ball we write $B_\rho^{n+k}(x_0)\equiv B_\rho(x_0)$.
	\item $\eta_{x_0,\rho}:\R^{n+k}\to \R^{n+k}$ given by $\eta_{x_0,\rho}(x):= \rho^{-1}(x-x_0)$.
	\item $\tau_{x_0}:\R^{n+k}\to \R^{n+k}$ given by $\tau_{x_0}(x):= \eta_{x_0,1}(x) = x-x_0$.
	\item For $s\geq 0$, $\H^s$ denotes the $s$-dimensional Hausdorff measure on $\R^{n+k}$.
        \item $\dim_\H$ denotes the Hausdorff dimension.
	\item $\w_n:= \H^n(B^n_1(0))$.
	\item For $A,B\subset\R^{n+k}$, $\dist_\H(A,B)$ denotes the Hausdorff distance between $A$ and $B$.
	\item For $A\subset\R^{n+k}$, $\dist(x_0,A):= \inf_{y\in A}|y-x_0|$.
	\item For $A\subset\R^{n+k}$, $B_\rho(A)\equiv (A)_\rho := \{x\in \R^{n+k}:\dist(x,A)<\rho\}$.
\end{itemize}
By a \emph{plane} we mean any $n$-dimensional affine subspace of $\R^{n+k}$. Given a plane $P$ we write:
\begin{itemize}
	\item $\pi_P$ for the orthogonal projection $\R^{n+k}\to P\subset\R^{n+k}$ onto $P$.
	\item $x^{\top_P} \equiv \pi_P(x)$ and $x^{\perp_P} \equiv x-x^{\top_{P}} \equiv \pi_{P^\perp}(x)$.
\end{itemize}
By a \emph{half-plane} we mean a closed half-plane, i.e.~any set which is the closure of one of the connected components of $P\setminus L$, where $P$ is any plane and $L$ is any $(n-1)$-dimensional affine subspace of $P$. For any half-plane $H$ we write $\pi_H$ for the orthogonal projection onto the unique plane (of the same dimension) containing $H$. The distinguished plane throughout the paper is
$$P_0 := \{0\}^k\times\R^n.$$
For $x_0\in P_0$, we will often abuse notation and write $B^n_\rho(x_0)$ to mean $\{x\in P_0:|x-x_0|<\rho\}$, with the superscript $n$ indicating that the ball is $n$-dimensional lying in the $n$-dimensional plane $P_0$. We also define cylinders $C_\rho(x_0):= \R^k\times B^n_\rho(x_0)$ and more generally for a plane $P$ and $x_0\in P$ we write
$$C_\rho(x_0;P):=\{x\in \R^{n+k}:|\pi_P(x)-x_0|<r\}.$$
This definition also makes sense for $x_0\not\in P$ by defining $C_\rho(x_0;P) := C_\rho(\pi_P(x_0);P)$. Note that $C_\rho(x_0) \equiv C_\rho(x_0;P_0)$.

\subsection{Varifolds and Cones}

Let $G(n,n+k)$ denote the Grassmannian of $n$-dimensional subspaces of $\R^{n+k}$. An $n$\emph{-varifold} $V$ on an open set $U\subset\R^{n+k}$ is a Radon measure on $U\times G(n,n+k)$. Recall then the following notation for $n$-varifolds $V$:
\begin{itemize}
	\item The \emph{weight measure} is denoted $\|V\|$.
	\item The \emph{support} is $\spt\|V\|$.
	\item If $M$ is a countably $n$-rectifiable set, it defines a (multiplicity one) $n$-varifold by $|M|(A):=\H^n(\{x:(x,T_xM)\in A\})$ for $A\subset U\times G(n,n+k)$ (here we insist that $n$-rectifiable sets have finite $\H^n$-measure).
	\item We write $T_xV$ for the \emph{approximate tangent plane} to $\spt\|V\|$ at $x$ when it exists.
\end{itemize}
We say that an $n$-varifold $V$ is \emph{integral} if there exist positive integers $\{c_j\}_{j=1}^\infty\subset\{1,2,\dotsc\}$ and $n$-rectifiable sets $\{M_j\}_{j=1}^\infty$ for which
$$V = \sum^\infty_{j=1}c_j|M_j|.$$
The \emph{regular part} of $V$, denoted $\reg(V)$, is the set of $x\in \spt\|V\|$ for which $B_\rho(x)\cap \spt\|V\|$ is a smooth, $n$-dimensional embedded $C^1$ submanifold of $B_\rho(x)$ for some $\rho>0$. We also write $\sing(V)$ for the (interior) \emph{singular set} of $V$, i.e.~$\sing(V):= (\spt\|V\|\setminus\reg(V))\cap U$.

Given an integral $n$-varifold $V$ in $B_2(0)$, we write
$$\hat{E}_V:=\left(\int_{C_1(0)}\dist^2(x,P_0)\, \ext\|V\|(x)\right)^{1/2}$$
for the $L^2$ \emph{height excess} of $V$ from $P_0$. We also call $\hat{E}_V$ the \emph{coarse} (\emph{height}) \emph{excess} of $V$. We can define the (scale-invariant) $L^2$ height excess over different cylinders $C_\rho(x_0)$, where $x_0\in P_0$ and $\rho>0$, by
$$\hat{E}_V(C_\rho(x_0)):= \hat{E}_{(\eta_{x_0,\rho})_\#V} \equiv \left(\frac{1}{\rho^{n+2}}\int_{C_\rho(x_0)}\dist^2(x,P_0)\, \ext\|V\|(x)\right)^{1/2}.$$

Let us now introduce notation for the various sets of cones we will be interested in:
\begin{itemize}
	\item Write $\mathcal{P}$ for the set of all integral $n$-varifolds $\Cbf$ in $\R^{n+k}$ of the form $\BC = |P_1| + |P_2|$, where $P_1$ and $P_2$ are \emph{distinct} $n$-dimensional planes intersecting along an affine subspace $S(\BC):= P_1\cap P_2\neq\emptyset$. We call $S(\BC)$ the \emph{spine} of $\BC$.
	\item For $j=0,\dotsc,n-1$, define $\mathcal{P}_j := \{\BC\in\mathcal{P}:\dim (S(\BC))=j\}$ and $\mathcal{P}_{\leq j}:=\{\BC\in\mathcal{P}:\dim(S(\BC))\leq j\}$.
	\item Write $\mathcal{P}_\emptyset$ for the set of all integral $n$-varifolds $\BC$ in $\R^{n+k}$ which are of the form $\BC = |P_1| + |P_2|$ with $P_1\cap P_2=\emptyset$.
	\item Write $\CC_{n-1}$ for the set of all integral $n$-varifolds $\BC$ in $\R^{n+k}$ which are of the form $\BC = \sum^4_{i=1}|H_i|$, where for each $i=1,\dotsc,4$ the $H_i$ are distinct half-planes which have a common boundary and only meet along their common boundary, which moreover is $S(\BC) := \bigcap^4_{i=1}H_i$, i.e.~the boundary of each half-plane is the same for each half-plane, which in this case is an affine $(n-1)$-dimensional subspace.
	\item Define $\CC:= \CC_{n-1}\cup \mathcal{P}_{\leq n-2}$.
\end{itemize}
Notice that $\mathcal{P}_{n-1}\subset\mathcal{C}_{n-1}$. The set of cones $\CC$ together with the set of (multiplicity 2) planes will be all the cones we are interested in for the present work. If we wish to include multiplicity two planes, we will use the notation $\widetilde{\CC}:=\CC\cup \{2|P|:P\text{ is an $n$-dimensional plane in $\R^{n+k}$}\}$. The subset of $\CC$ consisting of cones whose spine contains the origin (and so is a `cone' in the usual sense of geometric measure theory, namely being scale-invariant about $0$) is denoted by $\CC_0$. Similarly, we can define $\widetilde{\CC}_0\subset \widetilde{\CC}$ in an analogous manner.

For notational simplicity, if $V$ is a varifold and $x\in \R^{n+k}$, we write $\dist(x,V) \equiv \dist(x,\spt\|V\|)$.

\subsection{$Q$-Valued Functions and their Graphs}\label{sec:q-valued-notation}

There are many basic facts about multi-valued functions that we will not repeat here; standard references are \cite{Alm00, DLS11}. For an integer $Q\geq 1$, the space of unordered $Q$-tuples of points in $\R^k$ (identified with $\sum_{\alpha=1}^{Q} \llbracket a^{\alpha}\rrbracket$ where $\llbracket a^{\alpha} \rrbracket$ is the Dirac mass at $a^{\alpha} \in \R^{k}$) is denoted by $\A_Q(\R^k)$ which comes with the metric $\G$ defined by 
$$\G\left(\sum_{\alpha=1}^{Q} \llbracket a^{\alpha}\rrbracket\, , \sum_{\alpha=1}^{Q} \llbracket b^{\alpha}\rrbracket\right) := \min_{\pi} \, \sqrt{\sum_{\alpha=1}^{Q} \left|a^{\alpha} - b^{\pi(\alpha)}\right|^{2}},$$
where $\pi$ denotes a permutation of $\{1, 2, \ldots, Q\}$. If $g:B^n_1(0)\to \A_Q(\R^k)$ we write:
\begin{itemize}
	\item For $x\in B^n_1(0)$, $g(x) = \sum^Q_{\alpha=1}\llbracket g^\alpha(x)\rrbracket$.
	\item Letting $\{e_1,\dotsc,e_k\}$ denote the standard orthonormal basis of $\R^k$, we define $g^\kappa: B^n_1(0)\to \A_Q(\R)$ by $g^\kappa(x):= \sum^Q_{\alpha=1}\llbracket g^\alpha(x)\cdot e_\kappa\rrbracket$ for $\kappa\in \{1,\dotsc,k\}$.
	\item For each such $\kappa$, we can globally order the $Q$ values of $g^\kappa$, i.e.~there exist $Q$ functions $g^{\kappa,\alpha}:B^n_1(0)\to \R$ with $g^{\kappa,1}\leq\cdots\leq g^{\kappa,Q}$ such that $g^\kappa(x) = \sum^Q_{\alpha=1}\llbracket g^{\kappa,\alpha}(x)\rrbracket$.
	\item $g_a:= Q^{-1}\sum^Q_{\alpha=1}g^\alpha$ for the (single-valued) \emph{average-part} of $g$ (or just \emph{average} of $g$).
	\item $g_f:= g-g_a\equiv \sum^Q_{\alpha=1}\llbracket g^\alpha-g_a\rrbracket$ for the \emph{average-free part} of $g$.
	\item $|g(x)|^2:= \G(g(x),Q\llbracket 0\rrbracket)^2 = \sum^k_{\kappa=1}\sum^Q_{\alpha=1}|g^{\kappa,\alpha}(x)|^2$, where $|g(x)|\geq 0$.
\end{itemize}
If $g:B^n_1(0)\to \A_Q(\R^k)$ is Lipschitz, then we write
\begin{itemize}
	\item $\mathbf{v}(g)$ for the integral $n$-varifold naturally associated to the $n$-rectifiable set
	$$\graph(g):= \{(y,x)\in \R^k\times B^n_1(0):y=g^\alpha(x)\text{ for some }\alpha\in \{1,\dotsc,k\}\}$$
	whilst also taking into account multiplicity. Sometimes we will abuse notation and write $\graph(g)\equiv \mathbf{v}(g)$.
\end{itemize}
More precisely, $\mathbf{v}(g)$ is defined as the multi-valued push-forward of $B^n_1(0)$ under the map $F:x\mapsto (g(x),x)\in \A_Q(\R^{n+k})$ (cf.~\cite[Section 1.5]{Alm00}). We may then define
\begin{itemize}
	\item $J\equiv J(g):B^n_1(0)\to \A_Q(\R)$ to be the Jacobian of the map $F$, i.e.~for a.e. $x\in B^n_1(0)$, the $Q$ values of $J$ are given by
	$$J^\alpha(x):= \sqrt{\det[(\delta_{pq}+D_pg^\alpha(x)\cdot D_qg^\alpha(x))_{p,q=1,\dotsc,n}]} \qquad \text{for }\alpha\in \{1,\dotsc,Q\}.$$
	(Thus technically $J$ is only defined on $B^n_1(0)$ up to a set of measure $0$.)
\end{itemize}
Noting that $\sum^Q_{\alpha=1}J^\alpha(x)$ is a well-defined single-valued function, we also record here the following consequence of the area formula (cf.~\cite[Section 1.6 \& 1.7]{Alm00} or \cite[Corollary 1.11]{DLS13}): for any measurable set $A\subset B^n_\sigma(0)$ and any bounded compactly supported Borel function $f:\R^k\times B^n_1(0)\to \R$ we have
\begin{equation}\label{E:area-formula}
	\int_{\R^k\times A}f(x)\, \ext\|\mathbf{v}(g)\|(x) = \int_A\sum^Q_{\alpha=1}f(g^\alpha(x),x)J^\alpha(x)\, \ext x.
\end{equation}
Moreover, using the Cauchy--Binet formula and the Taylor expansion of $t\mapsto (1+t)^{1/2}$, we have the following standard estimates for the Jacobian:
\begin{equation}\label{E:J-bounds}
	1\leq J^\alpha(x)\leq 1+C|D^\alpha g(x)|^2
\end{equation}
where $C = C(n,k)\in (0,\infty)$.

\subsection{Stationary Integral Varifolds}

We now collect a few facts about stationary integral varifolds, all of which are well-known to experts (cf.~\cite{All72, Sim83} for more details). Recall that an integral $n$-varifold $V$ is said to be \emph{stationary} in $B_2(0)$ if
\begin{equation}\label{E:stationarity-1}
\int_{B_2(0)\times G(n,n+k)}\div_S\Phi(x)\, \ext V(x,S)=0
\end{equation}
for every $\Phi\in C^1_c(B_2(0);\R^{n+k})$. Such varifolds have approximate tangent spaces $\H^n$-a.e. The fact that $V$ is stationary in $B_2(0)$ implies (by taking $\Phi = \phi e_\kappa$) that the coordinate functions are harmonic on the support of $V$, in the sense that for any $\kappa\in\{1,\dotsc,n+k\}$ we have
\begin{equation}\label{E:stationarity-2}
\int_{B_2(0)}\nabla^V x^\kappa\cdot \nabla^V\phi(x)\, \ext\|V\|(x) = 0 \qquad \text{for all }\phi\in C^\infty_c(B_2(0)).
\end{equation}
The monotonicity formula then states that for any $x_0\in B_2(0)$ and $0<\sigma\leq\rho <2-|x_0|$, we have
\begin{equation}\label{E:monotonicity}
	\frac{\|V\|(B_\rho(x_0))}{\rho^n} - \frac{\|V\|(B_\sigma(x_0))}{\sigma^n} = \int_{B_\rho(x_0)\setminus B_\sigma(x_0)}\frac{|(x-x_0)^\perp|^2}{|x-x_0|^{n+2}}\, \ext\|V\|(x)
\end{equation}
where here (and throughout the paper) in an integrand depending on $x$, unless specified otherwise the notation $\perp$ denotes the projection onto the approximate tangent space of the measure being integrated over. Thus in \eqref{E:monotonicity}, $\perp\, \equiv\, \perp_{T_xV}$.

A standard consequence of the monotonicity formula is the existence everywhere of the \emph{density}
$$\Theta_V(x_0):= \lim_{r\downarrow 0}\frac{\|V\|(B_r(x_0))}{\w_n r^n}.$$
The density function is also upper semi-continuous in both the variables $x_0$ and $V$ (in their natural respective topologies).

It is elementary to check that
$$\sum^k_{\kappa=1}|\nabla^Vx^\kappa|^2 = \frac{1}{2}\|\pi_{P_0} - \pi_{T_xV}\|^2.$$
Thus, by taking $\Phi(x) = \eta(x)x^\kappa$ in \eqref{E:stationarity-2} for an appropriate choice of cut-off function $\eta$ we obtain for any $\sigma\in (0,1)$ the well-known reverse Poincaré inequality for stationary integral varifolds:
\begin{equation}\label{E:rpi}
	\int_{\R^k\times B^n_\sigma(0)}\|\pi_{P_0}-\pi_{T_xV}\|^2\, \ext\|V\|(x) \leq \frac{128}{(1-\sigma)^2}\hat{E}_V^2
\end{equation}
i.e.~we have a priori control on the $L^2$ tilt-excess in terms of the $L^2$ height excess. With more work Allard proved the following supremum estimate which we will need at various points:
\begin{theorem}[{\cite[Theorem 6, Section 7.5]{All72}}]\label{thm:allard-sup-estimate}
	If $V$ is a stationary integral $n$-varifold in $B_2(0)$ and $\sigma\in (0,1)$, then
	$$\sup_{x\in \spt\|V\|\cap C_\sigma(0)\cap B^{n+k}_{31/16}(0)}\dist(x,P_0) \leq \frac{C}{(1-\sigma)^n}\hat{E}_V$$
	for some $C = C(n,k)\in (0,\infty)$.
\end{theorem}

\subsection{Tangent Cones and Stratification}

Given a stationary integral $n$-varifold $V$ in $B_2(0)$ and $x\in \spt\|V\|$:
\begin{itemize}
	\item $\vartan_xV$ denotes the set of all varifold tangent cones to $V$ at $x$, i.e.~all limits of the form $\lim_{j\to\infty}(\eta_{x,\rho_j})_\#V$ for some sequence $\rho_j\downarrow 0$.
	\item Each varifold $\BC\in \vartan_x V$ is a stationary integral \emph{cone}, i.e.~for any $\lambda>0$ we have $(\eta_{0,\lambda})_\#\BC = \BC$.
	\item The \emph{spine} of $\BC$ is $S(\BC):= \{y\in \spt\|\BC\|: (\tau_y)_\#\BC = \BC\}$. It can be shown that $S(\BC)$ is a subspace and equals $S(\BC) = \{y\in\spt\|\BC\|:\Theta_{\BC}(y) = \Theta_{\BC}(0)\}$.
	\item For $j\in \{0,1,\dotsc,n\}$, we define the $j$\emph{-strata} of $V$ by
	$$\S_j:=\{x\in \sing(V):\dim(S(\BC))\leq j\,\text{ for all }\BC\in\vartan_xV\}.$$
\end{itemize}
The standard result regarding the strata is that for each $j$ (cf.~\cite[Theorem 2.26]{Alm00})
\begin{equation}\label{E:strata-bound}
	\dim_\H(\S_j)\leq j.
\end{equation}
We again note that $\BC\in \CC$ is not actually a cone in the above (usual) sense unless $0\in S(\BC)$.

\section{Main Results}

Our main results will concern the following class of stationary integral varifolds, indexed by $\beta\in (0,1)$; they are those which satisfy the topological structural condition in any $\beta$-coarse gap. Recall $P_0:= \{0\}^k\times\R^n$ and $C_\rho(x_0):=\R^k\times B^n_\rho(x_0)$.

\begin{defn}\label{defn:V-beta}
	Fix $\beta\in (0,1)$. Write $\mathcal{V}_\beta$ for the set of all stationary integral $n$-varifolds $V$ in $B_2(0)$ such that $V$ satisfies the topological structural condition in any $\beta$-coarse gap contained in $C_1(0)$. More precisely, whenever $x\in P_0$ and $\rho>0$ are such that $C_\rho(x)\subset C_1(0)$ and:
	\begin{enumerate}
		\item [(a)] $C_\rho(x)\subset\{\Theta_V<2\}$;
		\item [(b)] $\frac{3}{2}\leq (\w_n\rho^n)^{-1}\|V\|(C_\rho(x))\leq \frac{5}{2}$;
		\item [(c)] For some $z\in P_0^\perp\equiv \R^k\times\{0\}^n$,
		$$\frac{1}{\rho^{n+2}}\int_{C_\rho(x)}\dist^2(y,\{z\}\times\R^n)\, \ext\|V\|(y) < \beta^2;$$
	\end{enumerate}
	then $\spt\|V\|\cap C_{\rho/2}(x)$ is disconnected, and at least two of its connected components have non-empty intersection with $C_{\rho/4}(x)$.
\end{defn}

\begin{remark}\label{remark:inclusion}
	Notice that for $\beta_1<\beta_2$ we have $\mathcal{V}_{\beta_2}\subset\mathcal{V}_{\beta_1}$, i.e.~$\beta\mapsto \mathcal{V}_\beta$ increases as $\beta$ decreases. As such, when proving theorems for $\mathcal{V}_\beta$ we can always assume that $\beta$ is small.
\end{remark}
The topological structural condition guarantees that there are at least two `large' connected components of $\spt\|V\|$ within the half-cylinder of a $\beta$-coarse gap. This will, for instance, rule out `necks' joining sheets together, as well as current boundaries (if $V$ is the varifold associated to a current). Indeed, we will see how this implies various sheeting properties in Section \ref{sec:tsc}.

The reader will note that the topological structural condition in $\mathcal{V}_\beta$ \emph{only} refers to $\beta$-coarse gaps which are cylinders over the plane $P_0$. It does however not care about which plane parallel to $P_0$ the $L^2$ height excess of $V$ is small relative to in such a cylinder. Consequently, $\mathcal{V}_\beta$ does not remain closed under rotations which change $P_0$. It does have the following (positive) closure properties:
\begin{enumerate}
	\item [(i)] $\mathcal{V}_\beta$ \emph{is} closed under translations and rescalings, i.e.~if $x_0\in \R^k\times B_1^n(0)$ and $0<\rho<\dist(x,\del C_1(0))$, then $(\eta_{x_0,\rho})_\#V\res B_2(0)\in\mathcal{V}_\beta$ whenever $V\in\mathcal{V}_\beta$. Indeed, the key point here is that $\eta_{x_0,\rho}$ fixes $P_0$ and maps cylinders over $P_0$ to cylinders over $P_0$.
	\item [(ii)] $\mathcal{V}_\beta$ \emph{is} closed under rotations \emph{which leave $P_0$ invariant}, i.e.~if $\Gamma\in SO(n+k)$ has $\Gamma(P_0) = P_0$, then $\Gamma_\#V\in \mathcal{V}_\beta$ whenever $V\in\mathcal{V}_\beta$.
\end{enumerate}
For rotations which do not leave $P_0$ invariant, the impact is a change to the value of $\beta$, i.e.:
\begin{enumerate}
	\item [(iii)] Suppose $V\in\mathcal{V}_\beta$ and $\Gamma\in SO(n+k)$ is a small rotation, i.e.~$\|\Gamma-\id\|\ll 1$, and set $\eta:= \dist_\H(P_0\cap B_1, \Gamma(P_0)\cap B_1)$ (which is controlled by $\|\Gamma-\id\|$). Then, one can check that, if $\eta = \eta(n,k,\beta)$ is sufficiently small, then for some small $\delta = \delta(n,k)$ we have $(\eta_{0,1-\delta})_\#\Gamma_\#V\res B_2(0)\in \mathcal{V}_{\widetilde{\beta}}$ for some $\widetilde{\beta} = \widetilde{\beta}(n,k,\beta)\in (0,1)$. See Proposition \ref{prop:rotation-class} for further details.
\end{enumerate}
The reason for these observations is that usually in regularity arguments within geometric measure theory one works with a class of varifolds which is closed under homotheties and rotations. This provides certain freedoms which can be used to simplify arguments, such as when establishing properties of blow-ups or when iterating an excess decay lemma (in which, say, one can reduce to working with a fixed plane by rotating the new plane to the original one). The above observations show that we are allowed to do this, at the expense of decreasing $\beta$ by a fixed amount.

There are several other (more restrictive) classes one could work with instead of $\mathcal{V}_\beta$ by changing parts of the Definition \ref{defn:V-beta}, such as:
\begin{itemize}
	\item One could change the definition of a $\beta$-coarse gap to instead require that the $L^2$ \emph{tilt}-excess relative to $P_0$ in the cylinder is $<\beta$, rather than the minimum of the $L^2$ height excess over planes parallel to $P_0$. This class is still not naturally closed under rotations.
	\item The strongest assumption one could make would be to require the topological structural condition over \emph{any} cylinder over \emph{any} plane where the $L^2$ excess is small, rather than just over planes parallel to $P_0$. This class would then have natural closure properties under both homotheties and rotations. This would be sufficient for the application to $2$-valued stationary Lipschitz graphs, for example (see Theorem \ref{thm:main-2}).
\end{itemize}
There are several other ways in which one could change the definition of $\mathcal{V}_\beta$, allowing for different constants or different requirements in the $\beta$-coarse gaps. These all turn out to be equivalent to the current definition of $\mathcal{V}_\beta$, up to changing $\beta$ by a controlled amount; in particular, they are equivalent for the purpose of our main theorems. Let us quickly summarise some of these:
\begin{itemize}
	\item We will see that, up to changing the parameter $\beta$ by a controlled amount, the topological structural condition in $\mathcal{V}_\beta$ is equivalent to a sheeting condition in the cylinder $C_{\rho/2}(x)$ for any $\beta$-coarse gap $C_{\rho}(x)$; this observation is made precise in Proposition \ref{prop:sheeting-equivalence}. As such, one could alternatively work with this class of varifolds (denoted $\widetilde{\mathcal{V}}_\beta$ in Proposition \ref{prop:sheeting-equivalence}) and the results established would be equivalent to those established for $\mathcal{V}_\beta$. This explains the discussion in Part \ref{part:intro}, particularly in Remark \ref{remark:c}, of the topological structural condition being a sheeting property. Notice that a priori it appears that the requirement in the definition of $\mathcal{V}_\beta$ is weaker than requiring a sheeting property on a smaller cylinder, such as in $C_{\rho/4}(x)$, as indeed a sheeting property could still hold on $C_{\rho/4}(x)$ yet these graphs could belong to the same connected component in $C_{\rho/2}(x)$ (this is what happens for a catenoid with its neck not centred at the centre of the cylinder). This observation however shows that the two properties provide equivalent results in this paper.
	\item One possible generalisation of Definition \ref{defn:V-beta} is to introduce a second parameter $c\in (0,1)$ determining the scale at which $V$ has two disconnected components, i.e.~replacing $\spt\|V\|\cap C_{\rho/2}$ and $\spt\|V\|\cap C_{\rho/4}(x)$ in the topological structural condition by $\spt\|V\|\cap C_{c\rho}(x)$ and $\spt\|V\|\cap C_{c\rho/2}(x)$, respectively. However, if this class were denoted $\mathcal{V}_{\beta,c}$, then in fact we have $\mathcal{V}_{\beta,c}\subset\mathcal{V}_{\widetilde{\beta}}$ for suitable $\widetilde{\beta} = \widetilde{\beta}(\beta,c)$, simply by applying the definition condition to suitable translations and rescalings of $V$. Thus, it suffices to prove results for the classes $\mathcal{V}_{\beta}$ to deduce them for $\mathcal{V}_{\beta,c}$.
	\item Another possible modification of Definition \ref{defn:V-beta} is to introduce a parameter $\gamma\in (0,1)$ and only require the topological structural condition in $\beta$-coarse gaps of radius $\leq\gamma$; if this new class is denoted $\mathcal{V}_{\beta}^\gamma$, then clearly $\mathcal{V}_{\beta}\subset\mathcal{V}^\gamma_\beta$. One can then readily check (using the equivalence with a sheeting property mentioned two points above) that for appropriate $\widetilde{\beta} = \widetilde{\beta}(\beta,\gamma)$ we have $\mathcal{V}^\gamma_\beta\subset\mathcal{V}_{\widetilde{\beta}}$. Therefore, again it suffices to prove results for the classes $\mathcal{V}_{\beta}$ to get them for $\mathcal{V}_{\beta}^\gamma$.
\end{itemize}

\begin{remark}\label{remark:small-singular-set}
We briefly recall the following discussion from Remark \ref{remark:b} of Section \ref{sec:remarks}. Consider a stationary integral $n$-varifold $V$ in $B_2(0)$ which contains no triple junction singularities (or equivalently, by \cite{Sim93}, singular points $x\in \S_{n-1}\setminus \S_{n-2}$ with $\Theta_V(x) = \frac{3}{2}$). For instance, $V$ could be the varifold associated to a current with zero boundary. Then, we have $\dim_\H(\sing(V)\cap \{\Theta_V<2\})\leq n-3$, as $\sing(V)\cap \{\Theta_V<2\} = \S_{n-3}\cap \{\Theta_V<2\}$ by elementary tangent cone analysis and Simon's $\eps$-regularity theorem for the triple junction \cite{Sim93}. Therefore, in the case when $n=2$ the defining condition of $\mathcal{V}_\beta$ refers only to the regular part of $V$.
\end{remark}

To state our main regularity theorem, we need to recall the notion of a multi-valued function being \emph{generalised-}$C^1$ or \emph{generalised-}$C^{1,\alpha}$ (cf.~\cite{MW24}, where a similar notion was defined in codimension one).

\begin{defn}[Generalised-$C^1$ and $C^{1,\alpha}$ Regularity]\label{defn:gc-1}
	Fix a domain $U\subset B^n_2(0)$. We say that a continuous function $f:U\to \A_2(\R^k)$ is \emph{generalised-}$C^1$ in $U$ and belongs to $GC^{1}(U)$ if to each $x\in U$ we can assign a function $Af_x:\R^n\to \A_2(\R^k)$ such that:
	\begin{enumerate}
		\item [(1)] $\graph(Af_x)\in \widetilde{\CC}\cup\mathcal{P}_\emptyset$;
		\item [(2)] $\lim_{h\to 0}|h|^{-1}\G(f(x+h),Af_x(x+h))=0$;
		\item [(3)] We have $U = \mathcal{R}_f\cup \CC_f^* \cup \mathcal{B}_f$ is a disjoint union, where:
		\begin{itemize}
			\item $\mathcal{R}_f$ is the \emph{regular set} of $f$, namely the set of $x\in U$ for which there exists $\rho = \rho(f,x)>0$ such that $f|_{B_\rho(x)} = \llbracket f^1\rrbracket + \llbracket f^2\rrbracket$, where $f^1,f^2:B^n_{\rho}(x)\to \R^k$ are $C^1$ (and possibly intersect);
			\item $\CC_f^*$ is the set of \emph{twisted classical singularities} of $f$, namely the set of $x\in U\setminus \mathcal{R}_f$ for which there exists $\rho = \rho(f,x)>0$ such that $\graph(f)\res C_\rho(x)$ is equal to four embedded, $n$-dimensional, $C^1$ submanifolds-with-boundary with the same $C^1$ boundary, which contains $x$, and moreover only intersect at their common boundary. In particular, $\graph(Af_x)\in \CC_{n-1}\setminus\mathcal{P}_{n-1}$;
			\item $\mathcal{B}_f:= U\setminus (\mathcal{R}_f\cup\CC^*_f)$ is the \emph{branch set} of $f$; if $x\in\mathcal{B}_f$ then $\graph(Af_x) = 2|P|$ for some plane $P = P(f,x)$;
            \item [(4)] The assignment $x\mapsto Af_x$ is continuous on both $\mathcal{R}_f\cup \mathcal{B}_f$ and $\CC_f^*$ separately.
		\end{itemize}
		We also write $\CC_f$ for the set of classical singularities of $f$, namely the $x\in U$ for which either $x\in \CC_f^*$, or there exists $\rho = \rho(f,x)>0$ such that $f|_{B_\rho(x)} = \llbracket f^1\rrbracket + \llbracket f^2\rrbracket$, where $f^1,f^2:B_\rho^n(x)\to \R^k$ are $C^1$ and $f^1(x) = f^2(x)$.
	\end{enumerate}
	Furthermore, we say that $f\in GC^{1,\alpha}(U)$ and that $f$ is \emph{generalised-}$C^{1,\alpha}$ in $U$ if $f$ is generalised-$C^1$ in $U$ and additionally: (a) for each compact set $K\subset U$ and $y\in\mathcal{R}_f$, if $\rho_y>0$ is the largest radius on which $f|_{B_{\rho_y}(y)} = \llbracket f^1\rrbracket + \llbracket f^2\rrbracket$ with $f^1,f^2$ being $C^1$, the functions $f^1,f^2$ are actually in $C^{1,\alpha}(\overline{B}_{\rho_y}(y)\cap K)$; (b) each $y\in \CC_f^*$, is a $C^{1,\alpha}$ classical singularity of $\graph(f)$; (c) on each compact set $K\subset \mathcal{R}_f\cup \mathcal{B}_f$ we have $\sup_{x_1,x_2\in K:\,x_1\neq x_2}\frac{\G(Af_{x_1},Af_{x_2})}{|x_1-x_2|^\alpha}<\infty$, where for $x\in \mathcal{R}_f\cup\mathcal{B}_f$ we view $Af_x$ as being valued in $\A_2(\mathcal{M}_{k,n})$, where $\mathcal{M}_{k,n}$ is the set of $k\times n$ matrices.
\end{defn}

Our main regularity theorem is the following. We stress that $\hat{E}_V\equiv \hat{E}_{V,P_0}$ and so the smallness assumption at unit scale is over the same plane as the cylinders are taken over in the definition of $\mathcal{V}_\beta$.

\begin{theorem}\label{thm:main}
	Fix $\beta\in (0,1)$. Then, there exists $\eps_0 = \eps_0(n,k,\beta)\in (0,1)$ such that the following is true. Suppose that $V\in\mathcal{V}_\beta$ satisfies
	$$\frac{3}{2}\leq \frac{\|V\|(B_1(0))}{\w_n}\leq\frac{5}{2} \qquad \text{and} \qquad \hat{E}_V<\eps_0.$$
	Then, there is a $GC^{1,\alpha}$ function $f:B^n_{1/2}(0)\to \A_2(\R^k)$ such that
	$$\widetilde{V}\res C_{1/2}(0) = \mathbf{v}(f)$$
    where $\widetilde{V} := V\res B_{3/2}(0)$. Furthermore, the following conclusions hold:
	\begin{enumerate}
		\item [(1)] $\|f\|_{C^{0,1}}\leq C\hat{E}_V$;
		\item [(2)] At every point $x\in \spt\|\widetilde{V}\|\cap C_{1/2}(0)$, there is a unique tangent cone $\Cbf_x$ equal to either a plane of multiplicity $\in \{1,2\}$, a transversely intersecting pair of planes, or a (stationary) twisted union of $4$ half-planes, obeying
		$$\dist_\H(\Cbf_x\cap B_1(0), P_0\cap B_1(0))\leq C\hat{E}_V;$$
		\item [(3)] For every $x,y\in \sing(\widetilde{V})\cap C_{1/2}(0)$ we have
		$$\dist_\H(\Cbf_x\cap B_1(0), \Cbf_y\cap B_1(0))\leq C\hat{E}_V|x-y|^\alpha;$$
		\item [(4)] For all $x\in \sing(\widetilde{V})\cap C_{1/2}(0)$ we have
		$$r^{-n-2}\int_{B_r(x)}\dist^2(y,x+\BC_x)\, \ext\|V\|(y) \leq C\hat{E}_V^2r^{2\alpha} \qquad \text{for all }r\in (0, 1/4).$$
	\end{enumerate}
	Here, $C = C(n,k)\in (0,\infty)$ and $\alpha = \alpha(n,k)\in (0,1)$.
\end{theorem}

\textbf{Remark:} The fact that the constants $C,\alpha$ in Theorem \ref{thm:main} do not depend on $\beta$ follows once you have the conclusion that $V\res C_{1/2}(0)$ is a $2$-valued stationary Lipschitz graph, as then such graphs belong to $\mathcal{V}_{\beta}$ for \emph{any} $\beta\in (0,1)$ (cf.~Theorem \ref{thm:main-2}). So one may fix a choice of $\beta$ and reapply the theorem again to remove the constant dependencies on $\beta$ (up to the standard argument that one may change the radius of the cylinder in which the conclusion holds by applying the result to suitable homotheties of $V$).

\textbf{Note:} A trivial way in which $V\in\mathcal{V}_\beta$ for \emph{every} $\beta\in (0,1)$ is if $V$ has a point of density $\geq 2$ in \emph{every} cylinder; in this case, the result and proof essentially reduce to that of Allard's regularity theorem, with the conclusion being that $V\res C_{1/2}(0)$ is represented by a single-valued $C^{1,\alpha}$ graph with multiplicity $2$. Notice however that, a priori, just because $V$ has a point of density $\geq 2$ in every cylinder does not mean that $\Theta_V\geq 2$ $\H^n$-a.e; as such, one cannot simply apply Allard's regularity theorem to the varifold with the same support as $V$ and half the density function in order to deduce Theorem \ref{thm:main} in this case, as this might not be an integral varifold.

Thus, Theorem \ref{thm:main} says that for any $\beta\in (0,1)$ the only condition preventing a stationary integral varifold $V$ which is $\eps_0$-close to a multiplicity two plane (for $\eps_0$ as in Theorem \ref{thm:main}) from being a two-valued graph is if $V\not\in \mathcal{V}_\beta$. This gives as a corollary the following dichotomy for \emph{general} stationary integral varifolds close to a multiplicity $2$ plane:

\begin{corollary}\label{cor:main-2}
    Fix $\beta\in (0,1)$. Then there exists $\eps_0 = \eps_0(n,k,\beta)\in (0,1)$ such that the following is true. Suppose $V$ is a stationary integral varifold in $B^{n+k}_2(0)$ which obeys
    $$\frac{3}{2}\leq\frac{\|V\|(B_1(0))}{\w_n}\leq\frac{5}{2} \qquad \text{and} \qquad \hat{E}_V<\eps_0.$$
    Then, one of the following must hold:
    \begin{enumerate}
        \item [(1)] $V\in \mathcal{V}_\beta$, and so the conclusions of Theorem \ref{thm:main} hold;
        \item [(2)] $V\not\in \mathcal{V}_\beta$, i.e.~there is a ball $B_\rho(x_0)\subset B^n_1(0)$ such that on the cylinder $C_\rho(x_0)$ we have $\{\Theta_V<2\}\subset C_\rho(x_0)$, $\frac{3}{2}\leq (\w_n\rho^n)^{-1}\|V\|(C_\rho(x_0))\leq\frac{5}{2}$, and for some $z\in P_0^\perp$,
        $$\frac{1}{\rho^{n+2}}\int_{C_\rho(x_0)}\dist^2(x,\{z\}\times\R^n)\, \ext\|V\|(x) < \beta^2,$$
         yet $\spt\|V\|\cap C_{\rho/2}(x_0)$ is connected, or at most one of the connected components of $\spt\|V\|\cap C_{\rho/2}(x_0)$ has non-empty intersection with $C_{\rho/4}(x_0)$.
    \end{enumerate}
\end{corollary}

\begin{remark}
	In particular, if in Theorem \ref{cor:main-2} $V\res C_{1/2}(0)$ does \emph{not} resemble a $2$-valued graph over $B^n_{1/2}(0)$, we \emph{must} have $V\not\in\mathcal{V}_\beta$. This happens for example when for some $x\in B^n_{1/2}(0)\subset P_0$ we have
$$\sum_{y\,\in\,\pi_{P_0}^{-1}(x)\cap \spt\|V\|}\Theta_V(y) \neq 2.$$
One instance this happens is when $V\res C_{1/2}(0)$ contains a triple junction singularity. Another is when there is a ray $\R^k\times\{x\}\subset C_{1/2}(0)$ which \emph{does not} intersect $\spt\|V\|$. In both these scenarios, we must get an \emph{entire cylinder} in which $V$ is close to a plane of multiplicity $2$, $\Theta_V$ is $<2$ in the cylinder, yet $V$ has only one connected component of significant size in the cylinder. These would seem to be non-trivial facts, as the presence or absence of certain structures implies the existence of very controlled flat regions in $V$ which is failing an $\eps$-regularity theorem, which a priori the relationship between is unclear.
\end{remark}

We therefore see that \emph{if} one works with a class of varifolds $\mathscr{V}$ where one has already established the special case of the $\eps$-regularity theorem (for suitable $\eps = \eps(\mathscr{V})$) when one \emph{additionally assumes} that $\Theta_V<2$ everywhere (and thus in which one expects to have more information on the possible singularities of $V$), then $\mathscr{V}\subset \mathcal{V}_\beta$ for suitable $\beta = \beta(\mathscr{V})$ and hence, using Theorem \ref{thm:main}, one can deduce, for the class $\mathscr{V}$, the general $\eps$-regularity theorem close to a multiplicity $2$ plane. This is the key observation underpinning our applications of Theorem \ref{thm:main} to certain classes of stationary integral varifolds, such as to $2$-valued stationary Lipschitz graphs in the following:

\begin{theorem}\label{thm:main-2}
	Fix $L\in (0,\infty)$. Suppose $f:B_2^n(0)\to \A_2(\R^k)$ is a Lipschitz $2$-valued function with $\Lip(f)\leq L$ and such that the integral varifold $V = \mathbf{v}(f)$ associated to the graph of $f$ is stationary. Then, $V\in\mathcal{V}_\beta$ for any $\beta\in (0,1)$.
	
	In particular, Theorem \ref{thm:main} applies to $V$ (with e.g.~$\beta=1/2$), i.e.~there exists $\eps_0 = \eps_0(n,k)\in (0,1)$ such that if such $V$ obeys
	$$\frac{3}{2}\leq \frac{\|V\|(B_1(0))}{\w_n}\leq\frac{5}{2} \qquad \text{and} \qquad \hat{E}_V<\eps_0,$$
	then in fact $f|_{B_{1/2}^n(0)}$ is $GC^{1,\alpha}$, and conclusions (1) -- (4) of Theorem \ref{thm:main} hold.\end{theorem}

Let us now consider the special case where a stationary integral varifold $V$ has a singular point $x$ where one tangent cone is a multiplicity $2$ plane (which, after rotating, we can without loss of generality assume to be $2|P_0|$). Then, the only way $V$ is \emph{not} a $2$-valued stationary Lipschitz graph about $x$ is if there exist a \emph{sequence} of points $x_j\to x$ and radii $\rho_j\to 0$ such that $C_{\rho_j}(x_j)\subset\{\Theta_V<2\}$, $\frac{3}{2}\leq (\w_n \rho_j^n)^{-1}\|V\|(C_{\rho_j}(x_j))\leq \frac{5}{2}$, and for some $z_j\in P_0^\perp$ we have
$$\frac{1}{\rho_j^{n+2}}\int_{C_{\rho_j}(x_j)}\dist^2(x,\{z_j\}\times\R^n)\, \ext\|V\|(x) \to 0,$$
yet for each $j$ we have that either $\spt\|V\|\cap C_{\rho_j/2}(x_j)$ is connected, or at most one connected component of $\spt\|V\|\cap C_{\rho_j/2}(x_j)$ has non-empty intersection with $C_{\rho_j/4}(x_j)$. In particular, there are infinitely many cylinders contained in the set $\{\Theta_V<2\}$ for which a connected component of $V$ in each cylinder is close, as varifolds, to some plane of multiplicity $2$ parallel to $P_0$.

Using this, we are able to deduce structural information about the behaviour of a codimension one stationary integral varifold which has \emph{infinite Morse index} on its regular part about a singular point $x$ where one tangent cone is a plane with multiplicity $2$. Indeed, this corollary follows since a two-valued $C^{1,\alpha}$ graph in codimension one has stable (i.e.~Morse index $0$) regular part (see e.g.~\cite{Hie20}), and thus $V\res B_{\rho}(x)$ can never be a $2$-valued stationary Lipschitz graph for any $\rho>0$, meaning the above discussion applies. We therefore have:

\begin{corollary}\label{cor:main}
	Suppose $V$ is a stationary integral $n$-varifold in $B^{n+1}_2(0)$ with $0 \in {\rm spt} \, \|V\|$ such that the Morse index of $\reg(V)\cap B_\rho(0)$ is infinite for every $\rho>0$. Suppose also that one tangent cone to $V$ at 
    $0$ is $2|\{0\}\times \R^n|$. Then, there exist a sequence of points $x_j\to 0$ with $x_{j} \neq 0$ and radii $\rho_j\in (0,|x_j|)$ such that $C_{\rho_j}(x_j)\subset \{\Theta_V<2\}$, $\frac{3}{2}\leq (\w_n \rho_j^n)^{-1}\|V\|(C_{\rho_j}(x_j))\leq \frac{5}{2}$, and 
	$$\frac{1}{\rho_j^{n+2}}\int_{C_{\rho_j}(x_j)}\dist^2(x,\{z_j\}\times\R^n)\, \ext\|V\|(x) \to 0$$
	for some $z_j\in \R^k\times\{0\}^n$, and yet for each $j$ we have that either $\spt\|V\|\cap C_{\rho_j/2}(x_j)$ is connected, or at most one connected component of $\spt\|V\|\cap C_{\rho_j/2}(x_j)$ has non-empty intersection with $C_{\rho_j/4}(x_j)$.
\end{corollary}

This corollary appears to be non-obvious: for instance, a priori the index could be entirely contained within small smoothly embedded regions where the varifold has constant multiplicity two. Corollary \ref{cor:main} therefore says that, whilst the index might have such regions contributing towards it, there \emph{must} always be infinitely many regions of the above form. Each such region must either contain a triple junction singularity or must have index $\geq 1$ from the results in \cite{MW24} (as if there are no triple junctions and the index were zero, the varifold in such a cylinder would decompose as two minimal graphs by \cite[Theorem A]{MW24}, and thus be disconnected). Thus, if there are no triple junctions, there must be an infinite contribution to the index from the set $\{\Theta_V<2\}$; again, we remark that triple junctions cannot occur if the varifold arises as the varifold limit of smoothly embedded minimal surfaces. Also note that in this corollary the cylinders are \emph{all} over the plane $P_0\equiv \{0\}\times\R^n$ which was the tangent cone. If there were multiple different (multiplicity 2) planes arising as tangent cones, we would have the same conclusion on a sequence of cylinders over \emph{each} such plane.

\begin{remark}
    Suppose that $V$ is a stationary integral $n$-varifold in $B^{n+1}_2(0)$ with $0\in \spt\|V\|$ and such that $V$ has no triple junction singularities. Suppose that one tangent cone to $V$ at $0$ is $2|\{0\}\times\R^n|$. If $\reg(V)$ has finite Morse index on some ball $B_\rho(0)$, then a simple cut-off argument shows that there must exist $\rho^\prime$ such that $\reg(V)\cap B_{\rho^\prime}(0)$ has zero Morse index, at which point one may apply the main result of \cite{MW24} to deduce that the local structure about $0$ is given by the graph of a $C^{1,1/2}$ $2$-valued function. Thus, the only situation we do not understand is when $\reg(V)\cap B_\rho(0)$ has infinite Morse index for each $\rho>0$. Corollary \ref{cor:main} therefore gives an equivalent structural condition to this analytic condition, which can be viewed as a generalisation to arbitrary dimensions of the equivalence between infinite Morse index for minimal surfaces and genus.
\end{remark}

\begin{remark}
	Even if the varifold $V$ corresponds to a (multiplicity one) smoothly embedded minimal hypersurface with an isolated singularity with one tangent cone being a plane with multiplicity $2$, Theorem \ref{thm:main} still has non-trivial content. Indeed, suppose $V = |M|$, where $M$ is a smoothly embedded minimal hypersurface in $B^{n+1}_1(0)\setminus\{0\}$, and that $0\in\sing(V)\equiv\overline{M}\setminus M$ is such that one tangent cone is a multiplicity $2$ plane (which without loss of generality is $2|\{0\}\times\R^n|$). Then, even if we knew that $V\in \mathcal{V}_\beta$ for some $\beta\in (0,1)$, it is not a priori clear that $V$ corresponds to a two-valued graph about $0$ without utilising Theorem \ref{thm:main} (the argument partially simplifies in that the regularity of coarse blow-ups in Part \ref{part:blow-up-reg} is much easier, however as the new cone in an excess decay lemma can still be the union of two distinct planes, it seems that one cannot avoid the fine blow-up procedure in Part \ref{part:fine-reg}). Indeed, the only conclusion one can seemingly obtain directly from $V\in\mathcal{V}_\beta$ is that there is a sequence of annuli regions, $B_{\rho_i}(0)\setminus\overline{B}_{\tau_i}(0)$, where $\tau_i/\rho_i\to 0$, such that $(\eta_{0,\rho_i})_\#V$ is a two-valued graph over $B_1^n(0)\setminus \overline{B}^n_{\tau_i/\rho_i}(0)$; it is not clear that one can ever extend the two-valued graph structure across $B_{\tau_i/\rho_i}(0)$ without using Theorem \ref{thm:main}, as indeed cylinders centred at points in $B_{\tau_{i}/\rho_{i}}(0)$ that avoid $0$ (the unique density 2 point) must have small radii if the centre is close to $0$, but one does not necessarily know smallness of the $L^2$ excess of $V$ in such small cylinders in order to apply the defining property of $\mathcal{V}_\beta$.; indeed, this would be the case if one does not know uniqueness of the tangent cone (if one knew uniqueness of the tangent plane, a simple extension argument would determine the local structure as a $2$-valued graph). Moreover, after using Theorem \ref{thm:main} to get that $V$ is in fact a 2-valued graph about $0$, since the codimension is $1$ and $V$ is embedded away from $0$, we see that the $2$-valued graph globally decomposes as a sum of two single-valued minimal graphs, for which $0$ is removable, and then from the maximum principle we would get that $V$ coincides with a smoothly embedded minimal surface with multiplicity $2$, a contradiction. Thus, we see that in this setting we must also have $V\not\in\mathcal{V}_\beta$ for all $\beta\in (0,1)$, leading to the same conclusion as in Corollary \ref{cor:main} in this setting. Again, each such region must therefore contain index.
\end{remark}

Our proof of Theorem \ref{thm:main} will be the consequence of an excess decay-type lemma, reminiscent of that seen in the proof of Allard's regularity theorem, however more complicated due to the possibility that the new cone for which the excess has decayed towards might not be a single plane but instead a union of planes or even a (twisted) union of $4$ half-planes. In fact, our excess decay theorem will instead be phrased as a \emph{dichotomy}: either we get decay to a new \emph{plane}, or (in all other cases when the decay would instead be towards a cone in $\CC$) we actually automatically get all the conclusions of Theorem \ref{thm:main}. An excess decay dichotomy of this form was used in the work \cite{MW24} and will follow from an appropriate understanding of coarse blow-ups (see Part \ref{part:coarse-blow-ups} and Part \ref{part:blow-up-reg}). The advantage of such a dichotomy means that at each stage of the iteration process, if we are not in the second case (where iteration can stop), we are decaying to another multiplicity $2$ plane, from which it is simple to verify the hypotheses to iterate (up to $\beta$ changing by a controlled amount as seen in Proposition \ref{prop:rotation-class}). Establishing the second case of this dichotomy utilises a so-called \emph{fine $\eps$-regularity theorem} (this is the focus of Part \ref{part:fine-reg}, see Theorem \ref{thm:fine-reg}), which furthermore is also needed to establish the necessary regularity conclusion for coarse blow-ups.

Even more generally, our proof provides the following \emph{tri}chotomy for \emph{arbitrary} stationary integral varifolds close to a multiplicity $2$ plane. This is more similar in spirit to the main dichotomy in Simon's work on cylindrical tangent cones \cite[Lemma 1]{Sim93}, and is the more precise statement of Theorem \ref{thm:main-gap-intro}.

\begin{theorem}\label{thm:main-3}
	Fix $\beta>0$ and $\gamma>0$. Then, there exists $\eps = \eps(n,k,\beta,\gamma)\in (0,1)$, $\eta = \eta(n,k,\beta,\gamma)\in (0,1)$, and $\rho_0 = \rho_0(n,k,\beta,\gamma)\in (0,1)$ such that the following is true. Suppose that $V$ is a stationary integral $n$-varifold in $B^{n+k}_2(0)$ which obeys
	\begin{itemize}
		\item $\Theta_V(0)\geq 2$;
		\item $\frac{3}{2}\leq \w_n^{-1}\|V\|(B_1(0))\leq \frac{5}{2}$;
		\item $\hat{E}_V<\eps$.
	\end{itemize}
	Then, one of the following three alternatives must occur:
	\begin{enumerate}
		\item [(i)] there exists a cone $\BC \in \widetilde{\CC}_0$ and a scale $\theta = \theta(n,k,\beta,\gamma)\in (0,1/4)$ such that
		$$\frac{1}{\theta^{n+2}}\int_{C_\theta(0)}\dist^2(x,\BC)\, \ext\|V\|(x) \leq \frac{1}{4}\hat{E}_V^2;$$
		$$\dist_\H(\BC\cap B_1(0), P_0\cap B_1(0))\leq C\hat{E}_V;$$
		\item [(ii)] there is $x\in P_0$ and $\rho \geq \rho_0$ with $C_\rho(x)\subset C_1(0)$ obeying:
			\begin{itemize}
				\item $C_\rho(x)\subset\{\Theta_V<2\}$;
				\item $\frac{3}{2}\leq (\w_n\rho^n)^{-1}\|V\|(C_\rho(x))\leq \frac{5}{2}$;
				\item $\hat{E}_{V}(C_\rho(x))<\beta$;
				\item for any harmonic function $u:B_\rho(x)\to \R^k$ we have
		$$\frac{1}{\rho^{n+2}}\int_{C_\rho(x)}\dist^2(y,\graph(u))\, \ext\|V\|(y)>\eta\hat{E}_V^2;$$
				\item $V$ fails the topological structural condition in $C_\rho(x)$;
			\end{itemize}
		\item [(iii)] $V$ has a $(\beta,\gamma)$-fine gap (of possibly very small radius) in which $V$ fails the topological structural condition.
	\end{enumerate}
    Here, $C = C(n,k,\beta,\gamma)\in (0,\infty)$.
\end{theorem}
This should be viewed as a refinement of Corollary \ref{cor:main} in that we know much more about the $\beta$-coarse gap failing the topological structural condition. Indeed, it must either be (1) a $(\beta,\gamma)$-fine gap (of possibly very small radius), or (2) $\beta$-coarse gap which has a definite lower bound on its radius, is close to $2|P_0|$ at the scale of the cylinder (rather than just being close to a multiplicity two plane \emph{parallel} to $P_0$), \emph{and} within the coarse gap $V$ is far, relative to its $L^2$ height excess at scale $1$, from \emph{any} harmonic graph. Practically speaking, what this latter condition means is that any coarse blow-up must separate \emph{somewhere}, i.e.~it prevents the coarse blow-up being a single harmonic function with multiplicity $2$.

As mentioned in Remark \ref{remark:d}, our methods are also able to establish analogues of Theorem \ref{thm:main} and Theorem \ref{thm:main-3} when instead the stationary integral varifold is close to a cone in $\CC$, provided one makes a suitable adjustment to the notions of $\beta$-coarse gap and $(\beta,\gamma)$-fine gap. In fact, this argument is strictly simpler than that needed in the proof of Theorem \ref{thm:main} and Theorem \ref{thm:main-3} (indeed, it more or less follows the proof as in \cite{BK17}), as the main technical result needed (Theorem \ref{thm:fine-reg}) in the planar case may be viewed as a degenerate version of the situation when the varifold is close to a cone in $\CC$. We will discuss this situation in Section \ref{sec:non-degenerate} at the end of the paper to avoid the possibility of creating confusion by introducing specialised terminology to that setting which won't be used in the rest of the paper.

\section{Elementary Observations on the Structural Condition}\label{sec:tsc}

The aim of this section is two-fold: we first wish to give some lemmas concerning the behaviour of the varifolds in the cylinders as described in Definition \ref{defn:V-beta}. Then, we wish to prove Theorem \ref{thm:main-2} showing that $2$-valued stationary Lipschitz graphs belong to the classes $\mathcal{V}_\beta$ for $\beta\in (0,1)$.

\subsection{The Structural Condition}

We give two lemmas showing how the defining condition of $\mathcal{V}_\beta$ gives structural information on the varifold. These lemmas will be of vital importance in the paper.

\begin{lemma}\label{lemma:gap-1}
	There exists $\eps_0 = \eps_0(n,k)\in (0,1)$ such that the following is true. Suppose $V$ is a stationary integral $n$-varifold in $C_1(0)$ with:
	\begin{enumerate}
		\item [(1)] $\|V\|(C_1(0)) \leq \frac{5}{2}\w_n$;
		\item [(2)] $C_1(0)\subset\{\Theta_V<2\}$;
		\item [(3)] $\hat{E}_V<\eps_0$;
		\item [(4)] $\spt\|V\|\cap C_{1/2}(0)$ is (non-empty) disconnected, and at least two of its components have non-empty intersection with $C_{1/4}(0)$.
	\end{enumerate}
	Then, there are two smooth functions $f^1,f^2:B^n_{1/8}(x_0)\to \R^k$ solving the minimal surface system such that:
	\begin{enumerate}
		\item [(i)] $V\res C_{1/8}(0) = |\graph(f^1)| + |\graph(f^2)|$;
		\item [(ii)] $\graph(f^1)\cap \graph(f^2) = \emptyset$;
		\item [(iii)] $\|f^i\|_{C^3(B^n_{1/8}(0))}\leq C\hat{E}_V$ for some $C = C(n,k)$.
	\end{enumerate}
\end{lemma}

\begin{proof}
	We work by contradiction. If the lemma were false, then we could find sequences $\eps_j\downarrow 0$ and $(V_j)_j$ obeying the assumptions of the lemma with $V_j,\eps_j$ in place of $V$ and $\eps_0$, respectively, yet the conclusions fail. From the standard compactness theorem and interior height bound for stationary integral varifolds, and the constancy theorem, as well as assumptions (1),(3) and (4), we can pass to a subsequence to ensure that $V_j\weakly \theta|P_0|$ in $C_1(0)$ for some constant $\theta\in\{1,2\}$. However, if $\theta=1$ we could then apply Allard's regularity theorem to contradict (4). Thus, we must have $\theta=2$.
	
	Now for each $j$ let $(S_{j,i})_{i=1}^{N_j}$ be the connected components of $\spt\|V_j\|\cap C_{1/2}(0)$. By (4), we know $N_j\geq 2$, and by (1) and the monotonicity formula we know that $N_j<\infty$. Furthermore, let $(\tilde{S}_{j,i})_{i=1}^{M_j}$ be the connected components of $\spt\|V_j\|\cap C_{1/4}(0)$; again, our assumptions give $2\leq M_j<\infty$. Assumption (4) tells us that $S_{j,1}$ and $S_{j,2}$ have non-empty intersection with $C_{1/4}(0)$, and so, for instance, $S_{j,1}\cap C_{1/4}(0)$ must be comprised of a union of some of the $\tilde{S}_{j,i}$. Therefore, simply by relabelling, we may assume that $\tilde{S}_{j,1}\subset S_{j,1}$ and $\tilde{S}_{j,2}\subset S_{j,2}$.
	
	Now since $\spt\|V_j\|$ is closed, it is easy to check that $(|S_{j,i}|)_{j=1}^\infty$ is a sequence of stationary integral $n$-varifolds in $C_{1/2}(0)$ for each fixed $i$. From the convergence of $V_j$ to $\theta|P_0|$ in $C_1(0)$, we must have that $|S_{j,i}|\weakly Q_i|P_0|$ in $C_{1/2}(0)$, where necessarily $Q_i\in \{0,1,2\}$.
	
	Since $\tilde{S}_{j,i}\subset S_{j,i}$ for $i=1,2$, the monotonicity formula formula gives $\|V\|(S_{j,i})\geq \w_n(1/4)^n$ for $i=1,2$, and thus we deduce that $Q_i\geq 1$ for $i=1,2$. Since $\sum_i Q_i=2$, we must then have $Q_1=Q_2=1$ and $Q_i=0$ for $i\geq 3$. As $|S_{j,1}|$ and $|S_{j,2}|$ are stationary integral varifolds converging to a multiplicity $1$ plane, we can apply Allard's regularity theorem to these for $j$ sufficiently large in the region $C_{1/3}(0)$ (say). From this it follows that we must have $M_j=2$ for all $j$ large: indeed, if not a similar argument to the above based on the monotonicity formula, in the cylinder $C_{1/3}(0)$, would contradict the convergence to a plane of multiplicity $2$. Thus, we have reached the desired contradiction, completing the proof.	
\end{proof}

Now we may rephrase Lemma \ref{lemma:gap-1} for the classes $\mathcal{V}_\beta$ as follows.

\begin{lemma}\label{lemma:gap-2}
	Fix $\beta\in (0,1)$ and $\rho_0\in (0,1/4)$. Then, there exists and $\eps_0 = \eps_0(n,k,\beta,\rho_0)\in (0,1)$ such that the following is true. Suppose $V\in \mathcal{V}_\beta$ is such that
	$$\frac{3}{2}\leq\frac{\|V\|(B_1(0))}{\w_n}\leq \frac{5}{2} \qquad \text{and} \qquad \hat{E}_V<\eps_0.$$
	Then, if $C_\rho(x_0)\subset C_1(0)$ is a cylinder ($x_0\in P_0$, $\rho>0$) for which:
	\begin{itemize}
		\item $\Theta_V<2$ in $C_\rho(x_0)$;
		\item $\rho\geq\rho_0$.
	\end{itemize}
	Then, there are two smooth functions $f^1,f^2:B^n_{\rho/8}(x_0)\to \R^k$ solving the minimal surface system such that:
	\begin{enumerate}
		\item [(i)] $V\res C_{\rho/8}(x_0) = |\graph(f^1)| + |\graph(f^2)|$;
		\item [(ii)] $\graph(f^1)\cap \graph(f^2) = \emptyset$;
		\item [(iii)] $\|f^i\|_{C^3(B^n_{\rho/8}(x_0))}\leq C\hat{E}_V(C_\rho(x_0))$, where $C = C(n,k,\rho_0)\in (0,\infty)$.
	\end{enumerate}
\end{lemma}

\textbf{Remark:} Alternatively, in (iii) above we can take the scale-invariant $C^3$ norm on $B^n_{\rho/8}(x_0)$ and the constant $C$ to depend  only on $n,k$.

\begin{proof}
	We argue by contradiction. If the lemma were false, then for fixed $\beta\in (0,1)$ and $\rho_0\in (0,1/4)$, we can find $\eps_j\downarrow0$, varifolds $V_j\in \mathcal{V}_\beta$ obeying
	$$\frac{3}{2}\leq\frac{\|V_j\|(B_1(0))}{\w_n}\leq\frac{5}{2} \qquad \text{and} \qquad \hat{E}_{V_j}<\eps_j,$$
	and cylinders $C_{\rho_j}(x_j)\subset C_1(0)$ with $x_j\in P_0$, $\rho_j\geq\rho_0$, and $\Theta_{V_j}<2$ in $C_{\rho_j}(x_j)$, yet the conclusions with $\rho_{j}$ in place of $\rho_{0}$ and $x_{j}$ in place of $x_{0}$ fail. In particular, as $\rho_j\geq \rho_0$, we may pass to a subsequence so that $x_j\to x_0\in B^n_1(0)$ (the uniform lower bound on $\rho_j$ guarantees $x_0\in B_{(1+\rho_0)/2}^n(0)\subset B^n_1(0)$ rather than just in $\overline{B^n_1(0)}$). Note also that $V_j\weakly 2|P_0|$ in $B_1(0)$ by the conditions on $V_j$ above.
	
	Furthermore, we have
	$$\hat{E}_{V_j}(C_{3\rho_j/4}(x_j)) \leq (4/3)^{n+2}\rho_{0}^{-n-2}\int_{C_{\rho_j}(x_j)}\dist^2(x,P_0)\, \ext\|V_j\|(x) \leq (4/3)^{n+2}\rho_0^{-n-2}\eps_j$$
	and thus $\hat{E}_{V_j}(C_{3\rho_j/4}(x_j))\to 0$. Thus, for all $j$ sufficiently large we have, with $\tilde{\rho}_j:= 3\rho_j/4$:
	\begin{itemize}
		\item $\Theta_{V_j}<2$ in $C_{\tilde{\rho}_j}(x_j)$;
		 \item $\frac{3}{2}\leq (\w_n\rho_j^n)^{-1}\|V_j\|(C_{\tilde{\rho}_j}(x_j))\leq \frac{5}{2}.$
		\item $\hat{E}_{V_j}(C_{\tilde{\rho}_j}(x_j))<\beta$;
	\end{itemize}
	The second condition comes from the fact that $V_j\weakly 2|P_0|$ and there is a uniform positive lower bound on $\tilde{\rho}_j$. Hence, by definition of $\mathcal{V}_\beta$, we get that $\spt\|V_j\|\cap C_{\tilde{\rho}_j/2}(x_j)$ is (non-empty) disconnected, and at least two of its connected components have non-empty intersection with $C_{\tilde{\rho}_j/4}(x_j)$. So, if we consider $\tilde{V}_j:= (\eta_{x_j,\tilde{\rho}_j})_\#V_j$, we see that for all $j$ sufficiently large we can apply Lemma \ref{lemma:gap-1} to get the result (the mass ratio control coming from the fact that $V_j\weakly 2|P_0|$). 
	
	(Technically, this would give the conclusions over the ball $B_{3\rho/32}(x_0)$ instead of $B_{\rho/8}(x_0)$, however suitably changing the radii in the conclusion of Lemma \ref{lemma:gap-1} and the radius $3/4$ in the above argument allows us to get the conclusion as stated.)
\end{proof}

\textbf{Note:} In Lemma \ref{lemma:gap-2} we in fact only needed the defining property of $\mathcal{V}_\beta$ in $\beta$-coarse gaps when the $L^2$ height excess is small relative to $P_0$, rather than a plane parallel to $P_0$.

We can now make precise a previous remark made following Definition \ref{defn:V-beta} regarding the equivalence of the definition of $\mathcal{V}_\beta$ with a definition through a sheeting property. Indeed, for $\beta\in (0,1)$, define $\widetilde{\mathcal{V}}_\beta$ to be the set of all stationary integral $n$-varifolds $V$ in $B_2(0)$ obeying the following condition: whenever $x\in P_0$ and $\rho>0$ are such that $C_\rho(x)\subset C_1(0)$ and:
\begin{enumerate}
    \item [(a)] $C_\rho(x)\subset\{\Theta_V<2\}$;
    \item [(b)] $\frac{3}{2}\leq (\w_n\rho^n)^{-1}\|V\|(C_\rho(x))\leq \frac{5}{2}$; 
    \item [(c)] For some $z\in \R^k\times\{0\}^n$ we have
    $$\frac{1}{\rho^{n+2}}\int_{C_\rho(x)}\dist^2(y,\{z\}\times\R^n)\, \ext\|V\|(y) < \beta^2;$$
\end{enumerate}
then $V\res C_{\rho/2}(x) = \mathbf{v}(f^1) + \mathbf{v}(f^2)$, where $f^1,f^2:B^n_{\rho/2}(x)\to \R^k$ are smooth functions solving the minimal surface system.

We then have the following:

\begin{prop}\label{prop:sheeting-equivalence}
    Fix $\beta\in (0,1)$. Then for suitable $\widetilde{\beta} = \widetilde{\beta}(n,k,\beta)\in (0,\beta)$ we have
    $$\widetilde{\mathcal{V}}_\beta \subset \mathcal{V}_\beta \subset \widetilde{\mathcal{V}}_{\widetilde{\beta}}.$$
\end{prop}
\begin{proof}
    The inclusion $\widetilde{\mathcal{V}}_\beta\subset\mathcal{V}_\beta$ is trivial: indeed, if $V \in \widetilde{\mathcal{V}}_{\beta}$ and $C_\rho(x)$ is a cylinder as in Definition \ref{defn:V-beta}, then the defining property of $\widetilde{\mathcal{V}}_{\beta}$ implies the conditions for $V$ to be in $\mathcal{V}_\beta$.

    Now let $V\in\mathcal{V}_\beta$ and fix $\widetilde{\beta} = \widetilde{\beta}(n,k,\beta)\in (0,\beta)$ to be determined. Suppose that $C_\rho(x)$ is a cylinder as in the definition of $\widetilde{\mathcal{V}}_{\widetilde{\beta}}$, which without loss of generality by translating and rescaling we can assume is $C_1(0)$ and $z=0$. Lemma \ref{lemma:gap-2} with $\rho_0 = 1/8$ therein gives the existence of a constant $\eps_0 = \eps_0(n,k,\beta)\in (0,1)$ so that, if $\widetilde{\beta}<\eps_0$, then we have $V\res C_{1/64}(0) = \mathbf{v}(f^1) + \mathbf{v}(f^2)$ for two smooth functions $f^1,f^2:B^n_{1/64}(0)\to \R^k$ solving the minimal surface equation. Our assumptions then allow us to repeat this to $\widetilde{V}:= (\eta_{(0,y),1/4})_\#V$, for any $y\in B^n_{1/2}(0)$, provided $\widetilde{\beta} = \widetilde{\beta}(n,k,\beta)$ is sufficiently small (indeed, we need to ensure that $V$ is sufficiently close to a multiplicity $2$ plane as varifolds here to ensure mass bounds on the new cylinders). After patching together the functions acquired through this (using unique continuation on overlaps) we therefore get $V\res C_{1/2}(0) = \mathbf{v}(g^1) + \mathbf{v}(g^2)$ where $g^1,g^2:B^n_{1/2}(0)\to \R^k$ are smooth solutions to the minimal surface equation. This therefore shows that $V\in\widetilde{\mathcal{V}}_{\widetilde{\beta}}$, completing the proof.
\end{proof}

Next, we consider the impact of a small rotation on a varifold in $\mathcal{V}_\beta$. As included in the definition of $\mathcal{V}_\beta$, this requires the validity of the topological structural condition in $\beta$-coarse gaps when the $L^2$ height excess is small relative to some \emph{parallel} plane to $P_0$, rather than just $P_0$ itself.

\begin{prop}\label{prop:rotation-class}
        Fix $\beta\in (0,1)$ and $\delta\in (0,1)$. Then there exist $\eta = \eta(n,k,\beta,\delta)\in (0,1)$ and $\widetilde{\beta} = \widetilde{\beta}(n,k,\beta,\delta)\in (0,\beta)$ for which the following is true. Suppose $V\in \mathcal{V}_\beta$ and $\Gamma\in SO(n+k)$ is a rotation such that $\|\Gamma-\id\|<\eta$. Then we have
        $$(\eta_{0,1-\delta})_\#\Gamma_\#V\in\mathcal{V}_{\widetilde{\beta}}.$$
\end{prop}

\textbf{Note:} The homothety $\eta_{0,1-\delta}$ in Proposition \ref{prop:rotation-class} only plays the role of ensuring that the support of $V$ in cylinders over $\Gamma(P_0)$ centred at points in $B_{1-\delta}(0)$ lies within $C_1(0)$.

\begin{proof}
    Set $P:= \Gamma^{-1}(P_0)$. After performing a rotation and translation (recalling also that $\spt\|V\|\subset B_2(0)$ by definition, and so $\dist(x,P_0) <2$ for all $x\in \spt\|V\|$), the proposition is equivalent to showing that for suitably small $\eta = \eta(n,k,\beta,\delta)\in (0,1)$ and $\widetilde{\beta} = \widetilde{\beta}(n,k,\beta,\delta)\in (0,1)$, if $V\in \mathcal{V}_\beta$ and $\dist_\H(P\cap B_1(0),P_0\cap B_1(0))<\eta$, then for any cylinder $C_\rho(x_0;P)$ with $C_\rho(x_0;P)\cap \spt\|V\|\subset C_1(0)$ (for given $\delta>0$, this last inclusion is guaranteed for $x_0\in B_{1-\delta}(0)$ if $\eta = \eta(n,k,\delta)\in (0,1)$ is sufficiently small) obeying:
    \begin{enumerate}
        \item [(a)] $C_\rho(x_0;P)\subset \{\Theta_V<2\}$;
        \item [(b)] $\frac{3}{2}\leq (\w_n\rho^n)^{-1}\|V\|(C_\rho(x_0;P))\leq \frac{5}{2}$;
        \item [(c)] $\hat{E}_V(C_\rho(x_0;P))<\widetilde{\beta}$;
    \end{enumerate}
    then $\spt\|V\|\cap C_{\rho/2}(x_0;P)$ is disconnected and at least two of its connected components have non-empty intersection with $C_{\rho/4}(x_0;P)$. In turn, if $|y-x|<7\rho/8$, for suitably small $\eta$ and $\widetilde{\beta}$ (depending on $\beta$), by the triangle inequality we can guarantee that $C_{\rho/16}(y;P_0)$ is a $\min\{\beta,\eps_0\}$-coarse gap for $V$, where $\eps_0$ is as in Lemma \ref{lemma:gap-1}, and thus we can apply the defining property of $\mathcal{V}_\beta$, as well as Lemma \ref{lemma:gap-1} (as $\min\{\beta,\eps_0\}\leq\eps_0$) to obtain that $V\res C_{\rho/128}(y;P_0)$ is a sum of two (disjoint) minimal graphs. Repeating this process with different $y$ and patching together the resulting functions means that $V\res C_{\rho/2}(x_0;P)$ is represented by the sum of two disjoint minimal graphs, both of which have non-empty intersection with $C_{\rho/4}(x_0;P)$. This completes the proof.
\end{proof}

\subsection{Lipschitz $2$-Valued Graphs and Theorem \ref{thm:main-2}}

The main observation towards Theorem \ref{thm:main-2} is the following lemma:

\begin{lemma}\label{lemma:Lip-1}
	Let $f:B^n_1(0)\to \A_2(\R^k)$ be Lipschitz and $\Omega\subset B^n_1(0)$ an open and simply connected set. Suppose that $V = \mathbf{v}(f)$ is such that $\sing(V)\cap (\R^k\times \Omega)=\emptyset$. Then, either:
	\begin{enumerate}
		\item [(a)] \textnormal{(Complete Decomposition.)} There are two single-valued $C^1$ functions such that $f^1(x)\neq f^2(x)$ for all $x\in \Omega$ and $V = |\graph(f^1)| + |\graph(f^2)|$;
		\item [(b)] \textnormal{(Complete Coincidence.)} There is a single-valued $C^1$ function $f:\Omega\to \R^k$ such that $V = 2|\graph(f)|$.
	\end{enumerate}
\end{lemma}
\textbf{Note:} This is a fact about $C^1$ $2$-valued functions, as we make no assumption on whether $V$ is stationary. (Again, the singular set is just the usual singular set associated to an integral varifold, so the assumption in the above is that $\spt\|V\|\cap (\R^k\times\Omega)$ is locally a $C^1$ submanifold.)

\begin{proof}
	We first claim that there cannot be both coincidence points and non-coincidence points in $\Omega$. Indeed, suppose that there were. First, notice that the non-coincidence set
	$$U:= \{x\in B^n_1(0):f^1(x)\neq f^2(x)\}$$
	is open by continuity of $f$. We then need to consider (and contradict) two cases depending on whether $\overline{U}\cap U^c\cap \Omega$ is empty or not.
	
	Notice then that the case $\overline{U}\cap U^c\cap \Omega=\emptyset$ leads to an easy contradiction: indeed, if this intersection were empty, then $\Omega$ is the disjoint union of the relatively closed sets $U^c\cap \Omega$ and $\overline{U}\cap \Omega$, which by our original assumption are both non-empty. But $\Omega$ is connected, and so we get a contradiction.
	
	Thus we just need show that $\overline{U}\cap U^c\cap\Omega\neq\emptyset$ leads to a contradiction. In this case, there must be $x_0\in \Omega$ and $\rho>0$ with $B_\rho(x_0)\subset U$, $\overline{B_\rho(x_0)}\subset\Omega$, and $U^c\cap \del B_\rho(x_0)\neq\emptyset$ (indeed, one can simply choose a point $U$ closer to $U^c$ than to $\del\Omega$ and choose the largest ball about this point which lies in $U$). The third condition here guarantees a coincidence point on $\del B_\rho(x_0)$, i.e.~there is $x_1\in \del B_\rho(x_0)$ with $f(x_1) = 2\llbracket y_1\rrbracket$ for some $y_1\in \R^k$.
	
	By our assumption on $\Omega$, we know that $x_1$ is a regular point. This therefore gives a contradiction as at a regular point $V$ is locally a single sheet, whilst $x_1$ is chosen to be on the boundary of separation. More precisely, this can be argued as follows. By assumption on $\Omega$, we know $(y_1,x_1)\in \reg(V)$, and thus for any $\delta>0$ there is $r = r(\delta)>0$ sufficiently small and a $C^1$ function $g:T_{(y_1,x_1)}V\cap B_r(y_1,x_1)\to T^\perp_{(y_1,x_1)}V$ such that $\|g\|_{C^1}<\delta$ and
	$$V\res B_{r/2}(y_1,x_1) = |\graph(g)\cap B_{r/2}(y_1,x_1)|.$$
	In fact, as also $f$ is locally $C^1$ about $(y_1,x_1)$, there is an affine function $\ell:P_0\to P_0^\perp$ with $\text{Lip}(\ell)\leq\Lip(f)$ and $\graph(\ell) = T_{(y_1,x_1)}V$. Thus, we can reparameterise over $P_0$ instead, i.e.~we may instead choose $r = r(\delta)>0$ such that there is a (single-valued) $C^1$ function $h:B^n_r(x_1)\to P_0^\perp$ such that $V\res B_{r/2}(y_1,x_1) = |\graph(h)\cap B_{r/2}(y_1,x_1)|$.
	
	But now by choice of $x_1$ we may pick $x^\prime\in B^n_{r/2}(x_1)\cap U\cap \Omega$ sufficiently close to $x_1$ so that $\{(f^1(x^\prime),x^\prime),(f^2(x^\prime),x^\prime)\}\subset B_{r/2}(y_1,x_1)$ and (as $x^\prime\in U$) $|f^1(x^\prime)-f^2(x^\prime)|>0$. But then this is a contradiction, because the set
	$$\graph(f)\cap \{(y,x^\prime):y\in \R^k\}\cap B_{r/2}(y_1,x_1)$$
	would need to equal both $\{(f^1(x^\prime),x^\prime), (f^2(x^\prime),x^\prime)\}$ and $\{(h(x^\prime),x^\prime)\}$, which have different cardinalities.
	
	Thus, we have reached a contradiction in either case, and so there cannot be both coincidence and non-coincidence points, i.e.~we must have either (1) $\Omega\subset U$ or (2) $\Omega\subset B^n_1(0)\setminus U$. In case (2), then $f^1\equiv f^2$ in $\Omega$, i.e.~$V = 2|\graph(f)|$ for some function $f$, which by assumption must be $C^1$ as $\sing(V)\cap (\R^k\times\Omega)=\emptyset$. Hence in the case of (2), conclusion (b) of the lemma holds.
	
	We are now left with the case (1), i.e.~no point is a coincidence point. We will show that conclusion (a) of the lemma holds. Of course, we can assume $\text{Lip}(f)>0$ or else $f$ is constant and so the lemma is trivially true.
	
	Thus, fix $x_0\in \Omega$ and write $f(x_0) = \llbracket f^1(x_0)\rrbracket + \llbracket f^2(x_0)\rrbracket$. Set
	$$r_{x_0}:= \min\left\{\frac{|f^1(x_0)-f^2(x_0)|}{4\Lip(f)},\, \dist(x_0,\del\Omega)\right\}.$$
	Note that $r_{x_0}>0$. Then, for any $x\in B^n_{r_{x_0}}(x_0)$, let the values of $f(x) = \llbracket f^1(x)\rrbracket + \llbracket f^2(x)\rrbracket$ be ordered such that
	$$\G(f(x),f(x_0)) = |f^1(x)-f^1(x_0)| + |f^2(x)-f^2(x_0)|.$$
	This implies that
	$$|f^1(x)-f^1(x_0)| + |f^2(x)-f^2(x_0)|\leq \text{Lip}(f)|x-x_0|$$
	which, by definition of $r_{x_0}$, means that
	\begin{equation}\label{E:Lipschitz-1}
	x \in B^n_{r_{x_0}}(x_0)\ \ \ \Longrightarrow\ \ \
	\begin{cases}
	|f^1(x)-f^1(x_0)| \leq \frac{1}{4}|f^1(x_0)-f^2(x_0)|\ \ \text{and}\\
	|f^2(x)-f^2(x_0)| \leq \frac{1}{4}|f^1(x_0)-f^2(x_0)|.
	\end{cases}
	\end{equation}
	If $y,z\in B^n_{r_{x_0}}(x_0)$, then the triangle inequality gives
	$$|f^1(x_0)-f^2(x_0)| \leq |f^1(x_0)-f^1(y)| + |f^1(y)-f^2(z)| + |f^2(z)-f^2(x_0)|,$$
	$$|f^1(x_0)-f^2(x_0)| \leq |f^1(x_0)-f^1(z)| + |f^1(z)-f^2(y)| + |f^2(y)-f^2(x_0)|.$$
	Combining these inequalities with \eqref{E:Lipschitz-1} (with $y$ and $z$ in place of $x$ therein) we deduce
	\begin{equation}\label{E:Lipschitz-2}
		\frac{1}{2}|f^1(x_0)-f^2(x_0)| \leq |f^1(y)-f^2(z)| \ \ \ \text{and} \ \ \ \frac{1}{2}|f^1(x_0) - f^2(x_0)| \leq |f^1(z)-f^2(y)|.
	\end{equation}
	Again using the triangle inequality, \eqref{E:Lipschitz-1} and $y$ and $z$ in place of $x$ therein, and \eqref{E:Lipschitz-2}, we have:
	\begin{align*}
		|f^1(y)&-f^1(z)| + |f^2(y)-f^2(z)|\\
		& \leq |f^1(y)-f^1(x_0)| + |f^1(x_0)-f^1(z)| + |f^2(y)-f^2(x_0)| + |f^2(x_0)-f^2(z)|\\
		& \leq |f^1(x_0)-f^2(x_0)|\\
		& \leq |f^1(y) - f^2(z)| + |f^2(y)-f^1(z)|.
	\end{align*}
	This means that
	$$\G(f(y),f(z)) = |f^1(y)-f^1(z)| + |f^2(y)-f^2(z)|$$
	which implies
	$$|f^1(y)-f^1(z)| + |f^2(y)-f^2(z)| \leq \Lip(f)|y-z|.$$
	Hence, we see that $f^1$ and $f^2$ are Lipschitz (single-valued) functions on $B^n_{r_{x_0}}(x_0)$. Since there are no coincidence points or singular points of $\graph(f)$ in $\Omega$, we therefore deduce that $f^1$ and $f^2$ are both $C^1$ functions on $B^n_{r_{x_0}}(x_0)$. Let us call these functions $f^1_{x_0}$ and $f^2_{x_0}$. We have therefore established that for any $z\in \Omega$, there is a radius $r_z>0$ such that $\graph(f)\cap C_{r_z}(z)$ is the disjoint union of two $C^1$ graphs $G^a_z$, $G^b_z$, over $B^n_{r_z}(z)$.
	
	Fix again $x_0\in \Omega$ and take some other point $y\in \Omega$. Since $\Omega$ is open and path-connected, we can find a smooth, injective path $\gamma:[0,1]\to \Omega$ with $\gamma(0) = x_0$ and $\gamma(1) = y$. By the compactness of $\gamma:= \gamma([0,1])$, we have that $\gamma\subset \cup^N_{i=0} B^n_{r_i}(x_i)$, where $x_N = y$ and $\{x_i\}_{i=1}^{N-1}$ are some points on $\gamma$, and $r_i\leq r_{x_i}$ are such that $B^n_{r_i}(x_i)\cap B^n_{r_j}(x_j)\neq\emptyset$ if and only if $|i-j|\leq 1$ (i.e.~only neighbouring balls along the path intersect); we can also arrange so that $r_N = r_{x_N}\equiv r_y$.
	
	Now, for any $K\in \{1,2,\dotsc,N\}$ and smooth functions $f^1,f^2:\bigcup^{K-1}_{i=0}B^n_{r_i}(x_i)\to \R^k$ for which $\left. f^i\right|_{B^n_{r_0}(x_0)} = \left. f^i_{x_0}\right|_{B^n_{r_0}(x_0)}$, and for which $\graph(f)\cap \left(\R^k\times \cup^{K-1}_{i=0}B^n_{r_i}(x_i)\right)$ is the disjoint union of $\graph(f^1)$ and $\graph(f^2)$, we can uniquely define $C^1$ extensions $\tilde{f}^1$, $\tilde{f}^2$, of $f^1$ and $f^2$ defined on the domain $\bigcup^K_{i=0}B^n_{r_i}(x_i)$ for which $\graph(f)\cap\left(\R^k\times \cup^K_{i=0}B^n_{r_i}(x_i)\right)$ is the disjoint union of $\graph(\tilde{f}^1)$ and $\graph(\tilde{f}^2)$. By repeating this inductively on $K$, we establish that given any $y\in \Omega$ and a smooth, injective path $\gamma:[0,1]\to \Omega$ with $\gamma(0) = x_0$ and $\gamma(1) = y$, there are $C^1$ functions defined on $\mathcal{N}_y := \cup^N_{i=0}B^n_{r_i}(x_i)$ with:
	\begin{itemize}
		\item $\left. f^i\right|_{B^n_{r_{0}}(x_0)} = \left. f^i_{x_0}\right|_{B^n_{r_{0}}(x_0)}$ for $i=1,2$;
		\item $f^1(x)\neq f^2(x)$ for all $x\in \mathcal{N}_y$;
		\item $\graph(f)\cap (\R^k\times \mathcal{N}_y) = \graph(f^1)\cup \graph(f^2)$.
	\end{itemize}
	Without loss of generality, we may reorder $G^a_y$ and $G^b_y$ to assume that
	$$\graph(f^1)\cap C_{r_y}(y) = G^a_y \qquad \text{and} \qquad \graph(f^2)\cap C_{r_y}(y) = G^b_y$$
	(where here recall that $r_N\equiv r_y$). Recall also that here $y\in\Omega$ is arbitrary. We claim that this procedure globally defines two $C^1$ functions $f^1$ and $f^2$ on $\Omega$ for which $\graph(f)\cap (\R^k\times\Omega)$ is the disjoint union of $\graph(f^1)\cup \graph(f^2)$.
	
	To see this, let $\tilde{\gamma}:[0,1]\to\Omega$ be a different smooth, injective, path connecting $x_0$ to $y$. By performing the same procedure as above, we have another pair of $C^1$ functions $\tilde{f}^1$, $\tilde{f}^2$ defined on a neighbourhood $\tilde{N}$ or $\tilde{\gamma}$ with $\tilde{f}^i = f^i(x_0)$ on a neighbourhood of $x_0$ for $i=1,2$, $\tilde{f}^1(x)\neq \tilde{f}^2(x)$ for every $x\in\tilde{N}$, and such that
	$$\graph(f)\cap (\R^k\times \tilde{\N}) = \graph(\tilde{f}^1)\cup \graph(\tilde{f}^2).$$
	Assume for the sake of contradiction that this resulting in the opposite labelling of the values of $f(y)$, i.e.~we assume that
	$$\graph(\tilde{f}^1)\cap C_r(y) = G^b_y \qquad \text{and} \qquad \graph(\tilde{f}^2)\cap C_r(y) = G^a_y$$
	for some suitable $r>0$. Write $\Gamma:= \tilde{\gamma}^{-1}\ast\gamma$, where $\ast$ denotes the usual concatenation of curves. This is a piecewise smooth curve with $\Gamma(0) = \Gamma(1) = x_0$. Since $\Omega$ is simply connected, we can continuously contract $\Gamma$ whilst staying in $\Omega$, i.e.~we have a continuous family $\{\Gamma_s\}_{s\in [0,1]}$ of closed loops, all lying in $\Omega$ and with $\Gamma_0 = \Gamma$, such that $\Gamma_1\equiv z_*$ for some $z_*\in\Omega$. For $i\in \{1,2\}$, let $\xi^i:[0,1]\to \graph(f)$ be the $C^1$ curves for which $\graph(f)\cap (\R^k\times\gamma) = \xi^1\cup \xi^2$; in particular, $\xi^i(0) = (f^i(x_0),x_0)$ for $i=1,2$. Also, let $\tilde{\xi}^i:[0,1]\to \graph(f)$ be the $C^1$ curves for which $\graph(f)\cap (\R^k\times\tilde{\gamma}) = \tilde{\xi}^1\cup\tilde{\xi}^2$; again, $\tilde{\xi}^i(0) = (\tilde{f}^i(x_0),x_0)$ for $i=1,2$. By assumption, we have $\tilde{\xi}^1(1) = \xi^2(1)$ and $\tilde{\xi}^2(1) = \xi^1(1)$. Then, there is a piecewise smooth loop $\Xi:[0,1]\to \graph(f)$ given by:
	$$\Xi:= (\tilde{\xi}^1)^{-1}\ast\xi^2\ast (\tilde{\xi}^2)^{-1}\ast\xi^1.$$
	This obeys $\Xi(0) = \Xi(1) = (f^1(x_0),x_0)$ and moreover that $\graph(f)\cap (\R^k\times\Gamma) = \Xi$. By construction, the open sets $\mathcal{N}$, $\tilde{\mathcal{N}}$ (which satisfy $\Gamma\subset\mathcal{N}\cup\tilde{\mathcal{N}}$) and the established properties of the balls $B^n_{r_z}(z)$ give that there is some $\eps>0$ for which there exists a continuous family $\{\Xi_s\}_{s\in [0,\eps)}$ of continuous loops with
	\begin{enumerate}
		\item [(1)] $\Xi_0 = \Xi$;
		\item [(2)] $\Xi_s(0) = \Xi_s(1)$ for every $s\in [0,\eps)$; and
		\item [(3)] $\graph(f)\cap (\R^k\times\Gamma_s) = \Xi_s$ for every $s\in [0,\eps)$.
	\end{enumerate}
	Let $\eps_0\in (0,1]$ denote the supremum of the set of numbers $\eps>0$ for which there exists such a continuous family $\{\Xi_s\}_{s\in [0,\eps)}$ satisfying (1) -- (3) above. We first claim that $\eps_0<1$. Indeed, since there is $r_{z_*}>0$ such that $\graph(f)\cap C_{r_{z_*}}(z_*)$ is the disjoint union of $G^a_{z_*}\cup G^b_{z_*}$ of $C^1$ graphs $G^a_{z_*}$, $G^b_{z_*}$ over $B^n_{r_{z_*}}(z_*)$, we have that for $s$ sufficiently close to $1$, $\graph(f)\cap (\R^k\times\Gamma_s)$ is the disjoint union of two closed curves, one lying entirely in $G^a_{z_*}$ and one lying entirely in $G^b_{z_*}$. This however cannot equal the continuous image of $[0,1]$, and thus cannot be of the form $\Xi_s([0,1])$ for some continuous $\Xi_s$, violating (3). Therefore, we must have $\eps_0\in (0,1)$.
	
	Now take $\eps<\eps_0$. Then any finite covering of $\Gamma_\eps$ of the form $\{B^n_{\bar{r}_i}(y_i)\}_{i=1}^M$ with $y\in \Omega$, $\bar{r}_i<r_{y_i}$, must be such that $\Gamma_{\eps_0}\not\subset\cup^M_{i=1}B^n_{\bar{r}_i}(y_i)$, else we could increase the size of $\eps_0$ (as it is $\eps_0<1$), contradicting it being the supremum. Thus, if we consider a sequence $(\eps_j)_j$ with $\eps_j<\eps_0$ and $\eps_j\uparrow \eps_0$, we deduce that there must be a sequence of points $(y_j)_j$ with $y_j\in\Gamma_{\eps_j}$ and with $r_{y_i}\to 0$. After passing to a subsequence, we may assume that $y_j\to y_0\in \Gamma_{\eps_0}$. Since $r_{y_i}\to 0$, we must have (by definition of $r_{y_i}$) that $|f^1(y_j)-f^2(y_j)|\to 0$, and thus by continuity of $f$ we would get that $y_0$ is a coincidence point of $f$, i.e.~$y_0\not\in U$. But this is a contradiction, as $y_0\in\Omega$ and we were assuming $\Omega\subset U$. Thus, we have reached a contradiction, which shows that the procedure for defining $f^1(y)$ and $f^2(y)$ depends only on $y$ and not on the chosen path $\gamma$. Consequently, we have globally defined $C^1$ functions $f^1$ and $f^2$ for which conclusion (a) holds, completing the proof.
\end{proof}

Now we are in a position to show that Lipschitz minimal graphs always split into two disjoint, single-valued, Lipschitz minimal graphs in regions where the density is $<2$. Notice that this result immediately implies Theorem \ref{thm:main-2}, as the two graphs then give the desired connected components in the smaller cylinders.

\begin{lemma}\label{lemma:Lip-2}
	Let $V$ be a Lipschitz $2$-valued minimal graph over a plane $P$. Suppose that $U\subset P$ is open, simply connected, and that $\pi_P^{-1}(U)\cap \spt\|V\|\subset\{\Theta_V<2\}$. Then, there exist two stationary integral varifolds $V^1$ and $V^2$ which are (single-valued) Lipschitz graphs over $U$, have disjoint support, and are such that $V\res \pi^{-1}_P(U) = V^1+V^2$.
	
	In particular, Theorem \ref{thm:main-2} is true.
\end{lemma}

\textbf{Note:} The single-valued Lipschitz graphs in Lemma \ref{lemma:Lip-2} need not be smooth; however, they will be in the context of topological structural condition given in Definition \ref{defn:V-beta}, provided $\beta = \beta(n,k)\in (0,1)$ is sufficiently small, as an additional smallness-of-height assumption which would allow one to apply Allard regularity to each.

\begin{proof}
	It suffices to prove the case where $U = B$ is an open ball; once one proves this special case (which in fact is all that is needed for Theorem \ref{thm:main-2}), the general case follows by a simple path-uniqueness argument similar to that seen in Lemma \ref{lemma:Lip-1}.
	
	Let $f:P\to \A_2(P^\perp)$ be the $2$-valued Lipschitz graph with $V=\graph(f)$. Write $\S:= \sing(V)\cap \pi_P^{-1}(B)$. By Remark \ref{remark:small-singular-set} before Theorem \ref{thm:main}, as $V$ contains no triple junction singularities and $\Theta_V<2$ in the region of interest, we know that $\dim_\H(\S)\leq n-3$, and so in particular $\dim_\H(\pi_P(\S))\leq n-3$. Thus, we have that $B\setminus \pi_P(\S)$ is simply connected (see, e.g.~\cite[Appendix]{SW16}, for a proof of the fact that the complement in a ball in $\R^n$ of a closed set with $\H^{n-2}$-measure zero is simply connected, which can be applied here to deduce this).
	
	Write $\Omega:= B\setminus \pi_P(\S)$. Then we may apply Lemma \ref{lemma:Lip-1} here: Lemma \ref{lemma:Lip-1}(b) cannot happen as $\Theta_V<2$ everywhere, and so Lemma \ref{lemma:Lip-1}(a) must occur. Thus, there are two $C^1$ functions $f^1$, $f^2:\Omega\to P^\perp$ such that $\graph(f)\cap \pi^{-1}_P(\Omega)$ is the disjoint union of $\graph(f^1)$ and $\graph(f^2)$. Since we know $f$ is continuous, we can extend $f^1$ and $f^2$ to all of $B$ by continuity. Since $f$ is Lipschitz, it then follows that the extensions of $f^1$ and $f^2$ to $B$ are still Lipschitz as well.
	
	We claim that the graphs of $f^1$ and $f^2$ are both stationary in $B$. Indeed, we have that $V_{f^1}:= |\graph(f^1)|$ is an integral $n$-varifold that is stationary on $\pi^{-1}(\Omega)$ (as indeed $f^1$ is a $C^1$ solution to the minimal surface system over $\Omega$). However, as $\H^{n-2}(\S)=0$, we know that $\pi^{-1}(B\setminus\Omega)$ is a $\H^{n-1}$-null set. Moreover, as we have the following growth bound
	$$\|V_{f^1}\|(B_\rho(x))\leq C\rho^n \qquad \text{for all }x\in\pi^{-1}(B\setminus\Omega)$$
	(which follows from $f^1$ being Lipschitz, or the monotonicity formula for stationary varifolds), a standard cut-off argument implies that $V_{f^1}$ is stationary across $\pi^{-1}(B\setminus \Omega)$ also, and so $V_{f^1}$ is stationary in $\pi^{-1}(B)$. The same argument obviously also applies to the graph of $f^2$. This therefore completes the proof (indeed, the graphs must then be disjoint as a coincidence point has density $\geq 2$, cf.~Remark \ref{remark:coincidence-points} below).
\end{proof}

\begin{remark}\label{remark:coincidence-points}
In fact, we can show that if $V$ is a Lipschitz $2$-valued stationary graph and $x\in B_1(0)$ is a coincidence point of the graph, then $\Theta_V(x)\geq 2$. To see this, let $\BC$ be a tangent cone at $x$; this remains a $2$-valued Lipschitz graph, with the origin being a coincidence point. If $\Theta_\BC(y)\geq 2$ for any $y$, we would get $\Theta_V(x) = \Theta_{\BC}(0)\geq 2$ by upper semi-continuity of the density and so would be done. Thus, we may assume that $\Theta_{\BC}(y)<2$ for all $y$. But then again by Remark \ref{remark:small-singular-set} following Definition \ref{defn:V-beta}, this would imply that $\BC$ is embedded (with multiplicity $1$) away from a closed set of dimension $\leq n-3$. Furthermore, we may assume that the coincidence set of $\BC$ is just the origin, else we could take a tangent cone at any other coincidence point $p\neq 0$ to get a new cone which is a 2-valued Lipschitz graph, (after we quotient out the spine) of one lower dimension, and with the origin as a coincidence point; thus one could argue by induction on the dimension of the cone to show $\Theta_{\BC}(p)\geq 2$ (the low dimensional cases, $n\leq 2$, can be dealt with by hand, so assume $n\geq 3$). But then the argument in Lemma \ref{lemma:Lip-2} would imply that $\BC$ is the sum of two stationary single-valued graphs with the origin common to both of them, and thus as each has density $\geq 1$ at the origin, we see $\Theta_V(x) = \Theta_{\BC}(0)\geq 1+1 = 2$.
\end{remark}

\textbf{Remark:} Recalling Remark \ref{remark:special-case-stable}, we mention one other natural setting where Theorem \ref{thm:main} applies. Write $\S_2$ for the set of integral $n$-currents $T$ in $B^{n+1}_2(0)$ with $\del T \res B^{n+1}_{2}(0) = 0$ such that 
$|T|$ (i.e.\ the associated $n$-varifold) is stationary in $B^{n+1}_{2}(0)$ and has stable regular part. In particular, $T$ cannot have any triple junction singularities. Then, for suitably small $\beta = \beta(n,k)\in (0,1)$, we have that $|T| \in \mathcal{V}_\beta$. Indeed, this is a consequence of the regularity theory for codimension one stationary integral varifolds with stable regular part and no classical singularities, established by the third author \cite{Wic14}. So Theorem \ref{thm:main} can be applied to 
$|T|$. In fact, this special case of Theorem \ref{thm:main} has already been established in previous work by the second and third authors, cf.~\cite[Theorem D]{MW24}. Moreover, in this case one gets stronger conclusions, namely that $T$ is actually a $C^{1,1/2}$ $2$-valued graph, which is smoothly immersed away from the set of flat singular points which is countably $(n-2)$-rectifiable (see \cite{SW16, KW21}).

\part{\centering Blow-Ups of Stationary Integral Varifolds off a Plane}\label{part:coarse-blow-ups}

\section{Constructing Coarse Blow-Ups}\label{sec:constructing-coarse-blow-ups}

We now turn to constructing coarse blow-ups of stationary integral varifolds near multiplicity $Q$ planes for any $Q\in \{1,2,\dotsc\}$. For the main results of the present work we will only be interested in the case $Q=2$, however in this part we will look at the general $Q$ case.

\subsection{Lipschitz Graphical Approximation}

We begin by recording a version of Almgren's (weak) multi-valued Lipschitz approximation, \cite[Corollary 3.11]{Alm00}. To reach the version as stated below, we work only with stationary integral $n$-varifolds and apply the estimates from Theorem \ref{thm:allard-sup-estimate} and \eqref{E:rpi} to some of the original statements.

\begin{theorem}[{\cite[Corollary 3.11]{Alm00}}]\label{thm:Lipschitz-approx}
	Fix $\sigma\in (0,1)$, $L\in (0,1)$, and $Q\in \{1,2,\dotsc\}$. Then, there exists $\eps_0 = \eps_0(n,k,Q,\sigma,L)\in (0,1)$ such that the following holds. Let $V$ be a stationary integral $n$-varifold in $B^{n+k}_2(0)$ such that:
	\begin{itemize}
		\item $(\w_n2^n)^{-1}\|V\|(B^{n+k}_2(0))<Q+1/2$;
		\item $Q-1/2\leq \w_n^{-1}\|V\|(C_1(0))<Q+1/2$;
		\item $\hat{E}_V<\eps_0$.
	\end{itemize}
	Then, there exists a Lebesgue measurable set $\Sigma\subset B^n_\sigma(0)$ and a Lipschitz $Q$-valued function $u:B^n_\sigma(0)\to \A_Q(\R^k)$ with $\Lip(u)\leq L$ such that:
	\begin{enumerate}
		\item [(i)] $\H^n(\Sigma) + \|V\|(\R^k\times\Sigma) \leq C\hat{E}_V^2$ for some $C = C(n,k,Q,\sigma,L)\in (0,\infty)$;
		\item [(ii)] For each $x\in B^n_\sigma(0)\setminus \Sigma$ and each $y\in \spt\|V\|\cap \pi^{-1}_{P_0}(x)$, we have $\Theta_V(y)\in \{1,2,\dotsc\}$, $\sum_{y\in \spt\|V\|\cap \pi_{P_0}^{-1}(x)}\Theta_V(y) = Q$, and
		$$u(x) = \sum_{y\in\spt\|V\|\cap \pi^{-1}_{P_0}(x)}\Theta_V(y)\llbracket y\rrbracket;$$
		\item [(iii)] $V\res (\R^k\times (B^n_\sigma(0)\setminus\Sigma)) = \graph(u)\res (\R^k\times (B^n_\sigma(0)\setminus\Sigma))$;
		\item [(iv)] $\|u\|_{L^\infty}\leq C\hat{E}_V$, for some $C = C(n,k,\sigma)\in (0,\infty)$.
	\end{enumerate}
\end{theorem}

\textbf{Remarks:}
\begin{itemize}
	\item In the context of Theorem \ref{thm:Lipschitz-approx}, we write $\Omega:=B^n_\sigma(0)\setminus\Sigma$ and refer to $\Omega$ as the ``good set'' and $\Sigma$ as the ``bad set'' of the approximation.
	\item Instead of prescribing a \emph{constant} upper bound of $L$ on the Lipschitz constant in Theorem \ref{thm:Lipschitz-approx}, one can instead prescribe an upper bound on the Lipschitz constant of $\hat{E}_V^\beta$ for some choice of power $\beta\in (0,1)$. In turn, this changes the upper bound in Theorem \ref{thm:Lipschitz-approx}(i) to $C\hat{E}_V^{2-4\beta(n+1)Q}$, i.e.~we decrease the power on the measure of the bad set, allowing the bad set to be larger. This modification is in \cite[Corollary 3.12]{Alm00}. We could work with this formulation instead, as then the Lipschitz constant naturally $\to 0$ as $\eps_0\to 0$, however we will stick with the statement above and deal with any dependency on $L$ separately when it arises.
	\item If we choose $L\leq 1/2$ in Theorem \ref{thm:Lipschitz-approx}, then there exists $c = c(n,k)>0$ such that for $\H^n$-a.e. $x\in \graph(u)$,
	\begin{equation}\label{E:Lipschitz-approx-R1}
		|\pi_{T^\perp_x\graph(u)}(y)|\geq c|y| \qquad \text{for all }y\in \R^k\times\{0\}^n.
	\end{equation}
	\item Throughout, we will be using the notation $u^{\kappa,\alpha}$ to denote the $\alpha^{\text{th}}$-value of the $\kappa$-coordinate function of $u$, for $\kappa\in\{1,\dotsc,k\}$ and $\alpha\in \{1,\dotsc,Q\}$. However, we will make this choice in a fixed manner as follows: for each such $\kappa$, consider the function $u^\kappa\equiv u\cdot e_\kappa$, which is $\A_Q(\R)$-valued. Being $\R$-valued, we then can order the values, i.e.~we write $u^{\kappa,1}\leq u^{\kappa,2}\leq\cdots\leq u^{\kappa,Q}$ with $u\cdot e_\kappa = \sum_{\alpha=1}^Q\llbracket u^{\kappa,\alpha}\rrbracket$. This then has the property that for each $\alpha$ and $\kappa$, $u^{\kappa,\alpha}$ is a Lipschitz single-valued function. Indeed to see this, as $u$ is Lipschitz, for any $x,y$ we have for some choice of permutation $\sigma_{xy}\in S_Q$,
	$$\sum_{\alpha=1}^Q|u^\alpha(x)-u^{\sigma_{xy}(\alpha)}(y)|^2 = \G(u(x),u(y))^2 \leq \Lip(u)^2|x-y|^2,$$
	and thus for each $\kappa$, as $u^{\kappa,1}\leq \cdots \leq u^{\kappa,Q}$,
	\begin{align*}
	\sum^Q_{\alpha=1}|u^{\kappa,\alpha}(x)-u^{\kappa,\alpha}(y)|^2 = \G(u^\kappa(x),u^\kappa(y))^2 & \leq \sum_{\alpha=1}^Q|u^{\kappa,\alpha}(x)-u^{\kappa,\sigma_{xy}(\alpha)}(y)|^2\\
	& \leq \Lip(u)^2|x-y|^2,
	\end{align*}
	and so we see that each function $u^{\kappa,\alpha}$ as defined above is Lipschitz, with $\Lip(u^{\kappa,\alpha})\leq \Lip(u)$ (here, we have used that for an $\A_Q(\R)$-valued function, the $\G$-distance is attained by the permutation which orders the values). In particular, $u^{\kappa,\alpha}$ is differentiable almost everywhere and is a $W^{1,2}$ function for each $\kappa,\alpha$.  Notice however that (for $k>1$) it is \emph{not} necessarily the case that $\sum^Q_{\alpha=1}\llbracket u^\alpha(x)\rrbracket = \sum^Q_{\alpha=1}\llbracket\sum^k_{\kappa=1}e_\kappa u^{\kappa,\alpha}(x)\rrbracket$, as each of the $k$ permutations chosen for each coordinate function might differ. However, we will have for each $x$, $\sum_{\alpha=1}^Q \llbracket u^\alpha(x)\rrbracket = \sum^Q_{\alpha=1}\llbracket\sum_{\kappa=1}^k e_\kappa u^{\kappa,\sigma_{\kappa,x}(\alpha)}(x)\rrbracket$, for some permutation $\sigma_{\kappa,x}$ depending on both $\kappa$ and $x$. As there are only finitely many permutations which $\sigma_{\kappa,x}$ can be, we therefore also have that the derivative of $u$ as a multi-valued function is determined by, almost everywhere, the derivatives of the Lipschitz single-valued functions $u^{\kappa,\alpha}$. This is a key point when it comes to constructing blow-ups, as we will be taking the limits of rescalings of the single-valued functions $u^{\kappa,\alpha}$ via classical compactness results for single-valued $W^{1,2}$ functions, rather than for multi-valued functions. Notice that expressions which are invariant under permutations of $\alpha$, such as $\sum_{\alpha=1}^Q u^{\kappa,\alpha}$, are always well-defined and independent of the choice of ordering. It will always be the case that our estimates and identities involve such invariant expressions in the $Q$-values, and as such we can always decompose them in the above manner, as a sum of the Lipschitz functions $u^{\kappa,\alpha}$ (or their derivatives) constructed above. As such, the reader should keep in mind throughout the paper the difference between $u^\alpha$ and $u^{\kappa,\alpha}$ as described above. Often we will abuse notation and not explicitly write the permutations $\sigma_{\kappa,x}$ as above when the end result will not depend on this.
\end{itemize}
We will also need one further remark, for which we introduce some notation which we will use throughout the rest of the paper.
\begin{itemize}
	\item Whenever we are given a function $\zeta:B^n_1(0)\to \R^k$, we will denote by $\widetilde{\zeta}$ a choice of extension of $\zeta$ to $\R^k\times B^n_1(0)$ which is compactly supported in $\R^k\times B^n_1(0)$, and which satisfies
	$$\widetilde{\zeta}(x^1,\dotsc,x^{n+k}) = \zeta(x^{k+1},\dotsc,x^{n+k})$$
	when $|(x^1,\dotsc,x^k)|<1$, or at least in a neighbourhood of $\spt\|V\|$ (for the relevant $V$ for the context).
	\item The final remark for Theorem \ref{thm:Lipschitz-approx} is then the following comparison of derivatives on $\spt\|V\|$ or $\graph(u)$ with the derivatives on $P_0$ (the latter being denoted by $D$). Indeed, linear algebra gives for any functions $f_1,f_2\in C^1_c(B^n_1(0))$ and $\|V\|$-a.e. $x\in\spt\|V\|$,
	\begin{equation}\label{E:Lipschitz-approx-R2}
		\left|\nabla^V \widetilde{f}_1(x)\cdot\nabla^V\widetilde{f}_2(x) - D\widetilde{f}_1(x)\cdot D\widetilde{f}_2(x)\right| \leq \|\pi_{T_xV}-\pi_{P_0}\|^2\sup|Df_1|\sup|Df_2|.
	\end{equation}
\end{itemize}

We now prove some key estimates regarding Almgren's Lipschitz approximation, namely a priori energy bounds, as well as estimates on the Jacobian and various $L^2$ height quantities.

\begin{corollary}\label{cor:Lipschitz-approx}
	Fix $\sigma\in (0,1)$, $L\in (0,1)$, and let $V$ be a stationary integral $n$-varifold in $B_2(0)$ for which the conclusions of Theorem \ref{thm:Lipschitz-approx} hold for these choices of $\sigma,L$. Then the following estimates hold:
	\begin{equation}\label{E:Lipschitz-approx-1}
		\int_{B^n_\sigma(0)}\sum^k_{\kappa=1}\sum^Q_{\alpha=1}|Du^{\kappa,\alpha}|^2\leq C\hat{E}_V^2.
	\end{equation}
	\begin{equation}\label{E:Lipschitz-approx-2}
		\int_{B^n_\sigma(0)}\sum^Q_{\alpha=1}|J^\alpha(x)-1|\, \ext x \leq C\hat{E}_V^2.
	\end{equation}
	\begin{equation}\label{E:Lipschitz-approx-3}
		\left|\int_{C_\sigma(0)}\dist^2(x,P_0)\, \ext\|V\|(x) - \int_{B^n_\sigma(0)}|u|^2\right| \leq C\hat{E}_V^4.
	\end{equation}
	Here, $J$ is the Jacobian of $u$. Furthermore, for any $z\in \spt\|V\|\cap C_{3\sigma/4}(0)$ with $\Theta_V(z)\geq Q$, we have:
	\begin{equation}\label{E:Lipschitz-approx-4}
		\int_{B_{\sigma/4}(z)}\frac{\dist^2(x,z^{\perp_{P_0}}+P_0)}{|x-z|^{n+7/4}}\, \ext\|V\|(x) \leq C\hat{E}_V^2.
	\end{equation}
	\begin{equation}\label{E:Lipschitz-approx-5}
		\int_{B_{\sigma/4}(z)}\frac{\dist^2(x,P_0)}{|x-z|^{n-1/4}}\, \ext\|V\|(x) \leq C\hat{E}_V^2.
	\end{equation}
	Moreover, let $\tau\in (0,1/16)$. Then, if $L$ is an $(n-1)$-dimensional subspace of $\{0\}^k\times \R^n$ and $K\subset L$ is a closed subset such that
	$$K\cap B_{1/2}(0)\subset B_\tau(\{z\in \spt\|V\|:\Theta_V(z)\geq Q\}),$$
	then provided we allow $\eps_0$ in Theorem \ref{thm:Lipschitz-approx} to additionally depend on $\tau$, we have
	\begin{equation}\label{E:Lipschitz-approx-6}
	\int_{B_\tau(K)\cap C_{1/2}(0)}\|\pi_{T_xV}-\pi_{P_0}\|^2\, \ext\|V\|(x) \leq C_*\tau^{1/2}\hat{E}_V^2.
	\end{equation}
	In all the above, $C = C(n,k,Q,\sigma,L)\in (0,\infty)$, and $C_* = C_*(n,k,Q,L)\in (0,\infty)$.
\end{corollary}

\begin{proof}
	We first prove \eqref{E:Lipschitz-approx-1}. Since $|Du^{\kappa,\alpha}|\leq L$ for every $\kappa$ and $\alpha$, the integrands in \eqref{E:Lipschitz-approx-1} are all bounded; therefore, Theorem \ref{thm:Lipschitz-approx}(i) gives
	$$\int_{B^n_\sigma(0)}\sum^k_{\kappa=1}\sum^Q_{\alpha=1}|Du^{\kappa,\alpha}|^2 = \int_{\Omega}\sum^k_{\kappa=1}\sum^Q_{\alpha=1}|Du^{\kappa,\alpha}(x)|^2\, \ext x + E_1$$
	where $E_1$ satisfies $|E_1|\leq C\hat{E}_V^2$. Then, just focusing on the other term, using \eqref{E:J-bounds} and the area formula \eqref{E:area-formula}, we have
	\begin{align*}
		\int_\Omega \sum^k_{\kappa=1}\sum^Q_{\alpha=1}|Du^{\kappa,\alpha}|^2 & \leq \int_\Omega \sum^k_{\kappa=1}\sum^Q_{\alpha=1}|Du^{\kappa,\alpha}|^2 J^{\sigma^{-1}_{\kappa,x}(\alpha)}\\
		& = \sum^k_{\kappa=1}\sum^Q_{\alpha=1}\int_{(\R^k\times\Omega)\cap \graph(u^\alpha)}|D\widetilde{u}^{\kappa,\alpha}(x)|^2\, \ext\H^n(x).
	\end{align*}
	Now, since at each point $x\in \spt\|V\|\cap (\R^k\times \Omega)$ we have that $x^\kappa = \widetilde{u}^{\kappa,\alpha}(x)$ for some $\alpha\in \{1,\dotsc,Q\}$, we see from \eqref{E:Lipschitz-approx-R2} that for $\|V\|$-a.e. $x\in\spt\|V\|\cap (\R^k\times\Omega)\cap \{x^\kappa=\widetilde{u}^{\kappa,\alpha}(x)\}$, we have
	\begin{equation}\label{E:Lipschitz-approx-estimates-1}
		\left||\nabla^Vx^\kappa|^2 - |D\widetilde{u}^{\kappa,\alpha}(x)|^2\right| \leq \|\pi_{T_xV}-\pi_{P_0}\|^2\cdot L^2.
	\end{equation}
	The analogous estimate holds for every $\alpha\in \{1,\dotsc,Q\}$. So, using Theorem \ref{thm:Lipschitz-approx}(iii), we deduce that for each $\kappa\in \{1,\dotsc,k\}$,
	$$\int_{\R^k\times\Omega}|\nabla^V x^\kappa|^2\, \ext\|V\|(x) = \sum^Q_{\alpha=1}\int_{(\R^k\times\Omega)\cap\graph(u^\alpha)}|D\widetilde{u}^{\kappa,\alpha}(x)|^2\, \ext\H^n(x) + E_2$$
	where, using \eqref{E:Lipschitz-approx-estimates-1} followed by \eqref{E:rpi}, we can estimate $E_2$ as
	$$|E_2|\leq\int_{\R^k\times\Omega}\|\pi_{T_xV}-\pi_{P_0}\|^2 L^2\, \ext\|V\|(x) \leq C\hat{E}_V^2$$
	for $C = C(n,k,Q,\sigma,L)\in (0,\infty)$. Combining all the above, again with \eqref{E:rpi}, establishes \eqref{E:Lipschitz-approx-1}.
	
	To see \eqref{E:Lipschitz-approx-2}, notice that for $x\in \Omega$, \eqref{E:J-bounds} tells us that
	$$\sum^Q_{\alpha=1}|J^\alpha(x)-1|\leq C\sum^Q_{\alpha=1}|Du^\alpha(x)|^2 = C\sum_{\kappa=1}^k\sum_{\alpha=1}^Q|Du^{\kappa,\alpha}(x)|^2.$$
	Integrating this over $\Omega$ and using \eqref{E:Lipschitz-approx-1}, we have
	$$\int_\Omega\sum^Q_{\alpha=1}|J^\alpha(x)-1|\, \ext x \leq C\hat{E}_V^2.$$
	Then using Theorem \ref{thm:Lipschitz-approx}(i), we can add back in the integral over $\Sigma$ to the left-hand side at the cost of adjusting the constant on the right-hand side. This proves \eqref{E:Lipschitz-approx-2}.
	
	For \eqref{E:Lipschitz-approx-3}, write
	$$\int_{C_\sigma(0)}\dist^2(x,P_0)\, \ext\|V\|(x) = \int_{\R^k\times\Sigma}\dist^2(x,P_0)\, \ext\|V\|(x) + \int_{\R^k\times\Omega}\dist^2(x,P_0)\, \ext\|V\|(x).$$
	For the first integral on the right-hand side, we crudely estimate it by using Theorem \ref{thm:Lipschitz-approx}(i) and Theorem \ref{thm:allard-sup-estimate} to get
	$$\left|\int_{\R^k\times\Sigma}\dist^2(x,P_0)\, \ext\|V\|(x)\right| \leq \|V\|(\R^k\times\Sigma)\sup_{x\in C_\sigma(0)}\dist^2(x,P_0) \leq C\hat{E}_V^4.$$
	For the second integral on the right-hand side, we have
	\begin{align*}
		\int_{\R^k\times\Omega}\dist^2(x,P_0)\, \ext\|V\|(x) & = \int_\Omega\sum^Q_{\alpha=1}|u^\alpha(x)|^2 J^\alpha(x)\, \ext x\\
		& = \int_{B_\sigma^n(0)}\sum^Q_{\alpha=1}|u^\alpha(x)|^2 J^\alpha(x)\, \ext x + E_3
	\end{align*}
	where we can estimate $E_3$ using Theorem \ref{thm:Lipschitz-approx}(i) \& (iv) and the Lipschitz bound $\Lip(u)\leq L$:
	$$|E_3|\leq \int_\Sigma\sum^Q_{\alpha=1}|u^\alpha(x)|^2 J^\alpha(x)\, \ext x \leq C\H^n(\Sigma)\sup_{B^n_\sigma(0)}\max_\alpha |u^\alpha(x)|^2 \leq C\hat{E}_V^4.$$
	Finally, we can write
	$$\int_{B_\sigma^n(0)}\sum^Q_{\alpha=1}|u^\alpha(x)|^2 J^\alpha(x)\, \ext x = \int_{B^n_\sigma(0)}\sum^Q_{\alpha=1}|u^\alpha(x)|^2\, \ext x + E_4$$
	where
	$$|E_4| \leq \int_{B^n_\sigma(0)}\sum^Q_{\alpha=1}\sum^k_{\kappa=1}|u^{\kappa,\alpha}(x)|^2|J^\alpha(x)-1|\, \ext x \leq \sup_{x\in B^n_\sigma(0)}\max_{\alpha^\prime}|u^{\alpha^\prime}(x)|^2\int_{B_\sigma^n(0)}\sum^Q_{\alpha=1}|J^\alpha-1|$$
	so that using \eqref{E:Lipschitz-approx-2} and Theorem \ref{thm:Lipschitz-approx}(iv) we get
	$$|E_4|\leq C\hat{E}_V^4.$$
	Combining the above inequalities then gives \eqref{E:Lipschitz-approx-3}.
	
	Now we prove \eqref{E:Lipschitz-approx-4}. We start from the following consequence of the monotonicity formula:
	$$\frac{1}{\w_n}\int_{B_{\sigma/4}(z)}\frac{|(x-z)^{\perp_{T_xV}}|^2}{|x-z|^{n+2}}\, \ext\|V\|(x) = \frac{\|V\|(B_{\sigma/4}(z))}{\w_n(\sigma/4)^n} - \Theta_V(z).$$
	Now write $\bar{z}:= \pi_{P_0}(z)$. Using the inclusion $B_{\sigma/4}(z)\subset C_{\sigma/4}(\bar{z})$, Theorem \ref{thm:Lipschitz-approx}, the area formula \eqref{E:area-formula}, and \eqref{E:Lipschitz-approx-2}, we can estimate the right-hand side of this monotonicity formula by
	$$\frac{\|V\|(B_{\sigma/4}(z))}{\w_n(\sigma/4)^n}-Q\leq \frac{1}{\w_n(\sigma/4)^n}\int_{B^n_{\sigma/4}(\bar{z})}\sum^Q_{\alpha=1}|J^\alpha(x)-1|\, \ext x \leq C\hat{E}_V^2$$
	for some $C = C(n,k,Q,\sigma,L)\in (0,\infty)$. Thus, we have that
	\begin{equation}\label{E:Lipschitz-approx-estimates-2}
		\int_{B_{\sigma/4}(z)}\frac{|(x-z)^{\perp_{T_xV}}|^2}{|x-z|^{n+2}}\, \ext\|V\|(x) \leq C\hat{E}_V^2.
	\end{equation}
	Next, we use the first variation formula with the vector field
	$$\Phi(x) := \zeta(x-z)^2\cdot\frac{|(x-z)^{\perp_{P_0}}|^2}{|x-z|^{n+7/4}}(x-z),$$
	where $\zeta\in C^\infty_c(\R^{n+k})$ with $\zeta\equiv 1$ on $B_{\sigma/8}(0)$ and $\zeta\equiv 0$ on $\R^{n+k}\setminus B_{\sigma/4}(0)$, and $|D\zeta|\leq C\sigma^{-1}$. A cut-off function argument is needed to show that this vector field is valid in the first variation formula (due to it being singular at $x=z$) but this is straightforward. Now, a standard computation (cf.~\cite[Page 616]{Sim93}) gives that
	\begin{align*}
		\int_{B_{\sigma/4}(z)}\frac{|(x-z)^{\perp_{P_0}}|^2}{|x-z|^{n+7/4}}\, \ext\|V\|(x) & \leq C\int_{B_{\sigma/4}(z)}\zeta(x-z)^2\frac{|(x-z)^{\perp_{T_xV}}|^2}{|x-z|^{n+7/4}}\, \ext\|V\|(x)\\
		& \hspace{4em} + C\int_{B_{\sigma/4}(z)}|\nabla^V\zeta(x-z)|^2\frac{|(x-z)^{\perp_{P_0}}|^2}{|x-z|^{n-1/4}}\, \ext\|V\|(x).
	\end{align*}
	The first term on the right-hand side is then bounded by $C\hat{E}_V^2$ using \eqref{E:Lipschitz-approx-estimates-2}. For the second term on the right-hand side, since $\spt(D\zeta(\cdot-z))\subset B_{\sigma/4}(z)\setminus B_{\sigma/8}(z)$, the integrand is bounded pointwise by $C|(x-z)^{\perp_{P_0}}|^2$ for some $C = C(n,k,\sigma,L)\in (0,\infty)$. Using Allard's supremum estimate (Theorem \ref{thm:allard-sup-estimate}) then bounds this pointwise by $C\hat{E}_V^2$. Thus, the right-hand side of this expression is controlled by $C\hat{E}_V^2$, and so \eqref{E:Lipschitz-approx-4} follows by noting that $|(x-z)^{\perp_{P_0}}| = \dist(x,z^{\perp_{P_0}}+P_0)$.
	
	Next we prove \eqref{E:Lipschitz-approx-5}. For this, we will need the estimate
	\begin{equation}\label{E:Lipschitz-approx-estimates-3}
		\int_{B_{\sigma/4}(z)}\frac{1}{|x-z|^{n-1/4}}\, \ext\|V\|(x)\leq C = C(n,k,Q,\sigma,L)<\infty.
	\end{equation}
	Indeed, to see this note that for $m\in \{2,3,\dotsc\}$ we have
	$$\sigma^n 2^{-(m+1)n} \leq \|V\|(B_\rho(z)) \leq C\sigma^n 2^{-mn} \qquad \text{for all }\rho\in (\sigma 2^{-m-1},\sigma 2^{-m})$$
	by the monotonicity formula; here, $C = C(n,k,Q,1-\sigma)$. Moreover, for $x\in B_{\sigma 2^{-m}}(z)\setminus B_{\sigma 2^{-m-1}}(z)$ we have that $|x-z|^{-n+1/4} \leq \sigma^{-n+1/4}2^{(m+1)(n-1/4)}$. Therefore we have
	$$\int_{B_{\sigma 2^{-m}}(z)\setminus B_{\sigma 2^{-m-1}}(z)}\frac{1}{|x-z|^{n-1/4}}\, \ext\|V\|(x) \leq C\sigma^{1/4}2^{n-1/4 - m/4}$$
	where $C = C(n,k,Q,1-\sigma)\in (0,\infty)$. Summing this over $m\in \{2,3,\dotsc\}$ establishes \eqref{E:Lipschitz-approx-estimates-3}.
	
	Now, the triangle inequality gives
	$$\dist^2(x,P_0) \leq 2\dist^2(x,z^{\perp_{P_0}}+P_0) + 2|z^{\perp_{P_0}}|^2.$$
	If we divide this by $|x-z|^{n-1/4}$ and integrate with respect to $\|V\|$ over $B_{\sigma/4}(z)$, we have
	\begin{align*}
		\int_{B_{\sigma/4}(z)}\frac{\dist^2(x,P_0)}{|x-z|^{n-1/4}}\, \ext\|V\|(x) & \leq 2\int_{B_{\sigma/4}(z)}\frac{\dist^2(x,z^{\perp_{P_0}}+P_0)}{|x-z|^{n-1/4}}\, \ext\|V\|(x)\\
		& \hspace{4em} + 2|z^{\perp_{P_0}}|^2\int_{B_{\sigma/4}(z)}\frac{1}{|x-z|^{n-1/4}}\, \ext\|V\|(x).
	\end{align*}
	The first term on the right-hand side is controlled by $C\hat{E}_V^2$ using \eqref{E:Lipschitz-approx-4}. Using Theorem \ref{thm:allard-sup-estimate} and \eqref{E:Lipschitz-approx-estimates-3}, the second term on the right-hand side is controlled by $C\hat{E}_V^2$ as well. Thus, this proves \eqref{E:Lipschitz-approx-5}.
	
	Finally, we prove \eqref{E:Lipschitz-approx-6}. Notice that \eqref{E:Lipschitz-approx-4} gives (taking $\sigma = 3/4$ therein) that for any $z\in \spt\|V\|\cap C_{9/16}(0)$ with $\Theta_V(z)\geq Q$, we have
	$$\int_{B_{3\tau}(z)}\dist^2(x,z^{\perp_{P_0}}+P_0)\, \ext\|V\|(x) \leq C\tau^{n+7/4}\hat{E}_V^2.$$
	In view of the hypothesis $K\cap B_{1/2}\subset B_\tau(\{z\in\spt\|V\|:\Theta_V(z)\geq Q\})$
	the preceding estimate implies that for each $y\in K\cap B_{1/2}$, there exists $z\in \spt\|V\|\cap C_{9/16}(0)$ with $\Theta_V(z)\geq Q$ such that
	$$\int_{B_{2\tau}(y)}\dist^2(x,z^{\perp_{P_0}}+P_0)\, \ext\|V\|(x) \leq C\tau^{n+7/4}\hat{E}_V^2.$$
	This in turn implies, also using \eqref{E:rpi} applied to $(\eta_{z^{\perp_{P_0}},1})_\#V$ (more precisely, applying the argument leading to \eqref{E:rpi} using a suitable choice of test function therein), that for each $y\in K\cap B_{1/2}$,
	$$\int_{B_{3\tau/2}(y)}\|\pi_{T_xV}-\pi_{P_0}\|^2\, \ext\|V\|(x) \leq C\tau^{n-1/4}\hat{E}_V^2.$$
	Since $K\subset L$ is closed and $L$ is a subspace of dimension $n-1$, we may cover the set $B_\tau(K)\cap C_{1/2}(0)$ by balls $\{B_{3\tau/2}(y_j)\}_{j=1}^N$, with $y_j\in K\cap B_{1/2}$ and $N\leq C\tau^{1-n}$, $C = C(n)$, it follows that
	$$\int_{B_\tau(L)\cap C_{1/2}(0)}\|\pi_{T_xV}-\pi_{P_0}\|^2\, \ext\|V\|(x) \leq \tau^{3/4}\hat{E}_{V}^2$$
	which gives \eqref{E:Lipschitz-approx-6}. This completes the proof of the corollary.
\end{proof}

\subsection{Sequences of Stationary Integral Varifolds Converging to a Plane}\label{sec:coarse-blow-up-construction}

Let $(V_j)_j$ be a sequence of stationary integral $n$-varifolds in $B_2(0)$ such that $V_j\weakly Q|P_0|$ as varifolds in $C_2(0)$. Write $\hat{E}_j:= \hat{E}_{V_j}$, and fix $\sigma, L\in (0,1)$. Then, for sufficiently large $j$ (depending on $\sigma$ and $L$) we can invoke Theorem \ref{thm:Lipschitz-approx}; let $\Omega_j\subset B^n_\sigma(0)$ denote the measurable subset over which $V_j$ is graphical in the sense that for sufficiently large $j$ we have the conclusions of Theorem \ref{thm:Lipschitz-approx} hold with $V_j$ in place of $V$ and $\Sigma_j:= B^n_\sigma(0)\setminus\Omega_j$ in place of $\Sigma$, i.e.
$$V_j\res (\R^k\times\Omega_j) = \mathbf{v}(u_j)\res (\R^k\times\Omega_j)$$
where $u_j:B^n_\sigma(0)\to \A_Q(\R^k)$ is a Lipschitz function with $\Lip(u_j)\leq L$ and such that
\begin{equation}\label{E:bad-set-j}
	\|V_j\|(\R^k\times\Sigma_j) + \H^n(\Sigma_j)\leq C\hat{E}_j^2
\end{equation}
where $C = C(n,k,Q,\sigma,L)\in (0,\infty)$. We also write $v_j:= \hat{E}_j^{-1}u_j$ and $v^{\kappa,\alpha}_j := \hat{E}_V^{-1}u^{\kappa,\alpha}_j$.

To start, note that \eqref{E:Lipschitz-approx-3} gives
$$\int_{B^n_\sigma(0)}\sum^Q_{\alpha=1}\sum^{k}_{\kappa=1}|u^{\kappa,\alpha}_j|^2 = \int_{C_\sigma(0)}\dist^2(x,P_0)\, \ext\|V_j\|(x) + o(\hat{E}_j^2)$$
where $o(\bullet_j)$ denotes a term obeying $o(\bullet_j)/\bullet_j\to 0$ as $j\to\infty$. Since the right-hand side of the above is clearly at most $C\hat{E}_j^2$, this shows that $(v_j)_j$ is bounded uniformly in $L^2(B^n_\sigma(0))$. Moreover, \eqref{E:Lipschitz-approx-1} gives
$$\int_{B^n_\sigma(0)}\sum^Q_{\alpha=1}\sum^k_{\kappa=1}|Du_j^{\kappa,\alpha}|^2\leq C\hat{E}_V^2$$
showing that $(Dv_j^{\kappa,\alpha})_{j=1}^\infty$ is bounded in $L^2(B^n_\sigma(0))$ for each $\kappa,\alpha$. Therefore, for every $\sigma\in (0,1)$, we have that $(v^{\kappa,\alpha}_j)_j$ is bounded uniformly in $W^{1,2}(B^n_\sigma(0))$ (with a bound depending only on $n,k,Q,\sigma,L$; recall here that indeed $Dv^{\kappa,\alpha}_j$ \emph{is} the weak derivative of $v^{\kappa,\alpha}_j$).

Also, let $\Gamma_{j,\sigma}$ denote the compact set
$$\Gamma_{j,\sigma}:= \pi_{P_0}\left(\{\Theta_{V_j}\geq Q\}\right)\cap \overline{B^n_\sigma(0)}.$$
Then, using Rellich's compactness theorem (for \emph{single-valued} functions) and the sequential compactness of the Hausdorff metric on the space of compact subsets of a compact set, we deduce that there exist functions $v^{\kappa,\alpha}\in W^{1,2}(B^n_\sigma(0))$ for each $\kappa$ and $\alpha$, a closed subset $\Gamma_{v,\sigma}\subset \overline{B^n_\sigma(0)}$, and a subsequence $\{j^\prime\}\subset\{j\}$ along which we have:
\begin{itemize}
	\item $v^{\kappa,\alpha}_{j^\prime}\to v^{\kappa,\alpha}$ strongly in $L^2(B_\sigma^n(0))$;
	\item $Dv^{\kappa,\alpha}_{j^\prime}\weakly Dv^{\kappa,\alpha}$ weakly in $L^2(B^n_\sigma(0);\R^n)$;
	\item $\dist_\H(\Gamma_{j,\sigma},\Gamma_{v,\sigma})\to 0$.
\end{itemize}
Now, take a sequence $\sigma_p\uparrow 1$ and perform this construction with $\sigma = \sigma_p$ for $p\in \{1,2,\dotsc\}$. Then, by choosing an appropriate diagonal subsequence $p = p_j$, we get that there is a (relatively) closed set $\Gamma_v\subset B^n_1(0)$ and a subsequence of $\{j\}$ (which we again call $\{j^\prime\}$, abusing notation) along which we have:
\begin{itemize}
	\item $v^{\kappa,\alpha}_{j^\prime}\to v^{\kappa,\alpha}$ strongly in $L^2_{\text{loc}}(B_1^n(0))$;
	\item $Dv^{\kappa,\alpha}_{j^\prime}\weakly Dv^{\kappa,\alpha}$ weakly in $L^2_{\text{loc}}(B^n_1(0);\R^n)$;
	\item $\Gamma_{j,\sigma_{p_j}}\to\Gamma_{v}$ locally in Hausdorff distance in $B^n_1(0)$.
\end{itemize}
We also know by the above construction that $v^{\kappa,\alpha}\in L^2(B_1^n(0))\cap W^{1,2}_{\text{loc}}(B^n_1(0))$ (indeed, Fatou's lemma gives that $v^{\kappa,\alpha}\in L^2(B^n_1(0))$, as \eqref{E:Lipschitz-approx-3} gives for all $j$ sufficiently large that $\|v_j\|_{L^2(B_{\sigma_{p_j}}(0))}\leq 2$).

We can also take the weak $L^2$ limit of $v_j$ as a multi-valued function, giving rise to a function $v:B^n_1(0)\to \A_Q(\R^k)$ in $L^2(B^n_1(0);\A_Q(\R^k))$. Notice again, however, that we do not necessarily have $v^\alpha = \sum_{\kappa=1}^k\llbracket v^{\kappa,\alpha}\rrbracket$. What we do have is for $\H^n$-a.e $x\in B^n_1(0)$, there are permutations $\sigma_{\kappa,x}$ of $\{1,\dotsc,Q\}$ for which $v^\alpha(x) = \sum_{\alpha=1}^k e_\kappa v^{\kappa,\sigma_{\kappa,x}}(x)$.

\begin{defn}\label{defn:coarse-blow-up}
	We call $v$ (and the corresponding $v^{\kappa,\alpha}$) constructed in the above manner a \emph{coarse blow-up} of the sequence $(V_j)_j$ off $P_0$. We call $\Gamma_v$ the set of \emph{Hardt--Simon points} of $v$.
\end{defn}	
We write $\Bfrak_Q$ for the collection of all $v\in L^2(B_1^n(0);\A_Q(\R^k))$ that arise, in the manner described above, as coarse blow-ups of a sequence $(V_j)_j$ of stationary integral $n$-varifolds in $B_2(0)$ with $V_j\weakly Q|P_0|$ in $B_2(0)$.

\begin{remark}\label{remark:after-coarse-construction}
We note that a coarse blow-up does not depend on the choice of Lipschitz constant $L$ in the sequence of Lipschitz approximations given by Theorem \ref{thm:Lipschitz-approx}. Indeed, suppose $(u_j)_j$ and $(u_j^\prime)_j$ are two sequences of Lipschitz approximations given by Theorem \ref{thm:Lipschitz-approx} applied to a given sequence of stationary integral varifolds $(V_j)_j$ for different choices of Lipschitz constants $L,L^\prime$ therein. Then, it is still the case that the measure of the set where $u_j\neq u_j^\prime$ is going to zero (as the graphs of both $u_j$ and $u_j^\prime$ agree with $V_j$ on their respective good sets). This means that $v_j-v_j^\prime$ converges to zero in measure, and hence converges to zero almost everywhere along a subsequence. Thus, the coarse blow-ups will agree (up to redefining on sets of measure zero).
\end{remark}

We finish this subsection by recording the basic but important fact that the principle obstruction to having $L^2(B^n_1(0))$ convergence of the rescaled Lipschitz approximations to the blow-up is when the height excess of the varifolds concentrates near the boundary of the ball $B^n_1(0)$.

\begin{prop}[Non-Concentration of Height Excess implies $L^2$ Convergence]\label{prop:full-L2}
	Let $v$ be the coarse blow-up of a sequence of stationary integral varifolds $(V_j)_j$ off $P_0$. Suppose that for any $\eps>0$, we can find $\sigma\in (0,1)$ and $J$ such that
	$$\int_{C_1(0)\setminus C_\sigma(0)}\dist^2(x,P_0)\, \ext\|V_j\|(x)<\eps\hat{E}_{V_j}^2 \qquad \text{for all }j\geq J.$$
	Then we have $v^{\kappa,\alpha}_j\to v^{\kappa,\alpha}$ in $L^2(B_1^n(0))$ for every $\kappa$ and $\alpha$. In particular, we have
	$$\hat{E}_{V_j}^{-2}\int_{C_1(0)}\dist^2(x,\graph(\hat{E}_{V_j}v))\, \ext\|V_j\|(x)\to 0.$$
\end{prop}

\textbf{Note:} Despite $v$ only being a function in $L^2$, we can still make sense of $\dist(\cdot,\graph(\hat{E}_{V_j}v))$ as a function $\H^n$-a.e.

\begin{proof}
	Observe that
	$$\int_{B^n_1(0)}|v_j^{\kappa,\alpha}-v^{\kappa,\alpha}|^2 \leq \int_{B_1^n(0)\setminus B_\sigma^n(0)}|v_j^{\kappa,\alpha}|^2 + \int_{B^n_\sigma(0)}|v_j^{\kappa,\alpha}-v^{\kappa,\alpha}|^2 + \int_{B^n_1(0)\setminus B^n_\sigma(0)}|v^{\kappa,\alpha}|^2.$$
	For the terms on the right-hand side:
	\begin{itemize}
		\item the third term can be made small by choosing $\sigma$ sufficiently close to $1$ (as $v^{\kappa,\alpha}\in L^2(B^n_1(0))$;
	\end{itemize}
	and, fixing such a $\sigma$ to ensure the third term is $<\eps$ (say),
	\begin{itemize}
		\item the second term can be made small by using choosing $j$ sufficiently large, as $v^{\kappa,\alpha}_j$ converges strongly in $L^2$ to $v^{\kappa,\alpha}$ on $B_\sigma^n(0)$;
		\item the first term can be made small (for this choice of $\sigma$) by again choosing $j$ sufficiently large and using \eqref{E:Lipschitz-approx-3} and the assumption in the present proposition (as combined these give that, for $j$ large, the $L^2$ norm of $u_j$ is arbitrarily close in proportion to $\hat{E}_{V_j}$).
	\end{itemize}
	This therefore establishes the $L^2$ convergence in $B_1^n(0)$. With the $L^2$ convergence verified, the second claim of the proposition follows immediately from an argument essentially identical to that which established \eqref{E:Lipschitz-approx-3} using the local $L^2$ convergence along the blow-up sequence.
\end{proof}

\section{General Properties of Coarse Blow-Ups}

We now want to study the coarse blow-up class $\Bfrak_Q$ defined above and establish properties that it satisfies. Later in Part \ref{part:blow-up-reg} we will specialise to the case $Q=2$ and prove further, more refined, properties in that setting, assuming the topological structural condition holds in $\beta$-coarse gaps.

For general $Q$, the main properties and estimates are summarised in the following theorem.

\begin{theorem}\label{thm:blow-up-properties}
	The blow-up class $\Bfrak_Q$ satisfies the following properties:
	\begin{enumerate}
		\item [$(\Bfrak1)$] If $v\in \Bfrak_Q$, then $v\in L^2(B^n_1(0);\A_Q(\R^k))$, and $v^{\kappa,\alpha}\in W^{1,2}_{\textnormal{loc}}(B^n_1(0))$ for each $\kappa\in \{1,\dotsc,k\}$ and $\alpha\in\{1,\dotsc,Q\}$.
		\item [$(\Bfrak2)$] If $v\in\Bfrak_Q$, then $v^\kappa_{a}$ is smooth and harmonic on $B^n_1(0)$ for each $\kappa\in\{1,\dotsc,k\}$.
		\item [$(\Bfrak3)$] \textnormal{(Compactness and Closure Properties.)} If $v\in\Bfrak_Q$, then:
		\begin{enumerate}
			\item [$(\Bfrak3\text{I})$] For any $z\in B^n_1(0)$ and $\rho\in (0,\frac{3}{8}(1-|z|)]$, if $v\not\equiv 0$ on $B_\rho^n(z)$, then
			$$v_{z,\rho}(x):= \frac{v(z+\rho x)}{\|v(z+\rho(\cdot))\|_{L^2(B_1^n(0))}} \in \Bfrak_Q;$$
			\item [$(\Bfrak3\text{II})$] If $\gamma:\R^n\to\R^n$ is an orthogonal rotation of $\R^n$, then $v\circ\gamma\in\Bfrak_Q$;
			\item [$(\Bfrak3\text{III})$] If $\lambda\in\R^k$ and $L:\R^n\to \R^k$ is a linear function with $v\not\equiv \lambda+L$ on $B^n_1(0)$, then
			$$\frac{v-\lambda-L}{\|v-\lambda-L\|_{L^2(B^n_1(0))}}\in \Bfrak_Q;$$
			\item [$(\Bfrak3\text{IV})$] If $(v_p)_{p=1}^\infty\subset\Bfrak_Q$, then there exists a subsequence $\{p^\prime\}\subset\{p\}$ and $v\in \Bfrak_Q$ such that
			\begin{itemize}
				\item $v_{p^\prime}\to v$ strongly in $L^2_{\textnormal{loc}}(B^n_1(0);\A_Q(\R^k))$;
				\item $Dv_{p^\prime}^{\kappa,\alpha}\weakly Dv^{\kappa,\alpha}$ weakly in $L^2_{\textnormal{loc}}(B^n_1(0))$ for each $\kappa\in \{1,\dotsc,k\}$ and $\alpha\in \{1,\dotsc,Q\}$;
				\item $\Gamma_{v_{p^\prime}}\to \Gamma_v$ locally in Hausdorff distance in $B_1^n(0)$.
			\end{itemize}
		\end{enumerate}
		\item [$(\Bfrak4)$] \textnormal{(Squash Inequality.)} If $v\in \Bfrak_Q$ and $\kappa\in \{1,\dotsc,k\}$, we have
		$$\int_{B^n_1(0)}\sum^Q_{\alpha=1}|Dv^{\kappa,\alpha}|^2\zeta \leq -\int_{B^n_1(0)}\sum^Q_{\alpha=1}v^{\kappa,\alpha}Dv^{\kappa,\alpha}\cdot D\zeta$$
		and
        $$\int_{B^n_1(0)}\sum^Q_{\alpha=1}|Dv_f^{\kappa,\alpha}|^2\zeta \leq -\int_{B^n_1(0)}\sum^Q_{\alpha=1}v_f^{\kappa,\alpha}Dv_f^{\kappa,\alpha}\cdot D\zeta$$
        for every $\zeta\in C^1_c(B^n_1(0);\R)$ with $\zeta\geq 0$.
		\item [$(\Bfrak5)$] \textnormal{(Hardt--Simon Inequality.)} For each $v\in \Bfrak_Q$, $z\in \Gamma_v$, and $\rho\in (0,\frac{3}{8}(1-|z|)]$, we have
		$$\sum^k_{\kappa=1}\sum^Q_{\alpha=1}\int_{B^n_{\rho/2}(z)}R_z(x)^{2-n}\left|\frac{\del}{\del R_z}\left(\frac{v^{\kappa,\alpha}(x)-v^\kappa_a(z)}{R_z(x)}\right)\right|^2\, \ext x \leq C\rho^{-n-2}\int_{B^n_\rho(z)}|v|^2$$
		where $R_z(x) := |x-z|$ and $C = C(n,k,Q)\in (0,\infty)$.
		\item [$(\Bfrak6)$] \textnormal{(Gradient Non-Concentration.)} For each $v\in \Bfrak_Q$, $\delta>0$, and $\kappa\in \{1,\dotsc,k\}$, we have
		$$\sum^Q_{\alpha=1}\int_{B^n_{1/2}(0)\cap \{|v^{\kappa,\alpha}_f|<\delta\}}|Dv^{\kappa,\alpha}_f|^2 \leq C\delta$$
		where $C = C(n,k,Q)\in (0,\infty)$.
	\end{enumerate}
\end{theorem}

\begin{remark}\label{remark:after-blow-up-properties}
We note some simple consequences of these properties:
\begin{itemize}
	\item By taking in $(\Bfrak4)$ $\zeta$ such that $\zeta\equiv 1$ on $B^n_{1/2}(0)$ with $|D\zeta|\leq 4$ and then using Cauchy--Schwarz, we deduce the following energy estimates for $v\in \Bfrak_Q$:
$$\int_{B^n_{1/2}(0)}|Dv|^2\leq C\int_{B_1(0)}|v|^2, \qquad \text{and} \qquad \int_{B_{1/2}^n(0)}|Dv_f|^2 \leq C\int_{B_1(0)}|v_f|^2.$$
	\item Let $v\in \Bfrak_Q$, $z\in B_1^n(0)$, and $\rho_j\downarrow 0$. Consider the sequence $v_j:= v_{z,\rho_j}$ as defined in $(\Bfrak3\text{I})$; we know by $(\Bfrak3\text{I})$ that $v_j\in \Bfrak_Q$ for all $j$ sufficiently large, provided $v$ does not vanish on a neighbourhood of $z$. Thus, by $(\Bfrak3\text{IV})$, there is a subsequence $\{j^\prime\}\subset\{j\}$ and $\psi\in \Bfrak_Q$ such that $v_{z,\rho_{j^\prime}}\to \psi$ as described in $(\Bfrak3\text{IV})$. We call any such $\psi$ that arises in this way a \emph{tangent} (\emph{coarse}) \emph{blow-up of $v$ at $z$}. A priori, we know very little about tangent coarse blow-ups.
	\item If we knew sufficient regularity on $v$ at a point $z_0\in\Gamma_v$ (e.g.~that $v$ was Lipschitz) then $(\Bfrak5)$ would tell us that $v^\alpha(z_0) = v_a(z_0)$ for all $\alpha\in \{1,\dotsc,Q\}$, i.e.~that every point $z_0\in\Gamma_v$ is a coincidence point of $v$.
	\item In $(\Bfrak5)$, by using $(\Bfrak3\text{III})$ we can replace the right-hand side of the inequality with instead $C\rho^{-n-2}\int_{B^n_\rho(z)}|v-\ell_{v_a,z}|^2$, where $\ell_{v_a,z}(x):= v_a(z) + (x-z)\cdot Dv_a(z)$. Moreover, in $(\Bfrak6)$, again by applying $(\Bfrak3)$ (and letting $\sigma\uparrow 1$) we can replace the right-hand side of the inequality with $C\delta\left(\int_{B_1(0)}|v|^2\right)^{1/2}$.
\end{itemize}
\end{remark}

Let us now prove each of the properties in Theorem \ref{thm:blow-up-properties}.

\begin{proof}[Proof of $(\Bfrak1)$]
	This is immediate from the construction of $\Bfrak_Q$.
\end{proof}

\begin{proof}[Proof of $(\Bfrak2)$]
	Suppose that $v\in \Bfrak_Q$ arises as the coarse blow-up of $(V_j)_j$. Fix $\zeta\in C^1_c(B^n_1(0))$. From the first variation formula \eqref{E:stationarity-2}, we have
	\begin{equation}\label{E:FB2-1}
		\int_{C_1(0)}\nabla^{V_j}x^\kappa\cdot\nabla^{V_j}\widetilde{\zeta}\, \ext\|V_j\|(x) = 0
	\end{equation}
	for any $\kappa\in \{1,\dotsc,k\}$. After writing the part of the integral takes place over the graphical region $\R^k\times\Omega_j$ as an integral over the domain $B^n_1(0)$ by the area formula \eqref{E:area-formula}, and then adding
	$$\sum^Q_{\alpha=1}\int_{\Sigma_j}\nabla^{V_j}\widetilde{u}^{\kappa,\alpha}_j(x)\cdot\nabla^{V_j}\widetilde{\zeta}(x)\, J^\alpha_j(x)\, \ext x$$
	to both sides, \eqref{E:FB2-1} gives (note that here by $\nabla^{V_j}$ we mean the gradient operator on $V_j$ at the point $(u_j^\alpha(x),x)\in \spt\|V_j\|$)
	\begin{align*}
		\sum^Q_{\alpha=1}\int_{B^n_1(0)}&\nabla^{V_j}\widetilde{u}^{\kappa,\alpha}_j(x)\cdot\nabla^{V_j}\widetilde{\zeta}(x)\, J^\alpha_j(x)\, \ext x\\
		& = -\int_{\R^k\times \Sigma_j}\nabla^{V_j}x^\kappa\cdot\nabla^{V_j}\widetilde{\zeta}(x)\, \ext\|V_j\|(x) + \sum^Q_{\alpha=1}\int_{\Sigma_j}\nabla^{V_j}\widetilde{u}^{\kappa,\alpha}_j(x)\cdot\nabla^{V_j}\widetilde{\zeta}(x)\, J^\alpha_j(x)\, \ext x.
	\end{align*}
	Note that here we are writing $\Sigma_j:= B^n_1(0)\setminus\Omega_j$, which is a slight abuse of notation from earlier (as $\zeta$ has compact support, all the integrals are technically only happening over regions over $B_\sigma^n(0)$ for some $\sigma\in (0,1)$, so this does not impact anything for $j$ sufficiently large). This equation can be rewritten as
	$$\sum^Q_{\alpha=1}\int_{B^n_1(0)}Du^{\kappa,\alpha}_j\cdot D\zeta = -\int_{\R^k\times\Sigma_j}\nabla^{V_j}x^\kappa\cdot\nabla^{V_j}\widetilde{\zeta}\, \ext\|V_j\| + \sum^Q_{\alpha=1}\int_{\Sigma_j}\nabla^{V_j}\widetilde{u}^{\kappa,\alpha}_j\cdot\nabla^{V_j}\widetilde{\zeta}\, J^{\alpha}_j\, \ext x + F_j$$
	where
	$$|F_j|\leq \sum^Q_{\alpha=1}\left|\int_{B^n_1(0)}\nabla^{V_j}\widetilde{u}_j^{\kappa,\alpha}\cdot\nabla^{V_j}\widetilde{\zeta}\, J^{\alpha}_j - Du^{\kappa,\alpha}_j\cdot D\zeta\, \ext x\right|.$$
	Combining the above with \eqref{E:Lipschitz-approx-R2}, \eqref{E:J-bounds}, \eqref{E:rpi}, \eqref{E:Lipschitz-approx-3}, and \eqref{E:bad-set-j}, we get
	$$\left|\sum^Q_{\alpha=1}\int_{B^n_1(0)}Du_j^{\kappa,\alpha}\cdot D\zeta\right| \leq C\sup D\zeta\cdot\hat{E}_j^2$$
	where $C = C(n,k,Q,\sigma,L)\in (0,\infty)$; here, $\sigma\in (0,1)$ is such that $\spt(\zeta)\subset B^n_\sigma(0)$, which exists as $\zeta$ has compact support in $B^n_1(0)$. Therefore, dividing this by $\hat{E}_j$ and taking $j\to\infty$, using that $Dv_j^{\kappa,\alpha}\to Dv^{\kappa,\alpha}$ weakly in $L^2_{\text{loc}}(B^n_1(0))$, we get
	$$\int_{B^n_1(0)}Dv^\kappa_a\cdot D\zeta = 0.$$
	As $\zeta\in C^1_c(B^n_1(0))$ was arbitrary, this is saying $v^\kappa_a$ is weakly harmonic. Thus, by Weyl's lemma we get that $v^\kappa_a$ is smooth and harmonic.
\end{proof}

\begin{proof}[Proof of $(\Bfrak3)$]
	Let $v\in \Bfrak_Q$. Firstly, let $z\in B^n_1(0)$ and $\rho\in (0,\frac{3}{8}(1-|z|)]$ and suppose that $v\not\equiv 0$ in $B^n_\rho(z)$. Then, if $(V_j)_j$ is a sequence whose coarse blow-up off $P_0$ is $v$, it is straightforward to check that the coarse blow-up of $\tilde{V}_j:= (\eta_{(0,z),\rho})_\#V_j$ off $P_0$ is $v_{z,\rho}$ (here, we are using that the original coarse blow-up sequence $v_j$ converges strongly to $v$ in $L^2(B^n_\rho(z))$); note that here one also needs to check that $\tilde{V}_j$ converges to $Q|P_0|$ as varifolds. but this is straightforward. This establishes $(\Bfrak3\text{I})$. Similarly, if $\gamma:\R^n\to \R^n$ is a rotation of $\R^n$, then the coarse blow-up of $\tilde{\gamma}_\#V_j$ is $v\circ \gamma$, which establishes $(\Bfrak3\text{II})$.
	
	Next, we prove the compactness property $(\Bfrak3\text{IV})$ (this will be needed to prove $(\Bfrak3\text{III})$). Suppose that for each $p$, we have that $v_p$ is the coarse blow-up of the sequence $(V_j^{(p)})_{j=1}^\infty$. Then, for each $p$, we can choose $j_p$ such that the sequence $(j_p)_{p=1}^\infty$ satisfies:
	\begin{itemize}
		\item $(j_p)_{p=1}^\infty$ is strictly increasing;
		\item $\hat{E}_{V^{(p)}_{j_p}}\to 0$;
		\item The conclusions of Theorem \ref{thm:Lipschitz-approx} hold with $1-\frac{1}{p}$ in place of $\sigma$ and $V^{(p)}_{j_p}$ in place of $V$;
		\item For every $\kappa\in \{1,\dotsc,k\}$ and $\alpha\in \{1,\dotsc,Q\}$, we have
		$$\left\|(u^{(p)}_{j_p})^{\kappa,\alpha}/\hat{E}_{V_{j_p}^{(p)}} - v_p^{\kappa,\alpha}\right\|_{L^2(B^n_{1-\frac{1}{p}}(0))}<\frac{1}{p}$$
		where $u^{(p)}_{j_p}$ is the $Q$-valued Lipschitz approximation for $V^{(p)}_{j_p}$ given by Theorem \ref{thm:Lipschitz-approx} in the bullet point above.
	\end{itemize}
	That this is possible is clear from the construction of coarse blow-ups. Now, let $v$ be a blow-up of $(V^p_{j_p})_{p=1}^\infty$ (up to passing to a subsequence). Then, using the fourth bullet point above and the fact that $v_p^{\kappa,\alpha}$, $v^{\kappa,\alpha}$ are (uniformly) bounded in $W^{1,2}(B^n_\sigma(0))$ for each fixed $\sigma\in (0,1)$, we may pass to a subsequence $\{p^\prime\}\subset\{p\}$ for which $v_{p^\prime}\to v$ in the sense described in $(\Bfrak3\text{IV})$, completing its proof.
	
	Now we can prove $(\Bfrak3\text{III})$. Firstly, given $y\in\R^k$ for which $v-y\not\equiv 0$, we have that for all $\sigma\in (0,1)$ sufficiently close to $1$, the map $x\mapsto v(\sigma x)-y$ is not identically zero. So, if $\tau_j:\R^{n+k}\to \R^{n+k}$ is the translation map $x\mapsto x-\hat{E}_j y$, then we have that the coarse blow-up of the sequence $\tilde{V}_j:= (\tau_j\circ \eta_{0,\sigma_j})_\#V_j$ is
	$$\frac{v(\sigma(\cdot))-y}{\|v(\sigma(\cdot))-y\|_{L^2(B_1^n(0))}}.$$
	Thus this function belongs to $\Bfrak_Q$. Then, using $(\Bfrak3\text{IV})$, we can let $\sigma\uparrow 1$ and take a convergent subsequence, the limit of which belongs to $\Bfrak_Q$ and by construction is necessarily $\|v-y\|_{L^2(B_1^n(0))}(v-y)$; hence, we get $\|v-y\|_{L^2(B^n_1(0))}^{-1}(v-y)\in \Bfrak_Q$ as claimed.
	
	Secondly, given a linear function $L:\R^n\to \R^k$ for which $v-L\not\equiv 0$, let $R_j$ be a rotation of $\R^{n+k}$ that minimises $\|R-\id_{\R^{n+k}}\|$ over rotations object to the property that $R_j(\graph(\hat{E}_j L)) = P_0$. Then, as for $\sigma\in (0,1)$ sufficiently close to $1$ we have $v(\sigma(\cdot))-L\not\equiv 0$, for such $\sigma$ we have that the coarse blow-up of the sequence $\hat{V}_j:= (R_j\circ\eta_{0,\sigma})_\#V_j$ is
	$$\frac{v(\sigma(\cdot))-L}{\|v(\sigma(\cdot))-L\|_{L^2(B_1^n(0))}}\in \Bfrak_Q.$$
	Again, we can then let $\sigma\uparrow 1$ and use $(\Bfrak3\text{IV})$ to deduce that $\|v-L\|_{L^2(B^n_1(0))}^{-1}(v-L)\in\Bfrak_Q$.
	
	To deduce the general result from the above two cases, suppose that we are given $\lambda\in \R^k$ and a linear function $L:\R^n\to \R^k$ with $v-\lambda-L\not\equiv 0$. We can then apply the first fact with $y=\lambda$, and then the second with $\|v-\lambda\|_{L^2(B^n_1(0))}^{-1}L$ in place of $L$ and $\|v-y\|_{L^2(B^n_1(0))}^{-1}(v-y)$ in place of $v$ to deduce $(\Bfrak3\text{III})$. This completes the proof of all the properties in $(\Bfrak3)$.
\end{proof}

\begin{proof}[Proof of $(\Bfrak4)$]
	First notice that the second inequality follows from the first by writing $v_f = v-v_a$ and using the fact that by $(\mathfrak{B}2)$ $v_a$ is harmonic, and for harmonic functions the inequality is true and in fact is an equality. Thus, it suffices to prove the first claimed inequality. For this, fix $\kappa\in\{1,\dotsc,k\}$ and let $\zeta\in C^1_c(B^n_1(0))$ be such that $\zeta\geq 0$. By taking $\widetilde{\zeta}x^\kappa$ in the first variation formula \eqref{E:stationarity-2} we get
    $$\int_{C_1(0)}|\nabla^{V_j}x^\kappa|^2 \widetilde{\zeta}\, \ext\|V_j\|(x) = -\int_{C_1(0)}x^\kappa\nabla^{V_j}x^\kappa\cdot\nabla^{V_j}\widetilde{\zeta}\, \ext\|V_j\|(x).$$
	Omitting the non-graphical part in the integral on the left-hand side of the above (which we can do as $\zeta\geq 0$), using the bound $J_j^\alpha\geq 1$, and using the area formula \eqref{E:area-formula} to write the graphical parts over the domain $B^n_1(0)$ where possible, we get
	\begin{align*}
		\int_{\Omega_j}\sum^Q_{\alpha=1}&|\nabla^{V_j}\widetilde{u}^{\kappa,\alpha}_j(x)|^2\zeta\, \ext x\\
		& \leq -\int_{\Omega_j}\sum^Q_{\alpha=1}\widetilde{u}^{\kappa,\alpha}_j\nabla^{V_j}\widetilde{u}^{\kappa,\alpha}_j\cdot\nabla^{V_j}\widetilde{\zeta}\, J^\alpha_j\, \ext x - \int_{\R^k\times\Sigma_j}x^\kappa\nabla^{V_j}x^\kappa\cdot\nabla^{V_j}\widetilde{\zeta}\, \ext\|V_j\|.
	\end{align*}
	Using Theorem \ref{thm:allard-sup-estimate}, Cauchy--Schwarz, and \eqref{E:bad-set-j}, we have the bound (for all $j$ sufficiently large depending on $\sigma$, where $\sigma\in (0,1)$ is such that $\spt(\zeta)\subset B^n_\sigma(0)$)
	$$\left|\int_{\R^k\times\Sigma_j}x^\kappa\nabla^{V_j}x^\kappa\cdot\nabla^{V_j}\widetilde{\zeta}\, \ext\|V_j\|\right|\leq C\hat{E}_j^3$$
	where $C = C(n,k,Q,\zeta)$. Notice that we can also replace the $\nabla^{V_j}$ derivatives in the integral on the left-hand side by the derivatives in the $P_0$ direction via \eqref{E:Lipschitz-approx-R2}, up to incurring an error term which by \eqref{E:rpi} is controlled by (using the Lipschitz bound on $u_j$) $C_*L^2\hat{E}_j^2$, where $C_* = C_*(n,k,Q,\zeta)$ is independent of $L$. When we do the same thing in the first term on the right-hand side, using additionally Theorem \ref{thm:allard-sup-estimate}, we get an error which is at most $C_*L\hat{E}_j^3$. If we then divide the resulting inequality by $\hat{E}_j^2$ and take $j\to\infty$, we may use the lower semi-continuity of energy under weak convergence in $W^{1,2}$, as well as $\one_{\Omega_j}\to \one_{B^n_1(0)}$ a.e., to obtain
	$$\int_{B_1(0)}\sum^Q_{\alpha=1}|Dv^{\kappa,\alpha}|^2\zeta \leq C_*L^2 -\int_{B_1(0)}\sum^Q_{\alpha=1}v^{\kappa,\alpha}Dv^{\kappa,\alpha}\cdot D\zeta.$$
	This is almost $(\Bfrak4)$, except there is an additional $C_*L^2$ on the right-hand side. However, we can remove this by recalling that (cf.~Remark \ref{remark:after-coarse-construction} after the construction of blow-ups) that the final blow-up $v$ is independent of $L$. Thus, we can simply let $L\downarrow 0$ in the above to arrive at $(\Bfrak4)$. 
	
	(We also note that one can deal with the domain of integration on the left-hand side depending on $j$ by working over regions of the form $\bigcap_{\ell\geq J}\Omega_\ell$; these are independent of $j$, are contained in $\Omega_j$ for all $j\geq J$, and for suitable subsequences of $\{j\}$ will have measure converging to the measure of $B^n_1(0)$. One may then pass to the limit over a domain of this form for fixed $J$ and then take $J\to\infty$.)
\end{proof}

\begin{proof}[Proof of $(\Bfrak5)$]
	We fix $j\geq 1$ and for ease of notation will drop the subscript $j$ in the notation below, i.e.~we will just write $V = V_j$, $\Omega = \Omega_j$, and so forth.
	
	From the monotonicity formula, for any $z\in \spt\|V\|\cap C_{1/4}(0)$ we have
	$$\frac{1}{\w_n}\int_{B_{5/8}(z)}\frac{|(x-z)^{\perp}|^2}{|x-z|^{n+2}}\, \ext\|V\|(x) = \frac{\|V\|(B_{5/8}(z))}{\w_n(5/8)^n}-\Theta_V(z)$$
	where $\perp\,\equiv\, \perp_{T_xV}$ in integrals with respect to $\|V\|$. Suppose that $z$ is such that, for some constant $\tilde{C}\geq 0$, we have
	$$\Theta_V(z) \geq Q-\tilde{C}\hat{E}_V^2$$
	(in particular, if $\Theta_V(z)\geq Q$ this is satisfied with $\tilde{C}=0$; throughout the present discussion we include this slightly more general case, and allow constants to depend also on $\tilde{C}$). Arguing similar to the proof of \eqref{E:Lipschitz-approx-4} gives (writing $\bar{z}:= \pi_{P_0}(z)$):
	\begin{align*}
		\frac{\|V\|(B_{5/8}(z))}{\w_n(5/8)^n}-Q\leq \frac{1}{\w_n(5/8)^n}\int_{B^n_{5/8}(\bar{z})}\sum^Q_{\alpha=1}(J^\alpha-1) \leq C\int_{B^n_{5/8}(\bar{z})}\sum^Q_{\alpha=1}|Du^\alpha|^2 + C\hat{E}_V^2 \leq C\hat{E}_V^2
	\end{align*}
	for some $C = C(n,k,Q)$. Thus, we have that
	\begin{equation}\label{E:B5-1}
		\int_{B_{5/8}(z)}\frac{|(x-z)^{\perp}|^2}{|x-z|^{n+2}}\, \ext\|V\|(x) \leq C\hat{E}_V^2
	\end{equation}
	where now $C$ depends additionally on $\tilde{C}$. We will bound the left-hand side from below in two different ways. Clearly we have
	$$\int_{B_{5/8}(z)}\frac{|(x-z)^{\perp}|^2}{|x-z|^{n+2}}\, \ext\|V\|(x) \geq \int_{B_{5/8}(z)}|(x-z)^\perp|^2\, \ext\|V\|(x).$$
	Using the notation $\hat{x}:= x-\bar{x}$ where $\bar{x} = \pi_{P_0}(x)$ combined with the elementary inequality
	$$|(\hat{x}-\hat{z})^\perp|^2 \leq 2|(x-z)^\perp|^2 + 2n\sum^{n+k}_{i=k+1}|e_i^{\perp}|^2$$
	which holds for $|x-z|<5/8$, we can further bound this below by
	$$\frac{1}{2}\int_{B_{5/8}(z)}|(\hat{x}-\hat{z})^\perp|^2\, \ext\|V\| - n\int_{B_{5/8}(z)}\sum^{n+k}_{i=k+1}|e_i^\perp|^2\, \ext\|V\|.$$
	Using that $\text{trace}(\pi_{T_xV}) = n$, we have $\sum_{i=k+1}^{n+k}|e_i^{\perp}|^2 = \sum^k_{\kappa=1}|\nabla^V x^\kappa|^2 = \frac{1}{2}\|\pi_{P_0}-\pi_{T_xV}\|^2$. Thus, using \eqref{E:rpi} and the triangle inequality, the above expression is further bounded below by
	\begin{align*}
		&\frac{1}{2}\int_{B_{5/8}(z)}|(\hat{x}-\hat{z})^\perp|^2\, \ext\|V\|(x) - C\hat{E}_V^2\\
		& \geq \frac{1}{4}\int_{B_{5/8}(z)}|\hat{z}^\perp|^2\, \ext\|V\| - \frac{1}{4}\int_{B_{5/8}(z)}|\hat{x}^\perp|^2\, \ext\|V\|(x) - C\hat{E}_V^2\\
		& \geq \frac{1}{4}\int_{B_{5/8}(z)}|\hat{z}^\perp|^2\, \ext\|V\| - C\hat{E}_V^2.
	\end{align*}
	Using \eqref{E:Lipschitz-approx-R1}, \eqref{E:bad-set-j}, and the lower bound $J^\alpha\geq 1$ (cf.~\ref{E:J-bounds}), we have
	$$\int_{B_{5/8}(z)}|\hat{z}^\perp|^2\, \ext\|V\| \geq c|\hat{z}|^2\int_{B^n_{1/4}(\bar{z})\cap\Omega}\sum^Q_{\alpha=1}J^\alpha \geq c|\hat{z}|^2 - C\hat{E}_V^2.$$
	Combining all the above, we have therefore shown that for such $z$,
	\begin{equation}\label{E:B5-2}
		\dist(z,P_0) \equiv |\hat{z}| \leq C\hat{E}_V.
	\end{equation}
	Now let us start again with the left-hand side of \eqref{E:B5-1} and bound it from below in a different way. Firstly, we have that
	\begin{equation}\label{E:B5-3}
		\int_{B_{5/8}(z)}\frac{|(x-z)^{\perp}|^2}{|x-z|^{n+2}}\, \ext\|V\|(x)\geq \int_{\R^k\times (B^n_{1/2}(\bar{z})\cap \Omega)}\frac{|(x-z)^\perp|^2}{|x-z|^{n+2}}\, \ext\|V\|(x).
	\end{equation}
	Writing $R_{\bar{z}}\equiv R_{\bar{z}}(x):= |x-\bar{z}|$, we next claim that (again, we recall $\perp\, \equiv\, \perp_{T_xV}$ here)
	\begin{equation}\label{E:B5-4}
	-R_{\bar{z}}(\bar{x})^2\frac{\del}{\del R_{\bar{z}}}\left(\frac{u^\alpha(\bar{x})-\hat{z}}{R_{\bar{z}}(\bar{x})}\right)^{\perp} = (x-z)^{\perp}
	\end{equation}
	for $x = (u^\alpha(\bar{x}),\bar{x})$, $\bar{x}\in \Omega$. This is checked just as in \cite[Page 616]{Sim93}: indeed, let $\Phi(\bar{x}):= (u^\alpha(\bar{x}),\bar{x})$, so that for $x\in \spt\|\graph(u^\alpha)\|$ we have $(x-z)^{\perp} = (\Phi(\bar{x})-z)^\perp$. Since $\Phi$ always maps into $\graph(u^\alpha)$, we have
	$$\frac{\del\Phi}{\del R_{\bar{z}}} \in T_x\graph(u^\alpha)$$
	almost everywhere (as $u^\alpha$ is differentiable a.e.). By expanding the derivative and using this fact for $\bar{x}\in\Omega$ (so $T_x\graph(u^\alpha) = T_xV$), we readily get for such $\bar{x}$,
	$$\frac{\del}{\del R_{\bar{z}}}\left(\frac{\Phi(\bar{x})-z}{R_{\bar{z}}(\bar{x})}\right)^{\perp} = -\frac{(x-z)^\perp}{R_{\bar{z}}(\bar{x})^2}.$$
	But on the other hand, the homogeneity of $R_{\bar{z}}(\bar{x})$ easily leads to
	$$\frac{\del}{\del R_{\bar{z}}}\left(\frac{\Phi(\bar{x})-z}{R_{\bar{z}}(\bar{x})}\right)^\perp = \frac{\del}{\del R_{\bar{z}}}\left(\frac{u^\alpha(\bar{x})-\hat{z}}{R_{\bar{z}}(\bar{x})}\right)^{\perp}$$
	which gives the claimed identity \eqref{E:B5-4}. Therefore, using \eqref{E:B5-4} in \eqref{E:B5-3} gives
	\begin{align*}
		\int_{\R^k\times (B^n_{1/2}(\bar{z})\cap \Omega)}&\frac{|(x-z)^\perp|^2}{|x-z|^{n+2}}\, \ext\|V\|(x)\\
		& \geq \int_{B^n_{1/2}(\bar{z})\cap\Omega}\frac{R_{\bar{z}}(x)^4}{|\Phi(x)-z|^{n+2}}\sum^Q_{\alpha=1}\left|\frac{\del}{\del R_{\bar{z}}}\left(\frac{u^\alpha(x)-\hat{z}}{R_{\bar{z}}(x)}\right)^{\perp_{T_{\Phi(x)}V}}\right|^2J^\alpha(x)\, \ext x.
	\end{align*}
	Using \eqref{E:Lipschitz-approx-R1} again, together with $J^\alpha\geq 1$ (recalling \eqref{E:B5-1} and \eqref{E:B5-3})
	\begin{equation}\label{E:B5-5}
		\int_{B^n_{1/2}(\bar{z})\cap \Omega}\frac{R_{\bar{z}}(x)^{n+2}}{|\Phi(x)-z|^{n+2}}\cdot R_{\bar{z}}(x)^{2-n}\sum^Q_{\alpha=1}\left|\frac{\del}{\del R_{\bar{z}}}\left(\frac{u^\alpha(x)-\hat{z}}{R_{\bar{z}}(x)}\right)\right|^2\, \ext x \leq C\hat{E}_V^2.
	\end{equation}
	Now we prove the desired claim of $(\Bfrak5)$ by looking at what the above says along a sequence when we blow-up. Indeed, suppose that for some constant $\tilde{C}$, we have a sequence of points $z_j\in\spt\|V_j\|\cap C_{1}(0)$ with $z_j\to z_0$ as $j\to\infty$ with $\Theta_{V_j}(z_j)\geq Q-\tilde{C}\hat{E}_j^2$ (with $(\Bfrak5)$ corresponding to the case $\tilde{C}=0$ and $z_0\in\Gamma_v$). Thus, fixing $\rho\in (0,\frac{3}{8}(1-|z_0|)]$ and applying the above argument leading to \eqref{E:B5-5} to $(\eta_{z_j,\rho})_\#V_j$ in place of $V$ and then changing variables, we see that
	\begin{align}
	\nonumber \int_{B^n_{\rho/2}(\bar{z}_j)\cap\Omega_j}\frac{R_{\bar{z}_j}(x)^{n+2}}{|\Phi_j(x)-z_j|^{n+2}}\cdot R_{\bar{z}_j}(x)^{2-n}\sum^Q_{\alpha=1}&\left|\frac{\del}{\del R_{\bar{z}_j}}\left(\frac{u^\alpha_j(x)-\hat{z}_j}{R_{\bar{z}_j}(x)}\right)\right|^2\, \ext x\\
	& \leq C\rho^{-n-2}\int_{C_\rho(\bar{z}_j)}\dist^2(x,\hat{z}_j+P_0)\, \ext\|V_j\|(x)\label{E:B5-6}
	\end{align}
	for $C = C(n,k,Q,\tilde{C})\in (0,\infty)$. Now, for all sufficiently large $j$ depending on $\rho$, we have that $\|V_j\|(C_{\rho/16}(\bar{z}_j))\geq C\rho^n$ for some constant $C = C(n,k)\in (0,\infty)$. This means that there exists $y_j\in\spt\|V_j\|\cap C_{\rho/16}(\bar{z}_j)$ such that
	$$|\hat{y}_j|^2 \leq C\rho^{-n}\int_{C_{\rho/16}(\bar{z}_j)}\dist^2(x,P_0)\, \ext\|V_j\|(x).$$
	But now we can apply the estimate \eqref{E:B5-2} with $\tilde{V}:= (\eta_{y_j,\rho/2})_\#V_j$ in place of $V$ and $\tilde{z}:= (\rho/2)^{-1}(z_j-y_j)$ in place of $z$ in order to deduce that
	$$|\hat{z}_j-\hat{y}_j|^2 \leq C\rho^{-n}\int_{C_\rho(\bar{z}_j)}\dist^2(x,\hat{y}_j+P_0)\, \ext\|V_j\|(x)$$
	which gives
	$$|\hat{z}_j|^2 \leq C\rho^{-n}\int_{C_\rho(\bar{z}_j)}\dist^2(x,P_0)\, \ext\|V_j\|(x).$$
	This tells us that, up to passing to a subsequence, the limit $\lambda:= \lim_{j\to\infty}\hat{z}_j/\hat{E}_j\in\R^k$ exists. Now, dividing \eqref{E:B5-6} by $\hat{E}_j^2$ and taking $j\to\infty$, we get (we can replace $\dist(\cdot,\hat{z}_j+P_0)$ by $\dist(\cdot,P_0)$ in the right-hand side, up to increasing the constant, using the above height bound)
	$$\sum^k_{\kappa=1}\sum^Q_{\alpha=1}\int_{B^n_{\rho/2}(z_0)}R_{z_0}(x)^{2-n}\left|\frac{\del}{\del R_{z_0}}\left(\frac{v^{\kappa,\alpha}(x)-\lambda^\kappa}{R_{z_0}(x)}\right)\right|^2\, \ext x \leq C\rho^{-n-2}\int_{B^n_\rho(z_0)}|v|^2$$
	(here we note that $z_0\in P_0$). Here, $C = C(n,k,Q,\tilde{C})$. Now we just need to determine $\lambda$. The triangle inequality gives
	$$\int_{B_{\rho/2}^n(z_0)}R_{z_0}(x)^{2-n}\left|\frac{\del}{\del R_{z_0}}\left(\frac{v_a(x)-\lambda}{R_{z_0}(x)}\right)\right|^2\, \ext x <\infty$$
	which, as $v_a$ is a smooth function by $(\Bfrak2)$, implies that we must have $\lambda = v_a(z_0)$ (and thus $\lambda$ only depends on $z_0$ and not on the sequence $(z_j)_j$ chosen to converge to $z$). Thus, we now have for any $z\in P_0$ which is the limit point of a sequence $(z_j)_j$ with $\Theta_{V_j}(z_j)\geq Q-\tilde{C}\hat{E}_V^2$ and any $\rho\in (0,\frac{3}{8}(1-|z_0|)]$,
	$$\sum^k_{\kappa=1}\sum^Q_{\alpha=1}\int_{B^n_{\rho/2}(z_0)}R_{z_0}(x)^{2-n}\left|\frac{\del}{\del R_{z_0}}\left(\frac{v^{\kappa,\alpha}(x)-v^\kappa_a(z_0)}{R_{z_0}(x)}\right)\right|^2\, \ext x \leq C\rho^{-n-2}\int_{B^n_\rho(z_0)}|v|^2.$$
	where $C = C(n,k,Q,\tilde{C})$. Taking $z_0\in \Gamma_v$ and $\tilde{C}=0$ then gives $(\Bfrak5)$.
\end{proof}

\begin{proof}[Proof of $(\Bfrak6)$] Fix $\eps>0$ and $v\in \Bfrak_Q$. Since $|Dv^{\kappa,\alpha}|\in L^2(B^n_{3/4}(0))$ for every $\kappa\in\{1,\dotsc,k\}$ and $\alpha\in\{1,\dotsc,Q\}$, and moreover as $|Dv_a|\in L^2(B^n_{3/4}(0))$, we can use Egorov's theorem to find a Lebesgue measurable set $U\subset B^n_{3/4}(0)$ (with $\H^n(B^n_{3/4}(0)\setminus U)<\eps$) such that
\begin{equation}\label{E:B6-1}
	v_j^{\kappa,\alpha}\to v^{\kappa,\alpha} \qquad \text{uniformly on $U$}
\end{equation}
for every such $\kappa$ and $\alpha$, and moreover such that
\begin{equation}\label{E:B6-2}
	\int_{B^n_{3/4}(0)\setminus U}|Dv_{a}|^2 + \sum^k_{\kappa=1}\sum^Q_{\alpha=1}\int_{B^n_{3/4}(0)\setminus U}|Dv^{\kappa,\alpha}|^2<\eps.
\end{equation}
The proof is then comprised of three steps.

\textbf{Step 1:} \emph{Applying the first variation formula.} We define a cut-off function as follows. For $\delta\in (0,1)$, let $\eta_\delta:\R\to \R$ be a non-decreasing $C^2$ function satisfying
\begin{itemize}
	\item $\eta_\delta\equiv -2\delta$ on $(-\infty,-2\delta]$ and $\eta_\delta\equiv 2\delta$ on $[2\delta,\infty)$;
	\item $\eta_\delta(t) = t$ for $t\in (-\delta,\delta)$;
	\item $0\leq \eta_\delta^\prime\leq 4$ everywhere.
\end{itemize}
Now fix $\zeta\in C^1_c(B^n_1(0))$. Suppose that $v$ is the coarse blow-up of $(V_j)_j$ off $P_0$, and fix $\kappa\in\{1,\dotsc,k\}$. Using the first variation formula \eqref{E:stationarity-2} with the function $\phi(x):= \eta_\delta(x^\kappa - \hat{E}_j\widetilde{v}^\kappa_a(x))\widetilde{\zeta}(x)^2$ we get
\begin{equation}\label{E:B6-3}
	\int\nabla^{V_j}x^\kappa\cdot\nabla^{V_j}\left(\eta_\delta(x^\kappa-\hat{E}_j\widetilde{v}^\kappa_a(x))\widetilde{\zeta}(x)^2\right)\, \ext\|V_j\|(x) = 0.
\end{equation}
	We now expand the $\nabla^{V_j}$ derivative and analyse the resulting terms. Using the fact that $\sup|\eta_\delta|\leq 2\delta$ together with a simple upper bound for the term that involves $\nabla^{V_j}\widetilde{\zeta}$, from \eqref{E:B6-3} we get:
	\begin{align}
		\nonumber \int_{C_1(0)}\eta_\delta^\prime(x^\kappa-\hat{E}_jv_a^\kappa)\left(|\nabla^{V_j}x^\kappa|^2 - \hat{E}_j\nabla^{V_j}x^\kappa\cdot\nabla^{V_j}\widetilde{v}_a^\kappa\right)&\widetilde{\zeta}^2\, \ext\|V_j\|\\
		&\leq 2\delta\int_{C_1(0)}|\widetilde{\zeta}||\nabla^{V_j}\widetilde{\zeta}||\nabla^{V_j}x^\kappa|\, \ext\|V_j\|. \label{E:B6-4}
	\end{align}
	Now choose $\zeta$ so that $\zeta\equiv 1$ on $B^n_{1/2}(0)$, $\zeta\equiv 0$ on $B^n_1(0)\setminus B^n_{3/4}(0)$, and $|D\zeta|\leq 5$ everywhere. Let us also write $\delta = \bar{\delta}\hat{E}_j$. Thus, using the Cauchy--Schwarz inequality with \eqref{E:rpi}, we may bound the right-hand side of \eqref{E:B6-4} by $C\bar{\delta}\hat{E}_j^2$, so that \eqref{E:B6-4} becomes
	\begin{equation}\label{E:B6-5}
		 \int_{C_1(0)}\eta_\delta^\prime(x^\kappa-\hat{E}_jv_a^\kappa)\left(|\nabla^{V_j}x^\kappa|^2 - \hat{E}_j\nabla^{V_j}x^\kappa\cdot\nabla^{V_j}\widetilde{v}_a^\kappa\right)\widetilde{\zeta}^2\, \ext\|V_j\| \leq C\bar{\delta}\hat{E}_j^2
	\end{equation}
	where $C = C(n,k)\in (0,\infty)$.
	
	Now look at the integral on the left-hand side of \eqref{E:B6-5}. We begin by splitting the region of integration into the complementary regions $\R^k\times\Omega_j$ and $\R^k\times\Sigma_j$ (note that since $\spt(\zeta)\subset B^n_{3/4}(0)$, we are only really interested in the integral over $C_{3/4}(0)$). The integral that takes place over $\R^k\times\Sigma_j$ is equal to
	$$ \int_{\R^k\times\Sigma_j}\eta_\delta^\prime(x^\kappa-\hat{E}_jv_a^\kappa)\left(|\nabla^{V_j}x^\kappa|^2 - \hat{E}_j\nabla^{V_j}x^\kappa\cdot\nabla^{V_j}\widetilde{v}_a^\kappa\right)\widetilde{\zeta}^2\, \ext\|V_j\|(x).$$
	Since $\eta_\delta^\prime\geq 0$, the inequality \eqref{E:B6-5} is still true if we throw away the term here involving $|\nabla^{V_j}x^\kappa|^2$, i.e.~if we replace this by just
	$$-\hat{E}_j \int_{\R^k\times\Sigma_j}\eta_\delta^\prime(x^\kappa-\hat{E}_jv_a^\kappa)\nabla^{V_j}x^\kappa\cdot\nabla^{V_j}\tilde{v}_a^\kappa\widetilde{\zeta}^2\, \ext\|V_j\|(x).$$
	Using the supremum bound for the integrand together with \eqref{E:bad-set-j}, we immediately see that in absolute value this term is bounded by
	$$C\sup_{B^n_{3/4}(0)}|Dv_a|\hat{E}_j^3$$
	where $C = C(n,k,Q,L)\in (0,\infty)$, where $L$ is the prescribed upper bound on the Lipschitz constant for $u_j$ from Theorem \ref{thm:Lipschitz-approx}. (We stress that whilst the Lipschitz constant $L$ does appear in the proof in an important manner, ultimately the final statement in $(\Bfrak6)$ is independent of $L$, since as remarked after the construction of coarse blow-ups, the final coarse blow-up is independent on the choice of $L$.)
	
	Now, returning this all to \eqref{E:B6-5} and using the area formula to write the remaining part of the integral over $\R^k\times\Omega_j$ as instead an integral over $\Omega_j$, we get (and dividing through everything by $\hat{E}_j^2$):
	\begin{align}
		\nonumber \int_{\Omega_j\cap B^n_{3/4}(0)}\sum^Q_{\alpha=1}\eta^\prime_{\bar{\delta}\hat{E}_j}(u_j^{\kappa,\alpha}(x)-\hat{E}_jv_a^\kappa(x))\left(|\nabla^{V_j}v_j^{\kappa,\alpha}|^2 - \nabla^{V_j}v^{\kappa,\alpha}_j\cdot\nabla^{V_j}v^\kappa_a\right)&\zeta^2J_j^\alpha\, \ext x\\
		&\hspace{-2em} \leq C(\bar{\delta} + \sup_{B^n_{3/4}(0)}|Dv_a|\hat{E}_j). \label{E:B6-6}
	\end{align}
	For the integral here we use the pointwise identity
	$$|\nabla^{V_j}v_j^{\kappa,\alpha}|^2 - \nabla^{V_j}v^{\kappa,\alpha}_j\cdot\nabla^{V_j}v^\kappa_a = |\nabla^{V_j}v_j^{\kappa,\alpha}-\nabla^{V_j}v^\kappa_a|^2 + \nabla^{V_j}v_a^\kappa\cdot (\nabla^{V_j}v_j^{\kappa,\alpha}-\nabla^{V_j}v_a^\kappa).$$
	This splits the remaining integral into two term: the \emph{main term}
	\begin{equation}\label{E:B6-7}
		\int_{\Omega_j\cap B^n_{3/4}(0)}\sum^Q_{\alpha=1}\eta^\prime_{\bar{\delta}\hat{E}_j}(u^{\kappa,\alpha}_j(x)-\hat{E}_jv^\kappa_{a}(x))|\nabla^{V_j}v^{\kappa,\alpha}_j - \nabla^{V_j}v^\kappa_a|^2\zeta^2 J_j^\alpha\, \ext x
	\end{equation}
	and the \emph{error term}
	\begin{equation}\label{E:B6-8}
		\int_{\Omega_j\cap B^n_{3/4}(0)}\sum^Q_{\alpha=1}\eta^\prime_{\bar{\delta}\hat{E}_j}(u^{\kappa,\alpha}_j(x)-\hat{E}_jv^\kappa_a(x))\nabla^{V_j}v^\kappa_a\cdot (\nabla^{V_j}v_j^{\kappa,\alpha}-\nabla^{V_j}v^\kappa_a)\zeta^2 J^\alpha_j\, \ext x.
	\end{equation}
	In Step 2 we will estimate the error term \eqref{E:B6-8} and in Step 3 we will analyse the main term \eqref{E:B6-7} to conclude the result.
	
	\textbf{Step 2:} \emph{Controlling the error term \eqref{E:B6-8}.} We would first like to arrange that the domain of integration is $B^n_{3/4}(0)$. Doing this introduces another error term, namely
	$$\int_{\Sigma_j\cap B^n_{3/4}(0)}\sum^Q_{\alpha=1}\eta^\prime_{\bar{\delta}\hat{E}_j}(u^{\kappa,\alpha}_j(x)-\hat{E}_jv^\kappa_a(x))\nabla^{V_j}v^\kappa_a\cdot (\nabla^{V_j}v_j^{\kappa,\alpha}-\nabla^{V_j}v^\kappa_a)\zeta^2 J^\alpha_j\, \ext x.
$$
This latter error term is not difficult to estimate: using the supremum bounds on $\eta^\prime_{\bar{\delta}\hat{E}_j}$, $\zeta$, and $J^\alpha_j$ (cf.~\eqref{E:J-bounds}) as well as the fact that $\Lip(u_j)\leq L \leq 1$, we have that this is controlled by
$$\leq C\sup_{B^n_{3/4}(0)}|Dv_a| (\hat{E}_j^{-1} + \sup_{B^n_{3/4}(0)}|Dv_a|)\cdot\H^n(\Sigma_j)$$
and so using \eqref{E:bad-set-j} we see that this is at most
$$C(1+\sup_{B^n_{3/4}(0)}|Dv_a|\hat{E}_j)\sup_{B^n_{3/4}(0)}|Dv_a|\hat{E}_j$$
where $C = C(n,k,Q,L)\in (0,\infty)$. Therefore, we now just need to control
\begin{equation}\label{E:B6-9}
	\int_{B^n_{3/4}(0)}\sum^Q_{\alpha=1}\eta^\prime_{\bar{\delta}\hat{E}_j}(u^{\kappa,\alpha}_j(x)-\hat{E}_jv^\kappa_a(x))\nabla^{V_j}v^\kappa_a\cdot (\nabla^{V_j}v_j^{\kappa,\alpha}-\nabla^{V_j}v^\kappa_a)\zeta^2 J^\alpha_j\, \ext x.
\end{equation}
When we replace the $\nabla^{V_j}$ derivatives by the usual derivatives in the $P_0$ directions (by \eqref{E:Lipschitz-approx-R2}), we incur an error that can be controlled by $C\hat{E}_j$. We can also replace $J^\alpha_j$ in the integrand by $1$ which introduces an error term which using \eqref{E:Lipschitz-approx-2} can be controlled by $O(\hat{E}_j)$. Thus, this reduces controlling \eqref{E:B6-9} to controlling
\begin{equation}\label{E:B6-10}
	\int_{B^n_{3/4}(0)}\sum^Q_{\alpha=1}\eta^\prime_{\bar{\delta}\hat{E}_j}(u^{\kappa,\alpha}_j(x)-\hat{E}_jv^\kappa_a(x))Dv^\kappa_{a}\cdot(Dv^{\kappa,\alpha}_j-Dv^\kappa_a)\zeta^2\, \ext x.
\end{equation}
To bound this, note that it is in fact equal to
$$\hat{E}_j^{-1}\int_{B^n_{3/4}(0)}D\left(\sum^Q_{\alpha=1}\eta_{\bar{\delta}\hat{E}_j}(u^{\kappa,\alpha}_j(x)-\hat{E}_jv^\kappa_{a}(x))\right)\cdot Dv^\kappa_a\, \zeta^2\, \ext x$$
which we can integrate by parts to get that it is equal to
$$\hat{E}_j^{-1}\int_{B^n_{3/4}(0)}\sum^Q_{\alpha=1}\eta_{\bar{\delta}\hat{E}_j}(u^{\kappa,\alpha}_j(x)-\hat{E}_jv^\kappa_a(x))(\Delta v^\kappa_a + Dv^\kappa_a\cdot 2\zeta D\zeta)\, \ext x.$$
Recall from $(\Bfrak2)$ that $v^\kappa_a$ is a classically harmonic function, and so one part of the integrand vanishes. Using the fact that $\sup|\eta_{\bar{\delta}\hat{E}_j}|\leq 2\bar{\delta}\hat{E}_j$, the remaining part of the integral is bounded in absolute value by
$$20\bar{\delta}\sup_{B^n_{3/4}(0)}|Dv_a|.$$
Thus combining all of this in Step 2, we see that the total error term \eqref{E:B6-8} is controlled by
$$40\bar{\delta}\sup_{B^n_{3/4}(0)}|Dv_a| + C(1+\sup_{B^n_{3/4}(0)}|Dv_a|\hat{E}_j)\sup_{B^n_{3/4}(0)}|Dv_a|\hat{E}_j + C\hat{E}_j$$
i.e.~it is bounded by $40\bar{\delta}\sup_{B^n_{3/4}(0)}|Dv_a| + O(\hat{E}_j)$.

\textbf{Step 3:} \emph{Analyse the main term \eqref{E:B6-7}.} As in Step 2, we can replace $J^\alpha_j$ by $1$ up to introducing an error term which is $O(\hat{E}_j)$. We now want to replace the $\nabla^{V_j}$ by the usual derivatives, the issue here is that we only have the pointwise bound of $L\hat{E}_j^{-1}$ on $Dv_j^{\kappa,\alpha}$ and there is a quadratic term in $\nabla^{V_j}v_j^{\kappa,\alpha}$. Thus, \eqref{E:Lipschitz-approx-R2} allows us to replace $\nabla^{V_j}$ by the usual derivatives, up to introducing an error term of the form
$$C(n,k)(L^2+\sup_{B^n_{3/4}(0)}|Dv_a|^2\hat{E}_j^2).$$
If we now combine this with all the steps so far, we currently have
\begin{align}
\nonumber\int_{\Omega_j}\sum^Q_{\alpha=1}\eta^\prime_{\bar{\delta}\hat{E}_j}(u^{\kappa,\alpha}_j(x)-\hat{E}_jv^\kappa_a(x))&|Dv^{\kappa,\alpha}_j-Dv_a^\kappa|^2\zeta^2\, \ext x\\
& \leq C_*\bar{\delta}(1+\sup_{B^n_{3/4}(0)}|Dv_a|) + C_*L^2 + O(\hat{E}_j)\label{E:B6-11}
\end{align}
where here we stress that the constant $C_*$ only depends on $n,k,Q$ and \emph{not} on $L$ (this will be important as after we send $j\to\infty$, we will send $L\downarrow 0$).

Let us now focus on the remaining integral. Since $\eta^\prime_{\bar{\delta}\hat{E}_j}(u_j^{\kappa,\alpha}(x)-\hat{E}_jv^\kappa_a(x))\equiv 1$ on the set $\{x\in B^n_1(0):|v_j^{\kappa,\alpha}(x)-v^k_a(x)|<\bar{\delta}\}$, and $\zeta\equiv 1$ on $B^n_{1/2}(0)$, we have
\begin{align}
\nonumber\int_{\Omega_j}\sum^Q_{\alpha=1}\eta^\prime_{\bar{\delta}\hat{E}_j}(u_j^{\kappa,\alpha}(x)&-\hat{E}_jv^\kappa_a(x))|Dv_j^{\kappa,\alpha}(x)-Dv_a^\kappa(x)|^2\zeta^2\, \ext x\\
& \geq \int_{\Omega_j\cap B^n_{1/2}(0)}\sum^Q_{\alpha=1}\one_{\{|v_j^{\kappa,\alpha}-v^\kappa_a|<\bar{\delta}\}}(x)|Dv_j^{\kappa,\alpha}(x)-Dv_a^\kappa(x)|^2\, \ext x\nonumber\\
& \geq \int_{\Omega_j\cap B^n_{1/2}(0)\cap U}\sum^Q_{\alpha=1}\one_{\{|v_j^{\kappa,\alpha}-v^\kappa_a|<\bar{\delta}\}}(x)|Dv_j^{\kappa,\alpha}(x)-Dv_a^\kappa(x)|^2\, \ext x\label{E:B6-12}
\end{align}
where recall $U$ is as in \eqref{E:B6-1} and \eqref{E:B6-2}. Since $v^{\kappa,\alpha}_j\to v^{\kappa,\alpha}$ uniformly on $U$ we have, for all sufficiently large $j$ (depending on $\bar{\delta}$), that
$$U\cap \{|v^{\kappa,\alpha}-v_a^\kappa|<\bar{\delta}/2\}\subset U\cap \{|v^{\kappa,\alpha}_j-v^\kappa_a|<\bar{\delta}\}.$$
Then, for all $j$ sufficiently large, \eqref{E:B6-12} is at least
$$\int_{\Omega_j\cap B^n_{1/2}(0)\cap U}\sum^Q_{\alpha=1}\one_{\{|v^{\kappa,\alpha}-v^\kappa_a|<\bar{\delta}/2\}}(x)|Dv_j^{\kappa,\alpha}(x)-Dv_a^\kappa(x)|^2\, \ext x.$$
Now recall the following general fact: if
\begin{itemize}
	\item $(f_j)_j$ is uniformly bounded in $L^2(B_1^n(0))$, $f\in L^2(B^n_1(0))$, and $f_j\weakly f$ weakly in $L^2(B^n_1(0))$;
	\item $(h_j)_j$ is uniformly bounded in $L^\infty(B^n_1(0))$ with $h_j\to h$ almost everywhere in $B^n_1(0)$;
\end{itemize}
then we have $h_jf_j\weakly hf$ weakly in $L^2(B^n_1(0))$. Thus if we take $\liminf_{j\to\infty}$ in the above and apply this fact, using also that $\one_{\Omega_j}\to \one_{B^n_1(0)}$ a.e. and the lower semi-continuity of the $L^2$ norm under weak $L^2$ convergence, we get that
\begin{align*}
\liminf_{j\to\infty}\int_{\Omega_j\cap B^n_{1/2}(0)\cap U}\sum^Q_{\alpha=1}&\one_{\{|v^{\kappa,\alpha}-v^\kappa_a|<\bar{\delta}/2\}}(x)|Dv_j^{\kappa,\alpha}(x)-Dv_a^\kappa(x)|^2\, \ext x\\
& \geq \int_{B^n_{1/2}(0)\cap U}\sum^Q_{\alpha=1}\one_{\{|v^{\kappa,\alpha}-v^\kappa_a|<\bar{\delta}/2\}}(x)|Dv^{\kappa,\alpha}(x)-Dv^\kappa_a(x)|^2\, \ext x.
\end{align*}
Thus, currently we have
$$\int_{B^n_{1/2}(0)}\sum^Q_{\alpha=1}\one_{\{|v^{\kappa,\alpha}-v^\kappa_a|<\bar{\delta}/2\}}(x)|Dv^{\kappa,\alpha}(x)-Dv^\kappa_a(x)|^2\, \ext x \leq C_*\bar{\delta}(1+\sup_{B^n_{3/4}(0)}|Dv_a|) + C_*L^2.$$
We can add back in the part of the integral over $B^n_{1/2}(0)\setminus U$ and use \eqref{E:B6-2} to get
$$\int_{B^n_{1/2}(0)}\sum^Q_{\alpha=1}\one_{\{|v^{\kappa,\alpha}-v^\kappa_a|<\bar{\delta}/2\}}(x)|Dv^{\kappa,\alpha}(x)-Dv^\kappa_a(x)|^2\, \ext x \leq C_*\bar{\delta}(1+\sup_{B^n_{3/4}(0)}|Dv_a|) + C_*L^2 + C_*\eps.$$
Now, the choice of constant $\eps>0$ was arbitrary, so we can let $\eps\downarrow 0$ in the above. Also, as already remarked, $v$ is independent of the choice of Lipschitz constant $L$ and $C_*$ is independent of $L$, and so we can let $L\downarrow 0$ as well. Next, as $v$ is bounded in $W^{1,2}(B^n_{7/8}(0))$ by some $C(n,k,Q)$, it follows from elliptic estimates for harmonic functions that $\sup_{B^n_{3/4}(0)}|Dv_a|\leq C(n,k,Q)$. Thus, we end up with
$$\int_{B^n_{1/2}(0)}\sum^Q_{\alpha=1}\one_{\{|v^{\kappa,\alpha}-v^\kappa_a|<\bar{\delta}/2\}}(x)|Dv^{\kappa,\alpha}(x)-Dv^\kappa_a(x)|^2\, \ext x \leq C\bar{\delta}.$$
This completes the proof of $(\mathfrak{B}6)$.
\end{proof}

\textbf{Remark:} In fact, the argument for $(\mathfrak{B}6)$ shows a more general statement, namely that if $h$ is \emph{any} harmonic function $B^n_1(0)\to \R$, then for any $v\in \mathfrak{B}_Q$, $\delta>0$, and $\kappa\in\{1,\dotsc,k\}$, we have
$$\sum^Q_{\alpha=1}\int_{B^n_{1/2}(0)\cap \{|v^{\kappa,\alpha}-h|<\delta\}}|Dv^{\kappa,\alpha} - Dh|^2\leq C\delta\left(\int_{B_1(0)}|h|^2\right)^{1/2}$$
where $C = C(n,k,Q)\in (0,\infty)$, i.e.~we get an energy non-concentration estimate relative to \emph{any} harmonic function. Moreover, noting that $\{|v^\kappa-h|<\bar{\delta}/2\}\subset\{|v^{\kappa,\alpha}-h|<\bar{\delta}/2\}$ for all $\alpha=1,\dotsc,Q$, this also implies
$$\int_{B^n_{1/2}(0)\cap \{|v^\kappa - h|<\bar{\delta}/2\}}\sum^Q_{\alpha=1}|Dv^{\kappa,\alpha}-Dh|^2\, \ext x \leq C\bar{\delta}\left(\int_{B_1(0)}|h|^2\right)^{1/2}.$$

\part{\centering Fine Regularity near Classical Cones}\label{part:fine-reg}

For the remainder of the paper we will focus on the case $Q=2$.

In this part, we establish a regularity result (Theorem~\ref{thm:fine-reg} below) which is a special case of our main result Theorem~\ref{thm:main}. We shall refer this special case as a \emph{fine $\eps$-regularity theorem}. We will use Theorem~\ref{thm:fine-reg} first as a key ingredient in the analysis of coarse blow-ups of sequences of varifolds $V \in \mathcal{V}_\beta$ converging to a multiplicity $2$ plane, and again in a more direct way in the proof of Theorem~\ref{thm:main}. In Theorem~\ref{thm:fine-reg}, in addition to $\hat{E}_V<\eps$, we impose the hypothesis
\begin{equation}\label{E:special-case}\tag{$\star$}
\hat{F}_{V,\BC}<\gamma\hat{E}_V
\end{equation}
for suitably small $\gamma = \gamma(n,k,\beta)\in (0,1)$, where 
$\BC\in \CC$, that is, $\BC$ is a cone equal to either a transversely intersecting pair of planes or a twisted union of $4$ half-planes, and 
$\hat{F}_{V,\BC}$ is the two-sided fine excess of the varifold $V$ relative to $\BC$ (see \eqref{E:fine-excess} for a precise definition). Under this additional hypothesis, stronger conclusions than in Theorem~\ref{thm:main} will follow: it will turn out that Theorem~\ref{thm:fine-reg} will yield that \emph{all}  tangent cones to $V\res C_{1/2}(0)$ at singular points must be in $\CC$, that is, $V$ has no branch point singularities; and moreover that the singular set of 
$V \res C_{1/2}(0)$ is contained within an $(n-1)$-dimensional $C^{1,\alpha}$ submanifold, and that away from the singular set, $V \res C_{1/2}(0)$ is locally represented by the (disjoint) graphs of two smooth functions over an appropriate ball in the plane $\{0\}^{k} \times {\mathbb R}^{n}$, solving the minimal surface system and satisfying uniform estimates.

Analysis of the structure of a varifold under an assumption of the form \eqref{E:special-case} first appeared in \cite{Wic14} in a simpler, branch-point-free setting. Subsequently, this analysis was taken further and made more quantitative in the work \cite{MW24} (see \cite[Theorem 3.1]{MW24}), where it was used to prove a local structure theorem for certain, possibly branched, stable codimension 1 integral varifolds which, as a special case, include codimension 1 area minimising currents mod $p$. Since then, assumptions analogous to \eqref{E:special-case} have been used in various other contexts (e.g.~\cite{Min21a, KW21, KW23b, DLMS24}). However, in these works it has always been the case that one can work under the additional assumption that no ``significant density gaps'' occur, by which we mean that $S(\BC)$ is contained in a small neighbourhood of the set of singularities of $V$ with density $\geq \Theta_{\BC}(0)$. For us here, this is certainly not the case as density gaps can occur. Allowing for density gaps makes the analysis significantly more complicated, and forces us to argue in a manner similar to the work of the first author \cite{BK17}, utilising our topological structural assumption to understand the structure of the blow-up in density gap regions. The work \cite{BK17} however was in the context of $2$-valued Lipschitz graphs, and also, in more substantive contrast to the present setting, in a non-degenerate setting where $V$ is far from being planar (i.e.~$\hat{E}_V\geq \eps>0$ for suitable $\eps$ depending on $\BC$).

The precise statement of the fine $\eps$-regularity theorem is as follows.

\begin{thmx}[Fine $\eps$-Regularity Theorem]\label{thm:fine-reg}
	Fix $M\in [1,\infty)$ and $\beta\in (0,1)$. Then, there exists $\eps = \eps(n,k,M,\beta)\in (0,1)$ and $\gamma = \gamma(n,k,M,\beta)\in (0,1)$ such that the following is true. Suppose that $V\in\mathcal{V}_\beta$ and $\BC\in \CC$ satisfy:
	\begin{enumerate}
		\item [(a)] $\Theta_V(0)\geq 2$ and $(\w_n 2^n)^{-1}\|V\|(B_2^{n+k}(0))\leq \frac{5}{2}$;
		\item [(b)] $0\in S(\BC)\subset S_0 := \{0\}^{k+1}\times\R^{n-1}\subset P_0$;
		\item [(c)] $\nu(\BC)<\eps$ and $\hat{E}_V<\eps$;
		\item [(d)] $\hat{E}^2_V \leq M \inf_P \int_{C_1(0)}\dist^2(x,P)\, \ext\|V\|(x)$, where the infimum is taken over all planes $P$ with $S(\BC)\subset P$;
		\item [(e)] $\hat{F}_{V,\BC}<\gamma\hat{E}_V$.
	\end{enumerate}
	Then, there is a cone $\BC_0\in \CC$ with $0\in S(\BC_0)\subset S_0$ and an orthogonal rotation $\Gamma:\R^{n+k}\to \R^{n+k}$ with
	$$\dist_\H(\BC_0\cap C_1(0), \BC\cap C_1(0))\leq C\hat{F}_{V,\BC},$$
	$|\Gamma(e_\kappa)-e_\kappa|\leq C\hat{F}_{V,\BC}$ for $\kappa=1,\dotsc,k$ and $|\Gamma(e_{k+i})-e_{k+i}|\leq C\hat{E}_{V}^{-1}\hat{F}_{V,\BC}$ for $i=1,\dotsc,n$, such that $\BC_0$ is the unique tangent cone to $\Gamma^{-1}_\#V$ at $0$ and
	$$\sigma^{-n-2}\int_{C_\sigma(0)}\dist^2(x,\BC_0)\, \ext\|\Gamma^{-1}_\#V\|(x) \leq C\sigma^{2\mu}\hat{F}_{V,\BC}^2 \qquad \text{for all }\sigma\in (0,1/2).$$
	Furthermore, there exists a generalised-$C^{1,\mu}$ function $u:B^n_{1/2}(0)\to \A_2(P_0^\perp)$ such that:
	\begin{enumerate}
		\item [(1)] $V\res (C_{1/2}(0)\cap B_{3/2}^{n+k}(0)) = \mathbf{v}(u)$;
		\item [(2)] $\mathcal{B}_u\cap B_{1/2} = \emptyset$;
		\item [(3)] $\graph(u)\cap (\R^k\times \CC_u) = \sing(V)\cap C_{1/2}(0)\cap B^{n+k}_{3/2}(0)$; moreover,
		\begin{enumerate}
			\item [(i)] $\sing(V)\cap C_{1/2}(0)\cap B_{3/2}^{n+k}(0) \subset \graph(\phi)$, where $\phi: S_0\cap B_{1/2}^{n-1}(0)\to S_0^\perp$ is of class $C^{1,\mu}$ over $\overline{S_0\cap B_{1/2}^{n-1}(0)}$, where if $\phi = (\phi_1,\dotsc,\phi_k,\phi_{k+1})$, then $|\phi_\kappa|_{1,\mu;B_{1/2}^{n-1}(0)}\leq C\hat{F}_{V,\BC}$ for $\kappa=1,\dotsc,k$ and $|\phi_{k+1}|_{1,\mu;B_{1/2}^{n-1}(0)}\leq C\hat{E}_{V}^{-1}\hat{F}_{V,\BC}$;
			\item [(ii)] If $\Omega^\pm$ denote the two connected components of $B_{1/2}^n(0)\setminus \graph(\phi_{k+1})$, then on each of $\Omega^\pm$ we have that we can write $u|_{\Omega^\pm} = \llbracket u^1_\pm\rrbracket + \llbracket u^2_\pm\rrbracket$, where $u^\alpha_\pm \in C^{1,\mu}(\overline{\Omega^\pm};\R^k)$ and $|u^\alpha_\pm|_{1,\mu;\Omega^\pm}\leq C\hat{E}_V$ for each $\alpha=1,2$;
		\end{enumerate}
		\item [(4)] If $\BC_z\in \CC$ denotes the (unique) tangent cone to $V$ at $z\in \sing(V)\cap C_{1/2}(0)\cap B^{n+k}_{3/2}(0)$, then 
		$$C^{-1}\hat{E}_V\leq \dist_\H(\BC_z\cap C_1(0), P_0\cap C_1(0))\leq C\hat{E}_V.$$
		Moreover, if $z_1,z_2\in \sing(V)\cap C_{1/2}(0)\cap B^{n+k}_{3/2}(0)$ then
		$$\dist_\H(\BC_{z_1}\cap C_1(0), \BC_{z_2}\cap C_1(0))\leq C|\pi_{P_0}(z_1)-\pi_{P_0}(z_2)|^{\mu}\hat{F}_{V,\BC}.$$
	\end{enumerate}
	Here, $C = C(n,k,M,\beta)\in (0,\infty)$ and $\mu = \mu(n,k,M,\beta)\in (0,1)$.
\end{thmx}

\textbf{Remark:} The assumption that $\Theta_V(0)\geq 2$ can easily be removed: if there are no density $2$ points, the result follows immediately from Lemma \ref{lemma:gap-2} or Allard's regularity theorem, depending on whether $V$ is close to a multiplicity one or two copy of $P_0$.

\begin{remark}\label{remark:fine-reg-sufficiency}
Similarly to that seen in Section \ref{sec:non-degenerate}, one can in fact prove Theorem \ref{thm:fine-reg} under the weaker assumption that the topological structural condition is always satisfied in $(\beta,\tilde{\gamma})$-fine gaps, provided one allows the constants to depend on $\tilde{\gamma}$ as well.
\end{remark}

We briefly elaborate on the two ways in which we will use Theorem \ref{thm:fine-reg}. The first is in understanding coarse blow-ups of stationary integral varifolds in $\mathcal{V}_\beta$ converging to a plane of multiplicity $2$. Looking towards establishing an excess decay lemma, we would like to know that coarse blow-ups satisfy a $C^{1,\alpha}$ estimate of some type. How would one prove such an estimate? Qualitatively, if $v$ is a $C^{1,\alpha}$ function, we would expect to see points on $\graph(v)$ where the two values of $v$ intersect transversely, and so the tangent map is a pair of transversely intersecting planes. The regions in the sequence of stationary integral varifolds that produced such a region in the coarse blow-up are then regions which obey additionally \eqref{E:special-case}. Therefore, understanding the situation where \eqref{E:special-case} also holds intuitively allows us to understand regions in the coarse blow-up where transverse intersection of the sheets occurs, and once we understand these regions we can then hope to try to understand the coarse blow-up as a whole.

The second way in which we will use Theorem \ref{thm:fine-reg} is in establishing the relevant excess decay lemma used to prove Theorem \ref{thm:main}. It will be convenient to phrase this as an excess decay \emph{dichotomy}, essentially saying that, when one produces the coarse blow-up in the usual manner for a proof of an excess decay statement, if the tangent map to the coarse blow-up at the origin is a plane, then we get decay to a new plane at some smaller scale, or instead the tangent map is a transverse intersection of two planes (or a twisted cone), in which case we can invoke Theorem \ref{thm:fine-reg} to directly conclude a full structural result on a smaller ball. It is then sufficient to iterate such a dichotomy in order to conclude the full conclusions of Theorem \ref{thm:main}.

\section{Further Notation}\label{sec:further-notation}

Recall the notation $\mathcal{P}$, $\mathcal{P}_j$, $\mathcal{P}_{\leq j}$, $\CC_{n-1}$, $\CC$, and so forth from Section \ref{sec:prelim}. Recall that $\CC\in \CC$ is not technically a cone unless $0\in S(\BC)$.

For any $\BC\in \CC$, we define:
\begin{itemize}
	\item $r_\BC \equiv r_\BC(x):= \dist(x,S(\BC))$ for the distance of a point $x$ to $S(\BC)$.
	\item $\sfrak(\BC):= \dim(S(\BC))$ for the dimension of $S(\BC)$.
	\item If $\BC\in\CC$ is such that $0\in S(\BC)\subset P_0$, we say that $\BC$ is \emph{aligned}.
\end{itemize}
The integer $\sfrak(\BC)$ will play an important role, as it will be the parameter over which we perform a long and complex induction argument in order to prove the main $L^2$ estimates. A key point to note is that if two planes intersect transversely, small perturbations of the planes can only \emph{decrease} the dimension of the intersection.

If $\BC\in\CC$ is aligned, we typically use the following notation:
\begin{itemize}
	\item $e_1,\dotsc,e_n$ is an orthonormal basis of $P_0$ in such a way that $e_1,\dotsc,e_{\sfrak(\BC)}$ is an orthonormal basis of $S(\BC)$.
	\item $\ell:P_0\to \A_2(P_0^\perp)$ is the two-valued function for which $\BC = \graph(\ell)$.
	\item If $\BC = |P^1|+|P^2|\in \mathcal{P}$, then typically the planes $P^1$ and $P^2$ will be very close to $P_0$. For $\alpha=1,2$ we will write $p^\alpha:P_0\to P_0^\perp$ for the function with $P^\alpha = \graph(p^\alpha)$.
	\item If $\BC\in \CC_{n-1}$, we label the half-planes $H^\alpha$ in $\BC$ so that $H^1$ and $H^2$ are graphical over $\{x_n<0\}\cap P_0$ and $H^3$ and $H^4$ are graphical over $\{x^n>0\}\cap P_0$. We then write $h^\alpha$ for the relevant linear function whose graph over the relevant domain (either $\{x^n\leq 0\}\cap P_0$ or $\{x^n\geq 0\}\cap P_0$) is $H^\alpha$.
	\item If $\BC = |P^1| + |P^2|\in \mathcal{P}_{n-1}$ then, in addition to the bullet point above, we label the planes and half-planes so that $H^1\cup H^3 = P^1$ and $H^2\cup H^4 = P^2$. In particular, $p^1|_{\{x^n\leq 0\}} = h^1$, and so forth.
	\item We write $U^1_0 = U^2_0 = \{x^n<0\}\cap P_0$ and $U^3_0 = U^4_0 = \{x^n>0\}\cap P_0$. This notation allows for us to write expressions more simply using sums, for example,
	$$\sum^4_{\alpha=1}\int_{U^\alpha_0\cap B^n_{1/2}(0)}|h^\alpha(x)|^2\, \ext x.$$
\end{itemize}
It will be convenient to have a flexible language to discuss tubular regions around (subsets of) subspaces. For any subspace $S\subsetneq P_0$, write
\begin{itemize}
	\item $\Lambda(S):= \{(x,s)\in S\times [0,1): |x|^2+s^2<1\}$.
	\item For $x\in B_1(0)$, write $\Pi_S(x):= (x^{\top_S},|x^{\perp_S}|)\in \Lambda(S)$.
	\item For $U\subset \Lambda(S)$, write $\mathcal{T}^S(U):= \{x\in B_1(0):\Pi_S(x)\in U\}$ for the toroidal region about $S$ generated by $U$.
\end{itemize}
We also introduce notation for cubes centred at choices of parameters $(x_0,s_0)\in \Lambda(S)$:
\begin{itemize}
	\item Write $Q_r(x_0,s_0)$ for the \emph{open} $(\dim(S)+1)$-dimensional cube in $\Lambda(S)$ centred at $(x_0,s_0)$ with edge length $2r$ (i.e.~the ball of radius $r$ in the $\|\cdot\|_\infty$ metric).
	\item For a cube $Q$, write $e(Q)$ for its edge length.
	\item Given $\lambda>0$, write $\lambda Q$ for the dilation of $Q$ by a factor of $\lambda$ about its own centre, i.e.~the cube with the same centre but edge length $\lambda e(Q)$.
	\item We say that (open) cubes $Q_1$, $Q_2$ are \emph{adjacent} if they are disjoint yet $\del Q_1\cup \del Q_2\neq\emptyset$.
	\item The lattice of points in $S\times [0,\infty)$ whose coordinates are integral determines a collection $\mathscr{C}_0$ of cubes of unit edge length, the vertices of which are points of the lattice. This collection defines a sequence $\{\mathscr{C}_j\}_{j=1}^\infty$ of collections of cubes in the following way: for $p\geq 0$, each cube $Q\in \mathscr{C}_p$ determines $2^{\dim(S)+1}$ cubes in $\mathscr{C}_{p+1}$ by bisecting the edges of $Q$. We write $\mathscr{C}:=\cup^\infty_{j=0}\mathscr{C}_j$.
\end{itemize}
\textbf{Note:} Some of the above notation suppresses the dependence on the subspace $S$.

For $0<r\leq s_0$, we then define
$$\mathcal{T}^S_r(x_0,s_0):= \mathcal{T}^S(Q_r(x_0,s_0)).$$
We can unravel the definitions to see that
$$x\in \mathcal{T}_r^S(x_0,s_0)\ \ \Longleftrightarrow\ \ (x^{\top_S},|x^{\perp_S}|)\in Q_r(x_0,s_0).$$
Notice that because $0<r\leq s_0$, we always have $S\cap \mathcal{T}^S_r(x_0,s_0) = \emptyset$. It will also be convenient to extend the notation $\mathcal{T}_r^S(x_0,s_0)$ to allow $s_0=0$ but $r>0$; for this we define $\mathcal{T}_r^S(x_0,0) := B_r(x_0)$ (we stress that this is a \emph{ball}, not a toroidal region determined by a cube, and is also independent of $S$).

\textbf{Remark:} Notice that for $0<r\leq s_0$,
\begin{equation}\label{E:cube-R1}
	Q_r(x_0,s_0)\cap Q_{r^\prime}(x_0^\prime,s_0^\prime)\neq\emptyset\ \ \ \ \Longleftrightarrow\ \ \ \ \mathcal{T}_r^S(x_0,s_0)\cap \mathcal{T}_{r^\prime}^S(x_0^\prime,s_0^\prime) \neq\emptyset
\end{equation}
and if $x\in \mathcal{T}_r^S(x_0,s_0)$, then
\begin{align*}
	|x| & = |x^{\top_S}+x^{\perp_S}|\\
	& \leq \left|x^{\top_S} - x_0 + x^{\perp_S} - s_0\frac{x^{\perp_S}}{|x^{\perp_S}|}\right| + \left|x_0 + s_0\frac{x^{\perp_S}}{|x^{\perp_S}|}\right|\\
	& = \left(|x^{\top_S}-x_0|^2 + \left|x^{\perp_S}-s_0\frac{x^{\perp_S}}{|x^{\perp_S}|}\right|^2\right)^{1/2} + (|x_0|^2 + s_0^2)^{1/2}\\
	& = \left(|x^{\top_S}-x_0|^2 + (|x^{\perp_S}|-s_0)^2\right)^{1/2} + (|x_0|^2+s_0^2)^{1/2}\\
	& \leq r\sqrt{n} + (|x_0|^2+s_0^2)^{1/2}
\end{align*}
i.e.~we have
\begin{equation}\label{E:cube-R2}
	\mathcal{T}_r^S(x_0,s_0) \subset B_R(0) \qquad \text{for }R = r\sqrt{n}+(|x_0|^2+s_0^2)^{1/2}.
\end{equation}
Next, we recall and introduce notation for the various $L^2$ excess quantities we use. Let $V$ be a stationary integral $n$-varifold and $P$ be a plane through the origin.
\begin{itemize}
	\item The $L^2$ \emph{coarse excess} of $V$ relative to $P$ is $\hat{E}_{V,P}:=\left(\int_{B_1(0)}\dist^2(x,P)\, \ext\|V\|(x)\right)^{1/2}$.
	\item Write $\hat{E}_V:= \hat{E}_{V,P_0}$.
\end{itemize}
Often we will arrange or assume that $\hat{E}_V$ is comparable to $\mathcal{E}_V:=\inf_P\hat{E}_{V,P}$, where the infimum is taken over all $n$-dimensional planes $P$ passing through $0$, and so we informally refer to $\hat{E}_V$ as `the' coarse excess of $V$.
\begin{itemize}
	\item For $S\subsetneq P$, $(x_0,s_0)\in \Lambda(S)$, $0<r\leq s_0$, we write
	$$\hat{E}_{V,P}(\mathcal{T}^S_r(x_0,s_0)):= \left(r^{-n-2}\int_{\mathcal{T}^S_r(x_0,s_0)}\dist^2(x,P)\, \ext\|V\|(x)\right)^{1/2}$$
	for the scale-invariant excess in the region $\mathcal{T}^S_r(x_0,s_0)$. The excess $\hat{E}_{V,P}(\mathcal{T}^S_r(x_0,0))\equiv \hat{E}_{V,P}(B_r(x_0))$ is defined analogously.
	\item Analogously, we can define $\hat{E}_V(\mathcal{T}^S_r(x_0,s_0))$.
\end{itemize}
For $\BC\in \CC$, we then define
\begin{itemize}
	\item The (\emph{one-sided}) \emph{excess of $V$ relative to $\BC$} by $\hat{E}_{V,\BC}:=\left(\int_{B_1(0)}\dist^2(x,\BC)\, \ext\|V\|(x)\right)^{1/2}$, where for shorthand we write $\dist(x,\widetilde{V}) \equiv \dist(x,\spt\|\widetilde{V}\|)$ for a varifold $\widetilde{V}$.
	\item The (\emph{two-sided}) \emph{fine excess of $V$ relative to $\BC$} is then
	\begin{equation}\label{E:fine-excess}
		\hat{F}_{V,\BC}:=\left(\int_{B_1(0)}\dist^2(x,\BC)\, \ext\|V\|(x) + \int_{B_{1/2}\setminus \{r_{\BC}<1/8\}}\dist^2(x,V)\, \ext\|\BC\|(x)\right)^{1/2}.
	\end{equation}
	\item We can then define quantities $\hat{E}_{V,\BC}(\mathcal{T}^S_r(x_0,s_0))$, $\hat{F}_{V,\BC}(T^S_r(x_0,s_0))$.
    For instance,
	\begin{align*}
	\hat{F}^2_{V,\BC}(\mathcal{T}^S_r(x_0,s_0)) := r^{-n-2}\int_{\mathcal{T}_r^S(x_0,s_0)}&\dist^2(x,\BC)\, \ext\|V\|(x)\\
	& + r^{-n-2}\int_{\mathcal{T}^S_{r/2}(x_0,s_0)\setminus \{r_{\BC}<r/8\}}\dist^2(x,V)\, \ext\|\BC\|(x).
	\end{align*}
\end{itemize}
The most delicate parts of the proofs rely on comparisons between the fine excess relative to a fixed cone $\BC\in\CC$ and the fine excess relative to the best-approximating \emph{coarser} cone, i.e.~
\begin{itemize}
	\item If $\BC,\Dbf\in \CC$ and $S(\BC)\subsetneq S(\Dbf)$, we say that $\Dbf$ is \emph{coarser} than $\BC$.
\end{itemize}
Intuitively, coarser cones have fewer degrees of freedom, and thus provide a coarser approximation (hence their name). We then introduce the following notation:
\begin{itemize}
	\item Write $\mathscr{E}_{V,\BC}:=\inf_{P}\hat{E}_{V,P}$, where here the infimum is over all coarser \emph{planes} to $\BC$.
	\item Similarly we can define $\mathscr{E}_{V,\BC}(\mathcal{T}^S_r(x_0,s_0))$, etc.
	\item Define
	$$\mathscr{E}^*_{V,\BC}:=\begin{cases}
		\mathscr{E}_{V,\BC} & \text{if }\sfrak(\BC)=n-1;\\
		\inf_{\Dbf}\hat{E}_{V,\Dbf} & \text{if }\sfrak(\BC)<n-1,
	\end{cases}$$
	where here the infimum is taken over all coarser \emph{unions of planes} $\Dbf\in \mathcal{P}$ to $\BC$.
	\item If $0\in S(\BC)$, define
	$$\mathscr{F}^*_{V,\BC}:=\begin{cases}
		\mathscr{E}_{V,\BC} & \text{if }\sfrak(\BC)=n-1;\\
		\inf_{\Dbf}\hat{F}_{V,\Dbf} & \text{if }\sfrak(\BC)<n-1,
	\end{cases}$$
	where again the infimum is taken over all coarser \emph{unions of planes} $\Dbf\in\mathcal{P}$ to $\BC$.
	\item Similarly we can define $\mathscr{F}^*_{V,\BC}(\mathcal{T}^S_r(x_0,s_0))$, etc.
\end{itemize}
\textbf{Note:} Notice that $\mathscr{F}^*_{V,\BC}\leq C\mathscr{E}_{V,\BC}$ always holds for some $C = C(n,k)$, i.e.~the infimum over coarser unions of planes is always smaller than the infimum over single coarser planes. Of course, if $\sfrak(\BC)=n-1$ this is true by definition with $C=1$. When $\sfrak(\BC)<n-1$, we can approximate any single coarse plane $P$ by a coarse union of two planes converging to $P$, meaning that $\mathscr{F}^*_{V,\BC}$ is controlled by the infimum of the two-sided excess of $V$ relative to any coarse plane $P$ (for the two-sided excess relative to a plane, in the second term we take the integral over $B_{1/2}$, say). The second term in the two-sided excess relative to $P$ can however be controlled by the one-sided excess of $V$ relative to $P$ by using by Theorem \ref{thm:allard-sup-estimate} and Theorem \ref{thm:Lipschitz-approx}; this therefore gives $\mathscr{F}^*_{V,\BC}\leq C\mathscr{E}_{V,\BC}$ in this case for suitable $C = C(n,k)$. Also, clearly by the same argument we have $\mathscr{E}^*_{V,\BC}\leq \mathscr{E}_{V,\BC}$ always.

Thus, the quantities $\mathscr{E}^*_{V,\BC}$ and $\mathscr{F}^*_{V,\BC}$ represent the optimal excess quantities with respect to \emph{coarser} unions of planes. (We stress that, the asterisk $*$ signifies that the infimum is taken over \emph{unions} of planes, rather than single-planes or twisted cones, and the lack of an asterisk means the infimum is taken over \emph{single} planes which are coarser.) If, say, $\hat{F}_{V,\BC}\ll\mathscr{F}^*_{V,\BC}$, this would be saying that $V$ is significantly closer to $\BC$ than any coarser cone to $\BC$, and so we have gained something by using the cone $\BC$ rather than a cone with a larger spine.

We will also need notation for comparing cones to one another. For any $\BC,\Dbf\in \widetilde{\CC}_0$ cones or plane passing through $0$, define
\begin{itemize}
	\item $\nu(\BC):= \dist_\H(\BC\cap B_1(0),P_0\cap B_1(0))$;
	\item $\nu(\BC,\Dbf):=\dist_\H(\BC\cap B_1(0),\Dbf\cap B_1(0))$.
\end{itemize}
Here, we have used the shorthand $\BC\cap B_1(0)\equiv \spt\|\BC\|\cap B_1(0)$ in the context of Hausdorff distance. Of course, we can also make these definitions even if $\BC,\Dbf$ do not pass through $0$.

For an aligned cone $\BC\in\mathcal{P}$ close to $P_0$, we define the smallest and largest \emph{apertures} of $\BC$ by
$$a_*(\BC):= \inf_{\w\in P_0\cap S(\BC)^\perp:\, |\w|=1}|p^1(\w)-p^2(\w)|,$$
$$a^*(\BC):=\sup_{\w\in P_0\cap S(\BC)^\perp:\, |\w|=1}|p^1(\w)-p^2(\w)|,$$
where here we recall that $p^1,p^2$ denoted the linear functions over $P_0$ whose graphs gives the two planes in $\BC$. Note that when $\sfrak(\BC)=n-1$, $a^*(\BC) = a_*(\BC)$, as the set of $\w$ over which the infimum and supremum are taken reduces to two points, the values of which coincide for each. If instead $\BC\in\CC\setminus\mathcal{P}$ is aligned, write $\{x_+,x_-\} = \{\w\in P_0\cap S(\BC)^\perp:|\w|=1\}$, where $\pm x^n_\pm >0$. We then define
$$a_*(\BC) := \min\{|h^1(x_-)-h^2(x_-)|,\, |h^3(x_+)-h^4(x_+)|\}$$
$$a^*(\BC) := \max\{|h^1(x_-)-h^2(x_-)|,\, |h^3(x_+)-h^4(x_+)|\}.$$
\textbf{Note:} If we knew that $\BC\in \CC\setminus\mathcal{P}$ was \emph{stationary} as a varifold, then one can readily check that the stationary condition for such a cone implies $a^*(\BC) = a_*(\BC)$. This is of course not true for a general cone in $\CC\setminus\mathcal{P}$, however under a certain smallness of excess hypothesis (which will be true for the cones of this form we ultimately care about), we will get comparability of $a_*(\BC)$ and $a^*(\BC)$; cf.~\eqref{E:nu-a-1}.

\begin{remark}\label{remark:elementary-observations}
	We can deduce, using only elementary geometric arguments, that there are dimensional constants $c = c(n,k)\in (0,1)$ and $C = C(n,k)\in (1,\infty)$ such that for a aligned cone $\BC\in\CC$, we have:
\begin{enumerate}
	\item [(1)] $c\nu(\BC)\leq \max_{\alpha}|Dp^\alpha| \leq C\nu(\BC)$ if $\BC\in\mathcal{P}$ and $c\nu(\BC)\leq\max_{\beta}|Dh^\beta|\leq C\nu(\BC)$ if $\BC\not\in\mathcal{P}$.
	\item [(2)] There is an $(n-1)$-dimensional subspace $S$ with $S(\BC)\subset S\subset P_0$ such that for $\tau\in (0,1)$ and $x\in B^n_1(0)\setminus B_\tau(S)$ we have
	$$\hspace{3em}c\tau\nu(\BC)\leq \max_{\alpha=1,2}|p^\alpha(x)|\ \ \ \ \text{if }\BC\in\mathcal{P}\ \ \ \ \ \ \text{or}\ \ \ \ \ \ c\tau\nu(\BC)\leq\max_{\beta=1,\dotsc,4}|h^\beta(x)|\ \ \ \ \text{if }\BC\not\in\mathcal{P}.$$
	(Indeed, $S$ can be taken to be the orthogonal complement in $P_0$ of a unit vector $\w\in P_0$ realising $\sup_{|\w|=1}\max_\alpha |p^\alpha(\w)|$ or $\sup_{|\w|=1}\max_{\beta}|h^\beta(\w)|$, respectively.)
	\item [(3)] For $\BC\in \mathcal{P}$ and $x\in B^n_1(0)$, we have
	$$r_{\BC}(x) a_*(\BC)\leq|p^1(x)-p^2(x)|\leq r_{\BC}(x) a^*(\BC).$$
	\item [(4)] There is an $(n-1)$-dimensional subspace $S_*$ with $S(\BC)\subset S_*\subset P_0$ such that for $\tau\in (0,1)$ and $x\in B^n_1(0)\setminus B_\tau(S_*)$, we have
	$$c\tau a^*(\BC)\leq |p^1(x)-p^2(x)|\ \ \ \ \text{if }\BC\in\mathcal{P};$$
	$$c\tau a^*(\BC)\leq \begin{cases}
		|h^1(x)-h^2(x)| & \text{for }x^n<0;\\
		|h^3(x)-h^4(x)| & \text{for }x^n>0.
	\end{cases}
	\ \ \ \ \text{if }\BC\not\in \mathcal{P}.$$
	(Indeed, $S_*$ can be taken to be the orthogonal complement in $P_0$ of a unit vector attaining the infimum in the definition of $a^*(\BC)$.)
	\item [(5)] $a^*(\BC)\leq C\nu(\BC)$.
\end{enumerate}
\end{remark}

The last facts we need are for comparing a cone to a translate of itself. These remarks are analogous to those found in \cite[Page 908]{Wic14}. Fix $z\in B_1(0)$ and a cone $\BC\in\CC$. Write $\xi:= z^{\perp_{S(\BC)}}$. Then:
\begin{itemize}
	\item For any $x\in B_1(0)$, the triangle inequality for $\dist(x,\cdot)$ and $\dist_\H$ gives
	\begin{equation}\label{E:translates-1}
		\left|\dist(x,\BC) - \dist(x,(\tau_z)_\#\BC)\right| \leq \nu(\BC,(\tau_z)_\#\BC).
	\end{equation}
	\item If $\BC$ is aligned, then using Remark \ref{remark:elementary-observations}(1) above and arguing as in \cite[Page 908]{Wic14}, we have:
	\begin{equation}\label{E:translates-2}
		\nu(\BC,(\tau_z)_\#\BC) \leq C\left[|\xi^{\perp_{P_0}}| + \nu(\BC)|\xi^{\top_{P_0}}|\right]
	\end{equation}
	where $C = C(n,k)\in (0,\infty)$.
\end{itemize}

\section{Statements of Main $L^2$ Estimates}\label{sec:fine-L2}

In this section we state the main technical results. First, we introduce the main hypotheses and give some basic consequences. The first hypotheses are simply guaranteeing that $V$ is close to a multiplicity $2$ plane in the unit cylinder. 

\textbf{Hypothesis G.} Let $\eps\in (0,1)$. We say that an integral $n$-varifold in $B_2(0)$ satisfies Hypothesis $G(\eps)$ if 
$$\frac{3}{2}\leq (\w_n 2^n)^{-1}\|V\|(B_2(0)) \leq \frac{5}{2} \qquad \text{and} \qquad \hat{E}_V<\eps.$$

The more involved hypotheses guaranteeing that $V$ is much closer to a cone $\BC\in\CC$ in a given region than any coarser cone are the following.

\textbf{Hypothesis H.} Let $\eps\in (0,1)$ and $\gamma\in (0,1)$ be small, $M\in (1,\infty)$, and let $A\subset P_0$ be a subspace. Then, we say a stationary integral $n$-varifold $V$ in $B_2(0)$ and $\BC\in \CC$ satisfies Hypothesis $H(A,\eps,\gamma,M)$ in a region $\mathcal{T}:= \mathcal{T}^A_r(x_0,s_0)\subset B_1(0)$ if
\begin{enumerate}
	\item [(A1)] $0\in A\subset S(\BC)\subset S_0:= \{0\}^{k+1}\times\R^{n-1}$ (in particular, $\BC$ is aligned);
	\item [(A2)] $\nu(\BC)<\eps$ and $\hat{E}_V(\mathcal{T})<\eps$;
	\item [(A3)] $\hat{E}_V(\mathcal{T})\leq M\mathscr{E}_{V,\BC}(\mathcal{T})$;
	\item [(A4)] $\hat{F}_{V,\BC}(\mathcal{T})<\gamma\mathscr{F}^*_{V,\BC}(T)$.
\end{enumerate}

\textbf{Note:} If $s_0 = 0$, then Hypothesis $H(A,\eps,\gamma,M)$ is independent of $A$. This will be the case for the main fine blow-up class we construct. In this case, we will just say ``Hypothesis $H(\eps,\gamma,M)$ holds in $\mathcal{T} := \mathcal{T}^A_r(x_0,0) \equiv B_r(x_0)$''.

Hypothesis (A4) is analogous to \cite[Hypothesis 10.1(5), Page 895]{Wic14} for the present situation. We note that, whilst we can fix $S(\BC) \equiv S_0$ when $\sfrak(\BC)=n-1$, when $\sfrak(\BC)<n-1$ we do not know which $\sfrak(\BC)$-dimensional subspace of $S_0$ the spine is, as indeed a cone with $\sfrak(\BC)=n-1$ can be approximated by cones with lower dimensional spine meeting along any subspace of $S_0$.

\textbf{Remark:} Hypothesis H will guarantee that all points of density $\geq 2$ in $V$ must accumulate around $S(\BC)$. This will in turn (cf.~Theorem \ref{thm:graphical-rep}) give us the ability to represent $V$ away from $S(\BC)$ by single-valued graphs over the half-planes or planes forming $\BC$. However, we will actually need the freedom to \emph{also} represent $V$ away from $S(\BC)$ over other nearby cones $\BC_*\in \CC$ also, with $S(\BC_*)\subset S(\BC)$ \emph{and even} $S(\BC_*)\subsetneq S(\BC)$, i.e.~$\sfrak(\BC_*)<\sfrak(\BC)$. This freedom will lead to an easier understanding of the fine blow-ups we construct later on. Note that we are \emph{not} claiming anything about whether Hypothesis H holds for such a cone $\BC_*$, as we only wish to represent $V$ over $\BC_*$ away from the spine of the original cone $\BC$.

We first record some immediate consequences of these hypotheses.

\begin{lemma}\label{lemma:nu-a}
	Suppose a stationary integral $n$-varifold $V$ in $B_2(0)$ and $\BC\in \CC$ satisfy Hypothesis $H(A,\eps,\gamma,M)$ in the region $\mathcal{T}:= \mathcal{T}^A_r(x_0,s_0)\subset B_1(0)$ for some $\eps\in (0,1)$ and $\gamma\in (0,1/4)$. Suppose also that $(\w_n 2^n)^{-1}\|V\|(B_2(0))\leq 5/2$. Then, there exist constants $c = c(n,k)\in (0,1)$ and $C = C(n,k)\in (1,\infty)$ such that:
	\begin{enumerate}
		\item [(i)] $c\hat{E}_V(\mathcal{T})\leq \nu(\BC)\leq C\hat{E}_V(\mathcal{T})$.
		\item [(ii)] $\mathscr{F}^*_{V,\BC}(\mathcal{T})\leq Ca_*(\BC)$ if $\sfrak(\BC)<n-1$.
		\item [(iii)] $\mathscr{E}_{V,\BC}(\mathcal{T})\leq Ca^*(\BC)$.
		\item [(iv)] $\mathscr{E}_{V,\BC}(\mathcal{T})\leq Ca_*(\BC)$ if $\sfrak(\BC)=n-1$.
	\end{enumerate}
\end{lemma}
\textbf{Remark:} In light of (A3) we can replace $\mathscr{E}_{V,\BC}(\mathcal{T})$ by $M^{-1}\hat{E}_V(\mathcal{T})$ in Lemma \ref{lemma:nu-a}(iii). In particular, we have
\begin{equation}\label{E:nu-a-2}
	C^{-1}M^{-1}\hat{E}_V(\mathcal{T}) \leq a^*(\BC) \leq C\nu(\BC) \leq C^2\hat{E}_V(\mathcal{T}).
\end{equation}
Moreover, when $\sfrak(\BC)=n-1$, from (iv), (i), and Remark \ref{remark:elementary-observations}(5) of Section \ref{sec:further-notation}, we get
\begin{equation}\label{E:nu-a-1}
	C^{-1}M^{-1}\hat{E}_V(\mathcal{T}) \leq a_*(\BC) \leq a^*(\BC) \leq C\nu(\BC) \leq C^2\hat{E}_V(\mathcal{T})
\end{equation}
i.e.~when $\sfrak(\BC)=n-1$, Lemma \ref{lemma:nu-a} tells us that $\hat{E}_V(\mathcal{T})$, $\nu(\BC)$, $a_*(\BC)$, and $a^*(\BC)$ are all comparable.

\begin{remark}\label{remark:after-nu-a-second}
We will see from the proof that Lemma \ref{lemma:nu-a}(i) does not require the full strength of Hypothesis (A4), but just the estimate $\hat{F}_{V,\BC}(\mathcal{T})<\gamma\hat{E}_V(\mathcal{T})$. If $\sfrak(\BC)<n-1$, this is much weaker.
\end{remark}

\begin{proof}
	For $x\in \spt\|V\|\cap\mathcal{T}$, the triangle inequality gives
	$$\dist^2(x,P_0) \leq 2\dist^2(x,\BC) + 2\dist^2_\H(P_0\cap B_1,\BC\cap B_1)$$
	which integrating over $\mathcal{T}$ gives $\hat{E}_V(\mathcal{T})^2 \leq 2\hat{E}_{V,\BC}(\mathcal{T})^2 + C\nu(\BC)^2$. Hypothesis (A4) then allows us to absorb the $\hat{E}_{V,\BC}(\mathcal{T})^2$ term into the left-hand side, which gives $c\hat{E}_V(\mathcal{T})\leq \nu(\BC)$, which is the first inequality in (i).
	
	For the second inequality in (i), we will work with notation that assumes $\BC\in\mathcal{P}$: minor modifications deal with the case $\BC\not\in\mathcal{P}$.
	
	Write $K:= \mathcal{T}^A_{r/4}(x_0,s_0)\cap \{r_{\BC}\geq r/8\}$. For $x\in P^1\cap K$, we have
	$$\dist^2(x,P_0) \leq 2\dist^2(x,V) + 2\dist^2_\H(V\cap B_{r}(x_0), P_0\cap B_r(x_0)).$$
	After integrating this over $P^1\cap K$ we get
	\begin{align*}
	\int_{P^1\cap K}\dist^2(x,P_0)\, \ext\H^n(x) \leq 2\int_{P^1\cap K}\dist^2&(x,V)\, \ext\H^n(x)\\
	& + Cr^n\dist^2_\H(\spt\|V\|\cap B_r(x_0), P_0\cap B_r(x_0))
	\end{align*}
	where $C = C(n,k)\in (0,\infty)$. But since
	$$r^{-n-2}\int_{P^1\cap K}\dist^2(x,V)\, \ext\H^n(x)\leq C\hat{F}_{V,\BC}(\mathcal{T})^2,$$
	Hypothesis (A4) implies that
	$$\int_{P^1\cap K}\dist^2(x,V)\, \ext\H^n(x) \leq C\int_{\mathcal{T}} \dist^2(x,P_0)\, \ext\|V\|(x)$$
	(here, we are using Theorem \ref{thm:allard-sup-estimate} to control the second half of the two-sided excess relative to a single plane by the one-sided excess). Moreover, Theorem \ref{thm:allard-sup-estimate} also gives
	$$\dist^2(V\cap B_r(x_0),P_0\cap B_r(x_0)) \leq Cr^2\hat{E}_V(\mathcal{T})^2.$$
	Combining the previous inequalities, we therefore arrive at
	$$\int_{P^1\cap K}\dist^2(x,P_0)\, \ext\H^n(x) \leq C\int_{\mathcal{T}}\dist^2(x,P_0)\, \ext\|V\|(x)$$
	for some $C = C(n,k)\in (0,\infty)$. Notice that the left-hand side here does not depend on $V$, and this in fact gives $\nu(P^1\cap K,P_0\cap K) \leq Cr^2\hat{E}_V(\mathcal{T})$, which by scaling gives $\nu(P^1)\leq C\hat{E}_V(\mathcal{T})$. Performing the same calculations for $P^2$ then gives $\nu(\BC) \leq C\hat{E}_V(\mathcal{T})$, which proves the second inequality in (i).
	
	Next we prove (ii). We have $\BC\in \mathcal{P}_{\leq n-2}$. Take any coarser $\Dbf\in \CC$. By the triangle inequality and arguing similar to (i), we have
	$$\hat{F}_{V,\Dbf}(\mathcal{T}) \leq C\hat{F}_{V,\BC}(\mathcal{T}) + C\nu(\BC,\Dbf).$$
	Hypothesis (A4) means that we can absorb the term $C\hat{F}_{V,\BC}$ into the left-hand side, and deduce $\hat{F}_{V,\Dbf}(\mathcal{T})\leq C\nu(\BC,\Dbf)$. This implies that $\mathscr{F}^*_{V,\BC}(\mathcal{T})\leq C\nu(\BC,\Dbf)$ for every coarser $\Dbf\in \CC$. Now, let $\w_\BC\in P_0\cap S(\BC)^\perp$ be a direction vector realising the infimum in $a_*(\BC)$. Also, let $P_\alpha^\prime$ be a plane minimising $\nu(P^\alpha,L)$ over all coarser planes $L\subset P_0$ with $S(\BC)\subset L$ and $\w_{\BC}\in L$. Now, with $\Dbf:= |P_1^\prime| + |P_2^\prime|$ we have $\nu(\BC,\Dbf) \leq Ca_*(\BC)$, with $C = C(n,k)$ (one can argue by contradiction, for example, to see this). Combining this with the above establishes (ii).
	
	To prove (iii), we start with the inequality
	$$\hat{E}_{V,P}(\mathcal{T}) \leq 2\hat{E}_{V,\BC}(\mathcal{T}) + C\nu(\BC,P)$$
	for any plane $P$ with $P\supset S(\BC)$. We can use (A4) to absorb the excess on the right-hand side to deduce that $\mathscr{E}_{V,\BC}(\mathcal{T}) \leq C\nu(\BC,P)$ for all such planes $P$. If $\BC\in \mathcal{P}$, let $P$ be the plane defined by the graph of $\frac{1}{2}(p^1+p^2)$: then $\nu(\BC,P)\leq Ca^*(\BC)$ and $P\supset S(\BC)$, and so the above gives $\mathscr{E}_{V,\BC}(\mathcal{T})\leq Ca^*(\BC)$. If instead $\BC\in\CC_{n-1}\setminus\mathcal{P}$, then choose $P$ to be the plane determined by extending half-plane from the graph of $\frac{1}{2}(h^1+h^2)$ or $\frac{1}{2}(h^3+h^4)$, depending on which determines the value of $a^*(\BC)$. This also obeys $\nu(\BC,P)\leq a^*(\BC)$, and so the result follows.
	
	Finally we prove (iv); notice that if $\BC\in \mathcal{P}$, then $a^*(\BC) = a_*(\BC)$ and so the result is immediate from (iii). So, we may assume that $\BC\in \CC_{n-1}\setminus\mathcal{P}$. The basic idea here is that one of the two pairs of half-planes cannot be close to each other whilst the other two are not. We will prove this using a coarse blow-up, similar to \cite[Lemma 9.1]{Wic14}. To avoid repetition with some of the arguments in the proof of Theorem \ref{thm:graphical-rep} when $\BC\in \CC\setminus\mathcal{P}$, we will defer the proof of this case to the proof of Theorem \ref{thm:graphical-rep} when $\sfrak(\BC)=n-1$ in Section \ref{sec:base-case}.
\end{proof}

\subsection{Main $L^2$ Estimates and Graphical Representation}

We now state the three main technical results. The first of the main theorems is a graphical representation of $V$ away from $S(\BC)$. The key point is that under Hypothesis (A4), the set $\mathcal{D}(V)\cap T$ is forced to lie near $S(\BC)$, where
$$\mathcal{D}(V):= \{x\in \spt\|V\|:\Theta_V(x)\geq 2\}.$$

\begin{theorem}[Graphical Representation]\label{thm:graphical-rep}
	Fix $\tau\in (0,1/8)$, $\sigma\in (1/2,1)$, $a\in (0,1/2)$, and $M\in [1,\infty)$. Also fix $\beta\in(0,1)$. Then, there exist $\eps_{I}$ and $\gamma_{I}\in (0,1)$ depending on $n,k,\tau,\sigma,a,M, \beta$, such that the following is true. Suppose that $V\in \mathcal{V}_\beta$ and $\BC\in \CC$ satisfy Hypothesis $G(\eps_I)$ and $H(A,\eps_I,\gamma_I,M)$ in $\mathcal{T}:= \mathcal{T}^A_r(x_0,s_0)$, where either $s_0 = 0$ or $r/s_0\in (a,1-a)$. Then we have:
	\begin{enumerate}
		\item [(i)] If $\sfrak(\BC)=n-1$, we have
		\begin{align*}
			V\res (\R^k\times\pi_{P_0}(\mathcal{T}^A_{\sigma r}(x_0,s_0)))&\cap \{|x^n|>\tau r\}\\
			& = \sum^4_{\alpha=1}\left|\graph(h^\alpha|_{\{|x^n|>\tau r\}\cap U_0^\alpha} + u^\alpha)\right|\res (\R^k\times \pi_{P_0}(\mathcal{T}^A_{\sigma r}(x_0,s_0)))
		\end{align*}
		where $u^\alpha\in C^\infty(P_0\cap \mathcal{T}^A_{(1+\sigma)r/2}(x_0,s_0)\cap \{|x^n|>\tau r\}\cap U^\alpha_0;P_0^\perp)$. Moreover, $h^\alpha+u^\alpha$ is a smooth solution to the minimal surface system on its domain, and we have $\|u^\alpha\|_{L^2}\leq C\hat{E}_{V,\BC}(\mathcal{T})$ for each $\alpha=1,\dotsc,4$.
		\item [(ii)] If $\sfrak(\BC)<n-1$, we have
		\begin{align*}
			V\res (\R^k\times\pi_{P_0}(\mathcal{T}^A_{\sigma r}(x_0,s_0)))&\cap \{r_{\BC}(x)>\tau r\}\\
			& = \sum^2_{\alpha=1}\left|\graph(p^\alpha|_{\{r_{\BC}>\tau r\}} + u^\alpha)\right|\res (\R^k\times \pi_{P_0}(\mathcal{T}^A_{\sigma r}(x_0,s_0)))
		\end{align*}
		where $u^\alpha\in C^\infty(P_0\cap \mathcal{T}^A_{(1+\sigma)r/2}(x_0,s_0)\cap \{r_{\BC}\geq\tau r\};P_0^\perp)$ for $\alpha=1,2$. Moreover, $p^\alpha+u^\alpha$ is a smooth solution to the minimal surface system on its domain, and we have $\|u^\alpha\|_{L^2}\leq C\hat{E}_{V,\BC}(\mathcal{T})$ for each $\alpha=1,2$.
	\end{enumerate}
	Here, $C = C(n,k,\sigma,\tau)$.
\end{theorem}
We will regularly refer back to the functions $u^\alpha$ given in Theorem \ref{thm:graphical-rep}, often saying things such as \emph{``Let $\{u^\alpha\}_{\alpha}$ be the functions which represent $V$ graphically over $\BC$ in the sense of Theorem \ref{thm:graphical-rep}''}. Whether the index $\alpha$ ranges over $\{1,2\}$ or $\{1,2,3,4\}$ will be clear from context, and when convenient we will take the liberty to switch between the two index sets when $\BC\in \mathcal{P}_{n-1}$. 

\begin{remark}\label{remark:after-graphical-rep}
Notice that the conclusion of Theorem \ref{thm:graphical-rep} in either case implies that the two-sided excess of $V$ relative to $\BC$ and the one-sided excess are comparable, i.e.
$$\hat{F}_{V,\BC}(\mathcal{T}) \leq C\hat{E}_{V,\BC}(\mathcal{T})$$
for some $C = C(n,k)\in (0,\infty)$.
\end{remark}

The content of the second main result is certain key $L^2$ estimates which hold at points of `good' density, i.e.~points in $\mathcal{D}(V)$.

\begin{theorem}[Main $L^2$ Estimates]\label{thm:L2-estimates}
	Fix $\sigma\in (0,1)$, $a\in (0,1/2)$, $M\in [1,\infty)$, $\mu\in (0,1)$, and $\rho\in (0,(1-\sigma)/2]$. Also fix $\beta\in (0,1)$. Then, there exist $\eps_{II}$ and $\gamma_{II}\in (0,1)$ depending on $n,k,\sigma,a,M,\mu,\rho,\beta$ such that the following is true. Suppose that $V\in\mathcal{V}_\beta$ and $\BC\in \CC$ satisfy Hypothesis $G(\eps_{II})$ and $H(A,\eps_{II},\gamma_{II},M)$ in $\mathcal{T}:= \mathcal{T}^A_r(x_0,s_0)$, where either $s_0 = 0$ or $r/s_0\in (a,1-a)$. Then, for any $z\in \mathcal{D}(V)\cap \mathcal{T}^A_{\sigma r}(x_0,s_0)$ we have
	\begin{equation}\label{E:L2-1}
		(\sigma\rho r)^{-\mu}\int_{B_{\sigma\rho r}(z)}\frac{\dist^2(x,(\tau_z)_\#\BC)}{|x-z|^{n+2-\mu}}\, \ext\|V\|(x) \leq C^\prime\hat{E}_{V,(\tau_z)_\#\BC}^2(B_{\rho r}(z)).
	\end{equation}
	\begin{equation}\label{E:L2-2}
		(\sigma\rho r)^{-n}\int_{B_{\sigma\rho r}(z)}\sum^{\sfrak(\BC)}_{j=1}|e_j^{\perp_{T_xV}}|^2\, \ext\|V\|(x) \leq C\hat{E}_{V,(\tau_z)_\#\BC}(B_{\rho r}(z)).
	\end{equation}
	Furthermore in the case $\sfrak(\BC)=n-1$, if we fix $\tau\in (0,1/8)$ and allow $\eps_{II}$, $\gamma_{II}$ to additionally depend on $\tau$, then when Theorem \ref{thm:graphical-rep}(i) holds we have
	\begin{align}
		\nonumber\sum_{\alpha=1}^4\int_{B^n_{\sigma\rho r}(\bar{z})\cap \{|x^n|\geq\tau\}\cap U_0^\alpha}R_{\bar{z}}^{2-n}\left|\frac{\del}{\del R_{\bar{z}}}\left(\frac{u^\alpha(x)-(\xi^{\perp_{P_0}}-Dh^\alpha\cdot\xi^{\top_{P_0}})}{R_{\bar{z}}(x)}\right)\right|^2&\, \ext x\\
		&\hspace{-2em} \leq C\hat{E}^2_{V,(\tau_z)_\#\BC}(B_{\rho r}(z)). \label{E:L2-3}
	\end{align}
	Here, $\bar{z}:= \pi_{P_0}(z)$, $R_{\bar{z}}(x):= |x-\bar{z}|$, $\xi := z^{\perp_{S(\BC)}}$, and $C^\prime = C^\prime(n,k,\sigma,a,M,\mu)\in (0,\infty)$ and $C = C(n,k,\sigma,a,M)\in (0,\infty)$.
\end{theorem}

\textbf{Remark:} In fact, if one wishes to prove the estimates \eqref{E:L2-2} and \eqref{E:L2-3}, one does not need to assume the choice of $\eps_{II}$ or $\gamma_{II}$ depends on $\mu$.

The third and final main result is an estimate which controls the distance between $\BC$ and $(\tau_z)_\#\BC$. Ultimately, it enables us to replace $\hat{E}_{V,(\tau_z)_\#\BC}(B_{\rho r}(z))$ by $\hat{E}_{V,\BC}(\mathcal{T})$ on the right-hand side of the estimates in Theorem \ref{thm:L2-estimates}.

\begin{theorem}[Shift Estimate]\label{thm:shift-estimate}
	Fix $\sigma\in (0,1)$, $a\in (0,1/2)$, and $M\in [1,\infty)$. Also fix $\beta\in (0,1)$. Then there exist $\eps_{III}$ and $\gamma_{III}\in (0,1)$, depending on $n,k,\sigma,a,M,\beta$, such that the following is true. Suppose that $V\in\mathcal{V}_\beta$ and $\BC\in \CC$ satisfy Hypothesis $G(\eps_{III})$ and $H(A,\eps_{III},\gamma_{III},M)$ in $\mathcal{T}:=\mathcal{T}^A_r(x_0,s_0)$, where either $s_0=0$ or $r/s_0\in (a,1-a)$. Then, for any $z\in\mathcal{D}(V)\cap \mathcal{T}^A_{\sigma r}(x_0,s_0)$ we have
	\begin{equation}\label{E:shift}
		r^{-1}\left(|\xi^{\perp_{P_0}}| + \nu(\BC)|\xi^{\top_{P_0}}|\right) \leq C\hat{E}_{V,\BC}(\mathcal{T})
	\end{equation}
	where $\xi := z^{\perp_{S(\BC)}}$ and $C = C(n,k,\sigma,a,M)\in (0,\infty)$.
\end{theorem}

We will prove the above three theorems by induction on $\sfrak(\BC)\in \{n-1,n-2,\dotsc,1\}$. Along the way we will need a useful corollary, which we record and prove now. It provides, for example, certain $L^2$ non-concentration estimates.

\begin{corollary}\label{cor:L2-corollary}
	Fix $\sfrak\in \{n-1,n-2,\dotsc,1\}$ and suppose Theorems \ref{thm:graphical-rep}, \ref{thm:L2-estimates}, \ref{thm:shift-estimate} all hold when we additionally impose that $\sfrak(\BC)=\sfrak$. Fix $\sigma\in (0,1)$, $a\in (0,1/2)$, $M\in [1,\infty)$, and $\beta\in (0,1)$. Then, there exist constants $\eps_0$ and $\gamma_0\in (0,1)$, depending on $n,k,\sfrak,\sigma,a,M,\beta$, such that the following is true. Suppose that $V\in \mathcal{V}_\beta$ and $\BC\in \CC$ satisfy Hypothesis $G(\eps_0)$ and $H(A,\eps_0,\gamma_0,M)$ in $\mathcal{T}:= \mathcal{T}^A_r(x_0,s_0)$, where $\sfrak(\BC)=\sfrak$ and either $s_0=0$ or $r/s_0\in (a,1-a)$. Then for every $z\in\mathcal{D}(V)\cap \mathcal{T}^A_{\sigma r}(x_0,s_0)$:
	\begin{equation}\label{E:cor-1}
		\dist(z,S(\BC))\leq Cr\frac{\hat{F}_{V,\BC}(\mathcal{T})}{\mathscr{F}^*_{V,\BC}(\mathcal{T})};
	\end{equation}
	\begin{equation}\label{E:cor-2}
		\sup_{x\in \mathcal{T}^A_r(x_0,s_0)}|\dist(x,\BC)-\dist(x,(\tau_z)_\#\BC)| \leq \nu(\BC,(\tau_z)_\#\BC) \leq Cr\hat{E}_{V,\BC}(\mathcal{T}).
	\end{equation}
	If in addition we fix $\eta\in (0,\sigma(1-\sigma)/4)$ and allow $\eps_0$, $\delta_0$, and $\gamma_0$ to be small also depending on $\eta$, then
	\begin{equation}\label{E:cor-3}
		\int_{(\mathcal{D}(V)\cap \mathcal{T}^A_{\sigma r}(x_0,s_0))_{\eta r}}\dist^2(x,\BC)\, \ext\|V\|(x) \leq C\eta^{1/2}\int_{\mathcal{T}}\dist^2(x,\BC)\, \ext\|V\|(x);
	\end{equation}
	\begin{equation}\label{E:cor-4}
		\int_{(\mathcal{D}(V)\cap \mathcal{T}^A_{\sigma r}(x_0,s_0))_{\eta r}}\sum^{\sfrak}_{j=1}|e_j^{\perp_{T_xV}}|^2\, \ext\|V\|(x) \leq C(\eta r)^{-2}\int_{\mathcal{T}}\dist^2(x,\BC)\, \ext\|V\|(x).
	\end{equation}
	Here, $C = C(n,k,\sigma,a,M)\in (0,\infty)$.
\end{corollary}

\textbf{Remark:} A key observation is that when Corollary \ref{cor:L2-corollary} holds, in light of \eqref{E:cor-2} we can replace the $\hat{E}_{V,(\tau_z)_\#\BC}(B_\rho(z))$ terms appearing on the right-hand side of the estimates in Theorem \ref{thm:L2-estimates} by the `true' fine excess $\hat{E}_{V,\BC}(\mathcal{T})$.

To conclude this subsection, we will prove Corollary \ref{cor:L2-corollary}. In the next subsection we will begin the proof of the three main theorems.

\begin{proof}[Proof of Corollary \ref{cor:L2-corollary}]
	Let $\eps_0$ and $\gamma_0$ be small enough so that the conclusions of Theorem \ref{thm:graphical-rep}, \ref{thm:L2-estimates}, and \ref{thm:shift-estimate} all hold for $V$ and $\BC$ in the region $\mathcal{T}$ with this choice of $\sigma,a,M,\beta$ (we will be making certain choices of $\mu$ and $\rho$ in Theorem \ref{thm:L2-estimates} in our proof, so for Theorem \ref{thm:L2-estimates} we need $\eps_0,\gamma_0$ to be smaller than the minimum of the corresponding constants over these choices of $\mu$ and $\rho$).
	
	Now, for $z\in \mathcal{D}(V)\cap \mathcal{T}^A_{\sigma r}(x_0,s_0)$, writing $\xi:= z^{\perp_{S(\BC)}}$, we have
	\begin{align*}
		r^{-2}\dist^2(z,S(\BC)) & = r^{-2}(|\xi^{\top_{P_0}}|^2 + |\xi^{\perp_{P_0}}|^2)\\
		& \leq r^{-2}\hat{E}_V^{-2}(\mathcal{T})(|\xi^{\top_{P_0}}|^2 + \hat{E}_V^2(\mathcal{T})|\xi^{\perp_{P_0}}|^2)\\
		& \leq Cr^{-2}\hat{E}_V^{-2}(\mathcal{T})(|\xi^{\top_{P_0}}|^2 + \nu(\BC)^2|\xi^{\perp_{P_0}}|^2)\\
		& \leq C\hat{E}_V^{-2}(\mathcal{T})\hat{E}^2_{V,\BC}(\mathcal{T})\\
		& \leq C\mathscr{F}^*_{V,\BC}(\mathcal{T})^{-2}\hat{F}_{V,\BC}(\mathcal{T})^2
	\end{align*}
	where in the second inequality we have used Lemma \ref{lemma:nu-a}(i), in the third inequality we have used \eqref{E:shift} from Theorem \ref{thm:shift-estimate}, and the fourth inequality follows from the definitions of the excess quantities. This establishes \eqref{E:cor-1}.
	
	If we combine \eqref{E:translates-1} and \eqref{E:translates-2} with \eqref{E:shift}, we have for each $x\in\mathcal{T}^A_r(x_0,s_0)$,
	\begin{align*}
		|\dist(x,\BC)-\dist(x,(\tau_z)_\#\BC)| \leq \nu(\BC,(\tau_z)_\#\BC) \leq C(|\xi^{\top_{P_0}}|^2 + \nu(\BC)|\xi^{\perp_{P_0}}|^2)^{1/2} \leq Cr\hat{E}_{V,\BC}(\mathcal{T})
	\end{align*}
	which proves \eqref{E:cor-2}.
	
	Inequality \eqref{E:cor-3} is a non-concentration estimate similar to \cite[Corollary 3.2(ii)]{Sim93} or \cite[Lemma 10.8(b)]{Wic14}. To prove it, note that for any $z\in \mathcal{D}(V)\cap \mathcal{T}^A_{\sigma r}(x_0,s_0)$ we have
	\begin{align}
		(2\eta r)^{-n+1/4}\int_{B_{2\eta r}(z)}\dist^2(x,\BC)&\, \ext\|V\|(x) \leq \int_{B_{\sigma(1-\sigma)r/2}(z)}\frac{\dist^2(x,\BC)}{|x-z|^{n-1/4}}\, \ext\|V\|(x)\nonumber\\
		& \leq Cr^2\int_{B_{\sigma(1-\sigma)r/2}(z)}\frac{\dist^2(x,(\tau_z)_\#\BC)}{|x-z|^{n+7/4}}\, \ext\|V\|(x) + Cr^{2+1/4}\hat{E}_{V,\BC}^2(\mathcal{T})\nonumber\\
		& \leq Cr^{9/4}\hat{E}_{V,(\tau_z)_\#\BC}^2(B_{(1-\sigma)r/2}(z)) + Cr^{9/4}\hat{E}_{V,\BC}^2(\mathcal{T})\nonumber\\
		& \leq Cr^{9/4}\hat{E}_{V,\BC}^2(\mathcal{T}).\label{E:non-con-1}
	\end{align}
	Here, the first inequality is true for elementary reasons, the second inequality follows from \eqref{E:cor-2} just shown and the fact that $\int_{B_{\sigma(1-\sigma)r/2}(z)}|x-z|^{-n+1/4}\, \ext\|V\|(x) \leq Cr^{1/4}$ for some $C = C(n,k)\in (0,\infty)$, the third inequality follows from \eqref{E:L2-1} with $\rho = (1-\sigma)/2$ and $\mu = 1/4$, and the fourth inequality follows from \eqref{E:cor-2}. Here, $C = C(n,k,\sigma,a,M)$.
	
	Now, note that by \eqref{E:cor-1}, provided we choose $\gamma_0$ sufficiently small (depending on $\eta$) we can guarantee that for any $z\in \mathcal{D}(V)\cap \mathcal{T}^A_{\sigma r}(x_0,s_0)$ we have
	$$\dist(z,S(\BC))\leq \eta r/4.$$
	In particular, $\mathcal{D}(V)\cap \mathcal{T}^A_{\sigma r}(x_0,s_0)\subset (S(\BC))_{\eta r/4}$. (In fact, we only need the neighbourhood over region $S(\BC)\cap B_{\sqrt{n}\sigma r}(x_0)$.) Now cover $(S(\BC)\cap B_{\sqrt{n}\sigma r}(x_0))_{\eta r/4}$ by a collection of balls $\{B_{\eta r/2}(y_j)\}_{j=1}^N$ with $y_j\in S(\BC)$ and $N\leq C(\eta r/\sigma r)^{-\sfrak(\BC)} \leq C(\eta/\sigma)^{-n+1} \leq C\eta^{-n+1}$ (as $S(\BC)$ is a subspace of dimension at most $n-1$, and $\sigma\in (0,1)$); here $C = C(n)$. In particular, $\{B_{\eta r/2}(y_j)\}_{j=1}^N$ covers $\mathcal{D}(V)\cap \mathcal{T}^A_{\sigma r}(x_0,s_0)$. For each $j$, choose $z_j\in \mathcal{D}(V)\cap \mathcal{T}^A_{\sigma r}(x_0,s_0)\cap B_{\eta r/2}(y_j)$ (if it exists, otherwise we may discard this ball from the collection) and consider the ball $B_{\eta r}(z_j)$: note that $B_{\eta r}(z_j)\supset B_{\eta r/2}(y_j)$, and thus the collection of balls $\{B_{\eta r}(z_j)\}_{j=1}^N$ also covers $\mathcal{D}(V)\cap \mathcal{T}^A_{\sigma r}(x_0,s_0)$ (here we are abusing notation by reusing $N$ to denote the number of balls in this cover even though some balls may have been discarded).
	
	To summarise, we have found a collection of balls $\{B_{\eta r}(z_j)\}_{j=1}^N$ which covers $\mathcal{D}(V)\cap \mathcal{T}^A_{\sigma r}(x_0,s_0)$, $z_j\in\mathcal{D}(V)\cap \mathcal{T}^A_{\sigma r}(x_0,s_0)$ for all $j$, and obeys $N\leq C\eta^{-n+1}$. Hence, $\{B_{2\eta r}(z_j)\}_{j=1}^N$ covers $(\mathcal{D}(V)\cap \mathcal{T}^A_{\sigma r}(x_0,s_0))_{\eta r}$. So, using \eqref{E:non-con-1}, we get
	\begin{align*}
		\int_{(\mathcal{D}(V)\cap \mathcal{T}^A_{\sigma r}(x_0,s_0))_{\eta r}}\dist^2(x,\BC)\, \ext\|V\|(x) & \leq \sum^N_{j=1}\int_{B_{2\eta r}(z_j)}\dist^2(x,\BC)\, \ext\|V\|(x)\\
		& \leq N\cdot C(\eta r)^{n-1/4}\cdot r^{9/4}\hat{E}^2_{V,\BC}(\mathcal{T})\\
		& \leq C\eta^{3/4}\cdot r^{n+2}\hat{E}^2_{V,\BC}(\mathcal{T})\\
		& = C\eta^{1/2}\int_{\mathcal{T}}\dist^2(x,\BC)\, \ext\|V\|(x)
	\end{align*}
	which is exactly the claimed inequality in \eqref{E:cor-3} (as $\eta^{3/4}\leq \eta^{1/2}$).

	The proof of \eqref{E:cor-4} is similar to that of \eqref{E:cor-3}. First note that for $z\in \mathcal{D}(V)\cap \mathcal{T}^A_{\sigma r}(x_0,s_0)$,
	\begin{align}
		(\eta r)^2\int_{B_{2\eta r}(z)}\sum^{\sfrak}_{j=1}|e_j^{\perp_{T_xV}}|^2&\, \ext\|V\|(x) \leq C\int_{B_{4\eta r}(z)}\dist^2(x,(\tau_z)_\#\BC)\, \ext\|V\|(x)\nonumber\\
		& \leq C\int_{B_{4\eta r}(z)}\dist^2(x,\BC)\, \ext\|V\|(x) + C\|V\|(B_{4\eta r}(z))\cdot r^2\hat{E}_{V,\BC}^2(\mathcal{T}).\label{E:non-con-2}
	\end{align}
	The first inequality follows from \eqref{E:L2-2} with $\rho = 4\eta$ and $\sigma=1/2$ therein, and the second follows from \eqref{E:cor-2}; here $C = C(n,k,\sigma,a,M)$.
	
	Now cover $(\mathcal{D}(V)\cap \mathcal{T}^A_{\sigma r}(x_0,s_0))_{\eta r}$ with balls $\{B_{2\eta r}(z_j)\}_{j=1}^N$ centred at points $z_j\in \mathcal{D}(V)\cap\mathcal{T}^A_{\sigma r}(x_0,s_0)$ in such a way that $\{B_{4\eta r}(z_j)\}_{j=1}^N$ can be decomposed into at most $M\leq C(n,k)<\infty$ pairwise disjoint sub-collections. (This can be done arguing similarly to that in the proof of \eqref{E:cor-3}, except when choosing the covering centred on $S(\BC)$, we choose the cover $\{B_{\eta r/2}(y_j)\}_{j}$ such that the balls $\{B_{6\eta r}(y_j)\}_{j}$ have the property that they can be decomposed into at most $M\leq C(n,k)<\infty$ pairwise disjoint sub-collections, and then this property will be inherited by $\{B_{4\eta r}(z_j)\}_j$ by inclusion). Using this cover and the fact that $\|V\|(\mathcal{T})\leq Cr^n$ for suitable $C = C(n,k)$, we get
	\begin{align*}
		\int_{(\mathcal{D}(V)\cap \mathcal{T}^A_{\sigma r}(x_0,s_0))_{\eta r}}\sum^{\sfrak}_{j=1}|e_j^{\perp_{T_xV}}|^2&\, \ext\|V\|(x) \leq \sum_{j=1}^N\int_{B_{2\eta r}(z_j)}\sum^{\sfrak}_{j=1}|e_j^{\perp_{T_xV}}|^2\, \ext\|V\|(x)\\
		& \leq C(\eta r)^{-2}\cdot M\cdot\left[\int_{\mathcal{T}}\dist^2(x,\BC)\, \ext\|V\|(x) + r^{n+2}\hat{E}_{V,\BC}^2(\mathcal{T})\right]\\
		& \leq C(\eta r)^{-2}\int_{\mathcal{T}}\dist^2(x,\BC)\, \ext\|V\|(x).
	\end{align*}
	Here, the second inequality follows from summing \eqref{E:non-con-2} over the $M$ pairwise disjoint sub-collections of balls, as $B_{4\eta r}(z_j)\subset \mathcal{T}$ for each $j$. This therefore completes the proof of \eqref{E:cor-4} and hence of the corollary.
\end{proof}

\subsection{Structure of the Proofs of Theorem \ref{thm:graphical-rep}, \ref{thm:L2-estimates}, \ref{thm:shift-estimate}}

Theorem's \ref{thm:graphical-rep}, \ref{thm:L2-estimates}, and \ref{thm:shift-estimate} will be proved by induction on the spine dimension $\sfrak(\BC):=\dim(S(\BC))$. We stress now that the induction will be \emph{downwards}, i.e.~we use the results for larger spine dimensions to prove them for lower spine dimensions. Intuitively, this is because $\{\BC:\sfrak(\BC)=\sfrak\}$ is not necessarily a compact class, but $\{\BC:\sfrak(\BC)\geq \sfrak\}$ is, i.e.~cones with lower dimensional spines can limit onto those with higher dimensional spines.

The induction proof has three steps:
\begin{enumerate}
	\item [Step 1.] Base case: we will prove all three theorems when $\sfrak(\BC)=n-1$; this will be done in Section \ref{sec:base-case}. Thereafter, we can fix $\sfrak^*\in \{0,1,\dotsc,n-2\}$ and assume that the theorems all hold for $\sfrak(\BC)>\sfrak^*$.
	\item [Step 2.] To prove the induction step, we first need to construct the class of \emph{fine blow-ups} of suitable sequences of stationary integral varifolds $(V_j)_j$ relative to a sequence of cones $(\BC_j)_j$ with $\sfrak(\BC_j)\equiv \sfrak>\sfrak^*$. This construction will be done in Section \ref{sec:fine-blow-ups}, and will require the validity of all three of the theorems for spine dimension $>\sfrak^*$ (which at this stage is known).
	\item [Step 3.] We then complete the induction step, showing the validity of the three theorems for $\sfrak^*$ assuming their validity for all $\sfrak>\sfrak^*$. This completes the proof of the three theorems by induction. The inductive proof of Theorem \ref{thm:graphical-rep} is done in Section \ref{sec:graphical-rep}, whilst the inductive proof of Theorem \ref{thm:L2-estimates} and Theorem \ref{thm:shift-estimate} is done in Section \ref{sec:L2-shift}.
\end{enumerate}

We will break the proof of Theorem \ref{thm:L2-estimates} down into several smaller lemmas, as stated below. Indeed, Theorem \ref{thm:L2-estimates} will follow from a suitable application of Lemma \ref{lemma:L2-estimates-2} (stated below) to $(\eta_{z,\rho})_\#V$ in place of $V$.

\begin{lemma}\label{lemma:L2-estimates-2}
	Fix $\sigma\in (0,1)$ and $\beta\in (0,1)$. Then, there exists $\eps_0 = \eps_0(n,k,\sigma,\beta)\in (0,1)$ such that the following is true. Suppose we have $V\in \mathcal{V}_\beta$ and $\BC\in \CC$ aligned which satisfy:
	\begin{enumerate}
		\item [(i)] $\Theta_V(0)\geq 2$ and $(\w_n 2^n)^{-1}\|V\|(B_2(0))\leq 5/2$;
		\item [(ii)] $\nu(\BC)<\eps_0$ and $\hat{E}_{V,\BC}<\eps_0$.
	\end{enumerate}
	Then we have:
	\begin{equation}\label{E:initial-1}
		\int_{B_\sigma(0)}\frac{|x^{\perp_{T_xV}}|^2}{|x|^{n+2}}\, \ext\|V\|(x) \leq C\hat{E}_{V,\BC}^2.
	\end{equation}
	\begin{equation}\label{E:initial-2}
		\int_{B_\sigma(0)}\sum^{\sfrak(\BC)}_{j=1}|e_j^{\perp_{T_xV}}|^2\, \ext\|V\|(x) \leq C\hat{E}_{V,\BC}^2.
	\end{equation}
	Moreover, if $\mu\in (0,1)$ and we allow $\eps_0$ to also depend on $\mu$, then
	\begin{equation}\label{E:initial-3}
		\int_{B_\sigma(0)}\frac{\dist^2(x,\BC)}{|x|^{n+2-\mu}}\, \ext\|V\|(x) \leq C^\prime\hat{E}_{V,\BC}^2.
	\end{equation}
	Here, $C = C(n,k,\sigma,\beta)\in (0,\infty)$ and $C^\prime = C^\prime(n,k,\sigma,\mu,\beta)\in (0,\infty)$.
\end{lemma}

In turn, the proof of Lemma \ref{lemma:L2-estimates-2} relies on the following more technical lemma:

\begin{lemma}\label{lemma:Z}
	Fix $\sigma\in (0,1)$, $a\in (0,1/2)$. Also fix $\beta\in (0,1)$. Then there exist $\eps_0 = \eps_0(n,k,\sigma,a,\beta)\in (0,1)$ such that the following is true. Suppose that $V\in \mathcal{V}_\beta$ obeys
	$$(\w_n2^n)^{-1}\|V\|(B_2(0))\leq 5/2 \qquad \text{and} \qquad \hat{E}_V<\eps_0,$$
	$\BC\in \CC$ is aligned with $\nu(\BC)<\eps_0$, and write $\mathcal{T}:=\mathcal{T}_r^{S(\BC)}(x_0,s_0)$, where $s_0>0$ and $r/s_0\in (a,1-a)$. Then we have:
	\begin{equation}\label{E:Z-a}
		\int_{\mathcal{T}^{S(\BC)}_{\sigma r}(x_0,s_0)}|x^{\perp_{S(\BC)}\perp_{T_xV}}|^2\, \ext\|V\|(x)\leq C\int_{\mathcal{T}}\dist^2(x,\BC)\, \ext\|V\|(x).
	\end{equation}
	Moreover, for any function $\phi_{(0)}\in C^\infty_c(\mathcal{T}^{S(\BC)}_{\sigma r}(x_0,s_0))$ depending only on $x^{\top_{S(\BC)}}$ and $|x^{\perp_{S(\BC)}}|$ and satisfying:
	\begin{itemize}
		\item $|\phi_{(0)}|\leq \tilde{C}$ and $|D^\alpha\phi_{(0)}|\leq \tilde{C}r^{-|\alpha|}$ for each multi-index $\alpha$ with $|\alpha|\leq n-s(\BC)$, for some constant $\tilde{C}$ (we stress here $r$ is the length-scale of the domain $\mathcal{T}$),
		\item $|\phi_{(0)}|$ is a smooth function,
	\end{itemize}
	we have
	\begin{equation}\label{E:Z-b}
		\left|\int\phi_{(0)}(x)|x^{\perp_{S(\BC)}}|^2\, \ext\|V\| - \int\phi_{(0)}(x)|x^{\perp_{S(\BC)}}|^2\, \ext\|\BC\|\right| \leq C\int_{\mathcal{T}}\dist^2(x,\BC)\, \ext\|V\|(x).
	\end{equation}
	Here, $C = C(n,k,\sigma,a,\beta,\tilde{C})\in (0,\infty)$.
\end{lemma}

\textbf{Note:} We stress that the notation $x^{\perp_A\perp_B}$ is shorthand for $(x^{\perp_A})^{\perp_B}$.

The proof of both the above lemmas will also be by induction and fits into our previously described inductive proof. Indeed, for both lemmas the $\sfrak(\BC)=n-1$ case is established in Section \ref{sec:base-case}, whilst the $\sfrak(\BC)<n-1$ case is proved (inductively) in Section \ref{sec:L2-shift}.

\section{Base Case of Induction for Main $L^2$ Estimates}\label{sec:base-case}

In this section we will prove all three of the main theorems, as well as Lemma \ref{lemma:L2-estimates-2} and Lemma \ref{lemma:Z}, in the base case $\sfrak(\BC)=n-1$. The order will be as indicated before, namely: Theorem \ref{thm:graphical-rep}, Lemma \ref{lemma:Z}, Lemma \ref{lemma:L2-estimates-2}, Theorem \ref{thm:L2-estimates}, and finally Theorem \ref{thm:shift-estimate}. We therefore start with Theorem \ref{thm:graphical-rep}.

\begin{proof}[Proof of Theorem \ref{thm:graphical-rep} when $\sfrak(\BC)=n-1$.]
Fix $\tau\in (0,1/8)$, $\sigma\in (1/2,1)$, $a\in (0,1/2)$, and $M\in [1,\infty)$. Also fix $\beta\in (0,1)$. We will begin by supposing that we have sequences of the various objects and parameters in the statement of the theorem and, without loss of generality, we will perform some basic manipulations to simplify the set-up.

So, suppose we have sequences of numbers $\eps_j\downarrow 0$, $\gamma_j\downarrow 0$, as well as sequences of regions $\mathcal{T}^{A_j}_{r_j}(x_j,s_j)$ with either $s_j=0$ or $r_j/s_j\in (a,1-a)$, varifolds $V_j\in \mathcal{V}_\beta$, cones $\BC_j\in \CC$ with $\sfrak(\BC_j) \equiv n-1$ for all $j$, such that for each $j$, $V_j$ and $\BC_j$ satisfy Hypothesis $G(\eps_j)$ as well as $H(A_j,\eps_j,\gamma_j,M)$ in $\mathcal{T}^{A_j}_{r_j}(x_j,s_j)$. To prove the theorem, we simply need to show that Theorem \ref{thm:graphical-rep} (namely, Theorem \ref{thm:graphical-rep}(i)) holds along some subsequence (this would complete the proof by a contradiction argument).

Without changing notation, we first make some elementary observations to simplify the situation.
\begin{itemize}
	\item We may pass to a subsequence along which $\dim(A_j)$ is constant.
	\item Since $\sfrak(\BC_j)\equiv n-1$, Hypothesis (A1) of $H(A_j,\eps_j,\delta_j,M)$ gives that $S(\BC_j)\equiv S_0$, so is independent of $j$.
	\item By performing an ambient rotation under which $S_0$ and $P_0$ are invariant, we may assume without loss of generality that $A_j\equiv A$ for some fixed subspace $A\subset S_0$.
	\item We may also apply the homothety $\eta_{x_j,r_j}$ for each $j$ to push-forward the situation: since $x_j\in A\subset S_0$, this transformation fixes $A$, $S_0$, and $\BC_j$, and transforms the region $\mathcal{T}^{A_j}_{r_j}(x_j,s_j)$ to
	$$\mathcal{T}_j:= \eta_{x_j,r_j}(\mathcal{T}^A_{r_j}(x_j,s_j)) \equiv \mathcal{T}^A_1(0,s_j^\prime)$$
	where $s^\prime_j:= s_j/r_j$ is either $0$ or belongs to $(1/(1-a),1/a)$. We may then pass to a subsequence to assume that either $s_j^\prime\equiv 0$ for all $j$, or $s_j^\prime\to s_0\in [1/(1-a),1/a]$.
	\item Notice that under this homothety (and rotation), the push-forward of $V_j$ still satisfies Hypothesis $H(A,\eps_j,\gamma_j,M)$ in the region $\mathcal{T}_j$. Indeed, Hypotheses (A1) -- (A4) are scale-invariant and so are satisfied. Moreover, the push-forward of $V_j$ still lies in $\mathcal{V}_\beta$ as again the defining conditions of $\mathcal{V}_\beta$ are closed under homothety as well as rotations fixing $P_0$.
	\item After the homothety is applied, whilst we cannot maintain the condition given by Hypothesis $G(\eps_j)$, we can still maintain the upper bound of $5/2$ on the mass. Indeed, this is because before applying the homothety, the $V_j$ are converging to a multiplicity $2$ plane at scale $1$. Thus, to show $(\w_n (2r_j)^n)^{-1}\|V_j\|(B_{2r_j}(x_j))\leq 5/2$, we can use the monotonicity formula with the varifold convergence to a multiplicity $2$ plane.
	\item We may pass to a subsequence to assume that either $\BC_j\in\mathcal{P}$ for all $j$, or $\BC_j\in \CC\setminus\mathcal{P} \equiv \CC_{n-1}\setminus \mathcal{P}_{n-1}$ for all $j$.
\end{itemize}
Thus, we are now assuming that $V_j$ and $\BC_j$ satisfy Hypothesis $H(A,\eps_j,\gamma_j,M)$ in $\mathcal{T}_j \equiv \mathcal{T}^A_1(0,s_j^\prime)$. We then have two cases, depending on whether $\BC_j\in\mathcal{P}$ for infinitely many $j$ or not. Both cases are however extremely similar, so let us suppose first that $\BC_j\in \mathcal{P}$ for infinitely many $j$; recall that $p_j^1$, $p_j^2$ are then the linear functions representing the planes in $\BC_j$.

Write $\hat{E}_j:= \hat{E}_{V_j}(\mathcal{T}_j)$. First note that by Lemma \ref{lemma:nu-a}(i) and Remark \ref{remark:elementary-observations}(1) from Section \ref{sec:further-notation}, the sequences $(\hat{E}_{j}^{-1}p_j^\alpha)_{j=1}^\infty$ for $\alpha=1,2$ are uniformly bounded (on compact sets). As these are linear functions, this means that they have limits $p^\alpha_*$ as $j\to\infty$ (up to passing to a subsequence). Moreover, from (A3) of Hypothesis $H(A,\eps_j,\gamma_j,M)$ and Lemma \ref{lemma:nu-a}(iii) tell us that $p^1_*\neq p^2_*$, and so if
$$L:= \graph(p^1_*)\cup \graph(p^2_*)$$
then $L$ is the union of two transversely intersecting planes (indeed, they intersect along $S_0$, and so in fact $L\in \mathcal{P}_{n-1}$). In particular, the only coincidence points of $L$ lie on $S_0$.

Now let $v$ denote a coarse blow-up of $(V_j)_j$ off $P_0$ \emph{in the region} $\mathcal{T}^A_1(0,s_0)$; we can construct such an object due to Hypothesis $H(A,\eps_j,\gamma_j,M)$ holding here. To stress, by this we mean performing an analogous procedure to constructing coarse blow-ups off a plane but now only in the region $P_0\cap \mathcal{T}^A_1(0,s_0)$. Thus, all the graphical approximations are over increasing sub-regions, and we are controlling the $W^{1,2}$ norm on these subregions by the $L^2$ height excess $\hat{E}_j$, and so we construct blow-ups by renormalising by $\hat{E}_j\equiv \hat{E}_{V_j}(\mathcal{T}_j)$. These coarse blow-ups will obey many of the properties we saw in Part \ref{part:coarse-blow-ups}, in particular property $(\Bfrak5)$. In particular, since Hypothesis (A4) gives
$$\hat{E}_j^{-1}\hat{F}_{V_j,\BC_j}(\mathcal{T}_j)<\gamma_j\to 0$$
as $j\to\infty$, we deduce that $\graph(v) = L$, and hence $\graph(v)$ is the union of two transversely intersecting planes over the region $\mathcal{T}^A_1(0,s_0)\cap P_0$. (This therefore tells us that $V_j$ must be converging to $P_0$ with multiplicity two in $\mathcal{T}^A_1(0,s_0)$.)

Now, if we were to assume (for the sake of contradiction) that for infinitely many $j$ there is $z_j\in \mathcal{D}(V_j)\cap \{|x^n|\geq\tau/2\}\cap \mathcal{T}^A_{(1+\sigma)/2}(0,s^\prime_j)$, then by passing to a further subsequence so that $z_j\to z_0\in P_0\cap \{|x^n|\geq \tau/2\}\cap \mathcal{T}^A_1(0,s_0)$ and using the Hardt--Simon inequality (cf.~the third point of Remark \ref{remark:after-blow-up-properties} following Theorem \ref{thm:blow-up-properties}) we see that $z_0$ must be a coincidence point of $v$, i.e.~$v^1(z_0) = v^2(z_0)$. However, this would mean that $L$ has a coincidence point outside $S_0$, a contradiction. This contradiction therefore means that for all sufficiently large $j$ we must have
\begin{equation}\label{E:grep-1-1}
\mathcal{D}(V_j)\cap \{|x^n|\geq\tau/2\}\cap \mathcal{T}^A_{(1+\sigma)/2}(0,s_j) = \emptyset.
\end{equation}
If instead $\BC_j\in \CC_{n-1}\setminus \mathcal{P}_{n-1}$ for all $j$, then we can perform essentially an analogous argument to the above (except now the graph of the blow-up will consist of 4 half-planes meeting exactly along $S_0$) to also reach \eqref{E:grep-1-1}. The only additional ingredient we need for this is to know that the blow-up consists of $4$ \emph{distinct} half-planes, i.e.~that no two half-planes collapse onto a single half-plane and so all are multiplicity one (by Lemma \ref{lemma:nu-a}(iii), we know at least one pair does not). This follows from Lemma \ref{lemma:nu-a}(iv), which we will justify now. (In the justification of Lemma \ref{lemma:nu-a}(iv), we have the freedom to change $\tau$ and choose $\sigma = 3/4$, say, so that the constant in Lemma \ref{lemma:nu-a}(iv) only depends on $n,k$.)

If Lemma \ref{lemma:nu-a}(iv) were to fail, then we can further assume that we have $\hat{E}_{j}^{-1}a_*(\BC_j)\to 0$, and thus the blow-up $p_*$ would obey (up to rotation using $(\mathfrak{B}3\text{II})$ and labelling) $p^1_*\neq p^2_*$ are distinct linear functions $L^1,L^2$ over $\mathcal{T}^A_1(0,s_0)\cap P_0\cap \{x^n<0\}$, and $p^3_* \equiv p^4_*\equiv L$ are the same linear function on $\mathcal{T}^A_1(0,s_0)\cap P_0\cap \{x^n>0\}$, and $L^1\equiv L^2\equiv L = 0$ on $\{x^n=0\}$. In particular, the average of $p_*$ (which is harmonic by $(\mathfrak{B}2)$) must obey $(p_*)_a\equiv L$; by applying $(\mathfrak{B}3\text{III})$, we may then assume that $L\equiv 0$, i.e.~$p^1_* \equiv - p^2_*\not\equiv 0$ and $p^3_* \equiv p^4_* \equiv 0$, where $p^1_* = L_*$ is a linear function with $L_* = 0$ on $\{x^n=0\}$.

Since the coarse blow-up is zero on $\{x^n>0\}$, it is easy to check that for any $\tau>0$,
$$\hat{E}_j^{-2}\int_{\R^k\times (\mathcal{T}^A_{15/16}(0,s_0)\cap P_0\cap \{x^n>\tau\})}\dist^2(x,P_0)\, \ext\|V_j\|(x)\to 0$$
and consequently, \eqref{E:rpi} gives
\begin{equation}\label{E:grep-1-2}
\hat{E}_j^{-2}\int_{\R^k\times(\mathcal{T}^A_{7/8}(0,s_0)\cap P_0\cap \{x^n>\tau/2\})}\|\pi_{T_xV_j}-\pi_{P_0}\|^2\, \ext\|V_j\|(x)\to 0
\end{equation}
as $j\to\infty$. The same argument that we were using to justify \eqref{E:grep-1-1} now gives that for all sufficiently large $j$,
$$\Theta_{V_j}(x)<2 \qquad \text{for all }x\in \R^k\times (\mathcal{T}^A_{31/32}(0,s_0)\cap P_0\cap \{x^n<-\tau/8\}).$$
Thus, we may invoke Lemma \ref{lemma:gap-2} to see that for all sufficiently large $j$, we have that in the region $\R^k\times (\mathcal{T}^A_{15/16}(0,s_0)\cap P_0\cap \{x^n<-\tau/4\})$, $V_j$ is the sum of two smooth minimal graphs $u^1_j,u^2_j$ defined on $\mathcal{T}^A_{15/16}(0,s_0)\cap P_0\cap \{x^n<-\tau/4\}$, and which satisfy $\|u_j^i\|_{C^3}\leq C_\tau \hat{E}_j$ for $i=1,2$, where $C_\tau$ depends only on $n,k,\tau$.

We next claim that for all sufficiently large $j$ we must have
\begin{equation}\label{E:grep-1-3}
S_0\cap \mathcal{T}^A_{7/8}(0,s_0)\subset B_\tau(\mathcal{D}(V_j)).
\end{equation}
Indeed, if this were not the case then there would be a subsequence $(j^\prime)$ and a point $y\in S_0\cap\mathcal{T}^A_{7/8}(0,s_0)$ which obeys $\dist(y,\mathcal{D}(V_{j^\prime}))\geq 3\tau/4$ for all $j^\prime$. Subsequently, we could apply Lemma \ref{lemma:gap-2} to represent $V_{j^\prime}$ as minimal graphs in the cylinder $\R^k\times (B_{\tau/2}(y)\cap \mathcal{T}^A_{7/8}(0,s_0))$, which means after taking the coarse blow-up, the coarse blow-up here must be harmonic. But this is a contradiction, because the coarse blow-up needs to agree with $p_*$ here, yet $p^1_* \equiv - p^2_*\not\equiv 0$, yet $p^3_*\equiv p^4_*\equiv 0$ is not even smooth at $y\in S_0$. 

Thus we have now shown \eqref{E:grep-1-3}. We will now proceed to derive the contradiction needed to justify Lemma \ref{lemma:nu-a}(iv) in this case. Take in the first variation formula \eqref{E:stationarity-2} the test function $\widetilde{\zeta}$ for some $\zeta\in C^1_c(B_{1/4}(0)\cap P_0)$ to get
$$\int\nabla^{V_j}x^n\cdot\nabla^{V_j}\widetilde{\zeta}(x)\, \ext\|V_j\|(x) = 0.$$
By performing a computation which is analogous to those seen in later chapters, in particular those leading to \eqref{E:fine-b-10}, and arguing as in \cite[(9.12)]{Wic14} (which utilises the fact that the `bad set' of our Lipschitz approximation is contained with $\{x^n>-\tau/4\}$ as argued above), this gives
\begin{equation}\label{E:grep-1-4}
\int\sum_{\alpha=1,2}\left((J^\alpha_j(x)-1)D_{n}\zeta(x) - \sum^n_{p=1}D_pu_j^\alpha\cdot D_{n}u_j^\alpha D_p\zeta(x) J^\alpha_j(x)\right)\, \ext x = F_j
\end{equation}
where we have used that $\int D_n\zeta = 0$, and where
$$|F_j|\leq C\sup|D\zeta|\int_{\R^k\times (\mathcal{T}^A_{7/8}(0,s_0)\cap P_0\cap \{x^n>-\tau/4\})}\|\pi_{T_xV_j}-\pi_{P_0}\|^2\, \ext\|V_j\|(x).$$
If we split the integral controlling $F_j$ into two parts, one over $\{x^n>\tau/2\}$ and the other over $\{-\tau/4<x^n<\tau/2\}$, we can use \eqref{E:grep-1-2} for the first term and \eqref{E:Lipschitz-approx-6} of Corollary \ref{cor:Lipschitz-approx} with \eqref{E:grep-1-3} for the second term to see that, for all sufficiently large $j$ depending on $\tau>0$,
$$\hat{E}_j^{-2}|F_j| \leq C\sup|D\zeta|\tau^{1/2}.$$
We can then do similarly for the left-hand side of \eqref{E:grep-1-4}: the integral over the `bad set' $\Sigma_j$, is controlled by $C\sup|D\zeta|\int_{\R^k\times (\mathcal{T}^A_{7/8}(0,s_0)\cap P_0\cap \{x^n>-\tau/4\})}\|\pi_{T_xV_j}-\pi_{P_0}\|^2\, \ext\|V_j\|(x)$, which again is controlled by $C\hat{E}_j^2\sup|D\zeta|\tau^{1/2}$ for all $j$ large, and the integral over the `good set' in the region $\{x^n>\tau/4\}$ is controlled in an analogously, using \eqref{E:grep-1-2}. We therefore end up with, for any $\tau>0$ and all sufficiently large $j$ (depending on $\tau$):
\begin{align*}
&\left|\hat{E}_j^{-2}\int_{\{x^n<-\tau/4\}}\sum_{\alpha=1,2}\left((J^\alpha_j(x)-1)D_{n}\zeta(x) - \sum^n_{p=1}D_pu_j^\alpha\cdot D_{n}u_j^\alpha D_p\zeta(x) J^\alpha_j(x)\right)\, \ext x\right|\\
& \hspace{31em} \leq C\sup|D\zeta|\tau^{1/2}.
\end{align*}
On the region $\{x^n<-\tau/4\}$ we have $C^2$ convergence to the blow-up, and so if we take $j\to\infty$ in this we therefore get
$$\int_{x^n<-\tau/4}\sum_{\alpha=1,2}\left(|Dv^\alpha|^2D_n\zeta - 2\sum^n_{p=1}D_pv^\alpha\cdot D_n v^\alpha D_p\zeta\right) = C\sup|D\zeta|\tau^{1/2}.$$
Now letting $\tau\downarrow 0$, and noting that $v^\alpha = p^\alpha_*$ for $\alpha=1,2$, and that $D_ip^\alpha_* = 0$ for $i<n$, this therefore gives
$$\int_{\{x^n<0\}}\sum_{\alpha=1,2}|D_np^\alpha_*|^2D_n\zeta = 0.$$
As $|D_np^\alpha_*|^2$ is constant for $\alpha=1,2$, choosing any $\zeta$ with $\int_{\{x^n<0\}}D_n\zeta \neq 0$ gives that $|D_np^\alpha_*|=0$ for $\alpha=1,2$, i.e.~$p^1_* \equiv p^2_*\equiv 0$, which gives the desired contradiction. Thus, we have now proved Lemma \ref{lemma:nu-a}(iv), and in particular we have now shown \eqref{E:grep-1-1} holds in all cases.

Now we return to the proof of the theorem, given we now have \eqref{E:grep-1-1}. Consider a ball $B^n_{r_*}(x)$, where $r_* = \min\{\tau/4,(1-\sigma)/16\}$, centred at a point $x\in P_0\cap \{|x^n|\geq\tau\}\cap \pi_{P_0}(\mathcal{T}^A_{(1+3\sigma)/4}(0,s_j))$; then $B_{r_*}(x)\subset P_0\cap \{|x^n|\geq \tau/2\}\cap \pi_{P_0}(\mathcal{T}^A_{(1+\sigma)/2}(0,s_j))$, and so $C_{r_*}(x)\equiv \R^k\times B^n_{r_*}(x)\subset\{\Theta_{V_j}<2\}$ for all $j$ large. (One might be concerned about the regions $\mathcal{T}^A_{\sigma}(0,s_0)$ and so forth being curved, but in $\mathcal{T}_j$ we know the $V_j$ are converging to $P_0$, and so $V_j$ is very close to being flat. This with shrinking the radius is sufficient to ensure the correct inclusions on $\spt\|V_j\|$ hold in these regions.) But now we may use the defining property of the class $\mathcal{V}_\beta$ (given through Lemma \ref{lemma:gap-2}) to deduce that in such a region $C_{r_*}(x)$, $V_j\res (C_{r_*/8}(x)\cap \mathcal{T}^A_{(1+3\sigma)/4}(0,s_0))$ is equal to the disjoint union of two smooth minimal graphs over $P_0$. 

Repeating this over finitely many such points $x$ so that these balls cover the desired region, we deduce that $V_j\res \mathcal{T}^A_{(1+7\sigma)/8}(0,s_0)\cap \{|x^n|>\tau\}$ is equal to the disjoint union of four smooth minimal graphs, say $\{\tilde{u}^\alpha_j\}_{\alpha=1}^4$, where $\tilde{u}^1_j$, $\tilde{u}^2_j$ are graphs over a domain contained in $P_0\cap \mathcal{T}^A_{(1+3\sigma)/4}(0,s_0)\cap \{x^n<-\tau/2\}$ and $\tilde{u}^3_j$, $\tilde{u}^4_j$ are graphs over a domain contained in $P_0\cap \mathcal{T}^A_{(1+3\sigma)/4}(0,s_0)\cap \{x^n>\tau/2\}$. Notice that because $\hat{E}_j^{-1}\hat{F}_{V_j,\BC_j}(\mathcal{T}_j)\to 0$, it must be the case that (after relabelling if necessary) that $(\tilde{u}^\alpha_j-p^\alpha_j)/\hat{E}_{V_j}$ converges to zero in $L^2$ on the relevant region for $\alpha=1,2$, and similarly for $(\tilde{u}^{\alpha+2}_j-p^\alpha_j)/\hat{E}_{V_j}$. Then, if we set $u^\alpha:=\tilde{u}^\alpha-h^\alpha$ (where here in the case $\BC_j\in\mathcal{P}_{n-1}$, by $h^\alpha$ we mean $p^\alpha$ restricted to either $\{x^n>0\}$ or $\{x^n<0\}$ depending on the value of $\alpha$) and use Lemma \ref{lemma:nu-a}(iii), we deduce the rest of the specific conclusions of Theorem \ref{thm:graphical-rep} in the case $\sfrak(\BC)=n-1$ by, for example, Allard's regularity theorem applied to these 4 minimal graphs relative to the half-planes in $\BC_j$ they are close to.
\end{proof}

The proof of Theorem \ref{thm:graphical-rep} implies the following weak estimate (compare with the stronger Theorem \ref{thm:shift-estimate}) which will be needed in the proof of Theorem \ref{thm:L2-estimates} and Theorem \ref{thm:shift-estimate}.

\begin{corollary}\label{cor:weak-shift-estimate}
	Let $\eta>0$. Then under the assumptions of Theorem \ref{thm:graphical-rep}, if $\sfrak(\BC)=n-1$ and $\eps_I$, $\gamma_I$ are also allowed to depend on $\eta$, then for every $z\in\mathcal{D}(V)\cap \mathcal{T}^A_{\sigma r}(x_0,s_0)$ we have
	$$r^{-1}\left(|\xi^{\perp_{P_0}}| + \hat{E}_V(\mathcal{T})|\xi^{\top_{P_0}}|\right) < \eta \hat{E}_V(\mathcal{T})$$
	where $\xi = z^{\perp_{S(\BC)}}$.
\end{corollary}

Our next aim is to prove Theorem \ref{thm:L2-estimates} when $q=1$. For this, we will need a good understanding of the behaviour that can occur when the varifold exhibits small \emph{one-sided} fine excess in regions $\mathcal{T}_r^{S(\BC)}(x_0,s_0)$. The kind of arguments used in the proof of Theorem \ref{thm:graphical-rep} are key to this.

\begin{lemma}\label{lemma:tori-q=1}
	Fix $\sigma\in (1/2,1)$ and $a\in (0,1/2)$. Fix also $\beta\in (0,1)$. Then, there exists $\eps_0 = \eps_0(n,k,\sigma,a,\beta)\in (0,1)$ such that the following is true. Suppose that we have $V\in \mathcal{V}_\beta$, $\BC\in\CC_{n-1}$ is aligned, and $\mathcal{T}:=\mathcal{T}^{S(\BC)}_r(x_0,s_0)\subset B_1(0)$, with $s_0>0$ and $r/s_0\in (a,1-a)$, satisfy the following hypotheses:
	\begin{enumerate}
		\item [(a)] $(\w_n 2^n)^{-1}\|V\|(B^n_2(0))\leq 5/2$;
		\item [(b)] $\nu(\BC) + \hat{E}_V + \hat{E}_V(\mathcal{T})<\eps$.
	\end{enumerate}
	Then, $V\res \mathcal{T}^{S(\BC)}_{\sigma r}(x_0,s_0) = \sum^4_{i=1}V^i$, where for each $i\in \{1,\dotsc,4\}$ we have the following conclusion: either $V^i = 0$, or $V^i$ is a stationary integral $n$-varifold in $\mathcal{T}^{S(\BC)}_{\sigma r}(x_0,s_0)$ for which there exists $\alpha = \alpha(i)\in \{1,\dotsc,4\}$ such that either:
	\begin{enumerate}
		\item [$(1\star)$] $V^i$ is a smooth minimal graph, i.e.~there is a domain $\Omega^\alpha \supset H^\alpha\cap \mathcal{T}^{S(\BC)}_{\sigma r}(x_0,s_0)$ and a smooth function $u^\alpha:\Omega^\alpha\to (H^\alpha)^\perp$ solving the minimal surface system such that:
		\begin{enumerate}
			\item [(i)] $V^i = \graph(u^\alpha)\res \mathcal{T}^{S(\BC)}_{\sigma r}(x_0,s_0)$
			\item [(ii)] For $x\in \Omega^\alpha$ we have $\dist(x+u^\alpha(x),\BC) = |u^\alpha(x)|$.
		\end{enumerate}
	\end{enumerate}
	or alternatively,
	\begin{enumerate}
		\item [$(2\star)$] $V^i$ is well-approximated by a Lipschitz two-valued graph, i.e.~there is a domain $\Omega^\alpha\supset H^\alpha\cap\mathcal{T}^{S(\BC)}_{\sigma r}(x_0,s_0)$, a $\H^n$-measurable subset $\Sigma\subset\Omega^\alpha$, $u: \Omega^\alpha\to \A_2((H^\alpha)^\perp)$ Lipschitz with $\Lip(u)\leq 1/2$, and such that:
		\begin{enumerate}
			\item [(i)] $V^i\res ((\Omega^\alpha\setminus\Sigma)\times (H^\alpha)^\perp) = \graph(u)\res ((\Omega^\alpha\setminus\Sigma)\times (H^\alpha)^\perp)\cap \mathcal{T}^{S(\BC)}_{\sigma r}(x_0,s_0)$;
			\item [(ii)] $\H^n(\Sigma) + \|V^i\|(\Sigma\times (H^\alpha)^\perp)\leq Cr^n\hat{E}_{V,\BC}^2(\mathcal{T})$;
			\item [(iii)] For $x\in \Omega^\alpha\setminus\Sigma$ we have
			$$|u(x)|\leq Cr\hat{E}_{V,\BC}(\mathcal{T}).$$
		\end{enumerate}
	\end{enumerate}
	Here, $C = C(n,k,\sigma,a)\in (0,\infty)$.
\end{lemma}

\begin{proof}
	We start in a similar manner to the proof of Theorem \ref{thm:graphical-rep} in the $\sfrak(\BC)=n-1$ case. Indeed, fix $\sigma\in (1/2,1)$, $a\in (0,1/2)$, and $\beta\in (0,1)$. Then (supposing the lemma is false and looking for a contradiction) we may take sequences $\eps_j\downarrow 0$, $V_j\in \mathcal{V}_\beta$, $\BC_j\in \CC_{n-1}$ aligned, and $\mathcal{T}_j := \mathcal{T}_{r_j}^{S(\BC_j)}(x_j,s_j)$ satisfying the hypotheses of the present lemma, with $\eps_j,V_j,\BC_j,\mathcal{T}_j$ in place of $\eps_0, V,\BC,\mathcal{T}$, respectively.
	
	Notice that from (a) and (b), we may assume that $V_j\weakly \theta_0|P_0|$ as varifolds in $B_1(0)$, where $\theta_0\in \{0,1,2\}$. If $\theta_0 = 0$, the result is trivial (indeed by upper semi-continuity of density, $V^i=0$ for each $i$). Thus from now on we will assume that $\theta_0\in \{1,2\}$.
	
	Next, we make some simplifications analogous to those we did in the proof of Theorem \ref{thm:graphical-rep} in the $\sfrak(\BC)=n-1$ case. Indeed, arguing as in that proof we may rotate to assume that $S(\BC_j) \equiv S_0$ for every $j$, and by applying an appropriate homothety to the system we may assume that $\mathcal{T}_{r_j}^{S(\BC_j)}(x_j,s_j) = \mathcal{T}^{S_0}_1(0,s_j^\prime)=: \mathcal{T}_j^\prime$, where (as $s_j>0$ here) $s_j^\prime\in (1/(1-a),1/a)$ for all $j$ and $s_j^\prime\to s_0^\prime\in [1/(1-a),1/a]$. (Again, this homothety preserves all the assumptions for $j$ sufficiently large, except $\hat{E}_V<\eps$, which we won't need.) Now write
	$$\mathcal{T}^\prime:= \mathcal{T}^{S_0}_1(0,s_0^\prime).$$
	If we can show that the conclusions of the lemma hold along a subsequence, this will complete the proof.
	
	We start by passing to a subsequence along which $V_j\res \mathcal{T}^\prime\weakly V_0$. From (b), we know that $\spt\|V_0\|\subset P_0\cap \mathcal{T}^\prime$. Write $V_j^+:= V_j\res(\mathcal{T}^\prime\cap \{x^n>0\})$ and $V_0^+:= V_0\res \{x^n>0\}$, and also $V_j^-:= V_j\res (\mathcal{T}^\prime\cap\{x^n<0\})$ and $V_0^-:= V_0\res \{x^n<0\}$. We will focus on analysing $V_j^+$ and $V_0^+$, as the arguments that follow will also apply to $V_j^-$ and $V_0^-$.
	
	The constancy theorem gives that $V_0^+$ has constant multiplicity. From our mass upper bound at scale $1$, we know that $V_0^+$ has constant multiplicity either $0$, $1$, or $2$. If the multiplicity is $0$ then we know for sufficiently large $j$ that $\spt\|V_j^+\|\cap \mathcal{T}^{S_0}_\sigma(0,s_0^\prime) = \emptyset$, in which case the conclusions hold for $V_j^+$ (again, with the relevant $V^i$ being 0). Thus, we may assume this constant multiplicity is either $1$ or $2$. For each $j$ now pick a plane $P_j$ with $P_j\supset S_0$ and with
	$$\hat{E}_{V_j^+,P_j}(\mathcal{T}^\prime) = \mathcal{E}_{V_j^+,\BC_j}(\mathcal{T}^\prime).$$
	There is then a sequence of rotations $R_j\in SO(n+k)$ under which $S_0$, $\mathcal{T}_j^\prime$, and $\mathcal{T}^\prime$ are all invariant, and such that $R_j(P_j) = P_0$ and $R_j\to \id_{\R^{n+k}}$. Let us assume, without changing notation, that we rotate by $R_j$, and therefore may replace $\BC_j$ by $(R_j)_\#\BC_j$ and $V_j^+$ by $(R_j)_\#V_j^+$ (by Proposition \ref{prop:rotation-class}, will change the value of $\beta$ slightly, but we abuse notation and still write this as $\beta$). Thus, we have
	$$\hat{E}_{V_j^+}(\mathcal{T}^\prime) = \mathcal{E}_{V_j^+,\BC_j}(\mathcal{T}^\prime).$$
	Up to passing to a further subsequence, we can without loss of generality assume that one of the following holds:
	\begin{enumerate}
		\item [(1)] $\hat{E}_{V_j^+}(\mathcal{T}^\prime)^{-1}\hat{E}_{V_j^+,\BC_j}(\mathcal{T}^\prime)\to 0$;
		\item [(2)] $\hat{E}_{V_j^+,\BC_j}(\mathcal{T}^\prime)\geq C\hat{E}_{V_j^+}(\mathcal{T}^\prime)$ for some $C = C(n,k,\sigma,a,\beta)\in (0,\infty)$.
	\end{enumerate}
	We will consider the case where (1) holds first. Notice that the condition in (1) means that the proof of the first inequality in Lemma \ref{lemma:nu-a}(i) holds (as the proof goes through analogously), and thus we have
	\begin{equation}\label{E:tori-q=1-0}
	C_1\hat{E}_{V_j^+}(\mathcal{T}^\prime) \leq \nu(\BC_j\res \{x^n>0\}) \leq \nu(H^3_j) + \nu(H^4_j)
	\end{equation}
	for all $j$ sufficiently large, where here $\BC_j = \sum^4_{\alpha=1}|H^\alpha_j|$ and $C_1 = C_1(n,k)\in (0,\infty)$. If we now apply Lemma \ref{lemma:appendix-1} (proved in the Appendix to Part 3) with $V_j^+$ in place of $V$ and $\BC_j\res (\mathcal{T}^\prime\cap \{x^n>0\}$) in place of $P$, we guarantee the existence of a non-empty subset $I\subset\{3,4\}$ for which
	\begin{equation}\label{E:tori-q=1-1}
		\sum_{\alpha\in I}\nu(H^\alpha_j) \leq C(\hat{E}_{V_j^+,\BC_j}(\mathcal{T}^\prime) + \hat{E}_{V_j^+}(\mathcal{T}^\prime))
	\end{equation}
	where $C = C(n,k,a)$ (the $a$ dependence coming from the size of the domain $\mathcal{T}^\prime$). When combined with the condition in case (1), this implies that for all $j$ sufficiently large,
	\begin{equation}\label{E:tori-q=1-2}
		\sum_{\alpha\in I}\nu(H^\alpha_j) \leq C\hat{E}_{V_j^+}(\mathcal{T}^\prime).
	\end{equation}
	This tells us that for $\alpha\in I$, up to passing to a subsequence the limit $\hat{E}_{V_j^+}(\mathcal{T}^\prime)^{-1}h^\alpha_j$ exists; call it $h^\alpha_*$ (note that $h^\alpha_*$ is necessarily linear). Now let $v$ denote a coarse blow-up of $(V_j^+)_j$ off $P_0$ in the region $\mathcal{T}^\prime\cap \{x^n>0\}$. Write
	$$L:= \bigcup_{\alpha\in I}\graph(h^\alpha_*).$$
	The condition in (1) implies that we must have $\graph(v)\subset L$ (we may not have equality, as (1) only concerns one-sided excess quantities). This means that $\graph(v)$ is either a single plane, or the union of two planes (or, more precisely, it is formed of the subsets of planes in the region $(Q_1(0,s_0^\prime)\cap \{x^n>0\})\times P_0^\perp$).
	
	We now subdivide into two cases depending on whether $|I|=2$ or $|I|=1$. Up to relabelling, these two cases are (A) $I = \{3,4\}$, and (B) $I = \{3\}$.
	
	Let us first consider case (A), i.e.~$I = \{3,4\}$. We first prove that we have $L^2$ convergence to the blow-up by checking the hypothesis of Proposition \ref{prop:full-L2}. Indeed, we do this as follows. For fixed $\rho\in (0,1)$, the triangle inequality gives
	\begin{align*}
	\int_{(\mathcal{T}^\prime\setminus \mathcal{T}^{S_0}_\rho(0,s_0^\prime))\cap \{x^n>0\}}&\dist^2(x,P_0)\, \ext\|V^+_j\|(x)\\
	& \leq \int_{\mathcal{T}^\prime\cap \{x^n>0\}}\dist^2(x,\BC_j)\, \ext\|V_j^+\|(x) + C(1-\rho)[\nu(H^3_j)^2 + \nu(H^4_j)^2]
	\end{align*}
	for some $C = C(n,k,a)\in (0,\infty)$ Using \eqref{E:tori-q=1-2} (as $I=\{3,4\}$), this gives,
	$$\hat{E}_{V_j^+}(\mathcal{T}^\prime)^{-2}\int_{(\mathcal{T}^\prime\setminus \mathcal{T}^{S_0}_\rho(0,s_0^\prime))\cap \{x^n>0\}}\dist^2(x,P_0)\, \ext\|V_j^+\|(x) \leq \hat{E}_{V_j^+}(\mathcal{T}^\prime)^{-2}\hat{E}_{V_j,\BC_j}(\mathcal{T}^\prime)^2 + C(1-\rho).$$
	So now, taking $j$ sufficiently large we can make the first term on the right-hand side as small as desired (using the condition from (1)), and choosing $\rho$ sufficiently close to $1$ the second term on the right-hand side can be made as small as desired. Thus, Proposition \ref{prop:full-L2} applies, and hence we get strong $L^2$ convergence to the blow-up. Moreover, \eqref{E:tori-q=1-0} tells us that $v\not\equiv 0$.
	
	Suppose, for the sake of contradiction, that $\graph(v)$ was a (piece of a) single plane, i.e.~there is a plane $\tilde{P}$ (necessarily $\neq P_0$, as $v\not\equiv 0$) such that after extending $v$ to a linear function on $P_0$, we have $\graph(v) = \tilde{P}$. Write $\tilde{P}_j:= \graph(\hat{E}_{V_j^+}(\mathcal{T}^\prime)^{-1}v)$. Thus, we conclude from Proposition \ref{prop:full-L2} that
	$$\hat{E}_{V_j^+}(\mathcal{T}^\prime)^{-2}\int_{\mathcal{T}^\prime}\dist^2(x,\tilde{P}_j)\, \ext\|V_j^+\|(x)\to 0.$$
    Notice that we must have $S_0\subset \tilde{P}$, since $\graph(v)\subset L$ and $L$ is the limit of half-planes containing $S_0$. So, since $\tilde{P}_j\supset S_0$, we also have that $\tilde{P}_j$ belongs to the set which the infimum is taken over in the definition of $\mathcal{E}_{V_j^+,\BC_j}(\mathcal{T}^\prime)$ (which by construction equals $\hat{E}_{V_j^+}(\mathcal{T}^\prime)$).Therefore, we should have for all $j$,
	$$\hat{E}_{V_j^+}(\mathcal{T}^\prime)^{-2}\int_{\mathcal{T}^\prime}\dist^2(x,\tilde{P}_j)\, \ext\|V_j^+\|(x) \geq 1.$$
	Clearly these two statements contradict each other. Hence, we have shown that if (1) holds and $I=\{3,4\}$, then $\graph(v)$ is the union of two (pieces of) distinct planes, each of which is a subset of one of the planes of $L$. In this case we may then argue as in the final paragraphs of the proof of Theorem \ref{thm:graphical-rep} in the case $\sfrak(\BC)=n-1$ to conclude that $V_j^+$ is the union of two smooth graphs, each of which satisfies the conclusion $(1\star)$ of the present lemma, completing the proof of this case.
	
	Now suppose that (1) holds and $I = \{3\}$. This means that $\hat{E}_{V_j^+}(\mathcal{T}^\prime)^{-1}\nu(H^4_j) \to \infty$ (indeed if not we can reduce to the previous case). Since $\sfrak(\BC_j) = n-1$ and $s_0^\prime>0$, this in fact implies that
	$$\hat{E}_{V_j^+}(\mathcal{T}^\prime)^{-1}\inf_{x\in H_j^4\cap \mathcal{T}^\prime,\, y\in P_0\cap \mathcal{T}^\prime}|x-y|\to\infty.$$
	But from Allard's supremum estimate (Theorem \ref{thm:allard-sup-estimate}), we have that
	$$\sup_{x\in \spt\|V_j^+\|\cap \mathcal{T}^{S_0}_{(1+\sigma)/2}(0,s_0^\prime)}|\pi_{P_0^\perp}(x)| \leq C\hat{E}_{V_j^+}(\mathcal{T}^\prime)$$
	where $C = C(n,k,\sigma)\in (0,\infty)$. In light of these estimates and \eqref{E:tori-q=1-2}, we actually have
	$$\int_{\mathcal{T}_\sigma^{S_0}(0,s_0^\prime)}\dist^2(x,\BC_j)\, \ext\|V_j^+\|(x) = \int_{\mathcal{T}_\sigma^{S_0}(0,s_0^\prime)}\dist^2(x,H^3_j)\, \ext\|V_j^+\|(x).$$
	In this case, after applying Theorem \ref{thm:Lipschitz-approx} to approximate $V_j^+$ by a Lipschitz $2$-valued graph over $H^3_j$, we deduce the conclusions of $(2\star)$ hold in this case. This concludes the proof of the lemma under assumption (1).
	
	Finally, let us suppose that (2) holds. Applying the condition in (2) to \eqref{E:tori-q=1-1}, we get (up to relabelling, as $I$ is non-empty)
	$$\nu(H_j^3) \leq C\hat{E}_{V_j^+,\BC_j}(\mathcal{T}^\prime)$$
	where $C = C(n,k,\sigma,a,\beta)$. Combining this with the condition in (2) again, we see that
	$$\hat{E}_{V_j^+,H_j^3}(\mathcal{T}^\prime)\leq C\hat{E}_{V_j^+,\BC_j}(\mathcal{T}^\prime).$$
	Then, once again, for all sufficiently large $j$ we can use Theorem \ref{thm:Lipschitz-approx} to approximate $V_j^+$ by a Lipschitz $2$-valued graph over $H^3_j$ (or using Allard's regularity theorem, if the convergence is to a multiplicity one plane) and use the above inequality to deduce the conclusions of $(2\star)$ hold in this case. Applying all of the above arguments to $V_j^-$ also, we therefore complete the proof of the lemma.
\end{proof}

We are now in a position to prove Lemma \ref{lemma:Z} when $\sfrak(\BC)=n-1$.

\begin{proof}[Proof of Lemma \ref{lemma:Z} when $\sfrak(\BC)=n-1$]
Let $\eps_0$ be as in the statement of Lemma \ref{lemma:tori-q=1}. There are two cases to consider: (A) $\hat{E}_{V,\BC}(\mathcal{T})<\eps_0$ or (B) $\hat{E}_{V,\BC}(\mathcal{T})\geq \eps_0$.

First consider case (B). This means that
$$r^{n+2}\leq \eps_0^{-1}\int_{\mathcal{T}}\dist^2(x,\BC)\, \ext\|V\|(x).$$
We clearly in general have
\begin{equation}\label{E:Z-q=1-1}
	\sup_{x\in\mathcal{T}}\phi_{(0)}(x)|x^{\perp_{S(\BC)}}|^2 + \sup_{x\in\mathcal{T}}|x^{\perp_{S(\BC)}\perp_{T_xV}}|^2 \leq Cr^2
\end{equation}
and
$$\|\BC\|(\mathcal{T}) + \|V\|(\mathcal{T})\leq Cr^n$$
where $C = C(n,k,a,\tilde{C})\in (0,\infty)$. Therefore,
\begin{align*}
	\int_{\mathcal{T}_{\sigma r}^{S(\BC)}}|x^{\perp_{S(\BC)}\perp_{T_xV}}|^2\, \ext\|V\|(x) + \left|\int\phi_{(0)}|x^{\perp_{S(\BC)}}|^2\, \ext\|V\|(x)\right. & - \left.\int\phi_{(0)}|x^{\perp_{S(\BC)}}|^2\, \ext\|\BC\|(x)\right|\\
	& \leq Cr^{n+2}\\
	& \leq C\int_{\mathcal{T}}\dist^2(x,\BC)\, \ext\|V\|(x)
\end{align*}
where $C = C(n,k,\sigma,a,\beta,\tilde{C})$. This completes the proof of the lemma in the case of (B).

Now turn to case (A). Here, we are in a position to apply Lemma \ref{lemma:tori-q=1} (up to perhaps decreasing $\eps_0$ depending on $M$). Thus, using the notation from Lemma \ref{lemma:tori-q=1}, we have
$$V\res \mathcal{T}^{S(\BC)}_{(1+\sigma)r/2}(x_0,s_0) = \sum_{i\in I}V^i$$
where here we write $I\subset\{1,2,3,4\}$ for the subset of indices for which $V^i$ is non-zero.

Fix $i\in I$ such that $(2\star)$ of Lemma \ref{lemma:tori-q=1} holds; we will focus on this case as the case where $(1\star)$ holds is strictly simpler. Let $\gamma = \gamma(i)$ be such that we have $H^\gamma$, $\Omega^\gamma$, $\Sigma$, and $u:\Omega^\gamma\to \A_2((H^\gamma)^\perp)$ as described in Lemma \ref{lemma:tori-q=1}$(2\star)$. Now, for any $x\in \graph(u)\cap ((\Omega^\gamma\setminus\Sigma)\times (H^\gamma)^\perp)$, if we write $x^\prime := \pi_{H^\gamma}(x)$, then there is some $\alpha\in\{1,2\}$ for which
\begin{align*}
	x^{\perp_{S(\BC)}\perp_{T_xV}} & = (x^\prime + u^\alpha(x^\prime))^{\perp_{S(\BC)}\perp_{T_xV}}\\
	& = \pi_{T^\perp_xV}\left((x^\prime)^{\perp_{S(\BC)}}\right) + \pi_{T^\perp_{x}V}\left(u^\alpha(x^\prime)^{\perp_{S(\BC)}}\right)\\
	& = (\pi_{T^\perp_xV}-\pi_{T^\perp_{x^\prime}\BC})((x^\prime)^{\perp_{S(\BC)}}) + \pi_{T_{x^\prime}^\perp\BC}((x^\prime)^{\perp_{S(\BC)}}) + \pi_{T_x^\perp V}(u^\alpha(x^\prime)^{\perp_{S(\BC)}}).
\end{align*}
The second term here is clearly zero. Moreover, since $u^\alpha(x^\prime)\in (H^\gamma)^\perp\subset S(\BC)^\perp$, in the last term we have $u^\alpha(x^\prime)^{\perp_{S(\BC)}} = u^\alpha(x^\prime)$. Therefore, we get
$$|x^{\perp_{S(\BC)}\perp_{T_xV}}|^2 \leq 2\left[\|\pi_{T_xV}-\pi_{T_{x^\prime}\BC}\|^2|(x^\prime)^{\perp_{S(\BC)}}|^2 + |u^\alpha(x^\prime)|^2\right].$$
Using this, the area-formula \eqref{E:area-formula}, the Jacobian bounds \eqref{E:J-bounds}, the reverse Poincaré inequality \eqref{E:rpi} (indeed, note that here $x^\prime\in H^\gamma$, and so $\pi_{T_{x^\prime}\BC}$ is the projection onto a fixed plane), the inequality $\sup_{x\in\mathcal{T}}|(x^\prime)^{\perp_{S(\BC)}}|^2 \leq Cr^2$ where $C = C(n,k,a)$, as well as conclusions (i) and (iii) of Lemma \ref{lemma:tori-q=1}$(2\star)$, we have
\begin{align*}
	& \int_{\mathcal{T}^{S(\BC)}_{(1+\sigma)r/2}(x_0,s_0)\cap ((\Omega^\gamma\setminus \Sigma)\times (H^\gamma)^\perp)}|x^{\perp_{S(\BC)}\perp_{T_xV}}|^2\, \ext\|V^i\|(x)\\
	& \hspace{6em} \leq Cr^2\int_{\mathcal{T}^{S(\BC)}_{(1+\sigma)r/2}(x_0,s_0)\cap ((\Omega^\gamma\setminus \Sigma)\times (H^\gamma)^\perp)} \|\pi_{T_xV}-\pi_{T_{x^\prime}\BC}\|^2\, \ext\|V^i\|(x)\\
	& \hspace{15em} + C\int_{\int_{\mathcal{T}^{S(\BC)}_{(3+\sigma)r/4}(x_0,s_0)\cap (\Omega^\gamma\setminus \Sigma)}}\sum_{\alpha=1,2}|u^\alpha(x)|^2\, \ext\H^n(x)\\
	& \hspace{6em} \leq C\int_{\mathcal{T}}\dist^2(x,\BC)\, \ext\|V\|(x)
\end{align*}
where $C = C(n,k,\sigma,a)$. Moreover, using \eqref{E:Z-q=1-1} and conclusion (ii) of Lemma \ref{lemma:tori-q=1}$(2\star)$, we have
\begin{align*}
	\int_{\mathcal{T}^{S(\BC)}_{(1+\sigma)r/2}(x_0,s_0)\cap (\Sigma\times (H^\gamma)^\perp)}|x^{\perp_{S(\BC)}\perp_{T_xV}}|^2\, \ext\|V^i\|(x) & \leq Cr^2\|V^i\|(\Sigma\times (H^\gamma)^\perp)\\
	& \leq C\int_{\mathcal{T}}\dist^2(x,\BC)\, \ext\|V\|(x).
\end{align*}
Combining the above two inequalities and summing over $i$ (analogous, in fact simpler, arguments also give the above when $i\in I$ is such that Lemma \ref{lemma:tori-q=1}$(1\star)$ holds for $V^i$) completes the proof of the first estimate \eqref{E:Z-a} in Lemma \ref{lemma:Z}.

Now we turn our attention to proving the second estimate \eqref{E:Z-b} in Lemma \ref{lemma:Z}. Again, we will consider the case that $i\in I$ is such that Lemma \ref{lemma:tori-q=1}$(2\star)$ holds, as the case of Lemma \ref{lemma:tori-q=1}$(1\star)$ is strictly simpler.

Replacing $V^i$ by its graphical approximation $\graph(u)$ where possible, we see that
$$\int_{\mathcal{T}^{S(\BC)}_{\sigma r}(x_0,s_0)}\phi_{(0)}(x)|x^{\perp_{S(\BC)}}|^2\, \ext\|V^i\|(x) = \int_{\mathcal{T}^{S(\BC)}_{\sigma r}(x_0,s_0)}\phi_{(0)}(x)|x^{\perp_{S(\BC)}}|^2\, \ext\|\mathbf{v}(u)\|(x) + E_1$$
where we have, using \eqref{E:Z-q=1-1} and (ii) of Lemma \ref{lemma:tori-q=1}$(2\star)$,
\begin{align*}
	|E_1|& := \left|\int_{\mathcal{T}^{S(\BC)}_{\sigma r}(x_0,s_0)\cap (\Sigma\times (H^\gamma)^\perp)}\phi_{(0)}(x)|x^{\perp_{S(\BC)}}|^2\, \ext\|V^i\|(x)\right.\\
	& \hspace{5em}\left. -\int_{\mathcal{T}_{\sigma r}^{S(\BC)}(x_0,s_0)\cap (\Sigma\times (H^\gamma)^\perp)}\phi_{(0)}|x^{\perp_{S(\BC)}}|^2\, \ext\|\mathbf{v}(u)\|(x)\right|\\
	& \leq Cr^2\|V\|(\Sigma\times (H^\gamma)^\perp) + Cr^2\H^n(\Sigma)\\
	& \leq C\int_{\mathcal{T}}\dist^2(x,\BC)\, \ext\|V\|(x)
\end{align*}
where $C = C(n,k,\sigma,a)$. Now, use the area formula \eqref{E:area-formula} to write the integral over $\graph(u)$ over the domain $\Omega^\gamma$ in $H^\gamma$. Indeed, we have that
\begin{align*}
	\int_{\mathcal{T}^{S(\BC)}_{\sigma r}(x_0,s_0)}&\phi_{(0)}(x)|x^{\perp_{S(\BC)}}|^2\, \ext\|\mathbf{v}(u)\|(x)\\
	& = \int_F\sum_{\alpha=1,2}\phi_{(0)}(x+u^\alpha(x))\left(|x^{\perp_{S(\BC)}}|^2 + |u^\alpha(x)|^2\right)J^\alpha(x)\, \ext\H^n(x)
\end{align*}
where the region of integration $F\subset H^\gamma$ is
$$F:= \pi_{H^\gamma}\left(\graph(u)\cap \mathcal{T}^{S(\BC)}_{\sigma r}(x_0,s_0)\right).$$
Note that we can insist that $F\subset \mathcal{T}^{S(\BC)}_{(3+\sigma)r/4}(x_0,s_0)$ by choosing $\eps_0$ sufficiently small. Next, we wish to replace the integrand on the right-hand side of the above by the expression $2\phi_{(0)}(x)|x^{\perp_{S(\BC)}}|^2$, and control the error incurred in doing so. Indeed, we claim that
\begin{align}
\nonumber\left|\int_F\sum_{\alpha=1,2}\phi_{(0)}(x+u^\alpha(x))\left(|x^{\perp_{S(\BC)}}|^2 + |u^\alpha(x)|^2\right)\right.&J^\alpha(x)\, \ext\H^n(x)\\
& \left. - \int_{H^\gamma}2\phi_{(0)}(x)|x^{\perp_{S(\BC)}}|^2\, \ext\H^n(x)\right| \leq E_2\label{E:Z-q=1-2}
\end{align}
where
$$E_2 \leq C\int_{\mathcal{T}}\dist^2(x,\BC)\, \ext\|V\|(x)$$
with $C = C(n,k,\sigma,a,\tilde{C})$. To see this, one argues in the following manner:
\begin{itemize}
	\item Firstly, use (iii) of Lemma \ref{lemma:tori-q=1}$(2\star)$ to control the term containing $|u^\alpha(x)|^2$: since $\phi_{(0)}$ and $J^\alpha$ are uniformly bounded by some constant depending only on $n,k,\tilde{C}$, this enables us to control this contribution to the error term $E_2$ as above. After this, we are left with the integral
	$$\int_F \sum_{\alpha=1,2}\phi_{(0)}(x+u^\alpha(x))|x^{\perp_{S(\BC)}}|^2J^\alpha(x)\, \ext\H^n(x).$$
	\item Now we replace $\phi_{(0)}(x+u^\alpha(x))$ by $\phi_{(0)}(x)$. To do this, recall that by assumption $\phi_{(0)}(x)$ depends only on $x^{\top_{S(\BC)}}$ and $|x^{\perp_{S(\BC)}}|$. Notice that for fixed $x\in F$ and $\alpha\in\{1,2\}$, if $y = x+u^\alpha(x)$, then we have
	$$y^{\top_{S(\BC)}} = x^{\top_{S(\BC)}} \qquad \text{and} \qquad |y^{\perp_{S(\BC)}}| = \sqrt{|x^{\perp_{S(\BC)}}|^2 + |u^\alpha(x)|^2}.$$
	Thus, we are looking at
	$$\phi_{(0)}(y) = \phi_{(0)}\left(x^{\top_{S(\BC)}},\sqrt{|x^{\perp_{S(\BC)}}|^2 + |u^\alpha(x)|^2}\right).$$
	If we apply the mean value theorem to the function
	$$t\longmapsto \phi_{(0)}\left(x^{\top_{S(\BC)}},\sqrt{|x^{\perp_{S(\BC)}}|^2 + t}\right)$$
	over the interval $[0,|u^\alpha(x)|^2]$ (of course, if $|u^\alpha(x)|=0$ there is nothing to prove) we can deduce that
	$$|\phi_{(0)}(y)-\phi_{(0)}(x)| \leq |u^\alpha(x)|^2|y^{\perp_{S(\BC)}}|^{-1}\sup_F|D\phi_{(0)}|$$
	(here, we are also using the simple fact that the function $t\mapsto \frac{t}{\sqrt{b^2+t^2}}$ is increasing for $t>0$ for any $b\in \R$). By our assumptions on $\phi_{(0)}$ we know that $|D\phi_{(0)}|\leq \tilde{C}r^{-1}$, and since also $\sup_{x\in F}|y^{\perp_{S(\BC)}}|^{-1}\leq Cr^{-1}$ for suitable $C = C(n,k,\sigma,a)$, we therefore get
	$$|\phi_{(0)}(x+u^\alpha(x))-\phi_{(0)}(x)| \leq Cr^{-2}|u^\alpha(x)|^2.$$
	Thus, the error we get when replacing $\phi_{(0)}(x+u^\alpha(x))$ by $\phi_{(0)}(x)$ is controlled by
	$$Cr^{-2}\int_{F}\sum_{\alpha=1,2}|u^\alpha(x)|^2|x^{\perp_{S(\BC)}}|^2J^\alpha(x)\, \ext\H^n(x).$$
	Using the fact that $|x^{\perp_{S(\BC)}}|^2 \leq Cr^2$ on $F$, the powers of $r$ cancel in the above, and so using \eqref{E:J-bounds} and (iii) of Lemma \ref{lemma:tori-q=1}$(2\star)$, we can control this error term by $C\int_{\mathcal{T}}\dist^2(x,\BC)\, \ext\|V\|(x)$ as desired. Thus, the integrand we are left with now is
	$$\int_{F}\sum_{\alpha=1,2}\phi_{(0)}(x)|x^{\perp_{S(\BC)}}|^2J^\alpha(x)\, \ext\H^n(x).$$
	\item Finally, to replace $\sum_{\alpha=1,2}J^{\alpha}(x)$ by $2$, notice that for $\H^n$-a.e. $x\in F$, \eqref{E:J-bounds} gives
	$$2\leq \sum_{\alpha=1,2}J^\alpha(x) \leq 2 + C\sum_{\alpha=1,2}|Du^\alpha(x)|^2.$$
	Thus, the replacement gives an error term of the form
	$$\int_F\sum_{\alpha=1,2}\phi_{(0)}|x^{\perp_{S(\BC)}}|^2|Du^\alpha(x)|^2\, \ext\H^n(x).$$
	To control this error term, we argue similarly to how we have before, passing back from the graphical approximation to $V^i$ and using \eqref{E:rpi} coupled then with (iii) of Lemma \ref{lemma:tori-q=1}$(2\star)$. All the error terms incurred when doing this are also controlled by the desired quantity using (ii) of Lemma \ref{lemma:tori-q=1}$(2\star)$.
\end{itemize}
Thus, we have now established \eqref{E:Z-q=1-2}. Combining all the above, we have now shown that when $i$ is such that Lemma \ref{lemma:tori-q=1}$(2\star)$ holds, we have
$$\left|\int\phi_{(0)}(x)|x^{\perp_{S(\BC)}}|^2\, \ext\|V^i\|(x) - 2\int_{H^\gamma}\phi_{(0)}(x)|x^{\perp_{S(\BC)}}|^2\, \ext\H^n(x)\right| \leq E_3$$
where
$$E_3 \leq C\int_{\mathcal{T}}\dist^2(x,\BC)\, \ext\|V\|(x)$$
and here $C = C(n,k,\sigma,a)$. If instead $i$ is such that Lemma \ref{lemma:tori-q=1}$(1\star)$ holds, then an identical (and in fact simpler) argument to the above gives the same inequality, except there is no factor of $2$ on the second integral on the left-hand side (this factor corresponds exactly to the number of values in the graphical approximation, or equivalently the multiplicity of half-plane $V^i$ is close to). Notice that because the second integrand on the left-hand side of \eqref{E:Z-q=1-2} is independent of the choice of half-plane $H^\gamma$ in $\BC$, splitting the integral into a sum over $i$ we get
$$\left|\int\phi_{(0)}(x)|x^{\perp_{S(\BC)}}|^2\, \ext\|V\|(x) - \int\phi_{(0)}(x)|x^{\perp_{S(\BC)}}|^2\, \ext\|\BC\|(x)\right| \leq 4E_3$$
which gives \eqref{E:Z-b}. This therefore completes the proof of Lemma \ref{lemma:Z} when $\sfrak(\BC)=n-1$.
\end{proof}

We are now in a position to prove Lemma \ref{lemma:L2-estimates-2} when $\sfrak(\BC)=n-1$.

\begin{proof}[Proof of Lemma \ref{lemma:L2-estimates-2} when $\sfrak(\BC)=n-1$]
	Using the assumption that $\Theta_V(0)\geq 2$, a computation with the first variation formula (cf.~\cite[Lemma 3.4]{Sim93}) gives
	\begin{align}
		\nonumber\int_{B_\sigma(0)}\frac{|x^{\perp_{T_xV}}|^2}{|x|^{n+2}}&\, \ext\|V\|(x) + \frac{1}{2}\int_{B_\sigma(0)}\sum^{n-1}_{j=1}|e_j^{\perp_{T_xV}}|^2\, \ext\|V\|(x)\\
		& \leq C^\prime\int_{B_{(3+\sigma)/4}(0)}|x^{\perp_{S(\BC)}\perp_{T_xV}}|^2\, \ext\|V\|(x)\nonumber\\
		& \hspace{4em} + \int_{B_1(0)}\phi(|x|)|x^{\perp_{S(\BC)}}|^2\, \ext\|V\|(x) - \int_{B_1(0)}\phi(|x|)|x^{\perp_{S(\BC)}}|^2\, \ext\|\BC\|(x)\label{E:L2-2-q=1-1}
	\end{align}
	where here $C^\prime = C^\prime(n,k,\sigma)\in (0,\infty)$ and $\phi:[0,\infty)\to [0,\infty)$ is a smooth function with $\spt(\phi)\subset \left[\frac{1+\sigma}{2},\frac{3+\sigma}{4}\right)$ and such that $\sup|\phi| + \sup|\phi^\prime|\leq C$ for some absolute constant $C = C(\sigma)$. (Recall also that $e_1,\dotsc,e_{n-1}$ are a basis of $S(\BC)$ by assumption. We also note that, when compared with \cite[Lemma 3.4]{Sim93}, the above expression follows by setting $\phi(|x|)$ equal to $-2|x|^{-1}\psi(|x|)\psi^\prime(|x|)$ therein.) 
		
	Notice that the right-hand side of \eqref{E:L2-2-q=1-1} has the terms which are controlled in Lemma \ref{lemma:Z}. We will bring in Lemma \ref{lemma:Z} on a certain dyadic decomposition of $\Lambda(S(\BC))$. Indeed, we will construct a countable, pairwise disjoint, family $\mathscr{F}$ of dyadic cubes in $\Lambda(S(\BC))$ such that
	\begin{align}
		\nonumber\{(x,s)\in S(\BC)\times [0,1): |x|^2+s^2<\,&\, (3+\sigma)/4\}\\
		&\subset \bigcup_{Q\in\mathscr{F}}Q \subset \{(x,s)\in S(\BC)\times [0,1):|x|^2+s^2<(7+\sigma)/8\}\label{E:L2-2-q=1-2}
	\end{align}
	and satisfying certain properties. To insist on these inclusions, we need to ensure that the diameters of the cubes are sufficiently small. So, begin by denoting
	$$J:=\min\{j\geq 1: 2^{-j}\sqrt{\sfrak(\BC)+1}<(1-\sigma)/8\}.$$
	Then, let $\mathscr{F}_J$ denote all of the cubes $Q\in\mathscr{C}_J$ for which
	$$Q\cap \{(x,s)\in S(\BC)\times [0,1):|x|^2+s^2<(3+\sigma)/4\}\neq\emptyset.$$
	Then, for $p\geq J$, let $\mathscr{F}_{p+1}$ be obtained from $\mathscr{F}_p$ by bisecting every cube in $\mathscr{F}_p$ which is adjacent to $S(\BC)$ and keeping the rest, i.e.
	$$\mathscr{F}_{p+1}:=\{Q\in\mathscr{F}_p:Q\text{ not adjacent to }S(\BC)\}\cup\{Q\in\mathscr{C}_{p+1}:\hat{Q}\in\mathscr{F}_p\text{ and is adjacent to }S(\BC)\}$$
	where here $\hat{Q}$ denotes the parent of $Q$ in $\mathscr{C}$, i.e.~it is the cube for which $Q$ is one of the cubes which arises from bisecting $\hat{Q}$. Then set
	$$\mathscr{F}:=\bigcup_{p^\prime\geq J}\bigcap_{p\geq p^\prime}\mathscr{F}_p.$$
	Clearly the family $\mathscr{F}$ obeys the inclusions in \eqref{E:L2-2-q=1-2}. Moreover, note that the family $\mathscr{F}$ has the following properties:
	\begin{enumerate}
		\item [$(\mathscr{F}1)$] For every cube $Q = Q_r(x_0,s_0)\in\mathscr{F}$ we have $s_0>0$, and either $e(Q) = 2^{-J}$ (a constant depending on $\sfrak(\BC)$ and $\sigma$) or $e(Q)<2^{-J}$ with $\dist(Q,S(\BC))=2e(Q)$. This latter expression is equivalent to $\inf_{x\in Q,\, y\in S(\BC)}|x-y|=2e(Q)$. These conditions, along with the definition of $\mathscr{C}$, mean that in both cases there is a constant $a = a(n,\sigma)\in (0,1/2)$ for which $r/s_0\in (a,1-a)$.
		\item [$(\mathscr{F}2)$] If $Q_1,Q_2\in\mathscr{F}$ are adjacent, then $\frac{1}{2}e(Q_1)\leq e(Q_2)\leq 2e(Q_1)$.
		\item [$(\mathscr{F}3)$] There is a constant $N = N(n)\geq 1$ such that each point of $\{(x,s)\in S(\BC)\times [0,1):|x|^2+s^2<1\}$ is contained in at most $N$ cubes from the collection $\{\frac{5}{4}Q:Q\in\mathscr{F}\}$, i.e.
		$$\sum_{Q\in\mathscr{F}}\one_{\frac{5}{4}Q}\leq N \qquad \text{on }B_{1}(0).$$
		(This also follows in a standard manner from $(\mathscr{F}2)$, cf.~\cite[Chapter VI, Section 1.3]{Ste70}.)
	\end{enumerate}
	By construction, we clearly have
	$$B_{(3+\sigma)/4}(0)\subset\bigcup_{Q\in\mathscr{F}}\mathcal{T}^{S(\BC)}(Q)\subset B_{(7+\sigma)/8}(0).$$
	By $(\mathscr{F}3)$, we can also construct for each cube $Q\in\mathscr{F}$ a smooth function $\phi_Q:\Lambda(S(\BC))\to [0,1]$ such that:
	\begin{enumerate}
		\item [$(\mathscr{F}4)$] $\spt(\phi_Q)\subset\frac{9}{8}Q$;
		\item [$(\mathscr{F}5)$] $\sum_Q\phi_Q(y) = 1$ for all $y\in \cup_{Q\in\mathscr{F}}Q$;
		\item [$(\mathscr{F}6)$] $|D\phi_Q(x,s)|\leq Cs^{-1}$ for $(x,s)\in\Lambda(S(\BC))$, where $C = C(n,\sigma)\in (0,\infty)$. (And analogous bounds hold for higher derivatives.)
	\end{enumerate}
	Now for each $Q\in\mathscr{F}$, let $\tilde{\phi}_Q:B_1(0)\to [0,1]$ be the smooth function on $B_1(0)$ defined by:
	$$\tilde{\phi}_Q(x):= \phi_Q(x^{\top_{S(\BC)}},|x^{\perp_{S(\BC)}}|).$$
	Since $\tilde{\phi}_Q(x)$ depends only on $x^{\top_{S(\BC)}}$ and $|x^{\perp_{S(\BC)}}|$, from $(\mathscr{F}5)$ we get
	\begin{enumerate}
		\item [$(\mathscr{F}7)$] $\sum_{Q\in\mathscr{F}}\tilde{\phi}_Q(x) = 1$ for all $x\in B_{(3+\sigma)/4}(0)$.
	\end{enumerate} 
	From $(\mathscr{F}6)$ we have
	$$|D\tilde{\phi}_Q(x)|\leq C|x^{\perp_{S(\BC)}}|^{-1}$$
	for all $x\in B_1(0)$, where $C = C(n,\sigma)$. Hence we see that, for each $Q = Q_r(x_0,s_0)\in\mathscr{F}$, the function
	$$\phi_{(0)}(x):= \tilde{\phi}_Q(x)\phi(|x|)$$
	satisfies the hypotheses of Lemma \ref{lemma:Z}, with $\mathcal{T} = \mathcal{T}^{S(\BC)}_{5r/4}(x_0,s_0)$ and $\sigma = 9/10$ therein. Thus, provided $\eps_0$ is sufficiently small depending on $n,k,\sigma,\beta$ (as here $a = a(n,\sigma)$) are sufficiently small, we may apply Lemma \ref{lemma:Z} (which we know to be true when $\sfrak(\BC)=n-1$) to deduce
	$$\int_{\mathcal{T}^{S(\BC)}_{9r/8}(x_0,s_0)}|x^{\perp_{S(\BC)}\perp_{T_xV}}|^2\, \ext\|V\|(x) \leq C\int_{\mathcal{T}^{S(\BC)}_{5r/4}(x_0,s_0)}\dist^2(x,\BC)\, \ext\|V\|(x)$$
	and
	\begin{align*}
		\int\tilde{\phi}_Q(x)\phi(|x|)|x^{\perp_{S(\BC)}}|^2\, \ext\|V\|(x) - \int\tilde{\phi}_Q(x)\phi(|x|)&|x^{\perp_{S(\BC)}}|^2\, \ext\|\BC\|(x)\\
		& \leq C\int_{\mathcal{T}^{S(\BC)}_{5r/4}(x_0,s_0)}\dist^2(x,\BC)\, \ext\|V\|(x).
	\end{align*}
	where here $C = C(n,k,\sigma,\beta)\in (0,\infty)$. Now if we add these two inequalities together and sum the resulting inequality over $Q\in\mathscr{F}$, using $(\mathscr{F}3)$ and $(\mathscr{F}7)$, we can justify swapping the infinite sums with the integrals to deduce that (also using that by construction $\frac{5}{4}Q\subset B_1(0)$ for $Q\in\mathscr{F}$)
	\begin{align*}
	\int_{B_{(3+\sigma)/4}(0)}|x^{\perp_{S(\BC)}\perp_{T_xV}}|^2\, \ext\|V\|(x)& + \int\phi(|x|)|x^{\perp_{S(\BC)}}|^2\, \ext\|V\|(x)\\
	& - \int\phi(|x|)|x^{\perp_{S(\BC)}}|^2\, \ext\|\BC\|(x) \leq C\int_{B_1(0)}\dist^2(x,\BC)\, \ext\|V\|(x).
	\end{align*}
	Combining this with \eqref{E:L2-2-q=1-1}, we therefore deduce the first two estimates \eqref{E:initial-1} and \eqref{E:initial-2} of the lemma.
	
	Let us now turn to the third estimate \ref{E:initial-3}. Modifying slightly the construction in \cite[Page 615]{Sim93}, we see in light of Lemma \ref{lemma:nu-a} that there is a constant $C = C(n,k)$ and a Lipschitz function $d:B_1(0)\to [0,\infty)$ which is smooth away from $S(\BC)$ such that:
	\begin{enumerate}
		\item [(i)] $d$ is homogeneous of degree $1$;
		\item [(ii)] $C^{-1}\dist(x,\BC)\leq d(x)\leq C\dist(x,\BC)$ everywhere;
		\item [(iii)] $\Lip(d)\leq C$.
	\end{enumerate}
	(In fact, $d$ can be made to agree with $\dist(\cdot,\BC)$ in the region $\{x\in B_1(0)\setminus S(\BC):\dist(x,\BC)<C\hat{E}_{V,\BC}x^{\perp_{S(\BC)}}\}$.)
	
	Using the fact that $\left|\frac{x}{|x|}\cdot\nabla^V|x|\right| \leq 1$, and since by direct computation
	$$\div_{T_xV}(x|x|^{-n+\mu}) = \mu|x|^{-n+\mu} + n|x|^{-n+\mu}\left(1-\frac{x}{|x|}\cdot\nabla^V|x|\right)$$
	we see that
	\begin{equation}\label{E:L2-2-q=1-3}
	\div_{T_xV}(x|x|^{-n+\mu})\geq \mu|x|^{-n+\mu}.
	\end{equation}
	Then, applying the first variation formula with the vector field $\Phi(x):= \zeta(x)^2\frac{d(x)^2}{|x|^2}|x|^{-n+\mu}x$, where $\zeta\in C^\infty_c(B_1(0))$ and $\mu\in (0,1)$ (technically as $\Phi$ is singular at $0$ and $d$ is only Lipschitz along $S(\BC)$ one needs to justify its use in the first variation formula by an argument cutting off both $S(\BC)$ and $\{0\}$ separately, but this is a standard argument), we have using \eqref{E:L2-2-q=1-3},
	$$\int\mu|x|^{-n+\mu}\zeta(x)^2\frac{d(x)^2}{|x|^2}\, \ext\|V\|(x) \leq -\int |x|^{-n+\mu}x\cdot\nabla^V\left(\zeta(x)^2\frac{d(x)^2}{|x|^2}\right)\, \ext\|V\|(x).$$
	Using the product rule on the integrand on the right-hand side gives two terms. First, let us analyse the term which comes from the derivative hitting $\zeta(x)^2$. Using $2ab\leq a^2+b^2$ for suitable $a,b$, together with the fact that $d$ is comparable to $\dist(\cdot,\BC)$ by (ii) above, this term is estimate as follows:
	\begin{align*}
		-2\int|x|^{-n+\mu}\zeta(x)\frac{d(x)^2}{|x|^2}& x\cdot\nabla^V\zeta\, \ext\|V\|(x)\\
		& \leq C\int\frac{\dist^2(x,\BC)}{|x|^{n-\mu}}|\nabla^V\zeta|^2\, \ext\|V\|(x) + \frac{\mu}{4}\int_{B_1(0)}\zeta(x)^2\frac{d(x)^2}{|x|^{n+2-\mu}}\, \ext\|V\|(x).
	\end{align*}
	Now look at the second term from the product rule, namely when the derivative hits $d(x)^2/|x|^2$. Given that $d$ is homogeneous of degree one, we know that
	$$x\cdot D\left(\frac{d(x)}{|x|}\right)\equiv 0$$
	(here, $D$ is the ambient gradient operator on $\R^{n+k}$). Thus, we have
	$$x\cdot\nabla^V\left(\frac{d(x)^2}{|x|^2}\right) = -2\frac{d(x)}{|x|}x^{\perp_{T_xV}}\cdot D\left(\frac{d(x)}{|x|}\right).$$
	Additionally using the fact that $|D(d(x)/|x|)| \leq C/|x|$ (which follows from properties (i) and (iii) of $d$ above), as well as $2ab\leq a^2+b^2$, we estimate this term as follows:
	\begin{align*}
		\int 2\zeta(x)^2|x|^{-n+\mu}&\frac{d(x)}{|x|}x^{\perp_{T_xV}}\cdot D\left(\frac{d(x)}{|x|}\right)\, \ext\|V\|(x)\\
		& \leq C\int\zeta(x)^2\frac{|x^{\perp_{T_xV}}|^2}{|x|^{n+2-\mu}}\, \ext\|V\|(x) + \frac{\mu}{4}\int\zeta(x)^2\frac{d(x)^2}{|x|^{n+2-\mu}}\, \ext\|V\|(x).
	\end{align*}
	Combining all of the above, we therefore arrive at
	\begin{align*}
	\frac{\mu}{2}\int_{B_1(0)}\zeta^2(x)\frac{d(x)^2}{|x|^{n+2-\mu}}&\, \ext\|V\|(x)\\
	& \leq C\int\frac{\dist^2(x,\BC)}{|x|^{n-\mu}}|\nabla^V\zeta|^2\, \ext\|V\|(x) + C\int\zeta(x)^2\frac{|x^{\perp_{T_xV}}|^2}{|x|^{n+2-\mu}}\, \ext\|V\|(x).
	\end{align*}
	Using again property (ii) for $d$, this gives
	\begin{align*}
	\int_{B_1(0)}\zeta^2(x)\frac{\dist^2(x,\BC)}{|x|^{n+2-\mu}}&\, \ext\|V\|(x)\\
	& \leq C\int\frac{\dist^2(x,\BC)}{|x|^{n-\mu}}|\nabla^V\zeta|^2\, \ext\|V\|(x) + C\int\zeta(x)^2\frac{|x^{\perp_{T_xV}}|^2}{|x|^{n+2-\mu}}\, \ext\|V\|(x).
	\end{align*}
	Here, $C = C(n,k,\mu)$. Now choose $\zeta$ such that $\zeta\equiv 1$ on $B_\sigma(0)$, $|\zeta|\leq 1$ everywhere, $\zeta\equiv 0$ outside $B_{(1+\sigma)/2}(0)$, and $|D\zeta|\leq C\sigma^{-1}$. Notice that this means $|x|^{-n+\mu}|\nabla^V\zeta(x)| \leq C\sigma^{-n+\mu-1}$ pointwise everywhere, and so the first integral on the right-hand side above is controlled. For the second integral, we can appeal to \eqref{E:initial-1} with $(1+\sigma)/2$ in place of $\sigma$. Thus, we end up with
	$$\int_{B_\sigma(0)}\frac{\dist^2(x,\BC)}{|x|^{n+2-\mu}}\, \ext\|V\|(x) \leq C^\prime\int_{B_1(0)}\dist^2(x,\BC)\, \ext\|V\|(x)$$
	where $C^\prime = C^\prime(n,k,\sigma,\mu,\beta)$, which is exactly \eqref{E:initial-3}. This therefore completes the proof of Lemma \ref{lemma:L2-estimates-2} when $\sfrak(\BC)=n-1$.
\end{proof}

\begin{remark}\label{remark:L2-estimates-2-proof}
Up to changing some constants depending on the value of $\sfrak(\BC)$ (such as in \eqref{E:L2-2-q=1-1}), the only places where we used the fact that $\sfrak(\BC)=n-1$ in the above proof were in the applications of lemma \ref{lemma:Z}. Thus, once Lemma \ref{lemma:Z} has been established for a certain value of $\sfrak(\BC)$, the proof of Lemma \ref{lemma:L2-estimates-2} for that value of $\sfrak(\BC)$ follows in an analogous manner. We will therefore not repeat the proof of Lemma \ref{lemma:L2-estimates-2} when the time comes, and simply refer back to this remark.
\end{remark}

We are now in a position to prove Theorem \ref{thm:L2-estimates} when $\sfrak(\BC)=n-1$.

\begin{proof}[Proof of Theorem \ref{thm:L2-estimates} when $\sfrak(\BC)=n-1$]
	The key point is to show for $z\in \mathcal{D}(V)\cap \mathcal{T}^A_{\sigma r}(x_0,s_0)$ and $\rho\in (0,(1-\sigma)/2]$, the following is true: if $\eps_{II}$ and $\gamma_{II}$ are sufficiently small depending on $n,k,\sigma,a,M,\mu,\rho,\beta$, then the hypotheses of Lemma \ref{lemma:L2-estimates-2} hold for $\tilde{V}:= (\eta_{z,\rho r})_\#V \res B_1(0)$ and $\BC$ (note that $\tilde{V}$ still belongs to $\mathcal{V}_\beta$ as $\mathcal{V}_\beta$ is closed under this homothety).
	
	For this, notice that it is clear that $\Theta_{\tilde{V}}(0)\geq 2$ and $\nu(\BC)<\eps_0$ where $\eps_0$ is as in Lemma \ref{lemma:L2-estimates-2} (from the assumption on $z$ and (A2) of Hypothesis $H(A,\eps_{II},\gamma_{II},M)$, provided $\eps_{II}<\eps_0$). The mass ratio hypothesis for $\tilde{V}$ as in Lemma \ref{lemma:L2-estimates-2}(a) follows from Hypothesis $G(\eps_{II})$ and the monotonicity formula (again, provided $\eps_{II}$ is sufficiently small depending on $\sigma$). 
	
	Thus, it remains to check that $\hat{E}_{\tilde{V},\BC}<\eps_0$. Changing variables in the definition of $\hat{E}_{\tilde{V},\BC}$, we clearly have
	$$\hat{E}_{\tilde{V},\BC}^2 = (\rho r)^{-n-2}\int_{B_{\rho r}(z)}\dist^2(x,(\tau_z)_\#\BC)\, \ext\|V\|(x).$$
	Using the triangle inequality to compare the distance to $(\tau_z)_\#\BC$ to the distance to $\BC$, we therefore get
	$$\hat{E}_{\tilde{V},\BC}^2 \leq C\rho^{-n-2}\hat{E}_{V,\BC}^2(\mathcal{T}) + C\rho^{-2}\nu(\BC,(\tau_z)_\#\BC)^2$$
	where here $C = C(n,k)$ (note that we cannot remove the factor of $\rho^{-2}$ on the $\nu(\BC,(\tau_z)_\#\BC)$ term, as in the integral we are comparing the distance between $\BC$ and $(\tau_z)_\#\BC$ in $B_{\rho r}(z)$ and we only know that $z\in \mathcal{T}_{\sigma r}^A(x_0,s_0)$, and thus can only rescale by $r$). In light of Lemma \ref{lemma:nu-a}(i) and \eqref{E:translates-2}, this therefore gives
	$$\hat{E}_{\tilde{V},\BC}^2 \leq C\rho^{-n-2}\hat{E}^2_{V,\BC}(\mathcal{T}) + C\rho^{-2}\left(|\xi^{\perp_{P_0}}|^2 + \hat{E}_V^2(\mathcal{T})|\xi^{\top_{P_0}}|^2\right).$$
	Now choose $\eta>0$ so that $\rho^{-2}\eta<1$. Then, if $\eps_{II}$ and $\gamma_{II}$ are sufficiently small (depending on $\rho$) so that Corollary \ref{cor:weak-shift-estimate} holds with this choice of $\eta$, we therefore get
	$$\hat{E}_{\tilde{V},\BC}^2 \leq C\rho^{-n-2}\hat{E}_{V,\BC}^2(\mathcal{T}) + C\hat{E}_V^2(\mathcal{T}) \leq C(\rho^{-n-2}+1)\eps^2_{II}$$
	where in the second inequality here we have used the assumptions (A2) and (A4) of Hypothesis $H(A,\eps_{II},\gamma_{II},M)$. Hence, if we choose $\eps_{II}$ sufficiently small to make the last expression $<\eps_0$ from Lemma \ref{lemma:L2-estimates-2}, we have now shown that we may apply Lemma \ref{lemma:L2-estimates-2} to $\tilde{V}$ and $\BC$. From this, the first two estimates of Theorem \ref{thm:L2-estimates}, namely \eqref{E:L2-1} and \eqref{E:L2-2}, follow immediately.
	
	All that remains is to justify \eqref{E:L2-3}. From \eqref{E:initial-1} applied to $\tilde{V}$ (which we can do by the above), we have (with $\sigma$ replaced by $(1+\sigma)/2$ therein)
	$$\int_{B_{(1+\sigma)/2}(0)}\frac{|x^{\perp_{T_xV}}|^2}{|x|^{n+2}}\, \ext\|\tilde{V}\|(x) \leq C\hat{E}_{\tilde{V},\BC}^2.$$
	Changing variables to write this in terms of $V$, we get
	\begin{equation}\label{E:L2-q=1-1}
		\int_{B_{(1+\sigma)r\rho/2}(z)}\frac{|(x-z)^{\perp_{T_xV}}|^2}{|x-z|^{n+2}}\, \ext\|V\|(x) \leq C(r\rho)^{-n-2}\int_{B_{r\rho}(z)}\dist^2(x,(\tau_z)_\#\BC)\, \ext\|V\|(x).
	\end{equation}
	Write $\bar{x}:= \pi_{P_0}(x)$ and $\hat{x}:=x-\pi_{P_0}(x)$ for $x\in \R^{n+k}$. Looking at the left-hand side of the above inequality, for $x\in \spt\|V\|\cap \{|x^n|>\tau\}$ we then have the identity
	$$-R_{\bar{z}}(\bar{x})^2\frac{\del}{\del R_{\bar{z}}}\left(\frac{h^\alpha(\bar{x}) + u^\alpha(\bar{x})-\hat{z}}{R_{\bar{z}}(\bar{x})}\right)^{\perp_{T_xV}} = (x-z)^{\perp_{T_xV}}$$
	where $\alpha$ is such that $x\in \graph(h^\alpha+u^\alpha)$, and recall $R_{\bar{z}} = R_{\bar{z}}(x):= |x-\bar{z}|$; in particular, $x = \bar{x} + u^\alpha(\bar{x}) + h^\alpha(\bar{x})$. The proof of this is analogous to that seen in the proof of \eqref{E:B5-4} in the proof of $(\Bfrak5)$. Since $h^\alpha$ is linear on its domain and $D_ih^\alpha=0$ for $i<n$, we can write
	$$h^\alpha(\bar{x}) = (\bar{x}-\bar{z})\cdot Dh^\alpha(\bar{z}) + \xi^{\top_{P_0}}\cdot Dh^\alpha(\bar{z})$$
	where $\xi = z^{\perp_{S(\BC)}}$. Since $\frac{\del}{\del R_{\bar{z}}}\left(\frac{\bar{x}-\bar{z}}{R_{\bar{z}}(\bar{x})}\right)\equiv 0$, we therefore have that (as also $\xi^{\perp_{P_0}} = \hat{z}\equiv z^{\perp_{P_0}}$)
	$$(x-z)^{\perp_{T_xV}} = -R_{\bar{z}}(\bar{x})^2\frac{\del}{\del R_{\bar{z}}}\left(\frac{u^\alpha(\bar{x}) - (\xi^{\perp_{P_0}}-\xi^{\top_{P_0}}\cdot Dh^\alpha(\bar{z}))}{R_{\bar{z}}(\bar{x})}\right)^{\perp_{T_xV}}.$$
	Thus we have, using also \eqref{E:Lipschitz-approx-R1} and \eqref{E:J-bounds},
	\begin{align*}
	& \int_{B_{(1+\sigma)r\rho/2}(z)}\frac{|(x-z)^{\perp_{T_xV}}|^2}{|x-z|^{n+2}}\, \ext\|V\|(x) \geq \int_{(\R^k\times B^n_{\sigma \rho r}(\bar{z}))\cap \{|x^n|>\tau\}}\frac{|(x-z)^{\perp_{T_xV}}|^2}{|x-z|^{n+2}}\, \ext\|V\|(x)\\
	& \geq \sum_{\alpha=1}^4\int_{B^n_{\sigma\rho r}(\bar{z}))\cap \{|x^n|>\tau\}\cap U^\alpha_0}\frac{R_{\bar{z}}(x)^{n+2}}{|x+h^\alpha(x)+u^\alpha(x)-z|^{n+2}}R_{\bar{z}}(x)^{2-n}\\
	&\hspace{17em}\times\left|\frac{\del}{\del R_{\bar{z}}}\left(\frac{u^\alpha(x) - (\xi^{\perp_{P_0}}-\xi^{\top_{P_0}}\cdot Dh^\alpha(\bar{z}))}{R_{\bar{z}}(x)}\right)^{\perp_{T_xV}}\right|^2\, \ext x.
	\end{align*}
	where we have again used $x = \bar{x} + h^\alpha(\bar{x}) + u^\alpha(\bar{x})$ for $x\in\graph(h^\alpha+u^\alpha)$, for suitable $\bar{x}$. Now, as $|x+h^\alpha(x)+u^\alpha(x)-z| = |x-\bar{z} + (h^\alpha(x) + u^\alpha(x) - \hat{z})|$ for $x$ in the domain of $h^\alpha+u^\alpha$, ensuring that $R_{\bar{z}}(x)^{-1}|h^\alpha(x)+u^\alpha(x)-\hat{z}|$ is small by the choice of $\eps$ and $\gamma$ (depending on $\tau$; indeed, by Theorem \ref{thm:graphical-rep} we know that for sufficiently small $\eps$ and $\gamma$ we have $|\bar{z}|<\tau/2$, and so $R_{\bar{z}}(x)>\tau/2$ for $|x^n|>\tau$, and each of $h^\alpha(x), u^\alpha(x)$, and $\hat{z}$ approaches $0$ as $\eps,\gamma\downarrow 0$) means that we can insist that
	$$\frac{R_{\bar{z}}(x)^{n+2}}{|x+h^\alpha(x)+u^\alpha(x)-z|^{n+2}}\geq\frac{1}{2}$$
	in the latter integrand above. Substituting this into \eqref{E:L2-q=1-1} and using \eqref{E:Lipschitz-approx-R1}, we then deduce \eqref{E:L2-3}.
\end{proof}

The final step in the base case is to prove Theorem \ref{thm:shift-estimate} when $\sfrak(\BC)=n-1$.

\begin{proof}[Proof of Theorem \ref{thm:shift-estimate} when $\sfrak(\BC)=n-1$]
	We first claim the following: for every fixed $\sigma\in (1/2,1)$ and $\rho_0\in (0,(1-\sigma)/2]$, if $\eps_{III}$ and $\gamma_{III}$ are sufficiently small (depending on $\rho_0$ for the time being; later $\rho_0$ will be chosen to depend only on known parameters) then there are constants $\eta_1$ and $\tau_1$, depending only on $n$ and $k$ with the following property: for any $z\in \mathcal{D}(V)\cap \mathcal{T}^A_{\sigma r}(x_0,s_0)$ with $\xi:= z^{\perp_{S(\BC)}}\neq 0$, we have
	\begin{equation}\label{E:shift-q=1-1}
		\|V\|(A)\geq \eta_1(\rho_0 r)^n,
	\end{equation}
	where
	$$A := \Bigl\{x\in B_{\rho_0 r}(z)\cap \{r_{\BC}>\tau_1 r\}:|\xi^{\perp_{T_{x^\prime}\BC}}|\geq\eta_1 a^*(\BC)|\xi|\Bigr\},$$
	and where $x^\prime$ is the closest point projection of $x$ onto $\BC$. The proof of this claim is analogous to the proofs of the nearly identical claims in \cite[Lemma 6.21]{Wic04} and \cite{Wic14} (see also \cite[Lemma 3.9]{Sim93}). It relies on using Theorem \ref{thm:graphical-rep}, with $\tau$ therein depending on $n$, $k$ and $\rho_0$, which implies not only that $z$ must lie close to $S(\BC)$, but that for each half-plane $H^{\alpha}$ of $\BC$, $\spt\|V\| \cap B_{\rho_0 r}(z)$ contains plenty of points $x$ for which $x' \in H^{\alpha}$. 
	
	Since $|\xi^{\top_{P_0}}|\leq |\xi|$, we get from \eqref{E:shift-q=1-1} that
	\begin{equation}\label{E:shift-q=1-2}
	a^*(\BC)^2|\xi^{\top_{P_0}}|^2 \leq a^*(\BC)^2|\xi|^2 \leq \eta_1^{-3}(\rho_0 r)^{-n}\int_{A}|\xi^{\perp_{T_{x^\prime}\BC}}|^2\, \ext\|V\|(x).
	\end{equation}
	Now we will bound $|\xi^{\perp_{P_0}}|^2$ by the same integral as in the above. Indeed, for any $x\in \spt\|V\|\cap B_{\rho_0 r}(z)\cap \{r_{\BC}>\tau_1 r\}$ we have:
	\begin{equation}\label{E:shift-q=1-3}
	|\xi^{\perp_{P_0}}|^2 \leq 2|\xi^{\perp_{P_0}}-\xi^{\perp_{T_{x^\prime}\BC}}|^2 + 2|\xi^{\perp_{T_{x^\prime}\BC}}|^2 \leq C\|\pi_{P_0}-\pi_{T_{x^\prime}\BC}\|^2|\xi|^2 + 2|\xi^{\perp_{T_{x^\prime}\BC}}|^2
	\end{equation}
	where here $C = C(n,k)$. Now, since for all such $x$ Lemma \ref{lemma:nu-a} gives
	$$\|\pi_{P_0}-\pi_{T_{x^\prime}\BC}\| \leq C\nu(\BC) \leq C\hat{E}_{V}(\mathcal{T}) \leq Ca^*(\BC)$$
	when we integrate \eqref{E:shift-q=1-3} with respect to $\ext\|V\|$ over the region $A$, and use \eqref{E:shift-q=1-2} and \eqref{E:shift-q=1-1}, we obtain:
	\begin{align*}
		|\xi^{\perp_{P_0}}|^2 & \leq Ca^*(\BC)^2|\xi|^2 + C(\rho_0 r)^{-n}\int_{A}|\xi^{\perp_{T_{x^\prime}\BC}}|^2\, \ext\|V\|(x)  \leq C(\rho_0 r)^{-n}\int_{A}|\xi^{\perp_{T_{x^\prime}\BC}}|^2\, \ext\|V\|(x)
	\end{align*}
	where $C = C(n,k,\eta_1)$.
	
	So, currently we have
	\begin{equation}\label{E:shift-q=1-4}
		|\xi^{\perp_{P_0}}|^2 + a^*(\BC)^2|\xi^{\top_{P_0}}|^2 \leq C(\rho_0 r)^{-n}\int_{A}|\xi^{\perp_{T_{x^\prime}\BC}}|^2\, \ext\|V\|(x)
	\end{equation}
	where $C = C(n,k,\eta_1)$. We now aim to estimate this integral. By \eqref{E:translates-2}, Corollary \ref{cor:weak-shift-estimate}, and \eqref{E:nu-a-2}, we have
	$$\nu(\BC,(\tau_z)_\#\BC) \leq \eta a^*(\BC)$$
	for any choice of $\eta>0$, provided we choose the other constants sufficiently small. This means that by the definition of $A$, provided $\eps_{III}$ and $\gamma_{III}$ are sufficiently small depending on $\tau_1$, for every $x\in A$, the half-plane of $\BC$ closest to $x$ is the same as that which is closest to $x+z$, giving
	$$T_{x^\prime}\BC = T_{(x+z)^\prime}\BC.$$
	This observation therefore gives for $x\in A$,
	$$|\xi^{\perp_{T_{x^\prime}\BC}}| \leq |\xi^{\perp_{T_{x^\prime}\BC}}-(x-x^\prime)| + |x-x^\prime| = \dist(x,(\tau_z)_\#\BC) + \dist(x,\BC).$$
	So from this and \eqref{E:shift-q=1-4}, we get
	\begin{align*}
		|\xi^{\perp_{P_0}}|^2 & +a^*(\BC)^2|\xi^{\top_{P_0}}|^2\\
		& \leq C(\rho_0 r)^{-n}\left(\int_{B_{\rho_0 r}(z)\cap \{r>\tau_1 r\}}\dist^2(x,(\tau_z)_\#\BC)\, \ext\|V\| + \int_{B_{\rho_0 r}(z)\cap \{r>\tau_1 r\}}\dist^2(x,\BC)\, \ext\|V\|\right)\\
		& \leq C(\rho_0 r)^{-n}\int_{B_{\rho_0 r}(z)\cap \{r>\tau_1 r\}}\dist^2(x,(\tau_z)_\#\BC)\, \ext\|V\| + C\rho_0^{-n}r^2\hat{E}^2_{V,\BC}(\mathcal{T}).
	\end{align*}
	It remains to control the first term on the right-hand side of the above inequality. By choosing $\eps_{III}$ and $\gamma_{III}$ appropriately, we can use \eqref{E:L2-1} from Theorem \ref{thm:L2-estimates} with $\mu=1/4$ and with $\rho$ therein equal to $(1-\sigma)/2$. This then gives (assuming $\rho_0<\sigma(1-\sigma)/2$, as we may)
	\begin{align*}
		(\rho_0 r)^{-n}\int_{B_{\rho_0 r}(z)}\dist^2(x,(\tau_z)_\#\BC)\, \ext\|V\|(x) & \leq (\rho_0 r)^{7/4}\int_{B_{\rho_0 r}(z)}\frac{\dist^2(x,(\tau_z)_\#\BC)}{|x-z|^{n+7/4}}\, \ext\|V\|(x)\\
		& \leq C\rho_0^{7/4}r^{-n}\int_{B_{(1-\sigma)r/2}(z)}\dist^2(x,(\tau_z)_\#\BC)\, \ext\|V\|(x).
	\end{align*}
	where here $C = C(n,k,\sigma,a,M)$. Using now the triangle inequality with \eqref{E:translates-1} and \eqref{E:translates-2}, as well as Lemma \ref{lemma:nu-a}, we have
	\begin{align*}
	\rho_0^{7/4}r^{-n}\int_{B_{(1-\sigma)r/2}(z)}\dist^2(x,(\tau_z)_\#\BC)&\, \ext\|V\|(x)\\
	& \leq C\rho_0^{7/4}r^2\hat{E}^2_{V,\BC}(\mathcal{T}) + C\rho_0^{7/4}\left(|\xi^{\perp_{P_0}}|^2+a^*(\BC)^2|\xi^{\top_{P_0}}|^2\right)
	\end{align*}
	where here $C = C(n,k)$. Hence, if we choose $\rho_0$ sufficiently small depending on $n,k,\sigma,a,M$, we can combine the above inequalities and reabsorb the second term appearing on the right-hand side above into the left-hand side to get
	$$|\xi^{\perp_{P_0}}|^2 + a^*(\BC)^2|\xi^{\perp_{P_0}}|^2 \leq Cr^2\hat{E}_{V,\BC}^2(\mathcal{T})$$
	where $C = C(n,k,\sigma,a,M)$. Applying again Lemma \ref{lemma:nu-a}, namely \eqref{E:nu-a-2}, to replace $a^*(\BC)$ in the above by $\nu(\BC)$ up to a constant, this therefore completes the proof of Theorem \ref{thm:shift-estimate} when $\sfrak(\BC)=n-1$.
\end{proof}

\begin{remark}\label{remark:shift-estimate-proof}
The above proof of Theorem \ref{thm:shift-estimate} can be repeated in an analogous manner for general values of $\sfrak(\BC)$ once the other results have been established inductively. Because of this, we will not repeat the above proof in the more general inductive case when the time comes, and simply refer back to this remark.
\end{remark}

\section{Fine Blow-Ups Part I: Construction and Continuity}\label{sec:fine-blow-ups}

In this section, we construct a class of two-valued functions defined on $P_0$ that correspond to the linearisation of varifolds $V\in\mathcal{V}_\beta$ relative to a cone $\BC\in\CC$ under Hypothesis $G$ and $H$. These will be called \emph{fine blow-ups}. Understanding their properties will require knowing the validity of Theorem \ref{thm:graphical-rep}, \ref{thm:L2-estimates}, and \ref{thm:shift-estimate}. As such, we will work under an inductive hypothesis: we fix $\sfrak^*\in \{n-2,n-3,\dotsc,0\}$, and assume that Theorem \ref{thm:graphical-rep}, \ref{thm:L2-estimates}, and \ref{thm:shift-estimate} are known to hold when $\sfrak(\BC)> \sfrak^*$.

The following hypotheses describe the general situation in which we are able to construct fine blow-ups.

\begin{defn}\label{defn:dagger}
	We say that \emph{Hypothesis $\dagger$} holds when we have:
	\begin{itemize}
		\item A sequence $(V_j)_j\subset\mathcal{V}_\beta$ for some $\beta\in (0,1)$;
		\item A sequence $(\BC_j)_j\subset\CC$ such that $S(\BC_j)\equiv S$ is a fixed subspace with $\sfrak(\BC_j) \equiv \sfrak$ for all $j$, where $\sfrak\in \{n-1,\dotsc,\sfrak^*+1\}$;
		\item A subspace $A\subset S$ and regions $\mathcal{T}_j:= \mathcal{T}^A_{r_j}(x_j,s_j)$, where $r_j\to r_0$, $x_j\to x_0$, $s_j\to s_0$, where for some $a\in (0,1/2)$ we have $r_0/s_0\in [a,1-a]$ and furthermore $\mathcal{T}^A_{r_0}(x_0,s_0)\subset B_1(0)$. We also allow the case $\mathcal{T}_j\equiv B_1(0)$;
		\item A relatively closed set $\mathcal{D}\subset S\cap \mathcal{T}^A_{r_0}(x_0,s_0)$ such that $\mathcal{D}(V_j)\to \mathcal{D}$ locally in Hausdorff distance in $\mathcal{T}^A_{r_0}(x_0,s_0)$;
		\item Sequences of numbers $\eps_j,\gamma_j\downarrow 0$, and some $M\in [1,\infty)$ for which $V_j$ and $\BC_j$ satisfy Hypothesis $G(\eps_j)$ and Hypothesis $H(A,\eps_j,\gamma_j,M)$ in $\mathcal{T}_j$.
	\end{itemize}
\end{defn}

Furthermore, when $\sfrak=n-1$ and $\mathcal{D}\neq S\cap \mathcal{T}^A_{r_0}(x_0,s_0)$, we make the additional stipulation that
\begin{equation}\label{E:extra}
	\BC_j\in\mathcal{P} \text{ for all $j$.}
\end{equation}
The technical reason for \eqref{E:extra} is that when there is a gap, it will be convenient to use our topological structural condition via Lemma \ref{lemma:gap-2} to write $V_j$ as a graph over the \emph{planes} in $\BC_j$; if $\BC_j$ was formed of half-planes, this would not be possible. In fact, there is an alternative, equivalent assumption to \eqref{E:extra} (which one might find more natural in the context of the present work), namely that if $\BC_j\in \CC_{n-1}\setminus\mathcal{P}_{n-1}$, then instead of \eqref{E:extra} we add the assumption to Hypothesis $\dagger$ that
\begin{equation}\label{E:extra-2}
\left(\inf_{\Dbf\in \mathcal{P}_{n-1}}\hat{F}_{V_j,\Dbf}(\mathcal{T}_j)\right)^{-1}\hat{F}_{V_j,\BC_j}(\mathcal{T}_j) \to 0
\end{equation}
i.e.~we only ever form a fine blow-up off a union of $4$ half-planes when the fine excess relative to a twisted cone is significantly smaller than the fine excess relative to any union of planes with axis dimension $n-1$. Indeed, under this assumption, if we had $\mathcal{D}\neq S\cap \mathcal{T}^A_{r_0}(x_0,s_0)$, then any blow-up (fine or coarse) would need to be the union of two planes, rather than four half-planes, as the existence of a gap guarantees smoothness of the blow-up along the axis, ruling out $4$ half-planes as a possible blow-up. But taking a fine blow-up off a suitable choice of planes in $\mathcal{P}_{n-1}$ would contradict \eqref{E:extra-2}; thus, under \eqref{E:extra-2} it follows that $\mathcal{D} = S\cap \mathcal{T}^A_{r_0}(x_0,s_0)$. Moreover, if \eqref{E:extra-2} fails, we get $\hat{F}_{V_j,\BC_j}\geq C\inf_{\Dbf\in \mathcal{P}_{n-1}}\hat{F}_{V_j,\Dbf}(\mathcal{T}_j)$ for suitable $C = C(n,k,a,M)\in (0,\infty)$, and so we may instead replace $\BC_j$ with $\Dbf_j\in \mathcal{P}_{n-1}$. Of course, morally speaking, the reason for all of this is because if the $\BC_j$ are twisted cones, then all nearby cones are also twisted (as the spine cannot split due to the structure), and so one expects to have points of density $\geq 2$ converging all along $S$. As such, one could handle the case when $\BC_j$ is twisted somewhat independently of the case where $\BC_j\in \mathcal{P}$.

\begin{remark} The fine blow-up class will technically depend on the choice of the parameter $M\in [1,\infty)$, however we will suppress this dependence. We stress that the parameter $M$ is of importance, however, as the fine blow-up class for a given $M$ is \emph{not} necessarily closed under certain natural operations (for instance, those analogous to $(\mathfrak{B}3\text{I})$), but its closure under these operations \emph{will} be contained in the fine blow-up class for a fixed, larger choice of $M$.

Similarly, the construction of fine blow-ups has a dependence on the parameter $\beta\in (0,1)$. Later on when proving Theorem \ref{thm:fine-reg} by iterating an excess decay statement, we will need to perform small rotations to elements of $\mathcal{V}_\beta$, which as seen in Proposition \ref{prop:rotation-class} causes the value of $\beta$ to change by a controlled amount. Since all these rotations will be controlled, these will be contained within the fine blow-up class for a fixed, smaller choice of $\beta$. Since this still obeys all the same properties, this does not impact the argument. Alternatively, one could define the fine blow-up class for a given choice of $\beta$ by blowing-up sequences of the form $(\Gamma_j)_\#V_j$ in place of $V_j$, where $\Gamma_j\in SO(n+k)$ have $\Gamma_j\to \id$.
\end{remark}

We are now in a position to construct fine blow-ups. Suppose that Hypothesis $\dagger$ are satisfied. Write $\tilde{E}_j:= \hat{E}_{V_j,\BC_j}(\mathcal{T}_j)$ and $\hat{E}_j:= \hat{E}_{V_j}(\mathcal{T}_j)$. For simplicity, we sometimes drop the index $j$.

Now, if we fix any $\tau\in (0,1/8)$, $\sigma\in (1/2,1)$, and $\rho\in (0,1/4)$, we know that for sufficiently large $j$ we have
$$\eps_j<\min\{\eps_I(n,k,\tau,\sigma,a,M,\beta),\eps_{II}(n,k,\sigma,a,M,\rho,\beta),\eps_{III}(n,k,\sigma,a,M,\beta)\}$$
$$\gamma_j<\min\{\gamma_I(n,k,\tau,\sigma,a,M,\beta),\gamma_{II}(n,k,\sigma,a,M,\rho,\beta),\gamma_{III}(n,k,\sigma,a,M,\beta)\}$$
(with $\mu=1/4$ in Theorem \ref{thm:L2-estimates}). This means that for all sufficiently large $j$, the hypotheses of Theorems \ref{thm:graphical-rep}, \ref{thm:L2-estimates}, and \ref{thm:shift-estimate} will all be satisfied with these chosen parameters $\tau$, $\sigma$, and $\rho$ (and with $\mu = 1/4$ in Theorem \ref{thm:L2-estimates}). In light of the base case if $\sfrak=n-1$ from Section \ref{sec:base-case} or by appealing to the inductive hypothesis if $\sfrak < n-1$, we therefore have that all of the conclusions of Lemma \ref{lemma:nu-a}, Theorem \ref{thm:graphical-rep}, Theorem \ref{thm:L2-estimates}, Theorem \ref{thm:shift-estimate}, and Corollary \ref{cor:L2-corollary} hold with these parameters (of course, we need to also assume that $\eps_j$ and $\gamma_j$ are also smaller than the corresponding parameters in Corollary \ref{cor:L2-corollary}).

In light of Lemma \ref{lemma:nu-a}(i) and Remark \ref{remark:elementary-observations}(1) from Section \ref{sec:further-notation}, we can immediately pass to a subsequence along which we have the following limits:
\begin{itemize}
	\item if $\sfrak = n-1$, set $m^\alpha := \lim_{j\to\infty}\hat{E}_j^{-1}Dh^\alpha_j$ for $\alpha=1,\dotsc,4$;
	\item if $\sfrak < n-1$, set $m^\alpha := \lim_{j\to\infty}\hat{E}_j^{-1}Dp^\alpha_j$ for $\alpha = 1,2$.
\end{itemize}
Notice that when $\sfrak=n-1$, $m^\alpha$ is a $k\times n$ matrix, with zeroes in the first $(n-1)$-columns.

There are two slightly different situations to consider depending on whether or not points of density $\geq 2$ accumulate along all of $S_0\cap \mathcal{T}^A_{r_0}(x_0,s_0)$. When they do, the fine blow-up is determined by $4$ functions defined on subsets of half-planes, whilst when not the fine blow-up is determined by two functions defined on a subset of the full plane $P_0$. Of course, this accumulation can only happen along all of $S_0\cap \mathcal{T}_{r_0}^A(x_0,s_0)$ when $\sfrak=n-1$, and in particular, due to \eqref{E:extra} (or \eqref{E:extra-2}), it always happens when $\BC_j\in \CC\setminus \mathcal{P}$. Let us start by constructing blow-ups when this accumulation happens. 

So first suppose $\sfrak=n-1$ (so $S \equiv S_0$) and $\mathcal{D} = S\cap \mathcal{T}^A_{r_0}(x_0,s_0)$. Then for $\alpha=1,2$, extend the functions $u_j^\alpha$ given by Theorem \ref{thm:graphical-rep} to all of $P_0\cap \mathcal{T}^{A}_{r_j}(x_j,s_j)$ by defining their values to be zero on $P_0\cap \mathcal{T}_{r_j}^A(x_j,s_j)\cap \{-\tau r_j<x^n<0\}$, and similarly for $\alpha=3,4$ extend $u^\alpha_j$ to be zero on $P_0\cap \mathcal{T}^A_{r_j}(x_j,s_j)\cap \{0<x^n<\tau r_j\}$. For any compact subset $K\subset U^\alpha_0\cap \mathcal{T}^A_{r_0}(x_0,s_0)$, by first choosing $\tau$ and $\sigma$ depending on $K$, and then choosing $j$ sufficiently large, we have $K\subset U^\alpha_0\cap \mathcal{T}^A_{\sigma r_j}(x_j,s_j)\cap \{|x^n|>\tau r_j\}$. Thus, using Theorem \ref{thm:graphical-rep}, elliptic estimates, along with a compact exhaustion of $P_0\cap \mathcal{T}^A_{r_0}(x_0,s_0)$ and a diagonalisation argument, we deduce that there exist harmonic functions $\phi^\alpha:U^\alpha_0\cap \mathcal{T}^A_{r_0}(x_0,s_0)\to P_0^\perp$ for $\alpha=1,\dotsc,4$, such that
$$\tilde{E}_{j}^{-1}u^\alpha_j\to \phi^\alpha \qquad \text{as }j\to\infty,$$
where the convergence is in $C^3_{\text{loc}}(U^\alpha_0\cap \mathcal{T}^A_{r_0}(x_0,s_0)\setminus\{x^n=0\};P_0^\perp)$. Then, using the $L^2$ non-concentration inequality \eqref{E:cor-3} and the fact that $\mathcal{D} = S\cap \mathcal{T}^A_{r_0}(x_0,s_0)$, we deduce that for any $\eta\in (0,\sigma(1-\sigma)/4)$, for sufficiently large $j$ depending on $\eta$, we have
$$\int_{U^\alpha_0\cap (\mathcal{T}^A_{\sigma r_0}(x_0,s_0)\cap S)_{\eta r_0}}|\tilde{E}^{-1}_ju^\alpha_j|^2\, \ext\H^n \leq C\eta^{1/2}$$
where $C = C(n,k,\sigma,a,M)$. From this, combined with the convergence away from $S\equiv S_0\equiv \{x\in P_0 : x^n=0\}$, we deduce that the convergence $\tilde{E}_j^{-1}u_j^\alpha\to \phi^\alpha$ is in fact also in $L^2(U^\alpha_0\cap \mathcal{T}^A_{\sigma r_0}(x_0,s_0);P_0^\perp)$ for every $\sigma\in (1/2,1)$.

\begin{defn}
	When $\sfrak=n-1$ and $\mathcal{D} = S\cap \mathcal{T}^A_{r_0}(x_0,s_0)$, any collection of functions $\{\phi^\alpha\}_{\alpha=1}^4$ defined on $U^\alpha_0\cap \mathcal{T}^A_{r_0}(x_0,s_0)$ for each $\alpha=1,\dotsc,4$ constructed in the above manner is called a \emph{fine blow-up} of $(V_j)_j$ off $P_0$ relative to $(\BC_j)_j$ in $\mathcal{T}^A_{r_0}(x_0,s_0)$.
\end{defn}

Now we consider the second case, namely when either $\sfrak<n-1$ or we have $\sfrak=n-1$ and $\mathcal{D}\subsetneq S\cap \mathcal{T}^A_{r_0}(x_0,s_0)$. Write
$$r_{\mathcal{D}}(x):=\dist(x,\mathcal{D}) \qquad \text{and} \qquad \BC_j = |P^1_j| + |P^2_j|.$$
For all sufficiently large $j$, we have $\dist_\H(\mathcal{D}(V_j)\cap \mathcal{T}^A_{\sigma r_j}(x_j,s_j), \mathcal{D})<\tau r_0/2$. This means that for any fixed $y\in \mathcal{T}_{r_0}(x_0,s_0)\cap \{r_{\mathcal{D}}>2\tau r_0\}$, then for all sufficiently large $j$ (depending on $y$) we have
$$C_{\tau r_0}(y)\cap \mathcal{D}(V_j) = \emptyset.$$
By then appealing to Lemma \ref{lemma:gap-2}, we then get that $V^j\res C_{\tau r_0/8}(y) = V_j^1 + V_j^2$ for each such $j$, where each $V^\alpha_j$ is a (single-valued) smooth minimal graph and where $V^\alpha_j\weakly |\{0\}^k\times B^n_{\tau r_0/8}(y)|$ as varifolds, for $\alpha=1,2$. Now cover the set $\{r_{\BC}<\tau r_0/16\}\cap \{r_{\mathcal{D}}>3\tau r_0\}\cap \mathcal{T}^A_{\sigma r_0}(x_0,s_0)$ by a finite collection of balls $\{B^n_{\tau r_0/8}(y_p)\}_{p=1}^L$, where $L = L(n,\sigma)$ and $y_p\in \mathcal{T}^A_{r_0}(x_0,s_0)\cap \{r_{\mathcal{D}}>2\tau r_0\}$ for each $p$. Then, for all sufficiently large $j$ the above argument applies, and so applying Allard's regularity theorem on each of the single-valued minimal graphs just obtained, we conclude that for sufficiently large $j$ there is a pair $u^1_j, u^2_j$ of smooth functions $P_0\cap \{r_{\mathcal{D}}>3\tau r_0\}\cap \mathcal{T}^A_{\sigma r_0}(x_0,s_0)\to P_0^\perp$ such that
\begin{align*}
	V_j\res (\{x:r_{\mathcal{D}}(\pi_{P_0}(x)&)>3\tau r_0\}\cap \mathcal{T}^A_{\sigma r_0}(x_0,s_0))\\
	& = \sum_{\alpha=1,2}|\graph(u^\alpha_j + p^\alpha_j)|\res (\{x:r_{\mathcal{D}}(\pi_{P_0}(x))>3\tau r_0\}\cap \mathcal{T}^A_{\sigma r_0}(x_0,s_0)).
\end{align*}
Then, as before, we may extend each $u_j^\alpha$ to equal zero on $\{x\in P_0:r_{\mathcal{D}}\leq 3\tau r_0\}\cap \mathcal{T}^A_{\sigma r_0}(x_0,s_0)$. Then, using elliptic estimates and Lemma \ref{lemma:appendix-2} to bound the excess relative to a single plane $P^\alpha_j$ by the fine excess in the whole ball (as initially our estimates for $u^\alpha_j$ are in terms of their $L^2$ norm to a single plane, not the $L^2$ norm to $\BC_j$, which could be smaller in the region close to $S$), and a diagonal argument (letting $\tau\downarrow 0$, $\sigma\uparrow 1$) we have that there exist harmonic functions $\phi^1,\phi^2:\{r_{\mathcal{D}}>0\}\cap \mathcal{T}^A_{r_0}(x_0,s_0)\to P_0^\perp$ such that for $\alpha=1,2$,
$$\tilde{E}_j^{-1}u_j^\alpha\to \phi^\alpha \qquad \text{as }j\to\infty$$
where the convergence is in $C^3_{\text{loc}}(P_0\cap \mathcal{T}^A_{r_0}(x_0,s_0)\cap \mathcal{D}^c;P_0^\perp)$. We can then again appeal to \eqref{E:cor-3} to deduce that the convergence is also in $L^2(P_0\cap \mathcal{T}^A_{\sigma r_0}(x_0,s_0)\cap \mathcal{D}^c;P_0^\perp)$ for each $\sigma\in (1/2,1)$. 

\textbf{Remark:} We stress that when we say here that we are invoking elliptic estimates and Lemma \ref{lemma:appendix-2}, what we mean is that the any $C^{3,\gamma}$ estimate for $u^\alpha_j$ which we know holds away from the spine by Theorem \ref{thm:graphical-rep} can be extended across the spine (up to constants), as elliptic estimates for $u^\alpha_j$ allow one to control the $C^{3,\gamma}$ norm $u_j^\alpha$ on a ball by the $L^2$ norm on an annulus. This is what allows us to verify the assumptions of Lemma \ref{lemma:appendix-2}, to control the $C^{3,\gamma}$ norm of $u^\alpha_j$ by $\tilde{E}_j$, and thus gives the $C^3_{\text{loc}}$ convergence claimed above.

We now define:
\begin{defn}
	When $\sfrak<n-1$ or $\sfrak=n-1$ and $\mathcal{D}\subsetneq S\cap \mathcal{T}^A_{r_0}(x_0,s_0)$, any collection of functions $\phi^1,\phi^2:P_0\cap \mathcal{T}^A_{r_0}(x_0,s_0)\cap \mathcal{D}^c\to P_0^\perp$ constructed in the manner above is called a \emph{fine blow-up} of $(V_j)_j$ off $P_0$ relative to $(\BC_j)_j$ in $\mathcal{T}^A_{r_0}(x_0,s_0)$.
\end{defn}

\textbf{Remarks:}
\begin{enumerate}
	\item [(1)] Of course if $\sfrak=n-1$ and $\mathcal{D}\subsetneq S\cap \mathcal{T}^A_{r_0}(x_0,s_0)$, we can always define $\tilde{\phi}^\alpha:= \left.\phi^1\right|_{U^\alpha_0}$ for $\alpha=1,2$ and $\tilde{\phi}^\alpha := \left.\phi^2\right|_{U^\alpha_0}$ for $\alpha=3,4$, so that in either of the two cases where $\sfrak=n-1$, a fine blow-up is determined by the information of four smooth harmonic functions in the relevant regions, two on each of the half-planes $\{x^n<0\}$, $\{x^n>0\}$. These will then `pair up' across the boundary portion $\{x^n=0\}\setminus\mathcal{D}$ to form two smooth harmonic functions. We can similarly do this even in the $\sfrak<n-1$ case, and sometimes we will use this liberty to switch between an index set of $\{1,2\}$ (for functions defined on planes) and of $\{1,2,3,4\}$ (for functions defined on half-planes).
	\item [(2)] In light of (1), we will introduce the notation $\{\phi^\alpha\}_{\alpha\in I}$ for fine blow-ups, where it may not always be specified whether we are referring to a fine blow-up with $I = \{1,\dotsc,4\}$ or one with $I = \{1,2\}$.
	\item [(3)] If $\{\phi^\alpha\}_{\alpha\in I}$ is a fine blow-up of $(V_j)_j$ relative to a sequence $(\BC_j)_j$, we write
	$$\mathcal{D}(\phi):= \mathcal{D}.$$
	Of course, it should be noted that the set $\mathcal{D}$ is really determined by the sequence $(V_j)_j$, which as we passed to a subsequence does not necessarily uniquely determine the fine blow-up. We also write
	$$\sfrak(\phi):= \sfrak(\BC_j) \equiv \sfrak.$$
	\item [(4)] (Local boundedness of fine blow-ups.) Note that the elliptic estimates we used in the construction of fine blow-ups immediately gives that for any compact set $K\subset \mathcal{T}^A_{r_0}(x_0,s_0)\setminus\mathcal{D}$, there is a constant $C = C(n,k,K)$ such that $|\phi^\alpha|_{C^3(K)}\leq C$. Moreover, the argument based on the $L^2$ non-concentration estimate \eqref{E:cor-3} gives that for any $\sigma\in (0,1)$, there is a constant $C = C(n,k,\sigma)\in (0,\infty)$ such that
	$$\|\phi^\alpha\|_{L^2(P_0\cap \mathcal{T}^A_{\sigma r}(x_0,s_0))}\leq C.$$
\end{enumerate}

The main set of fine blow-ups we will be interested in for the proof of the fine $\eps$-regularity theorem are those defined on the unit ball. Thus, we define for each $\sfrak\in \{n-1,\dotsc,0\}$,

\begin{defn}
	Write $\mathfrak{b}_\sfrak$ for the collection of all fine blow-ups $\phi = \{\phi^\alpha\}_{\alpha\in I}$ with $\sfrak(\phi) = \sfrak$ and $0\in\mathcal{D}(\phi)$, and where $\mathcal{T}_j = B_1(0)$ for all $j$ and $A \equiv S$ ($\equiv S(\BC_j)$). We will also write $\mathfrak{b}:=\cup_{\sfrak=0}^{n-1}\mathfrak{b}_{\sfrak}$.
\end{defn}

A particularly important sub-class later will be the set of linear fine blow-ups. We introduce notation for this now.

\begin{defn}
	Let $\mathfrak{l}\subset\mathfrak{b}$ denote all fine blow-ups $\phi = \{\phi^\alpha\}_{\alpha\in I}$ of the following form: either $I = \{1,2\}$ and $\phi^1,\phi^2$ are both linear functions; or $I = \{1,\dotsc,4\}$ and $\phi^\alpha$ is linear on its domain for each $\alpha\in \{1,\dotsc,4\}$, and there is a linear function $\lambda:S_0\to S_0^\perp$ with $\lambda^\alpha(z) = \lambda^{\perp_{P_0}}(z) - m^\alpha\cdot\lambda(z)^{\top_{P_0}}$ for all $z\in S_0$.
\end{defn}

\textbf{Note:} When $\phi \in\mathfrak{l}$ is such that $\phi\in\mathfrak{b}_{n-1}$ and $\mathcal{D}(\phi)= S_0\cap B_1(0)$, then both index sets $I = \{1,2\}$ or $I = \{1,\dotsc,4\}$ are potentially valid, and so we allow both possibilities for $\phi\in\mathfrak{l}$.

An important remark is that each $\phi\in \mathfrak{l}$ is generated from the fine blow-up of a sequences of cones $\widetilde{\BC}_j$ relative to $\BC_j$.

\subsection{Continuity of Fine Blow-Ups}

To continue the proof of the induction step for the main technical theorems in Section \ref{sec:fine-L2}, we need to first observe some basic properties of fine blow-ups $\phi$ at points in $\mathcal{D}(\phi)$. In particular, we can establish a Hölder continuity estimate, as well as their values on $\mathcal{D}(\phi)$; in particular, they are continuous across $\mathcal{D}(\phi)$ always. When $\sfrak<n-1$, this will in fact suffice for establishing full regularity of fine blow-ups.

The key estimate for this continuity conclusion is the following lemma. For simplicity, here we use notation that views fine blow-ups as a quadruple of functions over half-planes via restriction onto $U^\alpha_0$, even if they can be viewed as two functions over planes (and similarly for the constants $m^\alpha$).

\begin{lemma}\label{lemma:fine-continuity-1}
	Let $\{\phi^\alpha\}_{\alpha\in I}$ be a fine blow-up in $\mathcal{T}^A_{r_0}(x_0,s_0)$ (where recall if $s_0\neq 0$, $r_0/s_0\in [a,1-a]$). Then, for every $z\in \mathcal{D}(\phi)\cap \mathcal{T}^A_{7r_0/8}(x_0,s_0)$, there exists $\lambda(z)\in S^\perp$ with $|\lambda(z)|\leq C$ and such that for every $\rho\in (0,r_0/8)$ we have
	\begin{align*}
		\sum_{\alpha\in I}\int_{B^n_{\rho/2}(z)\cap U^\alpha_0}&\frac{|\phi^\alpha(x)-(\lambda(z)^{\perp_{P_0}} - m^\alpha\cdot\lambda(z)^{\top_{P_0}})|^2}{|x-z|^{n+7/4}}\, \ext x\\
		& \hspace{2em} \leq C\sum_{\alpha\in I}\rho^{-n-7/4}\int_{B^n_\rho(z)}|\phi^\alpha(x)-(\lambda(z)^{\perp_{P_0}} - m^\alpha\cdot\lambda(z)^{\top_{P_0}})|^2\, \ext x.
	\end{align*}
	Here, $C = C(n,k,a,M)\in (0,\infty)$, and the notation $m^\alpha\cdot \lambda(z)^{\top_{P_0}}$ means multiplying the vector $\lambda(z)^{\top_{P_0}}$ by the matrix $m^\alpha$ to give a vector in $P_0^\perp$.
\end{lemma}

\textbf{Remark:} The function $\lambda:\mathcal{D}(\phi)\to S^\perp$ whose existence is asserted in Lemma \ref{lemma:fine-continuity-1} will play a key role in Section \ref{sec:fine-blow-ups-2}. If we want to stress its dependence on $\phi$, we write $\lambda\equiv \lambda_\phi$. We note that whilst Lemma \ref{lemma:fine-continuity-1} only gives that $\lambda$ is defined on $\mathcal{D}(\phi)\cap \mathcal{T}^A_{7r_0/8}(x_0,s_0)$, the argument in Lemma \ref{lemma:fine-continuity-1} shows that it is possible to define it at all $z\in \mathcal{D}(\phi)$. The constant $C$ then depends on the ratio of $r_0$ and the distance of $z$ from the boundary of $\mathcal{T}^A_{r_0}(x_0,s_0)$.

\begin{proof}
	Fix $z\in \mathcal{D}(\phi)\cap \mathcal{T}^A_{7r_0/8}(x_0,s_0)$ and $\rho\in (0,r_0/8)$. Suppose that $z_j\in\mathcal{D}(V_j)$ is such that $z_j\to z$. Pick $\rho^\prime\in (2\rho/3,\rho)$ so that for sufficiently large $j$ (depending on $\rho^\prime$) we have
	$$B_{\rho/2}(z)\subset B_{3\rho^\prime/4}(z_j)\subset B_{\rho^\prime}(z_j)\subset B_\rho(z).$$
	Write $\xi_j:= z_j^{\perp_S}$, and fix $\tau>0$. Now apply \eqref{E:L2-1} of Theorem \ref{thm:L2-estimates} with $z_j$ in place of $z$, $\sigma=3/4$, $\mu=1/4$, and the $\rho^\prime/r_0$ in place of $\rho$ therein (which is a valid choice). This gives
	$$\int_{B_{3\rho^\prime/4}(z_j)}\frac{\dist^2(x,(\tau_{z_j})_\#\BC_j)}{|x-z_j|^{n+7/4}}\, \ext\|V_j\|(x) \leq C(\rho^\prime)^{-n-7/4}\int_{B_{\rho^\prime}(z_j)}\dist^2(x,(\tau_{z_j})_\#\BC_j)\, \ext\|V_j\|(x)$$
	where $C = C(n,k,a,M)\in (0,\infty)$. If we now use the area formula to write the integral on the left-hand side in terms of the graphs $u^\alpha_j$ in the region $\{r_{\mathcal{D}}>\tau\}$ and bound the Jacobian below by $1$ (cf.~\eqref{E:J-bounds}), this gives for all $j$ sufficiently large
	\begin{align*}
	\sum_{\alpha\in I}\int_{\{r_{\mathcal{D}}>\tau\}\cap B_{\rho/2}(z)\cap U^\alpha_0}&\frac{|u_j^\alpha(x) - (\xi_j^{\perp_{P_0}}-Dp_j^\alpha\cdot\xi_j^{\top_{P_0}})|^2}{|x+u^\alpha_j(x)+p^\alpha_j(x)-z_j|^{n+7/4}}\, \ext x\\
	&\hspace{8em} \leq C(\rho^\prime)^{-n-7/4}\int_{B_{\rho^\prime}(z_j)}\dist^2(x,(\tau_{z_j})_\#\BC_j)\, \ext\|V_j\|(x)
	\end{align*}
	(indeed, the expression for $\dist(x,(\tau_{z_j})_\#\BC_j)$ in the region $\{r_{\mathcal{D}}>\tau\}$ follows from Theorem \ref{thm:shift-estimate} and Lemma \ref{lemma:nu-a}(ii), (iv), coupled with Hypothesis (A4)). Now we split the domain of integration in the integral on the right-hand side above into the complementary regions $B_{\rho^\prime}(z_j)\cap \{r_{\mathcal{D}}\leq \tau\}$ and $B_{\rho^\prime}(z_j)\cap \{r_{\mathcal{D}}>\tau\}$. For the integral over the latter domain, we again have the integral is controlled by (after using the area formula and controlling the Jacobian by a constant)
	$$C\sum_{\alpha\in I}\rho^{-n-7/4}\int_{\{r_{\mathcal{D}}>\tau\}\cap B_{\rho}(z)\cap U^\alpha_0}|u_j^\alpha(x)-(\xi_j^{\perp_{P_0}}-Dp_j^\alpha\cdot \xi_j^{\top_{P_0}})|^2\, \ext x.$$
	Here, we have also used $B_{\rho^\prime}(z_j)\subset B_\rho(z)$. For the other term, using \eqref{E:cor-2} and \eqref{E:cor-3}, the fact that $\|V_j\|(B_\rho(z)\cap \{r_{\mathcal{D}}\leq\tau\})<C\tau$ for suitable $C = C(n)$, we get
	\begin{align*}
		&\int_{B_{\rho^\prime}(z_j)\cap \{r_{\mathcal{D}}\leq \tau\}}\dist^2(x,(\tau_{z_j})_\#\BC_j)\, \ext\|V_j\|(x) \\
		& \hspace{3em} \leq \int_{B_{\rho}(z)\cap \{r_{\mathcal{D}}\leq \tau\}}\dist^2(x,\BC_j)\, \ext\|V_j\|(x) + C\hat{E}^2_{V_j,\BC_j}(\mathcal{T}_j)\cdot\|V_j\|(B_\rho(z)\cap \{r_{\mathcal{D}}\leq\tau\})\\
		& \hspace{3em} \leq C^\prime\tau^{1/2}\hat{E}^2_{V_j,\BC_j}(\mathcal{T}_j) + C\hat{E}_{V_j,\BC_j}^2(\mathcal{T}_j)\cdot\tau\\
		& \hspace{3em} \leq C^\prime\tau^{1/2}\tilde{E}_j^2
	\end{align*}
	where here $C^\prime$ depends on $n,k,\rho,a,M$ (although this term will be vanishing in the limit $\tau\downarrow 0$, so the dependence on $\rho$ here will not show up at the end of the proof). So, combining we currently have
	\begin{align}
		\nonumber\sum_{\alpha\in I}&\int_{\{r_{\mathcal{D}}>\tau\}\cap B_{\rho/2}(z)\cap U^\alpha_0}\frac{|u_j^\alpha(x) - (\xi_j^{\perp_{P_0}}-Dp_j^\alpha\cdot\xi_j^{\top_{P_0}})|^2}{|x+u^\alpha_j(x)+p^\alpha_j(x)-z_j|^{n+7/4}}\, \ext x\\
		& \hspace{2em} \leq C\sum_{\alpha\in I}\rho^{-n-7/4}\int_{\{r_{\mathcal{D}}>\tau\}\cap B_{\rho}(z)\cap U^\alpha_0}|u_j^\alpha(x)-(\xi_j^{\perp_{P_0}}-Dp_j^\alpha\cdot \xi_j^{\top_{P_0}})|^2\, \ext x + C^\prime\tau^{1/2}\tilde{E}_j^2.\label{E:fine-continuity-1-1}
	\end{align}
	We also know from Theorem \ref{thm:shift-estimate} (and Lemma \ref{lemma:nu-a}(i)) that $\tilde{E}_j^{-1}\xi_j^{\perp_{P_0}}$ and $\tilde{E}_j^{-1}\hat{E}_j\xi_j^{\top_{P_0}}$ are bounded sequences. Therefore, passing to a subsequence, we conclude that there exists $\lambda(z)\in S^\perp$ for which $|\lambda(z)|\leq C = C(n,k,a,M)$ and
	$$\tilde{E}_j^{-1}\xi_j^{\perp_{P_0}} \to \lambda(z)^{\perp_{P_0}} \qquad \text{and} \qquad \tilde{E}_j^{-1}\hat{E}_j\xi_j^{\top_{P_0}}\to \lambda(z)^{\top_{P_0}}.$$
	(Indeed, if $\lambda_1,\lambda_2$ are limits of these sequences, we set $\lambda(z) := (\lambda_1,\lambda_2)\in (P_0^\perp\times P_0^\top)\cap S^\perp$.)
	
	Now, if we divide \eqref{E:fine-continuity-1-1} by $\tilde{E}_j^2$ and take $j\to\infty$, then using that $\tilde{E}_j^{-1}u^\alpha_j\to \phi^\alpha$ in $L^2$, $\hat{E}_j^{-1}Dp_j^\alpha\to m^\alpha$, and the dominated convergence theorem (as $p_j^\alpha$ and $u^\alpha_j$ as converge to $0$ in $C^2$ locally away from $\mathcal{D}(\phi)$), we get
	\begin{align*}
		\sum_{\alpha\in I}&\int_{\{r_{\mathcal{D}}>\tau\}\cap B_{\rho/2}(z)\cap U^\alpha_0}\frac{|\phi^\alpha(x)-(\lambda(z)^{\perp_{P_0}} - m^\alpha\cdot \lambda(z)^{\top_{P_0}})|^2}{|x-z|^{n+7/4}}\, \ext x\\
		& \leq C\sum_{\alpha\in I}\rho^{-n-7/4}\int_{\{r_{\mathcal{D}}>\tau\}\cap B_\rho(z)\cap U^\alpha_0}|\phi^\alpha-(\lambda(z)^{\perp_{P_0}} - m^\alpha\cdot \lambda(z)^{\top_{P_0}})|^2\, \ext x + C^\prime \tau^{1/2}.
	\end{align*}
	Now if we let $\tau\downarrow 0$, we arrive at the desired conclusion.
\end{proof}

\begin{remark}\label{remark:fine-blow-up-is-0}
	Notice that if the points $z_j\in\mathcal{D}(V_j)$ were all on the spine $S(\BC_j)\equiv S$ to begin with, then $\xi_j := z^{\perp_S}_j \equiv 0$, and so we would get that $\lambda(z) = 0$ in the end. We will use this in our inductive proof of the main results in Section \ref{sec:fine-L2}.
\end{remark}

The desired continuity estimate for fine blow-ups now follows from Lemma \ref{lemma:fine-continuity-1}.

\begin{theorem}\label{thm:fine-continuity-2}
	Let $\phi = \{\phi^\alpha\}_{\alpha\in I}$ be a fine blow-up. Then:
	\begin{enumerate}
		\item [(1)] If $\sfrak=n-1$, then $\phi^\alpha\in C^{0,1/4}(\overline{(U^\alpha_0\cup S_0)\cap \mathcal{T}^A_{r_0/4}(x_0,s_0)};P_0^\perp)$ for $\alpha=1,\dotsc,4$, with the estimate
		$$\|\phi^\alpha\|_{C^{0,1/4}((U^\alpha_0\cup S_0)\cap \mathcal{T}^A_{r_0/4}(x_0,s_0))}\leq C\sum_{i=1}^4\|\phi^i\|_{L^2(U^i_0\cap \mathcal{T}^A_{r_0/2}(x_0,s_0))}.$$
		\item [(2)] If $\sfrak<n-1$ (so that $I = \{1,2\}$), then $\phi^\alpha \in C^{0,1/4}(P_0\cap \mathcal{T}^A_{r_0/4}(x_0,s_0);P_0^\perp)$ for $\alpha=1,2$, with the estimate
		$$\|\phi^\alpha\|_{C^{0,1/4}(P_0\cap \mathcal{T}^A_{r_0/4}(x_0,s_0))}\leq C\sum_{i=1,2}\|\phi^i\|_{L^2(P_0\cap \mathcal{T}^A_{r_0/2}(x_0,s_0))}.$$
	\end{enumerate}
	Here, $C = C(n,k,a,M)\in (0,\infty)$.
\end{theorem}

\begin{proof}
	Consider first the case $\sfrak=n-1$. Then for $z\in \mathcal{D}(\phi)\cap \mathcal{T}^A_{7r_0/8}(x_0,s_0)$ and $\alpha=1,\dotsc,4$, set $\phi^\alpha(z):= \lambda(z)^{\perp_{P_0}} - m^\alpha\cdot \lambda(z)^{\top_{P_0}}$, using the notation from Lemma \ref{lemma:fine-continuity-1}. Then, for any $0<\rho_1<\rho_2/2<1/16$, Lemma \ref{lemma:fine-continuity-1} implies that
	$$\sum^4_{\alpha=1}\rho_1^{-n}\int_{B^n_{\rho_1 r_0}(z)\cap U^\alpha_0}|\phi^\alpha(x)-\phi^\alpha(z)|^2\, \ext x \leq C\left(\frac{\rho_1}{\rho_2}\right)^{7/4}\rho_2^{-n}\sum_{\alpha=1}^4\int_{B_{\rho_2 r_0}^n(z)\cap U^\alpha_0}|\phi^\alpha(x)-\phi^\alpha(z)|^2\, \ext x.$$
	Moreover, for balls centred at points $y\in (S\setminus\mathcal{D})\cap \mathcal{T}^{A}_{7r_0/8}(x_0,s_0)$, by the construction of fine blow-ups we know that $\phi^\alpha$ are smooth harmonic functions in $B_{d(y)}(y)\cap U^\alpha_0$, where $d(y):= \dist(y,\mathcal{D}(\phi))$. Thus, we have standard elliptic estimates for harmonic functions in such balls. From here, we need to apply Campanato--Morrey-type arguments to deduce the conclusions of the lemma in this case, however this is sufficiently standard and so we omit the details (cf.~similar arguments in \cite[Lemma 4.3]{Wic14}, \cite[Lemma 5.3]{BK17}, and \cite{Min21b}). This completes the proof in this case.
	
	In the case $\sfrak <n-1$, again set $\phi^\alpha(z):= \lambda(z)^{\perp_{P_0}} - m^\alpha\cdot \lambda(z)^{\top_{P_0}}$ for $z\in \mathcal{D}(\phi)\cap\mathcal{T}^A_{7r_0/8}(x_0,s_0)$. Then Lemma \ref{lemma:fine-continuity-1} gives (as above, as now we can work on full planes rather than half-planes)
	$$\sum_{\alpha=1,2}\rho_1^{-n}\int_{B^n_{\rho_1 r_0}(z)}|\phi^\alpha(x)-\phi^\alpha(z)|^2\, \ext x \leq C\left(\frac{\rho_1}{\rho_2}\right)^{7/4}\rho_2^{-n}\sum_{\alpha=1,2}\int_{B_{\rho_2 r_0}^n(z)}|\phi^\alpha(x)-\phi^\alpha(z)|^2\, \ext x.$$
	Again, for balls centred at points $y\in (S\setminus\mathcal{D})\cap \mathcal{T}^{A}_{7r_0/8}(x_0,s_0)$, by the construction of fine blow-ups we know that $\phi^1,\phi^2$ are smooth harmonic functions in $B_{d(y)}(y)$, where $d(y):= \dist(y,\mathcal{D}(\phi))$. Thus, we may again use Campanato--Morrey-type arguments to conclude.
\end{proof}

\textbf{Remark:} When $\sfrak(\phi) = n-1$, the results of Theorem \ref{thm:fine-continuity-2} imply that $\phi^\alpha$ extends Hölder-continuously to the boundary portion of $S_0$ in $\mathcal{T}^A_{r_0/4}(x_0,s_0)$. In light of standard elliptic (interior) estimates for harmonic functions and the fact that $\left.\phi^\alpha\right|_{S_0}$ is bounded by a constant depending only on $n,k,a,M$ (which comes from the uniform bounds on $\lambda(z)$ provided by Lemma \ref{lemma:fine-continuity-1}) we in fact have
$$|D\phi^\alpha(y)| \leq C\dist(y,S_0)^{-7/8}\sum^4_{i=1}\|\phi^i\|_{L^2(U^i_0\cap \mathcal{T}^A_{r_0/4}(x_0,s_0))}$$
for any $y\in U^\alpha_0\cap \mathcal{T}^A_{3r_0/16}(x_0,s_0)$. This in turn implies that
\begin{equation}\label{E:fine-blow-up-L1}
\|D\phi^\alpha\|_{L^1(U^\alpha_0\cap\mathcal{T}^A_{r_0/8}(x_0,s_0))}\leq C.
\end{equation}
Here, $C = C(n,k,a,M)\in (0,\infty)$.

\section{Proof of Theorem \ref{thm:graphical-rep}}\label{sec:graphical-rep}

We now come to the inductive step of the proof of Theorems \ref{thm:graphical-rep}, \ref{thm:L2-estimates}, and \ref{thm:shift-estimate}. In the present section, we will prove Theorem \ref{thm:graphical-rep}.

The main goal in the proof of Theorem \ref{thm:graphical-rep} is to eliminate the possibility that there are points $z\in \mathcal{D}(V)\cap \{r_{\BC}>\tau r\}$. The main idea for doing this is to take a fine blow-up relative to a (sequence of) well-chosen coarser cones $\Dbf\in \mathcal{P}$. Loosely speaking, we will assume, for the sake of contradiction, that there is a point $z\in\mathcal{D}(V)\cap\{r_{\BC}>\tau r\}$ and we will choose $\Dbf$ to be a pair of planes which, up to a constant, best approximates $V$ subject to the condition $\sfrak(\Dbf)>\sfrak(\BC)$ whilst also having $z\in S(\Dbf)$. Then we will prove that the fine blow-up of $V$ relative to $\Dbf$ (which we can produce by our inductive assumptions) is equal to a pair of planes $\Dbf^\prime\neq \Dbf$ with $\sfrak(\Dbf^\prime)\geq \sfrak(\Dbf)$. Since $\Dbf$ was chosen to be the best such coarse approximation to $V$, this will give a contradiction.

\begin{proof}[Proof of Theorem \ref{thm:graphical-rep}]
	Fix $\beta\in (0,\infty)$, $\tau\in (0,1/8)$, $\sigma\in (1/2,1)$, $a\in (0,1/2)$, and $M\in [1,\infty)$. Suppose that we have sequences $\eps_j,\gamma_j\downarrow 0$, $(V_j)_j\subset\mathcal{V}_\beta$, $(\BC_j)_j\subset\mathcal{P}_{\sfrak}$, and regions $\mathcal{T}^{A_j}_{r_j}(x_j,s_j)$ such that for each $j$, we have that $V_j$ and $\BC_j$ obey Hypothesis $G(\eps_j)$ and $H(A_j,\eps_j,\gamma_j,M)$. We also assume that either $s_j = 0$ for all $j$, or $r_j/s_j\in (a,1-a)$.
	
	By making careful rotations as done previously, we may assume that $A_j\equiv A$ and $S(\BC_j)\equiv S$, where $S$ is a fixed $\sfrak$-dimensional subspace of $P_0$ and $A\subset S$. As in the proof of $\sfrak(\BC)=n-1$ case of the present theorem, we may apply a homothety to assume that $\mathcal{T}^A_{r_j}(x_j,s_j) = \mathcal{T}^A_1(0,s_j^\prime)=:\mathcal{T}_j$, where either $s_j^\prime \equiv 0$ for all $j$, or $s_j^\prime\in (1/(1-a),1/a)$ for all $j$ and $s_j^\prime\to s_0\in [1/(1-a),1/a]$. Write $\mathcal{T}_0 := \mathcal{T}^A_1(0,s_0)$. Again, when we apply this homothety, we know that Hypothesis $H(A,\eps_j,\gamma_j,M)$ is preserved in $\mathcal{T}_j$, however we lose the conditions in Hypothesis $G(\eps_j)$. However, we can still ensure that the mass upper bound of $(\w_n2^n)^{-1}\|V_j\|(B_2(0))\leq 5/2$ holds.
	
	By Hypothesis $H(A,\eps_j,\gamma_j,M)$, we know that $V_j\res \mathcal{T}_0\weakly 2|P_0|$ as varifolds in $\mathcal{T}_0$, and indeed we will restrict ourselves to $\mathcal{T}_0$ from now on. To complete the proof of Theorem \ref{thm:graphical-rep}, it suffices to show that Theorem \ref{thm:graphical-rep}(ii) holds along a subsequence.
	
	Suppose first (for the sake of contradiction) that for all $j$ sufficiently large there exists $z_j\in\mathcal{D}(V_j)\cap \{r_{\BC_j}\geq \tau\}\cap \mathcal{T}^A_{(1+\sigma)/2}(0,s_j^\prime)$. Note that as $S(\BC_j)\equiv S$ is independent of $j$, so is $r_{\BC_j}$; we will write $r_{\BC_j}\equiv r_S$ as a shorthand.
	
	\textbf{Step 1:} \emph{Constructing the sequence of coarser cones $(\Dbf_j)_j$.} For each $j$, we will pick a new cone $\Dbf_j$. The procedure for this choice is involved, and so we break it down into separate steps.
	
	We begin with an intermediate step constructing cones $\widetilde{\Dbf}_j$ which are essentially optimal coarser cones to $\BC_j$. This will involve working with an intermediary collection of cones $\{\Dbf_j^{(m)}\}_{m=0}^{\widetilde{m}}$, where $\widetilde{m}\geq 1$. This (finite) collection is produced by an algorithm inductively on $m$. Start by setting $\Dbf_j^{(0)}:= \BC_j$. Then, choose $\Dbf_j^{(1)}\in \mathcal{P}$ a coarser cone to $\Dbf_j^{(0)}$ such that
	$$\hat{F}_{V_j,\Dbf_j^{(1)}}(\mathcal{T}_j) \leq \frac{3}{2}\mathscr{F}^*_{V_j,\BD_j^{(0)}}(\mathcal{T}_j) \equiv \frac{3}{2}\inf\hat{F}_{V_j,\Dbf_j^{(0)}}(\mathcal{T}_j)$$
	where we recall that the infimum is taken over all coarser cones $\Dbf\in \mathcal{P}$ to $\Dbf_j^{(0)}$ (recall $\mathcal{P}$ consists only of transversely intersecting planes). We may pass to a subsequence to assume that $\sfrak(\Dbf_j^{(1)})\equiv \sfrak_1$ is constant. This is the $m=1$ case of our collection (which always happens). Now, inductively for each $m$ in turn, we consider the following trichotomy:
	\begin{enumerate}
		\item [(i)] There is a subsequence along which
		\begin{equation}\label{E:tri-1}
			\mathscr{F}^*_{V_j,\Dbf_j^{(m)}}(\mathcal{T}_j)^{-1}\hat{F}_{V_j,\Dbf_j^{(m)}}(\mathcal{T}_j)\to 0.
		\end{equation}
		\item [(ii)] (i) fails and $\sfrak_m=n-1$.
		\item [(iii)] (i) fails and $\sfrak_m<n-1$.
	\end{enumerate}
	In the case of (i) or (ii), we set $\widetilde{m} := m$ and the process terminates. In the case of (iii), the process continues in the following manner. 
	
	\textbf{Note:} In the case of (iii), by the failure of (i), this tells us that there is a constant $C = C(n,k,\tau,\sigma,a,M,\beta)>0$ such that for all $j$ sufficiently large,
	\begin{equation}\label{E:i-fail}
		\hat{F}_{V_j,\Dbf_j^{(m)}}(\mathcal{T}_j) \geq C\mathscr{F}^*_{V_j,\Dbf_{j}^{(m)}}(\mathcal{T}_j).
	\end{equation}
	Suppose that when $m=1$ we are in case (iii) of the above trichotomy (and so the process does not terminate). We then set $m=2$, and continue inductively in the same manner. More precisely, for each $j$ we choose a $\Dbf_j^{(m)}\in\mathcal{P}$ coarser to $\Dbf_j^{(m-1)}$ satisfying
	\begin{equation}\label{E:graphical-rep-0}
	\hat{F}_{V_j,\Dbf_j^{(m)}}(\mathcal{T}_j)\leq\frac{3}{2}\mathscr{F}^*_{V_j,\Dbf_j^{(m-1)}}(\mathcal{T}_j)
	\end{equation}
	and we pass to a subsequence to ensure that $\sfrak(\Dbf_j^{(m)})\equiv \sfrak_m$ ($>\sfrak_{m-1}$) is constant. Then once again we consider the trichotomy described in (i), (ii), and (iii) above. If either (i) or (ii) holds, the process terminates. In the case that (iii) holds, we increase the value of $m$ by $1$ and continue. This process must terminate in finitely many steps as $\sfrak_m$ is strictly increasing and bounded above by $n-1$. Hence, we must have $\widetilde{m}\leq n$. In particular, the process must terminate on either (i) or (ii).
	
	We then make the following definitions depending on whether the process terminates on (i) or (ii). Firstly, in either case we choose a plane $P_j^{*}$ coarser to $\Dbf_j^{(\widetilde{m})}$ (i.e.~so that $P_j^{*}\supsetneq S(\Dbf_j^{(\widetilde{m})})$) and such that
	$$\hat{E}_{V_j,P_j^{*}}(\mathcal{T}_j) = \mathscr{E}_{V_j,\Dbf_j^{(\widetilde{m})}}(\mathcal{T}_j).$$
	Then, we define cones $\widetilde{\Dbf}_j$ as follows:
	\begin{itemize}
		\item If the process terminates on (i), we set $\widetilde{\Dbf}_j:= \Dbf_j^{(\widetilde{m})}$ and pass to a subsequence along which \eqref{E:tri-1} holds with $m=\widetilde{m}$.
		\item If the process terminates on (ii), then we set $\widetilde{\Dbf}_j:= |P_j^*|$. Notice also that by the failure of (i), we have \eqref{E:i-fail} holds with $m=\widetilde{m}$, and so since $\sfrak_{\widetilde{m}}=n-1$ we get,
			\begin{equation}\label{E:ii-holds}
			\hat{F}_{V_j,\Dbf_j^{(\widetilde{m})}}(\mathcal{T}_j)\geq C\mathscr{F}^*_{V_j,\Dbf_j^{(\widetilde{m})}}(\mathcal{T}_j) \equiv C\mathscr{E}_{V_j,\Dbf_j^{(\widetilde{m})}}(\mathcal{T}_j) \equiv C\hat{E}_{V_j,P_j^*}(\mathcal{T}_j).
			\end{equation}
	\end{itemize}
	This completes the construction of $\widetilde{\Dbf}_j$ as well as of the planes $P_j^*$. Notice that $A\subset S \equiv S(\BC_j)\subsetneq S(\Dbf_j^{(\widetilde{m})})\subset S(\widetilde{\Dbf}_j)\subset P^*_j$.
		
	Now for each $j$, pick a rotation $R_j$ of $\R^{n+k}$ attaining the infimum
	$$|R_j-\id_{\R^{n+k}}| = \inf_{R}|R-\id_{\R^{n+k}}|$$
	where this infimum is taken over all rotations of $\R^{n+k}$ such that $R|_S = \id_S$, $R(S(\Dbf^{(\widetilde{m})}_j)) \subset S_0$, and $R(P_j^*) = P_0$. This infimum is always attained in this set, and note that $R_j(\mathcal{T}_j) = \mathcal{T}_j$ since $\left.R_j\right|_S = \id_S$ and $A\subset S$. Note also that $R_j\to \id_{\R^{n+k}}$. Now we rotate the system by $R_j$, and set
	$$\Dbf_j:= (R_j)_\#\widetilde{\Dbf}_j \qquad \text{and} \qquad W_j:= (R_j)_\#V_j.$$
	By definition of $\widetilde{\Dbf}_j$ based on the (i) and (ii) of the previous trichotomy, we either have $\Dbf_j = |P_0|$ for every $j$ (when the previous process terminated on (ii)), or there is $\sfrak^\prime\in \{n-1,\dotsc,\sfrak+1\}$ for which $\Dbf_j\in \mathcal{P}_{\sfrak^\prime}$. Also, based on the chosen rotation $R_j$ and the choice of $P_j^*$, in the latter case (i.e.~when the previous process terminated on (i)) we have
	\begin{equation}\label{E:graphical-rep-1}
	\hat{E}_{W_j}(\mathcal{T}_j) = \mathscr{E}_{W_j,\Dbf_j}(\mathcal{T}_j).
	\end{equation}
	Moreover, by our contradiction assumption we know that
	$$R_j(z_j)\in \mathcal{D}(W_j)\cap \{r_{(R_j)_\#\BC_j}\geq\tau\}\cap \mathcal{T}^A_{(1+\sigma)/2}(0,s_j^\prime).$$
	\textbf{Note:} Since $R_j$ is a rotation which might not fix $P_0$, we no longer know whether $W_j$ is in $\mathcal{V}_{\beta}$. However, since $R_j\to \id_{\R^{n+k}}$ we do know by Proposition \ref{prop:rotation-class} that $(\eta_{0,1-\eta_j})_\#W_j\in \mathcal{V}_{\beta_j}$ for some $\beta_j\uparrow \beta$ and $\eta_j\downarrow 0$. In particular, since for all $j$ sufficiently large $\mathcal{V}_{\beta_j}\subset \mathcal{V}_{\beta/2}$, these varifolds lie in the same class for some fixed parameter. Thus, up to replacing $\beta$ by $\beta/2$ and scaling by an amount arbitrarily close to $1$, we can without loss of generality assume $W_j\in \mathcal{V}_{\beta}$; for notational simplicity we will assume that this is done and the notation is understood. Similarly, any time throughout this proof where arbitrarily small rotations are made, we will abuse notation and write $\mathcal{V}_\beta$ for the class the rotated varifolds lie in.

\textbf{Step 2:} \emph{Optimality of the height excess relative to $\Dbf_j$.} We have two cases to consider, depending on whether the inductive construction in Step 1 terminated at (i) or (ii) of the trichotomy therein.

Suppose to begin with it terminated on (ii), so that $\Dbf_j = |P_0|$ for all $j$. This means that \eqref{E:ii-holds} gives
$$\hat{E}_{W_j}(\mathcal{T}_j) \leq C^{-1}\hat{F}_{V_j,\Dbf_{j}^{(\widetilde{m})}}(\mathcal{T}_j)$$
where $C = C(n,k,\sigma,a,M,\beta)$.
By the definition of $\Dbf_j^{(\widetilde{m})}$, namely the relation \eqref{E:graphical-rep-0} with $m=\widetilde{m}$ therein, combined with the failure of (i) in the trichotomy at $m=\widetilde{m}-1$, i.e.~\eqref{E:i-fail} with $m = \widetilde{m}-1$ therein, we get
$$\hat{F}_{V_j,\Dbf_j^{(\widetilde{m})}}(\mathcal{T}_j) \leq \frac{3}{2}\mathscr{F}^*_{V_j,\Dbf_j^{(\widetilde{m}-1)}}(\mathcal{T}_j) \leq \frac{3}{2}C^{-1}\hat{F}_{V_j,\Dbf_j^{(\widetilde{m}-1)}}(\mathcal{T}_j)$$
so that
$$\hat{E}_{W_j}(\mathcal{T}_j)\leq \frac{3}{2}C^{-2}\hat{F}_{V_j,\Dbf_j^{(m-1)}}(\mathcal{T}_j).$$
Repeating this another $\widetilde{m}-2$ times, we therefore get
$$\hat{E}_{W_j}(\mathcal{T}_j) \leq C\mathscr{F}^*_{V_j,\BC_j}(\mathcal{T}_j) \equiv C\mathscr{F}^*_{W_j,(R_j)_\#\BC_j}(\mathcal{T}_j)$$
for suitable $C = C(n,k,\sigma,a,M,\beta)$. Noting that the term on the right-hand side is an infimum over all coarser cones to $(R_j)_\#\BC_j$, but $R_j(S(\BC_j)) = S(\BC_j)$, this infimum is the same with $\BC_j$ in place of $(R_j)_\#\BC_j$, i.e.~we have
\begin{equation}\label{E:graphical-rep-2}
	\hat{E}_{W_j}(\mathcal{T}_j) \leq C\mathscr{F}^*_{W_j,\BC_j}(\mathcal{T}_j).
\end{equation}
This is the desired optimality of height excess in this case.

Now consider the second case, where the previous trichotomy terminated on (i). An identical argument to the case above based on the failure of (i) and the definition of $\Dbf_j^{(m)}$ gives
\begin{equation}\label{E:graphical-rep-3}
	\hat{F}_{W_j,\Dbf_j}(\mathcal{T}_j) \leq C\mathscr{F}^*_{W_j,\BC_j}(\mathcal{T}_j)
\end{equation}
where $C = C(n,k,\sigma,a,M,\beta)$. This is the desired inequality here also, giving that $\Dbf_j$ is optimal (up to a constant) over all coarser cones to $\BC_j$.

Before continuing to the next step, first notice that since $V_j\res \mathcal{T}_0\weakly 2|P_0\cap \mathcal{T}_0|$, we certainly have that $\Dbf_j$ converges in Hausdorff distance to $P_0$. In particular, $P_j^*$ converges to $P_0$, and thus $R_j\to \id_{\R^{n+k}}$ (as mentioned before), which in turn guarantees that $W_j\res \mathcal{T}_0\weakly 2|P_0\cap\mathcal{T}_0|$ and that $(R_j)_\#\BC_j$ converges in Hausdorff distance to $P_0$. Since in the second case above we have already arranged that $S\subset S(\Dbf_j)\subset S_0$ for every $j$, we can further perform a rotation (which fixes $P_0$, $S_0$, and $S$) to assume without loss of generality that $S(\Dbf_j)\equiv S^\prime$ is independent of $j$, with $S\subset S^\prime\subset S_0$.

There will be two cases to consider again, depending on whether $\Dbf_j$ is the single plane $P_0$ or is in $\mathcal{P}$. The latter case will require significant more work and so we start with it. The case of a single plane will be dealt with in Step 6.

\textbf{Step 3:} \emph{Moving the good density point to the spine.} As previously noted, we are currently focusing on the case where $\Dbf_j\in \mathcal{P}$, i.e.~the case where (i) of the trichotomy from Step 1 occurs. We first claim that $W_j$ and $\Dbf_j$ satisfy Hypothesis $H(A,\tilde{\eps}_j,\tilde{\gamma}_j,1)$ in $\mathcal{T}_j$ for appropriate choices of $\tilde{\eps}_j,\tilde{\gamma}_j\downarrow 0$. Indeed, conditions (A1) and (A2) are clear from the discussion at the end of Step 2. Furthermore, condition (A3) with $M=1$ is guaranteed by \eqref{E:graphical-rep-1}, and (A4) with a suitable sequence of $\tilde{\gamma}_j\downarrow 0$ is guaranteed by \eqref{E:tri-1}. 

In light of this and the induction hypothesis, we are now in a position to apply the main results of Section \ref{sec:fine-L2} to $W_j$ and $\Dbf_j$ (it should be noted that, if we undo the homothety done at the start of the proof, $W_j$ obeys Hypothesis $G(\tilde{\eps}_j)$ for suitable $\tilde{\eps}_j\downarrow 0$ since it is just a rotation of $V_j$ which obeyed this hypothesis $G(\eps_j)$; as such we can apply inductively the results from Section \ref{sec:fine-L2}).

In particular, using \eqref{E:cor-1}, we can assume that there is a relatively closed subset $\mathcal{D}\subset S^\prime\cap \mathcal{T}_0$ for which
\begin{equation}\label{E:graphical-rep-4}
	\mathcal{D}(W_j)\cap\mathcal{T}_0\to \mathcal{D}\qquad \text{locally in Hausdorff distance in $\mathcal{T}_0$.}
\end{equation}
We now wish to rotate $W_j$ slightly to place the point $R_j(z_j)\in \mathcal{D}(W_j)$ onto $S^\prime\equiv S(\Dbf_j)$. Indeed, for each $j$ pick a rotation $\phi_j$ of $\R^{n+k}$ attaining the infimum
$$|\phi_j - \id_{\R^{n+k}}| = \inf_\phi |\phi-\id_{\R^{n+k}}|$$
where the infimum is taken over all rotations $\phi$ of $\R^{n+k}$ which obey $\phi|_S = \id_S$ and $\phi(R_j(z_j))\in S^\prime$. Then, consider $\widetilde{W}_j:= (\phi_j)_\#W_j$. Note that if $y_j:= \phi_j(R_j(z_j))$, we now have
\begin{equation}\label{E:graphical-rep-5}
	y_j\in\mathcal{D}(\widetilde{W}_j)\cap S^\prime\cap \{r_{S}\geq \tau\}\cap \mathcal{T}^A_{(1+\sigma)/2}(0,s_j^\prime).
\end{equation}
We now claim that $\widetilde{W}_j$ and $\Dbf_j$ satisfy Hypothesis $H(A,\eps_j^\prime,\gamma_j^\prime,M^\prime)$ in $\mathcal{T}_j$ for some $\eps_j^\prime,\gamma_j^\prime\downarrow 0$, and some $M^\prime = M^\prime(n,k,\sigma,a)\in (1,\infty)$. Of course, (A1) and (A2) are immediate from the corresponding conclusions for $W_j$ and $\Dbf_j$, as $\phi_j$ fixes $A$ and so leaves $\mathcal{T}_0$ invariant. Hypothesis (A3) and (A4) will require more work to check. (We also note again that, because of \eqref{E:graphical-rep-4}, we know that $\phi_j\to \id_{\R^{n+k}}$, and so $\widetilde{W}_j$ (under an appropriate homothety) satisfies Hypothesis $G(\eps_j^\prime)$ for suitable $\eps_j^\prime$, since $W_j$ did. Again, this rotation will slightly change the value of $\beta$, but this can be controlled, as described before, since $\phi_j\to \id_{\R^{n+k}}$.)

Let us first verify (A3). For this, it will be crucial to use that the points $R_j(z_j)$ are uniformly bounded away from $S$ by $\tau>0$. Indeed, because $r_{S}(R_j(z_j))>\tau$, an elementary argument based on Euclidean geometry (essentially using that $|\sin(x)/x|$ is close to $1$ for $x$ close to $0$) gives that 
$$|\phi_j-\id_{\R^{n+k}}| \leq C\left|\frac{R_j(z_j)^{\perp_{S^\prime}}}{r_S(R_j(z_j))}\right| \leq C\tau^{-1}|R_j(z_j)^{\perp_{S^\prime}}|$$
for some $C = C(n,k)$ (notice that the distance to $S$ of $R_j(z_j)/r_S(R_j(z_j))$ is $1$). In turn, another simple argument of a similar nature gives
$$\nu(P_0,\phi_j(P_0)) \leq Cr_S(R_j(z_j))^{-1}\cdot\nu(P_0,(\tau_{R_j(z_j)^{\perp_{S^\prime}}})_\#P_0)$$
which therefore gives
\begin{equation}\label{E:graphical-rep-7}
	\nu(P_0,\phi_j(P_0)) \leq C\tau^{-1}|R_j(z_j)^{\perp_{P_0}}|.
\end{equation}
(note also that as $S^\prime \subset P_0$, $x^{\perp_{S^\prime}\perp_{P_0}}\equiv x^{\perp_{P_0}}$).
Hence, by Theorem \ref{thm:shift-estimate} applied to $W_j$ and $\Dbf_j$ we have
$$\nu(P_0,\phi_j(P_0)) \leq C_*\hat{E}_{W_j,\Dbf_j}(\mathcal{T}_j)$$
where $C_* = C_*(n,k,\tau,\sigma,a,M)$. Since $\hat{E}_{W_j,\Dbf_j}(\mathcal{T}_j)\leq C\tilde{\gamma}_j\hat{E}_{W_j}(\mathcal{T}_j)$ by (A4) for $W_j$ and $\Dbf_j$ in $\mathcal{T}_j$ (again, this uses Theorem \ref{thm:allard-sup-estimate} to control the two-sided excess relative to a plane by the one-sided excess), combining these we get from the triangle inequality,
$$\hat{E}_{\widetilde{W}_j}(\mathcal{T}_j) \leq C\hat{E}_{W_j}(\mathcal{T}_j) + C\nu(P_0,\phi_j(P_0)) \leq C\hat{E}_{W_j}(\mathcal{T}_j) + C_*\tilde{\gamma}_j\hat{E}_{W_j}(\mathcal{T}_j)$$
and thus after rearranging gives for all $j$ sufficiently large (also recalling \eqref{E:graphical-rep-1})
\begin{equation}\label{E:graphical-rep-8}
	\hat{E}_{\widetilde{W}_j}(\mathcal{T}_j) \leq C\hat{E}_{W_j}(\mathcal{T}_j) \equiv C\mathscr{E}_{W_j,\Dbf_j}(\mathcal{T}_j).
\end{equation}
where by choosing $j$ sufficiently large depending on $\tau$ so that $C_*\tilde{\gamma}_j<1/2$, we can ensure the constant $C$ above is independent of $\tau$.

Now pick any plane $P\supset S^\prime$ with $\hat{E}_{\widetilde{W}_j,P}(\mathcal{T}_j) = \mathscr{E}_{\widetilde{W}_j,\Dbf_j}(\mathcal{T}_j)$. For such a plane $P$ we must have $\nu(P)\leq C\hat{E}_{\widetilde{W}_j}(\mathcal{T}_j)$ where $C = C(n,k)$; this can be readily seen by a direct computation or a simple contradiction argument. So, \eqref{E:graphical-rep-8} gives
\begin{equation}\label{E:graphical-rep-9}
	\nu(P) \leq C\hat{E}_{W_j}(\mathcal{T}_j).
\end{equation}
Since $P\supset S^\prime$, by the same reasoning which lead us to \eqref{E:graphical-rep-7}, using also \eqref{E:translates-2}, we have
\begin{align}
	\nonumber\nu(P,\phi_j(P)) & \leq C\tau^{-1}\nu(P,(\tau_{R_j(z_j)^{\perp_{S^\prime}}})_\#P)\\
	& \leq C_*\left(|R_j(z_j)^{\perp_{S^\prime}\perp_{P_0}}| + \nu(P)|R_j(z_j)^{\perp_{S^\prime}\top_{P_0}}|\right).\label{E:graphical-rep-10}
\end{align}
where again $C_*$ also depends on $\tau$. Combining this with \eqref{E:graphical-rep-9}, Lemma \ref{lemma:nu-a}(i) (in order to replace $\nu(P)$ in the above with $\nu(\Dbf_j)$), and Theorem \ref{thm:shift-estimate} (applied to $W_j$ and $\Dbf_j$) we get
\begin{equation}\label{E:graphical-rep-11}
	\nu(P,\phi_j(P)) \leq C_*\hat{E}_{W_j,\Dbf_j}(\mathcal{T}_j).
\end{equation}
Now using the triangle inequality, $\phi_j(\mathcal{T}_j) = \mathcal{T}_j$, and the defining property of $P$, we have
\begin{align*}
	\mathscr{E}_{W_j,\Dbf_j}(\mathcal{T}_j) \leq \hat{E}_{W_j,P}(\mathcal{T}_j) & = \hat{E}_{\widetilde{W}_j,\phi_j(P)}(\mathcal{T}_j)^2\\
	& \leq 2\hat{E}_{\widetilde{W}_j,P}(\mathcal{T}_j) + C\nu(P,\phi_j(P))\\
	& \leq 2\mathscr{E}_{\widetilde{W}_j,\Dbf_j}(\mathcal{T}_j) + C_*\hat{E}_{W_j,\Dbf_j}(\mathcal{T}_j).
\end{align*}
Property (A4) for $W_j$ and $\Dbf_j$ allows us to absorb the second term on the last line here onto the left-hand side, provided $j$ is sufficiently large (again, using Theorem \ref{thm:allard-sup-estimate} to control the two-sided excess relative to a plane by the one-sided excess). Thus we get for all $j$ large enough,
$$\mathscr{E}_{W_j,\Dbf_j}(\mathcal{T}_j) \leq 4\mathscr{E}_{\widetilde{W},\Dbf_j}(\mathcal{T}_j).$$
Combining this with \eqref{E:graphical-rep-8} gives
$$\hat{E}_{\widetilde{W}_j}(\mathcal{T}_j) \leq M^\prime\mathscr{E}_{\widetilde{W},\Dbf_j}(\mathcal{T}_j)$$
for some $M^\prime = M^\prime(n,k,\sigma,a)\in (1,\infty)$. This establishes (A3) for $\widetilde{W}_j$ and $\Dbf_j$ for this choice of $M^\prime$.

Next we check (A4) for $\widetilde{W}_j$ and $\Dbf_j$. Once again, the reasoning which leads to \eqref{E:graphical-rep-7} gives
$$\nu(\Dbf_j,\phi_j(\Dbf_j)) \leq C_*\nu(\Dbf_j,(\tau_{R_j(z_j)})_\#\Dbf_j)$$
for some $C_* = C_*(n,k,\tau)$, and so combining this with \eqref{E:translates-2} and Theorem \ref{thm:shift-estimate} for $W_j$ and $\Dbf_j$, as well as (A4) for $W_j$ and $\Dbf_j$, we have
\begin{equation}\label{E:graphical-rep-12}
	\nu(\Dbf_j,\phi_j(\Dbf_j)) \leq C_*\hat{E}_{W_j,\Dbf_j}(\mathcal{T}_j) \leq C_*\tilde{\gamma}_j\mathscr{F}^*_{W_j,\Dbf_j}(\mathcal{T}_j).
\end{equation}
Since by the triangle inequality $\hat{F}_{\widetilde{W}_j,\Dbf_j}(\mathcal{T}_j) \leq 2\hat{F}_{W_j,\Dbf_j}(\mathcal{T}_j) + C\nu(\Dbf_j,\phi_j(\Dbf_j))$, this means that (again using (A4) for $W_j$ and $\Dbf_j$)
\begin{equation}\label{E:graphical-rep-13}
	\hat{F}_{\widetilde{W}_j,\Dbf_j}(\mathcal{T}_j) \leq C_*\tilde{\gamma}_j\mathscr{F}_{W_j,\Dbf_j}^*(\mathcal{T}_j).
\end{equation}
Now to get the appropriate upper bound on the right-hand side of \eqref{E:graphical-rep-13}, pick a cone $\Dbf\in\mathcal{P}$ with $S(\Dbf)\supsetneq S(\Dbf_j)$ and $\hat{F}_{\widetilde{W}_j,\Dbf}(\mathcal{T}_j) \leq \frac{3}{2}\mathscr{F}^*_{\widetilde{W}_j,\Dbf_j}(\mathcal{T}_j)$. Such a cone must satisfy an estimate of the form $\nu(\Dbf)\leq C\hat{E}_{\widetilde{W}_j}(\mathcal{T}_j)$ with $C = C(n,k)$ (again, one can verify this directly or with a simple contradiction argument). In light of \eqref{E:graphical-rep-8}, this means
\begin{equation}\label{E:graphical-rep-14}
	\nu(\Dbf) \leq C\hat{E}_{W_j}(\mathcal{T}_j).
\end{equation}
Next, using an estimate analogous to \eqref{E:graphical-rep-10}, followed by \eqref{E:graphical-rep-14}, Lemma \ref{lemma:nu-a}(i), Theorem \ref{thm:shift-estimate} (applied to $W_j$ and $\Dbf_j$), and then (A4) for $W_j$ and $\Dbf_j$, we establish that
\begin{equation}\label{E:graphical-rep-15}
	\nu(\Dbf,(\phi_j)_\#\Dbf) \leq C_*\tilde{\gamma}_j\mathscr{F}^*_{W_j,\Dbf_j}(\mathcal{T}_j)
\end{equation}
where $C_*$ again also depends on $\tau$.
Now since the triangle inequality gives
\begin{align*}
	\mathscr{F}^*_{W_j,\Dbf_j}(\mathcal{T}_j)\leq \hat{F}_{W_j,\Dbf}(\mathcal{T}_j) = \hat{F}_{\widetilde{W}_j,(\phi_j)_\#\Dbf}(\mathcal{T}_j) \leq 2\hat{F}_{\widetilde{W}_j,\Dbf}(\mathcal{T}_j) + C\nu(\Dbf,(\phi_j)_\#\Dbf),
\end{align*}
we can combine this with \eqref{E:graphical-rep-15} and then with \eqref{E:graphical-rep-13} and the defining property of $\Dbf$ to get
$$\hat{F}_{\widetilde{W}_j,\Dbf_j}(\mathcal{T}_j) \leq C\tilde{\gamma}_j\mathscr{F}^*_{\widetilde{W}_j,\Dbf_j}(\mathcal{T}_j)$$
which provides (A4) for $\widetilde{W}_j$ and $\Dbf_j$ for suitable $\gamma_j^\prime$. Thus, we have now shown that $\widetilde{W}_j$ and $\Dbf_j$ satsify Hypothesis $H(A,\eps_j^\prime,\gamma_j^\prime,M^\prime)$ in $\mathcal{T}_j$ for suitable $\eps_j^\prime,\gamma^\prime_j\downarrow 0$ and $M^\prime = M^\prime(n,k,\sigma,a,M)$. In particular, we can apply our inductive hypothesis to them.

\textbf{Step 4:} \emph{Blowing-up when $\Dbf_j\in \mathcal{P}$.} Continuing on from Step 3, it is now clear that $\widetilde{W}_j$ and $\Dbf_j$ satisfy Hypothesis $\dagger$ from Section \ref{sec:fine-blow-ups} in $\mathcal{T}_j$ where, crucially, $\sfrak(\Dbf_j) \equiv \sfrak^\prime>\sfrak$. Thus, we may take a fine blow-up of $\widetilde{W}_j$ off $\Dbf_j$ in $\mathcal{T}_j$: let $\psi = \{\psi^\alpha\}_{\alpha=1,2}$ denote this fine blow-up (where the indexing is viewing the fine blow-up as $\psi^1,\psi^2$, defined on the full plane, which we may do as $\Dbf_j\in\mathcal{P}$ for all $j$).

We claim that we can also take a (fine) blow-up of $\widetilde{\BC}_j:= (\phi_j\circ R_j)_\#\BC_j$ relative to $\Dbf_j$ in $\mathcal{T}_j$, but using the excess $\hat{E}_{\widetilde{W}_j,\Dbf_j}(\mathcal{T}_j)$ as the scaling factor (which is the same used when blowing up $\widetilde{W}_j$ off $\Dbf_j$). As these are unions of planes this is much easier to do, and indeed we just need to show
\begin{equation}\label{E:graphical-rep-16}
	\nu(\widetilde{\BC}_j,\Dbf_j)\leq C\hat{E}_{\widetilde{W}_j,\Dbf_j}(\mathcal{T}_j).
\end{equation}
To see this, it is not difficult to use Theorem \ref{thm:graphical-rep} for $\widetilde{W}_j$ and $\Dbf_j$ (which is possible by Step 3 and our inductive assumption) and argue similarly to the proof of the second inequality in Lemma \ref{lemma:nu-a}(i) to get
$$\nu(\widetilde{\BC}_j,\Dbf_j) \leq C\hat{F}_{\widetilde{W}_j,\widetilde{\BC}_j}(\mathcal{T}_j) + C\hat{F}_{\widetilde{W}_j,\Dbf_j}(\mathcal{T}_j)$$
Indeed, this is using that in a region a fixed distance away from the spine, the closest plane in $\Dbf_j$ or $\BC_j$ from a point in $V_j$ is locally constant (which can be seen from performing a coarse blow-up off $P_0$, as the coarse blow-up of both $\BC_j$ and $\Dbf_j$ must agree with the coarse blow-up of $V_j$, which furthermore is a union of two distinct planes here using Lemma \ref{lemma:nu-a}(ii)), and so we can make an analogous argument to that in the proof of the second inequality of  Lemma \ref{lemma:nu-a}(i), working over the (fixed) individual planes in $\BC_j$ and $\Dbf_j$ for which $V_j$ is close to in a given region. Now, using Hypothesis (A4) for $V_j$ and $\BC_j$ in the first term of the above (as it equals exactly $\hat{F}_{V_j,\BC_j}(\mathcal{T}_j)$; note this is the first time we are using this hypothesis on $V_j$ and $\BC_j$), and then using Remark \ref{remark:after-graphical-rep} following the statement of Theorem \ref{thm:graphical-rep} (applied to $\widetilde{W}_j$ and $\Dbf_j$) for the second (as then the two-sided and one-sided excesses of $\widetilde{W}_j$ to $\Dbf_j$ are comparable, cf.~\eqref{E:graphical-rep-19} below also), we get exactly \eqref{E:graphical-rep-16}.

Thus, this allows us to take the blow-up of $(\widetilde{\BC}_j)_j$ off $(\Dbf_j)_j$ with this scaling factor; denote it by $\{\chi^\alpha\}_{\alpha=1,2}$. As $\widetilde{\BC}_j$ and $\Dbf_j$ are unions of planes we know that $\chi^1$ and $\chi^2$ are linear. Notice that from (A4) for $V_j$ and $\BC_j$, it follows that (up to reordering)
\begin{equation}\label{E:graphical-rep-17}
	\psi^\alpha = \chi^\alpha \qquad \text{for }\alpha=1,2,
\end{equation}
and thus $\graph(\psi^\alpha)$ is a plane for $\alpha=1,2$. Since $S\subsetneq S^\prime\subsetneq P_0$, we have that $\chi^\alpha|_S\equiv 0$, and hence by \eqref{E:graphical-rep-17} that $\psi^\alpha|_S\equiv 0$ for $\alpha=1,2$. Hence, $\graph(\psi^1)\cap\graph(\psi^2)\supset S$. Moreover, in light of \eqref{E:graphical-rep-5} and Remark \ref{remark:fine-blow-up-is-0} following the proof of Lemma \ref{lemma:fine-continuity-1} (since $y_j\in S^\prime\in \mathcal{D}(\widetilde{W}_j)$), we know that, if we pass to a subsequence for which $y_j\to y_0$, then
$$y_0\in (S^\prime\setminus S)\cap \overline{T^A_{(1+\sigma)/2}(0,s_0)}$$
is a point such that $\psi^\alpha(y_0) = 0$ for $\alpha=1,2$. Hence if $S_*$ is the subspace spanned by $S$ and $y_0$, then we have $S_*\supsetneq S$ and
\begin{equation}\label{E:graphical-rep-18}
	\graph(\psi^1)\cap\graph(\psi^2)\supset S_*.
\end{equation}
To get useful information out of the fine blow-up $\psi$ we want to know that it is not trivial, i.e.~$\psi\not\equiv 0$. By \eqref{E:graphical-rep-17}, to show this is suffices to show that $\chi\not\equiv 0$. To see this, note from the triangle inequality we have
$$\hat{E}_{\widetilde{W}_j,\Dbf_j}(\mathcal{T}_j)\leq 2\hat{E}_{\widetilde{W}_j,\widetilde{\BC}_j}(\mathcal{T}_j) + C\nu(\widetilde{\BC}_j,\Dbf_j)$$
for some $C = C(n,k)$. The first term on the right-hand side of this inequality can be absorbed into the right-hand side using the fact that it equals $\hat{E}_{V_j,\BC_j}(\mathcal{T}_j)$, and so we can use Hypothesis (A4) for $V_j$ and $\BC_j$ and the fact that $S(\phi_j^{-1}(\Dbf_j))\supsetneq S(\BC_j)$ to get
\begin{equation}\label{E:graphical-rep-19}
	\hat{E}_{\widetilde{W}_j,\widetilde{\BC}_j}(\mathcal{T}_j) < \gamma_j\hat{F}_{W_j,\phi_j^{-1}(\Dbf_j)}(\mathcal{T}_j) = \gamma_j\hat{F}_{\widetilde{W}_j,\Dbf_j}(\mathcal{T}_j) \leq C\gamma_j\hat{E}_{\widetilde{W}_j,\Dbf_j}(\mathcal{T}_j)
\end{equation}
where in the last inequality we have used Remark \ref{remark:after-graphical-rep} following the statement of Theorem \ref{thm:graphical-rep} applied to $\tilde{W}_j$ and $\Dbf_j$. Thus, for all $j$ sufficiently large we get
$$\hat{E}_{\widetilde{W}_j,\Dbf_j}(\mathcal{T}_j) \leq C\nu(\widetilde{\BC}_j,\Dbf_j)$$
from which it follows that $\chi\not\equiv 0$.

\textbf{Step 5:} \emph{$L^2$ convergence to the fine blow-up and the contradiction when $\Dbf_j\in\mathcal{P}$.} Fix now $\tilde{\sigma}\in (0,1)$. Then, using the triangle inequality, \eqref{E:graphical-rep-19}, and \eqref{E:graphical-rep-16}, we have
\begin{align*}
	\int_{\mathcal{T}_j\setminus\mathcal{T}^A_{\tilde{\sigma}}(0,s_j^\prime)}\dist^2(x,\Dbf_j)\, \ext\|\widetilde{W}_j\|(x) & \leq 2\int_{\mathcal{T}_j\setminus\mathcal{T}^A_{\tilde{\sigma}}(0,s_j^\prime)}\dist^2(x,\widetilde{\BC}_j)\, \ext\|\widetilde{W}_j\|(x) + C(1-\tilde{\sigma})\nu(\widetilde{\BC}_j,\Dbf_j)^2\\
	& \leq C[\gamma^2_j+(1-\tilde{\sigma})]\int_{\mathcal{T}_j}\dist^2(x,\Dbf_j)\, \ext\|\widetilde{W}_j\|(x).
\end{align*}
Hence, we can therefore choose $\tilde{\sigma}$ close to $1$ and $j$ sufficiently large to ensure that the quantity
$$\hat{E}_{\widetilde{W}_j,\Dbf_j}(\mathcal{T}_j)^{-1}\int_{\mathcal{T}_j\setminus\mathcal{T}_{\tilde{\sigma}}^A(0,s_j^\prime)}\dist^2(x,\Dbf_j)\, \ext\|\widetilde{W}_j\|(x)$$
is as small as we wish; in turn, analogously to Proposition \ref{prop:full-L2}, this implies that we have global $L^2$ convergence to the blow-up on $\mathcal{T}_0$. This means that, if we write $\Dbf_j = \sum_{\alpha=1,2}|\graph(\tilde{p}_j^\alpha)|$ for some linear functions $\tilde{p}_j^1$, $\tilde{p}_j^2$ over $P_0$, and define new planes $L_j^1$, $L_j^2$ by
$$L_j^\alpha := \graph(\tilde{p}_j^\alpha + \hat{E}_{\widetilde{W}_j,\Dbf_j}(\mathcal{T}_j)\psi^\alpha) \qquad \text{for }\alpha=1,2,$$
then for the cones $\Dbf^*_j := |L_j^1| + |L_j^2|$, the above global $L^2$ convergence gives that (this is using Theorem \ref{thm:graphical-rep} for $\widetilde{W}_j$ and $\Dbf_j$)
\begin{equation}\label{E:graphical-rep-20}
	\frac{\hat{F}_{\widetilde{W}_j,\Dbf_j^*}(\mathcal{T}_j)}{\hat{E}_{\widetilde{W}_j,\Dbf_j}(\mathcal{T}_j)}\to 0.
\end{equation}
Notice also that since $\psi\not\equiv 0$ (from Step 4), we know that $\Dbf_j^*\neq \Dbf_j$ for all $j$. Moreover, as $S(\Dbf_j)\supset S_*$, \eqref{E:graphical-rep-18} tells us that $S(\Dbf_j^*)\supset S_*$ and so $\Dbf_j^*$ is a coarser cone than $\BC_j$. Thus, assuming for the moment that $\Dbf_j^*$ is \emph{not} a single plane, this means that $(\phi_j^{-1})_\#\Dbf_j^*$ is an admissible cone in the set over which the infimum is taken on the right-hand side of \eqref{E:graphical-rep-3}. This implies that
\begin{equation}\label{E:graphical-rep-21}
	\hat{E}_{W_j,\Dbf_j}(\mathcal{T}_j) \leq C\hat{F}_{W_j,(\phi_j^{-1})_\#\Dbf_j^*}(\mathcal{T}_j) \equiv C\hat{F}_{\widetilde{W}_j,\Dbf_j^*}(\mathcal{T}_j).
\end{equation}
Notice now that even if $\Dbf_j^*$ was a single plane, we can simply perturb it into a transversely intersecting pair of planes which still intersect along a subspace strictly containing $S$ and for which \eqref{E:graphical-rep-21} holds, or simply take a sequence of such transversely intersecting pairs of planes converging to $\Dbf_j^*$ in the right-hand side of \eqref{E:graphical-rep-3}: either way, we get \eqref{E:graphical-rep-21}, and so the form of $\Dbf_j^*$ is irrelevant. Then, combining using the triangle inequality and combining with \eqref{E:graphical-rep-21}, \eqref{E:graphical-rep-12}, and \eqref{E:graphical-rep-3}, we get
$$\hat{E}_{\widetilde{W}_j,\Dbf_j}(\mathcal{T}_j) \leq C\hat{F}_{\widetilde{W}_j,\Dbf_j^*}(\mathcal{T}_j).$$
However, this directly contradicts \eqref{E:graphical-rep-20}, giving the desired contradiction.

Currently, we have reached the desired contradiction when $\Dbf_j\in\mathcal{P}$ for all $j$, i.e.~when (i) of the original trichotomy from Step 1 held. To finish the contradiction, we need to deal with the remaining case, when $\Dbf_j = |P_0|$ for all $j$, i.e.~when (ii) of the original trichotomy held.

\textbf{Step 6:} \emph{Getting a contradiction in the planar case.} Now we are supposing that $\Dbf_j \equiv |P_0|$ for all $j$. In this case, we can take a coarse blow-up of $(W_j)_j$ off $P_0$ in $\mathcal{T}_0$; let $v:\mathcal{T}_0\cap P_0\to \A_2(P_0^\perp)$ denote this coarse blow-up.

Let $\tilde{c}_j:P_0\to \A_2(P_0^\perp)$ denote the two-valued linear function whose graph is exactly $(R_j)_\#\BC_j$. Using the second inequality in Lemma \ref{lemma:nu-a}(i) for $V_j$ and $\BC_j$ and Remark \ref{remark:elementary-observations}(1) from Section \ref{sec:further-notation}, up to passing to a subsequence the functions $\hat{E}_{W_j}(\mathcal{T}_j)^{-1}\tilde{c}_j$ converge to some two-valued linear function $\tilde{c}$. By the first inequality in Lemma \ref{lemma:nu-a}(i) (again applied to $V_j$ and $\BC_j$), we get that $\tilde{c}\neq 0$. Moreover, by (A4) for $V_j$ and $\BC_j$, it follows that we must have $v = \tilde{c}$; hence, we can write $v = \llbracket v^1\rrbracket + \llbracket v^2\rrbracket$, where $v^1,v^2$ are linear, and moreover we have $v\not\equiv 0$.

Since $S((R_j)_\#\BC_j)\equiv S\subset P_0$ for every $j$, we must then have $S\subset\graph(v^1)\cap \graph(v^2)$. Moreover, arguing similarly to as in Step 4, in light of the third point in Remark \ref{remark:after-blow-up-properties} following Theorem \ref{thm:blow-up-properties}, since there is a sequence of points $R_j(z_j)\in \mathcal{D}(W_j)\cap \mathcal{T}^A_{(1+\sigma)/2}(0,s_j^\prime)\cap \{r_S\geq\tau\}$, we find that there must be a point $z_0\in \overline{\mathcal{T}^A_{(1+\sigma)2}(0,s_0)}\setminus S$ for which $v^1(z_0) = v^2(z_0)$. Thus, if we define $L^\alpha_j:= \graph(\hat{E}_{W_j}(\mathcal{T}_j)v^\alpha)$ for $\alpha=1,2$ and $\Dbf_j^*:= |L_j^1| + |L_j^2|$, then we see that $S(\Dbf_j^*)\supsetneq S$. Moreover, as $v\not\equiv 0$ we know that $\Dbf^*_j\not\equiv 2|P_0|$ (we stress that $\Dbf_j^*$ could still be a multiplicity $2$ plane, it just cannot be $2|P_0|$).

We can then argue in an analogous manner to the argument in Step 5 used to verify \eqref{E:graphical-rep-20} to show that here we have
\begin{equation}\label{E:graphical-rep-22}
\hat{E}_{W_j}(\mathcal{T}_j)^{-1}\hat{E}_{W_j,\Dbf_j^*}(\mathcal{T}_j)\to 0.
\end{equation}
Consider first the case that $\Dbf_j^*$ is a single plane (with multiplicity $2$), i.e.~that $v^1\equiv v^2$. Then, by Theorem \ref{thm:allard-sup-estimate}, the two-sided and one-sided height excesses relative to $\Dbf_j^*$ are comparable, and so \eqref{E:graphical-rep-22} can be upgraded to
\begin{equation}\label{E:graphical-rep-23}
	\hat{E}_{W_j}(\mathcal{T}_j)^{-1}\hat{F}_{W_j,\Dbf_j^*}(\mathcal{T}_j)\to 0.
\end{equation}
If instead $v^1\not\equiv v^2$, then $\Dbf_j^*\in\mathcal{P}$. In this case, using the third point of Remark \ref{remark:after-blow-up-properties} following Theorem \ref{thm:blow-up-properties}, we must have that $\mathcal{D}(W_j)\cap \mathcal{T}_0$ converges to a closed subset which is contained within the coincidence set of $v$, which in this case is a subspace of dimension $\leq n-1$. In particular, for each $\eps>0$, for all $j$ large depending on $\eps$ we have $\mathcal{D}(W_j)\cap \mathcal{T}_0\subset B_{\eps}(S(\Dbf_j^*))\cap \mathcal{T}_0$. Hence, the defining property of the class $\mathcal{V}_\beta$ applies (through Lemma \ref{lemma:gap-2}) to give that again the two-sided and one-sided height excesses relative to $\Dbf_j^*$ are comparable (this is because now away from the spine, we know that the blow-up converges \emph{smoothly} to $v$, and thus $W_j$ must be close to each plane in $\Dbf_j^*$ and indeed can be expressed graphically over each plane in $\Dbf_j^*$). So, in either case we see that \eqref{E:graphical-rep-22} can be upgraded to \eqref{E:graphical-rep-23}.

But now just as in Step 5, since either $\Dbf_j^*$ is admissible (in the case $v^1\not\equiv v^2$) or is a limit of admissible cones (in the case $v^1\equiv v^2$) in the set over which the infimum is taken on the right-hand side of \eqref{E:graphical-rep-2}, we get from \eqref{E:graphical-rep-2}
$$\hat{E}_{W_j}(\mathcal{T}_j)\leq C\hat{F}_{W_j,\Dbf_j^*}(\mathcal{T}_j)$$
which directly violates \eqref{E:graphical-rep-23}. Thus, we get a contradiction in this case also.

\textbf{Step 7:} \emph{Concluding the proof of the theorem.} The result of Steps 1--6 is that we have now shown that for all $j$ sufficiently large we have
$$\mathcal{D}(V_j)\cap \{r_{\BC_j}\geq\tau\}\cap \mathcal{T}^A_{(1+\sigma)/2}(0,s_j^\prime) = \emptyset.$$
At this point, the desired conclusions of Theorem \ref{thm:graphical-rep} hold along a subsequence using the defining property of $\mathcal{V}_\beta$ (through Lemma \ref{lemma:gap-2}) and Allard's regularity theorem. Indeed, this follows in an analogous way to that seen in the proof of Theorem \ref{thm:graphical-rep} in the case $\sfrak(\BC)=n-1$ giving in Section \ref{sec:base-case}. Thus, we have now completed the inductive proof of Theorem \ref{thm:graphical-rep}.
\end{proof}

\begin{remark}\label{remark:after-graphical-rep-proof}
Hypothesis (A4) uses the two-sided fine excess on the left-hand side of the inequality; this was crucial in the proof of Theorem \ref{thm:graphical-rep}. In what follows, we will sometimes work under an assumption only on the \emph{one-sided} fine excess. We therefore remark now that the only places in the proof of Theorem \ref{thm:graphical-rep} in which we used the full strength of (A4) for $V_j$ and $\BC_j$ (i.e.~with the two-sided fine excess on the left-hand side rather than the one-sided excess) are in the verification of:
\begin{itemize}
	\item \eqref{E:graphical-rep-16}, in the case where $\Dbf_j\in \mathcal{P}$ for all $j$;
	\item When applying the second inequality in Lemma \ref{lemma:nu-a}(i) at the start of Step 6 to deduce $\nu((R_j)_\#\BC_j)\leq C\hat{E}_{W_j}(\mathcal{T}_j)$ in the case where $\Dbf_j \equiv |P_0|$ for all $j$.
\end{itemize}
Notice that both applications serve the same purpose, namely showing that we can blow-up the original sequence of cones off the coarser chosen cones (as indeed, we need to know that no pieces of the cone are far away from $V_j$, so we need to know closeness of the cone to $V_j$). Elsewhere, and in particular verifying \eqref{E:graphical-rep-19}, only requires that the \emph{one-sided} excess of $V_j$ relative to $\BC_j$ is much smaller than the one-sided fine excess relative the coarser cones. (One other place one might suspect a two-sided excess is needed is in verifying that the blow-up $\psi$ of the $V_j$ agrees with the blow-up $\chi$ of the cones (cf.~\eqref{E:graphical-rep-17}). However, under a one-sided assumption we would get that $\psi$ is contained in $\chi$, i.e.~both values could coincide with the same linear function, which is still sufficient for the rest of the argument.)  Thus, if we have by other means the inequalities in the bullet points above, then we can make do with a one-sided fine excess assumption rather than a two-sided one.
\end{remark}

Analogously to how the proof of Theorem \ref{thm:graphical-rep} in the case $\sfrak(\BC)=n-1$ gave us Corollary \ref{cor:weak-shift-estimate}, the inductive proof of Theorem \ref{thm:graphical-rep} gives us the following:

\begin{corollary}\label{cor:weak-shift-estimate-2}
	Let $\eta>0$. Then under the assumptions of Theorem \ref{thm:graphical-rep}, if $\eps_I$, $\gamma_I$ are also allowed to depend on $\eta$, then for every $z\in\mathcal{D}(V)\cap \mathcal{T}^A_{\sigma r}(x_0,s_0)$ we have
	$$r^{-1}\left(|\xi^{\perp_{P_0}}| + \mathscr{F}^*_{V,\BC}(\mathcal{T})|\xi^{\top_{P_0}}|\right) < \eta\mathscr{F}^*_{V,\BC}(\mathcal{T})$$
	where $\xi = z^{\perp_{S(\BC)}}$.
\end{corollary}

Indeed, this follows from the above inductive proof of Theorem \ref{thm:graphical-rep} because when taking the fine blow-up relative to the sequences of cones $\Dbf_j$, we not only know that $\hat{F}_{W_j,\Dbf_j}(\mathcal{T}_j)$ and $\hat{E}_{W_j,\Dbf_j}(\mathcal{T}_j)$ are comparable (from Remark \ref{remark:after-graphical-rep} after Theorem \ref{thm:graphical-rep}, which is applicable here due to the induction assumption), but also by choice of $\Dbf_j$ we know (cf.~\eqref{E:graphical-rep-2} and \eqref{E:graphical-rep-3}) that
$$\mathscr{F}^*_{W_j,\Dbf_j}(\mathcal{T}_j) \leq \hat{F}_{W_j,\Dbf_j}(\mathcal{T}_j) \leq C\mathscr{F}^*_{W_j,\BC_j}(\mathcal{T}_j).$$
Thus, when performing the relevant blow-up in the proof, the scale factor is comparable to $\mathscr{F}^*_{W_j,\BC_j}(\mathcal{T}_j)$, thus giving Corollary \ref{cor:weak-shift-estimate-2} as the good density points were forced to the spine during the blow-up.

\section{Proof of Theorem \ref{thm:L2-estimates} and Theorem \ref{thm:shift-estimate}}\label{sec:L2-shift}

Having now completed the inductive proof of Theorem \ref{thm:graphical-rep}, we move to proving inductively Theorem \ref{thm:L2-estimates} and Theorem \ref{thm:shift-estimate}. As already noted in Section \ref{sec:fine-L2}, these will rely on inductively proving Lemma \ref{lemma:L2-estimates-2} and Lemma \ref{lemma:Z}. As seen in Section \ref{sec:base-case}, a crucial step towards proving both these lemmas in the $\sfrak(\BC)=n-1$ case was Lemma \ref{lemma:tori-q=1}. Thus, we start with an analogue of Lemma \ref{lemma:tori-q=1} in the general setting.

\begin{lemma}\label{lemma:tori}
	Fix $\sigma\in (1/2,1)$, $a\in (0,1/2)$, $\beta\in (0,1)$, and $\tilde{\eps}\in (0,1)$. Then, there exists $\eps_0 = \eps_0(n,k,\sigma,a,\beta,\tilde{\eps})\in (0,1)$ such that the following is true. Suppose that we have $V\in\mathcal{V}_\beta$, an aligned cone $\BC = |P^1|+|P^2|\in\mathcal{P}_{\leq n-2}$, and $\mathcal{T}:= \mathcal{T}^{S(\BC)}_r(x_0,s_0)\subset B_1(0)$ with $s_0>0$ and $r/s_0\in (a,1-a)$, satisfying the following hypotheses:
	\begin{enumerate}
		\item [(a)] $(\w_n 2^n)^{-1}\|V\|(B^{n+k}_2(0)) \leq 5/2$;
		\item [(b)] $\nu(\BC) + \hat{E}_V + \hat{E}_V(\mathcal{T})<\eps_0$.
	\end{enumerate}
	Then one of the following conclusions holds for $\widetilde{V}:=V\res \mathcal{T}^{S(\BC)}_{\sigma r}(x_0,s_0)$:
	\begin{enumerate}
		\item [$(1\star)$] $\widetilde{V}$ is the union of at most two smooth minimal graphs, i.e.~$\widetilde{V} = V^1 + V^2$, where for each $\alpha\in \{1,2\}$, either $V^\alpha = 0$, or there exists a domain $\Omega^\alpha\supset P^\alpha\cap \mathcal{T}^{S(\BC)}_{\sigma r}(x_0,s_0)$ and a smooth function $u^\alpha:\Omega^\alpha\to (P^\alpha)^\perp$ solving the minimal surface system such that:
		\begin{enumerate}
			\item [(i)] $V^\alpha = \graph(u^\alpha)\res \mathcal{T}^{S(\BC)}_{\sigma r}(x_0,s_0)$;
			\item [(ii)] For $x\in\Omega^\alpha$ we have $\dist(x+u^\alpha(x),\BC) = |u^\alpha(x)|$;
		\end{enumerate}
	\end{enumerate}
	or
	\begin{enumerate}
		\item [$(2\star)$] $\widetilde{V}$ is well-approximated by a Lipschitz two-valued graph, i.e.~there is $\alpha\in \{1,2\}$, a domain $\Omega^\alpha\supset P^\alpha\cap \mathcal{T}^{S(\BC)}_{\sigma r}(x_0,s_0)$, a $\H^n$-measurable subset $\Sigma\subset\Omega^\alpha$, and $u:\Omega^\alpha\to \A_2((P^\alpha)^\perp)$ Lipschitz with $\Lip(u)\leq 1/2$ such that:
		\begin{enumerate}
			\item [(i)] $\widetilde{V}\res ((\Omega^\alpha\setminus\Sigma)\times (P^\alpha)^\perp) = \graph(u)\res ((\Omega\setminus\Sigma)\times (P^\alpha)^\perp)\cap \mathcal{T}_{\sigma r}^{S(\BC)}(x_0,s_0)$;
			\item [(ii)] $\H^n(\Sigma) + \|\widetilde{V}\|(\Sigma\times(P^\alpha)^\perp) \leq Cr^n\hat{E}_{V,\BC}(\mathcal{T})^2$;
			\item [(iii)] For each $x\in \Omega^\alpha\setminus\Sigma$ we have
                $$|u(x)|\leq Cr\hat{E}_{V,\BC}(\mathcal{T});$$
		\end{enumerate}
	\end{enumerate}
	or
	\begin{enumerate}
		\item [$(3\star)$] $V$ is equivalently approximated by a coarser transverse cone, i.e.~there exists $\Dbf\in\mathcal{P}$ coarser to $\BC$ with
		\begin{enumerate}
			\item [(i)] $\hat{E}_{V,\Dbf}(\mathcal{T}) \leq C\hat{E}_{V,\BC}(\mathcal{T})$; and
			\item [(ii)] There exists a rotation $R\in SO(n+k)$ with $R(\mathcal{T}) = \mathcal{T}$ such that $R_\#V$ and $R_\#\Dbf$ satisfy $\hat{E}_{R_\#V}(\mathcal{T})<\tilde{\eps}$ and $\nu(R_\#\Dbf)<\tilde{\eps}$.
		\end{enumerate}
	\end{enumerate}
	Here, $C = C(n,k,\sigma,a,\beta)$.
\end{lemma}

\begin{proof}
	We again start in a similar manner to that of the proof of Theorem \ref{thm:graphical-rep}. So, suppose that we have a sequence of numbers $\eps_j\downarrow 0$, varifolds $V_j\in\mathcal{V}_\beta$, cones $\BC_j = |P^1_j| + |P^2_j|\in \mathcal{P}_{\sfrak}$ with $\sfrak\leq n-2$, and regions $\mathcal{T}_j:= \mathcal{T}^{S(\BC_j)}_{r_j}(x_j,s_j)$ such that assumptions (a) and (b) of the present lemma are satisfied with $V_j$, $\BC_j$, $\eps_j$, and $\mathcal{T}_j$ in place of $V$, $\BC$, $\eps$, and $\mathcal{T}$ respectively. (We stress that we may inductively assume the validity of results in Section \ref{sec:fine-L2} for $\sfrak(\BC)>\sfrak$, as well as Theorem \ref{thm:graphical-rep} for $\sfrak(\BC) = \sfrak$.) Again, to prove the lemma we just need to show that we can always find a subsequence along which the conclusion holds.
	
	As we have done before, by performing appropriate rigid motions and dilations, we may assume that $S(\BC_j)\equiv S$ is a fixed $\sfrak$-dimensional subspace of $S_0$, and that $\mathcal{T}_j = \mathcal{T}^S_1(0,s_j^\prime)$, where $s_j^\prime\to s_0\in [1/(1-a),1/a]$. Write $\mathcal{T}_0:= \mathcal{T}_1^S(0,s_0)$.
	
	We now consider the following dichotomy:
	\begin{enumerate}
		\item [(A)] there is a subsequence of $\{j\}$ along which
		\begin{equation}\label{E:tori-1}
		\mathscr{E}_{V_j,\BC_j}^*(\mathcal{T}_j)^{-1}\hat{E}_{V_j,\BC_j}(\mathcal{T}_j) \to 0;
		\end{equation}
		\item [(B)] (A) fails, i.e.~for all sufficiently large $j$ we have
		\begin{equation}\label{E:tori-2}
		\hat{E}_{V_j,\BC_j}(\mathcal{T}_j) \geq C\mathscr{E}^*_{V_j,\BC_j}(\mathcal{T}_j)
		\end{equation}
		for some $C = C(n,k,\sigma,a,\beta)\in (0,\infty)$.
	\end{enumerate}
	The harder case is when (A) holds, so we will consider this first. So, suppose that (A) holds, and pass to a subsequence along which \eqref{E:tori-1} holds.
	
	We now pick an optimal coarser cone $\Dbf_j$ to $\BC_j$ in an identical manner to that done in Step 1 of the inductive proof of Theorem \ref{thm:graphical-rep} (with $A = S(\BC_j)\equiv S$ therein). The result is that, by performing appropriate rotations and passing to another subsequence, we can assume that there are coarser cones $\Dbf_j$ to $\BC_j$ where either:
	\begin{enumerate}
		\item [(A-i)] $\Dbf_j\in \mathcal{P}_{\sfrak^\prime}$ for every $j$, where $\sfrak^\prime>\sfrak$. Moreover, Hypothesis $H(S,\tilde{\eps}_j,\tilde{\gamma}_j,1)$ holds for $V_j$ and $\Dbf_j$ in $\mathcal{T}_j$ for suitable $\tilde{\eps}_j,\tilde{\gamma}_j\downarrow 0$. Furthermore, since by definition of $\mathscr{F}^*_{V_j,\BC_j}$ and Remark \ref{remark:after-graphical-rep} (following the statement of Theorem \ref{thm:graphical-rep}, applied inductively to $V_j$ and $\Dbf_j$, to give comparability of the two-sided and one-sided fine excess of $V_j$ relative to $\Dbf_j$) we have
		$$\mathscr{E}^*_{V_j,\BC_j}(\mathcal{T}_j) \leq \mathscr{F}^*_{V_j,\BC_j}(\mathcal{T}_j) \leq \hat{F}_{V_j,\Dbf_j}(\mathcal{T}_j) \leq C\hat{E}_{V_j,\Dbf_j}(\mathcal{T}_j)$$
		we get from \eqref{E:tori-1} that
		\begin{equation}\label{E:tori-4}
			\hat{E}_{V_j,\Dbf_j}(\mathcal{T}_j)^{-1}\hat{E}_{V_j,\BC_j}(\mathcal{T}_j) \to 0.
		\end{equation}
		In this case, we write $\tilde{u}^1_j$ and $\tilde{u}^2_j$ for the functions which graphically represent $V_j$ over $\Dbf_j$ (away from $S(\Dbf_j)$) in the sense of Theorem \ref{thm:graphical-rep}.
		\item [(A-ii)] $\Dbf_j = |P_0|$ for every $j$, $\hat{E}_{V_j}(\mathcal{T}_j) = \mathscr{E}_{V_j,\BC_j}(\mathcal{T}_j)$ (cf.~\eqref{E:graphical-rep-1}) and (cf.~\eqref{E:graphical-rep-2})
		\begin{equation}\label{E:tori-6}
			\hat{E}_{V_j}(\mathcal{T}_j) \leq C\mathscr{F}^*_{V_j,\BC_j}(\mathcal{T}_j).
		\end{equation}
	\end{enumerate}
	Here, $C = C(n,k,\sigma,a,\beta)$. Again, we will need to consider cases (A-i) and (A-ii) separately. We start with (A-i). In this case write $\Dbf_j = |L_j^1| + |L_j^2|$. Also, write $p_j^1$, $p_j^2$ for the linear functions over $P_0$ whose graphs are $P^1_j$, $P^2_j$, respectively, and similarly write $\ell^1_j$ and $\ell^2_j$ for the linear functions over $P_0$ whose graphs are $L_j^1$ and $L_j^2$, respectively.
	
	Firstly, we claim that after possibly relabelling the planes of $\Dbf_j$, we have
	\begin{equation}\label{E:tori-7}
		\nu(P_j^1,L_j^1) \leq C\hat{E}_{V_j,\Dbf_j}(\mathcal{T}_j)
	\end{equation}
	where $C = C(n,k,\sigma,a,\beta)$. To see this, suppose for contradiction that
	$$\hat{E}_{V_j,\Dbf_j}(\mathcal{T}_j)^{-1}\min_{\alpha,\alpha^\prime}\nu(P_j^\alpha,L_j^{\alpha^\prime})\to \infty.$$
	Then for $\tilde{\sigma}\in (7/8,1)$ we would have for all $j$ sufficiently large,
	\begin{align*}
		\hat{E}_{V_j,\BC_j}(\mathcal{T}_j)^2 & \geq \int_{\mathcal{T}^S_{\tilde{\sigma}}(0,s_0)\cap \{r_{\mathbf{D}_j}>1/8\}}\dist^2(x,\BC_j)\, \ext\|V_j\|(x)\\
		& \geq C^{-1}\sum_{\alpha=1,2}\int_{\mathcal{T}^S_{\tilde{\sigma}}(0,s_0)\cap \{r_{\mathcal{D}_j}>1/8\}\cap P_0}\min_{\alpha^\prime}|\ell_j^\alpha(x) + \tilde{u}_j^\alpha(x) - p_j^{\alpha^\prime}(x)|^2\, \ext x\\
		& \geq C^{-1}\sum_{\alpha=1,2}\int_{\mathcal{T}^S_{\tilde{\sigma}}(0,s_0)\cap \{r_{\mathcal{D}_j}>1/8\}\cap P_0}\min_{\alpha^\prime}(|\ell_j^\alpha-p_j^{\alpha^\prime}|^2 - |\tilde{u}^\alpha_j|^2)\, \ext x\\
		& \geq C^{-1}\min_{\alpha,\alpha^\prime}\nu(P_j^\alpha,L_j^\alpha)^2 - C\hat{E}_{V_j,\Dbf_j}(\mathcal{T}_j)^2.
	\end{align*}
	However, dividing both sides by $\hat{E}_{V_j,\Dbf_j}(\mathcal{T}_j)^2$, this gives a contradiction to \eqref{E:tori-4}. Thus, we have established \eqref{E:tori-7}.
	
	Next, we claim in addition to \eqref{E:tori-7} we also have
	\begin{equation}\label{E:tori-8}
		\nu(P_j^2,L_j^2) \leq C\hat{E}_{V_j,\Dbf_j}(\mathcal{T}_j).
	\end{equation}
	To prove this, suppose for contradiction that it was false, and that instead
	\begin{equation}\label{E:tori-9}
		\hat{E}_{V_j,\Dbf_j}(\mathcal{T}_j)^{-1}\nu(P_j^2,L_j^2)\to \infty.
	\end{equation}
	First note that Remark \ref{remark:elementary-observations}(4) from the end of Section \ref{sec:further-notation} (combined with a bound of the form $a^*(P^2_j\cup L^2_j)\geq c\nu(P^2_j,L_j^2)$) gives that there is some $(n-1)$-dimensional subspace $\tilde{S}\supsetneq S$ such that
	\begin{equation}\label{E:tori-10}
		|p_j^2(x)-\ell_j^2(x)|\geq c\nu(P^2_j,L_j^2) \qquad \text{for }x\in \mathcal{T}^S_{\sigma^\prime}(0,s_0)\cap P_0\setminus B_{1/8}(\tilde{S})
	\end{equation}
	for some constant $c = c(n,k)\in(0,1)$. Assume, without loss of generality, that by rotating we have $\tilde{S} = S_0$. Now, take $x\in \spt\|V\|$ such that $x\in \graph(\ell^2_j + \tilde{u}^2_j)\cap\mathcal{T}^{S}_{\sigma^\prime}(0,s_0)\cap \{|x^n|>1/8\}\cap \{r_{\Dbf_j}>1/8\}$, and let $x^\prime$ be a point in $\BC_j$ closest to $x$ (of course, $\{|x^n|>1/8\}\cap \{r_{\Dbf_j}>1/8\} \equiv \{|x^n|>1/8\}$, but for clarity we write both). Also, let $\bar{x}:= \pi_{P_0}(x)$, so that $x = \bar{x}+\ell_j^2(\bar{x}) + \tilde{u}^2_j(\bar{x})$. If $x^\prime\in P_j^2$, then there is a constant $c = c(n,k)\in (0,1)$ such that
	$$\dist(x,\BC_j) = |x-x^\prime| \geq c|x-(\bar{x}+p^2_j(\bar{x}))| \geq c|\ell^2_j(\bar{x})-p^2_j(\bar{x})| - |\tilde{u}^2_j(\bar{x})|.$$
	From the conclusions of Theorem \ref{thm:graphical-rep} (for $V_j$ and $\Dbf_j$) we also know that $\sup|\tilde{u}_j^\alpha|\leq C\hat{E}_{V_j,\Dbf_j}(\mathcal{T}_j)$. Combining all of this with \eqref{E:tori-10} and \eqref{E:tori-9}, we deduce that
	\begin{equation}\label{E:tori-11}
		\hat{E}_{V_j,\Dbf_j}(\mathcal{T}_j)^{-1}\dist(x,\BC_j)\to \infty \qquad \text{for all such }x.
	\end{equation}
	We now want to deduce \eqref{E:tori-11} also in the case that $x^\prime\in P^1_j$. Indeed, if on the other hand $x^\prime\in P^1_j$, then similarly
	$$\dist(x,\BC_j) = |x-x^\prime| \geq c|x-(\bar{x}+p^1_j(\bar{x}))| \geq c|\ell^2_j(\bar{x})-\ell^1_j(\bar{x})| - |\ell^1_j(\bar{x})-p^1_j(\bar{x})| - |\tilde{u}^2_j(\bar{x})|.$$
	From Lemma \ref{lemma:nu-a}(ii) (applied to $V_j$ and $\Dbf_j$) and the definition of $a_*(\Dbf_j)$ (as $r_{\Dbf_j}(x)>1/8$) we know that
	$$|\ell^2_j(\bar{x})-\ell^1_j(\bar{x})| > \frac{1}{8}a_*(\Dbf_j) \geq c\mathscr{F}^*_{V_j,\Dbf_j}(\mathcal{T}_j).$$
	But from \eqref{E:tori-7}, we have $|\ell_j^1(\bar{x})-p^1_j(\bar{x})|\leq C\hat{E}_{V_j,\Dbf_j}(\mathcal{T}_j)$. Therefore, using again the estimate $\sup|\tilde{u}^\alpha_j|\leq C\hat{E}_{V_j,\Dbf_j}(\mathcal{T}_j)$, we get that
	$$\dist(x,\BC_j) \geq c\mathscr{F}^*_{V_j,\Dbf_j}(\mathcal{T}_j) - C\hat{E}_{V_j,\Dbf_j}(\mathcal{T}_j).$$
	Using then Hypothesis (A4) for $V_j$ and $\Dbf_j$, we can then conclude that \eqref{E:tori-11} holds in this case also. Hence, we have now shown that \eqref{E:tori-11} holds for all $x\in G_j:=\graph(\ell^2_j+\tilde{u}^2_j)\cap \mathcal{T}^S_{\sigma^\prime}(0,s_0)\cap \{|x^n|>1/8\}\cap \{r_{\Dbf_j}>1/8\}$. (Note that in fact we actually proved a stronger statement than \eqref{E:tori-11}, namely that the infimum of this quantity over all such $x$ is $\to\infty$ as the lower bound was uniform in the $x$ variable.) However, since
	$$\|V_j\|(G_j)\geq c(n,k)>0$$
	this implies (by a simple infimum lower bound on the integral in $\hat{E}_{V_j,\BC_j}(\mathcal{T}_j)$) that
	$$\hat{E}_{V_j,\Dbf_j}(\mathcal{T}_j)^{-1}\hat{E}_{V_j,\BC_j}(\mathcal{T}_j)\to \infty.$$
	However, this is in direct violation of \eqref{E:tori-4}. This contradiction therefore establishes \eqref{E:tori-8}. Combining \eqref{E:tori-8} with \eqref{E:tori-7}, we have therefore shown
	\begin{equation}\label{E:tori-12}
		\nu(\BC_j,\Dbf_j) \leq C\hat{E}_{V_j,\Dbf_j}(\mathcal{T}_j).
	\end{equation}
	We are now in a position to use Remark \ref{remark:after-graphical-rep-proof} following the inductive proof of Theorem \ref{thm:graphical-rep} (i.e.~at the end of Section \ref{sec:graphical-rep}). Indeed, the inequality \eqref{E:tori-12} tells us that we can blow-up $\BC_j$ off $\Dbf_j$ rescaling by $\hat{E}_{V_j,\Dbf_j}(\mathcal{T}_j)$. As we can also take a fine blow-up of $V_j$ off $\Dbf_j$, recalling \eqref{E:tori-4}, we see that the fine blow-up of $V_j$ off $\Dbf_j$ is contained within the blow-up of $\BC_j$ off $\Dbf_j$; in particular, it is linear. However, we can now use this blow-up to argue in an analogous way to that seen in Steps 4, 5, and 6 of the inductive proof of Theorem \ref{thm:graphical-rep} to deduce that in this case, $(1\star)$ holds (indeed, we can assume in our argument, when looking for a contradiction, that for any $\tau>0$, there is $z_j\in \mathcal{D}(V_j)\cap \{r_{\BC_j}>\tau\}\cap \mathcal{T}^S_{(1+\sigma)/2}(0,s_0)$ lying on $S(\Dbf_j)$, as in Step 3 of the proof of Theorem \ref{thm:graphical-rep}. We can follow the argument therein to get a contradiction. Then, choosing $\tau = \sigma/4$ suffices for the present lemma). Thus, the proof is complete in the case of (A-i), with conclusion $(1\star)$ holding with $V^1,V^2$ both non-zero smooth minimal graphs (one close to each plane in $\BC_j$, due to the choice of $\Dbf_j$).
	
	Next suppose that we are in the case where (A-ii) holds. Then, similarly to what we saw in the case of (A-i), we may show that
	\begin{equation}\label{E:tori-13}
		\nu(P^1_j) \leq C\hat{E}_{V_j}(\mathcal{T}_j)
	\end{equation}
	(indeed, note that \eqref{E:tori-1} still gives the appropriate form of \eqref{E:tori-4} for this situation, namely that $\hat{E}_{V_j}(\mathcal{T}_j)^{-1}\hat{E}_{V_j,\BC_j}(\mathcal{T}_j)\to 0$, which is what is needed to replicate the argument leading to \eqref{E:tori-7}).
	
	There are then two cases:
	\begin{enumerate}
		\item [(a)] $\nu(P^2_j)\leq C\hat{E}_{V_j}(\mathcal{T}_j)$ for some $C = C(n,k,\sigma,a,\beta)$; or
		\item [(b)] $\hat{E}_{V_j}(\mathcal{T}_j)^{-1}\nu(P^2_j)\to \infty$.
	\end{enumerate}
	We consider each case in turn.
	
	In the case that (a) holds, then combining with \eqref{E:tori-13} we have $\nu(\BC_j)\leq C\hat{E}_{V_j}(\mathcal{T}_j)$. In this case, using Remark \ref{remark:after-graphical-rep-proof} following the inductive proof of Theorem \ref{thm:graphical-rep} at the end of Section \ref{sec:graphical-rep}, we can follow the proof as in Step 6 of the inductive proof of Theorem \ref{thm:graphical-rep} (note that the first inequality in Lemma \ref{lemma:nu-a}(i) only needs $\hat{E}_{V,\BC}(\mathcal{T})<\gamma\hat{E}_V(\mathcal{T})$, which we have in the present case); thus again enables us to rule out points of density $\geq 2$ in the desired region, and so the proof is completed as in the proof of Theorem \ref{thm:graphical-rep}, leading to conclusion $(1\star)$ of the present lemma (again, with both $V^1, V^2$ therein being non-zero).
	
	If instead (b) holds, this implies that the graph of the coarse blow-up coincides with that of the coarse blow-up of $p_j^1$. Thus, we have from the local $L^2$ convergence to the coarse blow-up,
	$$\hat{E}_{V_j}(\mathcal{T}_j)^{-1}\hat{E}_{V_j,P_j^1}(\mathcal{T}^S_{(1+\sigma)/2}(0,s_0^\prime))\to 0.$$
	By Allard's supremum estimate (Theorem \ref{thm:allard-sup-estimate}) we know that
	$$\sup_{x\in \spt\|V_j\|\cap \mathcal{T}^S_{(1+3\sigma)/4}(0,s_0^\prime)}\dist(x,P_j^1)\leq C\hat{E}_{V_j,P_j^1}(\mathcal{T}^S_{(1+\sigma)/2}(0,s_0^\prime)).$$
	We now claim that there is a fixed constant $c$ for which 
	\begin{equation}\label{E:tori-14}
	a_*(\BC_j)\geq c\hat{E}_{V_j}(\mathcal{T}_j).
	\end{equation}
	Once we have this, then for $j$ sufficiently large from the above it follows that for all $x\in \spt\|V_j\|\cap \mathcal{T}^S_{(1+3\sigma)/4}(0,s_0)$ we have (as such $x$ are a fixed distance away from $S$)
	$$\dist(x,\BC_j) = \dist(x,P_j^1).$$
	This shows that, depending on whether $V_j$ converges to $P_0$ with multiplicity one or two, we can conclude either $(1\star)$ (with one of the $V^i$ being $0$) or $(2\star)$ by graphically approximating $V_j$ over $P_j^1$ using either Allard's regularity theorem or Theorem \ref{thm:Lipschitz-approx}, as the above gives $\hat{E}_{V_j,P_j^1}(\mathcal{T}^S_{(1+3\sigma)/4}(0,s_0) \equiv \hat{E}_{V_j,\BC_j}(\mathcal{T}^S_{(1+3\sigma)/4}(0,s_0)) \leq C\hat{E}_{V_j,\BC_j}(\mathcal{T}_j)$. Thus, once we have shown \eqref{E:tori-14}, this would conclude the proof when (A) holds.
	
	To see \eqref{E:tori-14}, first note that the same argument as in Lemma \ref{lemma:nu-a}(ii), except now only using one-sided excesses instead of two-sided excesses in the argument therein, and using \eqref{E:tori-1}, gives that
	$$\mathscr{E}^*_{V_j,\BC_j}(\mathcal{T}_j) \leq Ca_*(\BC_j).$$
	We now wish to argue that in the present case we can improve \eqref{E:tori-6} to
	\begin{equation}\label{E:tori-15}
		\hat{E}_{V_j}(\mathcal{T}_j) \leq C\mathscr{E}^*_{V_j,\BC_j}(\mathcal{T}_j).
	\end{equation}
	Once we have this, we can combine the previous two inequalities to get $\hat{E}_{V_j}(\mathcal{T}_j) \leq Ca_*(\BC_j)$, giving \eqref{E:tori-14}. Thus, it suffices to prove \eqref{E:tori-15}. To see \eqref{E:tori-15}, first choose $\widetilde{\Dbf}_j\in \mathcal{P}$ coarser to $\BC_j$ for which 
	$$\hat{E}_{V_j,\widetilde{\Dbf}_j}(\mathcal{T}_j) \leq \frac{3}{2}\mathscr{E}^*_{V_j,\BC_j}(\mathcal{T}_j).$$
	Observe that $\widetilde{\Dbf}_j$ obeys
	\begin{equation}\label{E:tori-16}
		\nu(\widetilde{\Dbf}_j,P^1_j) \leq C\hat{E}_{V_j,\widetilde{\Dbf}_j}(\mathcal{T}_j).
	\end{equation}
	Indeed, \eqref{E:tori-16} follows because a sufficiently large region of $V_j$ is represented by a graph over $\widetilde{\Dbf}_j$ by the inductive use of the present lemma (Lemma \ref{lemma:tori}, using Lemma \ref{lemma:tori-q=1} instead in the base case), and as we can upper bound $\nu(\widetilde{\Dbf}_j,P_j^1)$ by the $L^2$ norm of the difference of the linear functions defining $\widetilde{\Dbf}_j$ and $P_j^1$ and use the triangle inequality, using also the fact that, on compact subsets of $C_1(0)$, we have $L^2$ convergence to the blow-up which is a single, linear function, and so $V_j$ is close to $P_j^1$ there, and that we have \eqref{E:tori-1}. But now notice that by the construction of $\Dbf_j$ leading to case (A-ii), we have
	\begin{equation}\label{E:tori-17}
	\frac{\hat{F}_{V_j,\widetilde{\Dbf}_j}(\mathcal{T}_j)}{\mathscr{F}^*_{V_j,\widetilde{\Dbf}_j}(\mathcal{T}_j)} \geq \frac{\mathscr{F}^*_{V_j,\BC_j}(\mathcal{T}_j)}{C\hat{E}_{V_j}(\mathcal{T}_j)} \geq C^{-1}
	\end{equation}
	where we have used the various definitions of the infimums here as well as \eqref{E:tori-6}. But we also have, writing $\widetilde{\Dbf}_j = |\widetilde{L}_1| + \widetilde{L}_2|$,
	\begin{align*}
		\left(\frac{\hat{F}_{V_j,\widetilde{\Dbf}_j}(\mathcal{T}_j)}{\mathscr{F}^*_{V_j,\widetilde{\Dbf}_j}(\mathcal{T}_j)}\right)^2 & \leq \frac{\hat{E}_{V_j,\widetilde{\Dbf}_j}(\mathcal{T}_j)^2 + \sum_{i=1,2}\int_{\mathcal{T}^S_{1/2}(0,s^\prime_j)}\dist^2(x,V_j)\, \ext\|\widetilde{L}_i\|(x)}{\mathscr{F}^*_{V_j,\widetilde{\Dbf}_j}(\mathcal{T}_j)^2}\\
		& \leq \frac{\hat{E}_{V_j,\widetilde{\Dbf}_j}(\mathcal{T}_j)^2 + 2\int_{\mathcal{T}^S_{1/2}(0,s^\prime_j)}\dist^2(x,V_j)\, \ext\|P_j^1\|(x) + C\nu(\widetilde{\Dbf}_j,P_j^1)^2}{\mathscr{F}^*_{V_j,\widetilde{\Dbf}_j}(\mathcal{T}_j)^2}\\
		& \leq C\frac{\hat{E}_{V_j,\widetilde{\Dbf}_j}(\mathcal{T}_j)^2}{\mathscr{F}^*_{V_j,\widetilde{\Dbf}_j}(\mathcal{T}_j)^2} + \frac{2\int_{\mathcal{T}^S_{1/2}(0,s^\prime_j)}\dist^2(x,V_j)\, \ext\|P_j^1\|(x)}{\mathscr{F}^*_{V_j,\widetilde{\Dbf}_j}(\mathcal{T}_j)^2}.
	\end{align*}
	Now using the bound $\mathscr{E}^*_{V_j,\widetilde{\Dbf}_j}(\mathcal{T}_j) \leq \mathscr{F}^*_{V_j,\widetilde{\Dbf}_j}(\mathcal{T}_j)$ in the first term and $\mathscr{F}^*_{V_j,\widetilde{\Dbf}_j} \geq \mathscr{F}^*_{V_j,\BC_j}(\mathcal{T}_j) \geq C^{-1}\hat{E}_{V_j}(\mathcal{T}_j)$ in the second, as well as Allard's supremum estimate (Theorem \ref{thm:allard-sup-estimate}) in the numerator of the second term, we get
	$$\left(\frac{\hat{F}_{V_j,\widetilde{\Dbf}_j}(\mathcal{T}_j)}{\mathscr{F}^*_{V_j,\widetilde{\Dbf}_j}(\mathcal{T}_j)}\right)^2 \leq C\frac{\hat{E}_{V_j,\widetilde{\Dbf}_j}(\mathcal{T}_j)^2}{\mathscr{E}^*_{V_j,\widetilde{\Dbf}_j}(\mathcal{T}_j)^2} + C\frac{\hat{E}_{V_j,P_j^1}(\mathcal{T}^S_{3/4}(0,s_j^\prime))^2}{\hat{E}_{V_j}(\mathcal{T}_j)^2}.$$
	The second term here $\to 0$ as $j\to\infty$ due to the strong convergence in $L^2$ on (fixed) compact subsets of $\mathcal{T}_j$ to the blow-up. Thus, combining the above with \eqref{E:tori-17}, we get that for all $j$ sufficiently large,
	$$C^{-1} \leq\frac{\hat{E}_{V_j,\widetilde{\Dbf}_j}(\mathcal{T}_j)}{\mathscr{E}^*_{V_j,\widetilde{\Dbf}_j}(\mathcal{T}_j)}$$
	and thus, using the definition of $\widetilde{\Dbf}_j$, we therefore get
	$$\mathscr{E}^*_{V_j,\widetilde{\Dbf}_j}(\mathcal{T}_j) \leq C\mathscr{E}^*_{V_j,\BC_j}(\mathcal{T}_j).$$
	We can now repeat this argument, choosing an optimal cone for $\mathscr{E}^*_{V_j,\widetilde{\Dbf}_j}(\mathcal{T}_j)$, and so forth, to inductively increase the spine dimension of the cone optimal cone chosen within, until we reach the spine dimension $n$, i.e.~planes. At this point, as $P_0$ is chosen to be the optimal plane, we get $\hat{E}_{V_j}(\mathcal{T}_j) \equiv \mathscr{E}_{V_j,\BC_j}(\mathcal{T}_j) \leq C\mathscr{E}^*_{V_j,\BC_j}(\mathcal{T}_j)$, which is exactly \eqref{E:tori-15}, thus completing the proof of the lemma in the case that (A) holds.
	
	We finally come to the case of (B). After passing to a subsequence along which \eqref{E:tori-2} holds, we then follow again the method of choosing an optimal coarser cone $\Dbf_j$ for each $j$ as was done in Step 1 of the inductive proof of Theorem \ref{thm:graphical-rep}, except now we choose an optimal cone for the \emph{one-sided} fine excess $\mathscr{E}^*_{V_j,\BC_j}$ (this is done analogously). There are then two cases depending on whether $\Dbf_j\equiv |P_0|$ for every $j$ or if $\Dbf_j\in \mathcal{P}$ for every $j$. Notice that if $\Dbf_j\in\mathcal{P}$, then by combining the corresponding property from \eqref{E:graphical-rep-3} (which would now be $\hat{E}_{W_j,\Dbf_j}(\mathcal{T}_j)\leq C\mathscr{E}^*_{W_j,\BC_j}(\mathcal{T}_j)$) with \eqref{E:tori-2}, we immediately get that conclusion $(3\star)$ of the present lemma holds. 
	
	So, we are left with the case that $\Dbf_j = |P_0|$ for all $j$. Combining the corresponding property from \eqref{E:graphical-rep-2} (which would now be $\hat{E}_{W_j}(\mathcal{T}_j)\leq C\mathscr{E}^*_{W_j,\BC_j}(\mathcal{T}_j)$) with \eqref{E:tori-2}, we get $\hat{E}_{V_j}(\mathcal{T}_j) \leq C\hat{E}_{V_j,\BC_j}(\mathcal{T}_j)$. If we now apply Lemma \ref{lemma:appendix-1}, we see that there is a plane in $\BC_j$, without loss of generality by relabelling say it is $P_j^1$, such that
	$$\nu(P_j^1) \leq C(\hat{E}_{V_j,\BC_j}(\mathcal{T}_j) + \hat{E}_{V_j}(\mathcal{T}_j)).$$
	So, using the above we get $\nu(P_j^1) \leq C\hat{E}_{V_j,\BC_j}(\mathcal{T}_j)$, and hence the triangle inequality gives $\hat{E}_{V_j,P_j^1}(\mathcal{T}_j) \leq 2\hat{E}_{V_j}(\mathcal{T}_j) + C\nu(P_j^1) \leq C\hat{E}_{V_j,\BC_j}(\mathcal{T}_j)$. Thus, if $V_j$ is close to $P_j^1$ with multiplicity $2$, we can apply Theorem \ref{thm:Lipschitz-approx} to represent $V_j$ as a two-valued graph over $P_j^1$ and thus immediately deduce that $(2\star)$ holds in this case, otherwise it must be close to $P_j^1$ with multiplicity one, and thus we may apply Allard's regularity theory to deduce that $(1\star)$ holds in this case (with $V^2=0$ therein). This then completes the proof of the lemma.
\end{proof}

Given Lemma \ref{lemma:tori}, we are now in a position to prove Lemma \ref{lemma:Z}.

\begin{proof}[Proof of Lemma \ref{lemma:Z}]
	Let $\eps_0$ be as in the statement of Lemma \ref{lemma:tori} for suitable choices of $\tilde{\eps}$ and $\tilde{\gamma}$ to be determined depending only on $n,k,\sigma,a,\beta$. Notice also that we may assume without loss of generality that $\phi_{(0)}$ is non-negative; indeed, the general case then follows by splitting $\phi_{(0)}$ as $\phi_{(0)}^+ - \phi_{(0)}^-$ where $\phi_{(0)}^+:= \max\{\phi_{(0)},0\}$ and $\phi^-_{(0)} := \max\{-\phi_{(0)},0\}$, which from our assumptions on $\phi_{(0)}$ still obey all the assumptions of Lemma \ref{lemma:Z} (which follows from, for instance, the identity $\phi_{(0)}^+ = \phi_{(0)}/2 + |\phi_{(0)}|/2$).
	
	There are two cases to consider: (A) $\hat{E}_{V,\BC}(\mathcal{T})<\eps_0$ or (B) $\hat{E}_{V,\BC}(\mathcal{T})\geq\eps_0$. Proving the lemma when (B) holds is identical to the corresponding proof in the $\sfrak(\BC)=n-1$ case seen in Section \ref{sec:base-case}. Thus, we focus on case (A). Recall that we are assuming $\BC\in \mathcal{P}_\sfrak$ with $\sfrak\leq n-2$.
	
	When (A) holds, we may apply Lemma \ref{lemma:tori}. In the case that either conclusion $(1\star)$ or $(2\star)$ of Lemma \ref{lemma:tori} holds, then we can follow completely analogous arguments to those seen in the proof of the base case to conclude the lemma. Thus, we just need to consider the case when $(3\star)$ of Lemma \ref{lemma:tori} holds.
	
	Let $\Dbf\in \mathcal{P}$ be the cone whose existence is asserted by $(3\star)$ of Lemma \ref{lemma:tori}. We first claim that
	\begin{equation}\label{E:Z-1}
		\int_{B_1(0)}\phi_{(0)}(x)|x^{\perp_{S(\BC)}}|^2\, \ext\|\BC\|(x) = \int_{B_1(0)}\phi_{(0)}(x)|x^{\perp_{S(\BC)}}|^2\, \ext\|\Dbf\|(x).
	\end{equation}
	Indeed, this follows from the fact that the integrand is a function of only $x^{\top_{S(\BC)}}$ and $|x^{\perp_{S(\BC)}}|$. For the sake of completeness, we will prove this by expanding into appropriate cylindrical coordinates. Let us assume that we have chosen coordinates so that, if $m = \sfrak(\BC)$ and $m^\prime = \sfrak(\Dbf)$,
	$$P_0 = \textnormal{span}\langle e_1,\dotsc,e_n\rangle,\ \ \ \ S(\BC) = \textnormal{span}\langle e_1,\dotsc,e_{m}\rangle,\ \ \ \ S(\Dbf) = \textnormal{span}\langle e_1,\dotsc,e_{m^\prime}\rangle,$$
	and let $g:\R^m\times [0,\infty)\to \R$ be the function such that
	$$g(x^{\top_{S(\BC)}},|x^{\perp_{S(\BC)}}|):= \phi_{(0)}(x)|x^{\perp_{S(\BC)}}|^2.$$
	Now let $L$ be any $n$-dimensional plane containing $S(\BC)$. Then the coarea formula and the geometry of the region $\mathcal{T}$ give (also recalling that $\phi_{(0)}$ has compact support in $\mathcal{T}$)
	\begin{align*}
		\int_L \phi_{(0)}(x)|x^{\perp_{S(\BC)}}|^2\, \ext\H^n(x) & = \int_{L\cap \mathcal{T}}g(x^{\top_{S(\BC)}},|x^{\perp_{S(\BC)}}|)\, \ext\H^n(x)\\
		& = \int_{\mathcal{T}\cap S(\BC)}\int_{\mathcal{T}\cap S(\BC)^\perp\cap L}g(y,|\xi|)\, \ext\H^{n-m}(\xi)\ext\H^m(y)\\
		& = 2\w_{n-m-1}\int_{[x_0-r,x_0+r]^m}\int^{s_0+r}_{s_0-r}g(y,\rho)\, \ext\rho\ext y
	\end{align*}
	where the factor of $2$ comes from the fact that $\mathcal{T}\cap S(\BC)^\perp\cap L$ has two connected components and the integral over each connected component is the same. Since this final expression does not depend on $L$, if we take $L$ to be the planes in $\BC$ and $\Dbf$, \eqref{E:Z-1} follows.
	
	Also note that using simple linear algebra based on the above projection maps (for example, the fact that the map $\pi_{S(\BC)^\perp}-\pi_{S(\Dbf)^{\perp}}$ is the projection onto $\textnormal{span}\langle e_{m+1},\dotsc,e_{m^\prime}\rangle$, which is contained within $S(\BC)^\perp$) we have
	\begin{align*}
		|x^{\perp_{S(\BC)}\perp_{T_xV}}|^2 & \leq 2|\pi_{T_x^\perp V}(\pi_{S(\Dbf)^\perp} (x))|^2 + 2|\pi_{T_x^\perp V}(\pi_{S(\BC)^\perp}-\pi_{S(\Dbf)^\perp})(x)|^2\\
		& \leq 2|x^{\perp_{S(\Dbf)}\perp_{T_xV}}|^2 + |x_{m+1}\pi_{T_x^\perp V}(e_{m+1})+\cdots + x_{m^\prime}\pi_{T_x^\perp V}(e_{m^\prime})|^2\\
		& \leq 2|x^{\perp_{S(\Dbf)}\perp_{T_xV}}|^2 + n2^n|x^{\perp_{S(\BC)}}|^2\sum^{m^\prime}_{i=1}|e_i^{\perp_{T_xV}}|^2
	\end{align*}
	Therefore, based on this, \eqref{E:Z-1}, and (i) of $(3\star)$ of Lemma \ref{lemma:tori}, to prove the present lemma it suffices to prove the estimates
	\begin{equation}\label{E:Z-2}
		\left|\int\phi_{(0)}(x)|x^{\perp_{S(\BC)}}|^2\, \ext\|V\| - \int\phi_{(0)}(x)|x^{\perp_{S(\BC)}}|^2\, \ext\|\Dbf\|\right| \leq C\int_{\mathcal{T}}\dist^2(x,\Dbf)\, \ext\|V\|(x);
	\end{equation}
	\begin{equation}\label{E:Z-3}
		\int_{\mathcal{T}^{S(\BC)}_{\sigma r}(x_0,s_0)}|x^{\perp_{S(\Dbf)}\perp_{T_xV}}|^2\, \ext\|V\|(x) \leq C\int_{\mathcal{T}}\dist^2(x,\Dbf)\, \ext\|V\|(x);
	\end{equation}
	\begin{equation}\label{E:Z-4}
		\int_{\mathcal{T}^{S(\BC)}_{\sigma r}(x_0,s_0)}|x^{\perp_{S(\BC)}}|^2\sum^{\sfrak(\Dbf)}_{i=1}|e_i^{\perp_{T_xV}}|^2\, \ext\|V\|(x) \leq C\int_{\mathcal{T}}\dist^2(x,\Dbf)\, \ext\|V\|(x);
	\end{equation}
	(note that we cannot simply apply Lemma \ref{lemma:Z} inductively to $V$ and $\Dbf$ to show \eqref{E:Z-3}, as the domain of integration in \eqref{E:Z-3} is over $\mathcal{T}_{\sigma r}^{S(\BC)}(x_0,s_0)$ and not $\mathcal{T}^{S(\Dbf)}_{\sigma r}(x_0,s_0)$ which is what is needed in Lemma \ref{lemma:Z}).
	
	To prove these, we first make a computation based on the first variation formula. For any $f\in C^\infty_c(\mathcal{T})$ (to be chosen shortly) consider the vector field
	$$\Phi(x) = f(x)^2x^{\perp_{S(\Dbf)}}.$$
	Then using $\div_{T_xV}\Phi(x) := \sum_{j=1}^{n+k} e_j\cdot\nabla^V\Phi^j(x) \equiv \sum_{j=1}^{n+k}e_j\cdot\pi_{T_xV}(D\Phi^j(x))$ and noting that $\Phi^i \equiv 0$ for $i\leq m^\prime\equiv \sfrak(\Dbf)$, we have
	\begin{align}
		\nonumber\div_{T_xV}\Phi(x) & = f(x)^2\sum_{j=m^\prime+1}^{n+k}e_j\pi_{T_xV}(e_j) + \sum_{j=m^\prime+1}^{n+k}x_j e_j\cdot\pi_{T_xV}(D(f(x)^2))\\
		&\nonumber = f(x)^2\textnormal{trace}_{S(\Dbf)^\perp}(\pi_{T_xV}) + \pi_{T_xV}(x^{\perp_{S(\Dbf)}})\cdot D(f(x)^2)\\
		&\nonumber = f(x)^2\text{trace}_{S(\Dbf)^\perp}(\pi_{T_xV})\\
		&\nonumber \hspace{5em} + \pi_{T_xV}(x^{\perp_{S(\Dbf)}})\cdot \pi_{S(\Dbf)}(D(f(x)^2)) + \pi_{T_xV}(x^{\perp_{S(\Dbf)}})\cdot \pi_{S(\Dbf)^\perp}(D(f(x)^2))\\
		& =: \text{I} + \text{II} + \text{III}.\label{E:Z-5}
	\end{align}
	where by $\text{trace}_S(\pi)$ for a subspace $S$ and projection map $\pi$, we simply mean the trace of $\pi$ but only over the directions in $S$. We analyse each of these three terms in turn. Firstly,
	\begin{align*}
	\text{I}:= f(x)^2\text{trace}_{S(\Dbf)^\perp}(T_xV) & = f(x)^2(n - \text{trace}_{S(\Dbf)}(\pi_{T_xV}))\\
	& = f(x)^2(n-\sfrak(\Dbf) + \text{trace}_{S(\Dbf)}(\pi_{T_x^\perp V})).
	\end{align*}
	Since $\text{trace}_{S(\Dbf)}(\pi_{T^\perp_xV}) = \sum_{i=1}^{m^\prime}e_i\cdot \pi_{T_x^\perp V}(e_i) = \sum_{i=1}^{m^\prime}|e_i^{\perp_{T_xV}}|^2$, we therefore have
	\begin{equation}\label{E:Z-6}
		\text{I} = f(x)^2\left(n-\sfrak(\Dbf) + \sum^{m^\prime}_{i=1}|e_i^{\perp_{T_xV}}|^2\right).
	\end{equation}
	For the second term, replacing $\pi_{T_xV}$ by $\id_{\R^{n+k}}-\pi_{T_x^\perp V}$ and then noting that $x^{\perp_{S(\Dbf)}}\cdot \pi_{S(\Dbf)}\equiv 0$ gives
	\begin{align*}
		\text{II} := \pi_{T_xV}(x^{\perp_{S(\Dbf)}})\cdot\pi_{S(\Dbf)}(D(f(x)^2)) = -\pi_{T_x^\perp V}(x^{\perp_{S(\Dbf)}})\cdot\pi_{T_x^\perp V}\pi_{S(\Dbf)}(D(f(x)^2)).
	\end{align*}
	Now since the operator norm of the map $\pi_{T_x^\perp V}:S(\Dbf)\to S(\Dbf)\cap T_x^\perp V$ obeys
	$$\|\pi_{T_x^\perp V}|_{S(\Dbf)}\| \leq \left(\sum^{m^\prime}_{j=1}|e_j^{\perp_{T_xV}}|^2\right)^{1/2}$$
	we have
	\begin{align}
		\nonumber\text{II} & \geq -|\pi_{T_x^\perp V}(x^{\perp_{S(\Dbf)}})|\cdot |\pi_{S(\Dbf)}(D(f(x)^2))|\cdot\left(\sum^{m^\prime}_{j=1}|e_j^{\perp_{T_xV}}|^2\right)^{1/2}\\
		\nonumber& = -2|x^{\perp_{S(\Dbf)}\perp_{T_xV}}|\cdot |f(x)|\cdot |\pi_{S(\Dbf)}(Df(x))|\cdot \left(\sum^{m^\prime}_{j=1}|e_j^{\perp_{T_xV}}|^2\right)^{1/2}\\
		& \geq -2|x^{\perp_{S(\Dbf)}\perp_{T_xV}}|^2|\pi_{S(\Dbf)}(Df(x))|^2 - \frac{1}{2}f(x)^2\sum^{m^\prime}_{j=1}|e_j^{\perp_{T_xV}}|^2\label{E:Z-7}
	\end{align}
	where in the last line we have used $2ab\leq 2a^2+\frac{1}{2}b^2$ for suitable $a,b\in \R$.
	
	For the third term, we will need to explicitly differentiate $f$ in the $S(\Dbf)^\perp$ directions. So at this point, let us choose the $f$ we will use. Let
	$$f(x)^2 := \phi_{(0)}(x)|x^{\perp_{S(\BC)}}|^2$$
	(note that this is a valid choice, as we are assuming $\phi_{(0)}$ is non-negative). In particular, we may think of $f$ as a function of $x^{\top_{S(\BC)}}$ and $|x^{\perp_{S(\BC)}}|$. For simplicity, let us write $\phi^\prime_{(0)}$ for the derivative of $\phi_{(0)}$ with respect to the (1-dimensional) $|x^{\perp_{S(\BC)}}|$ input variable. So, using the product rule and noting that $S(\Dbf)^\perp\subset S(\BC)^\perp$, we have
	\begin{align*}
		\text{III} & := \pi_{T_xV}(x^{\perp_{S(\Dbf)}})\cdot\pi_{S(\Dbf)^\perp}D(f(x)^2)\\
		& \equiv \pi_{T_xV}(x^{\perp_{S(\Dbf)}})\cdot\pi_{S(\Dbf)^\perp} D(\phi_{(0)}(x)|x^{\perp_{S(\BC)}}|^2)\\
		& = \pi_{T_xV}(x^{\perp_{S(\Dbf)}})\cdot\left(\phi^\prime_{(0)}|x^{\perp_{S(\BC)}}|x^{\perp_{S(\Dbf)}} + 2\phi_{(0)}x^{\perp_{S(\Dbf)}}\right).
	\end{align*}
	Writing $\pi_{T_xV}(x^{\perp_{S(\Dbf)}}) = x^{\perp_{S(\Dbf)}} - x^{\perp_{S(\Dbf)}\perp_{T_xV}}$, we deduce that
	\begin{equation}\label{E:Z-8}
		\text{III} = \left(|x^{\perp_{S(\Dbf)}}|^2 - |x^{\perp_{S(\Dbf)}\perp_{T_xV}}|^2\right)\psi_{(1)}(x)
	\end{equation}
	where $\psi_{(1)}(x):= \phi_{(0)}^\prime(x)|x^{\perp_{S(\BC)}}| + 2\phi_{(0)}(x)$. For future reference, notice that:
	\begin{enumerate}
		\item [(P1)] $\psi_{(1)}(x)$ is a smooth function of just $x^{\top_{S(\Dbf)}}$ and $|x^{\perp_{S(\Dbf)}}|$. To see this, note that because $S(\BC)\subsetneq S(\Dbf)$ we have
		$$|x^{\top_{S(\BC)}}| = \pi_{S(\BC)}(x^{\top_{S(\Dbf)}}),$$
		$$|x^{\perp_{S(\BC)}}|^2 = |x^{\perp_{S(\Dbf)}}|^2 + |\pi_{S(\BC)^\perp}(x^{\top_{S(\Dbf)}})|^2.$$
		Thus, away from $S(\BC)$ (which in $\mathcal{T}$ we are) we have that $x^{\top_{S(\BC)}}$ and $|x^{\perp_{S(\BC)}}|$ are smooth functions of $x^{\top_{S(\Dbf)}}$ and $|x^{\perp_{S(\Dbf)}}|$, and thus from the form of $\psi_{(1)}$ and the fact that $\phi_{(0)}$ is a smooth function of $x^{\top_{S(\BC)}}$ and $|x^{\perp_{S(\BC)}}|$, the result follows.
		\item [(P2)] $|\psi_{(1)}|\leq \tilde{C}$ and $|D^\alpha\psi_{(1)}|\leq \tilde{C}r^{-|\alpha|}$ for $|\alpha|\leq n-\sfrak(\Dbf)$. This follows from the definition of $\psi_{(1)}$ and the corresponding properties for $\phi_{(0)}$.
	\end{enumerate}
Thus, combining the first variation formula with \eqref{E:Z-5}, \eqref{E:Z-6}, \eqref{E:Z-7}, and \eqref{E:Z-8}, we get the inequality
\begin{align*}
	0\geq & \int f(x)^2\left(n-\sfrak(\Dbf)+\sum^{m^\prime}_{j=1}|e_j^{\perp_{T_xV}}|^2\right)\, \ext\|V\|(x)\\
	& + \int -2|x^{\perp_{S(\Dbf)}\perp_{T_xV}}|^2|\pi_{S(\Dbf)}(Df(x))|^2 - \frac{1}{2}f(x)^2\sum^{m^\prime}_{j=1}|e_j^{\perp_{T_xV}}|^2\, \ext\|V\|(x)\\
	& + \int \left(|x^{\perp_{S(\Dbf)}}|^2 - |x^{\perp_{S(\Dbf)}\perp_{T_xV}}|^2\right)\psi_{(1)}(x)\, \ext\|V\|(x).
\end{align*}
Rearranging this and substituting in the choice of $f$, we get
\begin{align}
	\nonumber (n-\sfrak(\Dbf))\int&\phi_{(0)}|x^{\perp_{S(\BC)}}|^2\, \ext\|V\|(x) + \frac{1}{2}\int\phi_{(0)}|x^{\perp_{S(\BC)}}|^2\sum^{m^\prime}_{j=1}|e_j^{\perp_{T_xV}}|^2\, \ext\|V\|(x)\\
	& \leq \int\left(\psi_{(1)} + 2|\pi_{S(\Dbf)}(Df(x))|^2\right)|x^{\perp_{S(\Dbf)}\perp_{T_xV}}|^2\, \ext\|V\|(x) - \int\psi_{(1)}|x^{\perp_{S(\Dbf)}}|^2\, \ext\|V\|(x).\label{E:Z-9}
\end{align}
Note now that $f$ is Lipschitz: Indeed, since $f = \phi_{(0)}^{1/2}|x^{\perp_{S(\BC)}}|$, to prove this claim we just need to analyse $\phi_{(0)}^{1/2}$. This is Lipschitz because $\phi_{(0)}$ is non-negative and smooth, so whenever $\phi_{(0)}$ vanishes so must $D\phi_{(0)}$. More precisely, a Glaeser-type inequality gives for each direction $e_j$,
$$|D_j(\phi_{(0)}^{1/2})| \leq C\sqrt{\sup|D^2\phi_{(0)}|} \leq Cr^{-1}$$
	and so combining this with the product rule we readily get $|Df|\leq C$. In particular, $|\pi_{S(\Dbf)}(Df)|\leq C$. Using this in \eqref{E:Z-9} as well as the bound $|\psi_{(1)}|\leq C$, we therefore have
	\begin{align}
		\nonumber(n-\sfrak(\Dbf))\int\phi_{(0)}&|x^{\perp_{S(\BC)}}|^2\, \ext\|V\|(x) + \frac{1}{2}\int\phi_{(0)}|x^{\perp_{S(\BC)}}|^2\sum^{m^\prime}_{j=1}|e_j^{\perp_{T_xV}}|^2\, \ext\|V\|(x)\\
		& \leq C\int_{\mathcal{T}_{\sigma r}^{S(\BC)}(x_0,s_0)}|x^{\perp_{S(\Dbf)}\perp_{T_xV}}|^2\, \ext\|V\|(x) - \int\psi_{(1)}(x)|x^{\perp_{S(\Dbf)}}|^2\, \ext\|V\|(x). \label{E:Z-10}
	\end{align}
	Next we try to find expressions for the relevant quantities appearing in \eqref{E:Z-10} for the cone $\Dbf$; this ultimately amounts to an integration by parts. Write $\Dbf = S(\Dbf)\times\Dbf_0 = \R^{m^\prime}\times\Dbf_0$, where $\Dbf_0$ is a pair of $q:= n-\sfrak(\Dbf)$ dimensional planes in $\R^{q+k}$ that intersect only at the origin. By the coarea formula we have (as $f(x) = f(x^{\top_{S(\BC)}},|x^{\perp_{S(\BC)}}|)$,
	\begin{align*}
		\int_{B_1(0)}f(x)^2\, \ext\|\Dbf\|(x) = \int_{S(\Dbf)}\int_{\{y\}\times S(\Dbf)^\perp}f\left(y^{\perp_{S(\BC)}},\sqrt{|y^{\perp_{S(\BC)}}|^2 + |\xi|^2}\right)^2\, \ext\|\Dbf_0\|(\xi)\, \ext\H^{m^\prime}(y).
	\end{align*}
	Switching the inner integral to spherical polar coordinates, the right-hand side of the above equals
	$$\int_{S(\Dbf)}\int_{\{y\}\times(S^{q+k-1}\cap\spt\|\Dbf_0\|)}\ext\w^{q-1}\int^\infty_0 f\left(y^{\top_{S(\BC)}},\sqrt{|y^{\perp_{S(\BC)}}|^2 + \rho^2}\right)^2\cdot\rho^{q-1}\, \ext\rho\ext\H^{m^\prime}(y).$$
	Now perform a one-dimensional integration by parts in the $\rho$ variable. For notational simplicity, we will now drop the arguments of $f$ and $f^\prime$ (we recall that $f^\prime$ means only the derivative in the $|x^{\perp_{S(\BC)}}|$ variable). Since
	$$\frac{\del}{\del\rho}f\left(y^{\top_{S(\BC)}},\sqrt{|y^{\perp_{S(\BC)}}|^2+\rho^2}\right)^2 = 2\rho(|y^{\perp_{S(\BC)}}|^2+\rho^2)^{-1/2}f^\prime f$$
	we have that
	$$\int^\infty_0 f\left(y^{\top_{S(\BC)}},\sqrt{|y^{\perp_{S(\BC)}}|^2 + \rho^2}\right)^2\cdot\rho^{q-1}\, \ext\rho = -\frac{2}{q}\int^\infty_0 \rho^2(|y^{\perp_{S(\BC)}}|^2 + \rho^2)^{-1/2}ff^\prime\cdot\rho^{q-1}\, \ext\rho$$
	(where we have rearranged so that the polar Jacobian factor $\rho^{q-1}$ still appears explicitly). Transforming back out of polar coordinates, recalling that $\rho = |x^{\perp_{S(\Dbf)}}|$ and $|x^{\perp_{S(\BC)}}| = \sqrt{|y^{\perp_{S(\BC)}}|^2+\rho^2}$, the result is that combining all the above we have now shown
	$$q\int_{B_1(0)}f(x)^2\, \ext\|\Dbf\|(x) = -\int_{B_1(0)}2f(x)f^\prime(x)|x^{\perp_{S(\BC)}}|^{-1}|x^{\perp_{S(\Dbf)}}|^2\, \ext\|\Dbf\|(x).$$
	Writing $\rho_1:= |x^{\perp_{S(\BC)}}|$ momentarily, as a function of just $\rho_1$ we have $f(\rho_1)^2 = \phi_{(0)}(\rho_1)\rho_1^2$. Thus, as $f^\prime$ denotes the derivative of $f$ in $\rho_1$, we have from differentiating this that $2ff^\prime = \phi_{(0)}^\prime\rho_1^2 + 2\phi_{(0)}\rho_1$. Hence, $2ff^\prime|x^{\perp_{S(\BC)}}|^{-1} = \psi_{(1)}$. Hence, the above becomes
	\begin{equation}\label{E:Z-11}
		(n-\sfrak(\Dbf))\int_{B_1(0)}\phi_{(0)}|x^{\perp_{S(\BC)}}|^2\, \ext\|\Dbf\|(x) = -\int_{B_1(0)}\psi_{(1)}|x^{\perp_{S(\Dbf)}}|^2\, \ext\|\Dbf\|(x).
	\end{equation}
	Combining \eqref{E:Z-11} with \eqref{E:Z-10}, we get
	\begin{align}
		\nonumber&(n-\sfrak(\Dbf))\left(\int\phi_{(0)}|x^{\perp_{S(\BC)}}|^2\, \ext\|V\|(x) - \int\phi_{(0)}|x^{\perp_{S(\BC)}}|^2\, \ext\|\Dbf\|(x)\right)\\
		&\label{E:Z-11.5} \hspace{2em} + \frac{1}{2}\int\phi_{(0)}|x^{\perp_{S(\BC)}}|^2\sum^{m^\prime}_{j=1}|e_j^{\perp_{T_xV}}|^2\, \ext\|V\|(x)\\
		&\nonumber\hspace{4em}\leq C\int_{\mathcal{T}_{\sigma r}^{S(\BC)}(x_0,s_0)}|x^{\perp_{S(\Dbf)}\perp_{T_xV}}|^2\, \ext\|V\|(x)\\
		&\nonumber\hspace{6em} -\left( \int\psi_{(1)}(x)|x^{\perp_{S(\Dbf)}}|^2\, \ext\|V\|(x) -\int\psi_{(1)}|x^{\perp_{S(\Dbf)}}|^2\, \ext\|\Dbf\|(x)\right).
	\end{align}
	Combining this with \eqref{E:Z-2}, \eqref{E:Z-3}, and \eqref{E:Z-4}, we see that it suffices to prove the estimates (writing $\phi_{(1)}:= -\psi_{(1)}$)
	\begin{equation}\label{E:Z-12}
		\left| \int\phi_{(1)}|x^{\perp_{S(\Dbf)}}|^2\, \ext\|V\|(x) -\int\phi_{(1)}|x^{\perp_{S(\Dbf)}}|^2\, \ext\|\Dbf\|(x)\right| \leq C\int_{\mathcal{T}}\dist^2(x,\Dbf)\, \ext\|V\|(x);
	\end{equation}
	\begin{equation}\label{E:Z-13}
		\int_{\mathcal{T}^{S(\BC)}_{\sigma r}(x_0,s_0)}|x^{\perp_{S(\Dbf)}\perp_{T_xV}}|^2\, \ext\|V\|(x) \leq C\int_{\mathcal{T}}\dist^2(x,\Dbf)\, \ext\|V\|(x);
	\end{equation}
	\begin{equation}\label{E:Z-14}
		\int_{\mathcal{T}^{S(\BC)}_{\sigma r}(x_0,s_0)}|x^{\perp_{S(\BC)}}|^2\sum^{\sfrak(\Dbf)}_{i=1}|e_i^{\perp_{T_xV}}|^2\, \ext\|V\|(x) \leq C\int_{\mathcal{T}}\dist^2(x,\Dbf)\, \ext\|V\|(x).
	\end{equation}
	Note that \eqref{E:Z-13} is the same as \eqref{E:Z-3}, and \eqref{E:Z-14} is the same as \eqref{E:Z-4}. Notice that once we prove \eqref{E:Z-12} and \eqref{E:Z-13}, then returning to \eqref{E:Z-11.5}, we see that \eqref{E:Z-14} actually follows; indeed, if one changes $\sigma$ to $(1+\sigma)/2$ in the above, and lets $\phi_{(0)}$ be non-negative such that $\phi_{(0)}\equiv 1$ on $\mathcal{T}_{\sigma}^{S(\BC)}(x_0,s_0)$ and $|D^\alpha\phi_{(0)}|\leq C(1-\sigma)^{-|\alpha|}r^{-|\alpha|}$ for all multi-indices $\alpha$, then \eqref{E:Z-14} follows immediately from \eqref{E:Z-11.5}. Thus, in fact it suffices to prove \eqref{E:Z-12} and \eqref{E:Z-13}.
	
	The two estimates \eqref{E:Z-12} and \eqref{E:Z-13} are ones we can tackle with our inductive assumption. Indeed, to prove these we follow a cubic decomposition analogous to that seen in the proof of the base case of Lemma \ref{lemma:Z}.
	
	For $\tilde{\sigma}\in (0,1]$, write $\Omega_\sigma:= \pi_{S(\Dbf)}(\mathcal{T}^{S(\BC)}_{\sigma r}(x_0,s_0))$. By repeating an essentially analogous construction to that which appears in the proof of the base case of Lemma \ref{lemma:L2-estimates-2}, we may find a countable and pairwise disjoint family $\mathscr{F}$ of dyadic cubes in $\Lambda(S(\Dbf))$ such that
	$$\Omega_\sigma\subset\bigcup_{Q\in\mathscr{F}}Q\subset\Omega_{(1+\sigma)/2},$$
	$$\mathcal{T}^{S(\BC)}_{\sigma r}(x_0,s_0)\setminus S(\Dbf)\subset\bigcup_{Q\in\mathscr{F}}\mathcal{T}^{S(\Dbf)}(\tfrac{5}{4}Q)\subset \mathcal{T}^{S(\BC)}_{r}(x_0,s_0),$$
	and satisfying properties analogous to $(\mathscr{F}1)$ -- $(\mathscr{F}7)$ therein. In particular, this means:
	\begin{enumerate}
		\item [(a)] There is a constant $c^\prime = c^\prime(n,\sfrak(\Dbf),\sigma)\in (0,1/2)$ such that if $Q = Q_\rho(y_0,t_0)\in\mathscr{F}$, then $\rho/t_0\in (c^\prime,1-c^\prime)$ (and $t_0>0$ always).
		\item Each point of $B_1(0)$ is contained in at most $N = N(\sfrak(\Dbf))<\infty$ cubes from the collection $\{\tfrac{5}{4}Q:Q\in\mathscr{F}\}$.
		\item For each cube $Q\in\mathscr{F}$, there is a smooth function $\phi_Q:\Omega_1\to [0,1]$ such that
		\begin{itemize}
			\item $\spt(\phi_Q)\subset\tfrac{9}{8}Q$;
			\item $\sum_Q \phi_Q\equiv 1$ on $\bigcup_{Q\in\mathscr{F}}Q$;
			\item $|D^\alpha\phi_Q(y,t)|\leq Ct^{-|\alpha|}$ for all multi-indices $\alpha$, where $C = C(n,|\alpha|,\sigma)\in (0,\infty)$.
		\end{itemize}
	\end{enumerate}
	Then for each $Q\in\mathscr{F}$, let $\tilde{\phi}_Q:B_1(0)\to [0,1]$ be the smooth function given by
	$$\tilde{\phi}_Q(x):= \phi_Q(x^{\top_{S(\Dbf)}},|x^{\perp_{S(\Dbf)}}|).$$
	By construction we then have $\sum_{Q\in\mathscr{F}}\tilde{\phi}_Q \equiv 1$ for all $x\in \mathcal{T}_{\sigma r}^{S(\BC)}(x_0,s_0)$, and $|D^\alpha\tilde{\phi}(x)|\leq C|x^{\perp_{S(\Dbf)}}|^{-|\alpha|}$ for $C = C(n,|\alpha|,\sigma)$ and $x\in \mathcal{T}^{S(\BC)}_{r}(x_0,s_0)$. We can now complete the proofs of \eqref{E:Z-12} and \eqref{E:Z-13} in an analogous way to the base case. Indeed, for \eqref{E:Z-13} simply note that
	\begin{align*}
		\int_{\mathcal{T}_{\sigma r}^{S(\BC)}(x_0,s_0)}|x^{\perp_{S(\Dbf)}\perp_{T_xV}}|^2\, \ext\|V\|(x) & = \sum_{Q\in\mathscr{F}}\int_{\mathcal{T}^{S(\BC)}_{\sigma r}(x_0,s_0)}\tilde{\phi}_Q|x^{\perp_{S(\Dbf)}\perp_{T_xV}}|^2\, \ext\|V\|(x)\\
		& \leq \sum_{Q\in\mathscr{F}}\int_{\mathcal{T}^{S(\Dbf)}(\frac{9}{8}Q)}|x^{\perp_{S(\Dbf)}\perp_{T_xV}}|^2\, \ext\|V\|(x).
	\end{align*}
	Now for each cube $Q\in\mathcal{F}$ we can apply the inductive hypothesis to apply Lemma \ref{lemma:Z} for $V$ and $\Dbf$: this can be done by property (ii) of $(3\star)$ from Lemma \ref{lemma:tori}, provided $\tilde{\eps}$ is sufficiently small. Thus, for each $Q\in\mathscr{F}$ we get
	$$\int_{\mathcal{T}^{S(\Dbf)}(\frac{9}{8}Q)}|x^{\perp_{S(\Dbf)}\perp_{T_xV}}|^2\, \ext\|V\|(x) \leq C\int_{\mathcal{T}^{S(\Dbf)}(\frac{5}{4}Q)}\dist^2(x,\Dbf)\, \ext\|V\|(x).$$
	Then using the bounded overlap property of the cubes $\mathscr{F}$ to get that when we sum over $Q\in\mathscr{F}$, the right-hand side is still controlled by
	$$C\int_{\bigcup_{Q\in\mathscr{F}}\mathcal{T}^{S(\Dbf)}(\frac{5}{4}Q)}\dist^2(x,\Dbf)\, \ext\|V\|(x) \leq C\int_{\mathcal{T}^{S(\BC)}_r(x_0,s_0)}\dist^2(x,\Dbf)\, \ext\|V\|(x).$$
	Combining the above, we deduce \eqref{E:Z-13}.
	
	For \eqref{E:Z-12}, notice that for each $Q = Q_\rho(y_0,t_0)\in\mathscr{F}$, in light of (P1) and (P2) previously, the function $\tilde{\phi}_{(0)}:= \phi_{(1)}\tilde{\phi}_Q$ is in $C^\infty_c(\mathcal{T}^{S(\Dbf)}_{9\rho/8}(y_0,t_0))$ and satisfies the hypotheses needed for its inductive use in Lemma \ref{lemma:Z} to $V$ and $\Dbf$. Thus, applying \eqref{E:Z-b} for each cube $Q\in\mathscr{F}$ with this function, we get
	$$\left|\int\phi_{(1)}\tilde{\phi}_Q|x^{\perp_{S(\Dbf)}}|^2\, \ext\|V\| - \int\phi_{(1)}\tilde{\phi}_Q|x^{\perp_{S(\Dbf)}}|^2\, \ext\|\Dbf\|\right| \leq C\int_{\mathcal{T}^{S(\Dbf)}(\frac{5}{4}Q)}\dist^2(x,\Dbf)\, \ext\|V\|(x).$$
	Now using the triangle inequality and summing this over $Q\in\mathscr{F}$ and again using the bounded overlap property, we reach \eqref{E:Z-12}. This therefore completes the inductive proof of Lemma \ref{lemma:Z}.
\end{proof}

We can now complete the inductive proofs of the remaining lemmas and theorems. However, given Lemma \ref{lemma:tori} and Lemma \ref{lemma:Z}, these are now very similar to the proofs in the $\sfrak(\BC)=n-1$ case.

\begin{proof}[Proof of Lemma \ref{lemma:L2-estimates-2}]
	Recalling Remark \ref{remark:L2-estimates-2-proof} after the proof of the base case of Lemma \ref{lemma:L2-estimates-2} in Section \ref{sec:base-case}, given Lemma \ref{lemma:Z} the proof of Lemma \ref{lemma:L2-estimates-2} follows in an analogous way to the base case. Thus, we have proved Lemma \ref{lemma:L2-estimates-2}.
\end{proof}

\begin{proof}[Proof of Theorem \ref{thm:L2-estimates}]
	Given the inductive proofs of Lemma \ref{lemma:L2-estimates-2} and Theorem \ref{thm:graphical-rep}, the proof of Theorem \ref{thm:L2-estimates} follows in an analogous manner to that of the base case seen in Section \ref{sec:base-case}. The only extra remark needed is that instead of using Corollary \ref{cor:weak-shift-estimate} within the proof, we instead use Corollary \ref{cor:weak-shift-estimate-2}.
\end{proof}

\begin{proof}[Proof of Theorem \ref{thm:shift-estimate}]
	Recalling Remark \ref{remark:shift-estimate-proof} following the proof of the base case of Theorem \ref{thm:shift-estimate} in Section \ref{sec:base-case}, given the inductive proofs of all the preceding results the proof of Theorem \ref{thm:shift-estimate} now follows in an analogous way to the proof of the base case. Thus we have now proved Theorem \ref{thm:shift-estimate}.
\end{proof}

Consequently, this concludes our inductive argument and thus we have now unconditionally proved Theorem \ref{thm:graphical-rep}, Theorem \ref{thm:L2-estimates}, Theorem \ref{thm:shift-estimate}, and Corollary \ref{cor:L2-corollary}. Our next aim will be to use these results to understand the higher regularity of the fine blow-ups constructed in Section \ref{sec:fine-blow-ups} (which we now know we can always construct).

\section{Fine Blow-Ups Part II: Regularity}\label{sec:fine-blow-ups-2}

With the main estimates from Theorem \ref{thm:graphical-rep}, Theorem \ref{thm:L2-estimates}, Theorem \ref{thm:shift-estimate}, and Corollary \ref{cor:L2-corollary} at our disposal, we now return to the discussion of fine blow-ups we started in Section \ref{sec:fine-blow-ups}. Before, in Theorem \ref{thm:fine-continuity-2}, we established a Hölder continuity result for general fine blow-ups. Here, we will use our estimates to improve this to a generalised-$C^{1,\alpha}$ regularity result.

The hardest type of fine blow-up to understand is when the fine blow-up is in $\mathfrak{b}_{n-1}$, i.e.~when we are taking a fine blow-up relative a sequence of cones with spine dimension $n-1$. In fact, the regularity of fine blow-ups in $\mathfrak{b}_\sfrak$, $\sfrak\leq n-2$, follows directly from removability results for continuous harmonic functions (which we know is the case for fine blow-ups via Theorem \ref{thm:fine-continuity-2}) across sets with vanishing $2$-capacity, such as a subspace of dimension $\leq n-2$. For $\mathfrak{b}_{n-1}$, fine blow-ups are instead comprised of $4$ harmonic functions $\phi^1,\dotsc,\phi^4$ each defined on half-planes rather than a full-plane, and so the regularity question becomes one of boundary regularity rather than that of removable singularities. As such, many of our results will be focused on the case of $\mathfrak{b}_{n-1}$.

\subsection{Additional Integral Identities for Fine Blow-Ups in $\mathfrak{b}_{n-1}$}

Our first lemma consists of two integral identities for fine blow-ups in $\mathfrak{b}_{n-1}$, which will be important steps to understanding their boundary behaviour.

\begin{lemma}\label{lemma:fine-b}
	Let $\phi = (\phi^\alpha)_{\alpha=1,\dotsc,4}\in\mathfrak{b}_{n-1}$. Then for any $\zeta\in C^\infty_c(B^n_{3/8}(0))$ with $D_{e_n}\zeta\equiv 0$ in a neighbourhood of $S\equiv S_0$ and $\int_{S\cap B^n_{3/8}(0)}\zeta\, \ext\H^{n-1}=0$, we have
	\begin{equation}\label{E:fine-b-a}
		\sum^4_{\alpha=1}\int_{B^n_{3/8}(0)\cap U^\alpha_0}\phi^{\alpha,\kappa}\Delta\zeta\, \ext x = 0 \qquad \text{for any }\kappa\in \{1,\dotsc,k\};
	\end{equation}
	\begin{equation}\label{E:fine-b-b}
		\sum^4_{\alpha=1}\int_{B^n_{3/8}(0)\cap U^\alpha_0}\left(\sum^k_{\kappa=1}m^\alpha_{\kappa,n}\phi^{\alpha,\kappa}\right)\Delta\zeta = 0.
	\end{equation}
\end{lemma}

\textbf{Remark:} We stress that the subscript $n$ in \eqref{E:fine-b-b} is the direction \emph{orthogonal} to the boundary $S_0$ of the half-planes in $P_0$. Here, we write $\phi^\alpha = (\phi^{\alpha,1},\dotsc,\phi^{\alpha,k})$, and we are always viewing the fine blow-up in $\mathfrak{b}_{n-1}$ as four functions over half-planes, even if it is possible to view them as two functions over $P_0$.

\begin{proof}
	Let $\zeta$ be as in the statement of the lemma, and let $\tau>0$ be such that $D_{e_n}\zeta(x)=0$ in $\{|x^n|<2\tau\}$. Recall that $\widetilde{\zeta}$ denotes the vertical extension of $\zeta$. Let $(V_j)_j$, $(\BC_j)_j$ be the sequences whose fine blow-up is $\phi$. The first variation formula \eqref{E:stationarity-2} then gives
	$$0 = \int_{C_1(0)}\nabla^{V_j}x^{n+\kappa}\cdot\nabla^{V_j}\widetilde{\zeta}(x)\, \ext\|V_j\|(x)$$
	for each $\kappa=1,\dotsc,k$ and each $j$. Using the notation from Theorem \ref{thm:graphical-rep} and writing $f_j^\alpha = h_j^\alpha + u_j^\alpha$, this becomes
	\begin{align}
		\nonumber \sum_{\alpha=1}^4\int_{\R^k\times(B^n_{3/8}(0)\cap \{|x^n|>2\tau\}\cap U^\alpha_0)}\nabla^{V_j}\widetilde{f}_j^{\alpha,\kappa}&\cdot\nabla^{V_j}\widetilde{\zeta}\, \ext\|V_j\|\\
		 & = -\int_{\R^k\times(B^n_{3/8}(0)\cap \{|x^n|<2\tau\})}\nabla^{V_j}x^{n+\kappa}\cdot\nabla^{V_j}\widetilde{\zeta}\, \ext\|V_j\|.\label{E:fine-b-1}
	\end{align}
	Consider the term on the right-hand side of \eqref{E:fine-b-1}. Using the fact that in a neighbourhood of $\spt\|V\|$ we have that $\widetilde{\zeta}$ is constant in the vertical directions, and that $D_{e_n}\tilde{\zeta}(x) = 0$ in $\{|x^n|<2\tau\}$, we have
	\begin{align*}
		&\left|\int_{\R^k\times(B^n_{3/8}(0)\cap \{|x^n|<2\tau\})}\nabla^{V_j}x^{n+\kappa}\cdot\nabla^{V_j}\widetilde{\zeta}\, \ext\|V_j\|(x)\right|\\
		& \leq \sup_{B^n_{3/8}(0)}|D\zeta|\left|\int_{\R^k\times(B^n_{3/8}(0)\cap\{|x^n|<2\tau\})}e_{n+\kappa}\cdot\sum^{n-1}_{i=1}\pi_{T_xV_j}(e_i)\, \ext\|V_j\|(x)\right|\\
		& \leq \sup|D\zeta|\sum^{n-1}_{i=1}\left(\int_{\R^k\times(B^n_{3/8}(0)\cap\{|x^n|<2\tau\})}|e_i^{\perp_{T_xV_j}}|^2\, \ext\|V_j\|(x)\right)^{1/2}\\
		& \hspace{18em}\times\left(\|V_j\|(\R^k\times (B^n_{3/8}(0)\cap\{|x^n|<2\tau\}))\right)^{1/2}\\
		& \leq C\sup|D\zeta|\tilde{E}_j\cdot\sqrt{\tau}
	\end{align*}
	where $C = C(n,k,M,\beta)$, where in the last line we have used \eqref{E:cor-4} from Corollary \ref{cor:L2-corollary} to control the first term in the region close to $\mathcal{D}(V_j)$ (with the choice of $\eta$ therein independent of $\tau$, e.g.~$\eta = 1/1000$) and the regularity away from $B_{2\tau}(\mathcal{D}(V_j))$ provided in the fine blow-up procedure to control the part of the integral away from $\mathcal{D}(V_j)$, and we have controlled the mass term using the monotonicity formula. Combining this with \eqref{E:fine-b-1} we therefore now have
	\begin{equation}\label{E:fine-b-2}
		\left|\sum_{\alpha=1}^4\int_{\R^k\times(B^n_{3/8}(0)\cap \{|x^n|>2\tau\}\cap U^\alpha_0)}\nabla^{V_j}\widetilde{f}_j^{\alpha,\kappa}\cdot\nabla^{V_j}\widetilde{\zeta}\, \ext\|V_j\|\right| \leq C\sup|D\zeta|\tilde{E}_j\sqrt{\tau}.
	\end{equation}
	Let us now turn to the left-hand side of \eqref{E:fine-b-2}. For this, we morally follow computations similar to \cite[(12.17) -- (12.21)]{Wic14}. Writing $J^\alpha_j$ for the Jacobian of the map $x\mapsto (f^\alpha_j(x),x)$ and $M^\alpha_j := \graph(f^\alpha_j)$, we can express this term as an integral over the domain $P_0$ as:
	$$\sum_{\alpha=1}^4\int_{B^n_{3/8}(0)\cap \{|x^n|>2\tau\}\cap U^\alpha_0}\pi_{TM_j^\alpha}(Dh_j^{\alpha,\kappa}(x)+Du_j^{\alpha,\kappa}(x))\cdot D\zeta(x) J^\alpha_j(x)\, \ext x.$$
	Here, we have suppressed the $x$-dependence in the notation $\pi_{TM^\alpha_j}$, which denotes the orthogonal projection onto the tangent space of $M^\alpha_j$ at the point $(f^\alpha_j(x),x)$. Recalling that $h_j$ is linear and depends only on $x^n$, this therefore becomes
	\begin{align}
	\nonumber\sum_{\alpha=1}^4D_nh^{\alpha,\kappa}_j&\int_{B^n_{3/8}(0)\cap \{|x^n|>2\tau\}\cap U^\alpha_0}\pi_{TM^\beta_j}(e_n)\cdot D\zeta J^\alpha_j\, \ext x\\
	& + \sum_{\alpha=1}^4 \int_{B^n_{3/8}(0)\cap \{|x^n|>2\tau\}\cap U^\alpha_0}\pi_{TM^\alpha_j}(Du_j^{\alpha,\kappa}(x))\cdot D\zeta(x) J^\alpha_j(x)\, \ext x.\label{E:fine-b-3}
	\end{align}
	We start by showing that the first term in \eqref{E:fine-b-3} is $o(\tilde{E}_j)$. For this, let us first show that we can ignore the Jacobian factor: indeed, using \eqref{E:J-bounds}, elliptic estimates (to control the energy of $u^\alpha_j$ by its $L^2$ norm on a larger set), and Theorem \ref{thm:graphical-rep}, we get
	\begin{align*}
		&\left|D_nh^{\alpha,\kappa}_j\int_{B^n_{3/8}(0)\cap \{|x^n|>2\tau\}\cap U^\alpha_0}\pi_{TM^\alpha_j}(e_n)\cdot D\zeta(x) (J^\alpha_j(x)-1)\, \ext x\right|\\
		& \hspace{22em} \leq C\sup|D\zeta|\int_{B^n_{1/2}(0)\cap \{|x^n|>\tau\}}|u^\alpha_j|^2\\
		& \hspace{22em} \leq C\sup|D\zeta|\tilde{E}_j^2.
	\end{align*}
	Hence, the first term in \eqref{E:fine-b-3} can be written as
	$$\sum_{\alpha=1}^4D_{e_n}h^{\alpha,\kappa}_j\int_{B^n_{3/8}(0)\cap \{|x^n|>2\tau\}\cap U^\alpha_0}\pi_{TM^\alpha_j}(e_n)\cdot D\zeta(x)\, \ext x + O(\tilde{E}_j^2).$$
	We now need a more detailed expression for $\pi_{TM^\alpha_j}(e_n)$. Write
	\begin{equation}\label{E:fine-b-4}
	\pi_{TM^\alpha_j}(e_n) = (\pi_{TM^\alpha_j} - \pi_{H^\alpha_j})(e_n) + \pi_{H^\alpha_j}(e_n).
	\end{equation}
	Now, since $e_1,\dotsc,e_{n-1}$, $(1+|D_{e_n}h^\alpha_j|^2)^{-1/2}(e_n+D_{e_n}h^\alpha_j)$ is an orthonormal basis for the plane containing $H^\alpha_j$, we can easily compute that
	$$\pi_{H^\alpha_j}(e_n)\cdot D\zeta = \frac{D_{e_n}\zeta}{(1+|D_{e_n}h^\alpha_j|^2)^{1/2}}.$$
	Thus, the contribution which comes from this term in \eqref{E:fine-b-4} is
	$$\sum_{\alpha=1}^4\frac{D_{e_n}h^{\alpha,\kappa}_j}{(1+|D_{e_n}h^{\alpha}_j|^2)^{1/2}}\int_{B^n_{3/8}(0)\cap \{|x^n|>2\tau\}\cap U^\alpha_0}D_{e_n}\zeta$$
	which equals zero since
	$$\int_{B^n_{3/8}(0)\cap \{|x^n|>2\tau\}\cap U^\alpha_0}D_{e_n}\zeta(x)\, \ext x = \int_{B^n_{3/8}(0)\cap \{|x^n|>0\}\cap U^\alpha_0}D_{e_n}\zeta(x)\, \ext x = \int_{B^n_{3/8}(0)\cap S}\zeta$$
	which from our initial assumptions on $\zeta$ is zero. The other contribution to the integral comes from the first term in \eqref{E:fine-b-4}, namely summing over $\alpha=1,\dotsc,4$ the terms
	$$D_{e_n}h^{\alpha,\kappa}_j\int_{B^n_{3/8}(0)\cap \{|x^n|>2\tau\}\cap U^\alpha_0}(\pi_{TM^\alpha_j}-\pi_{H^\alpha_j})(e_n)\cdot D\zeta(x)\, \ext x.$$
	For this, if we use Cauchy--Schwarz and then elliptic estimates for $u^\alpha_j$ (recall that $h^\alpha_j + u^\alpha_j$ is the object satisfying the minimal surface system, and so using elliptic estimates from the corresponding elliptic equation satisfied by $u^\alpha_j$) we can control this term by
	$$C\sup|D\zeta|\cdot |D_{e_n}h_j^\alpha|\cdot \tilde{E}_j$$
	which is $o(\tilde{E}_j)$ since $|D_{e_n}h^\alpha_j|\to 0$ as $j\to\infty$. Thus, returning to \eqref{E:fine-b-3} and \eqref{E:fine-b-2}, we now have
	\begin{equation}\label{E:fine-b-5}
		\sum_{\alpha=1}^4 \int_{B^n_{3/8}(0)\cap \{|x^n|>2\tau\}\cap U^\alpha_0}\pi_{TM^\alpha_j}(Du_j^{\alpha,\kappa}(x))\cdot D\zeta(x) J^\alpha_j(x)\, \ext x = o(\tilde{E}_j) + \sqrt{\tau}O(\tilde{E}_j).
	\end{equation}
	Dividing this by $\tilde{E}_j$ and taking $j\to\infty$, we can then conclude that (as $\pi_{TM^\alpha_j}\to \pi_{P_0}$)
	$$\sum_{\alpha=1}^4\int_{B^n_{3/8}(0)\cap \{|x^n|>2\tau\}\cap U^\alpha_0}D\phi^{\alpha,\kappa}\cdot D\zeta = O(\sqrt{\tau}).$$
	Now, in light of the fact that $D_{e_n}\zeta\equiv 0$ in $\{|x^n|<2\tau\}$, we can integrate by parts in the remaining integrals and then take $\tau\downarrow 0$ to arrive at
	$$\sum^4_{\alpha=1}\int_{B^n_{3/8}(0)\cap U^\alpha_0}\phi^{\alpha,\kappa}\Delta\zeta = 0$$
	for each $\kappa=1,\dotsc,k$, which is exactly \eqref{E:fine-b-a}. This establishes the first claim of the lemma.
	
	Now we turn our attention to \eqref{E:fine-b-b}. Again, letting $\zeta$ be as in the statement of the lemma and $\tau>0$ such that $D_{e_n}\zeta\equiv 0$ in $\{|x^n|<2\tau\}$, the first variation formula \eqref{E:stationarity-2} similarly gives
	\begin{equation}\label{E:fine-b-6}
		0 = \int_{C_1(0)}\nabla^{V_j}x^n\cdot \nabla^{V_j}\widetilde{\zeta}(x)\, \ext\|V_j\|(x)
	\end{equation}
	(i.e.~our variation is now within the plane $P_0$ rather than orthogonal to it). We will now perform a very similar analysis and blow-up procedure to that above leading to \eqref{E:fine-b-a}, except now dividing by $\tilde{E}_j\hat{E}_j$. Let
	$$K_\tau := \{|x^n|\leq 2\tau\}\cap \{x:\dist(\pi_{P_0}(x),\mathcal{D})\leq 2\tau\}.$$
	Again writing $M^\alpha_j = \graph(f^\alpha_j)$ where $f^\alpha_j = h^\alpha_j + u^\alpha_j$, and so on $M^\alpha_j$ we have
	$$\nabla^{V_j}x^n\cdot\nabla^{V_j}\widetilde{\zeta} = \pi_{TM^\alpha_j}(e_n)\cdot D\widetilde{\zeta}.$$
	We then have two cases, depending on whether $\mathcal{D}(\phi)\cap B^n_1(0) = S_0\cap B^n_1(0)$ or not. If $\mathcal{D}(\phi)\cap B^n_1(0)\neq S_0\cap B^n_1(0)$, then \eqref{E:fine-b-6} gives
	\begin{equation}\label{E:fine-b-7}
		\sum_{\alpha=1,2}\int_{C_{3/8}(0)\cap K_\tau^c}\pi_{TM^\alpha_j}(e_n)\cdot D\widetilde{\zeta}\, \ext\|V_j\| = -\int_{C_{3/8}(0)\cap K_{\tau}}\nabla^{V_j}x^n\cdot\nabla^{V_j}\widetilde{\zeta}\, \ext\|V_j\|,
	\end{equation}
	whilst if $\mathcal{D}(\phi)\cap B^n_1(0) = S_0\cap B^n_1(0)$ then \eqref{E:fine-b-6} gives
	\begin{align}
		\nonumber\sum_{\alpha=1,2}\int_{C_{3/8}(0)\cap \{x^n<-2\tau\}}\pi_{TM^\alpha_j}(e_n)\cdot D\widetilde{\zeta}\, \ext\|V_j\|+ \sum_{\alpha=3,4}&\int_{C_{3/8}(0)\cap \{x^n>2\tau\}}\pi_{TM^\alpha_j}(e_n)\cdot D\widetilde{\zeta}\, \ext\|V_j\|\\
		&  = -\int_{C_{3/8}(0)\cap \{|x^n|\leq 2\tau\}}\nabla^{V_j}x^n\cdot\nabla^{V_j}\widetilde{\zeta}\, \ext\|V_j\|.\label{E:fine-b-8}
	\end{align}
	In the case \eqref{E:fine-b-7}, one should recall \eqref{E:extra} in the definition of fine blow-ups, which means that when $\mathcal{D}(\phi)\cap S_0 \neq B_1^n(0)\cap S_0$, we are assuming that $\BC_j$ is a union of two planes, and so we can view the fine blow-up as defined over two planes instead of over 4 half-planes.
	
	From now on we will work in the case of \eqref{E:fine-b-7}: this is strictly more complicated than the case in \eqref{E:fine-b-8}, and as such once we have proved \eqref{E:fine-b-b} assuming \eqref{E:fine-b-7}, an analogous argument will prove \eqref{E:fine-b-b} assuming \eqref{E:fine-b-8} instead.
	
	Let us look at the right-hand side of \eqref{E:fine-b-7}. Since $D_{e_n}\zeta\equiv 0$ in the region $\{|x^n|<\tau\}$ and about $\spt\|V_j\|$ we have that $\widetilde{\zeta}$ is constant in the vertical directions, we have for $x\in K_{\tau}$,
	$$|\nabla^{V_j}x^n\cdot\nabla^{V_j}\widetilde{\zeta}| \leq \sup|D\zeta|\sum^{n-1}_{i=1}\left|e_n^{\perp_{T_xV_j}}\cdot e_i^{\perp_{T_xV_j}}\right|.$$
	Now since $\pi_{P_0^\perp}(e_n) = 0$, we therefore have
	$$\left|e_n^{\perp_{T_xV_j}}\right|^2 = |(\pi_{T_x^\perp V_j}-\pi_{P_0^\perp})(e_n)|^2 \leq \|\pi_{T_xV_j}-\pi_{P_0}\|^2$$
	we can estimate the right-hand side of \eqref{E:fine-b-7} using Cauchy--Schwarz as follows:
	\begin{align*}
		&\left|\int_{C_{3/8}(0)\cap K_{\tau}}\nabla^{V_j}x^n\cdot\nabla^{V_j}\widetilde{\zeta}\, \ext\|V_j\|\right|\\
		& \hspace{2em} \leq n\sup|D\zeta|\left(\int_{C_{3/8}(0)\cap K_\tau}\|\pi_{T_xV_j}-\pi_{P_0}\|^2\, \ext\|V_j\|\right)^{1/2}\left(\int_{B_{1/2}(0)}\sum^{n-1}_{i=1}|e_i^{\perp_{T_xV_j}}|^2\, \ext\|V_j\|\right)^{1/2}.
	\end{align*}
	Now using \eqref{E:cor-4} to control the second integral by $C\tilde{E}_j^2$ and the non-concentration of tilt-excess with respect to coarse excess to control the first integral by $\tau^{1/2}\hat{E}_j^2$ (namely \eqref{E:Lipschitz-approx-6} from Corollary \ref{cor:Lipschitz-approx}), we get that this is at most
	$$C\sup|D\zeta|\cdot \tau^{1/4}\hat{E}_j\cdot\tilde{E}_j$$
	where $C = C(n,k,M,\beta)$. Hence, \eqref{E:fine-b-7} becomes
	\begin{equation}\label{E:fine-b-9}
		\sum_{\alpha=1,2}\int_{C_{3/8}(0)\cap K_\tau^c}\pi_{TM^\alpha_j}(e_n)\cdot D\widetilde{\zeta}\, \ext\|V_j\| = \tau^{1/4}O(\hat{E}_j\tilde{E}_j).
	\end{equation}
	Now we focus on the remaining integral in \eqref{E:fine-b-9}. We start by computing $\pi_{TM^\alpha_j}$ explicitly. Let $N$ denote the $(n+k)\times n$ matrix whose columns are $\{e_i+D_if_j^\alpha(x)\}_{i=1}^n$. Since the columns of $N$ are a basis for $TM^\alpha_j$, it follows that $N(N^TN)^{-1}N^T$ is the $(n+k)\times (n+k)$ matrix representing $\pi_{TM^\alpha_j}$. As $D\widetilde{\zeta}$ is a vector parallel to $P_0$, we have from direct computation that $N^T D\widetilde{\zeta} = D\widetilde{\zeta}$, and similarly $N^Te_n = e_n$, and so if we write $B^\alpha_j$ for the $n\times n$ matrix $N^TN$, so that $(B^\alpha_j(x))_{pi} = \delta_{pi} + D_pf^\alpha_j(x) \cdot D_if^\alpha_j(x)$ for $p,i\in \{1,\dotsc,n\}$, and $A^\alpha_j := (B^\alpha_j)^{-1}$, then
	\begin{align*}
	\pi_{TM^\alpha_j}(e_n)\cdot D\widetilde{\zeta} = N(N^TN)^{-1}N^Te_n \cdot D\widetilde{\zeta} = (AN^Te_n)\cdot (N^TD\widetilde{\zeta}) & = A(e_n)\cdot D\widetilde{\zeta}\\
	& = \sum^n_{i=1}(A^\alpha_j(x))_{in}D_i\zeta(x).
	\end{align*}
	Then, expression the integral in \eqref{E:fine-b-9} as an integral over the domain in $P_0$, this integral equals
	\begin{equation}\label{E:fine-b-10}
		\sum_{\alpha=1,2}\int_{B^n_{3/8}(0)\cap \{r_{\mathcal{D}}>2\tau\}}\sum^n_{i=1}(A^\alpha_j(x))_{in}D_i\zeta(x)J^\alpha_j(x)\, \ext x.
	\end{equation}
	We will need to compute the integrand here rather carefully. To make things easier in the computations which follow, we will temporarily drop the indices $j$ and $\alpha$ and remove the explicit $x$-dependence from the notation; thus, we will write $f = f^\alpha_j(x)$, $A = A^\alpha_j(x)$, $J = J^\alpha_j(x)$, etc.
	
	Looking at the $n^{\text{th}}$ column of the matrix equation $BA = I_n$ tells us
	\begin{equation}\label{E:fine-b-11}
		\sum^n_{i=1}(\delta_{pi}+D_pf\cdot D_if)A_{in} = 0 \qquad \text{for }p=1,\dotsc,n-1;
	\end{equation}
	and
	\begin{equation}\label{E:fine-b-12}
		\sum^n_{i=1}(\delta_{ni} + D_n f\cdot D_i f) A_{in} = 1.
	\end{equation}
	Using the fact that $D_ph=0$ for $p\neq n$, \eqref{E:fine-b-11} gives that for $p\neq n$,
	\begin{align*}
		A_{pn} & = -\sum^n_{i=1}D_p f\cdot D_i f A_{in}\\
		& = - \sum_{i=1}^n\left[D_pu\cdot (D_i u + D_i h)\right]A_{in}\\
		& = -D_pu\cdot D_n h A_{nn} - \sum^n_{i=1}D_p u\cdot D_i u A_{in}.
	\end{align*}
	The second term on the last line above is quadratic in $Du$. Thus, by elliptic estimates and Theorem \ref{thm:graphical-rep}, its contribution in \eqref{E:fine-b-10} will be of order $\tilde{E}_j^2$, which is $o(\tilde{E}_j\hat{E}_j)$. We will record this by writing
	\begin{equation}\label{E:fine-b-13}
		A_{pn} = -D_p u \cdot D_n h A_{nn} + \boldsymbol{\eps}
	\end{equation}
	where we will use $\boldsymbol{\eps}$ to denote error terms which may different from line-to-line and whose contribution to \eqref{E:fine-b-10} after integration is $o(\tilde{E}_j\hat{E}_j)$.
	
	Now from \eqref{E:fine-b-12} (again using that $D_ph = 0$ for $p\neq n$) we have
	\begin{align*}
		A_{nn} & = 1-\sum^n_{i=1}(D_nu+D_nh)\cdot (D_iu+D_ih)A_{in}\\
		& = 1-\sum^n_{i=1}D_nu\cdot D_i u A_{in} - \sum^n_{i=1}D_nh\cdot D_iuA_{in} - D_nu\cdot D_nh A_{nn} - |D_nh|^2A_{nn}.
	\end{align*}
	The first term is quadratic in $Du$, and so once again will be $o(\tilde{E}_j\hat{E}_j)$ and so will be an $\boldsymbol{\eps}$ term. For the second, using that fact that $A\to I_n$ converges to the identity matrix as $j\to\infty$, and since $D_n h\cdot D_i u$ is $O(\tilde{E}_j\hat{E}_j)$, for $i\neq n$ we have that $A_{in}\to 0$ and so the terms $D_nh\cdot D_i u A_{in}$ with $i\neq n$ are therefore in fact $o(\tilde{E}_j\hat{E}_j)$. Hence, we can rewrite the above as
	$$A_{nn} = 1-2D_n u \cdot D_n h A_{nn} - |D_nh|^2A_{nn} + \boldsymbol{\eps}$$
	which after rearranging becomes
	$$A_{nn} = \left(1+2D_nu\cdot D_nh + |D_nh|^2\right)^{-1} + \boldsymbol{\eps}.$$
	Applying the Taylor expansion for $(1+t)^{-1}$ we therefore have that
	\begin{equation}\label{E:fine-b-14}
		A_{nn} = 1-2D_nu\cdot D_nh - |D_nh|^2(1+\mathbf{h}) + \boldsymbol{\eps}
	\end{equation}
	where here the term $\mathbf{h}$ denotes a constant only depending on $|D_nh|^2$ which $\to 0$ as $j\to\infty$.
	
	In \eqref{E:fine-b-13} and \eqref{E:fine-b-14} we have the desired expressions for the terms $A_{in}$ appearing in \eqref{E:fine-b-10}. We will also need an expression for the Jacobian factor $J$. Start by writing $B = P + Q$, where
	$$P_{pi} := \delta_{pi} + D_p h\cdot D_ih + D_ph\cdot D_i u + D_pu\cdot D_i h \qquad \text{and} \qquad Q_{pi}:= D_pu\cdot D_iu.$$
	Using the mean-value inequality for the matrix function $M\mapsto \sqrt{\det(M)}$, as $J = \sqrt{\det(B)}$, we get
	$$J = \sqrt{\det(P)} + \boldsymbol{\eps}.$$
	Now write $P = R+S$, where $R_{pi} := \delta_{pi} + D_p h \cdot D_i h$ and $S_{pi}:= D_ph\cdot D_i u + D_pu\cdot D_i h$. We may assume that $j$ is large enough so that $R$ is invertible. This means that by expanding the determinant about $R$, we have
	$$\det(P) = \det(R) + \det(R)\text{trace}(R^{-1}S) + O(\|S\|^2).$$
	By direct computation we have (again, using that $D_ph = 0$ for $p\neq n$)
	$$\det(R) = 1+|D_nh|^2 \qquad \text{and} \qquad \text{trace}(R^{-1}S) = \frac{2D_nh\cdot D_nu}{1+|D_nh|^2}.$$
	Using these together with the Taylor expansion for $t\mapsto (1+t)^{1/2}$ as well as the appropriate elliptic estimates for $Du$ and $Dh$, we get:
	\begin{align*}
		\sqrt{\det(P)} & = \left(1+|D_nh|^2 + 2D_nh\cdot D_nu + O(\|S\|^2)\right)^{1/2}\\
		& = \sqrt{1+|D_nh|^2}\left(1+\frac{2D_nh\cdot D_nu}{1+|D_nh|^2} + O(\|S^2\|\right)^{1/2}\\
		& = \sqrt{1+|D_nh|^2} + \frac{D_nh D_nu}{\sqrt{1+|D_nh|^2}} + \boldsymbol{\eps},
	\end{align*}
	and hence we therefore have
	\begin{equation}\label{E:fine-b-15}
		J = \sqrt{1+|D_nh|^2} + \frac{D_nh \cdot D_nu}{\sqrt{1+|D_nh|^2}} + \boldsymbol{\eps}.
	\end{equation}
	Thus we can compute, using \eqref{E:fine-b-14}, \eqref{E:fine-b-15}, and elliptic estimates (giving that any term involving $3$ or more factors of $u$ or $h$ combined is $\boldsymbol{\eps}$), that
	\begin{equation}\label{E:fine-b-16}
		A_{nn}J = -D_nu\cdot D_nh\left(\frac{1+2|D_nh|^2}{\sqrt{1+|D_nh|^2}}\right) + \sqrt{1+|D_nh|^2}(1-|D_nh|^2(1+\mathbf{h})) + \boldsymbol{\eps}.
	\end{equation}
	Then, writing the integrand in \eqref{E:fine-b-10} as
	$$\sum^{n-1}_{i=1}A_{in}D_i\zeta J + A_{nn}D_n\zeta J,$$
	and using \eqref{E:fine-b-13} and \eqref{E:fine-b-16}, the integral in \eqref{E:fine-b-10} becomes the sum over $\alpha=1,2$ of
	\begin{align*}
		-\int_{B^n_{3/8}(0)\cap\{r_{\mathcal{D}}>2\tau\}}& \sum^{n-1}_{i=1}D_pu\cdot D_nh A_{nn}D_i\zeta J\\
		&\hspace{-2em} - \int_{B^n_{3/8}(0)\cap\{r_{\mathcal{D}}>2\tau\}}D_nu\cdot D_nh\left(\frac{1+2|D_nh|^2}{\sqrt{1+|D_nh|^2}}\right)D_n\zeta\\
		& \hspace{3em} + \sqrt{1+|D_nh|^2}(1-|D_nh|^2(1+\mathbf{h}))\int_{B^n_{3/8}(0)\cap \{r_{\mathcal{D}}>2\tau\}}D_n\zeta\\
		& \hspace{8em} + o(\tilde{E}_j\hat{E}_j).
	\end{align*}
	As before in the proof of \eqref{E:fine-b-a}, by our choice of $\zeta$ the integral in the third term here is zero. We then divide \eqref{E:fine-b-9} by $\tilde{E}_j\hat{E}_j$ and let $j\to\infty$: using the definition of $m^\alpha$, the convergence properties to the fine blow-up in the region $\{r_{\mathcal{D}}>2\tau\}$, as well as the fact that $J^\alpha_j\to 1$, $A^\alpha_j\to I_n$ almost everywhere, and $D_nh\to 0$, combining \eqref{E:fine-b-9} with \eqref{E:fine-b-10} and the above computation gives
	\begin{align*}
		0 & = \sum_{\alpha=1,2}\int_{B^n_{3/8}(0)\cap \{r_{\mathcal{D}}>2\tau\}\cap \{x^n<0\}}\left(\sum^k_{\kappa=1}m^\alpha_{\kappa,n}D\phi^{\alpha,\kappa}(x)\right)\cdot D\zeta(x)\, \ext x\\
		& \hspace{2em} + \sum_{\alpha=3,4}\int_{B^n_{3/8}(0)\cap \{r_{\mathcal{D}}>2\tau\}\cap \{x^n>0\}}\left(\sum^k_{\kappa=1}m^\alpha_{\kappa,n}D\phi^{\alpha,\kappa}(x)\right)\cdot D\zeta(x)\, \ext x + C\tau^{1/4}.
	\end{align*}
	Now giving our assumptions on $\zeta$, we can integrate by parts in these integrals, and then let $\tau\downarrow 0$ and we get exactly \eqref{E:fine-b-b}. (Notice that here, when we integrate by parts in the directions $e_1,\dotsc,e_{n-1}$, there could be extra non-zero boundary terms introduced over $\{r_{\mathcal{D}}=2\tau\}\setminus\{|x^n|=2\tau\}$. However, these terms go to zero as $\tau\to 0$ (due to cancellation with each other) and so do not impact the limit. Alternatively, one could add back in the part of the integral over $B^n_{3/8}(0)\cap \{r_{\mathcal{D}}\leq 2\tau\}$ before integrating by parts, and note that this extra piece added $\to 0$ as $\tau\to 0$.) This completes the proof of the lemma.
\end{proof}

\subsection{Boundary Behaviour of Homogeneous Degree One Fine Blow-Ups}

Armed with Lemma \ref{lemma:fine-b}, we can now work towards establishing that the boundary values of fine blow-ups are smooth. By this, we primarily mean the regularity of the function $\lambda = \lambda_\phi:\mathcal{D}(\phi)\cap B^n_{1/8}(0)\to S^\perp$ constructed in Lemma \ref{lemma:fine-continuity-1} which gives the values of the fine blow-up $\phi$ along $\mathcal{D}(\phi)$, as away from $\mathcal{D}(\phi)$ we know the fine blow-up is smooth. This will be a crucial step to proving the full regularity result for fine blow-ups.

\begin{lemma}\label{lemma:fine-lambda}
	Let $\phi = \{\phi^\alpha\}_{\alpha\in I}\in\mathfrak{b}$ be a fine blow-up. Then:
	\begin{enumerate}
		\item [(1)] If $\phi \in \mathfrak{b}_{n-1}$, then:
		\begin{enumerate}
			\item [(a)] The function $\lambda = \lambda_\phi:\mathcal{D}(\phi)\cap B^n_{1/8}(0)\to S_0^\perp$ is smooth; furthermore, when $\mathcal{D}(\phi)\cap B_{1/8}^n(0)\subsetneq S_0\cap B^n_{1/8}(0)$, we can say that $\lambda$ is the restriction of a smooth function defined on $S_0$.
			\item [(b)] If $\phi$ is homogeneous of degree one, then $\lambda$ is the restriction to $\mathcal{D}(\phi)$ of a linear function.
		\end{enumerate}
		\item [(2)] If $\phi\in\mathfrak{b}_{\sfrak}$ with $\sfrak<n-1$ (in particular, $I = \{1,2\}$), then:
		\begin{enumerate}
			\item [(a)] For $\alpha=1,2$, $\phi^\alpha$ extends to a smooth (harmonic) function on $B^n_{1/8}(0)$. In particular, the function $\lambda = \lambda_\phi:\mathcal{D}(\phi)\cap B^n_{1/8}(0)\to S^\perp$ is the restriction of a smooth function defined on $B^n_{1/8}(0)$.
			\item [(b)] If $\phi$ is homogeneous of degree one, then $\phi^\alpha$ is linear for $\alpha=1,2$, and $\lambda$ is linear.
		\end{enumerate}
	\end{enumerate}
\end{lemma}

\textbf{Note:} If $\phi\in \mathfrak{b}_{n-1}$ is homogeneous of degree one, Lemma \ref{lemma:fine-lambda} says nothing about whether $\phi$ is linear or not, only that the function $\lambda$ is. Showing in this case that $\phi$ is in fact linear is the main step to proving regularity of fine blow-ups in $\mathfrak{b}_{n-1}$, and will only be achievable after these preliminary results (see Theorem \ref{thm:fine-classification} later).

\begin{proof}
	We will prove (1) first. Take $\zeta\in C^\infty_c(B^n_{3/8}(0))$ as in the statement of Lemma \ref{lemma:fine-b}. Thus, by \eqref{E:fine-b-a} we have
	$$\sum^4_{\alpha=1}\int_{B^n_{3/8}(0)\cap U^\alpha_0}\phi^\alpha\Delta\zeta = 0.$$
	If we further insist that $\zeta$ is an even function in the $x^n$ variable, then we can write this as
	\begin{equation}\label{E:fine-lambda-1}
	\int_{B^n_{3/8}(0)\cap \{x^n>0\}}\overline{\phi}(x)\Delta\zeta = 0
	\end{equation}
	where $\overline{\phi}(x):= \sum_{\alpha=1,2}\phi^\alpha(x-2x^n e_n) + \sum_{\alpha=3,4}\phi^\alpha(x)$. For $i\in \{1,\dotsc,n-1\}$ and $h\in (-1/64,1/64)$, define $\delta_{i,h}\zeta(x):= \zeta(x+he_i) - \zeta(x)$. Notice now that since $\spt(\zeta)$ is a compact subset of $B^n_{3/8}(0)$ and $\int_{S_0\cap B^n_{3/8}(0)}\zeta\, \ext\H^{n-1}=0$, for all $|h|$ sufficiently small we have $\int_{S_0\cap B^n_{3/8}(0)}\delta_{i,h}\zeta\, \ext\H^{n-1} = 0$. Thus, $\delta_{i,h}\zeta$ is even in $x^n$ and satisfies the hypotheses of Lemma \ref{lemma:fine-b}, and thus \eqref{E:fine-lambda-1} gives
	\begin{equation}\label{E:fine-lambda-2}
		\int_{B^n_{3/8}(0)\cap \{x^n>0\}}\overline{\phi}\Delta(\delta_{i,h}\zeta)\, \ext x = 0
	\end{equation}
	for any $i\in \{1,\dotsc,n-1\}$ and $|h|$ sufficiently small. Rearranging \eqref{E:fine-lambda-2} we get
	\begin{equation}\label{E:fine-lambda-3}
		\int_{B^n_{3/8}(0)\cap \{x^n>0\}}\delta_{i,-h}\overline{\phi}\Delta\zeta = 0.
	\end{equation}
	Next, by approximating any smooth function that is even in the $x^n$-variable by smooth functions for which $D_{e_n}\zeta\equiv 0$ in a neighbourhood of $S_0$, and by using the dominated convergence theorem, we can deduce that \eqref{E:fine-lambda-3} holds for any $\zeta\in C^\infty_c(B^n_{3/8}(0))$ which is even in $x^n$; indeed, we can also remove the condition $\int_{S_0\cap B^n_{3/8}(0)}\zeta\, \ext\H^{n-1} = 0$ by subtracting a constant from $\zeta$, which \eqref{E:fine-lambda-3} is insensitive to. Now let $\hat{\phi}$ denote the extension to $\overline{\phi}|_{B^n_{3/8}(0)\cap\{x^n>0\}}$ to $B^n_{3/8}(0)$ that is even in the $x^n$-variable. We therefore have
	\begin{equation}\label{E:fine-lambda-4}
		\int_{B^n_{3/8}(0)}\delta_{i,-h}\hat{\phi}\Delta\zeta = 0
	\end{equation}
	for any $\zeta\in C^\infty_c(B^n_{3/8}(0))$ which is even in the $x^n$-variable. However, since $\hat{\phi}$ is even in $x^n$, \eqref{E:fine-lambda-4} also trivially holds for any $\zeta$ that is odd in the $x^n$-variable. As \eqref{E:fine-lambda-4} is linear in $\zeta$, and any function is a linear combination of an odd and an even function, we deduce that \eqref{E:fine-lambda-4} holds for any $\zeta\in C^\infty_c(B^n_{3/8}(0))$. Thus, $\delta_{i,-h}\hat{\phi}$ is weakly harmonic in $B^n_{3/8}(0)$, and so Weyl's lemma implies that $\delta_{i,-h}\hat{\phi}$ is smooth and harmonic on $B^n_{3/8}(0)$.
	
	Now, away from $S_0$, it is clear that $h^{-1}\delta_{i,-h}\hat{\phi}$ converges pointwise to $D_i\hat{\phi}$ as $h\to 0$; this is simply because $\phi^\alpha$ is harmonic (and thus smooth) here for each $\alpha$. In light of \eqref{E:fine-blow-up-L1} and standard facts about difference quotients (namely that the $L^1$ norm of $h^{-1}\delta_{i,-h}\hat{\phi}$ over a compact subset of $B^n_{3/8}(0)$ is controlled by $\|D_i\hat{\phi}\|_{L^1(B^n_{3/8}(0))}$ for all $h$ sufficiently small; see \cite[Lemma 7.23]{gt}) and the dominated convergence theorem, it follows that also $h^{-1}\delta_{i,-h}\hat{\phi}\to D_i\hat{\phi}$ in $L^1(B^n_{1/8}(0))$ also. As the $L^1$ limit of harmonic functions must be harmonic, we deduce that $D_i\hat{\phi}$ is smooth and harmonic in $B^n_{1/8}(0)$ for each $i\in \{1,\dotsc,n-1\}$. The fundamental theorem of calculus then implies that $\hat{\phi}$ is smooth and harmonic in $B^n_{1/8}(0)$.
	
	Now, repeat this argument, but starting with \eqref{E:fine-b-b} of Lemma \ref{lemma:fine-b}, i.e.~if we set $\psi^\alpha := \sum^k_{\kappa=1}m^\alpha_{\kappa,n}\phi^{\alpha,\kappa}$, then we have the identity
	$$\sum^4_{\alpha=1}\int_{B^n_{3/8}(0)\cap U^\alpha_0}\psi^\alpha\Delta\zeta = 0.$$
	Then arguing exactly as above, we prove that the function $\hat{\psi}$ (analogously defined, i.e.~the extension of $x\mapsto \sum_{\alpha=1,2}\psi(x-2x^ne_n) + \sum_{\alpha=3,4}\psi(x)$ to the whole of $B^n_{3/8}(0)$ which is even in the $x^n$-variable) is smooth and harmonic on $B^n_{1/8}(0)$.
	
	Now for $z\in\mathcal{D}(\phi)$, write $\lambda = \lambda(z) = (\lambda^{\top_{P_0}},\lambda^{\perp_{P_0}}) = (\lambda^n,\lambda^{n+1},\dotsc,\lambda^{n+k})\in S_0^\perp$. Then, from the definitions of $\hat{\phi}$ and $\hat{\psi}$ above, and the forms of the values of $\phi(z)$ coming from Lemma \ref{lemma:fine-continuity-1}, we have
	\begin{equation}\label{E:fine-lambda-5}
		\hat{\phi}^\kappa(z) = 4\lambda(z)^{n+\kappa} - \sum_{\alpha=1}^4m^\alpha_{\kappa,n}\lambda^{n} \qquad \text{for }\kappa=1,\dotsc,k;
	\end{equation}
	\begin{equation}\label{E:fine-lambda-6}
		\hat{\psi}(z) = \sum_{\alpha=1}^4\sum^k_{\kappa=1}m^\alpha_{\kappa,n}(\lambda^{n+\kappa}-m^\alpha_{\kappa,n}\lambda^n).
	\end{equation}
	This gives $k+1$ equations for $k+1$ variables $\lambda^n,\dotsc,\lambda^{n+k}$. If we rearrange \eqref{E:fine-lambda-5} to solve for $\lambda^{n+\kappa}$ in terms of $\lambda^n$ and $\hat{\phi}^\kappa(z)$, and substitute this into \eqref{E:fine-lambda-6}, upon rearranging we get
	\begin{align*}
		\hat{\psi}(z) - \frac{1}{4}\sum^k_{\kappa=1}\sum^4_{\alpha=1}m^\alpha_{\kappa,n}\hat{\phi}^\kappa(z) & = \sum^4_{\alpha=1}\sum^k_{\kappa=1}m^\alpha_{\kappa,n}\left(\frac{1}{4}\sum^4_{\alpha^\prime=1}m^{\alpha^\prime}_{\kappa,n}-m^\alpha_{\kappa,n}\right)\lambda^n\\
		& = \frac{1}{4}\sum^k_{\kappa=1}\sum^4_{\alpha,\alpha^\prime=1}m^\alpha_{\kappa,n}(m^{\alpha^\prime}_{\kappa,n}-m^\alpha_{\kappa,n})\lambda^n\\
		& = \left[-\frac{1}{8}\sum^k_{\kappa=1}\sum^4_{\alpha,\alpha^\prime=1}(m^\alpha_{\kappa,n}-m^{\alpha^\prime}_{\kappa,n})^2\right]\lambda^n.
	\end{align*}
	However, it follows from Lemma \ref{lemma:nu-a} (in particular, \eqref{E:nu-a-1}) that $\sum^k_{\kappa=1}\sum^4_{\alpha,\alpha^\prime}(m^\alpha_{\kappa,n}-m^{\alpha^\prime}_{\kappa,n})^2 >0$ (else the coarse blow-up of the $\BC_j$ would be a single plane, which is impossible), and moreover we also know from the construction of fine blow-ups that $|m^\alpha|\leq C$. This means that $\lambda^n(z)$ can be written as a linear combination of $\hat{\phi}^{n+1}(z),\dotsc,\hat{\phi}^{n+k}(z)$ and $\hat{\psi}(z)$ with bounded coefficients. Then we can proceed to solve for $\lambda^{n+\kappa}$ as well, ultimately writing $\lambda$ as a linear function of $\hat{\phi}$ and $\hat{\psi}$. As $\hat{\psi}$ and $\hat{\phi}$ are smooth, this shows that indeed $\lambda$ is smooth, and it is the restriction of a smooth function to $\mathcal{D}(\phi)\cap B^n_{1/8}(0)$ (indeed, this linear combination of $\hat{\phi}$ and $\hat{\psi}$ was defined on all of $B^n_{1/8}(0)$). This proves (1)(a) of the lemma. To see (1)(b), notice that if $\phi$ is homogeneous of degree one, then $\hat{\phi}$ and $\hat{\psi}$ must also be homogeneous of degree one, and hence as they are themselves harmonic they must then be linear. Thus, $\lambda$ must be linear as it is a linear combination of linear functions, proving (1)(b).
	
	Now we prove (2). Note first that the smoothness of $\phi^1,\phi^2$ claimed in (2)(a) is immediate from Theorem \ref{thm:fine-continuity-2}: indeed, when $\sfrak(\phi)<n-1$, we have $\mathcal{D}(\phi)$ is contained within a (compact subset of a) subspace of dimension at most $n-2$, and so $\H^{n-2}(\mathcal{D}(\phi))<\infty$. However, for a bounded harmonic function defined on an $n$-dimensional plane, sets of finite $\H^{n-2}$-measure are removable singularities, and hence we see that $\phi^1,\phi^2$ must be smooth harmonic functions on all of $P_0$. From this the claim concerning $\phi$ in (2)(b) immediately follows: indeed, if $\phi$ is homogeneous of degree one, then $\phi^1,\phi^2$ are harmonic functions which are homogeneous of degree one, from which it follows they must be linear.
	
	Thus, to conclude the proof of (2) we only need to show the claims related to $\lambda$ therein. We will prove these by showing that $\lambda$ is a linear combination of $\phi^1,\phi^2$; once we have shown this both claims concerning $\lambda$ in (2) follow.
	
	Firstly, write $\hat{\phi}:=\sum_{\alpha=1,2}\phi^{\alpha,\kappa}$ for $\kappa=1,\dotsc,k$ and $\hat{\psi}^i := \sum_{\alpha=1,2}\sum^k_{\kappa=1}m^\alpha_{\kappa,i}\phi^{\alpha,\kappa}$ for $i=\sfrak+1,\dotsc,n$, where $\sfrak:=\sfrak(\phi)$. Then Lemma \ref{lemma:fine-continuity-1} gives for $z\in \mathcal{D}(\phi)$ the system of equations
	\begin{equation}\label{E:fine-lambda-7}
		\hat{\phi}^\kappa(z) = 2\lambda^{n+\kappa}-\sum_{\alpha=1,2}\sum^n_{i=\sfrak+1}m^\alpha_{\kappa,i}\lambda^i \qquad \text{for }\kappa=1,\dotsc,k;
	\end{equation}
	\begin{equation}\label{E:fine-lambda-8}
		\hat{\psi}^i(z) = \sum_{\alpha=1,2}\sum^k_{\kappa=1}m^\alpha_{\kappa,i}\left(\lambda^{n+\kappa}-\sum^n_{\ell = \sfrak+1}m^\alpha_{\kappa,\ell}\lambda^\ell\right) \qquad \text{for }\ell=\sfrak+1,\dotsc,n.
	\end{equation}
	Here, we are writing $\lambda = (\lambda^{\sfrak+1},\dotsc,\lambda^n,\lambda^{n+1},\dotsc,\lambda^{n+k})\in S^\perp$, so that $\lambda^{\perp_{P_0}} = (\lambda^{n+1},\dotsc,\lambda^{n+k})$ and $\lambda^{\top_{P_0}} = (\lambda^{\sfrak+1},\dotsc,\lambda^n)$. Thus, we have $n+k-\sfrak$ equations for $n+k-\sfrak$ unknowns. We will often drop the labelling on the sums for ease of notation.
	
	Firstly, \eqref{E:fine-lambda-7} implies for each $\kappa=1,\dotsc,k$ that
	\begin{equation}\label{E:fine-lambda-9}
		\lambda^{n+\kappa} = \frac{1}{2}\hat{\phi}^\kappa + \frac{1}{2}\sum_{\alpha,i}m^\alpha_{\kappa,i}\lambda^i.
	\end{equation}
	Plugging \eqref{E:fine-lambda-9} into \eqref{E:fine-lambda-8} gives for each $i=\sfrak+1,\dotsc,n$,
	$$\hat{\psi}^i = \sum_{\alpha,\kappa}m^\alpha_{\kappa,i}\left(\frac{1}{2}\hat{\phi}^\kappa + \frac{1}{2}\sum_{\alpha^\prime,\ell}m^{\alpha^\prime}_{\kappa,\ell}\lambda^\ell - \sum_\ell m^\alpha_{\kappa,\ell}\lambda^\ell\right)$$
	i.e.
	\begin{equation}\label{E:fine-lambda-10}
		\hat{\psi}^i - \frac{1}{2}\sum_{\alpha,\kappa}m^\alpha_{\kappa,i}\hat{\phi}^\kappa = \sum_\ell\left(\sum_{\alpha,\alpha^\prime,\kappa}\frac{1}{2}m^\alpha_{\kappa,i}m^{\alpha^\prime}_{\kappa,\ell} - \sum_{\alpha,\kappa}m^\alpha_{\kappa,i}m^\alpha_{\kappa,\ell}\right)\lambda^\ell.
	\end{equation}
	This is now $n-\sfrak$ equations for $(\lambda^{\sfrak+1},\dotsc,\lambda^n)$. We therefore turn our attention to the matrix determining solvability of \eqref{E:fine-lambda-10}: this is the $(n-\sfrak)\times (n-\sfrak)$ matrix $M$ whose $(i,j)^{\text{th}}$ entry is given by
	\begin{align*}
		M_{ij} & := \sum_{\alpha,\alpha^\prime,\kappa}\frac{1}{2}m^\alpha_{\kappa,\sfrak+i}m^{\alpha^\prime}_{\kappa,\sfrak+j} - \sum_{\alpha,\kappa}m^\alpha_{\kappa,\sfrak+i}m^\alpha_{\kappa,\sfrak+j}\\
		& = \frac{1}{2}\sum^k_{\kappa=1}\left(\sum_{\alpha,\alpha^\prime}m^\alpha_{\kappa,\sfrak+i}m^{\alpha^\prime}_{\kappa,\sfrak+j} - 2\sum_{\alpha}m^\alpha_{\kappa,\sfrak+i}m^\alpha_{\kappa,\sfrak+j}\right)\\
		& = -\frac{1}{2}\sum^k_{\kappa=1}(m^1_{\kappa,\sfrak+i}-m^2_{\kappa,\sfrak+i})(m^1_{\kappa,\sfrak+j}-m^2_{\kappa,\sfrak+j}).
	\end{align*}
	Now let $N = (N_{\kappa,j})_{\kappa,j}$ be the $k\times (n-\sfrak)$ matrix with entries given by $N_{\kappa,j}:= m^1_{\kappa,\sfrak+j} - m^2_{\kappa,\sfrak+j}$. The above shows that
	$$M_{ij} = -\frac{1}{2}\sum^k_{\kappa=1}N_{\kappa,i}N_{\kappa,j} = -\frac{1}{2}(N^TN)_{ij}$$
	i.e.
	\begin{equation}\label{E:fine-lambda-11}
		M = -\frac{1}{2}N^TN.
	\end{equation}
	From the construction of the fine blow-up, it follows that we must have $\text{rank}(N) = n-\sfrak$. Indeed, if this were to fail then the blow-up of the cones $\BC_j$ would be a cone with larger spine dimension (as $S(\BC_j)\equiv S$) which is impossible by Lemma \ref{lemma:nu-a}(ii); note here also that to even be in such a fine blow-up setting, it \emph{must} be that $n-\sfrak \leq k$ (as, after we quotient out by $S(\BC_j)$, we have two ($n-\sfrak$)-dimensional planes meeting only at the origin in $\R^{n-\sfrak+k}$, which is only possible if $n-\sfrak + k \geq 2(n-\sfrak)$, i.e.~$k\geq n-\sfrak$). It therefore follows from \eqref{E:fine-lambda-11} that $\text{rank}(M) = \text{rank}(N^TN) = \text{rank}(N) = n-\sfrak$. But this means that \eqref{E:fine-lambda-10} has a unique solution, i.e.~we can write $\lambda^{\sfrak+1},\dotsc,\lambda^n$, as a linear function of $\hat{\phi}$ and $\hat{\psi}$. Returning then to \eqref{E:fine-lambda-9}, this implies that we can also express $\lambda^{n+1},\dotsc,\lambda^{n+k}$ as linear functions of $\hat{\phi}$ and $\hat{\psi}$. From the definitions of $\hat{\phi}$ and $\hat{\psi}$, this means that we have written $\lambda$ as a linear function of $\phi^1,\phi^2$. From this, the remaining conclusions of (2) concerning $\lambda$ follow immediately from the already established conclusions in (2) related to $\phi^1,\phi^2$. Thus, the proof of the lemma is complete. 
\end{proof}

Notice at this point we have fully understood the regularity of fine blow-ups in $\mathfrak{b}_\sfrak$ with $\sfrak<n-1$: there are smooth, harmonic functions. Thus, the only remaining regularity question concerns blow-ups in $\mathfrak{b}_{n-1}$. This will be our focus for the rest of the section.

\subsection{Further Properties of Fine Blow-Ups}

In order to prove the regularity of fine blow-ups in $\mathfrak{b}_{n-1}$, we will need to know that the class $\mathfrak{b}_{n-1}$ has certain flexibility under operations corresponding to natural geometric variations of the fine blow-up process, in the same spirit as those seen for coarse blow-ups in Theorem \ref{thm:blow-up-properties}. Most importantly, we wish to have the ability to subtract elements of $\mathfrak{l}$ from general elements of $\mathfrak{b}_{n-1}$ (and even $\mathfrak{b}$). This corresponds to appropriately tilting both the blow-up sequence $V_j$ as well as the cones $\BC_j$, and we are even allowed to rotate the planes (or half-planes) in $\BC_j$ independently (in the case of half-planes, the axis still needs to be maintained), forming a new sequences of cones $\Dbf_j$. We stress: \emph{this new sequence of cones could have a lower spine dimension than the original sequence} (but will never have a larger spine dimension). At the fine blow-up level, this allows us to subtract \emph{different} linear functions for each component of the fine blow-up. Notice that this flexibility was not possible with coarse blow-ups, and indeed the corresponding property for coarse blow-ups, $(\mathfrak{B}3\text{III})$, only allows one to subtract the \emph{same} linear function from each component of the coarse blow-up.

\textbf{Remark:} An additional point which arises when performing these operations is the dependence on the fine blow-up procedure on the constant $M$ appearing in (A3), as well as the parameter $\beta$ in the class $\mathcal{V}_\beta$; these will change under certain operations. However, we can always guarantee that $M$ changes by at most a fixed amount independent of the operation, i.e.~there is some the new constant $M_*$ obeys $M_* \leq CM$ for some $C = C(n,k)$. For $\beta$, as seen previously we can guarantee that the rotated sequence of varifolds all belong to $\mathcal{V}_{\beta/2}$, say. These changes do not impact our arguments, as the primary property of the fine blow-ups we need are the estimates which they satisfy.

\begin{lemma}[Integrating $\mathfrak{l}$]\label{lemma:fine-integration}
	Let $\phi = \{\phi^\alpha\}_{\alpha\in I}\in\mathfrak{b}$ be a fine blow-up of $(V_j)_j$ relative to $(\BC_j)_j$. Let $\pi = \{\pi^\alpha\}_{\alpha\in I}\in \mathfrak{l}$ be realisable as a blow-up off of $(\BC_j)_j$ by a sequence of cones in $\CC$. Then, there exists a sequence of rotations $R_j\in SO(n+k)$ and a sequence of cones $(\Dbf_j)_j\subset \CC$ such that:
	\begin{itemize}
		\item $|R_j-\id_{\R^{n+k}}|\leq C\hat{E}_{V_j}^{-1}\hat{E}_{V_j,\BC_j}$; in particular $R_j\to \id_{\R^{n+k}}$;
		\item $S(\Dbf_j)\subset S(\BC_j)$;
		\item $\nu((R_j)_\#\BC_j, \Dbf_j) \leq C\hat{E}_{V_j,\BC_j}$;
	\end{itemize}
	and such that the following holds. If $\phi^\alpha\neq \pi^\alpha$ for some $\alpha\in I$, then we can take a blow-up of $(\eta_{0,\sigma_j}\circ R_j)_\#V_j$ relative to $\Dbf_j$ in the region determined by $\BC_j$, where $\sigma_j\uparrow 1$, resulting in the blow-up $\widetilde{\phi} = \{\widetilde{\phi}^\alpha\}_{\alpha\in I}$ given by
	$$\widetilde{\phi} := \frac{\phi^\alpha - \pi^\alpha}{\|\phi-\pi\|_{L^2(B_1^n(0))}}, \qquad \alpha\in I.$$
	Here, $\|\phi-\pi\|_{L^2(B_1(0))} \equiv \left(\sum_{\alpha\in I}\|\phi^\alpha-\pi^\alpha\|_{L^2(B_1(0))}^2\right)^{1/2}$. Furthermore, in this case we have strong convergence in $L^2_{\text{loc}}(B_1(0))$ to the blow-up. Moreover, if $\phi(0) = 0$, $\phi$ additionally satisfies the following improved Hardt--Simon inequality at $0$: we have
	$$\sum_{\alpha\in I}\int_{B_{1/2}^n(0)}R^{2-n}\left|\frac{\del}{\del R}\left(\frac{\phi^\alpha}{R}\right)\right|^2\, \ext x \leq C\sum_{\alpha\in I}\int_{B_1^n(0)}|\phi^\alpha - \pi^\alpha|^2\, \ext x.$$
	Alternatively, if $\phi^\alpha = \pi^\alpha$ for each $\alpha\in I$, then for every $\theta\in (0,1)$ we have
	$$\hat{E}_{V_j,\BC_j}^{-1}\hat{F}_{(R_j)_\#V_j,\Dbf_j}(B_\theta^{n+k}(0))\to 0.$$
	In the above, $C$ is a constant depending on $n,k,M,\beta$, and $\|\pi\|_{L^2(B_1(0))}$.
\end{lemma}

\textbf{Remark:} In applications of this lemma, we will know that $\|\pi\|_{L^2(B_1(0))}\leq C(n,k,M,\beta)$.

\textbf{Note:} We are \emph{not} claiming that the above ``blow-up'' is a \emph{fine} blow-up as previously defined, as indeed the spine of the cones we are blowing-up off may have dropped, and we are not claiming that these new cones obey any form of Hypothesis $\dagger$. We can still however perform a ``blow-up'' off them in the region away from the spine of the \emph{original} sequence of cones, however, which is enough to deduce the given Hardt--Simon inequality.

\begin{proof}
	Since $\pi\in \mathfrak{l}$, $\mathbf{\pi}$ has a spine $S^\prime$ which obeys $S^\prime\subset S\equiv S(\BC_j)$. By the linearity, we know that there is a linear function $\lambda_\pi:S^\prime\to (S^\prime)^\perp$ such that
	$$\left.\pi^\alpha\right|_S = \lambda^{\perp_{P_0}}_\pi - m^\alpha\cdot\lambda_\pi^{\top_{P_0}}$$
	for some $k\times (n-\sfrak^\prime)$ matrices $m^\alpha$, where $\sfrak^\prime := \dim(S^\prime)$. Since $\pi$ is linear and valued in $P_0^\perp$, this then means that there are constants $c_{\sfrak^\prime+1}^\alpha,\dotsc,c_n^\alpha\in P_0^\perp$ for which
	$$\pi^\alpha(x) = \lambda_\pi(\pi_{S^\prime}(x))^{\top_{P_0}} - m^\alpha\cdot\lambda_{\pi}(\pi_{S^\prime}(x))^{\top_{P_0}} + \sum^n_{i=\sfrak^\prime+1}c_i^\alpha x^i.$$
	Now, for $\alpha\in I$ and $j=1,2,\dotsc$, set
	$$\widetilde{P}^\alpha_j := \graph\left(p_j^\alpha + \hat{E}_{V_j,\BC_j}\sum^n_{\sfrak^\prime+1}c_i^\alpha x^i\right) \qquad \text{if }I = \{1,2\};$$
	$$\widetilde{H}_j^\alpha := \graph\left(h^\alpha_j + \hat{E}_{V_j,\BC_j}\sum^n_{\sfrak^\prime+1}c_i^\alpha x^i\right) \qquad \text{if }I = \{1,\dotsc,4\};$$
	and
	$$\Dbf_j := \begin{cases}
		|\widetilde{P}_1| + |\widetilde{P}_j^2| \qquad \text{if }I=\{1,2\};\\
		\sum^4_{\alpha=1}|\widetilde{H}_j^\alpha| \qquad \text{if }I=\{1,\dotsc,4\}.
	\end{cases}$$
	We stress that, if $\phi\in \mathfrak{b}_{\sfrak}$, where $\sfrak\leq n-2$, then we are always in the case where $I = \{1,2\}$; when $\phi\in\mathfrak{b}_{n-1}$, we are defaulting in the above to the index set $I = \{1,\dotsc,4\}$, where if $\BC_j\in \mathcal{P}$ for all $j$ then we could $\Dbf_j$ in either $\mathcal{P}$ or $\CC_{n-1}\setminus\mathcal{P}$.
	
	Clearly from these definitions we have
	$$\nu(\BC_j,\Dbf_j) \leq C\hat{E}_{V_j,\BC_j}$$
	where $C$ depends on $n,k$, and $\|\pi\|_{L^2(B_1)}$. Also, clearly from the definition we have $S(\Dbf_j) \equiv S^\prime\subset S\equiv S(\BC_j)$ for all $j$.
	
	Now, for each $j$ define the linear function $\ell_j:S^\prime\to (S^\prime)^\perp$ by
	$$\ell_j := \hat{E}_{V_j,\BC_j}\lambda_\pi^{\perp_{P_0}} + \hat{E}_{V_j}^{-1}\hat{E}_{V_j,\BC_j}\lambda_\pi^{\top_{P_0}} \equiv \ell_j^{\perp_{P_0}} + \ell_j^{\top_{P_0}}.$$
	Now let:
	\begin{itemize}
		\item $R_j^{(\perp)}\in SO(n+k)$ denote a rotation minimising $|R-\id_{\R^{n+k}}|$ over all $R\in SO(n+k)$ for which $R(\graph(\ell_j^{\perp_{P_0}})) = S^\prime$;
		\item $R_j^{(\top)}\in SO(n+k)$ denote a rotation minimising $|R-\id_{\R^{n+k}}|$ over all $R\in SO(n+k)$ for which $R(\graph(\ell_j^{\top_{P_0}})) = S^\prime$.
	\end{itemize}
	These definitions imply that
	\begin{equation}\label{E:b-4-3}
		|R_j^{(\perp)}-\id_{\R^{n+k}}|\leq C\hat{E}_{V_j,\BC_j};
	\end{equation}
	\begin{equation}\label{E:b-4-4}
		|R_j^{(\top)}-\id_{\R^{n+k}}| \leq C\hat{E}_{V_j}^{-1}\hat{E}_{V_j,\BC_j}.
	\end{equation}
	Notice also that because $\graph(\ell_j^{\top_{P_0}})\subset P_0$, $R_j^{(\top)}$ necessarily fixes $P_0^\perp$ and maps $P_0$ to $P_0$.
	
	Now write $R_j := R_j^{(\perp)}\circ R_j^{(\top)}$. Using \eqref{E:translates-2}, Lemma \ref{lemma:nu-a}(i), Theorem \ref{thm:shift-estimate} for $V_j$ and $\BC_j$, and that $\lambda_\pi$ is uniformly bounded, it follows that
	\begin{equation}\label{E:b-4-5}
		\nu(\BC_j,(R_j)_\#\BC_j) \leq C\hat{E}_{V_j,\BC_j}
	\end{equation}
	(note that the factor of $\hat{E}_{V_j}^{-1}$ in \eqref{E:b-4-4} cancels, up to a constant, with the factor of $\nu(\BC_j)$ coming from \eqref{E:translates-2} when analysing the effect of $R_j^{(\top)}$ on $\BC_j$). Now set $W_j:= (\eta_{0,\sigma}\circ R_j)_\#V_j$ for $\sigma\in (0,1)$ close to $1$ (note that for all $j$ sufficiently large, $W_j\in \mathcal{V}_{\beta/2}$). We now wish to check Hypothesis $\dagger$ for $W_j$ and the original sequence $\BC_j$. It is easy to verify all the hypotheses (with suitable $\tilde{\eps}_j\downarrow 0$) except for (A3) and (A4), however these can be verified in an analogous manner to \cite[Lemma 13.1]{Wic14}, using \eqref{E:b-4-5}.
	
	We can therefore take a fine blow-up of $(W_j)_j$ relative to $(\BC_j)_j$, however to get the result we need to blow-up $(W_j)_j$ off of $(\Dbf_j)_j$, which could have a lower spine dimension. We will give the details of this in the case where $\Dbf_j\in\mathcal{P}$ (in particular, $\BC_j\in\mathcal{P}$ also), with the argument in the other case being analogous. For notational simplicity we will also use notation for which $\sigma=1$ in the definition of $W_j$ (which we may do after rescaling $V_j$ and so forth appropriately).
	
	Notice that as we can apply Theorem \ref{thm:graphical-rep} to $W_j$ and $\BC_j$, we know that $\mathcal{D}(W_j)$ is accumulating to $S(\BC_j)$, and thus outside of a neighbourhood of $S(\BC_j)$, we can follow the argument in the proof of Theorem \ref{thm:graphical-rep} to represent $W_j$ by functions over $\Dbf_j$ as well away from $S(\BC_j)$; call these functions $(\tilde{u_j}^\alpha)_{\alpha=1,2}$, and let $(u_j^\alpha)_{\alpha=1,2}$ be the functions representing $V_j$ over $\BC_j$, also in the sense of Theorem \ref{thm:graphical-rep}. For $x\in B^n_{1-\tau}(0)\setminus (S)_\tau$ and all sufficiently large $j$ (depending on $\tau$) we will work out an expression for $\tilde{u}^\alpha_j(x)$ in terms of $u_j^\alpha$. Let us write
	$$f_j^\alpha(x):= p_j^\alpha(x) + u^\alpha_j(x) \qquad \text{and} \qquad \tilde{f}_j^\alpha(x) := \tilde{p}_j^\alpha(x) + \tilde{u}^\alpha_j(x),$$
	where $\tilde{p}_j^\alpha(x) := p_j^\alpha(x) + \hat{E}_{V_j,\BC_j}\sum^n_{i=\sfrak^\prime+1}c_i^\alpha x^i$ is the linear function whose graph is $\widetilde{P}_j^\alpha$. We then have that $x+\tilde{f}^\alpha_j(x) \in \spt\|W_j\| \equiv \spt\|(R_j)_\#V_j\|$, which means that $R_j^{-1}(x+\tilde{f}_j^\alpha(x))\in \spt\|V_j\|$. Thus, we have
	$$R_j^{-1}\left(x+p_j^\alpha(x)+\hat{E}_{V_j,\BC_j}\sum^n_{i=\sfrak^\prime+1}c_i^\alpha x^i + \tilde{u}^\alpha_j(x)\right) = x^\prime + p_j^\alpha(x^\prime) + u^\alpha_j(x^\prime)$$
	where
	$$x^\prime = (\pi_{P_0}\circ R_j^{-1})(x+\tilde{f}_j^\alpha(x)).$$
	Rearranging the above, as $R_j$ is linear, we get
	\begin{equation}\label{E:b-4-6}
		\tilde{u}^\alpha_j(x) = [R_j(x^\prime) - x] + [R_j(p^\alpha_j(x^\prime) - p^\alpha_j(x)] + R_j(u^\alpha_j(x^\prime)) - \hat{E}_{V_j,\BC_j}\sum^n_{i=\sfrak^\prime+1}c^\alpha_i x^i.
	\end{equation}
	We then claim that in fact this can be simplified to:
	\begin{align}\label{E:b-4-7}
		\tilde{u}^\alpha_j(x) = -\hat{E}_{V_j,\BC_j}\lambda_{\pi}(\pi_{S^\prime}(x))^{\perp_{P_0}} + \hat{E}_{V_j}^{-1}&\hat{E}_{V_j,\BC_j}Dp^\alpha \cdot\lambda_\pi (\pi_{S^\prime}(x))^{\top_{P_0}}\\
		& + u^\alpha_j(x) - \hat{E}_{V_j,\BC_j}\sum^n_{i=\sfrak^\prime+1}c^\alpha_i x^i + o(\hat{E}_{V_j,\BC_j}). \nonumber
	\end{align}
		Once we have shown \eqref{E:b-4-7}, it clearly gives that
		$$\hat{E}_{V_j,\BC_j}^{-1}\tilde{u}^\alpha_j \to \phi^\alpha - \pi^\alpha$$
		as $j\to\infty$. We also claim that \eqref{E:b-4-7} and the $L^2$ convergence to the fine blow-ups gives
		\begin{equation}\label{E:b-4-7.5}
		\frac{\hat{E}_{W_j,\Dbf_j}}{\hat{E}_{V_j,\BC_j}} \to \left(\sum_{\alpha\in I}\|\phi^\alpha-\pi^\alpha\|^2_{L^2(B^n_1(0))}\right)^{1/2}.
		\end{equation}
		Indeed, this would follow provided we could show a $L^2$ non-concentration bound for the $L^2$ excess of $W_j$ relative to $\Dbf_j$ in a neighbourhood of $S(\BC_j)$. However, this follows from the fact that $\nu((R_j)_\#\BC_j,\Dbf_j) \leq C\hat{E}_{V_j,\BC_j}$, as the triangle inequality gives for any $\tau>0$, $\sigma\in (0,1)$, and all $j$ sufficiently large,
		\begin{align*}
		\int_{B_{1-\sigma}(0)\cap B_\tau(S(\BC_j))}&\dist^2(x,\Dbf_j)\, \ext\|W_j\|(x)\\
		& \leq C\int_{B_{1-\sigma/2}(0)\cap B_{2\tau}(S(\BC_j))}\dist^2(x,\BC_j)\, \ext\|V_j\|(x) + C\tau\nu((R_j)_\#\BC_j,\Dbf_j)
		\end{align*}
		and so applying \eqref{E:cor-3} of Corollary \ref{cor:L2-corollary} to $V_j$ and $\BC_j$, we get
		$$\int_{B_{1-\sigma}(0)\cap B_\tau(S(\BC_j))}\dist^2(x,\Dbf_j)\, \ext \|W_j\|(x) \leq C\tau^{1/2}\hat{E}_{V_j,\BC_j}$$
		which implies \eqref{E:b-4-7.5}. Thus, writing
		$$\frac{\tilde{u}_j^\alpha}{\hat{E}_{W_j,\Dbf_j}} = \frac{\tilde{u}^\alpha_j}{\hat{E}_{V_j,\BC_j}}\cdot\frac{\hat{E}_{V_j,\BC_j}}{\hat{E}_{W_j,\Dbf_j}}$$
		the claimed blow-up result of the lemma follows. Notice also that, given the above, the claimed Hardt--Simon inequality follows from applying Lemma \ref{lemma:L2-estimates-2} to $W_j$ and $\Dbf_j$ and passing to the limit, namely by manipulating \eqref{E:initial-1} in an analogous manner to that seen in the proof of \eqref{E:L2-3} in Theorem \ref{thm:L2-estimates}. Thus, for this it remains to prove \eqref{E:b-4-7}.
		
		To prove \eqref{E:b-4-7}, given \eqref{E:b-4-6}, we just need to focus on the three terms:
		$$\textnormal{I}:= R_j(x^\prime)-x; \qquad \text{II} := R_j(p^\alpha_j(x^\prime))-p^\alpha_j(x); \qquad \text{and} \qquad \text{III}:= R_j(u^\alpha_j(x^\prime)).$$
		To simplify the computations which follow, we will drop the index $j$ (although notation $o(\bullet)$ will refer to limits taken as $j\to\infty$). We will therefore write $\hat{E}:= \hat{E}_{V_j}$ and $\tilde{E}:= \hat{E}_{V_j,\BC_j}$.
		
		Firstly, note that in general if $y\in P_0^\perp$ satisfies $|y| = o(1)$, then since $R^{(\top)}(y) = y$ and $|R^{(\perp)}-\id_{\R^{n+k}}|\leq C\tilde{E}_j$, we have
		\begin{equation}\label{E:b-4-8}
			R(y) = R^{(\perp)}(y) = y + (R^{(\perp)}-\id_{\R^{n+k}})(y) = y + o(\tilde{E})
		\end{equation}
		where we have used that $|y| = o(1)$ in the last equality and that $R^{(\top)}(y) = y$ in the first equality. Similar considerations (or just acting on \eqref{E:b-4-8} by $R^{-1}$) show that for such $y$ we also have
		\begin{equation*}\label{E:b-4-9}
			R^{-1}(y) = y + o(\tilde{E}).
		\end{equation*}
		Therefore, since $\tilde{f}^\alpha(x)\in P_0^\perp$ and $|\tilde{f}^\alpha(x)| \leq C\hat{E}$, we have that
		\begin{equation}\label{E:b-4-10}
			x^\prime = \pi_{P_0}(R^{-1}(x+\tilde{f}^\alpha(x))) = \pi_{P_0}(R^{-1}(x)) + o(\tilde{E}).
		\end{equation}
		Now let us compute $\text{I}$. From the definition of $\text{I}$ and \eqref{E:b-4-10}, the fact that $|R|$ is bounded, and the fact that $R\circ \pi_{P_0}\circ R^{-1} = \pi_{R(P_0)}$, we have that
		\begin{equation}\label{E:b-4-11}
			\text{I} = \pi_{R(P_0)}(x) - x + o(\tilde{E}).
		\end{equation}
		Let us compute $\pi_{R(P_0)}$ more precisely. As $\lambda_{\pi}:S^\prime\to (S^\prime)^\perp$ is linear, we may extend it by $0$ to a function (still denoted $\lambda_{\pi}$) $P_0\to (S^\prime)^\perp$ given by $\lambda_{\pi}(z):= \sum^n_{i=1}z_i\lambda_i$ where $\lambda\in (S^\prime)^\perp$ and $\lambda_i = 0$ for $i>\sfrak^\prime$. Now, a basis for $R(P_0)$ is $\{e_i-\tilde{E}\lambda_i^{\perp_{P_0}}\}_{i=1}^n$ (indeed, $R^{(\top)}$ maps $P_0$ to $P_0$, and so we just need to look at the effect on $R^{(\perp)}$), and so if $N$ is the $(n+k)\times n$ matrix whose column vectors are $\{e_i-\tilde{E}\lambda_i^{\perp_{P_0}}\}_{i=1}^n$, we therefore have
		\begin{equation}\label{E:b-4-12}
			\pi_{R(P_0)}(x) = N(N^TN)^{-1}N^T x.
		\end{equation}
		Since $x\in P_0$, from the definition of $N$ we clearly have $N^T x = x$. Now, $N^TN$ is equal to the $n\times n$ matrix $A = (A_{pq})^n_{p,q=1}$ given by
		$$A_{pq} = \delta_{pq} + \tilde{E}^2\sum^k_{\kappa=1}\lambda_p^{n+\kappa}\lambda_q^{n+\kappa},$$
		and so by Taylor expanding $A^{-1}$ about the identity we see that
		$$(N^TN)^{-1}N^T x = (N^TN)^{-1}x = x + o(\tilde{E}).$$
		Combining this with \eqref{E:b-4-12} and using the fact that $N$ is bounded we have $\pi_{R(P_0)}(x) = Nx + o(\tilde{E})$. Then by direct computation we deduce that
		$$\pi_{R(P_0)}(x) = x - \tilde{E}\lambda_{\pi}(\pi_{S^\prime}(x))^{\perp_{P_0}} + o(\tilde{E}).$$
		Combining this with \eqref{E:b-4-11} we conclude that
		\begin{equation}\label{E:b-4-13}
			\text{I} = -\tilde{E}\lambda_{\pi}(\pi_{S^\prime}(x))^{\perp_{P_0}} + o(\tilde{E})
		\end{equation}
		which is the desired expression for $\text{I}$.
		
		Next we compute $\text{II}$. Since $p^\alpha\in P_0^\perp$ with $|p^\alpha|_{B^n_1(0)}|\leq C\hat{E}$, we get from \eqref{E:b-4-8} and \eqref{E:b-4-10} that
		\begin{align*}
			\text{II} & = R(p^\alpha(x^\prime))-p^\alpha(x)\\
			& = p^\alpha(x^\prime)-p^\alpha(x) + o(\tilde{E})\\
			& = p^\alpha(\pi_{P_0}(R^{-1}(x))) - p^\alpha(x) + o(\tilde{E}).
		\end{align*}
		Since for $y\in P_0$ we have $p^\alpha(y) = (p^\alpha\circ\pi_{(S^\prime)^\perp})(y)$ (indeed, $p^\alpha$ is translation invariant along $S^\prime$), we can write the above as
		\begin{equation}\label{E:b-4-14}
			\text{II} = (p^\alpha\circ\pi_{(S^\prime)^\perp})(R^{-1}(x) - x) + o(\tilde{E})
		\end{equation}
		where here to remove the $\pi_{P_0}$ previously appearing, we are viewing $p^\alpha$ as a linear function on all of $\R^{n+k}$, extending it by $0$ in the $P_0^\perp$ directions. We of course have $R^{-1} = (R^{(\top)})^{-1}\circ (R^{(\perp)})^{-1}$, and so from \eqref{E:b-4-3} and the boundedness of $R^{(\top)}$, we get
		\begin{align*}
			R^{-1}(x) - x & = (R^{(\top)})^{-1}(x+(R^{(\perp)})^{-1}(x)-x) - x\\
			& = (R^{(\top)})^{-1}(x) - x + O(\tilde{E}).
		\end{align*}
		But since $|p^\alpha| = o(1)$ on $B^n_1(0)$ and $|\pi_{S^\perp}|\leq 1$, we deduce from this that \eqref{E:b-4-14} can be written as
		\begin{equation}\label{E:b-4-15}
			\text{II} = (p^\alpha\circ \pi_{(S^\prime)^\perp})\left((R^{(\top)})^{-1}(x) - x\right) + o(\tilde{E}).
		\end{equation}
		We then write $\text{II} = \text{II}_1 + \text{II}_2$, where
		$$\text{II}_1 := (p^\alpha\circ \pi_{(S^\prime)^\perp})\left((R^{(\top)})^{-1}(\pi_{S^\prime}(x)) - \pi_{S^\prime}(x)\right);$$
		$$\text{II}_2 := (p^\alpha\circ \pi_{(S^\prime)^\perp})\left((R^{(\top)})^{-1}(\pi_{(S^\prime)^\perp}(x)) - \pi_{(S^\prime)^\perp}(x)\right).$$
		To compute $\text{II}_1$, we get directly from the definition of $R^{(\top)}$ that
		$$(R^{(\top)})^{-1}(\pi_{S^\prime}(x)) = \|\pi_{S^\prime}(x)\|\cdot\frac{\pi_{S^\prime}(x) + \ell^{\top_{P_0}}(\pi_{S^\prime}(x))}{\|\pi_{S^\prime}(x) + \ell^{\top_{P_0}}(\pi_{S^\prime}(x))\|}.$$
		But since
		\begin{align*}
			\|\pi_{S^\prime}(x)+\ell^\top(\pi_{S^\prime}(x))\|^{-1} & = \left(\|\pi_{S^\prime}(x)\|^2 + \hat{E}^{-2}\tilde{E}^2\|\lambda_\pi(\pi_{S^\prime}(x))^{\top_{P_0}}\|^2\right)^{-1/2}\\
			& = \|\pi_{S^\prime}(x)\|^{-1}(1+O(\hat{E}^{-2}\tilde{E}^2))
		\end{align*}
		combining the above two expressions with \eqref{E:b-4-15} and again using that $|p^\alpha| \leq C\hat{E}$ on $B_1^n(0)$, and that $\lambda_\pi\in (S^\prime)^\perp$, we deduce that
		\begin{equation}\label{E:b-4-16}
			\text{II}_1 = \hat{E}^{-1}\tilde{E}Dp^\alpha\cdot\lambda_\pi(\pi_{S^\prime}(x))^{\top_{P_0}} + o(\tilde{E}).
		\end{equation}
		Now we compute $\text{II}_2$. Since $R^{(\top)}|_{P_0}$ (and hence also $(R^{(\top)})^{-1}|_{P_0}$) is an orthogonal transformation of $P_0\cong \R^n$, there are pairwise orthogonal $(R^{(\top)})^{-1}$-invariant subspaces $V_1,\dotsc,V_{m^\prime}$ with $P_0 = V_1\oplus\cdots\oplus V_{m^\prime}$ and such that $\dim(V_i) = 2$ for $1\leq i<m^\prime$ and $\dim(V_{m^\prime}) \in \{1,2\}$. We may therefore write
		\begin{equation}\label{E:b-4-17}
			\pi_{(S^\prime)^\perp}(x) = \sum^{m^\prime}_{i=1}\pi_{V_i}(\pi_{(S^\prime)^\perp}(x)).
		\end{equation}
		Note then the following:
		\begin{itemize}
			\item From \eqref{E:b-4-4}, if $\dim(V_{m^\prime})=1$ then we must have $R^{(\top)}|_{V_{m^\prime}} = \id_{V_{m^\prime}}$ (as the only non-trivial orthogonal $1$-dimensional transformation is the reflection, a distance $2$ from the identity), in which case we have
			$$(R^{(\top)})^{-1}(\pi_{V_{m^\prime}}(\pi_{(S^\prime)^\perp}(x))) - \pi_{V_{m^\prime}}(\pi_{(S^\prime)^\perp}(x)) = 0.$$
			\item If $i\in \{1,\dotsc,m^\prime\}$ is such that $V_i\subset S^\prime$, then $\pi_{V_i}(\pi_{(S^\prime)^\perp}(x)) = 0$, in which case we also have
			$$(R^{(\top)})^{-1}(\pi_{V_{i}}(\pi_{(S^\prime)^\perp}(x))) - \pi_{V_{i}}(\pi_{(S^\prime)^\perp}(x)) = 0.$$
			\item If $i\in \{1,\dotsc,m^\prime\}$ is such that $V_i\subset (S^\prime)^\perp$, then $(R^{(\top)})^{-1}|_{V_i^\perp}(S^\prime) = (R^{(\top)})^{-1}(S^\prime) = \graph(\ell^{\top_{P_0}})$. So, if we were to suppose (for the sake of contradiction) that $R^{(\top)}$ acted non-trivially on $V_i$, we would be able to define a rotation $R^\prime\in SO(P_0)$ with $R^\prime(\graph(\ell^{\top_{P_0}})) = S^\prime$ and $|R^\prime-\id|<|R^{(\top)}-\id|$, contradicting the definition of $R^{(\top)}$ (indeed, $R^\prime$ would simply equal $R^{(\top)}$ on $V_i^\perp$ and act trivially on $V_i$). Therefore, $(R^{(\top)})^{-1}|_{V_i} = \id_{V_i}$, and so again we have
			$$(R^{(\top)})^{-1}(\pi_{V_{i}}(\pi_{(S^\prime)^\perp}(x))) - \pi_{V_{i}}(\pi_{(S^\prime)^\perp}(x)) = 0.$$
		\end{itemize}
		The final case to consider is therefore if $i\in\{1,\dotsc,m^\prime\}$ is such that $\dim(V_i) = 2$, and $\dim(V_i\cap S^\prime) = \dim(V_i\cap (S^\prime)^\perp) = 1$. Hence, $(R^{(\top)})^{-1}|_{V_i}$ is a rotation by some angle $\alpha_i\in (0,2\pi)$ with $|\tan(\alpha_i)| \leq C\hat{E}^{-1}\tilde{E}$. Pick an orthonormal basis $v_1,v_2$ for $V$ with $v_1\in S^\prime$ and $v_2\in (S^\prime)^\perp$. Then clearly we have $\pi_{V_i}(\pi_{(S^\prime)^\perp}(x)) = av_2$ for some $a\in \R$. Therefore, if we represent the computation entirely in this basis, we have
		\begin{align*}
			\left|\pi_{S^\perp}\left((R^{(\top)})^{-1}(\pi_{V_i}(\pi_{(S^\prime)^\perp}(x))) - \pi_{V_i}(\pi_{(S^\prime)^\perp}(x))\right)\right| & = \left|\begin{pmatrix}0\\ 1 \end{pmatrix}\cdot\begin{pmatrix}\cos(\alpha_i)-1 & -\sin(\alpha_i)\\ \sin(\alpha_i) & \cos(\alpha_i)-1 \end{pmatrix}\begin{pmatrix}0\\ a\end{pmatrix}\right|\\
			& = |a|\cdot |1-\cos(\alpha_i)|\\
			& = |a| O(\alpha_i^2)\\
			& = |a| O(\hat{E}^{-2}\tilde{E}^2).
		\end{align*}
		Using again that $|p^\alpha| \leq C\hat{E}$, we deduce that
		\begin{equation*}
			\left|(p^\alpha\circ\pi_{(S^\prime)^\perp})\left((R^{(\top)})^{-1}(\pi_{V_i}(\pi_{(S^\prime)^\perp}(x))) - \pi_{V_i}(\pi_{(S^\prime)^\perp}(x))\right)\right| = o(\tilde{E}).
		\end{equation*}
		Combining all of the above with \eqref{E:b-4-17} and the definition of $\text{II}_2$, we get
		$$\text{II}_2 = o(\tilde{E}).$$
		Thus, combining this with \eqref{E:b-4-16}, we see that
		\begin{equation}\label{E:b-4-18}
			\text{II} = \hat{E}^{-1}\tilde{E}Dp^\alpha\cdot\lambda_\pi(\pi_{S^\prime}(x))^{\top_{P_0}} + o(\tilde{E}).
		\end{equation}
		Finally we turn to computing $\text{III}$. Since $u^\alpha\in P_0^\perp$ has $|u^\alpha| = o(1)$, \eqref{E:b-4-8} implies that $R(u^\alpha(x^\prime)) = u^\alpha(x^\prime) + o(\tilde{E})$. Then as $|x^\prime-x| = o(1)$ (as can be seen from \eqref{E:b-4-10} and \eqref{E:b-4-3}, \eqref{E:b-4-4}), the mean-value inequality and the bound $|Du|\leq C\tilde{E}$ imply that $u^\alpha(x^\prime) = u^\alpha(x) + o(\tilde{E})$, and so we conclude that
		\begin{equation}\label{E:b-4-19}
			\text{III} = u^\alpha(x) + o(\tilde{E}).
		\end{equation}
		Combining \eqref{E:b-4-13}, \eqref{E:b-4-18}, and \eqref{E:b-4-19} with \eqref{E:b-4-6} immediately gives \eqref{E:b-4-7}, thus completing the proof of the lemma.
\end{proof}

In addition to Lemma \ref{lemma:fine-integration}, we will also need several other closure properties for fine blow-ups. The first is a translation and dilation property analogous to $(\mathfrak{B}3\text{I})$ for coarse blow-ups. The second is the ability to subtract certain combinations of constants in $S^\perp$, analogous to the ability to subtract constants from coarse blow-ups seen in $(\mathfrak{B}3\text{II})$.

\begin{lemma}[Translation by $z\in\mathcal{D}(\phi)$ and Dilation]\label{lemma:fine-translation}
	For any $\sfrak\in \{0,1,\dotsc,n-1\}$, let $\phi = \{\phi^\alpha\}_{\alpha\in I}\in\mathfrak{b}_{\sfrak}$ be a fine blow-up of $(V_j)_j$ relative to $(\BC_j)_j$. Then for any $z\in \mathcal{D}(\phi)\in B^n_{1/4}(0)$ and $\rho\in (0,\frac{1}{2}(\frac{1}{4}-|z|))$ such that $\phi\not\equiv 0$ in $B_\rho(z)$, we have that $\phi_{z,\rho}\in \mathfrak{b}_{\sfrak}$ is a fine blow-up, where for $\alpha\in I$,
	$$\phi^\alpha_{z,\rho}(x):=\frac{\phi^\alpha(z+\rho x)}{\|\phi(z+\rho(\cdot))\|_{L^2(B_1^n(0))}}\, ,$$
	where here $\|\phi\|^2_{L^2(\bullet)} := \sum_{\alpha\in I}\|\phi^\alpha\|_{L^2(\bullet)}^2$. Indeed, we can find a sequence $z_j\to z$, $z_j\in \mathcal{D}(V_j)$, such that $\phi_{z,\rho}$ is the fine blow-up of $(\eta_{z_j,\rho})_\#V_j\res B_1(0)$ relative to $\BC_j$.
\end{lemma}

\begin{proof}
	Fix $z\in\mathcal{D}(\phi)$ and $\rho\in (0,\frac{1}{2}(\frac{1}{4}-|z|))$. Let $z_j\in\mathcal{D}(V_j)$ be such that $z_j\to z$. The result will follow if we can check Hypothesis $\dagger$ from Definition \ref{defn:dagger} holds (with $\beta/2$ in place of $\beta$, say) for the sequence $\widetilde{V}_j:= (\eta_{z_j,\rho})_\#V_j$ relative to the same sequence of cones $\BC_j$. In turn, this amounts to checking that Hypothesis $H(S(\BC_j),\tilde{\eps}_j,\tilde{\gamma}_j,\tilde{M})$ holds for suitable sequences $\tilde{\eps}_j,\tilde{\gamma}_j\downarrow 0$, and constant $\tilde{M}$ (which we will ensure is independent of $z$ and $\rho$). Hypothesis (A1) and (A2) are not difficult to check, and so we just need to verify (A3) and (A4).
	
	We start by verifying (A3). Firstly, we have by changing variables and the triangle inequality,
	\begin{align*}
		\hat{E}_{\widetilde{V}_j}^2 & = \rho^{-n-2}\int_{B_\rho(z_j)}\dist^2(x,z_j^{\perp_{P_0}}+P_0)\, \ext\|V_j\|(x)\\
		& \leq C\rho^{-n-2}\hat{E}_{V_j,\BC_j}^2 + C\rho^{-2}\cdot \rho^2\nu(\BC_j)^2 + C\rho^{-2}|z_j^{\perp_{P_0}}|^2\\
		& \leq C\rho^{-n-2}\gamma_j^2(\mathscr{F}^*_{V_j,\BC_j})^2 + C\hat{E}_{V_j}^2
	\end{align*}
	where in the last inequality we have used Theorem \ref{thm:shift-estimate} for $V_j$ and $\BC_j$ to control $|z^{\perp_{P_0}}|^2\equiv |z^{\perp_{S(\BC_j)}\perp_{P_0}}|^2$ by $C\hat{E}_{V_j,\BC_j}^2$, followed by (A4) for $V_j$ and $\BC_j$, and then finally Lemma \ref{lemma:nu-a}(i) for $V_j$ and $\BC_j$ to control $\nu(\BC_j)$ by $\hat{E}_{V_j}$. For $j$ sufficiently large (depending on $\rho$) the first term on the right-hand side can be absorbed into the left-hand side, and so we get
	$$\hat{E}_{\widetilde{V}_j} \leq C\hat{E}_{V_j} \leq CM\mathscr{E}_{V_j,\BC_j}$$
	where in the second inequality we have used (A3) for $V_j$ and $\BC_j$; here, $C = C(n,k)$. Thus, we will be done if we can prove that
	\begin{equation}\label{E:b-2-0}
		\mathscr{E}_{V_j,\BC_j} \leq C\mathscr{E}_{\widetilde{V}_j,\BC_j}
	\end{equation}
	for some $C = C(n,k)$. To see this, pick a plane $\widetilde{P}_j\supsetneq S(\BC_j)$ with $\hat{E}_{\widetilde{V}_j,\widetilde{P}_j} = \mathscr{E}_{\widetilde{V}_j,\BC_j}$. As we have seen in previous arguments, such a plane must satisfy $\nu(\widetilde{P}_j) \leq C\hat{E}_{\widetilde{V}_j}$ for some $C = C(n,k)$. Write $\tilde{p}_j$ for the linear function over $P_0$ whose graph is $\widetilde{P}_j$. Observe that
	$$\mathscr{E}_{V_j,\BC_j}^2 \leq \hat{E}_{V_j,\widetilde{P}_j}^2 \leq C\hat{E}_{V_j,\BC_j}^2 + C \nu(\BC_j,\widetilde{P}_j)^2.$$
	Again, using (A4) for $V_j$ and $\BC_j$, the first term on the right-hand side of the above can be absorbed into the left-hand side, provided $j$ is sufficiently large, giving that $\mathscr{E}_{V_j,\BC_j} \leq C\nu(\BC_j,\widetilde{P}_j)$. Using this, we have
	\begin{align*}
		\mathscr{E}_{V_j,\BC_j}^2 & \leq C\nu(\BC_j,\widetilde{P}_j)^2\\
		& \leq C\rho^{-n-2}\sum_{\alpha=1,2}\int_{B^n_{\rho/2}(z_j)\cap \{r_{\BC_j}>\rho/16\}}|p_j^\alpha-\tilde{p}_j|^2\\
		& \leq C\rho^{-n-2}\sum_{\alpha=1,2}\int_{B^n_{\rho/2}(z_j)\cap \{r_{\BC_j}>\rho/16\}}|u_j^\alpha|^2\, \ext x\\
		& \hspace{4em} + C\rho^{-n-2}\sum_{\alpha=1,2}\int_{B^n_{\rho/2}(z_j)\cap \{r_{\BC_j}>\rho/16\}}|p_j^\alpha+u^\alpha_j - (\tilde{p}_j + z_j^{\perp_{\widetilde{P}_j}})|^2\, \ext x + C\rho^{-2}|z_j^{\perp_{\widetilde{P}_j}}|^2\\
		& \leq C\rho^{-n-2}\hat{E}_{V_j,\BC_j}^2 + C\rho^{-n-2}\int_{B_\rho(z_j)}\dist^2(x,z_j+\widetilde{P}_j)\, \ext\|V_j\|(x)\\
		& \hspace{20em} + C\rho^{-2}|z_j^{\perp_{P_0}}|^2 + C\rho^{-2}\nu(\widetilde{P}_j)^2|z^{\perp_{S(\BC_j)}}|^2\\
		& C\rho^{-n-2}\hat{E}_{V_j,\BC_j}^2 + C\hat{E}_{\widetilde{V}_j,\widetilde{P}_j}^2 + C\rho^{-2}\hat{E}_{V_j,\BC_j}^2,
	\end{align*}
	where here, in the fourth inequality we have used that $S(\BC_j)\subset P_0\cap \widetilde{P}_j$ to get $|z_j^{\perp_{\widetilde{P}_j}}| \leq |z_j^{\perp_{P_0}}| + \|\pi_{P_0}-\pi_{\widetilde{P}_j}\|\cdot |z_j^{\perp_{S(\BC_j)}} \leq |z_j^{\perp_{P_0}}| + \nu(\widetilde{P}_j)|z_j^{\perp_{S(\BC_j)}}$, and then in the final inequality we have subsequently used that $\nu(\widetilde{P}_j) \leq C\hat{E}_{\widetilde{V}_j} \leq C\rho^{-n-2}\hat{E}_{V_j}$, and so combining with Lemma \ref{lemma:nu-a}(i) with $V_j$ and $\BC_j$ we get $\nu(\widetilde{P}_j) \leq C\nu(\BC_j)$, and so combining this with Theorem \ref{thm:shift-estimate} we see $|z_j^{\perp_{\widetilde{P}_j}}|\leq C\hat{E}_{V_j,\BC_j}$. Now, using (A4) for $V_j$ and $\BC_j$, we can absorb the first and third term on the right-hand side of the above into the left-hand side (for $j$ sufficiently large, depending on $\rho$), thus giving $\mathscr{E}_{V_j,\BC_j} \leq C\hat{E}_{\widetilde{V}_j,\widetilde{P}_j} \equiv C\mathscr{E}_{\widetilde{V}_j,\BC_j}$, which is \eqref{E:b-2-0}, as desired (with a constant depending only on $n,k$). This verifies that (A3) holds for $\widetilde{V}_j$ with a constant $\tilde{M} = CM$, where $C = C(n,k)$.
	
	Next we verify (A4) for $\widetilde{V}_j$ and $\BC_j$. We have for all $j$ sufficiently large, by changing variables and the triangle inequality,
	$$\hat{F}_{\widetilde{V}_j,\BC_j}^2\leq C\rho^{-n-2}\hat{F}_{V_j,\BC_j}^2 + C\rho^{-2}\cdot \rho^2\nu(\BC_j,(\tau_{z_j})_\#\BC_j)^2$$
	(we stress that, when undoing the transformation on the two-sided fine excess, the extra term when compared with the one-sided fine excess is not obviously controlled by the equivalent term in $\hat{F}_{V_j,\BC_j}$; however, using Theorem \ref{thm:graphical-rep}, for $j$ sufficiently large depending on $\rho$, it is comparable to the one-sided fine excess of $V_j$ relative to $\BC_j$; cf.~Remark \ref{remark:after-graphical-rep} following Theorem \ref{thm:graphical-rep}).
	
	Now using \eqref{E:cor-2} for $V_j$ and $\BC_j$, as well as Hypothesis (A4) for $V_j$ and $\BC_j$, we get
	\begin{equation}\label{E:b-2-1}
		\hat{F}_{\widetilde{V}_j,\BC_j} \leq \tilde{\gamma}_j\mathscr{F}^*_{V_j,\BC_j}
	\end{equation}
	where $\tilde{\gamma}_j\downarrow 0$. We now claim that
	\begin{equation}\label{E:b-2-2}
		\mathscr{F}^*_{V_j,\BC_j}\leq C\mathscr{F}^*_{\widetilde{V}_j,\BC_j}.
	\end{equation}
	Combining \eqref{E:b-2-2} with \eqref{E:b-2-1} would then give the desired property (A4) for $\widetilde{V}_j$ and $\BC_j$. To see \eqref{E:b-2-2}, pick a coarser cone $\Dbf_j\in \mathcal{P}$ to $\BC_j$ satisfying $\hat{F}_{\widetilde{V}_j,\Dbf_j} \leq \frac{3}{2}\mathscr{F}^*_{\widetilde{V}_j,\BC_j}$. It then follows that we have
	$$\nu(\BC_j,\Dbf_j) \leq C\hat{F}_{\widetilde{V}_j,\Dbf_j}.$$
	(Alternatively, one can specifically choose $\Dbf_j$ as in Step 1 of the proof of Theorem \ref{thm:graphical-rep} and then this follows from \eqref{E:graphical-rep-16} in the case $\Dbf_j$ therein is two planes, or otherwise Lemma \ref{lemma:nu-a}(i) if $\Dbf_j$ is a plane.)
	
	Thus, we have $\nu(\BC_j,\Dbf_j) \leq C\mathscr{F}^*_{\widetilde{V}_j,\BC_j}$. Then from the triangle inequality we get
	$$\mathscr{F}_{V_j,\BC_j}^* \leq \hat{F}_{V_j,\Dbf_j} \leq C\hat{F}_{V_j,\BC_j} + C\nu(\BC_j,\Dbf_j)$$
	and then using (A4) we may absorb the first term on the right-hand side into the left-hand side. This therefore gives $\mathscr{F}^*_{V_j,\BC_j} \leq C\nu(\BC_j,\Dbf_j) \leq C\mathscr{F}^*_{\widetilde{V}_j,\BC_j}$, which is exactly \eqref{E:b-2-2}. This shows that (A4) does hold for $\widetilde{V}_j$ and $\BC_j$, and consequently we can take a fine blow-up of $\widetilde{V}_j$ relative to $\BC_j$. Using the fact that
	$$\frac{\hat{E}_{(\eta_{z_j,\rho})_\#V_j,\BC_j}}{\hat{E}_{V_j,\BC_j}}\to \|\phi(z+\rho(\cdot))\|_{L^2(B_1(0))}$$
	as $j\to\infty$ (which follows from the $L^2_{\text{loc}}(B_1(0))$ convergence to the fine blow-up $\phi$, which is a consequence of the $L^2$ non-concentration estimate \eqref{E:cor-3}, we deduce the result of the lemma.
\end{proof}

\begin{lemma}[Translation perpendicular to $S$]\label{lemma:fine-perp-translation}
	Suppose $\phi = \{\phi^\alpha\}_{\alpha\in I}\in \mathfrak{b}$ is a fine blow-up of $(V_j)_j$ relative to $(\BC_j)_j$. Then, for any $\kappa\in S^\perp$ with $\phi^\alpha\neq \kappa^{\perp_{P_0}}-m^\alpha\cdot \kappa^{\top_{P_0}}$ for some $\alpha\in I$, we have that
	$$x\longmapsto \frac{\phi^\alpha(x) - (\kappa^{\perp_{P_0}} - m^\alpha\cdot\kappa^{\top_{P_0}})}{\left(\sum_{\alpha^\prime\in I}\|\phi^{\alpha^\prime}-(\kappa^{\perp_{P_0}}-m^{\alpha^\prime}\cdot\kappa^{\top_{P_0}})\|^2_{L^2(B^n_1(0))}\right)^{1/2}}, \qquad \alpha\in I$$
	is also a cone blow-up. Indeed, we can find a sequence $\sigma_j\uparrow 1$ such that it is the fine blow-up of $(\eta_{\tilde{E}_j\kappa,\sigma_j})_\#V_j\res B_1(0)$ relative to $\BC_j$.
\end{lemma}

\begin{proof}
	First, fix $\sigma\in (0,1)$ sufficiently close to $1$ so that $\phi^\alpha\neq \kappa^{\top_{P_0}}-m^\alpha\cdot\kappa^{\top_{P_0}}$ on $B_\sigma(0)$. We claim that we can take the fine blow-up of $\widetilde{V}_j:= (\eta_{\tilde{E}_j\kappa,\sigma})_\#V_j$ relative to $\BC_j$. For this, we need to check Hypothesis $\dagger$ (with $\beta/2$ in place of $\beta$, say), which as for the proof of Lemma \ref{lemma:fine-translation} amounts to checking Hypotheses (A3) and (A4) hold for appropriate $\tilde{\eps}_j,\tilde{\gamma}_j\downarrow 0$ and constant $\tilde{M}$. This is however analogous to those estimates involved in the proof of Lemma \ref{lemma:fine-translation} (they are even easier, in fact, as we are shifting $V_j$ by an amount whose size is the same order as $\hat{E}_{V_j,\BC_j}$), and so we omit the details. Given this for a fixed $\sigma\in (0,1)$, the general case follows by a diagonalisation argument analogous to that seen in the proof of $(\mathfrak{B}3\text{IV})$ from Theorem \ref{thm:blow-up-properties}, and as such we will not repeat the details.
\end{proof}

We now introduce the following useful terminology for fine blow-ups $\phi\in \mathfrak{b}$: for $z\in\mathcal{D}(\phi)\cap B_{1/4}(0)$ with $\phi(z)=0$ and $\rho\in (0,1]$, we say that $\pi\in\mathfrak{l}$ \emph{dehomogenises} $\phi$ in $B_\rho(z)$ if
$$\int_{B_\rho(z)}\left|\phi(x)-\pi\left(\frac{x-z}{\rho}\right)\right|^2\, \ext x = \inf_{\tilde{\pi}\in\mathfrak{l}}\int_{B_\rho(z)}\left|\phi(x) - \tilde{\pi}\left(\frac{x-z}{\rho}\right)\right|^2\, \ext x.$$
Note that the existence $\pi\in\mathfrak{l}$ which dehomogenises $\phi$ is immediate from the compactness of the space $\mathfrak{l}$ under a uniform $L^2$ bound, as indeed from the continuity estimate given by Theorem \ref{thm:fine-continuity-2} we know that for any fine blow-up $\phi$ and $z,\rho$ as above:
$$\rho^{-n}\int_{B_\rho(z)}|\phi(x)|^2\, \ext x \leq C$$
where $C = C(n,k,M,\beta)$. As such, any dehomogeniser satisfies (as $\pi\equiv 0$ is a valid choice of competitor in the infimum) $\int_{B_\rho(z)}|\pi((x-z)/\rho)|^2\, \ext x \leq 2\int_{B_\rho(z)}|\phi|^2 \leq 2C\rho^n$, and so by changing coordinates and using linearity of $\rho$ we get $\int_{B_1(0)}|\pi(x)|^2\, \ext x \leq 2C$. We will write $\pi^\phi_{z,\rho}(x):= \pi\left(\frac{x-z}{\rho}\right)$ for the (recentred) dehomogeniser of $\phi$ on $B_\rho(z)$, although often we will just write $\pi_{z,\rho}$ as $\phi$ will be clear from context.

By combining Lemma \ref{lemma:fine-integration}, Lemma \ref{lemma:fine-translation}, and Lemma \ref{lemma:fine-perp-translation}, we therefore see that fine blow-ups always obey an improved Hardt--Simon inequality with respect to their dehomogenisers:

\begin{lemma}[Hardt--Simon Inequality for Fine Blow-Ups]\label{lemma:fine-Hardt-Simon}
	Suppose that $\phi = \{\phi^\alpha\}_{\alpha\in I}\in \mathfrak{b}$ is a fine blow-up, and let $z\in \mathcal{D}(\phi)\cap B^n_{1/4}(0)$ and $\rho\in (0,1/4)$. Then, we have
	\begin{align*}
		\sum_{\alpha\in I}\int_{B^n_{\rho/2}(z)}&R_z(x)^{2-n}\left|\frac{\del}{\del R_z}\left(\frac{\phi^\alpha(x)-(\lambda(z)^{\perp_{P_0}}-m^\alpha\cdot\lambda(z)^{\top_{P_0}})}{R_z(x)}\right)\right|^2\, \ext x\\
		& C\rho^{-n-2}\sum_{\alpha\in I}\int_{B_\rho^n(z)}|\phi^\alpha(x) - (\lambda(z)^{\perp_{P_0}} - m^\alpha\cdot\lambda(z)^{\top_{P_0}} + \pi^\phi_{z,\rho}(x))|^2\, \ext x
	\end{align*}
	where $R_z\equiv R_z(x) := |x-z|$ and $C = C(n,k,M,\beta)\in (0,\infty)$.
\end{lemma}

\begin{proof}
	Apply Lemma \ref{lemma:fine-perp-translation} with $\kappa = \lambda(z)$ ($\equiv \phi(z)\in S^\perp$), then Lemma \ref{lemma:fine-translation} at $z$ and scale $\rho$, to get that
	$$x\longmapsto \frac{\phi(z+ \rho x) - (\lambda(z)^{\perp_{P_0}} - m^\alpha\cdot\lambda(z)^{\top_{P_0}})}{\left(\sum_{\alpha\in I}\|\phi^\alpha(z + \rho(\cdot)) -  (\lambda(z)^{\perp_{P_0}} - m^\alpha\cdot\lambda(z)^{\top_{P_0}})\|_{L^2(B_1(0))}^2\right)^{1/2}}$$
	is a fine blow-up. Then, apply Lemma \ref{lemma:fine-integration} to this fine blow-up with $\pi\in \mathfrak{l}$ therein being the dehomogeniser of this at $z=0$ and $\rho=1$ normalised by
	$$\left(\sum_{\alpha\in I}\|\phi^\alpha(z + \rho(\cdot)) -  (\lambda(z)^{\perp_{P_0}} - m^\alpha\cdot\lambda(z)^{\top_{P_0}})\|_{L^2(B_1(0))}^2\right)^{1/2}.$$
	Notice that, as we saw for general dehomogisers that the $L^2$ norm was controlled by twice the $L^2$ norm of the fine blow-up being dehomogised, we have a uniform upper bound on the $L^2$ norm of this dehomogeniser depending only on $n,k,M,\beta$. As such, if we apply the Hardt--Simon inequality given in Lemma \ref{lemma:fine-integration} and perform a change of variables, we arrive exactly at the claim of the lemma, with the constant $C$ depending only on $n,k,M,\beta$.
\end{proof}

The last observation we need is a compactness property not only for fine blow-ups, but also for fine blow-ups minus their dehomogenisers on the ball $B_1(0)$.

\begin{lemma}[Compactness for Fine Blow-Ups]\label{lemma:fine-compactness}
	Suppose that for $p\in \Z_{\geq 1}$, $\phi_p = \{\phi^\alpha_p\}_{\alpha\in I_p}\in\mathfrak{b}$ is a sequence of fine blow-ups such that the parameter $M = M_p$ used in their construction obeys $M:=\limsup_p M_p<\infty$, and that the parameter $\beta = \beta_p$ obeys $\beta:= \liminf_p\beta_p>0$. Then, there is a subsequence $\{p^\prime\}\subset\{p\}$ and a fine blow-up $\phi = \{\phi^\alpha\}_{\alpha=1}^4\in\mathfrak{b}_{n-1}$ such that:
		\begin{enumerate}
			\item [(i)] $\phi^\alpha_{p^\prime}\to \phi^\alpha$ in $L^2(K;P_0^\perp)$ for every compact $K\subset (U^\alpha_0\cup S_0)\cap B^n_1(0)$;
			\item [(ii)] $\dist_\H(\mathcal{D}(\phi_{p^\prime})\cap K,\mathcal{D}(\phi)\cap K)\to 0$ for every compact $K\subset B^n_1(0)$;
			\item [(iii)] If $\mathcal{D}(\phi)\subsetneq S_0\cap B^n_1(0)$, then $\mathcal{D}(\phi_{p^\prime})\subsetneq S_0\cap B^n_1(0)$ for all sufficiently large $p^\prime$, and $\phi^\alpha_{p^\prime}\to \phi^\alpha$ in $C^2_{\text{loc}}(B^n_1(0)\setminus D(\phi);P_0^\perp)$ for $\alpha=1,2$ (where here we are viewing the fine blow-ups are $\phi^1,\phi^2$, both defined over all of $P_0$).
		\end{enumerate}
		Indeed, if $v_p$ is the fine blow-up of $(V_{p,j})_{j=1}^\infty$ relative to $(\BC_{p,j})_{j=1}^\infty$, then there is an increasing $\{j_{p^\prime}\}_{p^\prime}$ such that $\phi$ is the fine blow-up of $(V_{p^\prime,j_{p^\prime}})_{p^\prime}$ relative to $(\BC_{p^\prime,j_{p^\prime}})_{p^\prime}$, with parameters $M$ and $\beta$ therein.
		
		Furthemore, if $\phi_p(0) = 0$ for all $p$, and $\pi_p:= \pi^{\phi_p}_{0,1}$ is the dehomogeniser of $\phi_p$ on the ball $B_1(0)$, then the same compactness conclusions hold for the sequence $\phi_p - \pi_p$, with the limit being also of the form $\phi - \pi$ where $\pi$ is the dehomogeniser of $\phi$ on $B_1(0)$.
\end{lemma}

\textbf{Remark:} We will refer to the convergence described in Lemma \ref{lemma:fine-compactness} by saying ``converge in the sense of fine blow-ups''.

\begin{proof}
	The proof follows from a diagonalisation argument analogous to that seen in the proof of $(\mathfrak{B}3\text{IV})$ in Theorem \ref{thm:blow-up-properties}, and as such we will not repeat the details. For the claim concerning $\phi_p-\pi_p$, notice that as seen in the proof of Lemma \ref{lemma:fine-translation}, we still have strong $L^2_{\text{loc}}$ convergence to the blow-up argument, and the convergence to $\phi-\pi$ with $\pi$ being the dehomogeniser can be seen by a simple contradiction argument, supposing that $\pi$ was not and using the strong $L^2$ convergence in some $B_\sigma(0)$ with $\sigma\in (0,1)$ suitably close to $1$.
\end{proof}

\textbf{Remark:} Notice that from the Hölder estimate given in Theorem \ref{thm:fine-continuity-2}, the $L^2_{\text{loc}}$ convergence in Lemma \ref{lemma:fine-compactness} can actually be upgraded to uniform convergence locally in $B_1(0)$. For the sequence of fine blow-ups this follows directly from Theorem \ref{thm:fine-continuity-2}, whilst for the case $\phi_p-\pi_p$, one needs the equivalent Hölder estimate. This follows from in the same manner as in Theorem \ref{thm:fine-continuity-2} once one has the corresponding $L^2$ decay estimate involving dehomogenisers, namely for $z\in\mathcal{D}(\phi)\cap B_{1/4}(0)$ and $\rho\in (0,1/4)$,
\begin{align*}
\sum_{\alpha\in I}\int_{B^n_{\rho/2}(z)}&\frac{|\phi^\alpha(x) - (\lambda(z)^{\perp_{P_0}} - m^\alpha\cdot\lambda(z)^{\top_{P_0}} + \pi^\phi_{z,\rho}(x))|^2}{|x-z|^{n+7/4}}\, \ext x\\
& \leq C\rho^{-n-7/4}\sum_{\alpha\in I}\int_{B_\rho^n(z)}|\phi^\alpha(x)-(\lambda(z)^{\perp_{P_0}}-m^\alpha\cdot\lambda(z)^{\top_{P_0}} + \pi^\phi_{z,\rho}(x))|^2\, \ext x
\end{align*}
where $C = C(n,k,M,\beta)\in (0,\infty)$. This estimate follows in a similar manner to how the Hardt--Simon inequality from Lemma \ref{lemma:fine-integration} did, except now using \eqref{E:initial-3} of Lemma \ref{lemma:L2-estimates-2}.

\subsection{Classification of Homogeneous Degree One Fine Blow-Ups}

We are now ready to completely classify homogeneous degree one fine blow-ups. Recall that the only case left we need to do this for is $\mathfrak{b}_{n-1}$. This classification is the last remaining step towards the final global regularity result for fine blow-ups.

\textbf{Remark:} When talking about homogeneous degree one fine blow-ups here, we will always work with their homogeneous extensions to the whole of $P_0$.

\begin{theorem}\label{thm:fine-classification}
	Suppose that $\phi = \{\phi^\alpha\}_{\alpha\in I}\in \mathfrak{b}_{n-1}$ is a fine blow-up which is homogeneous of degree one. Then, $\phi\in \mathfrak{l}$.
\end{theorem}

\begin{proof}[Proof of Theorem \ref{thm:fine-classification}]
	Our proof will closely follow that seen in \cite[Theorem 5.4]{BK17}.
	
	\textbf{Step 1.} \emph{Easy cases.} First suppose that $\mathcal{D}(\phi)$ obeys $\H^{n-2}(\mathcal{D}(\phi)\cap B^n_1(0))<\infty$. In particular, this implies $I = \{1,2\}$. In this case, $\mathcal{D}(\phi)$ is a removable set for the functions $\phi^1,\phi^2$, which are harmonic on $P_0\setminus\mathcal{D}(\phi)$ and continuous (by Theorem \ref{thm:fine-continuity-2}), meaning that $\phi^1,\phi^2$ are both smooth and harmonic on all of $P_0$. Consequently, they are linear, and the result follows. Thus, we can assume that $\H^{n-2}(\mathcal{D}(\phi)\cap B^n_1(0)) = \infty$.
	
	Notice that by Lemma \ref{lemma:fine-continuity-1}, we know that at every $z\in \mathcal{D}(\phi)$ we have $\phi^\alpha(z) = \lambda(z)^{\perp_{P_0}} - m^\alpha\cdot\lambda(z)^{\top_{P_0}}$, and by Lemma \ref{lemma:fine-lambda} we in fact know that $\lambda$ is a linear function. In particular, if we have $\mathcal{D}(\phi)\cap B_{1/4}(0) = S_0\cap B_{1/4}(0)$, then by standard boundary regularity theory for harmonic functions, we see that $\phi^\alpha$ is a smooth harmonic function on $\overline{U^\alpha_0\cap B^n_{1/4}(0)}$ for each $\alpha$. However, for homogeneous degree one harmonic functions, the half-space Dirichlet problem is unique up to the addition of a linear function of the form $x\mapsto x^n c$ for some $c\in P_0^\perp$, which therefore implies that $\phi^\alpha$ is linear, and so in particular $\phi\in \mathfrak{l}$. So, we can henceforth assume that $\mathcal{D}(\phi)\cap B^n_{1/4}(0)\subsetneq S_0\cap B_{1/4}^n(0)$; in particular, $I = \{1,2\}$ (although, at times it will be convenient to view $\phi$ as being indexed by $\{1,\dotsc,4\}$.
	
	\textbf{Step 2.} \emph{Setting up the induction for the harder case.} For any homogeneous degree one fine blow-up $\phi\in\mathfrak{b}_{n-1}$, we can write
	$$S(\phi):= \{z\in S_0:\phi(x+z) = \phi(x) \text{ for all }x\in B_1\}.$$
	Using the homogeneity of $\phi$ it is easy to verify that $S(\phi)$ is a linear subspace of $S_0$. Write $d:= \dim(S(\phi))\in \{0,1,\dotsc,n-1\}$. We will prove by (downwards) induction on $d$ that whenever $\H^{n-2}(\mathcal{D}(\phi)) = \infty$ and $\mathcal{D}(\phi)\cap B_{1/4}^n(0)\subsetneq S_0\cap B_{1/4}^n(0)$, we must have $\phi\in\mathfrak{l}$. Notice that when $d=n-1$, then $S(\phi) = S_0$, from which the homogeneity of $\phi$ immediately gives that $\phi\in\mathfrak{l}$ (and in fact as $\mathcal{D}(\phi)\neq S_0$, $\phi$ must be two linear functions). We may therefore fix $d\in\{n-2,n-3,\dotsc,0\}$ and assume inductively we have shown the result for all larger values of $d$; we may also assume the same for $\psi:=\phi - \pi^\phi_{0,1}$ (i.e.~if $S(\psi)>d$ and $\H^{n-2}(\mathcal{D}(\phi))=\infty$, $\mathcal{D}(\phi)\cap B_{1/4}^n(0)\subsetneq S_0\cap B^n_{1/4}(0)$, then $\psi\in \mathfrak{l}$). Notice in this case, as $\H^{n-2}(\mathcal{D}(\phi)) = \infty$ by assumption, we have that $\mathcal{D}(\phi)\setminus S(\phi)\neq \emptyset$.
	
	For $z\in\mathcal{D}(\phi)\cap B_{1/4}^n(0)$ where $\phi(z) = 0$ and $\rho\in (0,1/4)$, recall that $\pi^\phi_{z,\rho}\equiv \pi_{z,\rho}$ is the function dehomogenising $\phi$ in $B_\rho(z)$. If $\phi\equiv \pi_{z,\rho}$ then $\phi\in\mathfrak{l}$ and so there is nothing to prove. So, we may assume (looking for contradiction) that $\phi-\pi_{z,\rho}\not\equiv 0$. Fix any compact subset $K$ of $(S_0\setminus S(\phi))\cap B_{1/4}^n(0)$. We then claim the following:
	
	\textbf{Claim:} There exists $\eps = \eps(\phi,K)\in (0,1)$ such that for any $y\in \mathcal{D}(\phi)\cap K$ and any $\rho\in (0,\eps]$, we have the \emph{reverse Hardt--Simon inequality}:
	\begin{equation}\label{E:homo-fine-1}
		\int_{B_\rho(y)\setminus B_{\rho/2}(y)}R_y^{2-n}\left|\frac{\del}{\del R_y}\left(\frac{\phi - \phi(y)}{R_y}\right)\right|^2 \geq \eps\rho^{-n-2}\int_{B_\rho(y)}|\phi-(\phi(y) + \pi_{y,\rho})|^2
	\end{equation}
	where $R_y\equiv R_y(x) = |y-x|$ and $\phi^\alpha(z) \equiv \lambda(z)^{\perp_{P_0}} - m^\alpha\cdot\lambda(z)^{\top_{P_0}}$.
	
	Indeed, if this claim were false then there would exist a sequence of points $y_j\in \mathcal{D}(\phi)\cap K$ and radii $\rho_j\downarrow 0$ with $\phi- (\phi(y) + \pi_{y_j,\rho_j})\not\equiv 0$ (if it were, \eqref{E:homo-fine-1} is trivially satisfied), and a sequence $\eps_j\downarrow 0$ such that, setting $\pi^j:= \pi_{y_j,\rho_j}$ (and recalling that as $\pi^j$ is homogeneous of degree $1$ about $y_j$, $\frac{\del}{\del R_{y_j}}\left(\frac{\pi^j}{R_{y_j}}\right) = 0$),
	$$\int_{B_{\rho_j}(y_j)\setminus B_{\rho_j/2}(y_j)}R_{y_j}^{2-n}\left|\frac{\del}{\del R_{y_j}}\left(\frac{\phi-(\phi(y_j) + \pi^j)}{R_{y_j}}\right)\right|^2 < \eps_j\rho_j^{-n-2}\int_{B_{\rho_j}(y_j)}|\phi-(\phi(y_j) + \pi^j)|^2.$$
	By Lemma \ref{lemma:fine-translation} and Lemma \ref{lemma:fine-integration}, we then have that if
	$$w_j := \frac{\phi(y_j + \rho_j(\cdot)) - (\phi(y_j) + \pi^j(y_j+\rho_j(\cdot)))}{\|\phi(y_j+\rho_j(\cdot)) - (\phi(y_j) + \pi^j(y_j+\rho_j(\cdot)))\|_{L^2(B_1)}}$$
	then this is a blow-up in the sense of Lemma \ref{lemma:fine-translation}. In particular, we may apply the compactness result Lemma \ref{lemma:fine-compactness} to pass to a subsequence for which $w_j\to w$, where the convergence is as fine blow-ups as described in Lemma \ref{lemma:fine-compactness}; in particular, the function $w$ is dehomogenised by $0$ on the ball $B_1(0)$, as $w = \tilde{\phi} - \tilde{\pi}$ for some fine blow-up $\tilde{\phi}\in\mathfrak{b}_{n-1}$ and $\tilde{\pi}\in \mathfrak{l}$ dehomogenising $\tilde{\phi}$ on $B_1(0)$. Furthermore, by the compactness of $K$ we can also assume that $y_j\to y\in\mathcal{D}(\phi)\cap K$. Writing our contradiction assumption in terms of $w_j$, we have
	$$\int_{B_1(0)\setminus B_{1/2}(0)}R^{2-n}\left|\frac{\del}{\del R}\left(\frac{w_j}{R}\right)\right|^2 < \eps_j$$
	which from the convergence properties to $w$ implies
	$$\int_{B_1(0)\setminus B_{1/2}(0)}R^{2-n}\left|\frac{\del}{\del R}\left(\frac{w}{R}\right)\right|^2 = 0$$
	giving that $w$ is homogeneous of degree one in $B_1(0)\setminus B_{1/2}(0)$. By unique continuation of harmonic functions (as fine blow-ups and functions in $\mathfrak{l}$ are always harmonic away from $S_0$), it is therefore equal to its homogeneous degree one extension in $B_1(0)$. Since $w_j\to w$ strongly in $L^2(B_1)$, it also follows that $w$ is dehomogenised by $0$ in $B_1(0)$. Also, notice that $S(\phi)\subset S(w)$ by construction (noting that if $\phi$ is translation invariant along a given direction, the dehomogeniser can also be taken to be translation invariant in those directions).
	
	A well-known argument in this field now implies that $y\in S(w)$; see, for example, \cite[Lemma 3.12]{MW24} (the same argument works here). Since $y\not\in S(\phi)$, this gives $\dim(S(w))>\dim(S(\phi))$. Recalling that $w = \tilde{\phi}-\tilde{\pi}$, if $\H^{n-2}(\tilde{\phi})=\infty$ and $\mathcal{D}(\tilde{\phi})\cap B_{1/4}(0)\subsetneq S_0\cap B_{1/4}(0)$, then by the induction hypothessis we know that $w\in\mathfrak{l}$; if either of these conditions fail, then arguing as in Step 1 we also have that $w\in\mathfrak{l}$, and thus we necessarily have $w\in\mathfrak{l}$. But as $w$ is dehomogenised by $0$ on $B_1(0)$, this would imply that in fact $w\equiv 0$. However, this leads to a contradiction by a standard integration-along-rays argument (see \cite[Proof of Theorem 3.11]{MW24}), which would give that we should have $\|w\|_{L^2(B_{3/4})}>0$ in the present setting. This contradiction shows the validity of the claim.
	
	\textbf{Step 3.} \emph{$C^1$ regularity on half-planes away from $S(\phi)$.} Let $K$ again be any compact subset of $(S_0\setminus S(\phi))\cap B_{1/4}(0)$. Then, \eqref{E:homo-fine-1} gives
	\begin{equation}\label{E:homo-fine-2}
		\eps\rho^{-n-2}\int_{B_\rho(y)}|\phi-(\phi(y) + \pi_{y,\rho})|^2 \leq \int_{B_\rho(y)\setminus B_{\rho/2}(y)}R_y^{2-n}\left|\frac{\del}{\del R_y}\left(\frac{\phi - \phi(y)}{R_y}\right)\right|^2.
	\end{equation}
	Also, Lemma \ref{lemma:fine-Hardt-Simon} gives
	\begin{equation}\label{E:homo-fine-3}
		\int_{B_{\rho/2}(y)}R_y^{2-n}\left|\frac{\del}{\del R_y}\left(\frac{\phi - \phi(y)}{R_y}\right)\right|^2 \leq C\rho^{-n-2}\int_{B_\rho(y)}|\phi - (\phi(y) + \pi_{y,\rho})|^2.
	\end{equation}
	Combining these two inequalities gives
	$$\frac{\eps}{C}\int_{B_{\rho/2}(y)}R_y^{2-n}\left|\frac{\del}{\del R_y}\left(\frac{\phi - \phi(y)}{R_y}\right)\right|^2 \leq \int_{B_\rho(y)\setminus B_{\rho/2}(y)}R_y^{2-n}\left|\frac{\del}{\del R_y}\left(\frac{\phi - \phi(y)}{R_y}\right)\right|^2.$$
	Adding $\int_{B_{\rho/2}(y)}R_y^{2-n}\left|\frac{\del}{\del R_y}\left(\frac{\phi - \phi(y)}{R_y}\right)\right|^2$ to both sides (`hole-filling') we get
	$$\int_{B_{\rho/2}(y)}R_y^{2-n}\left|\frac{\del}{\del R_y}\left(\frac{\phi - \phi(y)}{R_y}\right)\right|^2 \leq \vartheta\int_{B_\rho(y)}R_y^{2-n}\left|\frac{\del}{\del R_y}\left(\frac{\phi - \phi(y)}{R_y}\right)\right|^2$$
	where $\vartheta := (1+\eps/C)^{-1}\in (0,1)$. Iterating this with $2^{-i}\rho$ in place of $\rho$ and using a standard argument to interpolate between these scales, we deduce that
	$$\int_{B_\sigma(y)}R_y^{2-n}\left|\frac{\del}{\del R_y}\left(\frac{\phi - \phi(y)}{R_y}\right)\right|^2 \leq C\left(\frac{\sigma}{\rho}\right)^{2\mu}\int_{B_\rho(y)}R_y^{2-n}\left|\frac{\del}{\del R_y}\left(\frac{\phi - \phi(y)}{R_y}\right)\right|^2$$
	for all $0<\sigma\leq \rho/2\leq \eps/4$, where here $\mu = \mu(\phi,K)\in (0,1)$ and $C = C(n,k)$. Using this with \eqref{E:homo-fine-2} and \eqref{E:homo-fine-3}, we readily deduce that
	\begin{equation*}
		\sigma^{-n-2}\int_{B_\sigma(y)}|\phi- (\phi(y) + \pi_{y,\sigma})|^2 \leq C\eps^{-1}\left(\frac{\sigma}{\rho}\right)^{2\mu}\cdot\rho^{-n-2}\int_{B_\rho(y)}|\phi- (\phi(y) + \pi_{y,\rho})|^2
	\end{equation*}
	for all $0<\sigma\leq\rho/2\leq \eps/8$. Using this and the triangle inequality, it is then straightforward to check that there exists a \emph{single} $\pi_y\in \mathfrak{l}$ for which
	\begin{equation}\label{E:homo-fine-5}
		\sigma^{-n-2}\int_{B_\sigma(y)}|\phi-(\phi(y) + \pi_y)|^2 \leq C\eps^{-1}\left(\frac{\sigma}{\rho}\right)^{2\mu}\cdot\rho^{-n-2}\int_{B_\rho(y)}|\phi- (\phi(y) + \pi_y)|^2
	\end{equation}
	for all $0<\sigma\leq\rho/2\leq\eps/8$. To summarise we have therefore shown that for each $y\in\mathcal{D}(\phi)\cap K$, there is $\pi_y\in \mathfrak{l}$ for which
	\begin{align}
		\nonumber\sigma^{-n-2}\int_{B_\sigma(y)}|\phi-(\lambda(y)^{\top_{P_0}}-\, &\, m\cdot\lambda(y)^{\top_{P_0}} - \pi_y)|^2\\
		& \leq \beta\left(\frac{\sigma}{\rho}\right)^{2\mu}\cdot\rho^{-n-2}\int_{B_\rho(y)}|\phi-(\lambda(y)^{\top_{P_0}}-m\cdot\lambda(y)^{\top_{P_0}} - \pi_y)|^2\label{E:homo-fine-6}
	\end{align}
	for all $0<\sigma\leq\rho/2\leq \eps/8$, where here we stress that $\eps = \eps(\phi,K)\in (0,1)$ and $\beta = \beta(\phi,K)\in (0,\infty)$.
	Away from $\mathcal{D}(\phi)$, we know that $\phi$ is harmonic. Thus, we know that for any $x_0 \not\in \mathcal{D}(\phi)$ and any affine function $\ell:\R^n\to \R^k$, we have
	\begin{align*}
		\sigma^{-n-2}\int_{B_\sigma(x_0)}|\phi(x)-(\phi(x_0) + (x-x_0)\cdot D\phi(x_0))|^2 \leq C\left(\frac{\sigma}{\rho}\right)^2\cdot\rho^{-n-2}\int_{B_\rho(x_0)}|\phi-\ell|^2
	\end{align*}
	for any $0<\sigma\leq\rho/2\leq \frac{1}{2}\min\{\frac{1}{4},\dist(x_0,\mathcal{D}(\phi))\}$. Possibly increasing the constant $C$ here, the above inequality still holds with $\ell$ on the right-hand side replaced by any $\pi\in\mathfrak{l}$. Then combining this estimate with \eqref{E:homo-fine-6} in a standard away using Campanato-style estimates (cf.~\cite[Lemma 4.3]{Wic14}) then shows that, for each $\alpha\in \{1,\dotsc,4\}$, $\phi^\alpha$ is $C^{1,\mu}$ on $K\cup (B^n_{1/4}(0)\setminus \mathcal{D}(\phi))$, and so by the arbitrariness of $K\subset (S_0\setminus S(\phi))\cap B^n_{1/4}(0)$, we see that $\phi^\alpha$ is $C^1$ on $U^\alpha_0\setminus S(\phi)$, and is $C^1$ up-to-the-boundary on $S_0\setminus S(\phi)$. Note that this is $C^1$ regularity \emph{only on half-planes}, due to the possibility of $\pi_y\in\mathfrak{l}$ having its graph a twisted cone. What we will need to do now is show that this regularity extends to $C^1$ regularity on $P_0$.
	
	\textbf{Step 4.} \emph{Concluding the argument.} We already have $C^1$ regularity away from $S(\phi)$ in the directions parallel to $S_0$. Next we will show for full $C^1$ regularity away from $S(\phi)$ by showing continuity of the derivative orthogonal to $S_0$.
	
	For $r>0$ and $y\in S_0$, consider the function $(r,y)\mapsto \phi^1(-r,y) - \phi^3(r,y)$: this is a homogeneous degree one harmonic function on a half-space, which by the continuity of $\phi$ on all of $P_0$ necessarily has $0$ boundary values. This means that $\phi^1(-r,y) - \phi^3(r,y) = ra$ for some constant $a\in P_0^\perp$. But by assumption $\mathcal{D}(\phi)\cap B^n_{1/4}(0)\subsetneq S_0\cap B^n_{1/4}(0)$, and so there is a point $(0,y)\in S_0$ where locally $\phi$ is a harmonic function, and so in particular $\phi^1$ and $\phi^3$ connect across $S_0$ as a smooth harmonic function on some open set. In particular, on this open set we would necessarily have
	$$-\frac{\del}{\del r}\phi^1 = \frac{\del}{\del r}\phi^3$$
	but since we just argued that $\phi^1(-r,y) = \phi^3(r,y) + ar$, differentiating in $r$ gives $a = 0$, i.e.~we have $\phi^1(-r,y) = \phi^3(r,y)$ everywhere. But this implies that the derivatives in the direction orthogonal to $S_0$ coincide across $S_0$, meaning that in fact $\phi^1$ and $\phi^3$ join in a $C^1$ manner across $S_0$. The same argument clearly applies to $\phi^2$ and $\phi^4$ In particular, we have now shown that $\phi$ is $C^1$ on all of $B^n_{1/4}(0)\setminus S(\phi)$, and from now on we will view $\phi$ as being indexed by $I = \{1,2\}$.
	
	We now wish to show that $\phi$ is harmonic away from $S(\phi)$; we already know that it is harmonic away from $\mathcal{D}(\phi)\subset S_0$, and thus we wish to show that it is harmonic on $(S_0\cap \mathcal{D}(\phi))\setminus S(\phi)$. From Lemma \ref{lemma:fine-lambda}, we know that there is a linear function $\ell^1$ such that $\phi^1-\ell^1$ vanishes on $\mathcal{D}(\phi)$. But then since a $C^1$ harmonic function which is harmonic away from its zero set is necessarily harmonic across its zero set as well (see for example \cite{JL05}), and we now know that $\phi^1$ is $C^1$ away from $S(\phi)$, we may apply this result to $\phi^1-\ell^1$, thus seeing that we can extend $\phi^1-\ell^1$ to a harmonic function on $B^n_{1/4}(0)\setminus S(\phi)$. Thus, $\phi^1$ is harmonic away from $S(\phi)$. The same applies to $\phi^2$, and thus we see that $\phi$ is harmonic away from $S(\phi)$.
	
	Finally, as $\dim(S(\phi))\leq n-2$, such a set is removable for a continuous harmonic function, and thus we see that $\phi = \{\phi^\alpha\}_{\alpha=1,2}$ is a pair of harmonic functions on $P_0$. As $\phi$ is homogeneous of degree one, this then implies that both $\phi^1,\phi^2$ are linear, and thus $\phi\in\mathfrak{l}$. This then completes the proof, as it would contradict $\phi$ not agreeing with its dehomogeniser from the beginning of Step 2.
	
	Note that the same argument also works in the case $\psi = \phi - \pi^\phi_{0,1}$, and thus we also establish the corresponding claim for such $\psi$ also. Thus, the induction is complete and the proof of the theorem is completed.
\end{proof}

We are now able to prove the regularity result for arbitrary fine blow-ups in $\mathfrak{b}_{n-1}$. For ease of referencing this statement later on, will include within it the case of any fine blow-up in $\mathfrak{b}$.

\begin{theorem}[Generalised-$C^{1,\alpha}$ Regularity Estimate for Fine Blow-Ups]\label{thm:fine-blow-up-reg}
	There exists $\mu = \mu(n,k,M,\beta)\in (0,1)$ such that the following is true. Suppose that $\phi = \{\phi^\alpha\}_{\alpha\in I}\in\mathfrak{b}$ is a fine blow-up relative to $(\BC_j)_j$. Then, there exists $\pi = \{\pi^\alpha\}_{\alpha\in I}$ arising as a fine blow-up relative to $(\BC_j)_j$ such that:
	\begin{enumerate}
		\item [(1)] $\phi_{0,\rho}\to \pi$ in the sense of fine blow-ups (cf.~Lemma \ref{lemma:fine-translation} and Lemma \ref{lemma:fine-compactness});
		\item [(2)] For each $\alpha\in I$ we have
		$$\sup_{B_1(0)\cap U^\alpha_0}|\pi^\alpha|\leq C\sum_{\alpha^\prime\in I}\|\phi^{\alpha^\prime}\|_{L^2(B_1(0)\cap U^{\alpha^\prime}_0)};$$
		\item [(3)] For any $\rho\in (0,3/16)$ we have
		$$\rho^{-n-2}\sum_{\alpha\in I}\int_{B^n_\rho(0)\cap U^\alpha_0}|\phi^\alpha- (\phi^\alpha(0) + \pi^\alpha)|^2 \leq C\rho^{2\mu}\sum_{\alpha\in I}\int_{B_1(0)\cap U^\alpha_0}|\phi^\alpha|^2.$$
	\end{enumerate}
	Here, $C = C(n,k,M,\beta)\in (0,\infty)$.
\end{theorem}

\begin{proof}
	The proof of this is analogous to a subset of the argument in Steps 2 and 3 of the proof of Theorem \ref{thm:fine-classification}. Indeed, one can prove now prove an analogous form of the reverse Hardt--Simon inequality except now at \emph{any} point $y\in\mathcal{D}(\phi)\cap B_{1/4}(0)$, with the contradiction coming from invoking Theorem \ref{thm:fine-classification}. Given the reverse Hardt--Simon inequality, one may then follow Step 3 of the proof of Theorem \ref{thm:fine-classification} in order to deduce the desired regularity claims.
\end{proof}

\section{Fine Excess Decay \& Regularity}\label{sec:fine-reg}

In this section, we use the regularity of fine blow-ups established in Theorem \ref{thm:fine-blow-up-reg} to complete the proof of the fine $\eps$-regularity theorem, Theorem \ref{thm:fine-reg}. We begin with the main excess improvement lemma.

\begin{lemma}[Fine Excess Improvement]\label{lemma:fine-improvement}
	Let $0<8^{n+1}\theta_{n+1}<\cdots<8\theta_{1}<\theta_0\in (0,1/8)$, $M\in (1,\infty)$. Fix also $\beta \in (0,1)$ and $\tilde{\eps}$, $\tilde{\gamma}\in (0,1)$. Then, there exist $\eps_0$, $\gamma_0\in (0,1)$ depending on $n,k,M,\beta, \{\theta_i\}_{i=0}^{n+1},\tilde{\eps},\tilde{\gamma}$ such that the following is true. Suppose that $V\in \mathcal{V}_\beta$, $\BC\in\CC$ is aligned, and they satisfy:
	\begin{enumerate}
		\item [(a)] $0\in\mathcal{D}(V)$ and $(\w_n 2^n)^{-1}\|V\|(B_2(0)) < \frac{5}{2}$;
		\item [(b)] $\nu(\BC) + \hat{E}_V<\eps_0$;
		\item [(c)] $\hat{E}_V<M\mathscr{E}_{V,\BC}$;
		\item [(d)] $\hat{F}_{V,\BC}<\gamma_0\hat{E}_V$.
	\end{enumerate}
	Then, there exists an aligned cone $\BC^\prime\in \CC$, a rotation $\Gamma\in SO(n+k)$, and $j\in \{0,\dotsc,n+1\}$ such that the following conclusions hold:
	\begin{enumerate}
		\item [(1)] $|\Gamma-\id_{\R^{n+k}}|\leq C\hat{E}_{V}^{-1}\hat{F}_{V,\BC}$.
		\item [(2)] $\nu(\BC,\Gamma_\#\BC^\prime) \leq C\hat{F}_{V,\BC}$;
		\item [(3)] $\hat{F}_{V,\Gamma_\#\BC^\prime}(B_{\theta_j}(0)) \leq C^*_j\theta_j^\mu\hat{F}_{V,\BC}$.
	\end{enumerate}
	Furthermore, there exists an aligned cone $\Dbf\in \CC$ and a rotation $R\in SO(n+k)$ such that:
	\begin{enumerate}
		\item [(4)] $S(\BC)\subset S(\Dbf)\subsetneq P_0$;
		\item [(5)] $\hat{F}_{V,R_\#\Dbf}\leq C\hat{F}_{V,\BC}$;
		\item [(6)] $\nu(\BC,R_\#\Dbf) \leq C\hat{F}_{V,\BC}$;
		\item [(7)] $R_\#^{-1}V$ and $\Dbf$ satisfy Hypothesis $H(\tilde{\eps},\tilde{\gamma},1)$ in $B_1(0)$;\end{enumerate}
	Finally, we also have:
	\begin{enumerate}
		\item [(8)] For any $n$-dimensional plane $P$ with $0\in P$, we have
		$$\nu(\BC,P) \leq C\left(\hat{E}_{\Gamma^{-1}_\#V,P}(B_{\theta_j}(0)) + \hat{E}_{V,\BC}\right).$$
	\end{enumerate}
	Here, $C = C(n,k,M,\beta,\{\theta_j\}_{j=0}^{n+1},\tilde{\eps},\tilde{\gamma})\in (0,\infty)$ and $\mu = \mu(n,k,M,\beta)\in (0,1)$. Furthermore, $C_0 = C_0(n,k,M,\beta)$, and for $j\geq 1$, $C^*_j = C^*_j(n,k,M,\beta,\theta_0,\dotsc,\theta_{j-1})$.
\end{lemma}

\begin{proof}
	We will prove the lemma by proving a few different versions of it with successively weaker hypotheses.
	
	\textbf{Step 1.} The first version of the lemma that we will prove is due to the technicality introduced when we made the stipulation \eqref{E:extra} in Hypothesis $\dagger$ (morally, the case where $\BC$ consists of 4 half-planes instead of two planes needs to be analysed separately). Indeed, for $\delta>0$ let us say that there is a $\delta$\emph{-gap} for $V\in \mathcal{V}_\beta$ and $\BC\in\CC$ if there exists $x_0\in S(\BC)\cap B_1(0)$ such that
	$$B_\delta(x_0)\cap\mathcal{D}(V) = \emptyset.$$
	Fix now $\theta \in (0,1/8)$. We then claim the following:
	
	\textbf{Claim 1:} \emph{There exist $\delta_0, \eps_1,\gamma_1\in (0,1)$ depending only on $n,k,M,\beta,\theta,\tilde{\eps},\tilde{\gamma}$ such that the following is true. If we have $V\in\mathcal{V}_\beta$ and $\BC\in\CC$ satisfy Hypothesis $H(\eps_1,\gamma_1,M)$ in $B_1(0)$, then either there is a $\delta_0$-gap, or all the conclusions of the present lemma hold with just the one scale $\theta$.}
	
	To prove this, as usual we will work by contradiction. Thus, we may take sequences $(V_j)_j\subset\mathcal{V}_\beta$, $(\BC_j)_j\subset\CC$, satisfying Hypothesis $H(\eps_j,\gamma_j,M)$ for some numbers $\eps_j,\gamma_j\downarrow 0$, and moreover we may assume that there are no $\delta_j$-gaps for $V_j$ and $\BC_j$, where $\delta_j\downarrow 0$. We may as usual pass to a subsequence along which $\sfrak(\BC_j) \equiv \sfrak$, and assume (by performing a rotation to $V_j$ and $\BC_j$ that maps $P_0$ to itself) that $S(\BC_j)\equiv S\subset S_0$. Since for every $y\in S\cap B_1$ we have $\mathcal{D}(V_j)\cap B_{\delta_j}(y) \neq\emptyset$, the relatively closed set $\mathcal{D}$ in Hypothesis $\dagger$ obeys $\mathcal{D} = S\cap B_1(0)$.
	
	Let us first verify (1) and (2). Let $(\phi^\alpha)_{\alpha\in I}$ denote a fine blow-up of $(V_j)_j$ relative to $(\BC_j)_j$ in $B_1(0)$. If $\phi\equiv 0$, then necessarily $\hat{E}_{V_j,\BC_j}^{-1}\hat{F}_{V_j,\BC_j}(B_{\theta}(0))\to 0$, and thus for all sufficiently large $j$ (1) and (2) hold with $\BC^\prime = \BC$ and $\Gamma = \id_{\R^{n+k}}$.
	
	Next, suppose that $\phi \in\mathfrak{l}$ and $\phi\not\equiv 0$. Let $\Dbf_j$ and $R_j$ be as in Lemma \ref{lemma:fine-integration}. Then, for all sufficiently large $j$ conclusions (1) and (2) of the present lemma follow from Lemma \ref{lemma:fine-integration} with $\BC_j^\prime := \Dbf_j$ and $\Gamma_j:= R_j^{-1}$.
	
	Finally, suppose that $\phi\not\in \mathfrak{l}$. Then, let $\pi\in\mathfrak{l}$ be as in Theorem \ref{thm:fine-blow-up-reg}, and apply Lemma \ref{lemma:fine-integration} with this $\pi$, again giving cones $\Dbf_j$ and rotations $R_j$ as in Lemma \ref{lemma:fine-integration}. Again, if we set $\Gamma_j := R_j^{-1}$ and $\BC_j^\prime:= \tilde{\BC}_j$, then conclusions (1) and (2) are immediate from Lemma \ref{lemma:fine-integration}, and conclusion (3) follows from combining the $L^2$ convergence to the blow-up given in Lemma \ref{lemma:fine-integration} with the decay estimate in Theorem \ref{thm:fine-blow-up-reg} (with $\rho = \theta$). This proves (1), (2), and (3) in this case. (Notice also that here we in fact get $\hat{F}_{V,\Gamma_\#\BC^\prime}(B_{\theta}(0))\leq C_1\theta^\mu \hat{E}_{V,\BC}$ in (3).)
	
	Conclusions $(4)-(7)$ all hold in this case for sufficiently large $j$, with $R_j := \id_{\R^{n+k}}$ and $\Dbf_j = \BC_j$, by our initial assumptions.
	
	Finally, let us first prove (8). Fix any $n$-dimensional plane $P$ with $0\in P$. Let $p:P_0\to P_0^\perp$ be the linear function for which $P = \graph(p)$. Observe the elementary fact that there is a constant $C = C(n,k,\theta)\in (0,\infty)$ such that for any subspace $\tilde{S}$ with $S\subset\tilde{S}\subsetneq P_0$, we have
	\begin{equation}\label{E:fine-imp-1}
		\nu(\BC_j,P)^2 \leq \begin{cases}
			C\sum^4_{\alpha=1}\int_{B_{\theta}(0)\cap U^\alpha_0\cap \{|x^n|>\theta/16\}}|h^\alpha_j-p|^2\, \ext x & \text{if }\sfrak=n-1;\\
			C\sum_{\alpha=1,2}\int_{B_{\theta}(0)\cap\{|x^{\perp_{\tilde{S}}}|>\theta/16\}}|p^\alpha_j-p|^2\, \ext x & \text{if }\sfrak<n-1.
		\end{cases}
	\end{equation}
	Let us now work with the $\sfrak=n-1$ notation, as the argument in both cases is analogous. Since Hypothesis $H(\eps_j,\gamma_j,M)$ are satisfied,
	we have functions $\{u^\alpha_j\}_{\alpha=1}^4$ representing $V_j$ over $\BC_j$ as in Theorem \ref{thm:graphical-rep}. Therefore, the triangle inequality gives
	\begin{align*}
		\sum^4_{\alpha=1}&\int_{B_{\theta}(0)\cap U^\alpha_0\cap \{|x^n|>\theta/16\}}|h^\alpha_j-p|^2\\
		& \leq C\sum^4_{\alpha=1}\int_{B_{\theta}(0)\cap U^\alpha_0\cap \{|x^n|>\theta/16\}}|u_j^\alpha|^2 + C\sum^4_{\alpha=1}\int_{B_{\theta}(0)\cap U^\alpha_0\cap \{|x^n|>\theta/16\}}|h^\alpha_j+u^\alpha_j - p|^2\\
		& \leq C\hat{E}_{V_j,\BC_j}^2 + C\hat{E}_{V_j,P}^2(B_{\theta}(0)).
	\end{align*}
	Combining this with \eqref{E:fine-imp-1}, we get (8). This finishes the proof of Claim 1.

	\textbf{Step 2.} Now we remove the dichotomy concerning the $\delta_0$-gap from Claim 1. Indeed, we claim the following. Fix now $\theta_0,\theta_1$ obeying $0<8\theta_0<\theta_1<1/8$.
	
	\textbf{Claim 2:} \emph{There exist $\eps_1,\gamma_1\in (0,1)$ depending only on $n,k,M,\beta,\theta_0,\theta_1,\tilde{\eps},\tilde{\gamma}$ such that the following is true. Suppose that $V\in\mathcal{V}_\beta$ and $\BC\in \CC$ satisfy Hypothesis $H(\eps_1,\gamma_1,M)$ in $B_1(0)$. Then, the conclusions of the present lemma hold with $j=0$ or $j=1$ therein.}
	
	To prove this, we begin exactly as in Step 1: if $V_j$ and $\BC_j$ have no $\delta_0$-gaps, where here $\delta_0 = \delta_0(n,k,M,\beta, \theta_1, \tilde{\eps}, \tilde{\gamma})$ is as in Claim 1 (in particular depending on $\theta_1$ but not $\theta_0$), then the conclusion follows from Claim 1. Thus, we may now suppose that there is a $\delta_0$-gap. If $\BC_j\in \mathcal{P}$ for infinitely many $j$ then, after passing to a subsequence, we can still check Hypothesis $\dagger$ and proceed analogously to the argument in Step 1. Thus, the last remaining case to consider is when for all $j$ sufficiently large we have $\BC_j\in \CC_{n-1}\setminus\mathcal{P}_{n-1}$; the complication of \eqref{E:extra} in Hypothesis $\dagger$ means that we need to make some adjustments. Let us pass to a subsequence so that $\BC_j\in \CC_{n-1}\setminus\mathcal{P}_{n-1}$.
	
	For each $j$, pick a new cone $\overline{\BC}_j\in\mathcal{P}_{n-1}$ with $S(\overline{\BC}_j) = S(\BC_j)\equiv S_0$ satisfying
	$$\hat{E}_{V_j,\overline{\BC}_j}\leq\frac{3}{2}\inf_{\overline{\BC}}\hat{E}_{V_j,\overline{\BC}}$$
	where the infimum is taken over all $\overline{\BC}\in \mathcal{P}_{n-1}$ with $S(\overline{\BC}) = S_0$. We claim that there is a constant $C_0 = C_0(n,k,\delta_0)\in (0,\infty)$ for which
	\begin{equation}\label{E:fine-imp-2}
		\hat{E}_{V_j,\overline{\BC}_j}\leq C_0\hat{E}_{V_j,\BC_j}.
	\end{equation}
	This follows from a contradiction argument analogous to \cite[Lemma 7.1]{BK17}; we sketch the details for completeness. Indeed, if \eqref{E:fine-imp-2} did not hold, we would instead have $\hat{E}_{V_j,\overline{\BC}_j}^{-1}\hat{E}_{V_j,\BC_j} \to 0$. We also know by assumption that for each $j$ there is $x_j\in S_0\cap B_1(0)$ for which
	$$B_{\delta_0}(x_j) \cap \mathcal{D}(V_j) = \emptyset.$$
	Pass to a subsequence for which $x_j\to x_0\in S_0\cap \overline{B}_1(0)$; in particular, we then have $\mathcal{D}(V_j)\cap B_{\delta_0/2}(x_0)\cap B_1(0) = \emptyset$ for all sufficiently large $j$. We also have $\nu(\BC_j,\overline{\BC}_j) \leq C\hat{E}_{V_j,\overline{\BC}_j}$. 
	
	Since Hypothesis $H(\eps_j,\gamma_j,M)$ holds for $V_j$ and $\BC_j$, we may apply Theorem \ref{thm:graphical-rep} to see that all points of density $\geq 2$ in $V_j$ must accumulate along $S_0$, which in turn implies (arguing in the same way as in the proof of Theorem \ref{thm:graphical-rep} in the case $\sfrak(\BC) = n-1$) that the conclusions of Theorem \ref{thm:graphical-rep} hold for $V_j$ and $\overline{\BC}_j$, and so in particular we may produce a blow-up of $V_j$ relative to $\overline{\BC}_j$; call this $v$. By analogous arguments to those seen in the proof of Lemma \ref{lemma:fine-integration}, we have strong $L^2_{\text{loc}}(B_1(0))$ convergence to this blow-up, and in fact a simple argument analogous to what we have done before shows that we can utilise Proposition \ref{prop:full-L2} to show that in fact the convergence to the blow-up is in fact in $L^2(B_1(0))$. By the bound $\nu(\BC_j,\overline{\BC}_j)\leq C\hat{E}_{V_j,\overline{\BC}_j}$, we may blow-up $\BC_j$ relative to $\overline{\BC}_j$, giving rise to a $2$-valued function $w$ which is linear on $B^n_1(0)\cap \{|x^n|>0\}$. By our contradiction assumption, we in fact have $w = v$. Notice that since
	$$\hat{E}_{V_j,\overline{\BC}_j} \leq C\hat{E}_{V_j,\BC_j} + C\nu(\BC_j,\overline{\BC}_j)$$
	and since $\hat{E}_{V_j,\overline{\BC}_j}^{-1}\hat{E}_{V_j,\BC_j}\to 0$, we can absorb the first term on the right-hand side into the left-hand side to see that $\hat{E}_{V_j,\overline{\BC}_j}^{-1}\nu(\BC_j,\overline{\BC}_j)\geq c>0$ for all $j$, which implies that $w\not\equiv 0$. However, since $\mathcal{D}(V_j)<2$ in $B_{\delta_0/2}(x_0)\cap B_1(0)$, in some fixed ball centred on $S_0$, we actually know that the blow-up $v$ of $V_j$ off $\overline{\BC}_j$ must consist of two smooth harmonic functions in this ball which, as $v = w$ is linear, implies that $v$ must be two linear functions on all of $B_1(0)$, i.e.~$\graph(v)$ is the union of two planes (which necessarily meet along $S_0$ from the form of $w$). But then arguing analogously to the proof of Lemma \ref{lemma:fine-integration}, utilising the strong $L^2$ convergence we have on all of $B_1(0)$ to the blow-ups, we see that if we form a new cone $\overline{\BC}_j^\prime$ by modifying $\overline{\BC}_j$ based on the blow-up $v$ just as in Lemma \ref{lemma:fine-integration}, we have that $\overline{\BC}_j^\prime\in \mathcal{P}_{n-1}$ has $S(\overline{\BC}_j^\prime) = S_0$ yet
	$$\hat{E}_{V_j,\overline{\BC}_j}^{-1}\hat{E}_{V_j,\overline{\BC}_j^\prime} \to 0$$
	which contradicts the choice of $\overline{\BC}_j$ in the infimum. This contradiction completes the proof of \eqref{E:fine-imp-2}.
	
	Now given \eqref{E:fine-imp-2}, it suffices to work with the cones $(\overline{\BC}_j)_j\in\mathcal{P}$ in order to show (1), (2), and (3) (notice that, in Step 1, the right-hand side of (3) was in terms of $\hat{E}_{V,\BC}$ instead of $\hat{F}_{V,\BC}$, and so if we can prove this with $\overline{\BC}_j$ in place of $\BC_j$ the result for $\BC_j$ follows from \eqref{E:fine-imp-2}). For this, we can simply follow the arguments in Step 1. (Notice that, due to \eqref{E:fine-imp-2}, $V_j$ and $\overline{\BC}_j$ can be shown to still satisfy Hypothesis $\dagger$ here, so the argument does pass through.) Notice however that the constant $C_0$ from \eqref{E:fine-imp-2} depends on $\delta_0$ and thus on $\theta_1$, and so when applying Step 1 here we are applying Step 1 with the choice of $\theta$ therein being $\theta_0$; the constants then will also depend on $\theta_1$. This proves (1), (2), and (3). For $(4)-(7)$, we can again just take $\Dbf_j := \BC_j$ and $R_j := \id_{\R^{n+k}}$, and (8) follows in an identical way to that seen in Step 1. This completes the proof of Claim 2.
	
	\textbf{Step 3.} By Claim 2, we have now eliminated the need to discuss $\delta$-gaps. The next step is to replace Hypothesis (A4) with the weaker hypothesis (d) of the present lemma. For this, we fix a set of scales $\{\theta_0,\theta_1,\dotsc,\theta_{n}\}$ obeying $0<8^{n+1}\theta_{n+1}<\cdots<8\theta_1<\theta_0<1/8$ (i.e.~as in the statement of the lemma). We then prove the following.
	
	\textbf{Claim 3:} \emph{For each $\sfrak\in \{0,1,\dotsc,n-1\}$, there exist constants $\eps_2^{(\sfrak)}, \gamma_2^{(\sfrak)}$ depending on $n,k,M,\beta$, $\theta_{0},\dotsc,\theta_{n+1-\sfrak},\tilde{\eps},\tilde{\gamma},\sfrak$, such that the following is true. If we have $V\in\mathcal{V}_\beta$ and $\BC\in \CC_{\sfrak}$ satisfy assumptions (a$)-($d) of the present lemma with $\eps_2^{(\sfrak)}$ and $\gamma_2^{(\sfrak)}$ in place of $\eps_0,\gamma_0$, respectively, then the conclusions of the lemma hold with $\theta = \theta_p$ for some $p\in \{0,\dotsc,n+1-\sfrak\}$, where the constants $C$ and $C^*_j$ additionally depend on $\sfrak$ (here, we will denote them by $C_\sfrak$ and $C^{(j,\sfrak)}$, respectively).} 
	
	We will prove this claim by (downwards) induction on $\sfrak$. The base case is when $\sfrak=n-1$. As in Step 1, consider appropriate sequences $(V_j)_j\subset\mathcal{V}_\beta$, $(\BC_j)_j\subset \CC$ with $\sfrak(\BC_j)\equiv \sfrak = n-1$ satisfying the hypotheses with $\eps_j$ and $\gamma_j$ in place of $\eps$ and $\gamma$, where $\eps_j,\gamma_j\downarrow 0$. Let $\gamma^*_1$ be the constant $\gamma_1$ from Claim 2. Since $\sfrak = n-1$, hypothesis (d) is identical to Hypothesis (A4) (this is using (c) also), and thus for all sufficiently large $j$ we can apply Claim 2 to get the required conclusions when $\sfrak=n-1$.
	
	So now suppose that $\sfrak < n-1$. If it is the case that there is a subsequence $\{j^\prime\}\subset\{j\}$ along which we have
	\begin{equation}\label{E:fine-imp-3}
		\frac{\hat{F}_{V_{j^\prime},\BC_{j^\prime}}}{\mathscr{F}^*_{V_{j^\prime},\BC_{j^\prime}}}<\gamma_1^*,
	\end{equation}
	then we can again verify (A4) and directly apply Claim 2 to get the desired conclusion (with decay at either scale $\theta_0$ or $\theta_1$). If there is no subsequence along which \eqref{E:fine-imp-3} holds, then we need to pick an optimal coarser cone. Indeed, in this case pass to a subsequence along which the negation of \eqref{E:fine-imp-3} holds. For each $j$, pick $\widetilde{\Dbf}_j\in\mathcal{P}$ coarser to $\BC_j$ and with $\hat{F}_{V_j,\widetilde{\Dbf}_j}\leq \frac{3}{2}\mathscr{F}^*_{V_j,\BC_j}$. Pass to a subsequence along which $\sfrak(\widetilde{\Dbf}_j) \equiv \sfrak^\prime>\sfrak$. Now let $P_j$ be a plane coarser to $\widetilde{\Dbf}_j$ for which
	\begin{equation*}
		\hat{E}_{V_j,P_j} = \mathscr{E}_{V_j,\widetilde{\Dbf}_j}.
	\end{equation*}
	Now let $\widetilde{R}_j\in SO(n+k)$ be a rotation minimising $|\widetilde{R}-\id_{\R^{n+k}}|$ over all $\widetilde{R}\in SO(n+k)$ subject to the condition that $\widetilde{R}(P_j) = P_0$. We claim that the hypotheses of the present claim are satisfied with $(\tilde{R}_j)_\#V_j$ and $(\tilde{R}_j)_\#\widetilde{\Dbf}_j$ in place of $V_j$ and $\BC_j$ (up to also replacing $\beta$ by $\beta/2$, say). Indeed, hypotheses (a) and (b) are simple to verify. For (c), note that $0\in S(\BC_j)\subset S((\widetilde{R}_j)_\#\widetilde{\Dbf}_j)\subsetneq P_0$, and
	$$\hat{E}_{(\widetilde{R}_j)_\#V_j} = \hat{E}_{V_j,P_j} = \mathscr{E}_{V_j,\widetilde{\Dbf}_j} = \mathscr{E}_{(\widetilde{R}_j)_\#V_j,(\widetilde{R}_j)_\#\widetilde{\Dbf}_j}$$
	which therefore shows that (c) holds with $M=1$ therein. The remaining thing to check if (d). For this, using that $\widetilde{R}_j$ is an isometry of $B_1(0)$ which fixes $S(\widetilde{\Dbf}_j)$, as well as the negation of \eqref{E:fine-imp-3}, we have
	\begin{equation}\label{E:fine-imp-5}
		\hat{F}_{(\widetilde{R}_j)_\#V_j,(\widetilde{R}_j)_\#\widetilde{\Dbf}_j} = \hat{F}_{V_j,\widetilde{\Dbf}_j} \leq \frac{3}{2}\mathscr{F}^*_{V_j,\BC_j} \leq \frac{3}{2\gamma^*_1}\hat{F}_{V_j,\BC_j} < \frac{3\gamma_j}{2\gamma^*_1}\hat{E}_{V_j}
	\end{equation}
	where in the last inequality we have used (d) for the original sequence $V_j$ and $\BC_j$. We now want to replace $\hat{E}_{V_j}$ on the right-hand side by $\hat{E}_{(\widetilde{R}_j)_\#V_j}$. From (c) for $V_j$ and $\BC_j$, the definition of $\mathscr{E}_{V_j,\BC_j}$ (using that $P_j\supset S(\BC_j)$) and the fact that $\widetilde{R}_j(P_j) = P_0$, we have
	$$\hat{E}_V \leq M\mathscr{E}_{V_j,\BC_j} \leq M\hat{E}_{V_j,P_j} = M \hat{E}_{(\widetilde{R}_j)_\#V_j}.$$
	Therefore, we see that \eqref{E:fine-imp-5} gives
	\begin{equation}\label{E:fine-imp-6}
		\hat{F}_{(\widetilde{R}_j)_\#V_j,(\widetilde{R}_j)_\#\widetilde{\Dbf}_j} \leq \frac{3M}{2\gamma^*_1}\cdot\gamma_j\cdot\hat{E}_{(\widetilde{R}_j)_\#V_j}.
	\end{equation}
	This verifies the appropriate from of assumption (d) needed. Now by appealing to the induction hypothesis, there exists $\eps^{(\sfrak^\prime)}_2$ and $\gamma^{(\sfrak^\prime)}_2$ depending on $n,k,M,\beta,\tilde{\eps},\tilde{\gamma},\theta_{0},\theta_1,\dotsc,\theta_{n+2-\sfrak^\prime}$ for which Claim 3 holds with $\sfrak^\prime$ in place of $\sfrak.$ Clearly by \eqref{E:fine-imp-6} for all sufficiently large $j$ we have $\frac{3M}{2\gamma^*_1}\cdot\gamma_j < \gamma_2^{(\sfrak^\prime)}$, and thus the hypotheses of Claim 3 are satisfied with $(\widetilde{R}_j)_\#V_j$ and $(\widetilde{R}_j)_\#\widetilde{\Dbf}_j$ in place of $V_j$ and $\BC_j$. Thus we can apply the present claim inductively to get a rotation $\Gamma_j$, a new aligned cone $\BC^\prime_j$, and $p\in \{2,3,\dotsc,n+2-\sfrak^\prime\}$, for which
	\begin{equation}\label{E:fine-imp-7}
		\hat{F}_{(\widetilde{R}_j)_\#V_j,(\Gamma_j)_\#\BC_j^\prime}(B_{\theta_p}(0)) \leq C^{(p,\sfrak^\prime)}\theta^\mu_p \hat{F}_{(\widetilde{R}_j)_\#V_j,(\widetilde{R}_j)_\#\widetilde{\Dbf}_j};
	\end{equation}
	\begin{equation}\label{E:fine-imp-8}
		\nu((\widetilde{R}_j)_\#\widetilde{\Dbf}_j,(\Gamma_j)_\#\BC_j^\prime) \leq C^{(p,\sfrak^\prime)}\hat{F}_{(\widetilde{R}_j)_\#V_j,(\widetilde{R}_j)_\#\widetilde{\Dbf}_j};
	\end{equation}
	\begin{equation}\label{E:fine-imp-8.5}
		|\Gamma_j - \id_{\R^{n+k}}|\leq C^{(p,\sfrak^\prime)}\hat{E}_{(\widetilde{R}_j)_\#V_j}^{-1}\hat{E}_{(\widetilde{R}_j)_\#V_j,(\widetilde{R}_j)_\#\widetilde{\Dbf}_j}.
	\end{equation}
	By combining \eqref{E:fine-imp-7} with \eqref{E:fine-imp-5}, we get conclusion (3) of the lemma holds for sufficiently large $j$ with $\tilde{R}_j^{-1}\circ\Gamma_j$ in place of $\Gamma$. Also, conclusion (1) of the lemma follows from \eqref{E:fine-imp-8.5} along with \eqref{E:fine-imp-5} and the previously established bound $\hat{E}_V \leq M\hat{E}_{(\widetilde{R}_j)_\#V_j}$. To prove (2), we need to verify that
	$$\nu(\BC_j,(\widetilde{R}_j^{-1}\circ\Gamma_j)_\#\BC_j^\prime) \leq C\hat{F}_{V_j,\BC_j}.$$
	As $\BC_j$, $\BC_j^\prime$, and $\widetilde{\Dbf}_j$ are all unions of planes, controlling the Hausdorff distance between any pair in terms the $L^2$ norm of the difference between the linear functions defining the planes, and using the triangle inequality, we get
	$$\nu(\BC_j,(\widetilde{R}_j^{-1}\circ\Gamma_j)_\#\BC_j^\prime) \leq C\nu(\widetilde{\Dbf}_j,(\widetilde{R}_j^{-1}\circ\Gamma_j)_\#\BC_j^\prime) + C\nu(\widetilde{\Dbf}_j,\BC_j).$$
	We now estimate each of these two terms separately. Indeed, using the fact that $\widetilde{R}_j$ is an isometry of $B_1(0)$, as well as \eqref{E:fine-imp-8} and \eqref{E:fine-imp-5}, we get
	\begin{align*}
	\nu(\widetilde{\Dbf}_j,(\widetilde{R}_j^{-1}\circ\Gamma_j)_\#\BC_j^\prime) & = \nu((\widetilde{R}_j)_\#\widetilde{\Dbf}_j,(\Gamma_j)_\#\BC_j^\prime)\\
	& \leq C^{(p,\sfrak^\prime)}\hat{F}_{(\widetilde{R}_j)_\#V_j,(\widetilde{R}_j)_\#\widetilde{\Dbf}_j}\\
	& \leq C^{(p,\sfrak^\prime)}\cdot\frac{3}{2\gamma^*_1}\hat{F}_{V_j,\BC_j}.
	\end{align*}
	This controls one term. To control the other, using a graphical representation in a similar way to the computation following \eqref{E:fine-imp-1} we have
	\begin{equation}\label{E:fine-imp-9}
	\nu(\widetilde{\Dbf}_j,\BC_j) \leq C\hat{F}_{V_j,\widetilde{\Dbf}_j} + C\hat{F}_{V_j,\BC_j}
	\end{equation}
	which again by \eqref{E:fine-imp-5} to controlled by $C\hat{F}_{V_j,\BC_j}$ for all sufficiently large $j$. (It should be noted that in order to get the appropriate graphical representation, by arguing as in Step 1 of the inductive proof of Theorem \ref{thm:graphical-rep} to choose an appropriate optimal cone for the two-sided fine excess, and using hypothesis (c) of the present lemma, in the present setting we always know that points of density $\geq 2$ in $V_j$ must accumulate around $S_0$, and so we always have a graphical representation of $V_j$ over the cones away from $S_0$, which is sufficient for the above bound.) Combining the above we therefore prove (2).
	
	Next we establish conclusions $(4)-(7)$. By the inductive application of the present claim above, we know that there exists an aligned cone $\Dbf_j$ and corresponding rotation $\Delta_j$ such that $(4)-(7)$ are satisfied with $(\widetilde{R}_j)_\#V_j$ and $(\widetilde{R}_j)_\#\widetilde{\Dbf}_j$ in place of $V_j$ and $\BC_j$, i.e.~we have
	\begin{enumerate}
		\item [$(4)^\prime$] $0\in S((\widetilde{R}_j)_\#\widetilde{\Dbf}_j)\subset S(\Dbf_j) \subsetneq P_0$;
		\item [$(5)^\prime$] $\hat{F}_{(\widetilde{R}_j)_\#V_j,(\Delta_j)_\#\Dbf_j}\leq C\hat{F}_{(\widetilde{R}_j)_\#V_j,(\widetilde{R}_j)_\#\widetilde{\Dbf}_j}$;
		\item [$(6)^\prime$] $\nu((\widetilde{R}_j)_\#\widetilde{\Dbf}_j, (\Delta_j)_\#\Dbf_j) \leq C\hat{F}_{(\widetilde{R}_j)_\#V_j,(\widetilde{R}_j)_\#\widetilde{\Dbf}_j}$;
		\item [$(7)^\prime$] $(\Delta_j)^{-1}_\#(\widetilde{R}_j)_\#V_j$ and $\Dbf_j$ satisfy Hypothesis $H(\tilde{\eps},\tilde{\gamma},1)$ in $B_1(0)$.
	\end{enumerate}
	We now claim that $(4)-(7)$ follow for $V_j$ and $\BC_j$ with $R_j:= \widetilde{R}_j^{-1}\circ\Delta_j$ and $\Dbf_j$. Since $S(\BC_j)\subset S((\widetilde{R}_j)_\#\widetilde{\Dbf}_j)$, we immediately have (4) from $(4)^\prime$. Now using \eqref{E:fine-imp-5} with $(5)^\prime$, $(6)^\prime$, and $(7)^\prime$, we deduce
	\begin{enumerate}
		\item [$(5)^{\prime\prime}$] $\hat{F}_{V_j,(R_j)_\#\Dbf_j} \leq C\hat{F}_{V_j,\BC_j}$;
		\item [$(6)^{\prime\prime}$] $\nu(\widetilde{\Dbf}_j,(R_j)_\#\Dbf_j) \leq C\hat{F}_{V_j,\BC_j}$;
		\item [$(7)^{\prime\prime}$] $(R_j^{-1})_\#V_j$ and $\Dbf_j$ satisfy Hypothesis $H(\tilde{\eps},\tilde{\gamma},1)$ in $B_1(0)$.
	\end{enumerate}
	Clearly $(5)^{\prime\prime}$ and $(7)^{\prime\prime}$ are exactly the desired forms of (5) and (7), respectively. For (6), another argument based on the linearity of the cones and the triangle inequality gives
	$$\nu(\BC_j,(R_j)_\#\Dbf_j) \leq C\nu(\BC_j,\widetilde{\Dbf}_j) + C\nu(\widetilde{\Dbf}_j,(R_j)_\#\Dbf_j)$$
	which combined with $(6)^{\prime\prime}$ and \eqref{E:fine-imp-9} gives (6). This therefore establishes $(4)-(7)$.
	
	Finally, one can verify (8) in an analogous way to how it was verified in Step 1 following \eqref{E:fine-imp-1}, again noting that our the present assumptions all points of density $\geq 2$ must accumulate around $S_0$, and thus away from $S_0$ we have appropriate graphical representations for the argument to hold. Thus, this completes the proof of Claim 3.
	
	\textbf{Step 4.} Now the version of the lemma in the actual statement can be deduced from Claim 3. Indeed, notice that the constants
	$$\eps_0:= \min_{\sfrak}\eps^{(\sfrak)}_2, \qquad \gamma_0:=\min_{\sfrak}\gamma^{(\sfrak)}_2, \qquad C:=\max_{\sfrak}C_{\sfrak}, \qquad \text{and}\qquad C_p^*:= \max_{\sfrak}C^{(p,\sfrak)}$$
	from Claim 3 have the correct dependencies, and thus the lemma follows with these choices from Claim 3.
\end{proof}

\begin{remark}\label{remark:after-fine-improvement}
	We record one additional consequence of Lemma \ref{lemma:fine-improvement}. For $z\in\mathcal{D}(V)\cap B_{1/2}(0)$, we have $R^{-1}(z)\in \mathcal{D}(R^{-1}_\#V)$. Then using (7) of Lemma \ref{lemma:fine-improvement} (provided $\tilde{\eps}$ and $\tilde{\gamma}$ are appropriately small depending on $n,k,M,\beta$) to invoke \eqref{E:cor-2} of Corollary \ref{cor:L2-corollary}, and then using (5), we have
	\begin{equation}\label{E:fine-imp-con-1}
	\nu(\Dbf,(\tau_{R^{-1}(z)})_\#\Dbf) \leq C\hat{F}_{V,\BC}.
	\end{equation}
Then, the triangle inequality gives (again, using the fact that for unions of planes Hausdorff distance can be controlled in terms of the $L^2$ norms of the differences of the linear functions defining the planes in the cones)
\begin{align}
	\nonumber\nu(\BC,(\tau_z)_\#\BC) & \leq C\nu(\BC,R_\#\Dbf) + C\nu(R_\#\Dbf,(\tau_z)_\#R_\#\Dbf) + C\nu((\tau_z)_\#R_\#\Dbf,(\tau_z)_\#\BC)\\
	& = C\nu(\BC,R_\#\Dbf) + C\nu(\Dbf,(\tau_{R^{-1}(z)})_\#\Dbf) + C\nu(R_\#\Dbf,\BC)\nonumber\\
	& \leq C\hat{F}_{V,\BC}\label{E:fine-imp-con-2}
\end{align}
where in the second line we have used that for $R$ a rotation matrix and $w$ any vector, the map $R^{-1}\circ \tau_w\circ R$ is simply $\tau_{R^{-1}(w)}$, and in the final line we have used \eqref{E:fine-imp-con-1} and (6) of Lemma \ref{lemma:fine-improvement}.
\end{remark}

Now we iterate Lemma \ref{lemma:fine-improvement} to deduce a uniform decay estimate at a point of density $\geq 2$.

\begin{corollary}\label{cor:fine-decay}
	Fix $M\in [1,\infty)$ and $\beta\in (0,1)$. Then, there exists $\eps_0 = \eps_0(n,k,M,\beta)\in (0,1)$ and $\gamma_0 = \gamma_0(n,k,M,\beta)\in (0,1)$ such that the following is true. Suppose that $V\in\mathcal{V}_\beta$, $\BC\in\CC$ is aligned, and they satisfy:
	\begin{enumerate}
		\item [(A)] $\Theta_V(0)\geq 2$ and $(\w_n2^n)^{-1}\|V\|(B_2(0))\leq \frac{5}{2}$;
		\item [(B)] $\nu(\BC) + \hat{E}_V<\eps$;
		\item [(C)] $\hat{E}_V\leq M\mathscr{E}_{V,\BC}$;
		\item [(D)] $\hat{F}_{V,\BC}<\gamma_0\hat{E}_V$.
	\end{enumerate}
	Then, there exists $\theta = \theta(n,k,M,\beta) \in (0,1/8^{n+3})$, an aligned cone $\BC^\prime\in\CC$, and a rotation $\Gamma\in SO(n+k)$ such that the following conclusions hold:
	\begin{enumerate}
		\item [(1)] $\nu(\BC,\Gamma_\#\BC^\prime) \leq C\hat{F}_{V,\BC}$;
		\item [(2)] $\hat{E}_{V,\Gamma_\#\BC^\prime}(C_\rho(0)) \leq C\rho^\mu \hat{F}_{V,\BC}$ for every $\rho\in (0,\theta)$, and $\Gamma_\#\BC^\prime$ is the unique tangent cone to $V$ at $0$;
		\item [(3)] $\nu(\Gamma_\#\BC^\prime, (\tau_z)_\#\Gamma_\#\BC^\prime) \leq C\rho^\mu \hat{F}_{V,\BC}$ for every $z\in B_{\rho/2}(0)\cap \mathcal{D}(V)$ and every $\rho\in (0,\theta/2)$.
	\end{enumerate}
\end{corollary}

\begin{proof}
	Let $\mu$ be as in Lemma \ref{lemma:fine-improvement}, and let $(C_j^*)_{j=0}^{n+1}$ be the constants from Lemma \ref{lemma:fine-improvement}; in particular, $C_0^* = C_0^*(n,k,M,\beta)$. So, we can choose $\theta_0 = \theta_0(n,k,M,\beta)\in (0,1/8)$ such that
	$$C_0^*\theta_0^\mu < 1/2.$$
	Now, inductively choose $\theta_j\in (0,\theta_{j-1}/8)$ such that $C_j^*\theta_j^\mu < 1/2$; by induction, we have that $C_j^* = C_j^*(n,k,M,\beta)$ and so $\theta_j = \theta_j(n,k,M,\beta)$. Thus, the choices of $\theta_j$ only depend on $n,k,M,\beta$, and so the corresponding choices of $\eps_0,\gamma_0$ in Lemma \ref{lemma:fine-improvement} only depend on $n,k,M,\beta,\tilde{\eps},\tilde{\gamma}$. Furthermore, we choose $\tilde{\eps}$ and $\tilde{\gamma}$ sufficiently small depending on $n,k,M,\beta$ so that one may apply Theorem's \ref{thm:graphical-rep}, \ref{thm:L2-estimates}, \ref{thm:shift-estimate}, and Corollary \ref{cor:L2-corollary} when Hypothesis $H(\tilde{\eps},\tilde{\gamma},1)$ are satisfied. As such, all of our constants now only depend on $n,k,M,\beta$.
	
	We then claim the following:
	
	\textbf{Claim:} \emph{There exists a sequence of aligned cones $(\BC_\ell)_{\ell=0}^\infty\subset \CC$ with $\BC_0 = \BC$ and a sequence of rotations $(\Gamma_\ell)_{\ell=0}^\infty\subset SO(n+k)$ with $\Gamma_0 = \id_{\R^{n+k}}$ such that for $\ell\geq 1$, the following conclusions hold:}
	\begin{equation}\label{E:fd-1}
		\nu(\BC_\ell,(\Gamma_\ell)_\#\BC_\ell) \leq C\hat{F}_{V,(\Gamma_{\ell-1})_\#\BC_{\ell-1}}(B_{\sigma_{\ell-1}}(0));
	\end{equation}
	\begin{equation}\label{E:fd-2}
		\nu((\Gamma_\ell)_\#\BC_\ell, (\Gamma_{\ell-1})_\#\BC_{\ell-1}) \leq C\hat{F}_{V,(\Gamma_{\ell-1})_\#\BC_{\ell-1}}(B_{\sigma_{\ell-1}}(0));
	\end{equation}
	\begin{equation}\label{E:fd-3}
		\hat{F}_{V,(\Gamma_\ell)_\#\BC_\ell}(B_{\sigma_\ell}(0)) \leq \frac{1}{2}\hat{F}_{V,(\Gamma_{\ell-1})_\#\BC_{\ell-1}}(B_{\sigma_{\ell-1}}(0));
	\end{equation}
	\begin{equation}\label{E:fd-4}
		\hat{F}_{V,(\Gamma_\ell)_\#\BC_\ell}(B_{\sigma_\ell}(0)) \leq 2^{-\ell}\hat{F}_{V,\BC};
	\end{equation}
	\emph{and moreover, for any plane $P$ with $0\in P$ we have}
	\begin{equation}\label{E:fd-5}
		\nu(\BC_{\ell-1},P) \leq C\hat{E}_{(\Gamma_{\ell}^{-1})_\#V,P}(B_{\sigma_{\ell}}(0)) + C\hat{E}_{V,(\Gamma_{\ell-1})_\#\BC_{\ell-1}}(B_{\sigma_{\ell-1}}(0));
	\end{equation}
	\emph{and for every $z\in B_{\sigma_{\ell-1}/2}(0)\cap \mathcal{D}(V)$ we have}
	\begin{equation}\label{E:fd-6}
		\nu((\Gamma_{\ell-1})_\#\BC_{\ell-1}, (\tau_{z})_\#(\Gamma_{\ell-1})_\#\BC_{\ell-1}) \leq C\hat{F}_{V,(\Gamma_{\ell-1})_\#\BC_{\ell-1}}(B_{\sigma_{\ell-1}}(0)).
	\end{equation}
	Here, $\sigma_\ell = \sigma_{\ell-1} \theta_{j_\ell}$ for some choice of $j_\ell\in \{0,\dotsc,n+1\}$, with $\sigma_0 = 1$, and $C = C(n,k,M,\beta)\in (0,\infty)$ is independent of $\ell$. (We will also have $|\Gamma_\ell-\id_{\R^{n+k}}|\leq C\hat{E}_{V,\BC}^{-1}\hat{F}_{V,\BC}$ for all $\ell$.)
	
	We prove this claim by induction on $\ell$. The $\ell=1$ case of \eqref{E:fd-1}--\eqref{E:fd-6} follows immediately from Lemma \ref{lemma:fine-improvement} applied to $V$ and $\BC$, which we can do from our assumptions, as well as Remark \ref{remark:after-fine-improvement} following the proof of Lemma \ref{lemma:fine-improvement}; notice that \eqref{E:fd-1} follows from the proof of Lemma \ref{lemma:fine-improvement} arguing analogously to that seen in \eqref{E:b-4-5}.
	
	Now suppose that we have produced cones $\{\BC_{\ell^\prime}\}_{\ell^\prime=1}^\ell$ and rotations $\{\Gamma_{\ell^\prime}\}_{\ell^\prime=1}^\ell$ satisfying \eqref{E:fd-1} -- \eqref{E:fd-6} and we wish to find $\BC_{\ell+1}$ and $\Gamma_{\ell+1}$. We will do this by verifying the hypotheses of Lemma \ref{lemma:fine-improvement} with $(\eta_{0,\sigma_\ell}\circ\Gamma_\ell^{-1})_\#V$ and $\BC_\ell$ in place of $V$ and $\BC$, respectively. First, note that hypothesis (a) follows from the fact that $\Theta_V(0)\geq 2$ together with the monotonicity formula, and hypotheses (b) is immediate from the induction hypothesis, namely through \eqref{E:fd-2}, \eqref{E:fd-4}, and the triangle inequality. So, we can now focus on proving assumptions (c) and (d) of Lemma \ref{lemma:fine-improvement}.
	
	First we check (c). This will enable us to determine a suitable value of $M$ for the work as well. Starting with the triangle inequality and using that $\Gamma_\ell\in SO(n+k)$ we have
	\begin{align}
		\nonumber\hat{E}_{(\eta_{0,\sigma_\ell}\circ\Gamma_\ell^{-1})_\#V} & \leq C\hat{E}_{(\eta_{0,\sigma_\ell}\circ\Gamma_\ell^{-1})_\#V,\BC_\ell} + C \nu(\BC_\ell)\\
		& \leq C\hat{F}_{V,(\Gamma_\ell)_\#\BC_\ell}(B_{\sigma_\ell}(0)) + C\nu(\BC_\ell)\nonumber\\
		& \leq C\cdot 2^{-\ell}\hat{F}_{V,\BC} + C\nu(\BC_\ell)\nonumber\\
		& \leq C\cdot 2^{-\ell}\cdot\gamma\hat{E}_V + C\nu(\BC_\ell)\label{E:fd-7}
	\end{align}
	where in the third inequality we have used \eqref{E:fd-4} and in the last we used the assumption (D) of the present lemma. Look now at the second term on the right-hand side. Using the triangle inequality, as well as \eqref{E:fd-1} and \eqref{E:fd-2} for each $\ell^\prime\leq \ell$, we have
	\begin{align*}
		\nu(\BC_\ell) & \leq C\nu (\BC_\ell,(\Gamma_\ell)_\#\BC_\ell) + C\sum^\ell_{\ell^\prime=1}\nu((\Gamma_{\ell^\prime})_\#\BC_{\ell^\prime},(\Gamma_{\ell^\prime-1})_\#\BC_{\ell^\prime-1}) + C\nu(\BC)\\
		& \leq C\hat{F}_{V,\BC} + C\hat{F}_{V,\BC}\sum^{\ell}_{\ell^\prime=1}2^{-\ell^\prime} + C\nu(\BC)\\
		& \leq C\hat{F}_{V,\BC} + C\nu(\BC).
	\end{align*}
	Following the argument in Lemma \ref{lemma:nu-a} we have the estimate $\nu(\BC)\leq C\hat{E}_V$ (cf.~Remark \ref{remark:after-nu-a-second} following Lemma \ref{lemma:nu-a}). Combining this with the above, and \eqref{E:fd-7}, we deduce that
	\begin{equation}\label{E:fd-8}
		\hat{E}_{(\eta_{0,\sigma_\ell}\circ\Gamma_\ell^{-1})_\#V} \leq C\hat{E}_V.
	\end{equation}
	Suppose now we are given a plane $P$ with $S(\BC_\ell)\subset P$. Using assumption (c) of the present corollary as well as \eqref{E:fd-5}, \eqref{E:fd-1}, and \eqref{E:fd-4}, we have
	\begin{equation}\label{E:fd-8.5}
	\hat{E}_V \leq M\hat{E}_{V,P} \leq CM(\hat{E}_{V,\BC_{\ell-1}} + \nu(P,\BC_{\ell-1})) \leq CM(\hat{F}_{V,\BC} + \hat{E}_{(\Gamma_\ell^{-1})_\#V,P}(B_{\sigma_\ell}(0))).
	\end{equation}
	Using assumption (D) of the present corollary, by choosing $\gamma$ of the present proof sufficiently small (depending on $n,k,M,\beta$) we can absorb the first term on the right-hand side here into the left-hand side. Then, after taking the infimum over all such planes and combining this with \eqref{E:fd-8}, this gives
	$$\hat{E}_{(\eta_{0,\sigma_\ell}\circ\Gamma_\ell^{-1})_\#V} \leq CM\mathscr{E}_{(\eta_{0,\sigma_\ell}\circ\Gamma_\ell^{-1})_\#V, \BC_{\ell}}.$$
	Now, with $C$ equal to the constant on the right-hand side here, set $M^\prime:= CM$. This completes the verification of (c) of Lemma \ref{lemma:fine-improvement}, with $M^\prime$ in place of $M$. We can then let $\eps_0 = \eps_0(n,k,M,\beta)$ and $\gamma_0 = \gamma_0(n,k,M,\beta)$ be as in Lemma \ref{lemma:fine-improvement}. Notice also that the above argument gives $\nu(\BC_\ell)\leq C\hat{E}_V$ for some fixed $C$ independent of $\ell$.
	
	Lastly, we need to verify assumption (d) of Lemma \ref{lemma:fine-improvement}. We have from \eqref{E:fd-4} and assumption (D) of the present corollary,
	\begin{equation}\label{E:fd-9}
		\hat{F}_{(\eta_{0,\sigma_\ell}\circ\Gamma_\ell^{-1})_\#V,\BC_\ell} \leq 2^{-\ell}F_{V,\BC} \leq 2^{-\ell}\cdot\gamma\hat{E}_V.
	\end{equation}
	Applying \eqref{E:fd-8.5} with $P=P_0$ gives
	$$\hat{E}_V \leq C\hat{F}_{V,\BC} + C\hat{E}_{(\eta_{0,\sigma_\ell}\circ\Gamma_\ell^{-1})V}.$$
	Using now assumption (D) of the present corollary, the first term on the right-hand side here can be absorbed in the left-hand side, leaving us with
	$$\hat{E}_V \leq C\hat{E}_{(\eta_{0,\theta^\ell}\circ\Gamma_\ell^{-1})_\#V}.$$
	Combining this with \eqref{E:fd-9} gives
	$$\hat{F}_{(\eta_{0,\theta^\ell}\circ\Gamma_\ell^{-1})_\#V,\BC_\ell} \leq C2^{-\ell}\gamma\hat{E}_{(\eta_{0,\theta^\ell}\circ\Gamma_\ell^{-1})_\#V}.$$
	By choosing $\gamma$ of the present corollary sufficiently small we can ensure that $C\gamma<\gamma_0$, and thus this establishes (d) of Lemma \ref{lemma:fine-improvement}. Thus, we can apply Lemma \ref{lemma:fine-improvement} with $(\Gamma_\ell^{-1})_\#V_\ell$ and $\BC_\ell$ in place of $V$ and $\BC$. This immediately verifies \eqref{E:fd-1}, \eqref{E:fd-2}, \eqref{E:fd-3}, \eqref{E:fd-4}, and \eqref{E:fd-5} inductively, and Remark \ref{remark:after-fine-improvement} following the proof of Lemma \ref{lemma:fine-improvement} verifies \eqref{E:fd-6} (as we have chosen $\tilde{\eps}, \tilde{\gamma}$ sufficiently small). This therefore completes the inductive step of the proof of the claim, and thus proves the claim.
	
	It then follows that $(\Gamma_j)_j$ and $(\BC_j)_j$ are Cauchy sequences\footnote{Indeed, we are getting that $(\eta_{0,\rho})_\#V$ has a unique limit as $\rho\to 0$ which is not a single plane, and thus $\mathscr{E}_{V,\BC}(B_\rho(0))$ has a uniform lower bound in terms of $\rho$. Thus, from conclusion (1) of Lemma \ref{lemma:fine-improvement}, we do know that $|\Gamma_j-\Gamma_{j-1}|\leq C\gamma_0\cdot 2^{-j}$ for all $j$ with a constant $C$ independent of $j$, giving that $(\Gamma_j)_j\subset SO(n+k)$ is a Cauchy sequence. This also shows that, for all $j$, $\Gamma_j$ is close to $\id_{\R^{n+k}}$, and so we can guarantee that $(\Gamma_j)_\#V\in \mathcal{V}_{\beta/2}$ for all $j$ by Proposition \ref{prop:rotation-class}.}, and so by completeness of their respective spaces they converge; thus we have $\Gamma_j\to \Gamma\in SO(n+k)$ and $\BC_\ell\to \BC^\prime\in \CC$ (the fact that $\BC^\prime$ is not a plane follows from assumptions (C) and (D)). It follows immediately that conclusion (1) of the corollary holds. Next, we interpolate between scales: pick $\rho\in (0,3\theta_{n+1}/4)$ and let $\ell\in \N$ be the unique positive integer for which $\sigma_{\ell+1}\leq\rho<\sigma_\ell$. Then, we ask that $2^{-\ell}\leq \rho^{\mu_*}$, i.e.~we pick $\mu_*\in (0,1)$ such that $-\ell\log(2) \leq \mu_*\log(\rho)$. This is satisfied if $\ell\log(1/2) \leq \mu_*(\ell+1)\log(1/\theta_{n+1})$,
	i.e.~for
	$$\mu_* \leq \frac{\ell}{\ell+1}\cdot\frac{\log(1/2)}{\log(1/\theta_{n+1})}$$
	and so in particular we can take $\mu_* = \frac{1}{2}\cdot\frac{\log(2)}{\log(\theta_{n+1})}$; notice that this only depends on $n,k,M,\beta$. Then, from \eqref{E:fd-4} we have
	$$\hat{E}_{V,\Gamma_\#\BC^\prime}(B_\rho(0)) \leq C\rho^{\mu_*}\hat{F}_{V,\BC}$$
	and this holds for all $\rho\in (0,3\theta/4)$; we stress that this is only decay at all scales for the \emph{one-sided} fine excess, as for the second part of $\hat{F}_{V,\Gamma_\#\BC^\prime}$ we only get the decay at the discrete set of scales $\{\sigma_\ell\}_{\ell= 1}^\infty$, as the domains of integration are not nested. Combining this with Allard's supremum estimate (Theorem \ref{thm:allard-sup-estimate}) it is straightforward to convert this into conclusion (2), in which the cylindrical $L^2$ height-excess appears on the left-hand side. Using the above and the triangle inequality in \eqref{E:fd-6} we then have
	$$\nu(\Gamma_\#\BC^\prime,(\tau_z)_\#\Gamma_\#\BC^\prime) \leq C\hat{F}_{V,\Gamma_\#\BC^\prime}(B_\rho(0))$$
	for any $\rho\in (0,\theta/2)$ and any $z\in B_{\rho/2}(0)\cap \mathcal{D}(V)$. Combining this with the estimate just above gives conclusion (3).
\end{proof}

We can now finally prove the fine $\eps$-regularity theorem, as stated in Theorem \ref{thm:fine-reg}. This will follow by applying Corollary \ref{cor:fine-decay} at different base points $z\in\mathcal{D}(V)\cap B_{1/2}(0)$ and utilising Lemma \ref{lemma:gap-2}.

\begin{proof}[Proof of Theorem \ref{thm:fine-reg}]
	To start, we pick an optimal coarser cone: by arguing the same way as in Step 1 of the (inductive) proof of Theorem \ref{thm:graphical-rep}, we have that, given $\tilde{\eps}$, $\tilde{\gamma}>0$, we can choose $\eps$ and $\gamma$ of the present theorem sufficiently small, depending on $n,k,M,\beta,\tilde{\eps},\tilde{\gamma}$, such that there exists a cone $\Dbf\in \CC$, a plane $\widetilde{P}$, and a rotation $R\in SO(n+k)$ such that:
	\begin{itemize}
		\item $S(\BC)\subset S(\Dbf)\subsetneq \widetilde{P}$;
		\item $R(\widetilde{P}) = P_0$;
		\item $|R-\id_{\R^{n+k}}|\leq C\hat{E}_V$;
		\item $R_\#V$ and $R_\#\Dbf$ satisfy Hypothesis $H(\tilde{\eps},\tilde{\gamma},1)$ in $B_1(0)$;
		\item $\nu(\BC,\Dbf)\leq C\hat{F}_{V,C} < C\gamma\hat{E}_{V,\widetilde{P}} \leq C\gamma\hat{E}_V$.
	\end{itemize}
	Indeed, in the case where the choice of $\Dbf_j$ within Step 1 of the proof of Theorem \ref{thm:graphical-rep} is a plane, this means that $\hat{E}_{R_\#V}\leq C\mathscr{F}^*_{V,\BC}$, and so by assumptions (d) and (e) of Theorem \ref{thm:fine-reg} we can simply take $\Dbf := \BC$ and $R = \id_{\R^{n+k}}$. Otherwise, we take the choice of $\Dbf_j$ to be $\Dbf$ and similarly for $R$.
	
	In particular, we can choose $\tilde{\eps}$ and $\tilde{\gamma}$ sufficiently small (depending on $n,k,M,\beta$) so that Theorem's \ref{thm:graphical-rep}, \ref{thm:L2-estimates}, \ref{thm:shift-estimate}, and Corollary \ref{cor:L2-corollary} all hold for $R_\#V$ and $R_\#\Dbf$ (with now $\beta/2$ in place of $\beta$); this now means $\eps$ and $\gamma$ only depend on $n,k,M,\beta$. Then by arguing as in the proof of Lemma \ref{lemma:fine-translation} for verifying (A3) within, we can check that for any $z\in \mathcal{D}(V)\cap B_{1/2}(0)$, the hypotheses of Corollary \ref{cor:fine-decay} are satisfied with $V_z$ and $\BC$ in place of $V$ and $\BC$, respectively, where $V_z := (\eta_{z,1/4})_\#R_\#V\res B_1(0)$. Moreover, we can check that we have the estimate (cf.~the argument leading to \eqref{E:b-2-1})
	\begin{equation}\label{E:fine-reg-1}
		\hat{F}_{V_z,R_\#\Dbf} \leq C\hat{F}_{R_\#V,R_\#\Dbf}
	\end{equation}
	for some $C = C(n,k,M,\beta)\in (0,\infty)$.
	As $R_\#V$ and $R_\#\Dbf$ satisfy Hypothesis $H(\tilde{\eps},\tilde{\gamma},1)$, Remark \ref{remark:after-graphical-rep} following Theorem \ref{thm:graphical-rep} gives that $\hat{F}_{R_\#V,R_\#\Dbf} \leq C\hat{E}_{R_\#V,R_\#\Dbf} \equiv C\hat{E}_{V,\Dbf}$, and thus we also know that by the construction of $\Dbf$ that $\hat{F}_{V_z,R_\#\Dbf} \leq C\hat{F}_{V,\BC}$. This means that the result of applying Corollary \ref{cor:fine-decay} to $V_z$ and $\BC$ is that there exist fixed numbers $\theta = \theta(n,k,M,\beta)\in (0,1/8^{n+3})$ and $\mu = \mu(n,k,M,\beta)\in (0,1)$ such that for any $z\in\mathcal{D}(V)\cap B_{1/2}(0)$, there exists an aligned cone $\BC_z^\prime\in \CC$ and a rotation $\Gamma_z\in SO(n+k)$ such that (incorporating the rotation $R$ into $\Gamma_z$):
	\begin{equation}\label{E:fine-reg-2}
		\nu(\BC,(\Gamma_z)_\#\BC^\prime_z) \leq C\hat{F}_{V,\BC};
	\end{equation}
	\begin{equation}\label{E:fine-reg-3}
		\hat{E}_{V,(\tau_z)_\#(\Gamma_z)_\#\BC^\prime_z}(C_\rho(\pi_{P_0}(z))) \leq C\rho^\alpha\hat{F}_{V,\BC} \qquad \text{for all }\rho\in (0,\theta/8);
	\end{equation}
	\begin{equation}\label{E:fine-reg-4}
		C^{-1}\hat{E}_V \leq \hat{E}_{(\eta_{0,\rho})_\#V} \leq C\hat{E}_V \qquad \text{for all }\rho\in (0,1);
	\end{equation}
	and for every $y\in \mathcal{D}(V_z)\cap B_{\rho/2}(0)$ and $\rho\in (0,\theta/2)$,
	\begin{equation}\label{E:fine-reg-5}
		\nu((\Gamma_z)_\#\BC_z^\prime,(\tau_y)_\#(\Gamma_z)_\#\BC_z^\prime) \leq C\rho^\mu\hat{F}_{V,\BC}.
	\end{equation}
	Here, $C = C(n,k,M,\beta)\in (0,\infty)$. Notice also that Corollary \ref{cor:fine-decay} also gives that the tangent cone to $V$ at such a point $z$ is unique and is equal to $\BC_z := (\Gamma_z)_\#\BC_z^\prime$.
	
	Next, consider a point $x_0\in (\spt\|V\|\setminus \mathcal{D}(V))\cap C_{1/2}(0)$ with $\dist(x,\mathcal{D}(V))\in (0,\theta/32)$. Pick $z\in \mathcal{D}(V)\cap C_{1/2}(0)$ with $|x_0-z| = \dist(x,\mathcal{D}(V))$. Then $\tilde{x}_0:= 4(R(x_0)-z)\in \spt\|V_z\|$ has $\dist(\tilde{x}_0,\mathcal{D}(V_z)) = |\tilde{x}_0| \in (0,\theta/8)$. Utilising \eqref{E:fine-reg-3}, we can then rescale and consider $(\eta_{0,|\tilde{x}_0|/2})_\#V_z$, and apply Lemma \ref{lemma:gap-2}, to see that (after writing everything back in terms of $V$) we have
	$$V\res C_{\dist(x,\mathcal{D}(V))/16}(\pi_{P_0}(x_0)) = |\graph(u^1_{x_0})| + |\graph(u^2_{x_0})|$$
	for two smooth functions $u^1_{x_0},u^2_{x_0}:B^n_{\theta/16}(x_0)\to P_0^\perp$ each of which solves the minimal surface equation and therefore satisfies standard elliptic estimates (this can be done assuming only that $\eps$ is sufficiently small depending on $n,k$). If instead $x_0\in (\spt\|V\|\setminus \mathcal{D}(V))\cap C_{1/2}(0)$ has $\dist(x,\mathcal{D}(V))>\theta/32$, then we can directly apply Lemma \ref{lemma:gap-2} to reach the same conclusion, provided $\eps$ is chosen to be sufficiently small depending only on $n,k,\beta,\theta$; since $\theta = \theta(n,k,M,\beta)\in (0,1)$, it follows that a \emph{single} choice of $\eps$ can be chosen, depending only on $n,k,M,\beta$, such that this conclusion holds on a neighbourhood of every $x_0 \in (\spt\|V\|\setminus \mathcal{D}(V))\cap C_{1/2}(0)$. In particular, we must have that
	$$\sing(V) \cap C_{1/2}(0) = \mathcal{D}(V)\cap C_{1/2}(0).$$
	From here, returning to Almgren's weak Lipschitz approximation, Theorem \ref{thm:Lipschitz-approx} (with the precise form of the bad set as in \cite[Corollary 3.11]{Alm00}), one can check that $\Sigma = \emptyset$, i.e.~$V\res C_{1/2}(0)$ agrees everywhere with a Lipschitz $2$-valued graph. Thus, conclusion (1) follows. Conclusion (3) then follows using Allard's supremum estimate from Theorem \ref{thm:allard-sup-estimate} as well as the estimates above. Notice we have already shown (2), from \eqref{E:fine-reg-3} and \eqref{E:fine-reg-4}.
	
	To show conclusion (4) requires a little more work. Pick two distinct points $z_1,z_2\in\mathcal{D}(V)\cap C_{1/4}(0)$, and write $\sigma:= |\pi_{P_0}(z_1)-\pi_{P_0}(z_2)|$ (which we know from (1) must obey $\sigma>0$). Notice that $\tilde{z}_1:= \Gamma_{z_2}^{-1}(4(R(z_1)-z_2))\in \mathcal{D}((\Gamma_{z_2}^{-1})_\#V_{z_2})$. Write $W:= (\eta_{\tilde{z}_1,4\sigma})_\#(\Gamma^{-1}_{z_2})_\#V_{z_2}$. Now, let us suppose that $\sigma<\theta/64$. Then, writing $\BC_{z_2}^* := (\Gamma_{z_2}^{-1}\circ R\circ \Gamma_{z_2})_\#\BC_{z_2}^\prime$ for simplicity,
	\begin{align*}
		\hat{F}_{W,\BC_{z_2}^\star} & = \hat{F}_{(\Gamma^{-1}_{z_2})_\#V_{z_2},(\tau_{\tilde{z}_1})_\#\BC_{z_2}^\star}(B_{4\sigma}(\tilde{z}_1))\\
		& \leq C\hat{F}_{(\Gamma^{-1}_{z_2})_\#V_{z_2},(\tau_{\tilde{z}_1})_\#\BC_{z_2}^\star}(B_{8\sigma}(0))\\
		& \leq C\hat{F}_{(\Gamma_{z_2}^{-1})_\#V_{z_2},\BC_{z_2}^\star}(B_{8\sigma}(0)) + C\nu((\tau_{\tilde{z}_1})_\#\BC_{z_2}^\star,\BC_{z_2}^\star)\\
		& \leq C\sigma^\mu\hat{F}_{V,\BC}
	\end{align*}
	where the last term is controlled by Theorem \ref{thm:shift-estimate} applied to $R_\#V$ and $R_\#\Dbf$, using the bounds then relating $\Dbf$ to $\BC$. We can then proceed to check all the hypotheses of Corollary \ref{cor:fine-decay} for $W$ and $(\Gamma_{z_2}^{-1}\circ R\circ \Gamma_{z_2})_\#\BC^\prime_{z_2}$ in place of $V$ and $\BC$. The result is that after applying Corollary \ref{cor:fine-decay}, we deduce the existence of another aligned cone $\widetilde{\BC}_{z_1}\in \CC$ and rotation $\widetilde{\Gamma}_{z_1}$ for which
	\begin{equation}\label{E:fine-reg-6}
		\nu((\Gamma_{z_2}^{-1}\circ R\circ\Gamma_{z_2})_\#\BC_{z_2}^\prime,(\widetilde{\Gamma}_{z_1})_\#\widetilde{\BC}_{z_1})\leq C\sigma^\mu\hat{F}_{V,\BC}.
	\end{equation}
	But unravelling the definition of $W$ and using the uniqueness of tangent cones already shown, we see that
	$$(\widetilde{\Gamma}_{z_1})_\#\widetilde{\BC}_{z_1} = (\Gamma_{z_2}^{-1}\circ R \circ \Gamma_{z_1})_\#\BC_{z_1}^\prime$$
	which if we substitute back into \eqref{E:fine-reg-6} we see
	$$\nu((\Gamma_{z_1})_\#\BC_{z_1}^\prime,(\Gamma_{z_2})_\#\BC_{z_2}^\prime) \leq C\sigma^\mu\hat{F}_{V,\BC}.$$
	Recalling that for general $z\in \mathcal{D}(V)\cap C_{1/2}(0)$ we have $\BC_{z} = (\Gamma_{z})_\#\BC_{z}^\prime$, conclusion (4) therefore follows. (The case $\sigma\geq \theta/64$ follows more easily from the previous estimates.) This completes the proof.
\end{proof}

This therefore completes the proof of the fine $\eps$-regularity theorem, which is the crucial technical tool to understanding coarse blow-ups of varifolds in $\mathcal{V}_\beta$.

\textbf{Remark:} We elaborate briefly on Remark \ref{remark:fine-reg-sufficiency} regarding the sufficiency of ruling out $(\beta,\gamma)$-fine gaps which fail the topological structural condition in order to establish Theorem \ref{thm:fine-reg}. The reader can verify that, during the proof of Theorem \ref{thm:fine-reg}, we would only use Lemma \ref{lemma:gap-2} at times when we knew the graph of a coarse blow-up off $P_0$ was a cone in $\CC$, and as such we would be able to verify that the $\beta$-coarse gap was in fact a $(\beta,\gamma)$-coarse gap along the blow-up sequence of varifolds. This is most evident in the proof of Theorem \ref{thm:graphical-rep}, as well as Lemma \ref{lemma:tori-q=1} and Lemma \ref{lemma:tori}. There is one place where one needs the topological structural condition in a $\beta$-coarse gap where the $(\beta,\gamma)$-fine gap condition is satisfied by a cone $\BC_*$ (say) which is \emph{not} in $\CC$, namely when verifying Lemma \ref{lemma:nu-a}(iv) in the case $\BC\in \CC\setminus\mathcal{P}$. One sees in the proof that the cone $\BC_*$ is the union of $4$ half-planes where two of the half-planes \emph{coincide}. However, since such a cone can be approximately to an arbitrary degree by a cone which is formed of the union of $4$ \emph{distinct} half-planes, this is still a $(\beta,\gamma)$-coarse gap as defined, and so the topological structural condition still applies.

\textbf{Remark:} We wish to remark an alternative approach to the proof of Theorem \ref{thm:fine-reg}, which establishes the classification of homogeneous degree one fine blow-ups (Theorem \ref{thm:fine-classification}) in a slightly different manner. This elaborates on Remark \ref{remark:cascade}. Of course, one must first establish the main estimates in Theorem's \ref{thm:graphical-rep}, \ref{thm:L2-estimates}, \ref{thm:shift-estimate}, and Corollary \ref{cor:L2-corollary} by downwards induction on the spine dimension, which in turn requires inductively constructing and proving the continuity estimate for fine blow-ups (albeit only if needed in a weak form). From this, the regularity result for fine blow-ups in $\mathfrak{b}_{\sfrak}$, with $\sfrak\leq n-2$, is immediate by removability theorems for continuous harmonic functions. At this point, one could establish a form of Theorem \ref{thm:fine-reg} for when the cone $\BC$ obeys $\sfrak(\BC)\leq n-2$, by replacing the hypothesis $\hat{F}_{V,\BC}<\gamma \hat{E}_V$ with the hypothesis 
\begin{equation}\label{E:fine-remark-1}
\hat{F}_{V,\BC}<\gamma \inf_{\Dbf}\hat{F}_{V,\Dbf},
\end{equation}
where the infimum here is taken over all coarser cones $\Dbf$ to $\BC$ with spine dimension $\sfrak(\Dbf)=n-1$. In fact, in this case the conclusion of Theorem \ref{thm:fine-reg} can be improved to $V$ decomposing as a sum of two minimal graphs. We will state this result (and others) in Section \ref{sec:additional-fine} (cf.~Theorem \ref{thm:fine-reg-extra-1}. Given this, the proof of Theorem \ref{thm:fine-reg} can be reduced to working with a cone of spine dimension precisely $n-1$. To classify homogeneous degree one fine blow-ups for $\mathfrak{b}_{n-1}$, and subsequently establish their regularity, one can first establish an $\eps$-regularity theorem for fine blow-ups $\phi$ in $\mathfrak{b}_{n-1}$ whose graphs are close to unions of planes which do not intersect \emph{exactly} along $S_0$; in particular, this would mean that the sequence of varifolds $(V_j)_j$ whose fine blow-up is $\phi$ fall into the regime to apply the special case of Theorem \ref{thm:fine-reg} just mentioned, under the assumption \eqref{E:fine-remark-1}. This allows one to change the proof of the reverse Hardt--Simon inequality within the proof of Theorem \ref{thm:fine-reg} to a dichotomy, saying that either we only need to subtract off $\pi\in \mathfrak{l}$ whose graph has spine dimension $n-1$, or we automatically get a regularity conclusion. The advantage of this approach is that it only requires ever taking fine blow-ups of cones of the same spine dimension, leading to a simpler form of Lemma \ref{lemma:fine-integration}. One can then complete the proof of the classification of fine blow-ups in a similar manner, and then prove the full form of Theorem \ref{thm:fine-reg}. In fact, this approach is morally the same as what we do in Part \ref{part:blow-up-reg} for the regularity of coarse blows-ups and the proof of Theorem \ref{thm:main}. The reason we cannot do an approach for coarse blow-ups as we did for fine blow-ups in Theorem \ref{thm:fine-classification} is because for coarse blow-ups we do not have the freedom to subtract off such a variety of linear functions, as we do not have the same flexibility as we had with the spine dimension in the fine blow-up procedure. For coarse blow-ups, as we saw in $(\mathfrak{B}3\text{III})$, the only flexibility we have is subtracting single, linear functions.

\subsection{Refined versions of Theorem \ref{thm:fine-reg}}\label{sec:additional-fine}

Here we wish to record two refined versions of the fine $\eps$-regularity theorem given in Theorem \ref{thm:fine-reg}. These results are not needed for the remainder of the proof of Theorem \ref{thm:main}, they simply provide refined conclusions in Theorem \ref{thm:fine-reg} under additional hypotheses. They should be compared with \cite[Theorem 1]{BK17} and \cite[Theorem 3]{BK17}, respectively (Theorem \ref{thm:fine-reg}, on the other hand, should be compared with \cite[Theorem 2]{BK17}).

The first refinement can be viewed as a ``small spine'' case, where one assumes closeness to a cone in $\mathcal{P}_{\leq n-2}$ compared to any cone in $\mathcal{P}_{n-1}$. This allows for a significantly improved regularity theorem.

\begin{theorem}[Fine $\eps$-Regularity Theorem: Small Spine Refinement]\label{thm:fine-reg-extra-1}
	Fix $M\in [1,\infty)$ and $\beta\in (0,1)$. Fix also $\alpha\in (0,1)$. Then, there exists $\eps_0 = \eps_0(n,k,M,\beta,\alpha)\in (0,1)$ and $\gamma = \gamma(n,k,M,\beta,\alpha)\in (0,1)$ such that the following is true. Suppose that $V\in\mathcal{V}_\beta$ and $\BC\in \mathcal{P}_{\leq n-2}$ satisfy:
	\begin{enumerate}
		\item [(a)] $\Theta_V(0)\geq 2$ and $(\w_n2^n)^{-1}\|V\|(B^{n+k}_2(0))\leq \frac{5}{2}$;
		\item [(b)] $0\in S(\BC)\subset S_0$;
		\item [(c)] $\nu(\BC)<\eps$ and $\hat{E}_V<\eps$;
		\item [(d)] $\hat{E}_V \leq M\mathscr{E}_{V,\BC}$;
		\item [(e)] $\hat{F}_{V,\BC} < \gamma\inf_{\Dbf\in\mathcal{P}_{n-1}}\hat{F}_{V,\Dbf}$.
	\end{enumerate}
	Then, there is a cone $\BC_0\in \mathcal{P}_{\leq n-2}$ with $0\in S(\BC_0)\subset S_0$ and an orthogonal rotation $\Gamma\in SO(n+k)$ with
	$$\dist_\H(\spt\|\BC_0\|\cap C_1(0),\spt\|\BC\|\cap C_1(0))\leq C\hat{F}_{V,\BC},$$
	$|\Gamma(e_\kappa)-e_\kappa|\leq C\hat{F}_{V,\BC}$ for $\kappa=1,\dotsc,k$ and $|\Gamma(e_{k+i})-e_{k+i}|\leq C\hat{E}_V^{-1}\hat{F}_{V,\BC}$ for $i=1,\dotsc,n$, such that $\BC_0$ is the unique tangent cone to $\Gamma_\#^{-1}V$ at $0$ and
	$$\sigma^{-n-2}\int_{C_\sigma(0)}\dist^2(x,\spt\|\BC_0\|)\, \ext\|\Gamma_\#^{-1}V\|(x) \leq C\sigma^{2\alpha}\hat{F}_{V,\BC}^2 \qquad \text{for all }\sigma\in (0,1/2).$$
	Furthermore, if $\BC_0 = |P_0^1| + |P_0^2|$ where $P_0^i$ are $n$-dimensional planes passing through the origin, then there exist two smooth (single-valued) functions $u_1,u_2$ with $u_i:P_0^i\cap B^{n+k}_{1/2}(0)\to (P_0^i)^\perp$, $i=1,2$, both solving the minimal surface system and such that
	\begin{enumerate}
		\item [(1)] $V\res C_{1/2}(0) = (\mathbf{v}(u_1) + \mathbf{v}(u_2))\res C_{1/2}(0)$ and $\sing(V)\cap C_{1/2}(0) = \graph(u_1) \cap \graph(u_2) \cap C_{1/2}(0)$;
        \item [(2)] $|u_i|_{1,\alpha;P_0^i\cap B^{n+k}_{1/2}(0)}\leq C\hat{F}_{V,\BC}$ for each $i=1,2$;
        \item [(3)] At each $x\in \sing(V)\cap C_{1/2}(0)$, the (unique) tangent cone to $V$ at $x$ belongs to $\mathcal{P}_{\leq n-2}$.
	\end{enumerate}
	Here, $C = C(n,k,M,\beta)\in (0,\infty)$.
\end{theorem}

The proof of this theorem follows the same logic as in Theorem \ref{thm:fine-reg}, utilising the fact that under the assumption (e) above, every tangent cone must remain in $\mathcal{P}_{\leq n-2}$ and as such the two-valued function decouples into two single-valued functions. Of course, in this case the regularity can be improved using standard elliptic regularity theory for the minimal surface system.

The second refinement is in the other extra case, where $V$ is significantly closer to a cone in $\CC_{n-1}\setminus\mathcal{P}$ than any cone in $\mathcal{P}$. In this situation, the conclusion can be refined to all tangent cones being in $\CC_{n-1}\setminus\mathcal{P}$, and the singularities forming a connected $C^{1,\alpha}$ submanifold.

\begin{theorem}[Fine $\eps$-Regularity Theorem: Twisted Refinement]\label{thm:fine-reg-extra-2}
	Fix $M\in [1,\infty)$ and $\beta\in (0,1)$. Fix also $\alpha\in (0,1)$. Then, there exists $\eps_0 = \eps_0(n,k,M,\beta,\alpha)\in (0,1)$ and $\gamma = \gamma(n,k,M,\beta,\alpha)\in (0,1)$ such that the following is true. Suppose that $V\in\mathcal{V}_\beta$ and $\BC\in \mathcal{C}_{n-1}\setminus\mathcal{P}$ satisfy:
	\begin{enumerate}
		\item [(a)] $\Theta_V(0)\geq 2$ and $(\w_n2^n)^{-1}\|V\|(B^{n+k}_2(0))\leq \frac{5}{2}$;
		\item [(b)] $S(\BC) = S_0$;
		\item [(c)] $\nu(\BC)<\eps$ and $\hat{E}_V<\eps$;
		\item [(d)] $\hat{E}_V \leq M\mathscr{E}_{V,\BC}$;
		\item [(e)] $\hat{F}_{V,\BC} < \gamma\inf_{\Dbf\in\mathcal{P}}\hat{F}_{V,\Dbf}$.
	\end{enumerate}
	Then, there is a cone $\BC_0\in \mathcal{C}_{n-1}\setminus\mathcal{P}$ with $S(\BC_0) = S_0$ and an orthogonal rotation $\Gamma\in SO(n+k)$ with
	$$\dist_\H(\spt\|\BC_0\|\cap C_1(0),\spt\|\BC\|\cap C_1(0))\leq C\hat{F}_{V,\BC},$$
	$|\Gamma(e_\kappa)-e_\kappa|\leq C\hat{F}_{V,\BC}$ for $\kappa=1,\dotsc,k$ and $|\Gamma(e_{k+i})-e_{k+i}|\leq C\hat{E}_V^{-1}\hat{F}_{V,\BC}$ for $i=1,\dotsc,n$, such that $\BC_0$ is the unique tangent cone to $\Gamma_\#^{-1}V$ at $0$ and
	$$\sigma^{-n-2}\int_{C_\sigma(0)}\dist^2(x,\spt\|\BC_0\|)\, \ext\|\Gamma_\#^{-1}V\|(x) \leq C\sigma^{2\alpha}\hat{F}_{V,\BC}^2 \qquad \text{for all }\sigma\in (0,1/2).$$
	Furthermore, there exists a generalised-$C^{1,\alpha}$ function $u:B^n_{1/2}(0)\to \A_2(P_0^\perp)$ such that:
	\begin{enumerate}
		\item [(1)] $V\res C_{1/2}(0) = \mathbf{v}(u)$;
		\item [(2)] $\mathcal{B}_u \cap B_{1/2} = \emptyset$;
		\item [(3)] $\graph(u)\cap (\R^k\times \CC_u) = \sing(V)\cap C_{1/2}(0)$; moreover,
		\begin{enumerate}
			\item [(i)] $\sing(V)\cap C_{1/2}(0) = \graph(\phi)$, where $\phi: S_0\cap B_{1/2}^{n-1}(0)\to S_0^\perp$ is of class $C^{1,\alpha}$ over $\overline{S_0\cap B_{1/2}^{n-1}(0)}$, where if $\phi = (\phi_1,\dotsc,\phi_k,\phi_{k+1})$, then $|\phi_\kappa|_{1,\mu;B_{1/2}^{n-1}(0)}\leq C\hat{F}_{V,\BC}$ for $\kappa=1,\dotsc,k$ and $|\phi_{k+1}|_{1,\mu;B_{1/2}^{n-1}(0)}\leq C\hat{E}_{V}^{-1}\hat{F}_{V,\BC}$;
			\item [(ii)] If $\Omega^\pm$ denote the two connected components of $B_{1/2}^n(0)\setminus \graph(\phi_{k+1})$, then on each of $\Omega^\pm$ we have that we can write $u|_{\Omega^\pm} = \llbracket u^1_\pm\rrbracket + \llbracket u^2_\pm\rrbracket$, where $u^\alpha_\pm \in C^{1,\mu}(\overline{\Omega^\pm};\R^k)$ solve the minimal surface system over their domains, and $|u^\alpha_\pm|_{1,\mu;\Omega^\pm}\leq C\hat{E}_V$ for each $\alpha=1,2$;
		\end{enumerate}
		\item [(4)] If $\BC_z\in \CC$ denotes the (unique) tangent cone to $V$ at $z\in \sing(V)\cap C_{1/2}(0)$, then $\BC_z\in \CC_{n-1}\setminus\mathcal{P}$, $C^{-1}\hat{E}_V\leq \dist_\H(\BC_z\cap C_1(0), P_0\cap C_1(0))\leq C\hat{E}_V$; moreover, if $z_1,z_2\in \sing(V)\cap C_{1/2}(0)$ then
		$$\dist_\H(\spt\|\BC_{z_1}\|\cap C_1(0), \spt\|\BC_{z_2}\|\cap C_1(0))\leq C|\pi_{P_0}(z_1)-\pi_{P_0}(z_2)|^{\alpha}\hat{F}_{V,\BC}.$$
	\end{enumerate}
	Here, $C = C(n,k,M,\beta)\in (0,\infty)$. In particular, $V\res C_{1/2}(0)$ is equal to the union of four, smooth, minimal submanifolds with boundary meeting only along an $(n-1)$-dimensional $C^{1,\alpha}$ submanifold, which is their common boundary.
\end{theorem}

The proof of this theorem again follows the same logic as in Theorem \ref{thm:fine-reg}, utilising the fact that under the assumption (e) above, every tangent cone must remain in $\mathcal{C}_{n-1}\setminus\mathcal{P}$. Of course, in this case the regularity can be improved using standard (boundary and interior) elliptic regularity theory for the minimal surface system.

\section{Appendix to Part \ref{part:fine-reg}}

We prove two elementary lemmas which were used in the inductive proofs of Theorem \ref{thm:graphical-rep}, \ref{thm:L2-estimates}, and \ref{thm:shift-estimate}. The first is a type of triangle inequality. The second gives a criterion for controlling the excess relative to a single plane by the union of two planes.

\begin{lemma}\label{lemma:appendix-1}
	Fix $Q\in \{1,2,\dotsc\}$. Then, there exists $\eps_0 = \eps_0(n,k,Q)\in (0,1)$ such that the following is true. Suppose that $V$ is a stationary integral $n$-varifold in $B^{n+k}_2(0)$, and $\BC = \sum^Q_{\alpha=1}|P^\alpha|$ is a sum of $Q$ planes for which the following hypotheses hold:
	\begin{enumerate}
		\item [(1)] $Q-\frac{1}{2} < (\w_n 2^n)^{-1}\|V\|(C_2(0)) \leq Q+\frac{1}{2}$;
		\item [(2)] $\hat{E}_V<\eps_0$;
		\item [(3)] $\nu(\BC)<\eps_0$.
	\end{enumerate}
	Then, there is a non-empty subset $I\subset\{1,\dotsc,Q\}$ such that
	$$\sum_{i\in I}\nu(P^i) \leq C(\hat{E}_{V,\BC} + \hat{E}_V)$$
	for some $C = C(n,k,Q)\in(0,\infty)$.
\end{lemma}

\begin{proof}
	Let $u:B^n_{15/16}(0)\to \A_Q(\R^k)$ denote the $Q$-valued Lipschitz approximation of $V$ given by Theorem \ref{thm:Lipschitz-approx} (with $\sigma=15/16$ and $L=1/2$ therein); this can be achieved provided $\eps_0 = \eps_0(n,k,Q)$ is sufficiently small. Let $p^1,\dotsc,p^Q:B^n_1(0)\to \R^k$ denote the linear functions for which $P^\alpha\cap C_1(0) = \graph(p^\alpha)$.
	
	Let $S_1,\dotsc,S_J$ be an ordered list of the set of $Q$-tuples $(s_1,\dotsc,s_Q)$ of non-negative integers obeying $\sum_i s_i = Q$. Clearly $J = J(Q)<\infty$. For each $i\in \{1,\dotsc,J\}$, define a subset $\Omega_{S_i}\subset B^n_{15/16}(0)\setminus\Sigma$ as follows (here, $\Sigma$ is as in Theorem \ref{thm:Lipschitz-approx}). If $S_i = (s_1,\dotsc,s_Q)$, say that $x\in B^n_{15/16}(0)\setminus\Sigma$ belongs to $\Omega_{S_i}$ if for every $\beta\in \{1,\dotsc,Q\}$,
	\begin{equation}\label{E:app-1-1}
		\#\{\alpha\in\{1,\dotsc,Q\}:|u^\alpha(x) - p^\beta(x)| = \min_\gamma|u^\alpha(x)-p^\gamma(x)|\} = s_\beta
	\end{equation}
	i.e.~if $s_1$ points in $\R^k\times\{x\}\cap \spt\|V\|$ are closest to $p_1$, $s_2$ points are closest to $p_2$, etc. (Here, to avoid repetitions, if there is a choice of planes realising the minimum for a given point we default to the one with the lowest index $\beta$.) Now, the sets $\{\Omega_{S_i}\}_{i=1}^J$ cover $B^n_{15/16}(0)\setminus\Sigma$ and furthermore are measurable by construction. Thus, there must be some $i_*\in \{1,\dotsc,J\}$ with $|\Omega_{S_{i_*}}|\geq c(n,Q)>0$. We will now entirely work with $\Omega_{S_{i_*}}$. Write $S_{i_*} = (s_1,\dotsc,s_Q)$. Now, by the triangle inequality and the defining property \eqref{E:app-1-1} of $S$ we have for each $x\in\Omega_{S_{i_*}}$
	\begin{equation}\label{E:app-1-3}
		\sum^Q_{\alpha=1}s_\alpha |p^\alpha(x)| \leq \sum^Q_{\alpha=1}\dist(u^\alpha(x),\{p^1(x),\dotsc,p^Q(x)\}) + \sum^Q_{\alpha=1}|u^\alpha(x)|.
	\end{equation}
	Also observe that since $\nu(\BC)<\eps_0$, there is a $C = C(n,k)$ for which
	$$\dist(u^\alpha(x),\{p^1(x),\dotsc,p^Q(x)\})\leq C\dist(u^\alpha(x),\BC).$$
	Combining this with \eqref{E:app-1-3}, squaring both sides and then integrating over $\Omega_{S_{i_*}}$, we get
	$$\sum_{\alpha:s_\alpha\neq 0}\int_{\Omega_{S_{i_*}}}|p^\alpha|^2 \leq C\int_{B_{15/16}(0)\setminus\Sigma}\sum^Q_{\alpha=1}\dist^2(u^\alpha(x),\BC)\, \ext x + C\int_{B^n_{15/16}(0)}|u^\alpha(x)|^2\, \ext x.$$
	Using then \eqref{E:Lipschitz-approx-3}, this gives
	$$\sum_{\alpha:s_\alpha\neq 0}\int_{\Omega_{S_{i_*}}}|p^\alpha|^2 \leq C\int_{B_{15/16}(0)\setminus\Sigma}\sum^Q_{\alpha=1}\dist^2(u^\alpha(x),\BC)\, \ext x + C\hat{E}_V^2.$$
	Now using the area formula to write the first integral on the right-hand side as an integral over $V$, using also \eqref{E:J-bounds}, we get
	\begin{equation}\label{E:app-1-4}
	\sum_{\alpha:s_\alpha\neq 0}\int_{\Omega_{S_{i_*}}}|p^\alpha|^2 \leq C\hat{E}_{V,\BC}^2 + \hat{E}_V^2.
	\end{equation}
	Set $I := \{\alpha:s_\alpha\neq 0\}$; clearly by definition we have $I\neq\emptyset$. As $|\Omega_{S_{i_*}}|\geq c(n,Q)>0$, \eqref{E:app-1-4} then implies the lemma.
\end{proof}

\textbf{Remark:} An analogous result is true if we replace the $Q$ planes $P^1,\dotsc,P^Q$ forming $\BC$ by a cone $\BC = \sum^{2Q}_{i=1}|H_i|$ formed of a union of $2Q$ half-planes $H_1,\dotsc,H_{2Q}$ meeting along a common $(n-1)$-dimensional space, with $Q$ half-planes on `each side' (i.e.~so as $\nu(\BC)\to 0$, $\BC$ is converging to $Q|P_0|$ as varifolds). Suppose we have ordered the half-planes so that $H_1,\dotsc,H_Q$ are on one side and $H_{Q+1},\dotsc,H_{2Q}$ are on the other. The conclusion would then be that there are non-empty subsets $I_1\subset\{1,\dotsc,Q\}$ and $I_2\subset\{Q+1,\dotsc,2Q\}$ for which
$$\sum_{i\in I_1\cup I_2}\nu(H_i) \leq C(\hat{E}_{V,\BC} + \hat{E}_V)$$
for some $C = C(n,k,Q)$. To prove this, we simply need to notice that on `each side' of the spine of $\BC$, one sees only a union of planes. More precisely, given any ball $B_\rho(x_0)\setminus\pi_{P_0}(S(\BC))$, we can apply Lemma \ref{lemma:appendix-1} with $\tilde{V} := (\eta_{x_0,\rho})_\#V \res C_1(0)$ in place of $V$ and with a collection of planes $\tilde{P}_1,\dotsc,\tilde{P}_Q$ in place of $P_1,\dotsc,P_Q$ which are defined (up to relabelling) by the stipulation that $\spt\|\BC\|\cap C_\rho(x_0) = \sum^Q_{\alpha=1}\tilde{P}^\alpha$.

Our second lemma provides a criterion for controlling the excess relative to one plane by the excess to a union of two planes.

\begin{lemma}\label{lemma:appendix-2}
	There exists $\delta_0 = \delta_0(n,k)\in (0,1)$ such that the following holds. Suppose that we have a smooth function $f:P_0\cap B^n_{2\sqrt{n}}(0)\to P_0^\perp$ and a cone $\BC = |P_0| + |P| \in \mathcal{P}$ with
	\begin{enumerate}
		\item [(1)] $\nu(\BC)<\delta_0$;
		\item [(2)] $\|f\|_{C^1(B^n_{2\sqrt{n}}(0))}<\delta_0 a^*(\BC)$.
	\end{enumerate}
	(Note that the domain of $f$ is one of the planes in $\BC$.) Then
	\begin{equation}\label{E:appendix-2}
		\int_{\graph(f)\cap (\R^k\times Q_0)}\dist^2(x,P_0)\, \ext\H^n(x) \leq C\int_{\graph(f)\cap (\R^k\times Q_0)}\dist^2(x,\BC)\, \ext\H^n(x),
	\end{equation}
	where $Q_0 = [-1,1]^n\subset P_0$. Here, $C = C(n,k)\in (0,\infty)$.
\end{lemma}

\begin{proof}
	We first reduce the lemma to a one-dimensional problem. Write $p:P_0\to P_0^\perp$ for the linear function with $\graph(p) = P$. Then, let $e_n\in P_0$ denote a unit vector for which $|p(e_n)| = a^*(\BC)$. Choosing coordinates so that $P_0$ is spanned by $e_1,\dotsc,e_n$, set $S:= \textnormal{span}\langle e_1,\dotsc,e_{n-1}\rangle$. Write $Q_0$ for the cube $[-1,1]^n\subset P_0$.
	
	Since the gradient of $f$ is bounded, the area formula and \eqref{E:J-bounds} give
	$$\int_{\graph(f)\cap (\R^k\times Q_0)}\dist^2(x,P_0)\, \ext\H^n(x) \leq C\int_{Q_0}|f|^2$$
	and, also using that $\nu(\BC)<\delta_0$,
	\begin{align*}
	\int_{\graph(f)\cap (\R^k\times Q_0)}\dist^2(x,\BC)\, \ext\H^n(x) & \geq \int_{Q_0}\dist^2((f(x),x),\BC)\, \ext x\\
	& \geq C^{-1}\int_{Q_0}\min\{|f(x)-p(x)|^2,|f(x)|^2\}\, \ext x
	\end{align*}
	where here $C = C(n,k)$. Then, by the coarea formula we have
	$$\int_{Q_0}|f(x)|^2\, \ext x = \int_{S\cap Q_0}\int^1_{-1}|f(y+te_n)|^2\, \ext t\,\ext\H^{n-1}(y)$$
	and
	\begin{align*}
	\int_{Q_0}\min\{&|f(x)-p(x)|^2,|f(x)|^2\}\, \ext x\\
	& = \int_{S\cap Q_0}\int^1_{-1}\min\{|f(y+te_n)-p(y+te_n)|^2,|f(y+te_n)|^2\}\, \ext t\,\ext\H^{n-1}(y).
	\end{align*}
	Combining all of the above, it would therefore suffice to prove that for every $y\in S\cap Q_0$, we have
	\begin{equation}\label{E:app-2-1}
		\int^1_{-1}|f(y+te_n)|^2\, \ext t \leq C\int^1_{-1}\min\{|f(y+te_n)-p(y+te_n)|^2\, |f(y+te_n)|^2\}\, \ext t.
	\end{equation}
	So fix a point $y_0\in S\cap Q_0$. As $p$ is linear we know that $p(y_0+te_n) = p(y_0) + tp(e_n)$. Set $w:= p(e_n)$; we know that $|w| = a^*(\BC)$. Now, since $w\in \R^k$, we may \emph{choose} an orthonormal basis of $\R^k$, say $\{e_{n+1},\dotsc,e_{n+k}\}$, for which
	$$w\cdot e_{n+\kappa} = a^*(\BC)/\sqrt{k}$$
	for each $\kappa\in\{1,\dotsc,k\}$ (i.e.~so that $w = (a^*(\BC)/\sqrt{k},\dotsc,a^*(\BC)/\sqrt{k})$ in this basis). Clearly we then have, for $f^\kappa = f\cdot e_{n+\kappa}$, $p^\kappa = p\cdot e_{n+\kappa}$,
	$$|f(y_0+te_n)|^2 = \sum^k_{\kappa=1}|f^\kappa(y_0+te_n)|^2,$$
	$$|f(y_0+te_n)-p(y_0+te_n)|^2 = \sum^k_{\kappa=1}|f^\kappa(y_0+te_n) - p^\kappa(y_0+te_n)|^2,$$
	and thus since it is always true for non-negative $a_1,\dotsc,a_k$ and $b_1,\dotsc,b_k$ that $\min\{\sum_i a_i,\sum_i b_i\}\geq \sum_i \min\{a_i,b_i\}$, we have that
	\begin{align*}
	\min\{|f(y_0+te_n)-p(y_0+te_n)|^2&,\, |f(y_0+te_n)|^2\}\\
	& \geq\sum^{k}_{\kappa=1}\min\{|f^\kappa(y_0+te_n)-p^\kappa(y_0+te_n)|^2,\, |f^\kappa(y_0+te_n)|^2\}.
	\end{align*}
	Therefore, to prove \eqref{E:app-2-1} it suffices to prove
	\begin{equation}\label{E:app-2-2}
		\int^1_{-1}|f^\kappa(y_0+te_n)|^2\, \ext t \leq C\int^1_{-1}\min\{|f^\kappa(y_0+te_n)-p^\kappa(y_0+te_n)|^2,\, |f^\kappa(y_0+te_n)|^2\}\, \ext t
	\end{equation}
	for each $\kappa\in \{1,\dotsc,k\}$, where by assumption we know that
	$$\|f^\kappa(y_0+(\cdot)e_n)\|_{C^1([-1,1])} < \delta a^*(\BC)$$
	and
	$$p^\kappa(y_0+te_n) = p^\kappa(y_0) + tw\cdot e_{n+\kappa}$$
	is such that $w\cdot e_{n+\kappa} = a^*(\BC)/\sqrt{k}$. Now fix $\kappa\in \{1,\dotsc,k\}$ and let us abuse notation by writing
	$$f(t):= f^\kappa(y_0+te_n),\ \ \ \ p(t):= p^\kappa(y_0+te_n),\ \ \ \ \text{and}\ \ \ \ \mu = w\cdot e_{n+\kappa}.$$
	Then, in this notation, to prove the lemma we just need to prove that if $f,p:[-1,1]\to \R$ are functions with
	\begin{itemize}
		\item $p(t) = p_0 + t\mu$;
		\item $\|f\|_{C^1([-1,1])}<\delta \mu$,
	\end{itemize}
	for some $\mu>0$, then provided we choose $\delta\in (0,1)$ sufficiently small, we have
	\begin{equation}\label{E:app-2-3}
		\int^1_{-1}|f(t)|^2\, \ext t \leq C\int^1_{-1}\min\{|f(t)-p(t)|^2,\, |f(t)|^2\}\, \ext t
	\end{equation}
	for some $C = C(n,k)\in (0,\infty)$. In fact, dividing both sides by $\mu^2$, we can assume that $\mu=1$ (replacing $p_0/\mu$ by $p_0$), and then dividing both sides by $\delta$, we can instead assume that $p(t) = p_0 + Mt$, $\|f\|_{C^1([-1,1])} < 1$, where $M = 1/\delta>0$. Thus to summarise again, we are now assuming that $p(t) = p_0 + Mt$, $\|f\|_{C^1([-1,1])}<1$, and we want to show \eqref{E:app-2-3} holds for $M>0$ sufficiently large.
	
	First, notice that $|f(t)-p(t)|^2 = |f(t)|^2$ if and only if $2f(t)p(t) = p(t)^2$, i.e.~if $p(t) = 0$ or $2f(t) = p(t)$. Being linear, $p(t) = 0$ has at most one solution. We claim that, if $M>2$, the second equation has at most one solution. Indeed, if $t_1,t_2$ are both solutions to $2f(t) = p(t)$, then by the mean value theorem
	$$M|t_1-t_2| = |p(t_1)-p(t_2)| = 2|f(t_1)-f(t_2)| \leq 2 \cdot \|f\|_{C^1}\cdot |t_1-t_2| \leq 2\cdot |t_1-t_2|$$
	which if $M> 2$ implies that $t_1=t_2$. Thus, the set of points where $|f(t)-p(t)|^2 = |f(t)|^2$ is at most $2$. It then follows (from $p$ being linear, for example) that the region where $|f(t)-p(t)|^2 < |f(t)|^2$ is a single, open, interval, say $I$. We therefore just need to show
	\begin{equation}\label{E:app-2-4}
	\int_I |f(t)|^2\, \ext t \leq C\int^1_{-1}\min\{|f(t)-p(t)|^2,|f(t)|^2\}\, \ext t.
	\end{equation}
	Of course, if $I=\emptyset$ then \eqref{E:app-2-4} follows trivially, and so we may suppose that at least one solution to $|f(t)-p(t)|^2 = |f(t)|^2$ occurs.
	
	Notice that we must have $|I|$ small, as $\sup|f|<1$ and from the form of $p$ the interval has length at most $<2/M$. The interval must therefore be either $[-1,a)$, $(a,b)$, or $(b,1]$ for some $a,b$. We will split the proof into two cases, depending into how many solutions to $|f(t)-p(t)|^2 = |f(t)|^2$ occur.
	
	\textbf{Case 1:} \emph{Both solutions occur.} In this case, we have $I = (a,b)$. Let us suppose that $a$ is such that $p(a) = 0$ (the other case will follow be an analogous argument). Therefore, $2f(b) = p(b)$ ($>0$ as $b>a$ as $p$ is increasing). In particular, from the form of $b$ we must have
	$$M(b-a) = 2f(b) \qquad \Longrightarrow \qquad (b-a) = \frac{2f(b)}{M}.$$
	Take any $t\in I$. Then we have by the mean value theorem, using $|f|_{C^1}<1$,
	$$|f(t) - f(b)| \leq 1\cdot |t-b| \leq b-a = \frac{2f(b)}{M}$$
	and so $|f(t)| \leq f(b) + \frac{2f(b)}{M} = f(b)\left(1+\frac{2}{M}\right)$. Hence,
	$$\int_I |f(t)|^2\, \ext t \leq (b-a)\cdot f(b)^2\left(1+\frac{2}{M}\right)^2 = f(b)^3\cdot\frac{2}{M}\left(1+\frac{2}{M}\right)^2.$$
	We now wish to provide a lower bound for $\int_I \min\{|f(t)-p(t)|^2, |f(t)|^2\}\, \ext t \equiv \int_I |f(t)-p(t)|^2\, \ext t$. Notice that since $p(b) = 2f(b)$ and $p$ is linear, we must have $p(t) = 2f(b) + M(t-b)$. We also know from $|f|_{C^1}<1$ that $f(t) \leq f(b) + (b-t)$ for all $t<b$; set $g(t) = f(b) + (b-t)$. Then, we have $g(t) = p(t)$ at $x=c$, where
	$$2f(b) + M(c-b) = f(b) + (b-c) \qquad \Longrightarrow \qquad c = b - \frac{f(b)}{M+1}.$$
	Therefore, on the interval $[c,b]$, we have that $|f(t)-p(t)| = p(t) - f(t) \geq p(t) - [f(b)+(b-t)] \equiv p(t)-g(t)$. Hence,
	\begin{align*}
	\int_I |f(t) - p(t)|^2\, \ext t \geq \int^b_c |p(t)-g(t)|^2\, \ext t & = \int^b_c |(M+1)(b-t)-f(b)|^2\, \ext t\\
	& = \int^{\frac{f(b)}{M+1}}_0 |(M+1)t - f(b)|^2\, \ext t\\
	& = (M+1)^2\int^{\frac{f(b)}{M+1}}_0\left|t-\frac{f(b)}{M+1}\right|^2\, \ext t\\
	& = \frac{f(b)^3}{3(M+1)}.
	\end{align*}
	Thus, to prove \eqref{E:app-2-4} in this case it suffices choose $C$ such that
	$$f(b)^3\cdot\frac{2}{M}\left(1+\frac{2}{M}\right)^2 \leq C\cdot\frac{f(b)^3}{3(M+1)}$$
	i.e.
	$$\frac{(M+1)(M+2)^2}{M^3} \leq \frac{C}{6}.$$
	Of course, as the right-hand side has a limit of $1$ as $M\to\infty$ (in fact, for $M>7$ it is $<2$), and thus this is clearly possible if $M>7$ and $C>12$. This concludes the proof of \eqref{E:app-2-4} in Case 1.
	
	\textbf{Note:} In this case we actually prove $\int_I |f(t)|^2\, \ext t \leq C\int_I |f(t)-p(t)|^2$.
	
	\textbf{Case 2:} \emph{Only one solution occurs.} In this case, $I = [-1,a)$ or $I = (b,1]$. We will look at the case $I = (b,1]$, as the proof is analogous in both cases. Furthermore, we have two choices depending on whether $b$ is the zero of $p$ or the solution of $2f = p$; we will focus first on the former, as again the latter case is analogous.
	
	First, we claim that if $M>6$, then we must have $f(1)\geq 3(1-b)$. Indeed, if not, then since $p(b)=0$, we have $M(1-b)=p(1) - p(b) = p(1)$, and as there is no solution to $2f = p$ we must have $p(1)<2f(1)$. Thus, we get $M(1-b) < 2f(1)$. But then if $f(1)<3(1-b)$, this would imply $M<6$, giving the contradiction.
	
	Thus, we must have $f(1)\geq 3(1-b)$, i.e.~the size of $f$ is larger than a fixed proportion of the size of $I$. On $[b,1]$ we necessarily have (as $|f|_{C^1}<1$) $f(t) \leq f(1) + (1-b)$, and so
	$$\int_I |f(t)|^2\, \ext t \leq (1-b)\cdot|f(1)+(1-b)|^2.$$
	But on $[b-(1-b),b] \equiv [2b-1,b]$, we also know (again, from $|f|_{C^1}<1$ and the mean value theorem) that $f(t)\geq f(1) - 2(1-b)$. So,
	$$\int_{[-1,1]\setminus I}|f(t)|^2 \geq \int^b_{2b-1}|f(1)-2(1-b)|^2 = (1-b)|f(1)-2(1-b)|^2.$$
	Note that on $[-1,1]\setminus I$, the minimum of $|f-p|^2$ and $|f|^2$ is $|f|^2$ by definition. So, to prove \eqref{E:app-2-4} it suffices to show
	$$|f(1)+(1-b)|^2 \leq C|f(1)-2(1-b)|^2.$$
	Now, since we have $f(1)\geq 3(1-b)$ and $M(1-b)<2f(1)$, we see that
	$$|f(1)+ (1-b)|^2 \leq \left|f(1) + \frac{2}{M}f(1)\right|^2 = f(1)^2\left|1+\frac{2}{M}\right|^2$$
	and
	$$|f(1)-2(1-b)|^2 \geq f(1)^2\left|1-\frac{4}{M}\right|^2$$
	and thus in fact we just need $(1+2/M)^2 \leq C(1-4/M)^2$; for $M\geq 8$, we can simply take $C = 16$ and thus we would get \eqref{E:app-2-4} in this case also. Thus, this completes the proof of the lemma.
	
	\textbf{Note:} In Case 2 we actually prove $\int_I |f(t)|^2\, \ext t \leq C\int_{[-1,1]\setminus I}|f(t)|^2\, \ext t$.
\end{proof}

\textbf{Remark:} We see from the above proof that we actually only need to assume that $f$ is Lipschitz and that $\|f\|_{C^{0,1}}<\delta_0 a^*(\BC)$. Furthermore, the same proof also works for multi-valued Lipschitz functions: indeed, in the reduction step when we reduce to the coordinate functions, arguing similarly as in Section \ref{sec:constructing-coarse-blow-ups}, we can work with an ordering of the given ($\A_Q(\R)$-valued) coordinate functions where each function is Lipschitz, at which point the problem reduces to the same one-dimensional, single-valued, version as seen in the above proof.

\part{\centering Regularity of Coarse Blow-Ups and Conclusion}\label{part:blow-up-reg}

In this section, we return to the discussion of coarse blow-ups started in Part \ref{part:coarse-blow-ups}. There, we constructed and discussed properties of coarse blow-ups of arbitrary sequences of stationary integral varifolds converging to a plane with any multiplicity. We will now restrict ourselves to the case where the plane has multiplicity $2$, as well as only to coarse blow-ups which arise from blowing-up sequences of varifolds in $\mathcal{V}_{\beta}$, for some $\beta\in (0,1)$. Since for a fixed $\beta$ the class $\mathcal{V}_\beta$ is not closed under rotations, and thus coarse blow-ups of varifolds in $\mathcal{V}_\beta$ are not closed under the operation in $(\mathfrak{B}3\text{III})$ (all other properties in Theorem \ref{thm:blow-up-properties} do hold for this restricted class of coarse blow-ups, however) it will be convenient to work with a slightly larger class of coarse blow-ups, namely those which arise from blowing-up a sequence of varifolds in $\mathcal{V}_{\beta_j}$, where $\beta_j\uparrow \beta$. This is then a subclass of the full (multiplicity $2$) coarse blow-up class $\mathfrak{B}_2$; we will abuse notation slightly and denote this subclass by $\mathfrak{B}_\beta$. The reader should note two things:
\begin{itemize}
	\item \emph{All} the properties $(\mathfrak{B}1)-(\mathfrak{B}6)$ seen in Theorem \ref{thm:blow-up-properties} hold for the subclass $\mathfrak{B}_\beta$, and so in particular the coarse blow-ups whose existence is claimed in $(\mathfrak{B}3)$ will all belong to $\mathfrak{B}_\beta$. This follows from the proof of Theorem \ref{thm:blow-up-properties}.
	\item $\mathcal{V}_\beta$ is contained within the class of coarse blow-ups generated by sequences of varifolds in $\mathcal{V}_{\beta/2}$ converging to a multiplicity $2$ plane. In particular, we may apply results from the rest of the paper which are true for $\mathcal{V}_{\beta/2}$, such as Theorem \ref{thm:fine-reg}, along a blow-up sequence.
\end{itemize}
Using the analysis in Part \ref{part:fine-reg}, which concluded with the fine $\eps$-regularity theorem given in Theorem \ref{thm:fine-reg}, we will be able to conclude a full regularity result for coarse blow-ups in $\mathfrak{B}_\beta$. In turn, the regularity of these coarse blow-ups will then allow us to prove Theorem \ref{thm:main}.

Throughout this part of the paper, $\beta\in (0,1)$ will be fixed. As some extra notation, for $v\in \mathfrak{B}_\beta$, write $\CC_v$ for the set of $x\in B_1(0)$ such that $v(x)$ is a classical singularity of $v$ (by which we mean $\graph(v)$ is locally a $C^1$ perturbation of a cone in $\CC$), and write $\mathcal{B}_v:= \Gamma_v\setminus \CC_v$.

\section{Further Properties of Coarse Blow-Ups in $\mathfrak{B}_\beta$}

We start by deducing further properties of the coarse blow-up class $\mathfrak{B}_\beta$ in addition to $(\mathfrak{B}1)-(\mathfrak{B}6)$, which hold specifically due to the blow-ups arising from sequences of varifolds which obey the topological structural condition in certain coarse gaps. Recall that $\Gamma_v$ is the (relatively closed) subset of $B^n_1(0)$ consisting of the Hardt--Simon points of the coarse blow-up $v$, as defined in Definition \ref{defn:coarse-blow-up}.

\begin{theorem}\label{thm:further-blow-up-properties}
	Fix $\beta\in (0,1)$. Then the blow-up class $\mathfrak{B}_\beta$ satisfies the following properties:
	\begin{enumerate}
		\item [$(\mathfrak{B}7)$] \textnormal{(Decomposition away from $\Gamma_v$.)} Let $v\in \mathfrak{B}_\beta$. Then for any $x_0\in B^n_1(0)\setminus\Gamma_v$ and $0<r<\min\{\dist(x_0,\Gamma_v),1-|x_0|\}$, we have that $v|_{B^n_r(x_0)}$ decomposes as the sum of two smooth harmonic functions $v_1,v_2:B^n_r(x_0)\to \R^k$.
		\item [$(\mathfrak{B}8)$] \textnormal{(Hölder continuity.)} Let $v\in\mathfrak{B}_\beta$. Then (after redefining $v$ on a set of measure zero) we have that for any $\mu\in (0,1)$, $v|_{B_\sigma(0)}\in C^{0,\mu}(\overline{B_\sigma(0)};\A_2(\R^k))$ for each $\sigma\in (0,1)$, with the estimate
		$$\|v^{\kappa,\alpha}\|_{C^{0,\mu}(B_\sigma(0))}\leq C\|v\|_{L^2(B_1(0))}$$
		for each $\kappa\in \{1,\dotsc,k\}$ and $\alpha\in \{1,2\}$, where $C = C(n,k,\beta,\sigma,\mu)\in (0,\infty)$. Moreover, for each $z\in \Gamma_v$, we have $v^\alpha(z) = v_a(z)$ for $\alpha=1,2$.
		\item [$(\mathfrak{B}9)$] \textnormal{(Homogeneity continuation property.)} Let $v\in\mathfrak{B}_\beta$, and suppose that $v$ is homogeneous of degree $1$ on the annulus $B_1(0)\setminus \overline{B_r(0)}$ for some $r\in (0,1)$, i.e.~$\frac{\del(v/R)}{\del R}=0$ a.e.~in $B_1(0)\setminus \overline{B_r(0)}$, where $R(x) = |x|$. Then, $v$ is homogeneous of degree $1$ in $B_1(0)$.
		\item [$(\mathfrak{B}10)$] \textnormal{(Transverse $\eps$-regularity.)} Let $\gamma\in (0,1)$. Then, there exists $\eps = \eps(n,k,\beta,\gamma)\in (0,1)$ such that whenever $v\in\mathfrak{B}_\beta$ satisfies $0\in\Gamma_v$, $\|v\|_{L^2(B_1)}=1$, and
		$$\int_{B_1}\G(v,\psi)^2 < \eps$$
		for some $\psi:\R^n\to \A_2(\R^k)$ with $\psi_a\equiv 0$ and $\mathbf{v}(\psi)\in \CC$, then we have the following:
		\begin{enumerate}
			\item [(i)] $v$ is generalised-$C^{1,\mu}$ in $B_{1/2}(0)$. Moreover, $\mathcal{B}_v\cap B_{1/2}(0) = \emptyset$, $\CC_v\cap B_{1/2}\subset \graph(\phi)\cap B_{1/2}(0)$ for some function $\phi:S_0\cap B_{1/2}(0)\to S_0^\perp\subset \R^{k+1}$ of class $C^{1,\mu}$ with $|\phi|_{1,\mu;S_0\cap B_{1/2}(0)}\leq C\gamma$. Furthermore, if $\Omega^\pm$ denote the two connected components of $B_{1/2}^n(0)\setminus\graph(\phi_{k+1})$, then $v|_{\Omega^{\pm}}\in C^{1,\mu}(\overline{\Omega^{\pm}})$ and $|v|_{1,\mu;\Omega^{\pm}}\leq C$. 
			\item [(ii)] There is a function $\widetilde{\psi}:\R^n\to \A_2(\R^k)$ with $\mathbf{v}(\widetilde{\psi})\in \CC$ such that
			$$\dist_\H(\graph(\widetilde{\psi})\cap C_1, \graph(\psi)\cap C_1) \leq C\gamma$$
			and
			$$\sigma^{-n-2}\int_{B_\sigma(0)}\G(v(x),\widetilde{\psi}(x))^2\, \ext x \leq C\gamma\sigma^{2\mu} \qquad \text{for all }\sigma\in (0,1/2).$$
		\end{enumerate}
		Here, $C = C(n,k,\beta)\in (0,\infty)$ and $\mu = \mu(n,k,\beta)\in (0,1)$.
		\item [$(\mathfrak{B}11)$] \textnormal{(Squeeze identity at transverse points.)} Let $v\in \mathfrak{B}_\beta$ and suppose $x_0\in\Gamma_v$ is such that there exists a tangent map $\psi$ to $v$ at $x_0$ with $\mathbf{v}(\psi)\in \CC$. Then, there exists $r = r(v,x_0)\in (0,\frac{1}{2}(1-|x_0|))$ such that the following is true.
		\begin{enumerate}
			\item [(i)] $\Gamma_v\cap B^n_r(x_0)$ is contained within an $(n-1)$-dimensional $C^1$ submanifold of $B_r^n(x_0)$;
			\item [(ii)] $\Gamma_v\cap B^n_r(x_0) = \CC_v\cap B_r(x_0)$;
			\item [(iii)] For any $\zeta\in C^1_c(B_r(x_0);\R^n)$,
			$$\int_{B_r(x_0)}\sum_{\alpha=1,2}\sum_{\kappa=1}^k\sum_{i,j=1}^n\left(|Dv^{\kappa,\alpha}|^2\delta_{ij} - 2D_iv^{\kappa,\alpha}D_jv^{\kappa,\alpha}\right)D_i\zeta^j = 0,$$
			and
			$$\int_{B_r(x_0)}\sum_{\alpha=1,2}\sum_{\kappa=1}^k\sum_{i,j=1}^n\left(|Dv_f^{\kappa,\alpha}|^2\delta_{ij}-2D_iv_f^{\kappa,\alpha}D_jv_f^{\kappa,\alpha}\right)D_i\zeta^j = 0.$$
		\end{enumerate}
	\end{enumerate}
\end{theorem}

\textbf{Remark:} The homogeneity in $(\mathfrak{B}9)$ can also be understood using the continuity from $(\mathfrak{B}8)$.

We will prove each property in turn.

\begin{proof}[Proof of $(\mathfrak{B}7)$]
	Fix $x_0\in B^n_1(0)\setminus\Gamma_v$ and pick $r,r^\prime\in (0,\min\{\dist(x_0,\Gamma_v),1-|x_0|\})$ with $r^\prime>r$. By the construction of coarse blow-ups and the definition of $\Gamma_v$, we see that for all $j$ sufficiently large, if $(V_j)_j\subset\mathcal{V}_\beta$ is the sequence whose blow-up is $v$, we must have $\Theta_{V_j}<2$ in $C_{r^\prime}(x_0)$. Thus, again for all $j$ sufficiently large (depending on $r,r^\prime$) we can apply Lemma \ref{lemma:gap-2} to see that $V_j\res C_{(r+r^\prime)/2}(x_0)$ is equal to the sum of two single-valued smooth minimal graphs, with elliptic estimates which control their $C^{4}$ norm by $\hat{E}_{V_j}$. Thus, when we blow-up $V_j$ in this region, Arzelà--Ascoli can be applied to see that $v|_{B^n_r(x_0)}$ is equal to the sum of two smooth harmonic functions, proving $(\mathfrak{B}7)$.
\end{proof}

\begin{proof}[Proof of $(\mathfrak{B}8)$]
	This proof follows that seen in \cite[Lemma 3.2]{MW24}. For $\delta>0$, let $\eta_\delta\in C^1([0,\infty))$ be a non-decreasing function such that $\eta_\delta(t) = 0$ for $t\in [0,\delta/2]$, $\eta_\delta(t)=1$ for $t\in [\delta,\infty)$, and $|D\eta_\delta|\leq 3/\delta$ for all $t\in [0,\infty)$. Let $\phi\in C^1([0,\infty))$ be such that $\phi(t) = 1$ for $t\in[0,1/4]$, $\phi(t) = 0$ for $t\in [3/8,\infty)$, and $|\phi^\prime(t)|\leq 10$ for all $t\in [0,\infty)$. Note that since $v^{\kappa,\alpha}\in W^{1,2}_{\text{loc}}(B_1(0))$ for each $\kappa\in \{1,\dotsc,k\}$ and $\alpha\in \{1,2\}$, we have that $|v-v_a(z)|^2\in W^{1,1}_{\text{loc}}(B_1(0))$ for any fixed $z\in B_1(0)$. For any fixed $\gamma>0$, fixed $z\in B_{1/2}(0)\cap \Gamma_v$, and any fixed $\rho\in (0,1/4)$, set $X^i(x):=\phi(R_z/\rho)^2\eta_\delta(R_z)R_z^{-n+\gamma-2}|v(x)-v_a(z)|^2(x^i-z^i)$, where $R_z = |x-z|$. Then by the divergence theorem we have $\int_{\R^n}D_iX^i = 0$, whence
	\begin{align}
		\nonumber\gamma\int\phi(R_z/\rho)^2\eta_\delta(R_z)R_z^{-n+\gamma-2}|v-v_a(z)|^2 & = -\int\phi(R_z/\rho)^2\eta_\delta(R_z)R_z^{1-n+\gamma}\frac{\del(R_z^{-2}|v-v_a(z)|^2)}{\del R_z}\\
		& -2\int\rho^{-1}\phi(R_z/\rho)\phi^\prime(R_z/\rho)\eta_\delta(R_z)R_z^{-n+\gamma-1}|v-v_a(z)|^2\nonumber\\
		& -\int\phi(R_z/\rho)^2\eta_\delta^\prime(R_z)R_z^{-n+\gamma-1}|v-v_a(z)|^2.\label{E:B8-1}
	\end{align}
	Noting that $\frac{\del(R_z^{-2}|v-v_a(z)|^2)}{\del R_z} = \sum_{\kappa=1}^k\sum_{\alpha=1,2}2R_z^{-1}(v^{\kappa,\alpha}-v_a^\kappa(z))\frac{\del((v^{\kappa,\alpha}-v_a^\kappa(z))/R_z)}{\del R_z}$ a.e. in $B_1(0)$, we see, using the Cauchy--Schwarz inequality (in the form $2ab\leq \eps a^2 + \eps^{-1}b^2$) and $(\mathfrak{B}5)$, that this combined with \eqref{E:B8-1} implies
	\begin{align*}
		\int&\phi(R_z/\rho)^2\eta_\delta(R_z)R_z^{-n+\gamma-2}|v-v_a(z)|^2\\
		& \leq C\int\left(\phi(R_z/\rho)^2R_z^{2-n+\gamma}\left|\frac{\del}{\del R_z}\left(\frac{v-v_a(z)}{R_z}\right)\right|^2 + \rho^{-2}\phi^\prime(R_z/\rho)^2R_z^{-n+\gamma}|v-v_a(z)|^2\right)\\
		& \leq C\rho^{-n-2+\gamma}\int_{B_\rho(z)}|v-v_a(z)|^2
	\end{align*}
	where $C = C(n,k,\gamma)\in (0,\infty)$, and where we have discarded the last integral on the right-hand side of \eqref{E:B8-1} (note that we have used $(\mathfrak{B}3\text{III})$ with $\lambda = v_a(z)$ and $L\equiv 0$ to get this right-hand side in $(\mathfrak{B}5)$). Letting $\delta\downarrow 0$ and using the monotone convergence theorem, we deduce from this that
	$$\int_{B_{\rho/4}(z)}\frac{|v-v_a(z)|^2}{R_z^{n+2-\gamma}}\leq C\rho^{-n-2+\gamma}\int_{B_\rho(z)}|v-v_a(z)|^2$$
	for every $\rho\in (0,1/4)$ and every $z\in B_{1/2}\cap \Gamma_v$, whence
	$$\sigma^{-n}\int_{B_\sigma(z)}|v-v_a(z)|^2 \leq C\left(\frac{\sigma}{\rho}\right)^{2-\gamma}\rho^{-n}\int_{B_\rho(z)}|v-v_a(z)|^2$$
	for every $\sigma,\rho$ with $0<\sigma\leq\rho/4\leq 1/16$ and every $z\in B_{1/2}(0)\cap \Gamma_v$. Since by $(\mathfrak{B}7)$ we know that $v$ is harmonic in $B_1(0)\setminus\Gamma_v$, the desired conclusion, with $\sigma=1/2$, with the values of $v$ for each $z\in B_{1/2}(0)\cap \Gamma_v$ given by $v^\alpha(z) = v_a(z)$ for $\alpha=1,2$, follows from this, with $\mu = 1-\frac{\gamma}{2}$, and an appropriate version of the Campanato lemma (e.g.~as in \cite[Lemma 4.3]{Wic14} or \cite{Min21b}). The conclusion for $\sigma\in (1/2,1)$ then follows by applying the conclusion for $\sigma=1/2$ to appropriately translated and rescaled $v$ and using a covering argument, which is possible by $(\mathfrak{B}3\text{I})$.
\end{proof}

\begin{proof}[Proof of $(\mathfrak{B}9)$]
	This proof follows that seen in \cite[Lemma 3.3]{MW24}. First note the following general fact: if $\Omega\subset\R^n$ is a connected, open set, and if $u:\Omega\to \R$ is a harmonic function which is homogeneous of degree $1$ on some (non-empty) open subset of $\Omega$, then $u$ is homogeneous of degree $1$ on $\Omega$. Indeed, for $x\in\Omega\setminus\{0\}$, define $w:=\frac{\del(u/R)}{\del R} \equiv \frac{x\cdot Du - u}{R^2}$, where $R(x) = |x|$, and note that $\Delta(R^2w) = \Delta(x\cdot Du) - \Delta u = 0$ in $\Omega\setminus\{0\}$. Since $w\equiv 0$ on some open subset of $\Omega$, it follows from unique continuation for harmonic functions that $R^2w\equiv 0$ on $\Omega\setminus\{0\}$, and hence $w\equiv 0$ on $\Omega\setminus\{0\}$. Thus, $u$ is homogeneous of degree $1$ on $\Omega$.
	
	Now to prove the lemma, note that since $v_a$ is harmonic on $B_1$ (by $(\mathfrak{B}2)$) and homogeneous of degree $1$ in $B_1\setminus B_r$, we have that $v_a$ is homogeneous of degree $1$ in $B_1$ (and hence in fact is linear). Set $U:= \{x\in B_r: |v(x)-v_a(x)|\neq 0\}$. By continuity of $v$ (provided by $(\mathfrak{B}8)$) we know that $U$ is open. Moreover, for $x\in U$ we must have $v^1(x)\neq v^2(x)$, and so by $(\mathfrak{B}8)$ it follows that $x\not\in \Gamma_v$, and so by $(\mathfrak{B}7)$ we know that there is $\rho_x>0$ such that $v|_{B_{\rho_x}(x)}$ is harmonic. If $U = \emptyset$, then $v^\alpha = v_a$ on $B_r$ for $\alpha=1,2$ and so the lemma follows (as $v_a$ is homogeneous of degree $1$). So suppose that $U\neq\emptyset$, and let $\Omega$ be any connected component of $U$. If $v^\alpha(x) = v_a(x)$ for $\alpha=1,2$ and each $x\in\del \Omega$ (boundary taken in $B_1$), then since $v^\alpha$ is harmonic on $\Omega$ for each $\alpha$, it follows from the weak maximum principle that $v^\alpha\equiv v_a$ on $\overline{\Omega}$ for each $\alpha$, which is impossible since $\Omega\subset U$. So, there must be a point $x_0\in\del\Omega$ such that $v^\alpha(x_0)\neq v_a(x_0)$ for some $\alpha=1,2$. Now, it is not possible that $x_0\in B_r$ (for if $x_0\in B_r$, then $x_0\in U$ and hence, since $\Omega$ is a component, $x_0\not\in\del\Omega$). Thus, $x_0\in \del B_r\cap \del\Omega$. Then since $\Omega_0\equiv \Omega\cup B_{\rho_{x_0}}(x_0)$ is connected, $v^\alpha$ is harmonic in $\Omega_0$ for each $\alpha$, and (by assumption) $v^\alpha$ is homogeneous of degree $1$ in $B_{r_{x_0}}(x_0)\setminus B_r$, it follows from the above general fact that $v^\alpha$ is homogeneous of degree $1$ in $\Omega$ for each $\alpha$. This is then true for every component of $U$, and so $v$ is homogeneous of degree $1$ in $U$. We also have, since $v^\alpha = v_a$ on $B_r\setminus U$ for each $\alpha=1,2$, that $\frac{\del(v^\alpha/R)}{\del R} = 0$ a.e. on $B_r\setminus U$ for each $\alpha=1,2$. Thus, $\frac{\del(v/R)}{\del R} = 0$ a.e. on $B_r$. But for any function $h\in C^1(\overline{B}_r)$, any $s\in (0,r)$, and any $\rho,\sigma\in [s,r]$, we have that
	$$\left|\frac{h(\rho\w)}{\rho}-\frac{h(\sigma\w)}{\sigma}\right|\leq s^{1-n}\int^r_s t^{n-1}\left|\frac{\ext(h(t\w)/t)}{\ext t}\right|\, \ext t$$
	for any $\w\in S^{n-1}$, whence integrating this with respect to $\w$,
	$$\int_{S^{n-1}}\left|\frac{h(\rho\w)}{\rho} - \frac{h(\sigma\w)}{\sigma}\right|\, \ext\w \leq s^{1-n}\int_{B_r\setminus B_s}\left|\frac{\del(h/R)}{\del R}\right|.$$
	By an approximation argument, this extends to any $h\in C^0(\overline{B_r})\cap W^{1,2}(B_r)$, and so taking $h = v^{\kappa,\alpha}$, we conclude that $\frac{v(\rho\w)}{\rho} = \frac{v(\sigma \w)}{\sigma}$ for every $\w\in S^{n-1}$ and every $\rho,\sigma\in (0,r]$, i.e.~$v$ is homogeneous of degree $1$ on $B_r$ and hence on $B_1$.
\end{proof}

\begin{remark}\label{remark:holder-continuity}
	Up until this point, other than stationarity the only additional property we have used along the blow-up sequence is that on a $(\beta/2)$-coarse gap of \emph{fixed} size (possibly being very small), we have the topological structural condition. In particular, we can prove Hölder continuity of such coarse blow-ups always without the fine $\eps$-regularity. Thus, the place where we need the fine $\eps$-regularity analysis is in controlling the energy about parts of the coarse blow-up which look like cones in $\CC$, in particular establishing $(\mathfrak{B}11)$ about such points. If one could do this, then one could possibly prove the regularity for coarse blow-ups of this type.
\end{remark}

\begin{proof}[Proof of $(\mathfrak{B}10)$]
	The proof follows that seen in \cite[Theorem 3.4]{MW24}, except phrased as a contradiction argument instead. Suppose that $(\mathfrak{B}10)$ were false. Then, we could find $\gamma\in (0,1)$ as well as sequences $v_j\in\mathfrak{B}_\beta$ with $0\in \Gamma_{v_j}$, $\|v_j\|_{L^2(B_1(0))}=1$, and $\psi_j:\R^n\to \A_2(\R^k)$ with $(\psi_j)_a\equiv 0$ and $\mathbf{v}(\psi_j)\in \CC$, with $\int_{B_1}\G(v_j,\psi_j)\to 0$, yet the conclusions fail. We may pass to a subsequence such that $\graph(\psi_j)$ converges locally in Hausdorff distance to $\graph(\psi)$, where $\psi:\R^n\to \A_2(\R^k)$ is such that $\graph(\psi)\in \CC$ or $\graph(\psi)$ is a plane; since $(\psi_j)_a\equiv 0$ we would also need $\psi_a\equiv 0$, and so in the latter case we would have $\psi\equiv 0$. But notice that since $\|v_j\|_{L^2(B_1)}=1$ for all $j$ and $v_j\to \psi$ in $L^2(B_1(0))$, it cannot be the case that $\psi\equiv 0$, and so we must have $\psi\in \CC$. Applying now $(\mathfrak{B}3\text{IV})$, we can also suppose that $v_j\to \psi$ in the sense described in $(\mathfrak{B}3\text{IV})$. Of course, by $(\mathfrak{B}3\text{II})$ we can rotate everything to assume that $0\in S(\psi) := S(\mathbf{v}(\psi))\subset S_0$.
	
	Suppose now that $v_j$ is the coarse blow-up of the sequence $(V_{j,\ell})_\ell\subset\mathcal{V}_{\beta_j}$, where $\beta_j\uparrow \beta$ (a priori, the parameter $\beta_{j}$ will depend on $\ell$ also, but using Remark \ref{remark:inclusion} we can ensure this parameter is independent of $\ell$, up to passing to a subsequence in $\ell$). In particular, for all $j$ sufficiently large we have $V_{j,\ell}\in \mathcal{V}_{\beta/2}$ for all $\ell$.
	
	Let $u_{j,\ell}$ be the Lipschitz approximation of $V_{j,\ell}$ whose blow-up is $v_j$. By assumption, we know that for every $\delta>0$, for all $j$ sufficiently large we have
	$$\int_{B_1(0)}\G(v_j(x),\psi(x))^2\, \ext x < \delta^2.$$
	Additionally, given $\sigma\in (0,1)$ we can, for each $j$, choose $\ell_j$ sufficiently large (depending on $\delta$ and $\sigma$) so that for all $\ell\geq \ell_j$ we have
	$$\int_{B_\sigma(0)}\G(\hat{E}_{V_{j,\ell}}^{-1}u_{j,\ell}(x),v_j(x))^2\, \ext x < \delta^2.$$
	Let us write $W_j:= V_{j,\ell_j}$ and $u_j:= u_{j,\ell_j}$, and $\BC_j:= \graph(\hat{E}_{W_j}\psi)$. By the triangle inequality we have
	\begin{equation}\label{E:B10-0}
	\int_{B_\sigma(0)}\G(\hat{E}_{W_j}^{-1}u_j(x),\psi(x))^2\, \ext x <2\delta^2.
	\end{equation}
	In particular, we can assume that the coarse blow-up of $(W_j)_j$ is $\psi$. By the Hardt--Simon inequality, this means that $\mathcal{D}(W_j)\cap C_1(0)$ must converge in Hausdorff distance to a (closed) subset of $S(\psi)$.
	
	We claim that for all $j$ sufficiently large we have
	\begin{equation}\label{E:B10-1}
		\hat{E}_{W_j}\leq\frac{3}{2}\mathscr{E}_{W_j,\BC_j}.
	\end{equation}
	Indeed, if this were to fail then there is a subsequence we can pass to (without changing notation) and a sequence of planes $P_j\supsetneq S(\BC_j)$ such that
	\begin{equation}\label{E:B10-2}
	\hat{E}_{W_j,P_j}<\frac{2}{3}\hat{E}_{W_j}.
	\end{equation}
	In particular, $P_j\neq P_0$. Moreover, arguing as we have done many times throughout the paper, controlling the Hausdorff distance between two planes by the $L^2$ norm of the difference between the linear functions representing them over a suitable subset and then using the triangle inequality along with an appropriate graphical representation, it is clear that there must be a constant $C = C(n,k)\in (0,\infty)$ for which $\nu(P_j)\leq C\hat{E}_{W_j}$, and so we may further assume that $\hat{E}_{W_j}^{-1}p_j\to L$ for some linear function $L:\R^n\to \R^k$; here, $p_j$ is the linear function over $P_0$ with graph $P_j$, and so in particular $p_j\equiv 0$ on $S(\BC_j)$, meaning that $S(\BC_j)\subset\{L=0\}$. Dividing \eqref{E:B10-2} by $\hat{E}_{W_j}$ and taking $j\to\infty$, we see that
	$$\int_{B_1(0)}|\psi(x)-L(x)|^2\, \ext x \leq \frac{2}{3}.$$
	(Technically, we pass to the limit the integral on the left-hand side over $B_\sigma(0)$, for some $\sigma\in (1/2,1)$, and then let $\sigma\to 1$ to arrive at this.)
	Since $\psi_a\equiv 0$, we have $|\psi(x)-L(x)|^2 = |\psi(x)|^2 + 2|L(x)|^2  - 2L(x)\sum_{\alpha=1,2}\psi^\alpha(x) = |\psi(x)|^2 + 2|L(x)|^2$, and thus this becomes
	$$\int_{B_1}|\psi|^2 + 2\int_{B_1}|L|^2 \leq \frac{2}{3}$$
	but since $\|\psi\|_{L^2(B_1)}=1$, this is a contradiction. Hence we have established \eqref{E:B10-1}, which is assumption (d) of Theorem \ref{thm:fine-reg} with $M=\frac{3}{2}$ (and $\beta/2$ in place of $\beta$).
	
	By arguing analogously to the justification of \cite[(14.33)]{Wic14} or \cite[(3.40)]{MW24} also shows that given $\gamma\in (0,1)$, if we choose $\sigma\in (0,1)$ sufficiently close to $1$ and then $j$ sufficiently large (depending also on $\gamma$ and $\sigma$) then we have
	$$\hat{E}_{W_j,\BC_j} \leq \frac{\gamma}{2}\hat{E}_{W_j}.$$
	Moreover, as $\mathcal{D}(W_j)\cap C_1(0)$ converges in Hausdorff distance to a subset of $S(\psi)$ (as mentioned above), we may apply Lemma \ref{lemma:gap-2} away from any fixed neighbourhood of $S(\psi)$ for all sufficiently large $j$. In particular, we know that $W_j$ decomposes as two single-valued smooth minimal graphs in the region $B_{1/2}(0)\setminus\{r_{\BC_j}<1/8\}$, and we can conclude from these observations that 
	$$\int_{B_{1/2}(0)\setminus\{r_{\BC_j}<1/8\}}\dist(x,W_j)\, \ext\|\BC_j\|(x) = \int_{B_{1/2}(0)\setminus\{r_{\BC_j}<1/8\}}\dist^2(x,\graph(u_j))\, \ext\|\BC_j\|(x).$$
	By elliptic estimates, we know that $Du_j\to 0$ uniformly in $B^n_{1/2}(0)\setminus\{r_{\BC_j}<1/8\}$, which ultimately means that the integral above is bounded by
	$$C\int_{B_{1/2}(0)}\G(\hat{E}_{W_j}\psi(x),u_j(x))^2\, \ext x.$$
	Invoking \eqref{E:B10-0}, this is then controlled by $C\delta^2\hat{E}_{W_j}^2$, and thus we have therefore proved that (choosing $\delta$ sufficiently small depending on $\gamma$)
	$$\hat{F}_{W_j,\BC_j}\leq \gamma\hat{E}_{W_j}.$$
	Clearly $W_j,\BC_j\to 2|P_0|$ in $C_1(0)$, which means that all the other hypotheses of Theorem \ref{thm:fine-reg} are easily verifiable for $W_j$ and $\BC_j$ and all sufficiently large $j$. Thus, for suitable large $j$ we can apply Theorem \ref{thm:fine-reg} to $W_j$ and $\BC_j$. This means that we can in fact that $\graph(u_j)$ to be exactly equal to $V_{j,\ell_j}$, and indeed $\graph(u_{j,\ell})$ exactly equal to $V_{j,\ell}$ for all $\ell\geq \ell_j$. Passing the conclusions of Theorem \ref{thm:fine-reg} to the blow-up then gives the desired conclusions of $(\mathfrak{B}10)$ hold for $v_j$ for all sufficiently large $j$, providing the desired contradiction and thus proving $(\mathfrak{B}10)$.
\end{proof}

\begin{proof}[Proof of $(\mathfrak{B}11)$]
	Given our assumptions, after translating, `tilting', and rescaling (via $(\mathfrak{B}3\text{I})$ and $(\mathfrak{B}3\text{II})$) we can assume that we are in the setting of $(\mathfrak{B}10)$ for suitable $\psi$ therein, which we can then apply to see that the conclusions of $(\mathfrak{B}10)$ hold in some neighbourhood $B_r(x_0)$ of $x_0$, where $r = r(v,x_0)>0$. In particular, claims (i) and (ii) are immediate. Now let $(\widetilde{V}_j)_j\subset \mathcal{V}_{\beta_j}\subset\mathcal{V}_{\beta/2}$ be a sequence of varifolds generating $v$ (here, without loss of generality $\beta/2\leq\beta_j\uparrow\beta$). Now let $\Gamma\in SO(n+k)$ be the rotation $\Gamma(x):= \gamma(\pi_{P_0}(x))$, where $\gamma:\R^n\to\R^n$ is a rotation taking $S(\psi)$ to a subset of $S_0$. The proof of $(\mathfrak{B}10)$ then in fact shows that following: after possibly passing to a subsequence of $(\widetilde{V}_j)_j$, for any sequence of points $(z_j)_j$ with $z_j\to 0$ and $z_j\in\mathcal{D}(\widetilde{V}_j)$, if we set $V_j:= (\Gamma\circ\eta_{z_k,1-|z_k|})_\#\widetilde{V}_j$ and $\BC_j:= \mathbf{v}(\hat{E}_{V_j}\psi\circ\gamma)\in \CC$, then $v$ is the coarse blow-up of $(V_j)_j$ and for all sufficiently large $j$ the hypotheses of Theorem \ref{thm:fine-reg} are satisfied with $V_j$ and $\BC_j$ in place of $V$ and $\BC$, respectively. In particular, $V_j\res C_{1/2}(0)$ is a two-valued graph over $B_{1/2}^n(0)$ for all sufficiently large $j$.
	
	Now fix $\zeta\in C^1_c(B_{1/4}(0);\R)$. Taking the vector field $X = \widetilde{\zeta}e_p$, for $p\in\{1,\dotsc,n\}$, in the first variation formula for $V_k$, we get
	$$\int\nabla^{V_j}x^p\cdot\nabla^{V_j}\widetilde{\zeta}(x)\, \ext\|V_j\|(x) = 0.$$
	Since $\|V_j\|$-a.e. on $\spt\|V_j\|$ we have $\nabla^{V_j}x^p\cdot \nabla^{V_j}\widetilde{\zeta}(x) = \pi_{T_xV_j}(e_p)\cdot D\widetilde{\zeta}(x)$, a direct computation analogous to that leading to \eqref{E:fine-b-10} shows that
	$$\int \nabla^{V_j}x^p\cdot \nabla^{V_j}\widetilde{\zeta}(x)\, \ext\|V_j\|(x) = \int_{B_1^n(0)}\sum_{\alpha=1,2}\sum^n_{q=1}(A^\alpha_j(x))_{qp}D_q\zeta(x) J^\alpha_j(x)\, \ext x$$
	where $A^\alpha_j(x)$ is the $n\times n$ matrix that is the inverse of the matrix $B^\alpha_j$ given by $(B^\alpha_j)_{pq} = \delta_{pq}+\sum^k_{\kappa=1}D_pu^{\kappa,\alpha}_jD_q u^{\kappa,\alpha}_j$, i.e.~we have that (suppressing the dependence on $\alpha$ and $j$ momentarily): 
	$$\sum^n_{i=1}(\delta_{qi}+D_q u^\alpha\cdot D_i u^\alpha)A_{ip} = \delta_{qp}.$$
	This leads to:
	$$A_{qp} = \delta_{qp} - D_qu^\alpha\cdot D_pu^\alpha + \boldsymbol{\eps}$$
	where the contribution to the integral given by the term denoted by $\boldsymbol{\eps}$ is $o(\hat{E}_j^2)$. Using this computation and the fact that $\int D_p\zeta(x)\, \ext x = 0$, we get
	$$\int\sum_{\alpha=1,2}\left((J^\alpha_j(x)-1)D_p\zeta(x) - \sum^n_{q=1}D_qu^\alpha\cdot D_pu^\alpha D_q\zeta(x) J^\alpha_j(x)\right)\, \ext x = o(\hat{E}_j^2).$$
	If we divide this by $\hat{E}_j^2$ and take $j\to\infty$, using the $C^{1,\mu}$ estimates for $u^\alpha$ on $\Omega^\pm$ provided in Theorem \ref{thm:fine-reg} in order to pass the derivatives pointwise to the limit, we end up with
	$$\int\sum_{\alpha=1,2}\sum^n_{q=1}\left(|Dv^\alpha|^2\delta_{pq} - 2D_qv^\alpha\cdot D_p v^\alpha\right)D_q\zeta(x)\, \ext x = 0.$$
	Now, given any $\zeta = (\zeta^1,\dotsc,\zeta^n)\in C^1(B_{1/4}(0);\R^n)$, if we set $\zeta = \zeta^p$ in the above and then sum over $p$, we get the first claimed identity of $v$. The second identity for $v_f$ follows from this, using that $v_f = v-v_a$ and that $v_a$ is harmonic by $(\mathfrak{B}2)$, as for harmonic functions the claimed identity is trivially satisfied by an integration by parts.
\end{proof}

Next, we will extend property $(\mathfrak{B}11)$ to show that in fact the squeeze identity holds for coarse blow-ups in $\mathfrak{B}_\beta$ within any region in which the coarse blow-up is known to be $GC^1$. This will require the energy non-concentration estimate, $(\mathfrak{B}6)$. We remark that, in fact in the present setting a version of this estimate is implied by the continuity estimate $(\mathfrak{B}8)$, the fact that $v\in\mathfrak{B}_\beta$ is harmonic whenever $|v(x)-v_a(x)|>0$ by $(\mathfrak{B}7)$, as well as the energy estimate which follows from $(\mathfrak{B}4)$. Indeed, we have:

\begin{lemma}[Alternative energy non-concentration estimate]\label{lemma:non-con-alt}
	Let $v\in \mathfrak{B}_\beta$ and $\delta>0$. Then for every $\rho\in (0,1)$ and $\kappa\in\{1,\dotsc,k\}$ we have
	$$\int_{B_\rho(0)}\sum_{\alpha=1,2}\one_{\{|v^{\kappa,\alpha}_f|<\delta\}}(x)|Dv^{\kappa,\alpha}_f(x)|^2\, \ext x \leq C(1-\rho)^{-2}\delta\left(\int_{B_1}|v^\kappa_f|^2\right)^{1/2}$$
	where $C = C(n,k)\in (0,\infty)$.
\end{lemma}

\begin{proof}
	We follow \cite[Lemma 3.7]{MW24}. Fix $\delta>0$. For $\eps<\delta/2$, define $\eta_\eps:\R\to \R$ to be the odd extension to $\R$ of the function $t\mapsto \delta_{(\eps,\infty)}(t)\cdot\min\left\{\frac{t-\eps}{\delta-\eps},1\right\}$ for $t>0$. Let $\sigma\in (0,1)$ and $\zeta\in C^1_c(B_\sigma)$. By $(\mathfrak{B}8)$, $v$ is continuous on $B_1$ and so for each $\alpha=1,2$ we have by $(\mathfrak{B}7)$ that $v^{\kappa,\alpha}_f$ is harmonic on the open set $\{|v^\alpha_f|>0\}$, and so by an approximation argument we may take $\eta_\eps(v^{\kappa,\alpha}_f)\zeta$ as a test function in the weak harmonicity of $v^{\kappa,\alpha}_f$, i.e.
	$$\int_{B_1}Dv_f^{\kappa,\alpha}\cdot D(\eta_\eps(v^{\kappa,\alpha}_f)\zeta)=0.$$
	Upon rearranging this becomes
	$$\int_{B_1}\zeta\cdot\eta^\prime_\eps(v^{\kappa,\alpha}_f)\cdot|Dv^{\kappa,\alpha}_f|^2 = -\int_{B_1}\eta_\eps(v^{\kappa,\alpha}_f)Dv^{\kappa,\alpha}_f\cdot D\zeta$$
	from which we get
	$$\frac{\delta}{\delta-\eps}\int_{B_1}\zeta\cdot\one_{\{\eps<|v^{\kappa,\alpha}_f|<\delta\}}\cdot|Dv^{\kappa,\alpha}_f|^2 \leq \w_n^{1/2}\sup|\eta_\eps|\sup_{B_\sigma}|D\zeta|\left(\int_{B_\sigma}|Dv_f^{\kappa,\alpha}|^2\right)^{1/2}.$$
	By $(\mathfrak{B}4)$ (cf.~the first point in Remark \ref{remark:after-blow-up-properties} after Theorem \ref{thm:blow-up-properties}) we have $\int_{B_\sigma}|Dv^\kappa_f|^2 \leq C(1-\sigma)^{-2}\int_{B_1}|v^\kappa_f|^2$, where $C = C(n,k)$, and so we get from the above, after summing over $\alpha$ and letting $\eps\to 0$ and noting that $Dv^{\kappa,\alpha}_f = 0$ a.e. on $\{v^{\kappa,\alpha}_f=0\}$,
	$$\int_{B_1}\sum_{\alpha=1,2}\zeta\cdot\one_{\{|v^{\kappa,\alpha}_f|<\delta\}}\cdot|Dv^{\kappa,\alpha}_f|^2 \leq C\delta(1-\sigma)^{-1}\sup_{B_\sigma}|D\zeta|\left(\int_{B_1}|v^\kappa_f|^2\right)^{1/2}.$$
	For given $\rho\in (0,1)$, choosing $\sigma = (1+\rho)/2$ and $\zeta$ such that $\zeta=1$ on $B_\rho$ and $|D\zeta|\leq 2/(\sigma-\rho)$ everywhere, the desired conclusion follows.
\end{proof}

Let us now prove the squeeze identity for $GC^1$ coarse blow-ups in $\mathfrak{B}_\beta$.

\begin{lemma}[Squeeze identity for $GC^1$ coarse blow-ups]\label{lemma:squeeze}
	Let $v\in \mathfrak{B}_\beta$ and let $\Omega\subset B_{1/2}(0)$ be open. Suppose that $v$ is of class $GC^1$ in $\Omega$. Then, for every $\zeta\in C^1_c(\Omega;\R^n)$,
	$$\int_{\Omega}\sum^Q_{\alpha=1}\sum^n_{i,j=1}\left(|Dv^\alpha|^2\delta_{ij}-2D_iv^\alpha\cdot D_jv^\alpha\right)D_i\zeta^j = 0$$
	and
	$$\int_{\Omega}\sum^Q_{\alpha=1}\sum^n_{i,j=1}\left(|Dv^\alpha_f|^2\delta_{ij}-2D_iv^\alpha_f\cdot D_jv^\alpha_f\right)D_i\zeta^j = 0.$$
\end{lemma}

\begin{proof}
	We follow \cite[Lemma 3.9]{MW24}. The first assertion follows from the second by writing $v^\alpha = v^\alpha_f + v_a$, noting that $v^1_f+v^2_f = 0$ and that the claimed identity holds for $v_a$ (since it is harmonic by $(\mathfrak{B}2)$. Thus, we will focus on proving the second identity.
	
	Let $\zeta\in C^1_c(\Omega;\R^n)$. For $\delta>0$ small, consider a smooth function $\gamma_\delta:\R\to \R$ such that $\gamma_\delta(t) = t+\frac{3}{4}\delta$ for $t<-\delta$, $\gamma_\delta(t)\equiv 0$ for $|t|<\delta/2$, and $\gamma_\delta(t) = t-\frac{3}{4}\delta$ for $t>\delta$, in such a way that $|\gamma_\delta^\prime(t)|\leq 1$ and $|\gamma_\delta^{\prime\prime}(t)|\leq 3\delta^{-1}$ for all $t\in \R$.
	
	We know from $(\mathfrak{B}8)$ that $\Gamma_v\subset\Omega\subset\{|v_f|=0\}$ and that $\gamma_\delta(v^{\kappa,\alpha}_f)|_\Omega$ is a smooth function for each $\alpha=1,2$ and $\kappa=1,\dotsc,k$. Thus by direct calculation and the dominated convergence theorem we have (using summation convention)
	\begin{align*}
		\lim_{\delta\to 0}\int_{B_1(0)}\sum_{\alpha=1,2}\sum_{\kappa=1}^k(|D\gamma_\delta(v^{\kappa,\alpha}_f)|^2\delta_{ij}& - 2D_i\gamma_\delta(v^{\kappa,\alpha}_f)D_j\gamma_\delta(v^{\kappa,\alpha}_f))D_i\zeta^j\\
		& = \int_{B_1(0)}\sum_{\alpha=1,2}(|Dv^\alpha_f|^2\delta_{ij}-2D_iv^\alpha_f\cdot D_jv^\alpha_f)D_i\zeta^j.
	\end{align*}
	We now compute the limit on the left-hand side in an alternative way. Let $\eps>0$. Since $v$ is generalised-$C^1$ in $\Omega$, we know that $|Dv_f|=0$ on $\mathcal{B}_v\cap\Omega$. For each $x\in \CC_v\cap \Omega$, let $\rho_x>0$ denote the radius from $(\mathfrak{B}11)$. Since $\spt(\zeta)$ is compact, we can pick finitely many points $x_1,\dotsc,x_N\in \CC_v\cap \Omega$ for which
	$$\spt(\zeta)\cap \CC_v \cap B_\eps(\mathcal{B}_v) \subset\bigcup^N_{i=1}B_{\rho_{x_i}}(x_i).$$
	Then choose an open set $\mathcal{O}$ such that
	$$\dist(\mathcal{O},\CC_v\cup\mathcal{B}_v)>0 \qquad \text{and} \qquad \spt(\zeta)\subset \mathcal{O}\cup B_\eps(\mathcal{B}_v)\cup\bigcup^N_{i=1}B_{\rho_{x_i}}(x_i)$$
	(e.g.~take $\mathcal{O} = \Omega\setminus \overline{B_\eps(\mathcal{B}_v)\cup\cup^N_{i=1}B_{\rho_{x_i}}(x_i)}$.) Write $\mathcal{U}:= \{\mathcal{O}\}\cup \{B_\eps(\mathcal{B}_v)\}\cup \{B_{\rho_{x_i}}(x_i):i=1,\dotsc,N\}$; this is an open cover of $\spt(\zeta)$. Now let $(\phi_\beta)_{\beta\in A}$ be a smooth partition of unity subordinate to $\mathcal{U}$ (which of course depends on $\eps$ but is independent of $\delta$). Note that:
	\begin{enumerate}
		\item [(I)] The indexing set $A$ has cardinality $|A| = N+2<\infty$ (i.e.~a finite number depending on $\eps$);
		\item [(II)] $\mathcal{O}$ is contained in the set where $v_f$ is harmonic;
		\item [(III)] For each $\beta\in A$, there is some $B\in\mathcal{U}$ for which $\spt(\phi_\beta)$. This means that we can write $A$ as a disjoint union $A_\CC\cup A_{\mathcal{B}}\cup A_{\mathcal{O}}$ depending on the set in which the support of $\phi_\beta$ lies; indeed, $A_{\mathcal{O}} := \{\beta\in A:\spt(\phi_\beta)\subset\mathcal{O}\}$, $A_{\CC} := \{\beta\in A\setminus A_{\mathcal{O}}:\spt(\phi_\beta)\subset B_{\rho_{x_i}}(x_i)\text{ for some }i\in\{1,\dotsc,N\}\}$, and $A_{\mathcal{B}}:= A\setminus(A_{\mathcal{O}}\cup\mathcal{A}_{\CC})$.
		\item [(IV)] $\sum_{\beta\in A}\phi_\beta = 1$.
	\end{enumerate}
	Since $v$ is generalised-$C^1$ and $|Dv_f|=0$ on $\mathcal{B}_v\cap\Omega$, it follows that
	\begin{equation}\label{E:squeeze-1}
		\sup_{B_{\eps}(\mathcal{B}_v)\cap \spt(\zeta)\cap \CC^c_v}|Dv_f|\to 0\qquad \text{as }\eps\to 0.
	\end{equation}
	Now consider the integral:
	\begin{equation}\label{E:squeeze-2}
		\sum_{B_1(0)}\sum_{\alpha=1,2}\sum_{\kappa=1}^k(|D\gamma_\delta(v^{\kappa,\alpha}_f)|^2\delta_{ij} - 2D_i\gamma_\delta(v_f^{\kappa,\alpha})D_j\gamma_\delta(v_f^{\kappa,\alpha}))D_i(\phi_\beta\zeta^j)\, \ext x.
	\end{equation}
	If $\beta\in A_{\mathcal{B}}$, then by integrating by parts (which we can do as $\gamma_\delta(v^{\kappa,\alpha}_f)$ is smooth) we have that this is equal to
	\begin{align*}
	-\int_{B_1(0)}\sum_{\alpha=1,2}\sum_{\kappa=1}^k\left[2D_\ell\gamma_\delta(v^{\kappa,\alpha}_f)\right.&D_{i\ell}\gamma_\delta(v^{\kappa,\alpha}_f)\delta_{ij}\\
	&\left. - 2\Delta\gamma_\delta(v^{\kappa,\alpha}_f)D_j\gamma_\delta(v^{\kappa,\alpha}_f) - 2D_i\gamma_\delta(v^{\kappa,\alpha}_f)D_{ij}\gamma_\delta(v^{\kappa,\alpha}_f)\right]\phi_\beta\zeta^j.
	\end{align*}
	Notice that the first and third terms in the above cancel (as we are summing over the indices). Since pointwise on $\spt(\gamma_\delta^\prime)$ we have $\Delta\gamma_\delta(v^\alpha_f) = \gamma_{\delta}^{\prime\prime}(v^{\kappa,\alpha}_f)|Dv^{\kappa,\alpha}_f|^2$, we get:
	\begin{align*}
		\sum_{\beta\in A_{\mathcal{B}}}\int_{B_1(0)}\sum_{\alpha=1,2}\sum_{\kappa=1}^k&(|D\gamma_\delta(v^{\kappa,\alpha}_f)|^2\delta_{ij} - 2D_i\gamma_\delta(v^{\kappa,\alpha}_f)D_j\gamma_\delta(v^{\kappa,\alpha}_f))D_i(\phi_\beta\zeta^j)\\
		& = 2\int_{B_1(0)}\sum_{\alpha=1,2}\sum_{\kappa=1}^k\gamma_{\delta}^{\prime\prime}(v^{\kappa,\alpha}_f)|Dv^{\kappa,\alpha}_f|^2D\gamma_\delta(v^{\kappa,\alpha}_f)\cdot\zeta\left(\sum_{\beta\in A_{\mathcal{B}}}\phi_\beta\right).
	\end{align*}
	Now as $\spt(\gamma_\delta^{\prime\prime})\subset [-\delta,\delta]$ and $\spt(\phi_\beta)\subset B_\eps(\mathcal{B}_v)$ for $\beta\in A_{\mathcal{B}}$, the integral on the right hand side above (for each $\alpha$ and $\kappa$) only takes place over the set $\{|v^{\kappa,\alpha}_f|<\delta\}\cap B_\eps(\mathcal{B}_v)$. Thus, since $|\gamma_\delta^\prime|\leq 1$, we can estimate it, in absolute value, by:
	$$\leq 2\sup_\R|\gamma_\delta^{\prime\prime}|\sup_{\spt(\zeta)\cap B_\eps(\mathcal{B}_v)\cap \CC_v^c}\sum_{\alpha=1,2}|Dv^\alpha_f|\cdot\sup|\zeta|\cdot\int_{B_{1/2}(0)}\sum_{\alpha=1,2}\sum_{\kappa=1}^k\one_{\{|v^{\kappa,\alpha}_f|<\delta\}}|Dv^{\kappa,\alpha}_f|^2.$$
	So, using the energy non-concentration estimate (either from $(\mathfrak{B}6)$ or Lemma \ref{lemma:non-con-alt}), and the fact that $\sup|\gamma_\delta^{\prime\prime}|\leq 2\delta^{-1}$, we see that for $\delta$ sufficiently small,
	\begin{align*}
	&\left|\sum_{\beta\in A_{\mathcal{B}}}\int_{B_1(0)}\sum_{\alpha=1,2}\sum_{\kappa=1}^k(|D\gamma_\delta(v^{\kappa,\alpha}_f)|^2\delta_{ij} - 2D_i\gamma_\delta(v^{\kappa,\alpha}_f)D_j\gamma_\delta(v^{\kappa,\alpha}_f))D_i(\phi_\beta\zeta^j)\right|\\
	&\hspace{29em} \leq C\sup_{\spt(\zeta)\cap B_{\eps}(\mathcal{B}_v)\cap \CC_v^c}|Dv_f|
	\end{align*}
	where $C$ is independent of $\delta$ and $\eps$.
	
	Now suppose $\beta\in A_{\mathcal{O}}$. Then \eqref{E:squeeze-2} is equal to
	$$\int_{B_1(0)}\sum_{\alpha=1,2}\sum^k_{\kappa=1}|\gamma_\delta^\prime(v^{\kappa,\alpha}_f)|^2(|Dv_f^{\kappa,\alpha}|^2\delta_{ij} - 2D_iv^{\kappa,\alpha}_fD_jv_f^{\kappa,\alpha})D_i(\phi_\beta\zeta^j)$$
	which as $\delta\downarrow 0$ converges to
	$$\int_{B_1(0)}\sum_{\alpha=1,2}\sum^{k}_{\kappa=1}(|Dv^{\kappa,\alpha}_f|^2\delta_{ij} - 2D_iv^{\kappa,\alpha}_fD_jv^{\kappa,\alpha}_f)D_i(\phi_\beta\zeta^j).$$
	By integrating by parts and making a pointwise calculation for $x\in\mathcal{O}$, this is equal to $0$ (since the $v^{\kappa,\alpha}_f$ are harmonic on $\mathcal{O}$).
	
	Finally, suppose $\beta\in A_{\CC}$ in \eqref{E:squeeze-2}. Then we have
	\begin{align*}
		&\int(|D\gamma_\delta(v^{\kappa,\alpha}_f)|^2\delta_{ij}-2D_i\gamma_\delta(v^{\kappa,\alpha}_f)D_j\gamma_\delta(v^{\kappa,\alpha}_f))D_i(\zeta^j\phi_\beta)\\
		& \hspace{10em} = \int|\gamma^\prime_\delta(v^{\kappa,\alpha}_f)|^2(|Dv^{\kappa,\alpha}_f|^2\delta_{ij}-2D_iv^{\kappa,\alpha}_fD_jv^{\kappa,\alpha}_f)D_i(\zeta^j\phi_\beta)\\
		& \hspace{10em} \to \int(|Dv^{\kappa,\alpha}_f|^2\delta_{ij}-2D_iv^{\kappa,\alpha}_fD_jv^{\kappa,\alpha}_f)D_i(\zeta^j\phi_\beta)
	\end{align*}
	as $\delta\to 0$, and moreover by $(\mathfrak{B}11)$ this final expression equals zero.
	
	Thus, summing the above three calculations over the different $\beta\in A$, then letting $\delta\to 0$ first, and then letting $\eps\to 0$, we get the claimed identity.
\end{proof}

For $v\in\mathfrak{B}_\beta$ that is generalised-$C^1$ on an open set $\Omega\subset B_1$, and for $B_\rho(y)\subset\Omega$, we define the \emph{frequency function} by:
$$N_{v,y}(\rho):= \frac{\rho^{2-n}\int_{B_\rho(y)}|Dv|^2}{\rho^{1-n}\int_{\del B_\rho(y)}|v|^2}.$$
This is well-defined whenever $v$ is not identically zero on $B_\rho(y)$; indeed, if $\int_{\del B_\rho(y)}|v|^2=0$, then the squash inequality from $(\mathfrak{B}4)$ would imply that $|Dv|=0$ on $B_\rho(y)$, and hence $v\equiv 0$ on $B_\rho(y)$.

By combining the squash inequality $(\mathfrak{B}4)$ and the squeeze identity (Lemma \ref{lemma:squeeze}) we can prove in a standard fashion (see e.g.~\cite{Alm00} or \cite{SW16}) that, where defined, $N_{v,y}(\rho)$ is a non-decreasing function of $\rho$. In particular for each $y\in \Omega$ either there is a $\rho>0$ for which $v|_{B_\rho(y)}\equiv 0$, of the limit
$$N_v(y):= \lim_{\rho\downarrow0}N_{v,y}(\rho)$$
exists in $[0,\infty)$. When defined, we call $N_v(y)$ the \emph{frequency} of $v$ at $y$. It is also a standard consequence of the monotonicity of $N_{v,y}(\cdot)$ that if $\Omega$ is connected, then either $v\equiv 0$ in $\Omega$ or $v\not\equiv 0$ on any ball $B_\rho(y)\subset\Omega$. Hence, if $\Omega$ is connected, then $N_{v,y}$ is well-defined for each $y\in \Omega$ and $\rho\in (0,\dist(y,\del\Omega))$, unless $v\equiv 0$ in $\Omega$. As the squash inequality and squeeze identity also hold for $v_f$, we can also analogously define $N_{v_f,y}(\rho)$ and $N_{v_f}(y)$, and the preceding facts hold with $v_f$ in place of $v$. Thus, we have the following:

\begin{theorem}[Monotonicity of the Frequency Function]\label{thm:frequency}
	Let $v\in\mathfrak{B}_\beta$ and suppose that $v$ is generalised-$C^1$ on a connected open set $\Omega\subset B_{1/2}$, and that $v\not\equiv 0$ in $\Omega$. Then, $v\not\equiv 0$ on each ball $B_\rho(y)\subset \Omega$, $N_{v,y}(\rho)$ is well-defined for each $y\in \Omega$ and $\rho\in (0,\dist(y,\del\Omega))$, and $N_{v,y}(\rho)$ is a monotonically non-decreasing function of $\rho$; in particular, the frequency $N_v(y)=\lim_{\rho\to 0}N_{v,y}(\rho)$ exists in $[0,\infty)$ for each $y\in\Omega$, and $N_v(y)$ is an upper semi-continuous function of $y$. Furthermore:
	\begin{enumerate}
		\item [(i)] for each $y\in\Omega$ and $0<\sigma\leq\rho<\dist(y,\del\Omega)$ we have
		$$\left(\frac{\sigma}{\rho}\right)^{2N_{v,y}(\rho)}\rho^{-n}\int_{B_\rho(y)}|v|^2 \leq \sigma^{-n}\int_{B_\sigma(y)}|v|^2 \leq \left(\frac{\sigma}{\rho}\right)^{2N_v(y)}\rho^{-n}\int_{B_\rho(y)}|v|^2;$$
		\item [(ii)] if $N_{v,y}(\rho)$ is constant for $\rho\in (\rho_1,\rho_2)$, then $v$ is a homogeneous function with respect to the variable $|x-y|$ on the interval $(\rho_1,\rho_2)$, with degree of homogeneity equal to the constant value of $N_{v,y}(\rho)$ on this interval.
	\end{enumerate}
	Moreover, the same conclusions hold with $v_f$ in place of $v$ whenever $v_f\not\equiv 0$ in $\Omega$.
\end{theorem}

\begin{remark}\label{remark:squash}
    Technically, the left-most inequality in Theorem \ref{thm:frequency}(i) requires one to have the squash \emph{equality}, rather than the squash inequality. However, in the setting of \ref{thm:frequency} one can prove the equality version of the squash inequality in an analogous manner to how we proved the squeeze identity. We will use this later.
\end{remark}

\section{Regularity of Coarse Blow-Ups in $\mathfrak{B}_\beta$}

We are now in a position to understand the regularity of coarse blow-ups in $\mathfrak{B}_\beta$. We start as we did for fine blow-ups, by classifying homogeneous degree one elements of $\mathfrak{B}_\beta$.

\begin{theorem}[Classification of Homogeneous Degree One Blow-Ups]\label{thm:classification}
	Suppose that $v\in\mathfrak{B}_\beta$ is homogeneous of degree $1$ in $B_1(0)$. Then, either $v$ is linear (i.e.~$v(x) = 2\llbracket \ell(x)\rrbracket$ for some linear function $\ell:\R^n\to \R^k$), or $\graph(v)\in \CC$.
\end{theorem}

\begin{proof}
	We follow \cite[Theorem 3.11]{MW24}. Note first that $v_a$ is harmonic (by $(\mathfrak{B}2)$) and homogeneous of degree $1$, and so is linear. Hence if $v=2\llbracket v_a\rrbracket$ then the conclusion holds. If not, then by $(\mathfrak{B}3\text{III})$ we have that $\|v-v_a\|_{L^2(B_1)}^{-1}(v-v_a)\in \mathfrak{B}_\beta$. So, it suffices to establish the theorem (with the conclusion that $\mathbf{v}(v)\in\CC$) for $v\in\widetilde{\mathfrak{B}}_\beta$, where
	$$\widetilde{\mathfrak{B}}_\beta := \left\{v\in\mathfrak{B}_\beta:\frac{\del(v/R)}{\del R}=0\text{ a.e. in }B_1,\ v_a\equiv 0,\, \|v\|_{L^2(B_1(0))} = 1\right\}.$$
	For each $v\in\widetilde{\mathfrak{B}}_\beta$, let $\widetilde{v}:\R^n\to\A_2(\R^k)$ denote the homogeneous degree $1$ extension of $v$ to $\R^n$. Denote by $S(\widetilde{v})$ the set of points $z\in \R^n$ such that $\widetilde{v}$ is invariant under translation by $z$ (i.e.~$\widetilde{v}(x+z) = \widetilde{v}(x)$ for every $x\in\R^n$), and note that by homogeneity of $\widetilde{v}$, we have that $S(\widetilde{v})$ is a linear subspace of $\R^n$. So, we can write $\widetilde{\mathfrak{B}}_\beta = \cup^n_{j=0}\mathcal{H}_j$, where $\mathcal{H}_j := \{v\in\widetilde{\mathfrak{B}}_\beta:\dim(S(\widetilde{v})) = j\}$. Note that $\mathcal{H}_n = \emptyset$, since $\|v\|_{L^2(B_1(0))} = 1$ and $v_a\equiv 0$ for each $v\in\widetilde{\mathfrak{B}}_\beta$. If $v\in\mathcal{H}_{n-1}$ then by homogeneity we have $\mathbf{v}(v)\in \CC$, and so the conclusion follows. So, to prove the theorem, we just need to verify this when $v\in \mathcal{H}_d$ for $d\leq n-2$. We will do this by (downward) induction.
	
	Fix an element $v\in\mathcal{H}_d$, where we are now assuming the conclusion for any $\mathcal{H}_{d_*}$ with $d_*>d$. We claim the following:
	
	\textbf{Claim:} \emph{For any compact subset $K\subset B_1\setminus S(\widetilde{v})$, there is $\eps = \eps(v,K)\in (0,\dist(K,\del B_1\cup S(\widetilde{v})))$ such that the following holds.} For each $z\in K\cap \Gamma_v$ and each $\rho\in (0,\eps]$, either:
	\begin{enumerate}
		\item [(a)] The conclusions of $(\mathfrak{B}10)$ hold on $B_{3\rho/8}(z)$; in particular, $v$ is a generalised-$C^{1,\mu}$ function in $B_{3\rho/8}(z)$, and there is a function $\psi_z:\R^n\to \A_2(\R^k)$ with $\graph(\psi_z)\in \CC$ such that for all $0<\sigma\leq 3\rho/8$,
		\begin{equation}\label{E:classification-1}
		\sigma^{-n-2}\int_{B_\sigma(z)}\G(v(x),\psi_z(x))^2\, \ext x \leq C\left(\frac{\sigma}{\rho}\right)^{2\mu}\cdot\rho^{-n-2}\int_{B_\rho(z)}|v|^2;
		\end{equation}
		\item [(b)] The \emph{reverse Hardt--Simon inequality} holds, i.e.
		\begin{equation}\label{E:classification-2}
			\sum_{\alpha=1,2}\int_{B_\rho(z)\setminus B_{\rho/2}(z)}R_z^{2-n}\left|\frac{\del}{\del R_z}\left(\frac{v^\alpha}{R_z}\right)\right|^2\geq \eps\rho^{-n-2}\int_{B_\rho(z)}|v|^2,
		\end{equation}
		where $R_z:=|x-z|$.
	\end{enumerate}
	Here, $C = C(n,k,\beta)\in (0,\infty)$ and $\mu = \mu(n,k,\beta)\in (0,1)$ are fixed independent of $v$, $K$.
	
	To prove this, we argue by contradiction analogously to that seen in the proof of Theorem \ref{thm:fine-classification}; the key technical difference (which is why we need now a dichotomy, utilising $(\mathfrak{B}10)$ through the fine $\eps$-regularity theorem) is that we can only subtract off linear functions for coarse blow-ups. For completeness, let us sketch the proof. Supposing the claim is not true (with $C$ to be chosen depending only on $n,k$), then for each $i=1,2,\dotsc$, there are numbers $\eps_i\downarrow 0$, points $z,z_i\in K\cap \Gamma_v$ with $z_i\to z$, and radii $\rho_i>0$ with $\rho_i\to 0$ such that assertion (a) with $\rho = \rho_i$ and $z=z_i$ fails for each $i$, and also
	\begin{equation}\label{E:classification-3}
		\sum_{\alpha=1,2}\int_{B_{\rho_i}(z_i)\setminus B_{\rho_i/2}(z_i)}R_{z_i}^{2-n}\left|\frac{\del}{\del R_{z_i}}\left(\frac{v^\alpha}{R_{z_i}}\right)\right|^2 < \eps_i\rho_i^{-n-2}\int_{B_{\rho_i}(z_i)}|v|^2.
	\end{equation}
	Set $w_i:= v_{z_i,\rho_i}$; note that by $w_i\in \mathfrak{B}_\beta$ by $(\mathfrak{B}3\text{I})$, and by $(\mathfrak{B}8)$ $w_i$ is continuous on $B_1$. By $(\mathfrak{B}3\text{IV})$ and $(\mathfrak{B}8)$, we can find a subsequence (which we pass to) and an element $w_*\in \mathfrak{B}_\beta$ such that $w_i\to w_*$ locally uniformly and locally weakly in $W^{1,2}$ on $B_1$. Moreover, as $v_a\equiv 0$, we have that $(w_*)_a\equiv 0$. The same one-dimensional integration argument as used in Theorem \ref{thm:fine-classification} gives that $w_*\not\equiv 0$, and \eqref{E:classification-3} implies that $w_*$ is homogeneous of degree $1$ in $B_1\setminus B_{1/2}$, and hence $(\mathfrak{B}9)$ gives that $w_*$ is homogeneous of degree $1$ in $B_1$. Thus, $w_*/\|w_*\|_{L^2(B_1(0))}$ (which belongs to $\mathfrak{B}_\beta$ by $(\mathfrak{B}3\text{III})$) is a non-zero element of $\widetilde{\mathfrak{B}}_\beta$, and thus must be in $\mathcal{H}_{d^*}$, for some $d^*\geq 1$. But by construction, we necessarily have that $S(\widetilde{w}_*)\supset S(\widetilde{v})\cup \{z\}$, and thus $d_*>d$; thus, our induction hypothesis gives that $\mathbf{v}(w_*)\in \CC$. But then we can apply $(\mathfrak{B}9)$, with $\psi = w_*$ and $v = (w_i)_{0,3/4}$, for all sufficiently large $i$ to conclude that alternative (a) above must hold with $z = z_i$, $\rho = \rho_i$, with the constant $C = C(n,k)\in (0,\infty)$ as in $(\mathfrak{B}10)$. This is contrary to our assumption, and thus we have shown that the claim holds.
	
	Combining \eqref{E:classification-2} with $(\mathfrak{B}5)$, we then get the following dichotomy: if $z\in K\cap \Gamma_v$ and $\rho\in (0,\eps]$, then either:
	\begin{enumerate}
		\item [(i)] The conclusions of $(\mathfrak{B}10)$ hold in $B_{3\rho/8}(z)$; in particular, $v$ is generalised-$C^{1,\mu}$ on $B_{3\rho/8}(z)$ and there is a function $\psi_z:\R^n\to \A_2(\R^k)$ with $\graph(\psi_z)\in \CC$ for which we have the estimate \eqref{E:classification-1}; or
		\item [(ii)] We have that \eqref{E:classification-2} holds and that
		$$\sum_{\alpha=1,2}\int_{B_{\rho/2}(z)}R_z^{2-n}\left|\frac{\del}{\del R_z}\left(\frac{v^\alpha}{R_z}\right)\right|^2 \leq \theta\sum_{\alpha=1,2}\int_{B_{\rho}(z)}R_z^{2-n}\left|\frac{\del}{\del R_z}\left(\frac{v^\alpha}{R_z}\right)\right|^2 $$
	\end{enumerate}
	where $\theta = \theta(v,K)\in (0,1)$. We claim that from this, the following dichotomy holds for each $z\in K\cap \Gamma_v$: either
	\begin{enumerate}
		\item [(I)] The conclusions of $(\mathfrak{B}10)$ hold in some neighbourhood of $z$, and moreover there is a function $\psi_z:\R^n\to \A_2(\R^k)$ with $\graph(\psi)\in \CC$ for which we have the estimate
		$$\rho^{-n-2}\int_{B_\rho(z)}\G(v(x),\psi_z(x))^2\, \ext x \leq C\rho^{2\gamma}\int_{B_{\eps}(z)}|v|^2$$
		for some $C = C(v,K)\in (0,\infty)$ and all $\rho\in (0,3\eps/8]$; or
		\item [(II)] We have that \eqref{E:classification-2} holds with $\rho = 2^{-i}\eps$ for each $i=1,2,\dotsc$, and hence
		$$\sum_{\alpha=12}\int_{B_\sigma(z)}R_z^{2-n}\left|\frac{\del}{\del R_z}\left(\frac{v^\alpha}{R_z}\right)\right|^2 \leq C\left(\frac{\sigma}{\rho}\right)^{2\gamma}\sum_{\alpha=1,2}\int_{B_{\rho}(z)}R_z^{2-n}\left|\frac{\del}{\del R_z}\left(\frac{v^\alpha}{R_z}\right)\right|^2$$
		for all $0<\sigma\leq\rho/2\leq\eps/2$.
	\end{enumerate}
	Here, $C = C(v,K)\in (0,\infty)$ and $\gamma = \gamma(v,K)\in (0,1)$. Indeed, for each fixed $z\in K\cap \Gamma_v$, the dichotomy (i) or (ii) above holds for $\rho = 2^{-i}\eps$, $i=0,1,2,\dotsc$. Let $I$ be the first time (i) holds, i.e.~the smallest integer $i\geq 0$ such that (i) holds with $\rho = 2^{-i}\eps$. If $I=0$ we have alternative (I); also, if $I\geq 1$, then iterating (ii) for $i=0,1,\dotsc,I-1$ and combining with the estimate provided in (i) as well as $(\mathfrak{B}5)$ and \eqref{E:classification-2}, we again have (I). If no such $I$ exists, i.e.~if (ii) always holds, then iterating (ii) for all $i$ gives alternative (II).
	
	Finally, in the case that (II) holds, we can again use $(\mathfrak{B}5)$ and \eqref{E:classification-2} in conjunction with the estimate in (II) to replace (II) with:
	\begin{enumerate}
		\item [(II)$^\prime$] For all $0<\sigma\leq\rho/2\leq\eps/4$ we have
		\begin{equation}\label{E:classification-4}
			\sigma^{-n-2}\int_{B_{\sigma}(z)}|v|^2\leq C\left(\frac{\sigma}{\rho}\right)^{2\gamma}\rho^{-n-2}\int_{B_\rho(z)}|v|^2.
		\end{equation}
	\end{enumerate}
	Note that the set of points $z\in \Gamma_v\in \text{int}(K)$ where alternative (I) holds is $\CC_v\cap \text{int}(K)$, and hence we have shown that for each point $z\in \Gamma_v\cap \text{int}(K)\setminus \CC_v$, the estimate \eqref{E:classification-4} holds. From this, it is straightforward to check (e.g.~using the Campanato lemma \cite[Lemma 4.3]{Wic14}) that $v|_{\text{int}(K)}$ is generalised-$C^{1,\gamma}$ in $\text{int}(K)$, where $\gamma = \gamma(v,K)\in (0,1)$. In particular, as $K\subset B_1\setminus S(\widetilde{v})$ was an arbitrary compact set, we see that $v$ is generalised-$C^1$ in $B_1\setminus S(\widetilde{v})$.
	
	We now claim that $\Gamma_v\subset S(\widetilde{v})$. Our method for showing this will rely on the frequency function, which now can be brought into play in view of the generalised-$C^1$ regularity of $v$ on $B_1\setminus S(\widetilde{v})$. Indeed, by Lemma \ref{lemma:squeeze} we have that the squeeze identity
	$$\int_{B_1(0)}\sum_{\alpha=1,2}(|Dv|^2\delta_{ij}-2D_iv^\alpha\cdot D_jv^\alpha)D_i\zeta^j\, \ext x = 0$$
	holds for all $\zeta\in C^1_c(B_1\setminus S(\widetilde{v});\R^k)$. Since $S(\widetilde{v})$ is a linear subspace of dimension at most $n-2$, it has zero $2$-capacity, and hence, since we also have that $Dv$ is bounded in $B_1\setminus S(\widetilde{v})$ by the homogeneity of $\widetilde{v}$ (recalling that $\widetilde{v}$ is translation invariant under $S(\widetilde{v})$), we may perform a standard excision argument to see that we can in fact take any $\zeta\in C^1_c(B_1(0);\R^k)$ in the squeeze identity. Armed now with the squeeze identity and the squash inequality $(\mathfrak{B}4)$ (in fact, by Remark \ref{remark:squash} we have the squash equality here), we can use standard arguments (see e.g.~\cite{SW16}) to show that (since $\widetilde{v}\not\equiv 0$) the frequency $N_{\widetilde{v}}(y)$ is well-defined at every point $y\in B_1(0)$, and in fact all the conclusions of Theorem \ref{thm:frequency} hold with $\widetilde{v}$ in place of $v$. In particular, as $\widetilde{v}$ is homogeneous of degree $1$, we have that $N_{\widetilde{v}}(0)=1$. It follows from upper semi-continuity of $N_{\widetilde{v}}$ and the homogeneity of $\widetilde{v}$ that $N_{\widetilde{v}}(y)\leq N_{\widetilde{v}}(0) = 1$ for all $y\in B_1(0)$. Moreover, frequency monotonicity and homogeneity of $\widetilde{v}$ give that if $N_{\widetilde{v}}(y) = N_{\widetilde{v}}(0)=1$, then $\widetilde{v}$ is translation invariant along directions parallel to $y$, i.e.~$y\in S(\widetilde{v})$, and so $S(\widetilde{v}) = \{y\in \R^n:N_{\widetilde{v}}(y) = N_{\widetilde{v}}(0)=1\}$. We also have that for each $y\in B_1(0)$ and $0<\sigma\leq\rho$:
	\begin{equation}\label{E:classification-5}
		\left(\frac{\sigma}{\rho}\right)^{2N_{\widetilde{v},y}(\rho)}\rho^{-n}\int_{B_\rho(y)}|\widetilde{v}|^2\leq\sigma^{-n}\int_{B_\sigma(y)}|\widetilde{v}|^2 \leq \left(\frac{\sigma}{\rho}\right)^{2N_{\widetilde{v}}(y)}\rho^{-n}\int_{B_\rho(y)}|\widetilde{v}|^2.
	\end{equation}
	Now, if $y\in \mathfrak{B}_v\cap (B_1\setminus S(\widetilde{v}))$, then (II)$^\prime$ must hold for some $\eps = \eps(y,v)>0$ (which can be taken to be the $\eps$ corresponding to some fixed compact set $K = K(y)\subset B_1\setminus S(\widetilde{v})$ with $y\in K$ in the argument leading to (I) and (II)$^\prime$). In particular, combining (II)$^\prime$ with \eqref{E:classification-5}, we have for fixed $\rho\in (0,\eps/2]$ and all $\sigma\in (0,\rho/2]$:
	$$\left(\frac{\sigma}{\rho}\right)^{2N_{\widetilde{v},y}(\rho)}\rho^{-n}\int_{B_\rho(y)}|\widetilde{v}|^2 \leq \sigma^{-n}\int_{B_\sigma(y)}|\widetilde{v}|^2 \leq \sigma^2\cdot C\left(\frac{\sigma}{\rho}\right)^{2\gamma}\rho^{-n-2}\int_{B_\rho(y)}|\widetilde{v}|^2$$
	and thu $\sigma^{1+\gamma-N_{\widetilde{v},y}(\rho)}\geq \widetilde{C}>0$ for some $\widetilde{C} = \widetilde{C}(v,y,\rho)$ and $\gamma = \gamma(v,y)\in (0,1)$. As we may take $\sigma\downarrow 0$, we see that $N_{\widetilde{v},y}(\rho)\geq 1+\gamma$. As $\rho\in (0,\eps/2)$ was arbitrary, we can then take $\rho\downarrow 0$ to see that $N_{\widetilde{y}}(y)\geq 1+\gamma>1 = N_{\widetilde{v}}(0)$, a contradiction. Thus, we must have that $\mathcal{B}_{\widetilde{v}}\subset S(\widetilde{v})$. So, if we have $\Gamma_v\not\subset S(\widetilde{v})$, then we can find $z\in B_1\setminus S(\widetilde{v})$ for which (I) holds. In particular, as each cone in $\CC$ is determined by linear functions, we readily deduce from the estimate in (I) that $\rho^{-n}\int_{B_\rho(z)}|v|^2\leq C\rho^2$ for all $\rho\in (0,\eps]$, and some $C$ independent of $\rho$, and so in the same way as above, we have $N_{\widetilde{v}}(z)\geq 1$, and thus $N_{\widetilde{v}}(z)=1$, meaning that $z\in S(\widetilde{v})$, which is a contradiction. Thus, we must have $\Gamma_v\subset S(\widetilde{v})$.
	
	Hence, from $(\mathfrak{B}7)$ we know that $\widetilde{v}$ is harmonic on $\R^n\setminus S(\widetilde{v})$. In particular, as $\widetilde{v}$ is continuous on $\R^n$ by $(\mathfrak{B}8)$, and $\dim(S(\widetilde{v}))\leq d \leq n-2$, thus giving that $S(\widetilde{v})$ has vanishing $2$-capacity, we see that $S(\widetilde{v})$ is removable for $\widetilde{v}$ and so $\widetilde{v}$ is harmonic on all of $\R^n$. As it is homogeneous of degree $1$, it then follows that it must be linear, i.e.~$v^1,v^2$ are both linear functions, whose graphs are planes $P_1$ and $P_2$, say. As $\widetilde{v}(0) = 0$ from the homogeneity and continuity, it must be that $P_1\cap P_2 \supset \{0\}$, and in fact we must have $P_1\cap _2 = S(\widetilde{v})$. Therefore, $\mathbf{v}(v)\in \mathcal{P}_{d}\subset \CC$. Thus, we have shown by induction that if $v\in \mathcal{H}_d$ for $d\leq n-2$, then $\mathbf{v}(v)\in \mathcal{P}_{d}$, which completes the proof of the theorem.
\end{proof}

Just as for fine blow-ups, having classified the homogeneous degree one elements of the coarse blow-up class $\mathfrak{B}_\beta$, we now deduce generalised-$C^{1,\mu}$ regularity for coarse blow-ups.

\begin{theorem}[Regularity of Coarse Blow-Ups]\label{thm:coarse-regularity}
	Fix $\beta\in (0,1)$. Then, there exists $\mu = \mu(n,k,\beta)\in (0,1)$ such that $\mathfrak{B}_\beta\subset GC^{1,\mu}(B_{1/2}(0);\A_2(\R^k))$. Moreover, if $v\in \mathfrak{B}_\beta$ and $0\in \Gamma_v$, then there exists $\phi:\R^n\to \A_2(\R^k)$ with $\graph(\phi)\in \CC$ or $\phi = 2\llbracket L\rrbracket$ for a linear function $L:\R^n\to \R^k$ such that for every $0<\sigma\leq\rho/2\leq 3/16$, we have
	$$\sigma^{-n-2}\int_{B_\sigma(0)}\G(v(x)-v_a(0),\phi(x))^2\, \ext x \leq C\left(\frac{\sigma}{\rho}\right)^{2\mu}\cdot\rho^{-n-2}\int_{B_\rho(0)}|v|^2$$
	where $C = C(n,k,\beta)\in (0,\infty)$.
\end{theorem}

\begin{proof}
	 In view of property $(\mathfrak{B}3\text{I})$ and the fact that $\|v\|_{L^2(B_1(0))}\leq 1$ for each $v\in\mathfrak{B}_\beta$, it suffices to prove the claimed estimate for $\rho=3/8$, so assume $\rho=3/8$. First note that we can repeat the first part of the argument of Theorem \ref{thm:classification}, leading to the dichotomy that (a) or (b) holds, to establish the existence of $\eps = \eps(n,k,\beta)\in (0,1/2)$ such that the same dichotomy ((a) or (b)) must hold with $z=0$ for any $v\in\mathfrak{B}_\beta$ such that $0\in \Gamma_v$, provided we additionally assume that $v_a(0) = 0$ and $Dv_a(0)=0$. More precisely, there is $\eps = \eps(n,k,\beta)\in (0,1/2)$ such that for each $v\in \mathfrak{B}_\beta$ with $0\in \Gamma_v$, $v_a(0) = 0$, and $Dv_a(0) = 0$, either:
	 \begin{enumerate}
	 	\item [(a)] The conclusions of $(\mathfrak{B}10)$ hold on $B_{9/64}(0)$; in particular, $v$ is a generalised-$C^{1,\mu}$ graph on $B_{9/64}(0)$, and there is a function $\psi_0:\R^n\to \A_2(\R^k)$ with $\graph(\psi_0)\in \CC$ such that
		\begin{equation}\label{E:coarse-reg-1}
			\rho^{-n-2}\int_{B_\rho(0)}\G(v(x),\psi_0(x))^2\, \ext x \leq C\rho^{2\mu}\int_{B_{3/8}(0)}|v|^2
		\end{equation}
		for all $0<\rho\leq 9/64$, where $C = C(n,k,\beta)\in (0,\infty)$; or
		\item [(b)] We have
		\begin{equation}\label{E:coarse-reg-2}
			\sum_{\alpha=1,2}\int_{B_{3/8}(0)\setminus B_{3/16}(0)}R^{2-n}\left|\frac{\del}{\del R}\left(\frac{v^\alpha}{R}\right)\right|^2\geq \eps\int_{B_{3/8}(0)}|v|^2,
		\end{equation}
		where $R = |x|$.
	 \end{enumerate}
	 To prove this, we argue by contradiction exactly as in the proof of Theorem \ref{thm:classification}, but with $\|v_i\|^{-1}_{L^2(B_1(0))}v_i$ taking the place of $w_i = v_{z_i,\rho_i}$ appearing in the argument, where $(v_i)_i\subset\mathfrak{B}_\beta$ is a general sequence such that both options (a) and (b) with $v = v_i$ and $\eps = \eps_i$ are assumed to fail, with $\eps_i \to 0$; note that this argument utilises the classification of homogeneous degree $1$ elements provided by Theorem \ref{thm:classification} to reach a contradiction. (To contradict the case where the limit is a plane, note that we know the averages converge locally uniformly to the average of the limit, but harmonic estimates from $(\mathfrak{B}2)$ allow us to upgrade this convergence to $C^2$ convergence, so the limit would need to be zero, providing the contradiction.)
	 
	 Still subject to the conditions $v_a(0) = 0$, $Dv_a(0)=0$, this then leads to the final dichotomy as in the proof of Theorem \ref{thm:classification} (and by the same argument), namely that either the statement (I) holds with $z=0$ (and with the estimate $\rho^{-n-2}\int_{B_\rho(0)}\G(v(x),\psi_0(x))^2\, \ext x \leq C\rho^{2\gamma}\int_{B_{3/8}(0)}|v|^2$ for all $\rho\in (0,9/64)$), or the statement (II)$^\prime$ holds (with the estimate $\sigma^{-n-2}\int_{B_\sigma(0)}|v|^2\leq C\left(\frac{\sigma}{\rho}\right)^{2\gamma}\rho^{-n-2}\int_{B_\rho(0)}|v|^2$ for all $0<\sigma\leq\rho/2\leq 3/32$), where the constants $C,\gamma$ both depend only on $n$, $k$, $\beta$. In either case, this provides the desired estimate in the present theorem in the special case $v_a(0) = 0$ and $Dv_a(0) = 0$. The claimed estimate in the general case (i.e.~without the assumption $v_a(0) = 0$ and $Dv_a(0) = 0$) follows immediately from this special case in view of the property $(\mathfrak{B}3\text{III})$ and a standard derivative estimate for harmonic functions (applied to $v_a$).
	 
	 Finally, to see the claim that $v$ is of class $GC^{1,\mu}$ in $B_{1/2}$, note that we can apply, for any $z\in \Gamma_v\cap B_{1/2}$, the estimate just proved with $v_{z,1/2}$ in place of $v$ to obtain the corresponding estimate at the base point $z$ and with some $\phi_z:\R^n\to \A_2(\R^k)$ in place of $\phi$ where $\graph(\phi_z)\in \CC$ or $\graph(\phi_z)$ a multiplicity two plane $L_z$. Then, setting $\CC_v = \{z\in\Gamma_v\in B_{1/2}:\graph(\phi_z)\in\CC\}$, $\mathcal{B}_v = \Gamma_v\cap B_{1/2}\setminus \CC_v$, and $\mathcal{R}_v = B_{1/2}\setminus\Gamma_v$, we can employ property $(\mathfrak{B}7)$, $(\mathfrak{B}10)$, standard pointwise estimates for harmonic functions, together with the decay estimate of the present theorem with a Campanato estimate to show that $v$ is $GC^{1,\mu}$ in $B_{1/2}(0)$. This completes the proof.
\end{proof}

\section{Proof of Main Theorem}

We now prove the main result, Theorem \ref{thm:main}. This is done by utilising the asymptotic decay estimates for coarse blow-ups provided by Theorem \ref{thm:coarse-regularity}, and the fine $\eps$-regularity theorem Theorem \ref{thm:fine-reg}.

\begin{proof}[Proof of Theorem \ref{thm:main}]
	We first claim the following excess decay-like lemma.
	
	\textbf{Claim:} \emph{There exists $\eps = \eps(n,k,\beta)\in (0,1)$ and $\theta = \theta(n,k,\beta)\in (0,1)$ such that if $P$ is a plane in $\R^{n+k}$ with $\dist_\H(P\cap C_1(0),P_0\cap C_1(0))<\eps$, and if $V\in\mathcal{V}_\beta$ is such that $\Theta_V(0)\geq 2$, $(\w_n2^n)^{-1}\|V\|(B^{n+k}_2(0))<5/2$, and $\hat{E}_{V,P}<\eps$, then either:
	\begin{enumerate}
		\item [(i)] there is a plane $\widetilde{P}$ with $\dist_\H(\widetilde{P}\cap C_1(0), P\cap C_1(0))<C\hat{E}_{V,P}$ and $\hat{E}_{V,\widetilde{P}}(B_\theta(0))\leq \frac{1}{2}\hat{E}_{V,P}$; or,
		\item [(ii)] there is a cone $\BC\in \CC$ with $\dist_\H(\BC\cap C_1(0),P\cap C_1(0))<C\hat{E}_{V,P}$ and $\hat{E}_{V,\BC}(B_\rho(0))\leq C\rho^{\mu}\hat{E}_{V,P}$ for all $\rho\in (0,\theta/8]$.
	\end{enumerate}
	Here, $C = C(n,k,\beta)\in (0,\infty)$ and $\mu = \mu(n,k,\beta)\in (0,1)$.}
	
	To prove this, we argue by contradiction. So suppose we have sequences of varifolds $(V_j)_j\subset\mathcal{V}_\beta$ and a sequence of planes $(P_j)_j$ with $P_j\to P_0$ in Hausdorff distance in $C_1(0)$, such that $\Theta_{V_j}(0)\geq 2$, $(\w_n 2^n)^{-1}\|V_j\|(B^{n+k}_2(0))<5/2$, and $\hat{E}_{V_j,P_j}\to 0$. It suffices to prove, with $\theta = \theta(n,k,\beta)\in (0,1)$, $\mu = \mu(n,k,\beta)\in (0,1)$, and $C = C(n,k,\beta)\in (0,\infty)$ to be chosen, that for infinitely many $j$ we have either:
	\begin{enumerate}
		\item [(I)] there is a plane $\widetilde{P}_j$ with $\dist_\H(\widetilde{P}_j\cap C_1, P_j\cap C_1)<C\hat{E}_{V_j,P_j}$ and $\hat{E}_{V_j,\widetilde{P}_j}(B_\theta(0))\leq\frac{1}{2}\hat{E}_{V_j,P_j}$; or,
		\item [(II)] there is $\BC_j\in \CC$ with $\dist_\H(\BC_j\cap C_1, P_j\cap C_1)<C\hat{E}_{V_j,P_j}$ and $\hat{E}_{V_j,\BC_j}(B_\theta(0))\leq C\rho^\mu\hat{E}_{V_j,P_j}$ for all $\rho\in (0,\theta/8]$.
	\end{enumerate}
	Let $\Gamma_j:\R^{n+k}\to \R^{n+k}$ be a rotation such that $\Gamma_j(P_j) = \{0\}^k\times\R^n$ and $\|\Gamma_j-\id_{\R^{n+k}}\|\to 0$. Let $\widetilde{V}_j := (\Gamma_j)_\#V_j$; clearly by Proposition \ref{prop:rotation-class}, up to rescaling a small amount (which we ignore for notational simplicity, as we can take the rescaling to $1$ as $j\to\infty$), we have $\widetilde{V}_j\in\mathfrak{B}_{\beta_j}$ for some $\beta_j\uparrow\beta$. Then $\hat{E}_j:=\hat{E}_{\widetilde{V}_j} \equiv \hat{E}_{V_j,P_j}\to 0$. Let $v\in\mathfrak{B}_\beta$ be a coarse blow-up of $(\widetilde{V}_j)_j$. Since $\Theta_{\widetilde{V}_j}(0)\geq 2$, we have that $0\in \Gamma_v$ and $v_a(0) = 0$ (indeed, this is because $|v(0)|=0$ as it is the limit of this point when blown-up). So by Theorem \ref{thm:coarse-regularity}, there is $\phi:\R^n\to \A_2(\R^k)$ with $\mathbf{v}(\phi)\in \CC$ or a multiplicity $2$ plane such that for every $\sigma\in (0,3/16]$,
	\begin{equation}\label{E:main-1}
		\sigma^{-n-2}\int_{B_\sigma(0)}\G(v(x),\phi(x))^2\, \ext x \leq \widetilde{C}\sigma^{2\mu}\int_{B_{3/8}(0)}|v|^2
	\end{equation}
	where $\widetilde{C} = \widetilde{C}(n,k,\beta)\in (0,\infty)$ and $\mu = \mu(n,k,\beta)\in (0,1)$. Since $\int_{B_1}|v|^2\leq 1$, \eqref{E:main-1} and the homogeneity of $\phi$ imply that
	\begin{equation}\label{E:main-2}
		\int_{B_1(0)}|\phi|^2\leq \widetilde{C}_1
	\end{equation}
	for some $\widetilde{C}_1 = \widetilde{C}_1(n,k,\beta)\in (0,\infty)$. Also, note that \eqref{E:main-1} implies that $\sigma^{-n-2}\int_{B_\sigma(0)}|v_a(x)-\phi_a(x)|^2\, \ext x \leq \frac{1}{2}\widetilde{C}_1\sigma^{2\mu}\int_{B_{3/8}(0)}|v|^2$ for all $\sigma\in (0,3/16]$, and so since $v_a$ is harmonic, $v_a(0) = 0$, and $\phi$ is homogeneous of degree $1$, it follows that
	\begin{equation}\label{E:main-3}
		\phi_a(x) = Dv_a(0)\cdot x \qquad \text{for }x\in \R^n.
	\end{equation}
	Now choose $\theta = \theta(n,k,\beta)\in (0,1)$ such that
	\begin{equation}\label{E:main-4}
		\max\{\widetilde{C}_1(2\theta)^{\mu},(2\theta)^{\mu}\}<\min\{3/16,\eps\}
	\end{equation}
	where $\eps = \eps(n,k,\beta)\in (0,1)$ is as in $(\mathfrak{B}10)$ (where here $\mu$ is now the minimum power from the above with that in $(\mathfrak{B}10)$). We now divide into two cases:
	
	\emph{Case 1:} $(2\theta)^{-n-2}\int_{B_{2\theta}(0)}|v(x)-\phi_a(x)|^2\, \ext x <(2\theta)^\mu$. In this case, set $\widetilde{P}_j:= \Gamma_j^{-1}(\graph(\hat{E}_j\phi_a))$. Noting that by \eqref{E:main-3} that $\phi_a$ is linear, it is then straightforward to verify that option (I) holds for infinitely many $j$, with $C = C(n,k,\beta)$.
	
	\emph{Case 2: Case 1 fails.} In this case we have
	\begin{equation}\label{E:main-5}
		(2\theta)^{-n-2}\int_{B_{2\theta}(0)}|v(x)-\phi_a(x)|^2\, \ext x \geq (2\theta)^\mu.
	\end{equation}
	It follows that we must also have that $\mathbf{v}(\phi)\in \CC$. Indeed, if not then $\mathbf{v}(\phi) = 2|L|$, where $L = \graph(\phi_a)$, and so $\phi^\alpha(x) = \phi_a(x)$ for every $x\in \R^n$ and $\alpha\in \{1,2\}$. Thus, \eqref{E:main-1} would give $\sigma^{-n-2}\int_{B_\sigma(0)}|v(x)-\phi_a(x)|^2\, \ext x \leq \widetilde{C}_1\sigma^{2\mu}\int_{B_{3/8}(0)}|v|^2$ for all $\sigma\in (0,3/16]$, and thus taking $\sigma = 2\theta$ and using \eqref{E:main-4} we would see that Case 1 must, contrary to our assumption. So, we know that $\mathbf{v}(\phi)\in \CC$.
	
	Next note that there is a sequence of rotations $\widetilde{\Gamma}_j:\R^{n+k}\to \R^{n+k}$ with $\widetilde{\Gamma}_j(\graph(\hat{E}_j\phi_a)) = P_0$ and $\|\widetilde{\Gamma}_j-\id_{\R^{n+k}}\|\to 0$ and such that
	$$\widetilde{v}(x):= \frac{v(2\theta x)-\phi_a(2\theta x)}{\|v(2\theta(\cdot))-\phi_a(2\theta(\cdot))\|_{L^2(B_1)}}$$
	is the coarse blow-up of a subsequence of $((\eta_{0,2\theta}\circ\widetilde{\Gamma}_j)_\#\widetilde{V}_j)_j$. Then we have by \eqref{E:main-1},
	$$\int_{B_1}\G(\widetilde{v},\widetilde{\phi})^2 < \widetilde{C}_1(2\theta)^{2\mu}\cdot (2\theta)^2\|v(2\theta(\cdot))-\phi_a(2\theta(\cdot))\|_{L^2(B_1)}^{-1} \leq \widetilde{C}_1 \leq \widetilde{C}_1(2\theta)^{\mu} < \eps$$
	where $\widetilde{\phi}(x):=\|v(2\theta(\cdot))-\phi_a(2\theta(\cdot))\|_{L^2(B_1)}^{-1}\phi_f(2\theta x)$. Since $\|\widetilde{v}\|_{L^2(B_1)}=1$ and $\widetilde{\phi}_a\equiv 0$ we see, by the choice of $\eps$ and $(\mathfrak{B}10)$, that for sufficiently large $j$ the hypothees of Theorem \ref{thm:fine-reg} are satisfied with $W_k = (\eta_{0,2\theta}\circ\widetilde{\Gamma}_j)_\#\widetilde{V}_j$ in place of $V$ and $\hat{\BC}_j = \mathbf{v}(\hat{E}_{W_j}\widetilde{\phi})$ in place of $\BC$. Hence by applying Theorem \ref{thm:fine-reg} and noting that by \eqref{E:main-5} and \eqref{E:main-2} we have that $\int_{B_1}|\widetilde{\phi}|^2 \leq (2\theta)^{-\mu}\widetilde{C}_1 \equiv C = C(n,k,\beta)$, we see that option (II) must hold for infinitely many $j$ (with $\BC_j = (\widetilde{\Gamma}_j\circ\Gamma_k)^{-1}_\#\widetilde{\BC}_j$, where $\widetilde{\BC}_j$ is the cone $\widetilde{\BC}_j$ provided by Theorem \ref{thm:fine-reg} when that theorem is applied with $W_k$ in place of $V$ and $\hat{\BC}_j$ in place of $\BC$). This establishes the claim asserted at the beginning of the proof.
	
	Now we can apply the claim iteratively to deduce that for any plane with $\dist(P\cap C_1(0), P_0\cap C_1(0))<\eps$, one of the following must hold:
	\begin{enumerate}
		\item [(i)$^\prime$] there is a sequence of planes $(P_j)_j$ with $P_1 = P$, $\dist_\H(P_{j+1}\cap C_1(0),P_j\cap C_1(0))<C\hat{E}_{(\eta_{0,\theta^j})_\#V,P_j}$, and $\hat{E}_{(\eta_{0,\theta^j})_\#V,P_{j+1}}\leq \frac{1}{2}\hat{E}_{(\eta_{0,\theta^{j-1}})_\#V,P_{j}}$ for all $j\geq 1$; or,
		\item [(ii)$^\prime$] there is an integer $I\geq 1$ and a finite sequence of planes $P_1 = P,P_2,\dotsc,P_I$ such that (i)$^\prime$ holds for $j=1,\dotsc,I-1$ (if $I\geq 2$) and there is $\BC\in \CC$ with $\dist_\H(\BC\cap C_1(0),P_I\cap C_1(0))<C\hat{E}_{(\eta_{0,\theta^{I-1}})_\#V,P_I}$ and
		$$(\rho\theta^{I-1})^{-n-2}\int_{C_{\rho\theta^{I-1}}(0)}\dist^2(x,\BC)\, \ext\|V\|(x) \leq C\rho^{2\mu}\hat{E}^2_{(\eta_{0,\theta^{I-1}})_\#V,P_{I}},$$
		for all $\rho\in (0,\theta/8]$.
	\end{enumerate}
	From these, it is standard to deduce that there are constants $\mu_* = \mu_*(n,k,\beta)\in (0,1)$ and $C = C(n,k,\beta)\in (0,\infty)$ such that for any plane $P$ with $\dist_\H(P\cap C_1(0),P_0\cap C_1(0))<\eps$ we have either:
	\begin{enumerate}
		\item [(A)] there is a (unique) plane $P_*$ with $\dist_\H(P_*\cap C_1(0),P\cap C_1(0))<C\hat{E}_{V,P}$ and $\hat{E}_{(\eta_{0,\rho})_\#V,P_*}\leq C\rho^{\mu_*}\hat{E}_{V,P}$ for all $\rho\in (0,\theta/8]$; or,
		\item [(B)] there is a (unique) cone $\BC_*\in \CC$ with $\dist_\H(\BC_*\cap C_1(0),P\cap C_1(0))<C\hat{E}_{V,P}$ and
		$$\rho^{-n-2}\int_{C_\rho(0)}\dist^2(x,\BC_*)\, \ext\|V\|(x) \leq C\rho^{2\mu_*}\hat{E}_{V,P}^2$$
		for all $\rho\in (0,\theta/8]$.
	\end{enumerate}
	Indeed, (A) holds if (i)$^\prime$ holds and (B) holds if (ii)$^\prime$ holds. In particular, $V$ has a unique tangent cone at $0$ which is $\BC_*$ if (B) holds, and $2|P_*|$ if (A) holds (indeed, the multiplicity must be at least $2$ as $\Theta_V(0)\geq 2$, but then must be at most $2$ due to the mass upper bound).
	
	To complete the proof, note that if the hypotheses of Theorem \ref{thm:main} are satisfied with $\eps = \eps(n,k,\beta)\in (0,1)$ sufficiently small, then for any $z\in C_{3/4}(0)$ with $\Theta_V(z)\geq 2$, we may repeat the argument leading to (A) and (B) with $(\eta_{z,1/4})_\#V$ in place of $V$. This gives that $\{z\in C_{3/4}(0):\Theta_v(z)\geq 2\} = \{z\in C_{3/4}(0):\Theta_V(z) = 2\}$. At this point, the proof is then completed in an analogous manner to that seen for Theorem \ref{thm:fine-reg}, and so we are done.
\end{proof}

We remark that Theorem \ref{thm:main-3} follows from our arguments above leading to Theorem \ref{thm:main}. Indeed, recall Remark \ref{remark:holder-continuity}: assuming that Theorem \ref{thm:main-3}(ii)--(iii) fail along a blow-up sequence (for some fixed constants $\beta$, $\gamma$, and with corresponding sequences $\eta_j\to 0$, $(\rho_0)_j\to 0$, we know that we can prove that such coarse blow-ups are Hölder continuous. Indeed, when we don't have the topological structural condition in a $\beta$-coarse gap, the claim in Theorem \ref{thm:main-3}(ii) regarding harmonic graphs still gives that both components of the coarse blow-up will coincide with a (single) harmonic function in this region, giving $(\mathfrak{B}7)$. However, the assumption that Theorem \ref{thm:main-3}(iii) fails is sufficient for us to apply Theorem \ref{thm:fine-reg} (see Remark \ref{remark:fine-reg-sufficiency}) in order to deduce that these coarse blow-ups satisfy all the properties in Theorem \ref{thm:further-blow-up-properties}. The proof of the generalised-$C^{1,\mu}$ regularity of the coarse blow-ups then goes through as above, leading to the claimed decayed statement in Theorem \ref{thm:main-3}(i) holding along some subsequence of the contradiction sequence, which completes the proof of Theorem \ref{thm:main-3}.

\part{\centering Additional Observations and Extensions}\label{part:extension}

In this final part of the paper, we wish to prove some extensions of Theorem \ref{thm:main}. The first is an a posteriori technical improvement to the coarse blow-up procedure: we can in fact show that we have have energy convergence locally to the coarse blow-up, thus improving our previous convergence (which was strong in $L^2_{\text{loc}}(B_1(0))$ and weak in $W^{1,2}_{\text{loc}}(B_1(0))$) to strong convergence in $W^{1,2}_{\text{loc}}(B_1(0))$. Moreover, we in fact have that the blow-ups convergence in $C^{0,\alpha}_{\text{loc}}(B_1(0))$ to the coarse blow-up for any $\alpha\in (0,1)$, providing also local uniform convergence. If we knew that there were now twisted classical singularities (which is true for example in codimension one) then the conclusion of Theorem \ref{thm:main} would in fact be that the graph is $C^{1,\alpha}$ in the multi-valued sense, and so one could apply Arzelà--Ascoli to in fact get local uniform convergence in $C^{1,\beta}$ for any $\beta<\alpha$, and so in particular the derivatives also converge locally uniformly to the blow-up. These observations are needed in order to utilise the planar frequency function of B.~Krummel and the third author \cite{KW23a} in the context of Theorem \ref{thm:main} in order to establish control on the size of the singular set, which we intend to do in later work.

Secondly, we give an extension of Theorem \ref{thm:main} to integral varifolds which have non-vanishing mean curvature lying in $L^p$ for some $p>n$. Then we end with providing an extension of these results from the present setting, where the varifolds lie in an ambient Euclidean space, to one where they lie in a Riemannian manifold.

\section{Energy Convergence to Coarse Blow-Ups}

We wish here to give a technical improvement to the conclusions of Theorem \ref{thm:main} showing that, a posteriori, we in fact have local energy convergence in the coarse blow-up procedure. The reader should recall that this was the key issue which was preventing us from passing information from the stationarity condition at the varifold level to the coarse blow-up. This energy convergence will be key for showing approximate monotonicity of the planar frequency function in the forthcoming work \cite{KMW25}.

\begin{theorem}
	Fix $\beta\in (0,1)$ and $\eta>0$. Then, there exists $\eps = \eps(n,k,\beta,\eta)\in (0,1)$ such that the following holds. Suppose $V\in \mathcal{V}_\beta$ obeys:
	\begin{itemize}
		\item $(\w_n 2^n)^{-1}\|V\|(B_2^{n+k}(0)) < \frac{5}{2}$;
		\item $\frac{3}{2} \leq \w_n^{-1}\|V\|(C_1(0))< \frac{5}{2}$;
		\item $\hat{E}_V<\eps$.
	\end{itemize}
	Then, the conclusions of Theorem \ref{thm:main} hold and, if $f:B^n_{1/2}(0)\to \A_2(\R^k)$ is the function given in Theorem \ref{thm:main}, then there is a coarse blow-up $v\in \mathfrak{B}_\beta$ such that
	$$\sum_{\kappa=1}^k\sum_{\alpha=1,2}\int_{B_{1/2}^n(0)}|\hat{E}_V^{-1}Df^{\kappa,\alpha}-Dv^{\kappa,\alpha}|^2 \leq \eta.$$
	In particular, the convergence to the coarse blow-up is strongly in $W^{1,2}_{\textnormal{loc}}(B_1(0))$.
\end{theorem}

\textbf{Remark:} Note that as a consequence of Theorem \ref{thm:coarse-regularity}, we know that any function in $\mathfrak{B}_\beta$ is $GC^{1,\alpha}$ and so admits a monotone frequency function by Theorem \ref{thm:frequency}.

\textbf{Remark:} The energy convergence would be a simple consequence of Arzelà--Ascoli (for $C^{1,\alpha}$ functions) and Theorem \ref{thm:main} if the functions where $C^{1,\alpha}$ instead of generalised-$C^{1,\alpha}$. We therefore give an argument which avoids the need for Arzelà--Ascoli for the derivatives and only uses a $C^{0,\alpha}$ version.

\begin{proof}
	Of course provided $\eps = \eps(n,k,\beta)\in (0,1)$ is sufficiently small the conclusions of Theorem \ref{thm:main} hold, and so we are only interested in the second claim. For this, we will argue by contradiction. If it is false for a given $\beta$ and $\eta$, we may therefore take a sequence of varifolds $(V_j)_j\subset \mathcal{V}_\beta$ converging as varifolds to the plane $P_0$ with multiplicity two in $C_1(0)$. Let $f_j$ be the functions representing $V_j$ over $C_{\sigma_j}(0)$ where $\sigma_j\uparrow 1$, and take the coarse blow-up, i.e.~we have $v_j:= \hat{E}_{V_j}^{-1}f_j \to v\in\mathfrak{B}_\beta$ where the convergence is strongly in $L^2(B_{7/8}(0))$ and weakly in $W^{1,2}(B_{7/8}(0))$ (in the sense described in Section \ref{sec:coarse-blow-up-construction}). Moreover, by Theorem \ref{thm:main}, we have a uniform bound on the $C^{0,1}$ norm of the multi-valued functions $v_j$, and thus by Arzelà--Ascoli we can also assume that the convergence is also uniform in $B_{7/8}(0)$.
	
	It is then standard that weak convergence in $W^{1,2}$ gives for each $\kappa,\alpha$ that
	$$\|Dv^{\kappa,\alpha}\|_{L^2(B_{1/2}^n(0)} \leq \liminf_{j\to\infty}\|Dv_j^{\kappa,\alpha}\|_{L^2(B_{1/2}^n(0))}.$$
	Thus, to establish the claimed strong convergence it suffices to show for each $\kappa\in \{1,\dotsc,k\}$,
	\begin{equation}\label{E:energy-conv-1}
		\sum_{\alpha=1,2}\limsup_{j\to\infty}\|Dv_j^{\kappa,\alpha}\|_{L^2(B_{1/2}^n(0))} \leq \sum_{\alpha=1,2}\|Dv^{\kappa,\alpha}\|_{L^2(B_{1/2}^n(0))}.
	\end{equation}
	For this, we will use an improved version of the energy non-concentration estimate $(\mathfrak{B}6)$ but \emph{before} passing to a coarse blow-up; this improved estimate is possible due to the regularity conclusion provided by Theorem \ref{thm:main}, in particular the uniform convergence to the coarse blow-up. Indeed, if one analyses the proof of $(\mathfrak{B}6)$ from Theorem \ref{thm:blow-up-properties}, in light of Theorem \ref{thm:main} there are significant simplifications. Indeed, in \eqref{E:B6-5} there is now no `bad set' on the integral on the left-hand side, and so the set $\Omega_j$ therein can always be taken to be $B_{3/4}^n(0)$. Moreover, the Lipschitz constant is now $O(\hat{E}_j)$, and thus any terms involving the Lipschitz constant be absorbed into other error terms. Furthermore, in \eqref{E:B6-12} there is no need to introduce the set $U$, as the uniform convergence is now on all of $B_{1/2}^n(0)$. Thus, combining \eqref{E:B6-11} and \eqref{E:B6-12} now gives the following:
	
	\textbf{Claim:} \emph{Fix $\beta\in (0,1)$ and $\bar{\delta}>0$. Then, there exists $\eps = \eps(n,k,\beta,\bar{\delta})\in (0,1)$ such that under the assumptions of the current lemma with this choice of $\eps$, we have for any $\delta>\bar{\delta}$ and $\kappa\in \{1,\dotsc,k\}$,
	$$\sum_{\alpha=1,2}\int_{B^n_{1/2}(0)}\one_{\{|v^\alpha-v_a|<\delta\}}\cdot |\hat{E}^{-1}_{V}Df^{\kappa,\alpha} - Dv_a^\kappa|^2 \leq C\delta$$
	where $C = C(n,k)\in (0,\infty)$ and $f,v$ are as in the conclusions of the present lemma.}

	Fix now $\kappa\in \{1,\dotsc,k\}$, and also fix $\eps>0$. We may then choose $\delta = \delta(n,k,\eps)\in (0,1)$ in the above claim so that $C\delta<\eps$, we see that for all sufficiently large $j$ we must have
	$$\sum_{\alpha=1,2}\int_{B_{1/2}^n(0)}\one_{\{|v^\alpha-v_a|<\delta\}}\cdot|Dv_j^{\kappa,\alpha}-Dv_a^\kappa|^2 <\eps.$$
	From the uniform convergence of $v_j$ to $v$, we know that for all sufficiently large $j$ (depending on $j$) we have
	$$B_{1/2}(0)\cap \{|v^\alpha-v_a|\geq\delta\}\subset B_{1/2}(0)\cap\{|v^\alpha_j-(v_j)_a|\geq 3\delta/4\}.$$
	This latter conclusion tells us that on the set $\R^k\times (B_{1/2}(0)\cap \{|v^\alpha_j-v_a|\geq\delta\})$, $V_j$ can no have points of density $\geq 2$ for all $j$ sufficiently large. Thus, by Lemma \ref{lemma:gap-2} (and the uniform bounds given by Theorem \ref{thm:main}), for all sufficiently large $j$ we know that $V_j$ splits as a sum of two smooth minimal graphs here, and thus $f_j$ has $C^2$ estimates in this region which guarantees that $Dv_j\to Dv$ strongly in $L^2$ (in fact, uniformly) here. Thus, we get
	\begin{align*}
		\sum_{\alpha=1,2}\int_{B_{1/2}^n(0)}|Dv_j^{\kappa,\alpha}-Dv^\kappa_a|^2 & = \sum_{\alpha=1,2}\int_{B_{1/2}^n(0)}\one_{\{|v^\alpha-v_a|<\delta\}}|Dv_j^{\kappa,\alpha}-Dv^\kappa_a|^2\\
		& \hspace{8em} + \sum_{\alpha=1,2}\int_{B_{1/2}^n(0)}\one_{\{|v^\alpha-v_a|\geq\delta\}}|Dv^{\kappa,\alpha}_j-Dv^\kappa_a|^2\\
		& \leq \eps + \eps + \sum_{\alpha=1,2}\int_{B_{1/2}^n(0)\cap \{|v^\alpha-v_a|\geq\delta\}}|Dv^{\kappa,\alpha}-Dv^\kappa_a|^2\\
		& \leq 2\eps + \sum_{\alpha=1,2}\int_{B_{1/2}^n(0)}|Dv^{\kappa,\alpha}-Dv^\kappa_a|^2
	\end{align*}
	for all $j$ sufficiently large. Hence, taking $\limsup_{j\to\infty}$ and then $\eps\downarrow 0$, we see that
	$$\limsup_{j\to\infty}\sum_{\alpha=1,2}\int_{B_{1/2}^n(0)}|Dv_j^{\kappa,\alpha}-Dv^\kappa_a|^2 \leq \sum_{\alpha=1,2}\int_{B_{1/2}^n(0)}|Dv^{\kappa,\alpha}-Dv^\kappa_a|^2.$$
	Finally, notice that since $v^{\kappa,\alpha}_j\to v^{\kappa,\alpha}$ weakly in $W^{1,2}(B^n_{1/2}(0))$, we have $\int_{B_{1/2}^n(0)}Dv_j^{\kappa,\alpha}\cdot Dv^\kappa_a \to \int_{B_{1/2}^n(0)}Dv^{\kappa,\alpha}\cdot Dv_a^\kappa$, and thus the above is actually equivalent to
    $$\limsup_{j\to\infty}\sum_{\alpha=1,2}\int_{B^n_{1/2}(0)}|Dv_j^{\kappa,\alpha}|^2 \leq \sum_{\alpha=1,2}\int_{B_{1/2}^n(0)}|Dv^{\kappa,\alpha}|^2.$$
	This gives \eqref{E:energy-conv-1}, thus proving the theorem.
\end{proof}

\section{Generalisation to Non-Zero Mean Curvature and Riemannian Manifolds}\label{sec:generalisation}

In this section we will extend Theorem \ref{thm:main} first to the setting where the (generalised) mean curvature of the varifold lies in $L^p$, for some $p>n$ (in particular, it need not be zero), and then we subsequently extend this result to the case of more general ambient Riemannian manifolds.

\subsection{Non-Zero Generalised Mean Curvature}

Recall that the \emph{generalised} \emph{mean curvature} $H$ of an integer rectifiable $n$-varifold $V$ on an open set $U\subset\R^{n+k}$ is determined by the identity
$$\int \div_{T_xV}(\Phi)\, \ext\|V\|(x) = -\int \Phi\cdot H\, \ext\|V\|(x) \qquad \text{for all }\Phi\in C^1_c(U;\R^{n+k}).$$
In particular, $H$ is an $\R^{n+k}$-valued function defined $\H^n$-a.e. on $\spt\|V\|$, and the case that $V$ is stationary is equivalent to $H\equiv 0$.

Suppose now that $V$ is an integer rectifiable $n$-varifold in $B^{n+k}_2(0)\subset \R^{n+k}$. Allard's regularity theorem \cite{All72} (see also \cite{Sim83}) applies to $V$ provided the mean curvature of $V$ lies in $L^p(B_2^{n+k}(0))$ for some $p>n$. Let us briefly explain why this is true in light of the stationary case. Allard's regularity theorem assumes that, for some $\eps>0$,
$$\frac{1}{2}\leq (\w_n 2^n)^{-1}\|V\|(B^{n+k}_2(0))\leq \frac{3}{2},\qquad \hat{E}_V<\eps,\qquad \text{and} \qquad \left(\int_{B_1(0)}|H|^p\, \ext\|V\|\right)^{1/p} < \eps.$$
One then defines a modified \emph{excess} on a ball $B_\rho(x)\subset B^{n+k}_1(0)$ by:
$$\hat{E}^*_{V,\eps}(B_\rho(x)) := \max\left\{\hat{E}_V(B_\rho(x)),\, \eps^{-1/2}\left(\rho^{p-n}\int_{B_\rho(x)}|H|^p\, \ext\|V\|\right)^{1/p}\right\},$$
where as a shorthand $\hat{E}_{V,\eps}^*\equiv \hat{E}_{V,\eps}^*(B_1(0))$. Notice that for $p>n$, the $L^p$ norm of $H$ \emph{automatically} decays at a rate of $\rho^{1-\frac{n}{p}}$; indeed,
$$\left(\rho^{p-n}\int_{B_\rho(0)}|H|^p\, \ext\|V\|\right)^{1/p} \leq \rho^{1-\frac{n}{p}}\left(\int_{B_1(0)}|H|^p\, \ext\|V\|\right)^{1/p}.$$
The entire proof of Allard's regularity theorem in the stationary setting then passes through utilising the excess $\hat{E}^*_{V,\eps}$ in place of the usual height excess $\hat{E}_{V}$, provided one can control the additional error terms arising from the mean curvature term in the first variation formula by $\hat{E}_{V,\eps}^*$. Thus is done by using Hölder's inequality, which for example gives for $\Phi\in C^1_c(U;\R^{n+k})$ with $\spt(\Phi)\subset B_\rho(x)$,
\begin{align*}
\left|\int\Phi\cdot H\, \ext\|V\|\right| \leq \sup|\Phi|\int_{B_\rho(x)}|H|\, \ext\|V\| & \leq \sup|\Phi|\left(\int_{B_\rho(x)}|H|^p\, \ext\|V\|\right)^{1/p}\cdot\|V\|(B_\rho(x))^{1-\frac{1}{p}}\\
& \leq C\sup|\Phi|\cdot \rho^{n-1}\cdot\left(\rho^{p-n}\int_{B_\rho(x)}|H|^p\, \ext\|V\|\right)^{1/p}\\
& \leq C\sup|\Phi|\rho^{n-1}\cdot\hat{E}_{V,\eps}^*(B_\rho(x))\cdot \eps^{1/2}.
\end{align*}
In particular, as the left-hand side of this expression scales like $\rho^{n-1}$, this tells us that such terms are $o(\hat{E}_{V,\eps}^*(B_\rho(x)))$ during the blow-up procedure, as we will be sending $\eps\to 0$ along the blow-up sequence. The result of the above computation and similar ones (depending on the choice of $\Phi$) means that the additional terms arising from the generalised mean curvature can be incorporated into the error terms involving now the modified excess $\hat{E}_{V,\eps}^*$, and so they both do not impact the blow-up process as well as vanish in the blow-up limit. Moreover, the coarse blow-up (which in the proof of Allard regularity is single-valued) is still a harmonic function, by the same argument as in $(\mathfrak{B}2)$. One therefore still proves the analogous form of excess decay here, now for $\hat{E}_{V,\eps}^*$ (and with $\eps$ now fixed depending on $n,k,p$), and deduces the desired regularity conclusion from the Lipschitz approximation. The optimal regularity of the graph is then determined by the decay rate of $\hat{E}_{V,\eps}^*$, which from the incorporation of the generalised mean curvature term has a decay rate in general of at most $1-\frac{n}{p}\in (0,1)$.

Let us now describe the generalisation of Theorem \ref{thm:main} to the setting of generalised mean curvature in $L^p$, $p>n$.
Fix parameters $\beta\in (0,\infty)$ and $p>n$. Write $\mathcal{V}_{\beta,p}$ for the set of all integral $n$-rectifiable varifolds $V$ in $B_2(0)$ with generalised mean curvature $H$ in $L^p(B_2(0))$ which obey the following condition: whenever $x\in P_0$ and $\rho>0$ are such that $C_\rho(x)\subset C_1(0)$ and:
\begin{enumerate}
	\item [(a)] $C_\rho(x)\subset\{\Theta_V<2\}$;
	\item [(b)] $\frac{3}{2}\leq (\w_n\rho^n)^{-1}\|V\|(C_\rho(x))\leq \frac{5}{2}$;
	\item [(c)] For some $z\in P_0^\perp\equiv \R^k\times \{0\}^n$,
	$$\left(\frac{1}{\rho^{n+2}}\int_{C_\rho(x)}\dist^2(y,\{z\}\times\R^n)\, \ext\|V\|(y)\right)^{1/2} + \left(\rho^{n-p}\int_{C_\rho(x)}|H|^p\, \ext\|V\|\right)^{1/p}<\beta;$$
\end{enumerate}
then $\spt\|V\|\cap C_{\rho/2}(x)$ is disconnected, and at least two of its connected components have non-empty intersection with $C_{\rho/4}(x)$.

Theorem \ref{thm:main} then has the following generalisation to the class $\mathcal{V}_{\beta,p}$.
\begin{theorem}\label{thm:main-gen-mc}
	Fix $\beta\in (0,1)$ and $p>n$. Then, there exists $\eps_0 = \eps_0(n,k,\beta,p)\in (0,1)$ such that the following is true. Suppose that $V\in \mathcal{V}_{\beta,p}$ satisfies:
	$$\frac{3}{2}\leq\frac{\|V\|(B_1(0))}{\w_n}\leq \frac{5}{2} \qquad \text{and} \qquad \hat{E}_V + \|H\|_{L^p(C_1(0))}<\eps_0.$$
	Then, there is a $GC^{1,\alpha}$ function $f:B^n_{1/2}(0)\to \A_2(\R^k)$ with
	$$\widetilde{V}\res C_{1/2}(0) = \mathbf{v}(f)$$
	where $\widetilde{V} := V\res B^{n+k}_{3/2}(0)$. Moreover, conclusions (1) -- (4) of Theorem \ref{thm:main} hold, with $\hat{E}_V$ therein replaced by $\hat{E}_V + \|H\|^{1/2}_{L^p(B_1(0))}$, where the constants therein obey $C = C(n,k,p)\in (0,\infty)$ and $\alpha = \alpha(n,k,p)\in (0,\frac{1}{2}(1-\frac{n}{p}))$.
\end{theorem}

Since our original proof of Theorem \ref{thm:main} utilises only first variation arguments, the same proof goes through, making the same modifications as described above, to prove Theorem \ref{thm:main-gen-mc}. The reader will note that the estimates concluded in Theorem \ref{thm:main-gen-mc} are in terms of $\hat{E}_V + \|H\|^{1/2}_{L^p(B_1(0))}$, rather than $\hat{E}_V + \|H\|_{L^p(B_{1}(0))}$. The need for an additional square root on the mean curvature term is explained by the fact that our argument is more involved than that in Allard's regularity theorem described above, in particular through the need to use more involved first variation arguments. For instance, the proof of the Hardt--Simon inequality for coarse blow-ups (Theorem \ref{thm:blow-up-properties}(5)), as well as the proof of the squeeze identity, use first variation arguments which require dividing by the \emph{square} of the excess, however the extra terms coming from the mean curvature are only \emph{linear} in $\|H\|_{L^p}$. This is resolved by taking in the definition of the excess the \emph{square root} of $\|H\|_{L^p}$, which still naturally decays (at a rate now of $\frac{1}{2}(1-\frac{n}{p})$) and leads to the bound above.

It should also be noted that one also must modify the definition of fine excess in Part \ref{part:fine-reg} to also incorporate the generalised mean curvature in an analogous way to $\hat{E}^*_V$. Again, this does not cause any issues due to it decaying at a fixed rate and that the mean curvature terms disappear in the fine blow-up.

\subsection{Generalisation to Riemannian Manifolds}

In this section we note how to extend Theorem \ref{thm:main-gen-mc} (and thus Theorem \ref{thm:main}) to smooth ambient Riemannian manifolds $(M^m,g)$. To begin with, recall that by the Nash embedding theorem, we can without loss of generality assume that $M$ is (isometrically) embedded in some Euclidean space $\R^{n+k}$, which without loss of generality passes through the origin and with the tangent space to $M$ at $0$ being $\R^m\cong \R^m\times\{0\}^{n+k-m}$.

By scaling, we may then work under the assumption that the second fundamental form $A_M$ of $M$ in the unit ball is small, i.e.~$\|A_M\|_{L^\infty(B^{n+k}_1(0))} \leq \eps$, for some small $\eps\in (0,1)$. In particular, the injectivity radius of $M$ is $>1$. Write $B^M_1(0)$ for the geodesic ball in $M$ of radius $1$ centred at $0$. We then write $\mathcal{V}^M_{\beta,p}$ for the set of all integral $n$-rectifiable varifolds in $M$ whose generalised mean curvature $H$ within $M$ is in $L^p(B^M_1(0))$, $p>n$, and which obey the analogous structural condition as in the definition of $\mathcal{V}_{\beta,p}$ in the previous section. The generalisation of Theorem \ref{thm:main-gen-mc} is then the following:

\begin{theorem}\label{thm:main-gen-manifold}
	Fix $\beta\in (0,1)$ and $p>n$, and suppose $(M^m,g)$ is a Riemannian manifold isometrically embedded in some Euclidean space $\R^{n+k}$ and whose second fundamental form $A_M$ obeys $|A_M|<1$ in $B_{1}^{n+k}(0)$. Then, there exists $\eps_0 = \eps_0(n,k,\beta,p)\in (0,1)$ such that the following is true. Suppose that $V\in\mathcal{V}^M_{\beta,p}$ satisfies:
	$$\frac{3}{2}\leq\frac{\|V\|(B_1(0))}{\w_n}\leq\frac{5}{2} \qquad \text{and} \qquad \hat{E}_V + \|H\|_{L^p(B^M_1(0))} + \|A_M\|_{L^\infty(B_1^M(0))} < \eps_0.$$
	Then, there is a $GC^{1,\alpha}$ function $f:B^n_{1/2}(0)\to \A_2(\R^k)$ such that
	$$\widetilde{V} \res C_{1/2}(0) = \mathbf{v}(f),$$
	where $\widetilde{V} := V\res B_{3/2}^{n+k}(0)$. More precisely, conclusions (1) -- (4) of Theorem \ref{thm:main} hold, with $\hat{E}_V$ therein replaced by $\hat{E}_V + \|H\|^{1/2}_{L^p(B_1^M(0))} + \|A_M\|^{1/2}_{L^\infty(B_1^M(0))}$, and the constants satisfy $C = C(n,k,p)\in (0,\infty)$, $\alpha = \alpha(n,k,p)\in (0,\frac{1}{2}(1-\frac{n}{p}))$.
\end{theorem}

The proof of Theorem \ref{thm:main-gen-manifold} follows in an analogous manner to Theorem \ref{thm:main} and Theorem \ref{thm:main-gen-mc}, except now we have an additional error term in the first variation coming from the curvature of $M$. However, once again this has a natural decay rate (as the curvature scales as the radius), and so once again by appropriately modifying the excess quantities one can absorb any error terms into the new excess. As the additional term in the excess quantity has a natural decay rate, this does not impact the argument substantially and so one can deduce Theorem \ref{thm:main-gen-manifold}.

\section{Stationary Integral Varifolds Close to the Union of Two Planes or $4$ Half-Planes}\label{sec:non-degenerate}

We saw in Theorem \ref{thm:fine-reg} that one can establish an $\eps$-regularity theorem for stationary integral varifolds $V$ in the regime $\hat{E}_V\ll 1$, $\hat{F}_{V,\BC}\ll \hat{E}_V$ for some $\BC\in \CC$, \emph{provided} $V$ does not admit any $(\beta,\gamma)$-fine gaps (cf.~Remark \ref{remark:fine-reg-sufficiency}). In this section, we wish to record some analogous results when instead we assume that $\inf_{P}\hat{E}_{V,P}\geq \eps$ (the infimum taken over planes $P$) and $\hat{F}_{V,\BC}<\tilde{\eps}$, for suitably small $\tilde{\eps} = \tilde{\eps}(\eps,\beta,\gamma)$, i.e.~$V$ is close to a cone in $\CC$ but far from any plane. In fact, we shall allow our parameters to depend on the cone $\BC$. These results should be viewed as direct analogues of those in \cite{BK17}.

Thus, we fix $\BC\in \CC$ throughout this section. We then define gap regions adapted to the cone $\CC$ in an analogous way to how $\beta$-coarse gaps where defined relative to $P_0$.

\begin{defn}\label{defn:cone-gaps}
    Fix $\beta\in (0,1)$ and $V$ a stationary integral $n$-varifold in $B^{n+k}_2(0)$. We call a ball $B_\rho(x)\subset B^{n+k}_2(0)$ a $\beta$\emph{-coarse $\BC$-gap} (for $V$) if the following hold:
    \begin{enumerate}
        \item [(1)] $B_\rho(x)\subset\{\Theta_V<2\}$;
        \item [(2)] $\frac{3}{2} \leq (\w_n\rho^n)^{-1}\|V\|(B_\rho(x))\leq \frac{5}{2}$;
        \item [(3)] We have
        $$\frac{1}{\rho^{n+2}}\int_{B_\rho(x)}\dist^2(y,\BC)\, \ext\|V\| + \frac{1}{\rho^{n+2}}\int_{B_\rho(x)\setminus\{r_{\BC}<\rho/8\}}\dist^2(y,V)\, \ext\|\BC\| < \beta^2.$$
    \end{enumerate}
\end{defn}
Notice that for $\eps>0$, there is $\beta_* = \beta_*(\BC,\eps)\in (0,1)$ sufficiently small such that if $B_\rho(x)$ is a $\beta_*$-coarse $\BC$-gap, then the $L^2$ height excess of $(\eta_{x,\rho})_\#V$ in $B_1(0)$ to any plane is $>\eps$. As such, taking $\beta\in (0,\beta_*)$, we can guarantee that in a $\beta$-coarse $\BC$-gap $V$ is far from any plane. We will fix this choice of $\beta_*$.

\begin{defn}
    We say that $V$ satisfies the topological structural condition in a $\beta$-coarse $\BC$-gap $B_\rho(x)$ if $\spt\|V\|\cap B_{\rho/2}(x)$ has at least two connected components which intersect $B_{\rho/4}(x)$.
\end{defn}

The main $\eps$-regularity results can then be loosely described as follows, depending on the form that $\BC$ takes.

\begin{theorem}
    Fix $\BC\in \mathcal{P}_{\leq n-2}$ a cone passing through $0$ and $\beta\in (0,\beta_*)$. Then, there exists $\eps = \eps(\BC,\beta)\in (0,1)$ such that the following is true. Suppose that $V$ is a stationary integral $n$-varifold in $B^{n+k}_2(0)$ obeying:
    \begin{itemize}
        \item $\frac{3}{2}\leq \w_n^{-1}\|V\|(B_1(0))\leq \frac{5}{2}$;
        \item $\hat{F}_{V,\BC}<\eps$.
    \end{itemize}
    Then, provided $V$ has no $\beta$-coarse $\BC$-gaps which fail the topological structural condition, we have that $V\res B_{1/2}(0)$ is the union of two smooth minimal submanifolds, each of which lies close to one of the two planes in $\BC$, and which intersect only along a subset of a smooth $(n-2)$-dimensional submanifold.
\end{theorem}

\begin{theorem}
    Fix $\BC\in \CC_{n-1}\setminus\mathcal{P}$ a cone passing through $0$ and $\beta\in (0,\beta_*)$. Then, there exists $\eps = \eps(\BC,\beta)\in (0,1)$ such that the following is true. Suppose that $V$ is a stationary integral $n$-varifold in $B^{n+k}_2(0)$ obeying:
    \begin{itemize}
        \item $\frac{3}{2}\leq \w_n^{-1}\|V\|(B_1(0))\leq \frac{5}{2}$;
        \item $\hat{F}_{V,\BC}<\eps$.
    \end{itemize}
    Then, provided $V$ has no $\beta$-coarse $\BC$-gaps which fail the topological structural condition, we have that $V\res B_{1/2}(0)$ is the union of four smooth minimal submanifolds with boundary, meeting only along an $(n-1)$-dimensional $C^{1,\alpha}$ submanifold, which is their common boundary.
\end{theorem}

\begin{theorem}
    Fix $\BC\in \mathcal{P}_{n-1}$ a cone passing through $0$ and $\beta\in (0,\beta_*)$. Then, there exists $\eps = \eps(\BC,\beta)\in (0,1)$ such that the following is true. Suppose that $V$ is a stationary integral $n$-varifold in $B^{n+k}_2(0)$ obeying:
    \begin{itemize}
        \item $\frac{3}{2}\leq \w_n^{-1}\|V\|(B_1(0))\leq \frac{5}{2}$;
        \item $\hat{F}_{V,\BC}<\eps$.
    \end{itemize}
    Then, provided $V$ has no $\beta$-coarse $\BC$-gaps which fail the topological structural condition, then the singular set of $V\res B_{1/2}(0)$ is contained within an $(n-1)$-dimensional submanifold, and at each singular point there is a unique tangent cone, which belongs to $\CC$.
\end{theorem}

In turn, these results follow from a excess decay dichotomy, either giving decay of the excess to a nearby cone, or there being a gap of significant size.

\begin{theorem}\label{thm:cone-gap-or-decay}
    Fix $\BC\in \CC$ a cone passing through $0$ and $\beta\in (0,\beta_*)$. Then, there exists $\eps = \eps(\BC,\beta)\in (0,1)$ and $\rho_0 = \rho_0(\BC,\beta)\in (0,1)$ such that the following is true. Suppose $V$ is a stationary integral varifold in $B^{n+k}_2(0)$ which obeys:
    \begin{itemize}
        \item $\Theta_V(0)\geq 2$;
        \item $\frac{3}{2}\leq \w_n^{-1}\|V\|(B_1(0))\leq \frac{5}{2}$;
        \item $\hat{F}_{V,\BC}<\eps$.
    \end{itemize}
    Then, one of the following alternatives occurs:
    \begin{enumerate}
        \item [(i)] there is a cone $\BC^\prime\in \CC$ passing through the origin with $\dist_\H(\BC\cap B_1,\BC^\prime\cap B_1)\leq C\hat{E}_{V,\BC}$ such that for some $\theta = \theta(n,k,\beta)\in (0,1/4)$, $\hat{F}_{V,\BC^\prime}(B_\theta) < \frac{1}{2}\hat{E}_{V,\BC}$.
        \item [(ii)] There is a $\beta$-coarse $\BC$-gap for $V$ with radius $\rho\geq \rho_0$ which fails the topological structural condition.
    \end{enumerate}
    Here, $C = C(n,k,\beta)\in (0,\infty)$.
\end{theorem}

\bibliographystyle{alpha} 
\bibliography{references}

\end{document}